Ilgar Giyas oglu Aliyev

# Technological foundations of management decision-making in the reconstruction of complex gas pipeline system

MONOGRAPH

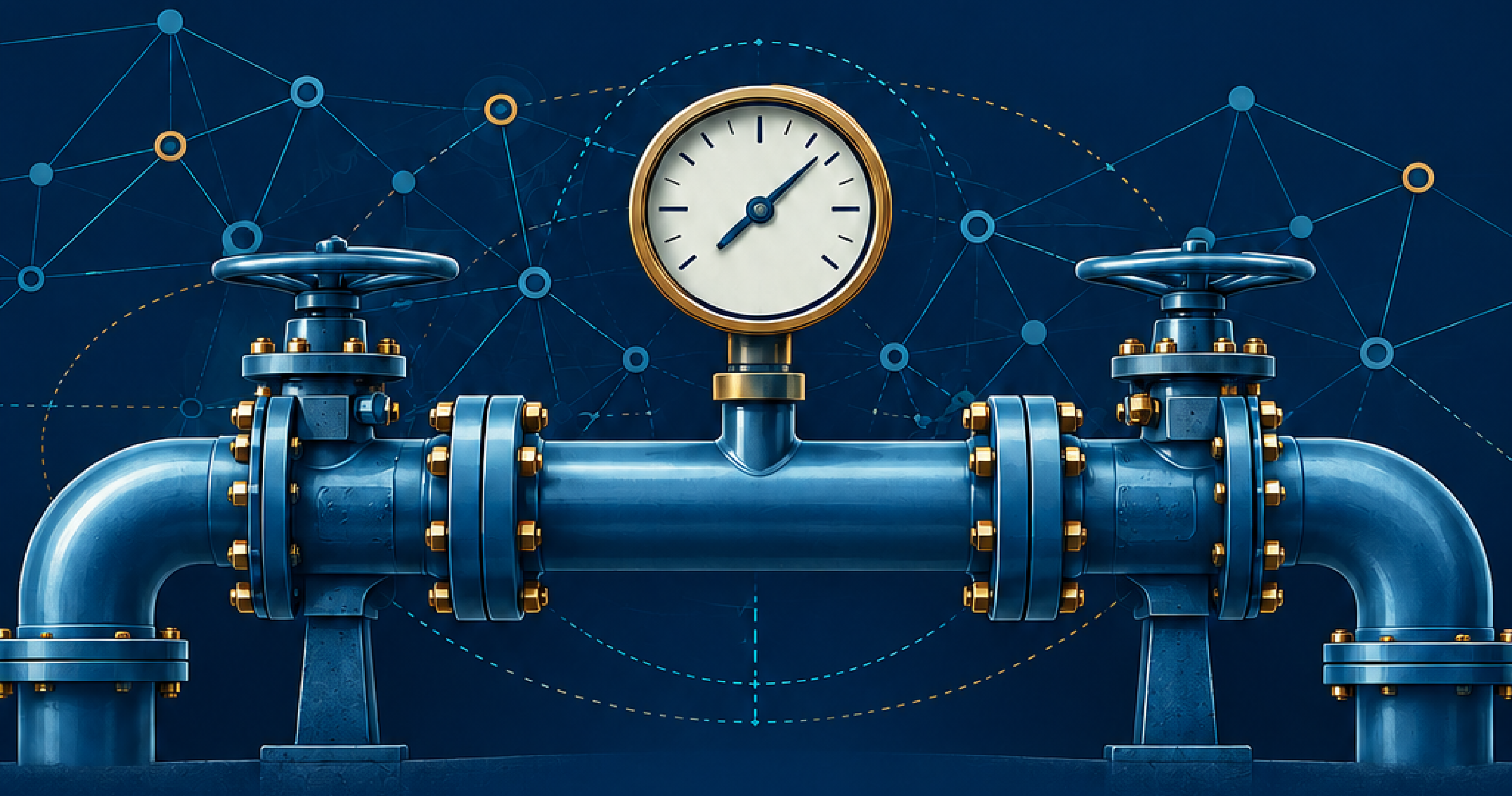

**Reviewers:**

**Nadir Bafadin oglu Agayev**
Professor, Department of Computer Systems and Programming, National Aviation Academy, Doctor of Technical Sciences

**Fagan Ganbar oglu Aliyev**
President of the International Academy of Ecoenergy, Head of the Department of Ecology, Azerbaijan University of Architecture and Construction, Doctor of Technical Sciences, Professor

# TABLE OF CONTENTS

# Preface

Analysis of global trends, issues in the development of the gas industry, and possible levels of gas consumption in various countries indicates that natural gas is increasingly taking on a more significant role compared to other fuel types in the developed countries of the world. However, with the increase in consumers adapted to using gas, there is a need for a substantial increase in natural gas reserves, expansion of networks, and regulatory complexes. Growing gas demand and resource consumption can be met not only through the construction of new networks and structures but also through the reconstruction of gas pipelines using potential reserves and efficient management. Ensuring the reliable operation of existing gas supply systems requires ongoing operational and capital expenditures, meaning that the continuous updating of morally and physically outdated equipment and networks is necessary.

The gas supply system, as an energy transmission system, continues to develop successfully compared to other fuel types. Along with the development of the resource base, the management system related to the reliable operation of the fuel and energy complex is also continually evolving.

This reliability is primarily ensured through high-tech safety measures, considering the current trends in the globalization of the energy sector.

However, this process is accompanied by reconstruction in existing gas pipelines. Currently, the reconstruction of main gas pipelines is a priority area for the innovative development of the management system, taking into account the use of advanced technologies and the operational principles of the fuel and energy complex. This also facilitates the prompt resolution of automation issues.

An important feature of the monograph is the extensive analysis and examination of the parameters of non-stationary gas flow in complex gas pipelines, with a key characteristic being the in-depth study of the problem. Specifically, obtaining mathematical expressions for non-stationary gas flow allows for targeted control of gasodynamic processes. To this end, effective solutions for the reconstruction of gas transport systems have been investigated based on an information database about pipeline failure conditions, which was derived from changes in the technological parameters of gas flow in the peripheral sections of pipelines during emergency modes. Simultaneously, the resolution of dynamic issues in complex main gas pipelines has led to the discovery of previously unknown relationships and regularities, a more comprehensive analysis of problems, and the provision of well-founded conclusions and recommendations. Based on the results of theoretical studies, a calculation scheme has been developed that enables the creation of normative-technical documents for justifying organizational and technological decisions in the operation of main gas pipelines.

Ensuring a high level of guaranteed gas supply to consumers is typically achieved through the reserve capacity of the pipeline system. Therefore, in the reconstruction of main gas pipeline systems, networks are generally proposed as parallel multi-line configurations. In the event of an emergency, due to the operational management of pipeline sections, there is no significant interruption in gas supply at the output (end) of the multi-line system, and consumers receive a continuous gas supply.

During the justification stage of the reconstruction of main gas pipelines, to achieve the required levels of reliability and productivity, pipelines are usually arranged in multi-line configurations (a minimum of two lines) that are interconnected. Proposals have been made for the theoretical and technical foundations of operating parallel gas pipelines to ensure the efficient placement of distances between these connecting pipelines.

The high-level operational reliability and efficiency of the reconstructed gas transmission scheme are achieved under several specific conditions:

- During the reconstruction phase of the existing gas pipeline system, solutions used should generally be standard and supported by innovative management technologies.
- Analytical calculations and methods have been used to justify the use of optimal schemes and to find ways to reduce costs and gas losses at all stages of the development and operation of the gas transmission system.
- Reliable operation of main gas pipeline systems requires extensive information for managing the elements according to their nature. The information base is obtained through a systematic approach based on the joint consideration of complex concepts such as system theory and feedback.

The reconstruction strategy for the gas supply system involves innovative approaches, new technical solutions, technical regulation, and the application of modern equipment due to the ongoing digital transformation of the oil and gas industry.

The management of the gas transmission system utilizes data, processes, and the capabilities and complexity of the organizational structures of the enterprise. The reconstruction of the existing gas supply system requires significant changes in normative documents to establish a new management system. The transition to intellectual technologies during reconstruction will ensure competitiveness, technological and environmental safety, as well as reduced operational costs during the renewal of gas supply system facilities.

The implementation of feasible options for reconstruction in the gas transmission process is ensured by:

- Comprehensive automation of technological equipment based on the creation of a management system for main gas pipelines;
- Planning of repair and restoration works at failure sites considering the actual condition of technological equipment;
- Ensuring automated management of facilities using complex algorithms for process management and safety.

The reconstruction of the complex gas pipeline ensures automatic control and remote management of technological equipment in the linear section. Simultaneously, it requires the collection of data related to ensuring reliable and safe operation of gas transmission facilities and the creation of an automated information system.

This monograph presents the methodology, model, and theoretical foundations developed for ensuring the technological foundations of decision-making in the management system for the reconstruction of the gas supply system. It mainly discusses methods that have practical applications in the reconstruction of gas pipelines with various configurations.

The text and presentation style of the monograph are intended for a broad audience, primarily project designers, specialists from the linear operation department of the gas supply system, as well as master's and doctoral students. Chapter 1 discusses the fundamental principles of gas pipeline system reconstruction. Chapter 2 examines the reconstruction methods of gas pipelines with various configurations for technological process improvement. Chapter 3 presents theoretical studies on the optimization methods of the parameters of reconstructed gas pipelines. Chapter 4 provides data on the technological foundations of management decisions in gas pipeline reconstruction. The monograph is based on work prepared by us since 2004 and research results published in journals from local and foreign countries.

# Chapter 1. Key Principles of Gas Pipeline System Reconstruction

A new concept for the operation of complex gas pipeline systems has been developed, defining the operational limitations of the equipment, and deriving mathematical expressions for managing gas-dynamic processes. Based on the obtained theoretical foundation, technological contradictions in the gas pipeline reconstruction process have been clarified and resolved within the proposed operational principle of emergency automatic valves installed on the gas pipelines. Consequently, a calculation method has been developed to determine the calibration indicator for regulating the valves of emergency automatic cranes in the improvement of gas transmission system management. This method is focused on optimizing the closing time of the valve in case of a gas pressure drop in the main pipeline during an emergency and obtaining an analytical expression for determining this time.

The reconstruction of gas supply systems primarily involves modernizing existing network equipment to increase gas consumption capacity in the future, ensure high reliability in gas supply, and improve equipment efficiency. One of the reconstruction options for gas pipelines that have been in operation for many years includes replacing outdated and physically worn-out equipment with modern, highly efficient, reliable, and safe equipment that meets current standards, thereby improving the system's management. This reconstruction option is mainly implemented by increasing the capacity of gas pipelines and introducing new equipment. One of the most critical pieces of equipment requiring further research is the emergency automatic shut-off valves. These valves are used in main gas pipelines to prevent catastrophic escalation in the event of a breach in pipeline integrity.

The working principle of these valves is based on the rapid drop in gas pressure within the main gas pipeline. If the pressure drop rate exceeds 0.1 MPa/min, the sensors detect the signal, and the valve automatically activates to close the pipeline. The pressure drop rate is determined by the difference in the stable gas pressure at the valve location in the main pipeline and the pressure decline during an emergency. This corresponds to the pressure drop caused by a complete rupture of the main gas pipeline. Theoretically and practically, the pressure drop rate can reach 0.2-0.5 MPa per minute in the event of a full pipeline rupture [3,14]. Therefore, during the reconstruction phase of gas pipelines, developing a new concept for the operation of automatic equipment that closes the pipeline in an emergency and operates efficiently at lower gas pressure drop rates is critical.

The paper provides a comprehensive overview of the methods for solving special differential equations describing the fundamental laws of gas movement, the calculation of non-stationary operating conditions of gas pipelines, and the theoretical investigation of gas dynamic parameters during non-stationary processes in the pipeline. Additionally, the conformity of the law governing the variation of pressure along the length of the pipeline over time under various operating conditions in a parallel pipeline system has been determined. Furthermore, recommendations for boundary conditions ensuring continuous operation in the pipeline section are proposed in relation to the occurrence of accidents in the gas pipeline.

The efficient operation regime of the pipeline system is determined to prevent the sudden cessation of gas supply to consumers through the theoretical investigation of the obtained analytical expressions of non-stationary gas dynamics. The purpose of the research is based on the premise that the strategy of reconstruction for managing non-stationary gas transportation processes with parallel pipeline systems relies on the characteristics of operating conditions to minimize the potential damage accidents could cause.

## 1.1. Study of Effective Solutions for the Improvement of Complex Gas Supply System Management

Until recently, these valves were exclusively mechanical. The advantage of a mechanical automatic valve is that it is entirely independent of the power supply system, as the regulation of the automatic valve's parameters is carried out using a special calibration valve. The downside is that the operation, maintenance, calibration, and adjustment of these complex mechanical systems require high costs. Currently, several foreign companies produce electronic-type automatic shut-off valves with various designs for main gas pipelines. These electronic automatic valves fulfill all the functions of mechanical ones while being simpler and more cost-effective to operate.

Electronic automatic shut-off valves use a pressure sensor to continuously measure the gas pressure in the pipeline before or after the valve, compare it with the predefined point, and close the valve if necessary. These modern valves are designed to completely eliminate gas leaks and save energy. However, as previously mentioned, their operating principle has not changed: like before, each valve automatically activates based on the gas pressure drop rate in the pipeline. The main difference is that several pressure sensors are now connected to a single point, meaning the valve closest to the gas leak will automatically shut off.

It is clear that in subsequent moments, starting from the first activated automatic valve, the automatic valves installed downstream along the pipeline will close either sequentially or randomly, depending on the ratio between the gas consumption of the consumers and the amount of leaking gas. If the ratio is less than one, they will close sequentially; if greater than one, they will close randomly. In the section of the pipeline from the starting point to the first closed valve, no further closures will occur since the pressure will increase over time.

Of course, this operational principle can be considered acceptable for linear gas pipelines. However, for complex gas pipelines, adopting this principle for efficiently managing non-stationary gasdynamic processes is not feasible. Therefore, theoretical studies are needed to determine a new method for activating automatic emergency valves to improve the management of complex gas supply systems and provide effective solutions for their reconstruction.

Many studies related to this issue have devoted, to varying degrees, to examining relatively localized problems. In other words, the calibration of valves in cases of a breach (rupture) in linear gas pipelines is conducted based on the pressure drop rate at the right and left sections of each valve [15,16,90,96].

The key difference between the current investigation we are conducting and the existing solved research lies in the development of a new operational concept for the equipment of complex gas pipelines. This is achieved through mathematical modeling of non-stationary gas flow, determining operational time limits for the equipment, and deriving mathematical expressions for controlling gas-dynamic processes. Thus, based on the theoretical foundation obtained, we clarify and resolve technological contradictions within the framework of the future operational principles of the devices installed in the gas pipeline during the reconstruction of the gas transmission system. Moreover, a calculation scheme is developed to improve the management of the gas transmission system.

To achieve this, we propose the following scheme to define the effective operational principles of devices installed on complex gas pipelines (Scheme 1.1.1).

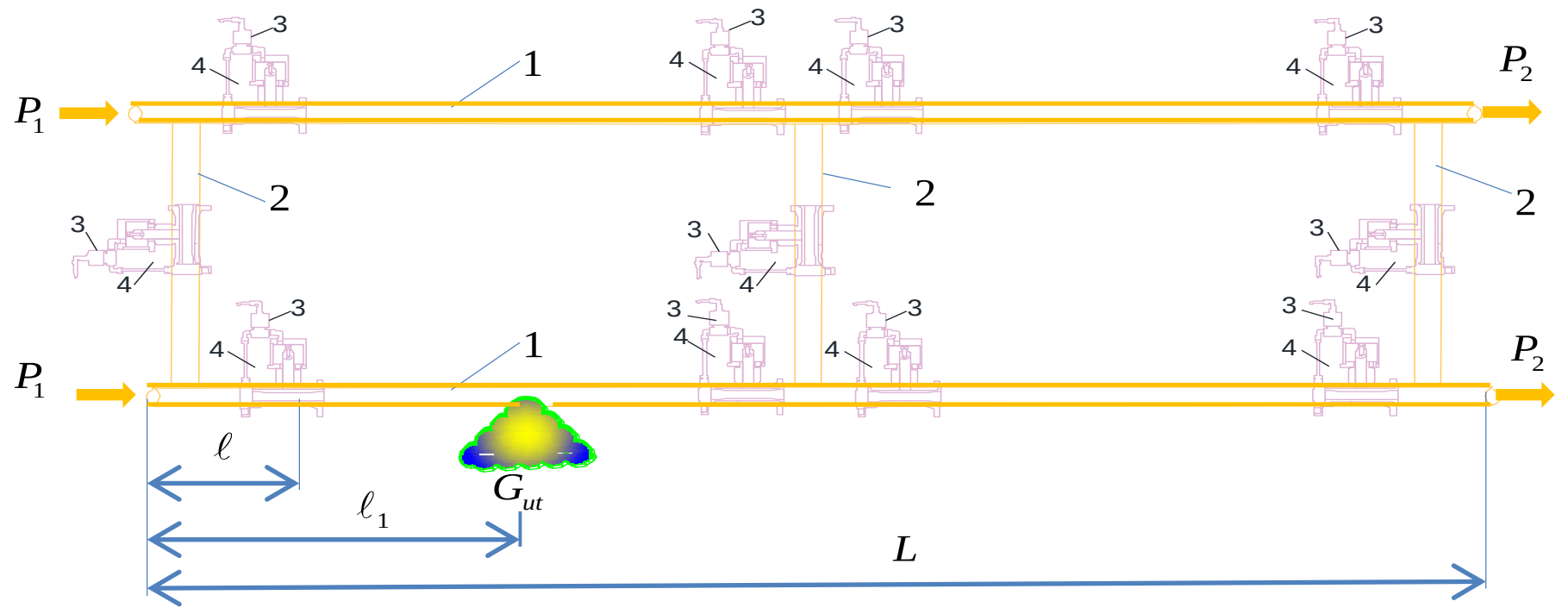

**Scheme 1.1.1. Efficient Placement Scheme of Emergency Automatic Valves in Complex Gas Pipelines**

**1- Pipeline of the parallel gas pipeline; 2- Connecting line of the parallel pipeline; 3-Valves that activate the automatic shut-off valves; 4- Automatic shut-off valves; L – The length of the gas pipeline; $G_{ut}$ - $x=\ell_1$ point where gas is leaking from the breached pipeline; $P_1$ and $P_2$ - The pressure of the gas flow in steady-state at the beginning and end of the pipeline.**

Using the Scheme for Installing Emergency Automatic Shut-off Valves in Complex Gas Pipelines, the new operational principle for the automatic device is examined for isolating the damaged section of the parallel gas pipeline from the main pipeline. This method focuses on optimizing the closing time of the valve during an emergency when the gas pressure drops in the main gas pipeline, and deriving an analytical expression for determining this time.

As shown in Scheme 1.1.1, except for the connecting pipes installed at the beginning and end of the gas pipeline, two automatic valves are placed beside each connecting pipe—one on each side. However, in the event of an emergency, only one of the valves next to the connecting pipe should operate. On the other hand, to prevent gas loss, the automatic valves on both sides of the breach ($x=\ell_1$) should close simultaneously ($t=t_1$). Therefore, the traditional principle of each valve automatically closing based on the pressure drop rate individually is not suitable. The pressure sensors of the automatic shut-off valves should be connected to a single point in the pipeline. In other words, the valves on either side of the gas leak should automatically close at the same time.

Thus, according to Scheme 1.1.1, to ensure the effective operation of the emergency automatic valves, it is necessary to coordinate the simultaneous operation of the valves located on the right side of the first connecting pipe and the left side of the next connecting pipe starting from the beginning of the gas pipeline. This is because, to continuously supply gas to the consumers fed by the parallel gas pipeline, the automatic valves installed on the connecting pipes on both sides of the breach must open simultaneously. The determination of the activation time for these automatic valves has been extensively studied [3]. It has been noted that these valves are activated remotely by automated control systems as soon as the breach location is identified at the dispatch centers [3,95].

Therefore, improving the management of the gas transmission system involves determining the activation time of the automatic valves installed on the gas pipeline and limiting this time based on the breach location, which is a key aspect of our research. At the same time, analytically determining the pressure drop rate at the end of the gas pipeline depending on the amount of leaking gas is crucial. This will enhance the reliability of the gas pipeline, ensure continuous gas supply to consumers, and guarantee reliable sealing of the breach. Moreover, timely activation of the automatic shut-off devices in case of a breach in the main gas pipeline prevents the catastrophic development of the emergency situation and effectively controls gas loss. This prevents significant gas leakage into the atmosphere and provides a basis for continuous gas supply to consumers. Therefore, applying this new operational principle of modern devices in the reconstruction of long-used gas pipelines will contribute to improving their reliability.

## Materials and Methods

To develop the new concept for automatic shut-off valves, it is essential to understand the nature of pressure changes in non-stationary gas flow during emergency situations (e.g., breach at the point $x = \ell_1$ in the gas pipeline). Furthermore, obtaining mathematical expressions for non-stationary gas flow enables purposeful management of gas-dynamic processes. Thus, solving dynamic problems in complex main gas pipelines will reveal previously unknown relationships and patterns, allow for a more comprehensive analysis of the problem, and provide well-founded conclusions and recommendations.

Non-stationary processes involve changes in gas flow parameters over time at each cross-section of the pipeline. For non-stationary gas flow, I.A. Charny proposed a system of differential equations related to gas flow parameters, such as pressure P(x,t) and gas mass flow G(x,t), as functions of two variables. Let us consider the proposed model [78].

$$\begin{cases} -\dfrac{\partial P}{\partial x} = \lambda \dfrac{\rho V^2}{2d} \\ -\dfrac{1}{c^2}\dfrac{\partial P}{\partial t} = \dfrac{\partial G}{\partial x} \end{cases} \tag{1}$$

Here, $G=\rho v$; $P=\rho \cdot c^2$ ; $\rho = \dfrac{P}{zRT}$

**T** – Average temperature of the gas;
**z = z(p,T)** – Compressibility factor, dependent on pressure and temperature;
**c** – Speed of sound in the gas for isothermal processes, m/sec;
**P** – Pressure in the cross-sections of the gas pipeline, Pa;
**x** – Coordinate along the axis of the gas pipeline, m;
**t** – Time coordinate, sec;
**λ** – Hydraulic resistance coefficient;
**ρ** – Average density of the gas, kg/m$^3$;
**v** – Average velocity of the gas flow, m/sec;
**d** – Diameter of the gas pipeline, m;
**R** – Universal gas constant, J/(mol·K);
**G** – Mass flow rate of the gas, $\dfrac{Pa \times \sec}{m}$

It is evident that solving the system of nonlinear equations (1) is not feasible because such equations are not integrable. In engineering calculations, approximate methods are employed. In other words, the equations are linearized. The linearization method proposed by I.A. Charny is more convenient for solving these equations. The expression characterizing this linearization is as follows [78].

$$2a = \lambda \frac{V}{2d} \tag{2}$$

Here, 2a is the linearization coefficient proposed by I.A. Charny, then the system of equations (1) and (2) will be transformed as follows

$$\begin{cases} -\dfrac{\partial P}{\partial X} = 2aG \\ -\dfrac{1}{c^2}\dfrac{\partial P}{\partial t} = \dfrac{\partial G}{\partial X} \end{cases} \qquad (3)$$

If we take the derivative with respect to x from the first expression of the system of equations (3) and determine $\frac{\partial G}{\partial X}$ from the second expression, then the equation will be as follows.

$$\begin{cases} -\dfrac{\partial^2 P}{\partial X^2} = 2a\dfrac{\partial G}{\partial X} \\ \dfrac{\partial G}{\partial X} = -\dfrac{1}{c^2}\dfrac{\partial P}{\partial t} \end{cases} \qquad (4)$$

Then, by substituting the expression for $\frac{\partial G}{\partial X}$ from the first expression into the second expression in equation (4), we obtain the heat transfer equation for the mathematical solution of the problem.

$$\frac{\partial^2 P}{\partial X^2} = \frac{2a}{c^2}\frac{\partial P}{\partial t} \qquad (5)$$

Thus, for the uniqueness of the solution to the equation (5) reflecting the studied process and the complete determination of the given process, the initial and boundary conditions must be specified correctly. The initial conditions provide the distribution of the sought function (pressure) along the pipeline length at t=0 when inertia forces are not considered. The boundary conditions ensure the predetermined variations of pressure and flow at the endpoints of the pipeline during the process under observation. Properly specifying the boundary conditions completes the mathematical model of the process and allows for a thorough and accurate investigation of the occurring physical phenomena.

In stationary regimes in parallel pipeline lines, pressure distribution depending on the coordinate is assumed to be linear. The diameter of the pipeline lines in the examined gas pipeline is uniform throughout. The supply intensity and the total mass flow rate of the gas, $G_0$, are given, and the pressure at the starting point is $P_1$. Then,

$$t=0;\; P_i(x,0)=P_1-2aG_0x_i \quad i=1,2$$

For the main sections of the gas pipeline, the boundary conditions for the examined process can be obtained as follows (Scheme 1.1.1). In previous research, the adjustment of valves in the operation of automatic shut-off valves was based on the rate of pressure drop [15,90]. These studies confirmed that the pressure change in the section where the pipeline is breached follows an exponential dependency. Therefore, the gas mass flow rate in that section can be represented as follows.

$$\frac{P_1(\ell,t)}{\partial x} = -P_1\frac{2a\beta}{g}e^{-\beta t}$$

It is clear that at the moment the automatic valve closes ($t=t_1$), the gas mass flow rate $\frac{P_1(\ell,t)}{\partial x} = 0$ will be zero. To express this physical process at the boundary condition $x=\ell$, we consider the Heaviside function $\sigma(t)$. The Heaviside function is a piecewise constant function, where the value is zero for negative arguments and one for positive arguments.

Thus, for the complete determination of the examined physical process under the condition of maintaining a constant mass flow rate at the starting point of the gas pipeline (x=0), we accept the boundary conditions as follows.

At the point $x=0 \quad \frac{P_1(x,t)}{\partial X} = 2aG_0$

At the point $x=\ell \quad \frac{P_1(x,t)}{\partial X} = -P_1 \frac{2a\beta}{g} e^{-\beta t}\left[\sigma(t)-\sigma(t-t_1)\right]$

At the point $x=\ell_1 \quad \begin{cases} P_1(x,t) = P_2(x,t) \\ \frac{P_1(x,t)}{\partial X} - \frac{P_2(x,t)}{\partial X} = -2aG_{ut}\left[\sigma(t)-\sigma(t-t_1)\right] \end{cases}$

It is clear that by applying the Laplace transform [P(x,s)= $\int_0^\infty P(x,t)e^{-st}dt$ ] to equation (5), it is converted into a second-order ordinary differential equation. The boundary conditions will then lead to general solutions as follows.

At the point $x=0 \quad \frac{P_1(x,s)}{\partial x} = \frac{2aG_0}{s}$

At the point $x=\ell \quad \frac{P_1(x,s)}{\partial x} = -P_1 \frac{2a\beta}{g} \frac{\left(1-e^{-st_1}\right)}{(s+\beta)}, \qquad t<t_1$

At the point $x=\ell_1 \quad \begin{cases} P_1(x,s) = P_2(x,s) \\ \frac{P_1(x,s)}{\partial x} - \frac{P_2(x,s)}{\partial x} = \frac{-2aG_{ut}}{s}\left(1-e^{-st_1}\right) \end{cases} \qquad t<t_1$

$$P_1(x,s) = \frac{P_1 - 2aG_0 x}{s} + A_1 sh\mu x + A_2 ch\mu x \qquad 0 \le x \le \ell_1 \tag{6}$$

$$P_2(x,s) = \frac{P_1 - 2aG_0 x}{s} + A_3 sh\mu x + A_4 ch\mu x \qquad \ell_1 \le x \le L \tag{7}$$

Here, $\mu = \sqrt{\frac{2as}{c^2}}$

By considering the initial and boundary conditions in equations (6) and (7), we determine the coefficients $A_1$, $A_2$, $A_3$ and $A_4$ .

$$A_1 = 0 \;; \qquad A_2 = \frac{2aG_0}{\mu s \cdot sh\mu\ell} - P_1 \frac{2a\beta}{g\mu} \frac{\left(1-e^{-st_1}\right)}{(s+\beta)sh\mu\ell} \;; \qquad A_3 = \frac{2aG_{ut}}{\mu s}\left(1-e^{-st_1}\right)ch\mu\ell_1$$

$$A_4 = \frac{2aG_0}{\mu s \cdot sh\mu\ell} - P_1 \frac{2a\beta}{g\mu} \frac{\left(1-e^{-st_1}\right)}{(s+\beta)sh\mu\ell} - \frac{2aG_{ut}}{\mu s}\left(1-e^{-st_1}\right)sh\mu\ell_1$$

Considering the values of the coefficients $A_1$, $A_2$, $A_3$ and $A_4$ in equations (6) and (7), we obtain the transformed equation of the examined process, or in other words, the dynamic state of the complex gas pipeline. The resulting equation will express the distribution of pressure along the length of the gas pipeline in the transformed form for the non-stationary flow regime of the complex gas pipeline.

$$P_1(x,s)=\frac{P_1-2aG_0x}{s}+\frac{2aG_0}{\mu s}\frac{ch\mu x}{sh\mu\ell}-P_1\frac{2a\beta}{g\mu}\frac{\left(1-e^{-st_1}\right)}{(s+\beta)}\frac{ch\mu x}{sh\mu\ell} \qquad 0\le x\le \ell_1 \tag{8}$$

$$P_2(x,s)=\frac{P_1-2aG_0x}{s}+\frac{2aG_0}{\mu s}\frac{ch\mu x}{sh\mu\ell}-P_1\frac{2a\beta}{g\mu}\frac{\left(1-e^{-st_1}\right)}{(s+\beta)}\frac{ch\mu x}{sh\mu\ell}- \qquad \ell_1\le x\le L \tag{9}$$
$$-\frac{2aG_{ut}}{\mu s}\left(1-e^{-st_1}\right)sh\mu\left(\ell_1-x\right)$$

To ensure the possibility of applying the inverse Laplace transform to each boundary, we express equations (8) and (9) in the following form.

$$P_1(x,s)=\frac{P_1-2aG_0x}{s}+\frac{2aG_0}{s\cdot b\sqrt{s}}\frac{chb\sqrt{s}x}{shb\sqrt{s}\ell}-P_1\frac{2a\beta}{gb\sqrt{s}}\frac{\left(1-e^{-st_1}\right)}{(s+\beta)}\frac{chb\sqrt{s}x}{shb\sqrt{s}\ell} \qquad 0\le x\le \ell_1 \tag{10}$$

$$P_2(x,s)=\frac{P_1-2aG_0x}{s}+\frac{2aG_0}{s\cdot b\sqrt{s}}\frac{chb\sqrt{s}x}{shb\sqrt{s}\ell}-P_1\frac{2a\beta}{g\cdot b\sqrt{s}}\frac{\left(1-e^{-st_1}\right)}{(s+\beta)}\frac{chb\sqrt{s}x}{shb\sqrt{s}\ell}- \qquad \ell_1\le x\le L$$
$$-\frac{2aG_{ut}}{2s}\left(1-e^{-st_1}\right)\left[\frac{shb\sqrt{s}\left(2\ell_1-x\right)}{b\sqrt{s}chb\sqrt{s}\ell_1}-\frac{shb\sqrt{s}x}{b\sqrt{s}chb\sqrt{s}\ell_1}\right] \tag{11}$$

Here, $b=\sqrt{\dfrac{2a}{c^2}}$

Based on the solution of the examined problem, it is more appropriate to determine the original solution of the mathematical expression for the non-stationary flow of gas in the pipeline. Therefore, the inverse Laplace transform of each boundary is applied as specified in equations (10) and (11). Consequently, the original expressions of equations (10) and (11) will be equal to the following functions.

$$0\le x\le \ell_1$$

$$P_1(x,t)=P_1-2aG_0x+\frac{2aG_0t}{b^2\ell}+\frac{4aG_0}{b^2\ell}\sum_{n=1}^{\infty}(-1)^n\,\mathrm{Cos}\frac{\pi nx}{\ell}\frac{\left(1-e^{-\alpha n^2t}\right)}{\alpha n^2}-$$
$$-\begin{cases}0\to t\succ t_1\\ P_1\dfrac{2a}{gb^2\ell}\left(1-e^{-\beta t}\right)-P_1\dfrac{4a\beta}{gb^2\ell}\displaystyle\sum_{n=1}^{\infty}(-1)^n\,\mathrm{Cos}\frac{\pi nx}{\ell}\frac{\left(e^{-\beta t}-e^{-\alpha n^2t}\right)}{\alpha n^2-\beta}\to t\prec t_1\end{cases} \tag{12}$$

$$\ell_1\le x\le L$$

$$P_2(x,t)=P_1-2aG_0x+\frac{2aG_0t}{b^2\ell}+\frac{4aG_0}{b^2\ell}(-1)^n\sum_{n=1}^{\infty}\mathrm{Cos}\frac{\pi nx}{\ell}\frac{\left(1-e^{-\alpha n^2t}\right)}{\alpha n^2}-$$
$$-\begin{cases}0\to t\succ t_1\\ P_1\dfrac{2a}{gb^2\ell}\left(1-e^{-\beta t}\right)-P_1\dfrac{4a\beta}{gb^2\ell}\displaystyle\sum_{n=1}^{\infty}(-1)^n\,\mathrm{Cos}\frac{\pi nx}{\ell}\frac{\left(e^{-\beta t}-e^{-\alpha n^2t}\right)}{\alpha n^2-\beta}-\\ -2aG_{ut}\left[\left(\ell_1-x\right)-\dfrac{1}{b^2\ell_1}\displaystyle\sum_{n=1}^{\infty}\left[(-1)^{n+1}-1\right]\sin\frac{\pi x(2n-1)}{2\ell_1}\frac{e^{-\alpha_1(2n-1)^2t}}{\alpha_1(2n-1)^2}\right]\to t\prec t_1\end{cases} \tag{13}$$

Here, $b^2 = \frac{2a}{c^2}$, $\alpha = \frac{\pi^2 c^2}{2a\ell^2} . \alpha_1 = \frac{\pi^2 c^2}{8a\ell_1^2}$

As a result of deriving the mathematical expressions (12) and (13) for unsteady gas flow, it is possible to achieve the purposeful solution of the problems set forth. First, we determine the value of the pressure function at the point x= $\ell$ in equation (12) to analytically determine the operation time of automatic valves.

$$P(\ell,t) = P_1 - 2aG_0\ell + \frac{c^2 G_0 t}{\ell} + \frac{2c^2 G_0}{\ell}\sum_{n=1}^{\infty}\frac{(1-e^{-\alpha n^2 t})}{\alpha n^2} - \\ -\begin{Bmatrix} 0 \to t \succ t_1 \\ P_1 \frac{c^2}{g\ell}\left(1-e^{-\beta t}\right) - P_1 \frac{2c^2\beta}{g\ell}\sum_{n=1}^{\infty}\frac{\left(e^{-\beta t} - e^{-\alpha n^2 t}\right)}{\alpha n^2 - \beta} \to t \prec t_1 \end{Bmatrix} \tag{14}$$

To obtain a more concise mathematical expression for the solution of the problem under consideration and to study the process of pressure change at the point x= $\ell$ in the gas pipeline requiring investigation, we accept the simplification n=2, which has been theoretically and experimentally tested in known areas of application for expression (14) [28]. Then,

$$P(\ell,t) = P_1 - 2aG_0\ell + \frac{c^2 G_0 t}{\ell} + \frac{2c^2 G_0}{\ell}\frac{(1-e^{-4\alpha t})}{4\alpha} - \\ -\begin{Bmatrix} 0 \to t \succ t_1 \\ P_1 \frac{c^2}{g\ell}\left(1-e^{-\beta t}\right) - P_1 \frac{2c^2\beta}{g\ell}\frac{\left(e^{-\beta t} - e^{-4\alpha t}\right)}{4\alpha - \beta} \to t \prec t_1 \end{Bmatrix} \tag{15}$$

If we take two series from the Maclaurin expansion of exponential functions, as a result of several mathematical operations, equation (15) for the cases when $t \prec t_1$ will be as follows.

$$P(\ell,t) = P_1 - 2aG_0\ell + 3\frac{c^2 G_0}{\ell}t - 3P_1\frac{c^2\beta}{g\ell}t \tag{16}$$

It is known that, according to the conditions of the study, both automatic emergency valves located to the right and left of the gas leakage point are designed to automatically shut off the damaged sections of the gas pipeline during accidents and emergencies at the same time. Additionally, the valves installed on the cranes are connected to gas pressure sensors at the same point (the end section) of the gas pipeline and close simultaneously. For this reason, in order to properly analyze the results of the study, we consider the pressure impulse of the emergency automatic valve as the pressure change occurring at the end section of the pipeline. Therefore, the calibration of the emergency automatic valves is based on the pressure drop indicator of the gas flow in the pipeline. If the pressure in the section where the shut-off valve is installed decreases by 20% within a short period of time, the valves automatically activate both valves and completely isolate the damaged section of the pipeline from the main part of the gas pipeline. We assume the closing time of the valves to be $t = t_1$ . Then, if we accept $P(\ell,t_1)$ = 0,8 $P(\ell,0)$ in equation (16), we derive the following formula to determine the operating time of the valves through the cranes.

$$t_1 = \frac{0,2\ell\left(P_1 - 2aG_0\ell\right)}{3c^2\left(P_1\dfrac{\beta}{g} - G_0\right)} \tag{17}$$

From the analysis of the expressions in formula (17), it is evident that after determining the value of the expression β, which characterizes the pressure drop rate in addition to the parameters of the gas pipeline, we can calculate the operating time of the valves. For this purpose, by substituting x = L into equation (13), we obtain the following expression to determine the time-dependent change in pressure at the end section of the gas pipeline under the condition $t \prec t_1$.

$$\begin{aligned} P_2(L,t) &= P_1 - 2aG_0L + \frac{c^2G_0t}{\ell} + \frac{2c^2G_0}{\ell}\sum_{n=1}^{\infty}(-1)^n \operatorname{Cos}\frac{\pi nL}{\ell}\frac{\left(1-e^{-\alpha n^2 t}\right)}{\alpha n^2} - \\ &- P_1\frac{c^2}{g\ell}\left(1-e^{-\beta t}\right) - P_1\frac{2c^2\beta}{g\ell}\sum_{n=1}^{\infty}(-1)^n \operatorname{Cos}\frac{\pi nL}{\ell}\frac{\left(e^{-\beta t}-e^{-\alpha n^2 t}\right)}{\alpha n^2-\beta} - \\ &- 2aG_{ut}\left[\left(\ell_1 - L\right) + \frac{c^2}{2a\ell_1}\sum_{n=1}^{\infty}\left[(-1)^{n+1}-1\right]\sin\frac{\pi L(2n-1)}{2\ell_1}\frac{e^{-\alpha_1(2n-1)^2 t}}{\alpha_1(2n-1)^2}\right] \end{aligned} \tag{18}$$

Based on the analysis of the study, the automatic valve closest to the leakage point and the valve on the other side of the gas leakage point must operate simultaneously. This means that the pressure sensors of both valves are connected to a single point. To this end, it is important to analyze the pressure drop rate at the end section of the pipeline using equation (18), depending on the location of the gas pipeline's integrity breach, the placement of the valves, and the amount of leaking gas.This analysis is broadly reflected in section 1.2.

### 1.2. Analysis of a New Concept for the Effective Operation Principle of Emergency Automatic Valves in the Reconstruction of Existing Parallel Gas Pipelines

In complex gas pipelines equipped with emergency automatic shut-off valves, the time-dependent change function of pressure at the end section of a cross-sectional area during an emergency regime (18) depends on the parameters in the stationary regime ($P_1, G_0, L$, etc.), as well as the coefficient characterizing the pressure drop rate in the section where the valves are installed, the leakage point of the pipeline, and the values of the mass of the leaking gas. Among these parameters, the leakage point ( $\ell_1$ ), the location of the operating automatic valve ( $\ell$ ), and the coefficient characterizing the pressure drop rate (β) are interdependent quantities.

From this perspective, we assume a distance of $\ell \prec \ell_1$ = 5000 m between the location of the automatic valve to be activated and the leakage point. Using expression (18), in order to study the time-dependent change in pressure at the end section of the gas pipeline (for determining the operating time of the shut-off valve) at different leakage points and different values of leaking gas, we accept the following data.

$P_1$ = 55 $\times 10^4$ Pa : $G_0$= 30 Pa ×sec/m : $2a = 0.1\frac{1}{\sec}$ ; $c = 383.3\frac{m}{\sec}$ ; d= 0,7 m: L = $10\times 10^4$ m;

1-st version $\ell_1 = 35000m$; 2-nd version $\ell_1 = 45000m$; 3-rd version $\ell_1 = 75000m$

If we assume that the distance between the location of the automatic valve and the leakage point is 5000 m, then for the first case $\ell = 30000m$, for the second case $\ell = 40000m$, and for the third case $\ell = 70000m$. It is evident that in equation (18), the coefficient β characterizing the rate of pressure drop can vary along the length of the gas pipeline. Therefore, according to the requirements of the unsteady gas flow gasodynamic process, we accept the values of β. Specifically, the values of β will differ depending on the distance between the pressure measurement point (pressure sensor) and the points where the gas leakage occurs. As the distance between the coordinates decreases, the value of β will increase based on the known gasodynamic process. On the other hand, for fixed distances between coordinates, the values of β should increase as the mass of gas leaking into the environment at the breach point of the pipeline increases. Thus, in accordance with theoretical and practical research [1,8], we accept the values of β and use equation (18) to determine the values of gas pressure over time at the end section of the gas pipeline. The results of the report are noted in the following table (Table 1.2.1).

| t, sec | P(L,t), $10^4$ Pa. | | | | | | | | |
|---|---|---|---|---|---|---|---|---|---|
| | $\ell = 30000m$, $\ell_1 = 35000m$ | | | $\ell = 40000m$, $\ell_1 = 45000m$ | | | $\ell = 70000m$, $\ell_1 = 75000m$ | | |
| | β= $0{,}9\cdot 10^{-3}$ | β= $1{,}6\cdot 10^{-3}$ | β= $2{,}5\cdot 10^{-3}$ | β= $1{,}3\cdot 10^{-3}$ | β= $2{,}5\cdot 10^{-3}$ | β=$4{,}1\cdot 10^{-3}$ | β= $2{,}1\cdot 10^{-3}$ | β= $3{,}7\cdot 10^{-3}$ | β= $5{,}8\cdot 10^{-3}$ |
| | $G_{ut}$= 0,1$G_0$ | $G_{ut}$= 0,3$G_0$ | $G_{ut}$= 0,5$G_0$ | $G_{ut}$= 0,1$G_0$ | $G_{ut}$= 0,3$G_0$ | $G_{ut}$= 0.5$G_0$ | $G_{ut}$= 0.1$G_0$ | $G_{ut}$= 0.3$G_0$ | $G_{ut}$= 0.5$G_0$ |
| 250 | 24,87 | 24,83 | 24,61 | 24,09 | 23,98 | 23,87 | 23,23 | 22,71 | 22,26 |
| 300 | 24,66 | 24,21 | 23,76 | 23,66 | 23,08 | 22,76 | 22,52 | 21,78 | 21,18 |
| 350 | 24,50 | 23,72 | 23,12 | 23,28 | 22,35 | 21,96 | 22,11 | 21,06 | 20,46 |
| 400 | 24,38 | 23,32 | 22,66 | 22,96 | 21,77 | 21,41 | 21,69 | 20,53 | 20,02 |
| 450 | 24,30 | 23,01 | 22,33 | 22,69 | 21,33 | 21,06 | 21,35 | 20,16 | 19,81 |
| 500 | 24,26 | 22,78 | 22,13 | 22,48 | 20,99 | 20,88 | 21,09 | 19,94 | 19,78 |

**Table 1.2.1. Pressure Values in the Last Section of Gas Pipelines During Non-Stationary Regimes Depending on the Location of the Accident and the Volume of Leaking Gas.**

Using the values from Table 1.2.1, we create a graph depicting the time-dependent changes in pressure in the last section of the pipeline for different accident locations and volumes of leaking gas (Graph 1.2.1).

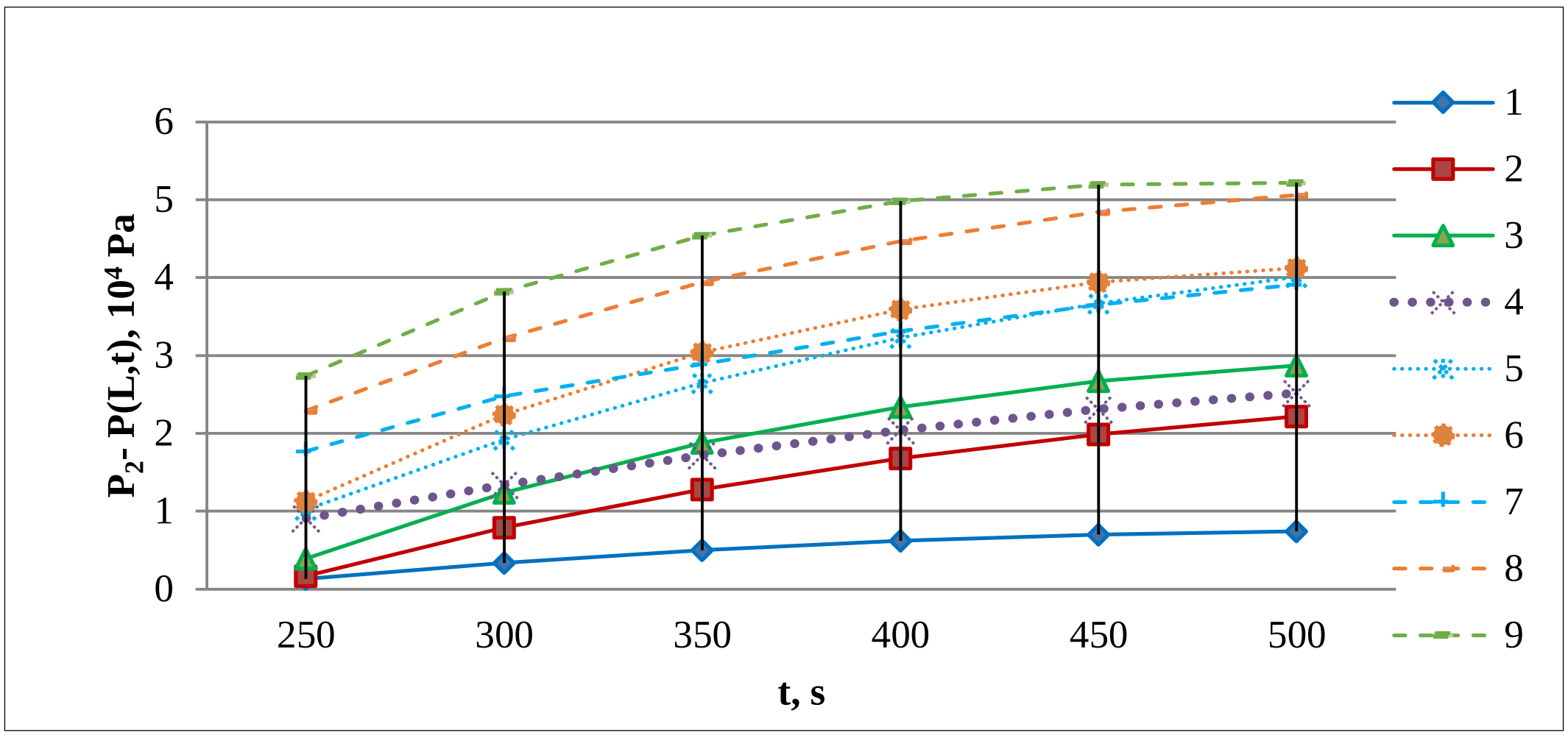

**Graph 1.2.1. Shows the change in the rate of pressure drop at the end section of the gas pipeline for different values of leaking gas mass and different leakage points. (1,2,3-$\ell_1 = 35000m$; 4,5,6-$\ell_1 = 45000m$; 7,8,9-$\ell_1 = 75000m$; 1,4,7-$G_{ut}$=0,1$G_0$; 2,5,8-$G_{ut}$=0,3$G_0$; 3,6,9- $G_{ut}$=0,5$G_0$).**

In the analysis of Figure 1.2.1, the first notable gasodynamic process is that, regardless of the different leakage points, the rate of pressure drop initially reaches a minimum value at the early moments of time, and this value increases as time progresses. Secondly, an increase in the volume of gas leakage also raises the value of the pressure drop rate.

If the pipeline's integrity breach occurs at a distance of $\ell_1 = 35000m$, and under the condition of $G_{ut}$=0,1$G_0$, the rate of pressure drop at the end section of the pipeline is 0.13 × $10^4$ Pa at t=250 seconds. Under the condition of $G_{ut}$=0,3$G_0$, the pressure drop rate increases to 0.17 × $10^4$ Pa at t=250 seconds, and under $G_{ut}$=0,5$G_0$, it further increases to 0.39 × $10^4$ Pa at t=250 seconds. Similarly, under the condition of $G_{ut}$=0,1$G_0$, the pressure drop rate at the end section of the pipeline is 0.62 × $10^4$ Pa at t=400 seconds, under $G_{ut}$=0,3$G_0$ it is 1.68 × $10^4$ Pa, and under $G_{ut}$=0,5$G_0$, it increases to 2.34 × $10^4$ Pa at t=400 seconds (Table 1.2.1 and Graph 1.2.1 ).

Based on the analysis, the third observed gasodynamic process indicates that as the leakage point of the gas pipeline approaches the end of the pipeline, the rate of pressure drop at the end section increases. For example, when
$\ell_1 = 45000m$, the increase in pressure drop rate compared to values observed at a distance of $\ell_1 = 35000m$ will be as follows:

Under the condition of $G_{ut}$=0,1$G_0$, if the pressure drop rate at the end section of the pipeline increases from 0.62 × $10^4$ Pa to 2.04 × $10^4$ Pa at t=400 seconds, then under $G_{ut}$=0,3$G_0$ it increases from 1.68 × $10^4$ Pa to 3.23 × $10^4$ Pa, and under $G_{ut}$=0,5$G_0$ it increases from 2.34 × $10^4$ Pa to 3.59 × $10^4$ Pa. For a distance of $\ell_1 = 75000m$, compared to the values observed at a distance of $\ell_1 = 35000m$, the increase will be as follows: Under $G_{ut}$=0,1$G_0$, the pressure drop rate at the end section of the pipeline increases from 0.62 × $10^4$ Pa to 4.29 × $10^4$ Pa at t=400 seconds; under $G_{ut}$=0,3$G_0$ it increases from 1.68 × $10^4$ Pa to 4.84 × $10^4$ Pa, and under $G_{ut}$=0,5$G_0$ it increases from 2.34 × $10^4$ Pa to 4.98 × $10^4$ Pa. However, as the leakage point of the gas pipeline approaches the end of the pipeline, the increase in pressure drop rate exceeds the increase in pressure due to the volume of leaking gas.

Thus, for cases where the leakage point is at distances of $\ell_1 = 35000m$ and $\ell_1 = 45000m$, the difference in the measured pressure drop rates at the end section of the pipeline is as follows: at t=250 seconds, for a gas leakage with $G_{ut}$=0,1$G_0$, the difference is 0.78 × $10^4$ Pa. Despite keeping the volume of gas leakage constant, the difference increases to 1.22 × $10^4$ Pa at t=350 seconds and further increases to 1.78 × $10^4$ Pa at t=500 seconds.

However, with an increase in the mass of gas leakage, specifically under the condition $G_{ut}=0,3G_0$ , the difference in the pressure drop rates at t=500 seconds is approximately $1.79 \times 10^4$ Pa, similar to the value observed under $G_{ut}=0,1G_0$ . Under the condition $G_{ut}=0,5G_0$ , the difference in pressure drop rates at t=500 seconds decreases to $1.25 \times 10^4$ Pa.

Similarly, this process can be observed in cases where the leakage point is at distances of $\ell_1 = 35000m$ and $\ell_1 = 75000m$. Under the condition $G_{ut}=0,1G_0$, the difference in the pressure drop rates at t=500 seconds is $1.39 \times 10^4$ Pa. Under $G_{ut}=0,3G_0$ , the difference is $1.05 \times 10^4$ Pa at t=500 seconds, and under $G_{ut}=0,5G_0$ , the difference decreases to $1.10 \times 10^4$ Pa at t=500 seconds. Therefore, considering the volume of the leaking gas is more important for the rapid activation of automatic valves.

Based on the analysis of the graph, it can be noted that if traditional emergency automatic valves are installed in the pipeline of the analyzed gas pipeline, they will only activate within **t=400** seconds (6.7 minutes) when a gas leakage occurs at $\ell_1 = 75000m$ distance and $G_{ut}=0,5G_0$ volume, by receiving a valve signal to open the valve. This is because at the point where the valve is located, the stationary gas pressure decreases by 20% within the specified time. At this time, the rate of pressure drop will be 0.08 MPa/min. According to previous studies by authors who investigated this problem, when the rate of pressure drop in a main gas pipeline exceeds 0.1 MPa/min, the valves that receive the signal automatically open the valve. Therefore, based on the traditional operation principle of valves, when the integrity of the pipeline of the analyzed gas pipeline is compromised (broken) and the mass of the leaking gas exceeds 50% of the gas mass transported by the pipeline, and the distance between the rupture point and the location of the automatic valves installed on the pipeline is 25,000 m, the valves will close. In other cases, they will either not close or will close only after a long period. As a result, a significant amount of gas fuel will be wasted into the environment, and consumers will be deprived of gas supply for a long time. Such operation of gas pipelines under current conditions is unacceptable. The first reason is the increasing demand for energy and the relevance of energy conservation due to energy shortages. On the other hand, the significant damage caused to consumers due to the prolonged lack of gas supply is unacceptable in today's market economy. Therefore, there is a great need for the development of a new concept for the operation of emergency automatic valves in the design and reconstruction of complex gas pipelines.

As mentioned above, the use of pressure impulse for calibration of automatic valves at each installed valve's coordinates, based on pressure changes over time, does not meet the modern operational requirements for gas pipelines. The parameter **β**, which characterizes the rate of pressure drop, will have different values for each gas pipeline according to its parameters, based on equation (18). This means that when the integrity of high-pressure gas pipelines is compromised, the pressure will drop to the required value (20%) within a short period. In low-pressure gas pipelines, the activation time of the automatic valves will be longer. Moreover, the ratio of pressures between the beginning and end of a gas pipeline can differ by 2 to 3 times depending on the pipeline's length. Consequently, when the integrity is compromised at the beginning of the pipeline, the pressure drop rate at the end will be several times lower compared to the beginning, leading to prolonged valve closure and undesired gas loss. Conversely, if a pipeline break occurs at the end of the pipeline, the pressure drop rate at the end will be more than twice the value at the beginning (see Table 1.2.1, Graph 1.2.1). To avoid these issues, it is crucial to determine the value of **β** for each gas pipeline with specific parameters. Analytically determining **β** involves sufficiently complex limits (infinite sums) as per equation (18). Therefore, using predefined theoretical and experimental values, equation (18) can be expressed in an empirical form.

**Analysis of Theoretical and Experimental Data**

The analysis of theoretical and experimental data is of particular importance for gas-dynamic processes. By using the values of certain function parameters obtained theoretically, simple mathematical expressions can be derived through the solution of the inverse problem. This leads to the determination of the sought parameter. Recently, the interest in solving inverse problems, the collection of more data, and the existence of advanced digital equipment and software for processing large volumes of theoretical data efficiently and rapidly are all related to the development of numerical devices and modern reporting programs.

The most commonly used method for parameter determination is the least squares method or minimizing the least squares function. The drawbacks of applying the least squares method in nonlinear gas flow stimulate the search for and development of linear approaches for parameter estimation. First, we use a simple yet common approach by linearizing the problem with the use of the inverse function. This approach is generally limited to relatively simple problems where the sought parameters appear as arguments of the function. Moreover, this approach is applied when the argument of the multidimensional measured function is the only monotonically dependent argument. Thus, using the data in Table 1.2.1 and employing methods for researching theoretical and experimental data, we obtain the following empirical equation.

$$P(L,t) = P(L,0) - P(L,0)\left[1 - \left(1{,}1 - \frac{G_{ut}}{10G_0}\right)\right]\exp\left(-\frac{P_2}{1{,}3P_1}\beta t\right) \quad (19)$$

Using the parameters provided above and the empirical formula (19), we determine the time-dependent values of pressure changes at the end of the gas pipeline. We then compare these values with those obtained from solving the heat transfer equation, as given by expression (18), and listed in Table 1.2.1 (see Table 1.2.2).

| t, sec | P(L,t), $10^4$ Pa. | | | | | | | | |
|---|---|---|---|---|---|---|---|---|---|
| | $\ell = 30000m$, | | | $\ell = 40000m$, | | | $\ell = 70000m$, | | |
| | Values Calculated Using Equation (18) | Values Calculated Using Equation (19) | Relative Error of Values Calculated Using Equations (18) and (19) | Values Calculated Using Equation (18) | Values Calculated Using Equation (19) | Relative Error of Values Calculated Using Equations (18) and (19) | Values Calculated Using Equation (18) | Values Calculated Using Equation (19) | Relative Error of Values Calculated Using Equations (18) and (19) |
| | $G_{ut} = 0{,}1G_0$ | | | $G_{ut} = 0{,}1G_0$ | | | $G_{ut} = 0{,}1G_0$ | | |
| 250 | 24,87 | 25,19 | 0,0128 | 24,09 | 24,33 | 0,0097 | 23,23 | 22,70 | 0,0226 |
| 300 | 24,66 | 24,79 | 0,0053 | 23,66 | 23,78 | 0,0052 | 22,52 | 21,89 | 0,0281 |
| 350 | 24,50 | 24,41 | 0,0039 | 23,28 | 23,25 | 0,0015 | 22,11 | 21,11 | 0,0454 |
| 400 | 24,38 | 24,02 | 0,0147 | 22,96 | 22,72 | 0,0103 | 21,69 | 20,35 | 0,0616 |
| 450 | 24,30 | 23,65 | 0,0269 | 22,69 | 22,21 | 0,0211 | 21,35 | 19,62 | 0,0810 |
| 500 | 24,26 | 23,28 | 0,0404 | 22,48 | 21,72 | 0,0339 | 21,09 | 18,92 | 0,1032 |
| t, sec | $G_{ut} = 0{,}3G_0$ | | | $G_{ut} = 0{,}3G_0$ | | | $G_{ut} = 0{,}3G_0$ | | |
| 250 | 24,83 | 24,72 | 0,0041 | 23,98 | 23,88 | 0,0042 | 22,71 | 22,29 | 0,0187 |
| 300 | 24,21 | 24,34 | 0,0051 | 23,08 | 23,34 | 0,0113 | 21,78 | 21,49 | 0,0133 |
| 350 | 23,72 | 23,96 | 0,0101 | 22,35 | 22,82 | 0,0208 | 21,06 | 20,72 | 0,0160 |
| 400 | 23,32 | 23,58 | 0,0113 | 21,77 | 22,31 | 0,0246 | 20,53 | 19,98 | 0,0268 |
| 450 | 23,01 | 23,21 | 0,0089 | 21,33 | 21,81 | 0,0226 | 20,16 | 19,26 | 0,0447 |
| 500 | 22,78 | 22,85 | 0,0030 | 20,99 | 21,32 | 0,0154 | 19,94 | 18,57 | 0,0687 |

| t, sec | $G_{ut} = 0,5G_0$ | | | $G_{ut} = 0,5G_0$ | | | $G_{ut} = 0,5G_0$ | | |
|---|---|---|---|---|---|---|---|---|---|
| 250 | 24,61 | 24,26 | 0,0139 | 23,87 | 23,43 | 0,0181 | 22,26 | 21,87 | 0,0173 |
| 300 | 23,76 | 23,88 | 0,0050 | 22,76 | 22,91 | 0,0063 | 21,18 | 21,09 | 0,0046 |
| 350 | 23,12 | 23,51 | 0,0167 | 21,96 | 22,39 | 0,0196 | 20,46 | 20,33 | 0,0063 |
| 400 | 22,66 | 23,14 | 0,0215 | 21,41 | 21,89 | 0,0224 | 20,02 | 19,60 | 0,0210 |
| 450 | 22,33 | 22,78 | 0,0202 | 21,06 | 21,40 | 0,0159 | 19,81 | 18,90 | 0,0461 |
| 500 | 22,13 | 22,42 | 0,0133 | 20,88 | 20,92 | 0,0017 | 19,78 | 18,22 | 0,0786 |

**Table 1.2.2. Comparison of Theoretical and Empirical Pressure Values at the End of the Complex Gas Pipeline Over Time**

From Table 1.2.2, it is evident that the relative error between the theoretical values calculated using function (18) and the values obtained using the empirical formula (19) for pressure at the end of the gas pipeline over time ranges from 0.0030 to 0.08, except for point 1 (highlighted in red in Table 1.2.2). The presence of relative error is due to the linearization and simplification of the complex mathematical expression obtained from the nonlinear system of equations.

Now, let us examine the graphical representation of the theoretical function (18) and the empirical function (19). For this purpose, using Table 1.2.2, we prepare the following graph (Graph 2), which illustrates the time-dependent pressure values at the end of the complex gas pipeline for the case where the rupture location is at a distance of $\ell = 40000m$. This graph is based on the theoretical values and those calculated using the empirical formula derived from experimental data (Graph 1.2.2)

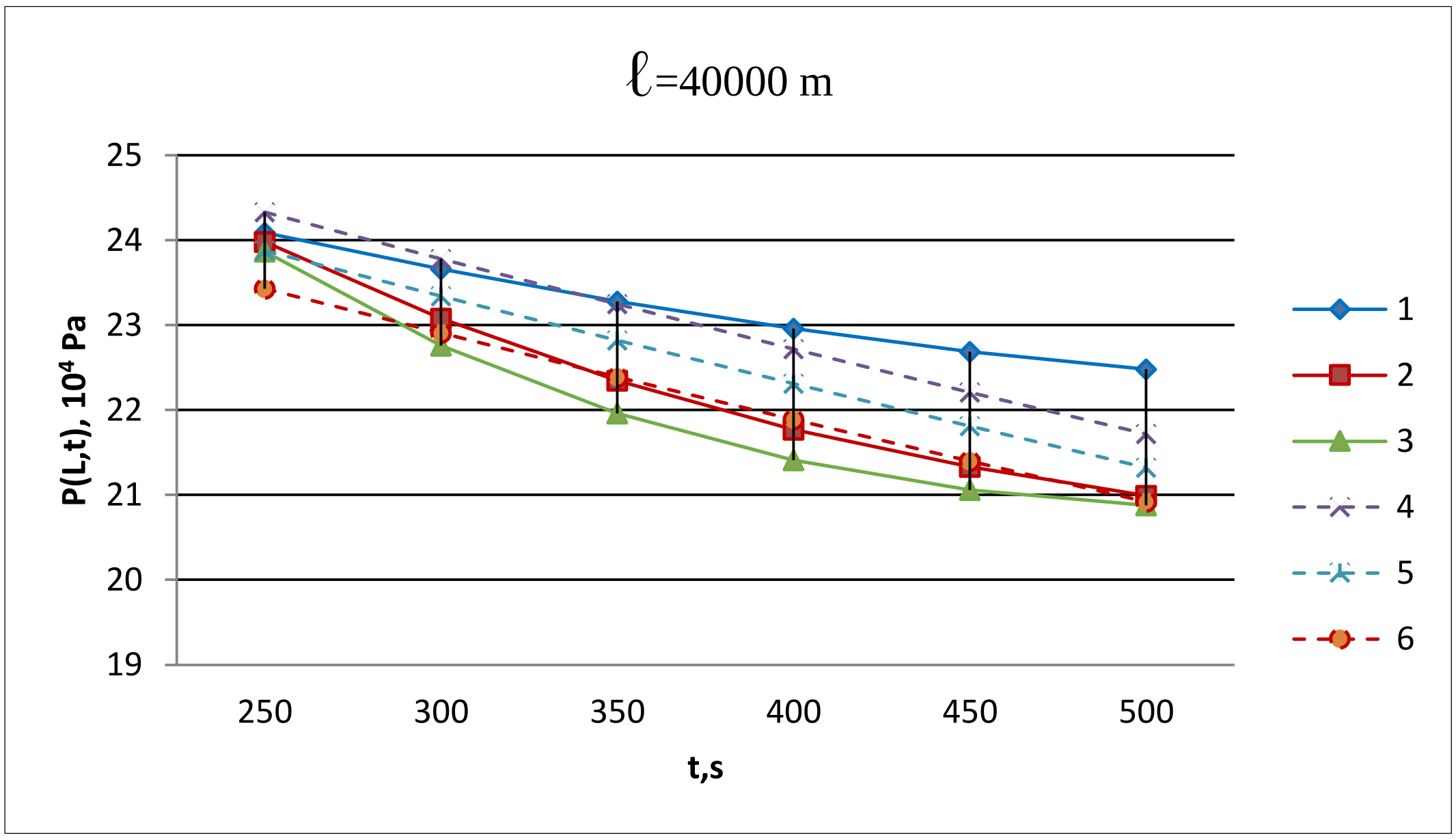

**Graph 1.2.2. Comparison of Theoretical and Experimental Pressure Values Over Time at the End of the Complex Gas Pipeline**

Analysis of the data in Graph 1.2.2 leads to the conclusion that the pressure changes when the gas pipeline ruptures follow an exponential dependency. Thus, both the theoretical and experimental functions are consistent with each other in terms of both their values and graphical representation. This indicates that the empirical equation (19) is suitable for studying physical processes and determining the sought parameter.

Therefore, using equation (19), we obtain the following formula to determine the parameter β, which characterizes the pressure drop rate in gas pipelines at time $t=t_1$ :

$$\beta = \frac{P_2}{P_1 t_1} \ell n \left[ \frac{\varphi}{\psi} \right] \quad (20)$$

Here, $\varphi = 1 - \frac{P(L,t_1)}{P(L,0)}$ ; $\psi = 0,1\left(1 - \frac{G_{ut}}{G_0}\right)$

To verify the accuracy of the values obtained for β based on the empirical formula derived from the experimental data, we use Table 1.2.2 to determine the value of β in the final section of the gas pipeline using equation (19). Next, using the determined value of β and the data provided for the gas pipeline, we calculate the operation time of the automatic valves at different leak locations and for different gas masses using equation (20). Based on these indicators, we determine the calibration parameter value for adjusting the automatic valve actuators. The results are summarized in the following table (Table 1.2.3).

| $\frac{G_{ut}}{G_0}$ | $\beta$ , $1/\text{sec}$ | $\ell = 30000m$ | | $\ell = 40000m$ | | $\ell = 70000m$ | |
|---|---|---|---|---|---|---|---|
| | | $t_1$, s | Calibration Indicator (%) | $t_1$, s | Calibration Indicator (%) | $t_1$, s | Calibration Indicator (%) |
| 0,1 | 0,0057 | 224 | 28 | 226 | 28 | 237 | 26 |
| 0,3 | 0,0075 | 164 | 29 | 168 | 28 | 175 | 27 |
| 0,5 | 0,0099 | 122 | 29 | 125 | 28 | 131 | 27 |

Table 1.2.3. Values of the activation time and calibration indicators of the isolation valves based on different quantities and locations of leaking gas.

The value of β shown in Table 1.2.3 represents the pressure drop rate at time $t=t_1$ in non-stationary gas flow. Therefore, the values may not necessarily match the values shown in Table 1.2.1 for different locations of the isolation valves. This is because, based on the problem setup, the values noted in Table 1.2.1 are average values of the gas pressure drop rate.

Now, considering the installation step of the automatic valves on the gas pipeline as 5000 meters, we determine the operating time of the valves every 5000 meters using the formula (17) at different values of pressure drop rate in non-stationary gas flow. We then create the following graph (Graph 1.2.3).

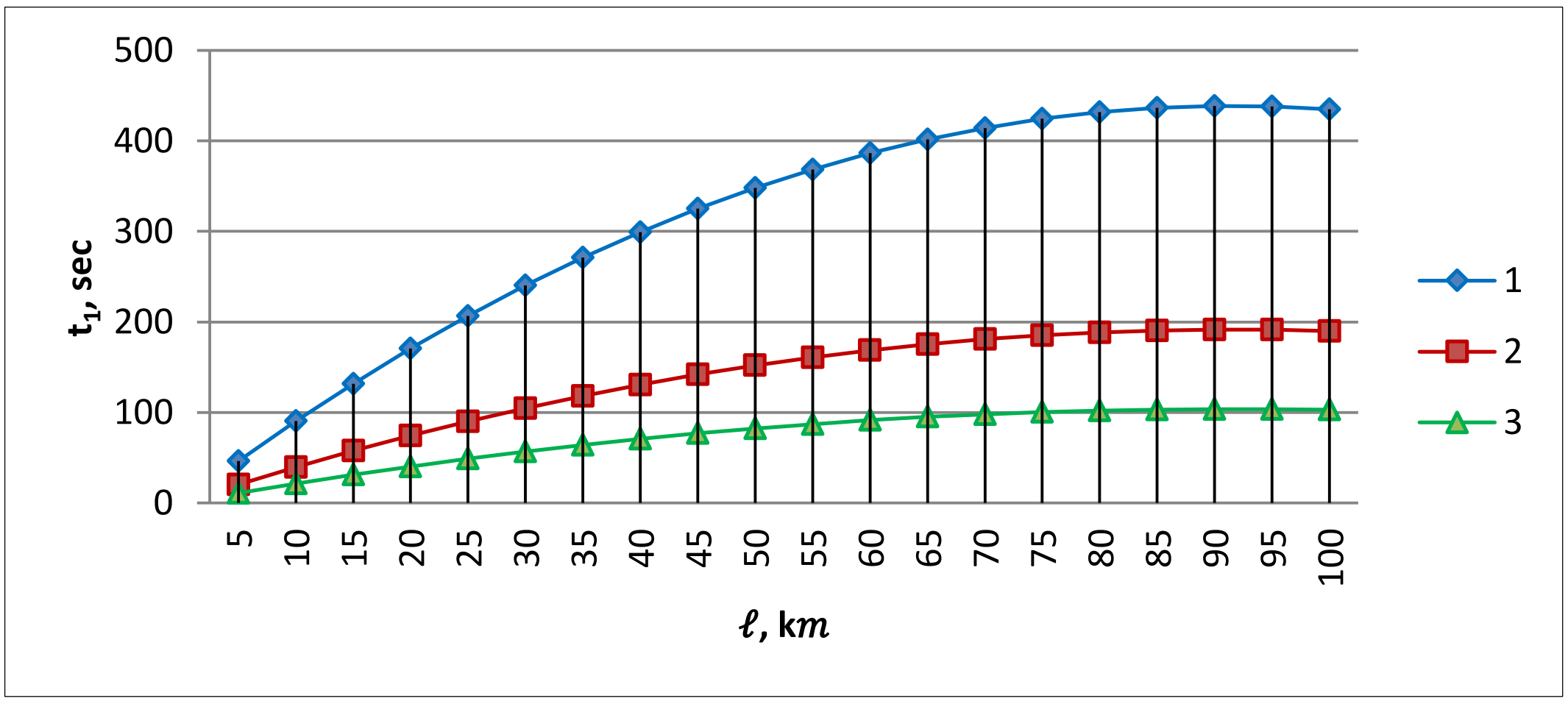

**Graph 1.2.3. The variation of activation time of the isolation valves depending on the installation location of the valves for different quantities of leaking gas. (1-$G_{ut}$=0,1$G_0$; 2-$G_{ut}$=0,3$G_0$; 3-$G_{ut}$=0,3$G_0$ ).**

Based on the analysis of Figure 1.2.3, we can state that when the leaking gas constitutes 10% of the gas pipeline's initial flow rate, the activation time of the isolation valves ranges from 47 to 435 seconds. If the leaking gas constitutes 30% of the flow rate, the activation time ranges from 20 to 190 seconds. For 50% of the flow rate, the activation time ranges from 11 to 103 seconds. The variation slightly increases from the beginning to the end of the gas pipeline, which aligns with the requirements of the gas dynamic process. It is known that as you move from the beginning to the end of the pipeline, the hydraulic resistance affects the pressure of the gas flow, causing the pressure drop rate to decrease accordingly. The analysis also indicates that in the sections closer to the end of the gas pipeline ($\ell = 85000 \div 100000m$), this increase stabilizes. Therefore, this variability cannot be considered a reliable basis for the operational management algorithm of gas transportation systems based on changes in technological parameters. The duration for the valve activation time should be set in such a way that it is sufficient for recording the leaking gas mass in a non-stationary gas flow at the dispatcher's station.

An important finding from the reports is that the calibration indicator value for adjusting the valves, regardless of the valve's location, the leaking gas mass, or the leak point, is almost the same (Table 1.2.3). The analysis showed that the adjustment of the automatic valves in emergency situations may vary for each gas pipeline with specific parameters. Based on the proposed scheme and the gas pipeline parameters, the calibration indicator for the automatic valves is set at an average of 28% according to Table 1.2.3. This means that in the event of an emergency in the gas pipeline, if the stationary pressure drop rate in the section where the valves are located decreases by 28%, the automatic valves on either side of the leak point in the pipeline will automatically close. From this moment, the pressure at the beginning of the damaged pipeline will start to increase, while the pressure at the end will start to decrease. This gas dynamic process is recorded at the dispatcher's station from that moment onward, and based on the research conducted in [5], the automatic valves on the left side of the closed section and the automatic valves on the connecting pipeline on the right side will be remotely controlled to open (Scheme 1.1.1)

The study of these processes and the consideration of the results of the analytical solutions in dispatcher control stations are crucial for the reconstruction of existing complex gas pipelines. This approach helps prevent large volumes of gas from leaking into the atmosphere and ensures the continuous supply of gas to consumers.

## Conclusion

The solution of the differential equations system describing the non-stationary gas flow along the length of a complex gas pipeline has been obtained, considering the operating principles of automatic shut-off valves. The boundary conditions included in the solution, related to Heaviside functions, allow for compact expressions of the mathematical model, but do not restrict the application scope of the model to changes in gas flow or pressure drop rates. Thus, the methodology for calculating the operational time of automatic valves and the conditions for their closure under non-stationary gas flow conditions has been developed.

An empirical calculation formula has been developed based on the analysis of theoretical and experimental data to determine the rate of pressure drop for managing and regulating gas transmission in emergency modes. The calculation results showed that the error in comparing the obtained values of algorithms does not exceed 8%. This method improves the calibration of automatic shut-off valves in gas transmission systems.

An analytical formula for determining the shutdown time and pressure drop rate of automatic shut-off valves in complex pipelines has been developed. The proposed methods can be widely used to improve the management of main pipeline systems. Using the developed algorithms and calculation methods, taking into account the non-stationary regime, effective solutions for the reconstruction of complex main gas pipelines can be achieved.

### 1.3. Methodology for Selecting Efficient Operating Conditions for Non-Stationary Gas Transmission During the Reconstruction Phase

The selection of optimal strategies for reconstruction, primarily ensuring the efficient exploitation regime of gas pipelines, and enhancing the reliability of elements installed in the pipeline's section under operating conditions are associated with minimizing the potential damage from accidents to an acceptable level. To address such issues, it is necessary to develop comprehensive control schemes for gas pipeline operations and utilize modern tools.

Control over the operation of main pipelines is primarily implemented through the calculation of transient processes during accidents. This allows for the calculation of technological alternatives and the making management decisions before the entire hydraulic system is fully deployed. To achieve this, it is necessary to use technical tools based on the simplest principles of exponential form, which have undergone theoretical and experimental tests in known areas of application. One of the significant challenges encountered during the calculation of gas dynamics issues is determining the appropriate method for solving processes that undergo changes depending on time and coordinates. This type of solution should be characteristic of modern applications of science, engineering, and technology problems. The establishment and selection of methods for constructing and analyzing different schemes of gas dynamics are crucial for solving the given problem. The main objective of the research is to determine and select optimal operating regimes during the system's reconstruction stage.

**Analysis and Generalized Methods of Existing Research:**

For various technological conditions, especially those limited to a specific scope, individual researchers have developed specific methods for the operational management of pipelines [24,42,61,65]. Certain areas of gas dynamics have been studied for quite some time and are undergoing significant development intensively. As a result, many important and interesting findings have been obtained. However, general methods for solving gas dynamic issues are still not available. Moreover, it should be noted that the existence and uniqueness of a general solution have not yet been proven. This is attributed to the complexity of gas dynamics equations and, above all, their non-linearity. Considering such events complicates the mathematical formalization of the issues and leads to the emergence of independent problems in their solution. However, the essence of the problem still lies in the classical equations of gas dynamics. Therefore, developing and preparing effective methods for solving these equations is an important issue for many areas of modern science.

Practically applying reliable and effective schemes and algorithms that have been tested in experiments to solve complex problems and utilizing a variety of programs to implement numerical algorithms on computers should be adapted to the requirements of today. Modern software allows for the calculation of parameters of transported gas in real-time by considering changes occurring in its physical processes [12].

Learning the approximation of the difference scheme (estimated calculations) for piecewise functions is not difficult in both linear and non-linear cases. Proving the stability of the scheme actually leads to obtaining certain apriori evaluations that express the continuous dependence of the solution on the initial data of the problem. Unlike linear situations, establishing such evaluations for non-linear equations is associated with significant challenges, and such evaluations are completely absent for gas dynamics equations. Therefore, the verification of scheme stability is usually carried out in some linear analogs of the initial problem.

It should be noted that for the class of one-dimensional non-stationary gas dynamics problems under consideration, other numerical solution methods also exist, such as the method of characteristics [30], the method of discontinuous decomposition [23], the "large particles" method, and so on.
A detailed description of these can be found in the referenced works. Instead of exploring all possible methods, we will implement boundary conditions, one of the simplest approaches, to address the issue. The majority of research in the field of gas pipeline dynamics is associated with the method of solving gas dynamic partial differential equations with various initial and boundary conditions.

It is possible to obtain solutions to the problem of determining the distribution of pressures and flow rates along the pipeline. Non-stationary gas flows in the main pipeline system are described by a system of non-linear partial differential equations. One of the captivating areas of research is the use of analytical methods [32].

### Analysis of Nonlinear Differential Equation Systems Using Analytical Methods.

The system of nonlinear differential equations typically relies on simplifying the initial equations and their subsequent manipulation. However, the difficulty of analytical solutions doesn't necessarily outweigh the benefits of simpler methods. This approach's disadvantage lies in its reliance on external similarity to significantly simplify the equations. Simplification is not only employed to simplify the initial equations for gas flow but also during the modeling of other physical processes, which is not always justified. On the other hand, the possibility of visually observing variable pressure fields and gas flow consumption in a complex network of gas pipelines stimulates research efforts. In this regard, we adopt the scheme of the following parallel gas pipeline, which accommodates various operating conditions (Scheme 3.1.1) .

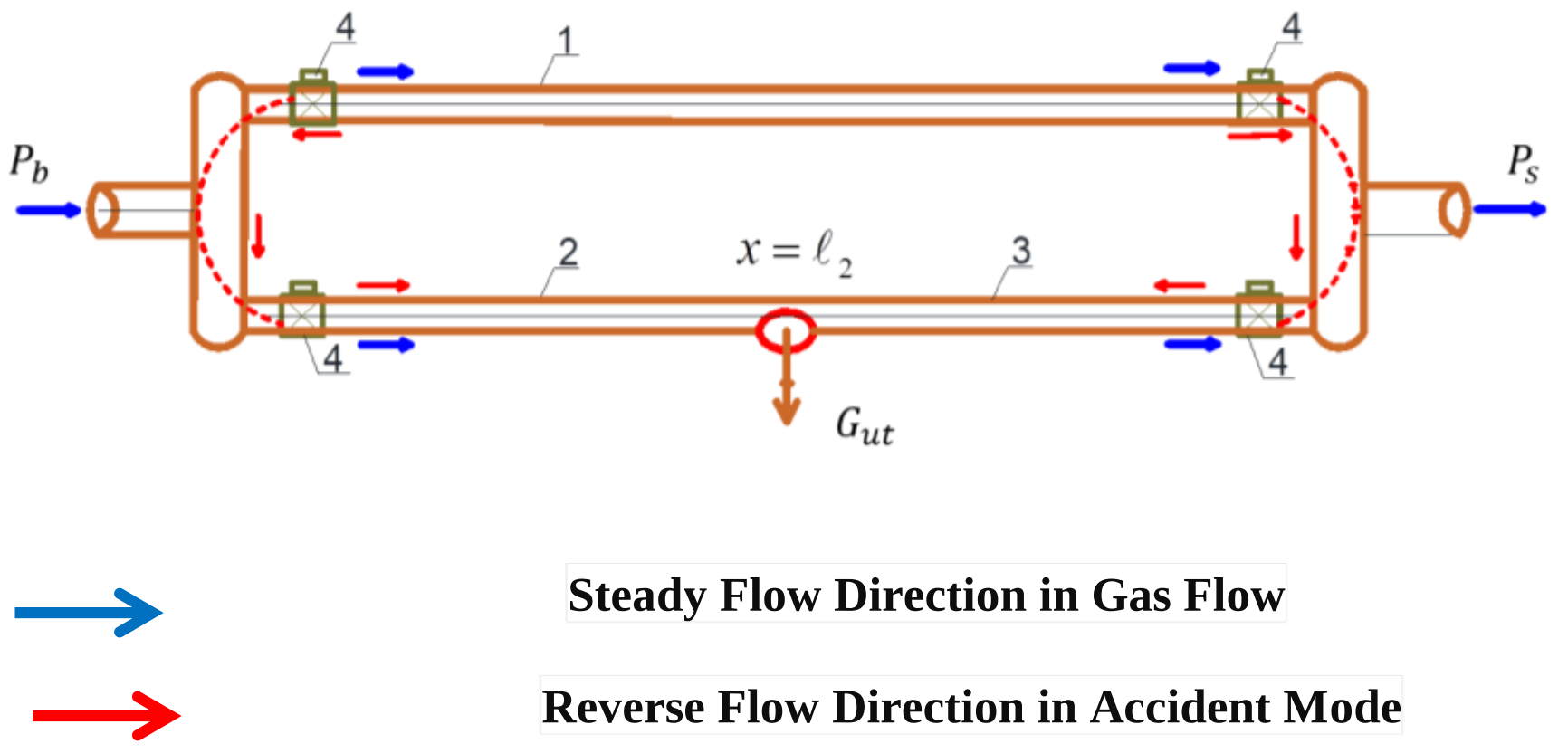

**Scheme 3.1.1. The prinsipial scheme of a parallel gas pipeline operating in unfieldt hydraulic regime.**

1- **Intact pipeline section of the gas pipeline;**
2- **Portion of the damaged pipeline section from the starting point of the damaged pipeline to the leakage point ($x = \ell_2$);**

3- **Portion of the damaged pipeline section from the leakage point to the end point; $P_b$ - Pressure at the starting point of the gas pipeline in steady-state regime; $P_s$ - Pressure at the end point of the gas pipeline in steady-state regime; $G_{ut}$ - Mass flow rate of gas leakage, $\frac{Pa\cdot sec}{m}$.**

The valve fittings (such as valves and taps) installed on parallel gas pipelines are positioned at a certain distance, taking into account the ease of their technical operation and repair as well as the requirements of gas supply facilities [8]. Let's assume that an accident occurs at point $x=\ell_2$ on one of the lines of the parallel gas pipeline due to a rupture in the pipeline. As shown in Scheme 1.3.1, both the damaged and intact sections of the gas pipelines will operate under non-steady-state gas flow conditions. It should be noted that selecting the parameters based on the operation of the main gas pipelines is required for calculating the occurrence of non-stationary processes caused by any accident situation [22]. This requirement should be considered one of the highest priorities when selecting the operational regime of the gas pipeline.Technologically-based information has been analyzed to determine the physical laws governing accident situations in parallel pipeline transportation processes. For any technological scheme, there are certain parameters and indicators that characterize it. All the technological issues mentioned above are resolved by mathematical modeling of gas flow parameters related to the solution of the following system of differential equations [8].
To calculate the gas flow in the pipeline, we utilize the dimensionless continuity and momentum equations recommended by I.A. Charney in the literature cited in the first reference [32,78].

$$\begin{cases} \dfrac{\partial P(x,t)}{\partial x} = -2aG(x,t) \\ \dfrac{\partial G(x,t)}{\partial x} = -\dfrac{1}{c^2}\dfrac{\partial P(x,t)}{\partial t} \end{cases} \tag{1}$$

Here, $2a=\lambda\frac{V}{2d}$; - I.A. Charney theorem; $G=\rho v$; $P=\rho\cdot c^2$;
c- Sound propagation speed of gas for an isothermal process, m/sec.;
d - The internal diameter of the pipe, m;
x- A coordinate aligned with the gas flow and coinciding with the axis of the pipe., m;
P -Absolute average gas pressure in sections, Pa;
v - Average gas velocity at the cross-section, m/sec;
ρ- Average density of the gas, κq/ $m^3$ ;
t -time, sec;
λ- Hydraulic resistance coefficient of the gas pipeline section, dimensionless;
G- Mass flow rate of the gas flow, $\frac{Pa\cdot sec}{m}$;

If we differentiate the first equation of the system with respect to x, and substitute $\frac{\partial G}{\partial X}$ in the second equation, then after taking the Laplace transform, we obtain the heat transfer equation for the mathematical solution of the problem.

$$\frac{d^2P_i(x,S)}{dx^2} = \frac{2a}{c^2}\left[SP_i(x,S) - P_i(x,0)\right] \text{ : i=1,2,3} \tag{2}$$

Using this equation, we need to determine boundary conditions to properly formulate the problems. The set of boundary and initial conditions establishes the connection between models, and the proper formulation of boundary conditions is often achieved through the modification of at least a general source of information database.

*Analysis of non-steady-state gas dynamic processes under operating conditions with fixed pressure maintained at the ends of parallel pipeline lines.*

Let's assume that a fixed pressure is maintained at the end of the parallel gas pipeline. Initial conditions are homogeneous for the system under consideration. For solving the mathematical model of equation (2), we accept the initial condition of the initial stationary distribution of gas flow pressures in the calculation domain.

$t=0;\ P_i(x,0)=P_b-2aG_0x:\ i=1,2,3$

Using the obtained general model, the boundaries of the calculation domain, conditions for pressure variation at the boundaries, and changes in gas flow are determined. Considering the fixed pressure at the end of the gas pipeline, we take into account pressure functions that describe the flow rate G(t) and its abrupt changes at the beginning of the pipeline. Thus, the following boundary conditions are accepted for solving equation (2):

At the point x=0 $\begin{cases} P_1(x,t)=P_2(x,t) \\ \dfrac{\partial P_1(x,t)}{\partial X}+\dfrac{\partial P_2(x,t)}{\partial X}=-2aG_s(t) \end{cases}$

At the point $x=\ell_2$ $\begin{cases} P_2(x,t)=P_3(x,t) \\ \dfrac{\partial P_3(x,t)}{\partial X}-\dfrac{\partial P_2(x,t)}{\partial X}=-2aG_{ut}(t) \end{cases}$

At the point x= L $\begin{cases} P_1(x,t)=P_s \\ P_3(x,t)=P_s \end{cases}$

In this case, the general form of the solution for equation (2) will be as follows:

$$P_1(x,s)=\frac{P_i-2aG_0x}{S}+c_1Sh\lambda x+c_2Ch\lambda x \qquad 0\le x\le L \tag{3}$$

$$P_2(x,S)=\frac{P_i-2aG_0x}{S}+c_3Sh\lambda x+c_4Ch\lambda y \qquad 0\le x\le \ell_2 \tag{4}$$

$$P_3(x,S)=\frac{P_i-2aG_0x}{S}+c_5Sh\lambda x+c_6Ch\lambda x \qquad \ell_2\le x\le L \tag{5}$$

Here, $\lambda=\sqrt{\dfrac{2as}{c^2}}$ , $\beta=\sqrt{\dfrac{2ac^2}{s}}$

(3), (4), and (5) functions represent the determination of pressure distribution along the axial sections of the gas pipeline in a general form. To assign physical meaning to these equations, it is necessary to determine the constants $c_1$, $c_2$, $c_3$, $c_4$, $c_5$ and $c_6$. By using the Gauss method to determine the mentioned constants, and substituting them into functions (3), (4), and (5), and after performing simple operations, we obtain the following expressions to determine the pressure distribution in converted form:

$$0\le x\le L$$

$$P_1(x,S)=\frac{P_1(x,0)}{S}+\frac{\beta}{2}\left[G_0(S)-\frac{G_0}{S}\right]\frac{sh\lambda(L-x)}{ch\lambda L}-$$
$$-\beta G_{ut}(S)\frac{sh\lambda(L-\ell_2)sh\lambda(L-x)}{sh\lambda 2L} \qquad (6)$$

$$0\le x\le \ell_2$$

$$P_2(x,S)=\frac{P_2(x,0)}{S}+\frac{\beta}{2}\left[G_0(S)-\frac{G_0}{S}\right]\frac{sh\lambda(L+x)}{ch\lambda L}-,$$
$$-\beta G_{ut}(S)\frac{sh\lambda(L-\ell_2)sh\lambda(L-x)}{sh\lambda 2L} \qquad (7)$$

$$\ell_2\le x\le L$$

$$P_3(x,S)=\frac{P_3(x,0)}{S}+\frac{\beta}{2}\left[G_0(S)-\frac{G_0}{S}\right]\frac{sh\lambda(L+x)}{ch\lambda L}-$$
$$-\beta G_{ut}(S)\frac{sh\lambda(L+\ell_2)sh\lambda(L-x)}{sh\lambda 2L} \qquad (8)$$

The original general solutions for functions (6), (7), and (8) will be written as follows:

$$P_1(x,t)=P_s+8aG_oL\sum_{n=1}^{\infty}\frac{e^{-\alpha_1 t}}{[\pi(2n-1)]^2}Cos\frac{\pi(2n-1)x}{2L}+$$
$$+\frac{c^2}{L}\sum_{n=1}^{\infty}Cos\frac{\pi(2n-1)x}{2L}\int_0^t[G_0(\tau)]\cdot e^{-\alpha_1(t-\tau)}d\tau- \qquad 0\le x\le L \qquad (9)$$
$$-\frac{c^2}{L}\sum_{n=1}^{\infty}(-1)^n Sin\frac{\pi n(L-\ell_2)}{2L}Sin\frac{\pi n(L-x)}{2L}\int_0^t G_{ut}(\tau)\cdot e^{-\alpha_2(t-\tau)}d\tau$$

$$P_{2,3}(x,t)=P_s+8aG_oL\sum_{n=1}^{\infty}\frac{e^{-\alpha_1 t}}{[\pi(2n-1)]^2}Cos\frac{\pi(2n-1)x}{2L}+$$
$$\frac{c^2}{L}\sum_{n=1}^{\infty}Cos\frac{\pi(2n-1)}{2L}\int_0^t[G_0(\tau)]\cdot e^{-\alpha_1(t-\tau)}d\tau- \qquad (10)$$
$$-\frac{c^2}{L}\sum_{n=1}^{\infty}(-1)^n\int_0^t G_{ut}(\tau)\cdot e^{-\alpha_2(t-\tau)}d\tau\begin{cases}Sin\frac{\pi n(L-\ell_2)}{2L}Sin\frac{\pi n(L+x)}{2L}, 0\le x\le \ell_2\\ Sin\frac{\pi n(L+\ell_2)}{2L}Sin\frac{\pi n(L-x)}{2L}, \ell_2\le x\le L\end{cases}$$

Here, $\alpha_1=\frac{\pi^2(2n-1)^2c^2}{8aL^2}$, $\qquad \alpha_2=\frac{\pi^2n^2c^2}{8aL^2}$

The non-steady-state operating conditions of gas pipelines lead to significant fluctuations in pressure and can compromise the normal functioning of gas supply to consumers. Studying these processes under various operating conditions of gas pipelines and considering the results of the analytical solutions obtained at dispatcher control points will significantly increase the reliability of gas supply and the efficiency of gas transportation processes. This function should be considered one of the highest priorities for the reconstruction of gas pipelines. When selecting the optimal operating conditions for gas pipelines, it should be noted that the task of calculating non-steady-state processes caused by any accident situation is fundamental to the operation of trunk gas pipelines. Selecting parameters for information measurement and control systems for gas distribution networks should minimize the time of transfer to dispatcher control points. During operation, information-measurement and control systems for gas distribution networks should ensure reliable collection and transmission of data and, using modern methods, monitor accident situations and efficiently manage non-steady-state gas dynamic processes to ensure uninterrupted gas supply to consumers.

During accidents in trunk gas pipelines, temporary processes must be determined at least one order of magnitude higher than the calculation time [4]. For this purpose, it is necessary to adopt theoretical and experimental simplification methods based on known areas of application and to use calculation programs based on simple rules. Currently, such programs are available and are specialized simplifications for the operational management of pipeline networks [22].

$$G_0(t)=G_0, \quad G_{ut}(t)=G_{ut}=constat, \ n=1,2,3.....10.$$

Adopting the simplifications mentioned above, we consider the following expressions (9) and (10) and the given data to determine the law of pressure change of gas along the axial length of the complex parallel pipeline in different scenarios of gas leakage locations (1st scenario - $\ell_2 = 2{,}5\times 10^4$ m; 2nd scenario - $\ell_2 = 7{,}5\times 10^4$ m).

$P_b = 55\times10^4$ Pa: $P_s = 4\times10^4$ Pa: $G_0 = G_1+G_2 = 30$ Pa×sec/m: 2a=0.1 1/sec :c= 383,3 m/sec: L = $10^5$ m: d= 0,7 m: 1st scenario $\ell_2 = 2{,}5 \times\cdot 10^4$ m; 2nd scenario $\ell_2 = 7{,}5\times 10^4$ m:

First, considering the fixed pressure at the end of the parallel gas pipeline under operating conditions, we assume a rupture location at a distance of $\ell_2 = 2{,}5\times 10^4$ m from the end, with gas leakage occurring every 300 seconds, and calculate the values of gas pressure along the gas pipeline for different values of leaked gas mass at point x = 25000 m. We record the calculated values in Table 1.3.1 below.

| **$\ell_2$=2500 m** | **$P_1(25000,t)$, $10^4$ Pa** | | | | **$P_2(25000,t)$, $10^4$ Pa** | | | |
|---|---|---|---|---|---|---|---|---|
| t, sec | $G_{ut}$ = 0,4 $G_0$ | $G_{ut}$ = 0,8 $G_0$ | $G_{ut}$ = 1,2 $G_0$ | $G_{ut}$ = 1,6 $G_0$ | $G_{ut}$ = 0,4 $G_0$ | $G_{ut}$ = 0,8 $G_0$ | $G_{ut}$ = 1,2 $G_0$ | $G_{ut}$ = 1,6 $G_0$ |
| 300 | 49,94 | 48,64 | 47,33 | 46,03 | 44,17 | 37,09 | 30,01 | 22,94 |
| 600 | 48,36 | 45,46 | 42,57 | 39,67 | 41,48 | 31,71 | 21,94 | 12,17 |
| 900 | 47,07 | 42,89 | 38,72 | 34,54 | 39,66 | 28,07 | 16,48 | 4,89 |
| 1200 | 46,05 | 40,84 | 35,64 | 30,44 | 38,32 | 25,39 | 12,45 | 0 |
| 1500 | 45,22 | 39,19 | 33,15 | 27,12 | 37,28 | 23,3 | 9,33 | 0 |

Table 1.3.1. Pressure values over time for different leaked gas masses at the point of gas leakage at x = 25000 m along the gas pipeline (1st scenario).

Using Table 1.3.1, let's create a graph (Graph 1.3.1) showing the variation of pressure at the point l=25000 m of the complex gas pipeline over time and depending on the leaked gas mass.

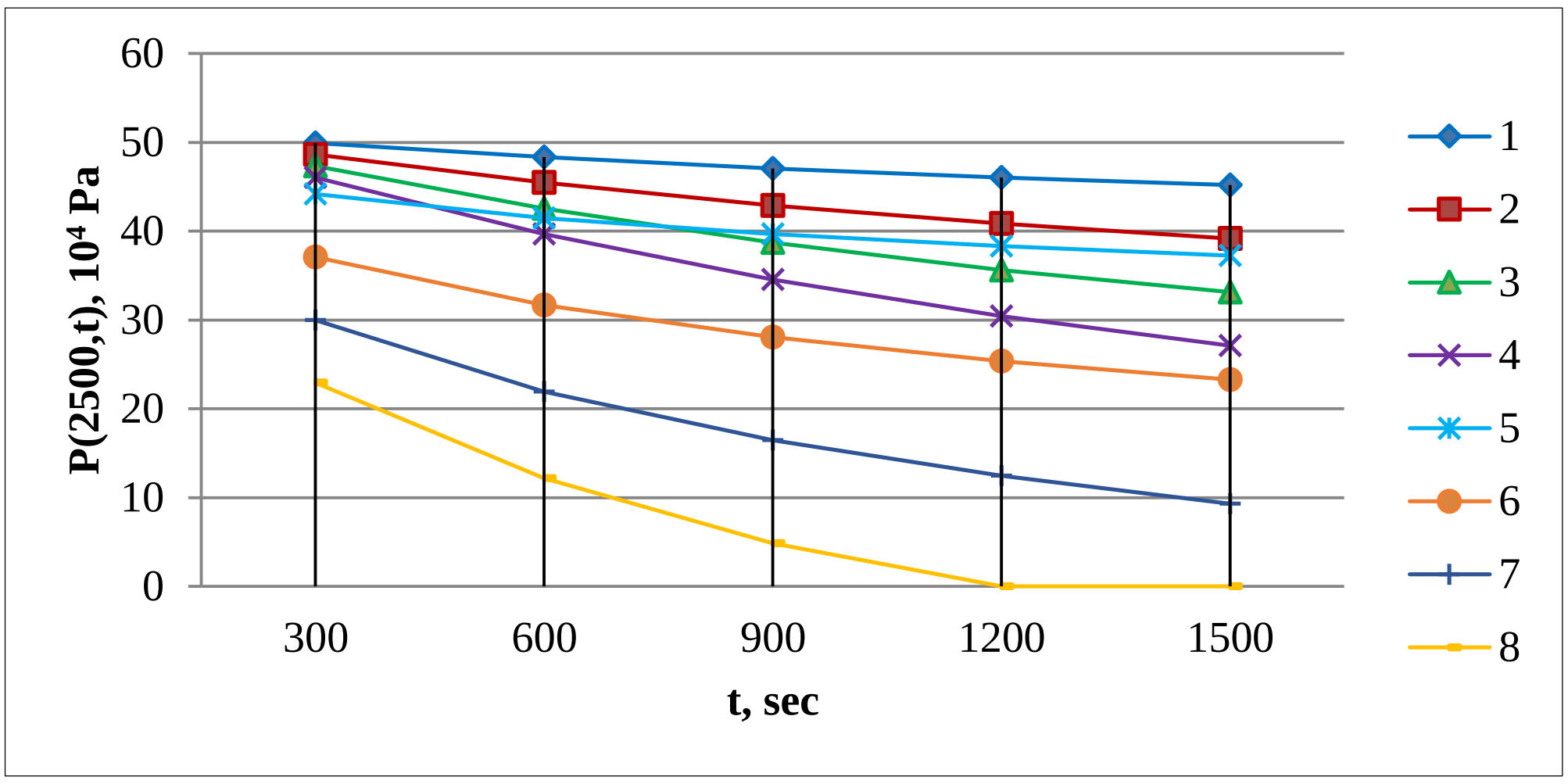

**Graph 1.3.1. Variation of pressure over time at the point x = 25000 m of the gas pipeline in the first scenario of gas leakage. 1, 2, 3, and 4 - Unaffected pipeline section (1-$G_{ut}$=0,4 $G_0$, 2-$G_{ut}$=0,8 $G_0$, 3-$G_{ut}$=1,2 $G_0$, 4-$G_{ut}$ =1,6 $G_0$ ); 5, 6, 7, and 8 - Affected pipeline section (at the point of rupture), (5-$G_{ut}$=0,4 $G_0$, 6-$G_{ut}$=0,8 $G_0$, 7-$G_{ut}$=1,2 $G_0$, 8-$G_{ut}$ =1,6 $G_0$ ).**

Table 1.3. and Graph 1.3.1 demonstrate that the non-steady state caused by the rupture of one pipeline in a parallel gas pipeline operating under unit hydraulic conditions disrupts the steady state of the other unaffected pipeline. This process can be explained by the fact that, in the event of an accident resulting from the rupture of one pipeline, the gas flow of the main sections of the pipeline tends towards the damaged section. In other words, in the operating conditions where the pressure remains constant at the end of the parallel gas pipeline and the mass flow at the beginning changes over time, the mass flows at the beginning of the pipeline are directed towards the damaged section of the pipeline. In this case, the main section of the pipeline operates in a non-stationary gas movement regime due to the extraction of gas from the starting sections of the undamaged pipeline. Consequently, the length of the undamaged pipe, its pressure, and the gas flow consumption vary depending on both coordinates and time due to the extraction of non-uniform gas mass from the starting section and the continuous supply of gas to the consumer institutions from the final section. As a result, a decrease in pressure along the length of both pipelines occurs. Based on the provided data, in the stationary regime, the pressure at x=25000 m for both the damaged and undamaged pipelines is P(25000, 0) =51,25×$10^4$ Pa . Analysis of Graph 1 reveals that the magnitude of pressure decrease in the damaged pipeline is several times greater than that in the undamaged pipeline.

Thus, for a duration of t=300 seconds and when the volume of the leaked gas mass is 40% of the volume of the pipeline's final section ($G_{ut}$ =0,4 $G_0$), the magnitude of pressure decrease in the damaged pipeline is [P(25000,0)-$P_2$ (25000,t)] = 51,25- 44,17)=7,08×$10^4$Pa Meanwhile, the magnitude of pressure decrease in the undamaged pipeline will be [P(25000,0)-$P_1$(25000,t)] =(51,25-49,94)= 1,31×$10^4$ Pa. The value of the pressure drop in the damaged pipeline is 5.4 times greater than the pressure drop occurring in the undamaged pipeline. The value of this ratio varies depending on time for a constant value of gas leakage. For instance, at t=600 seconds, it is 3.37; at t=900 seconds, it is 2.77; at t=1200 seconds, it is 2.48, and at t=1500 seconds, it decreases to 2.32.

It appears from this that the value of this ratio reaches its maximum at t=300 seconds for a constant value of gas leakage ($G_{ut}$ =0,4 $G_0$) and then decreases over time (Graph 1.3.2).

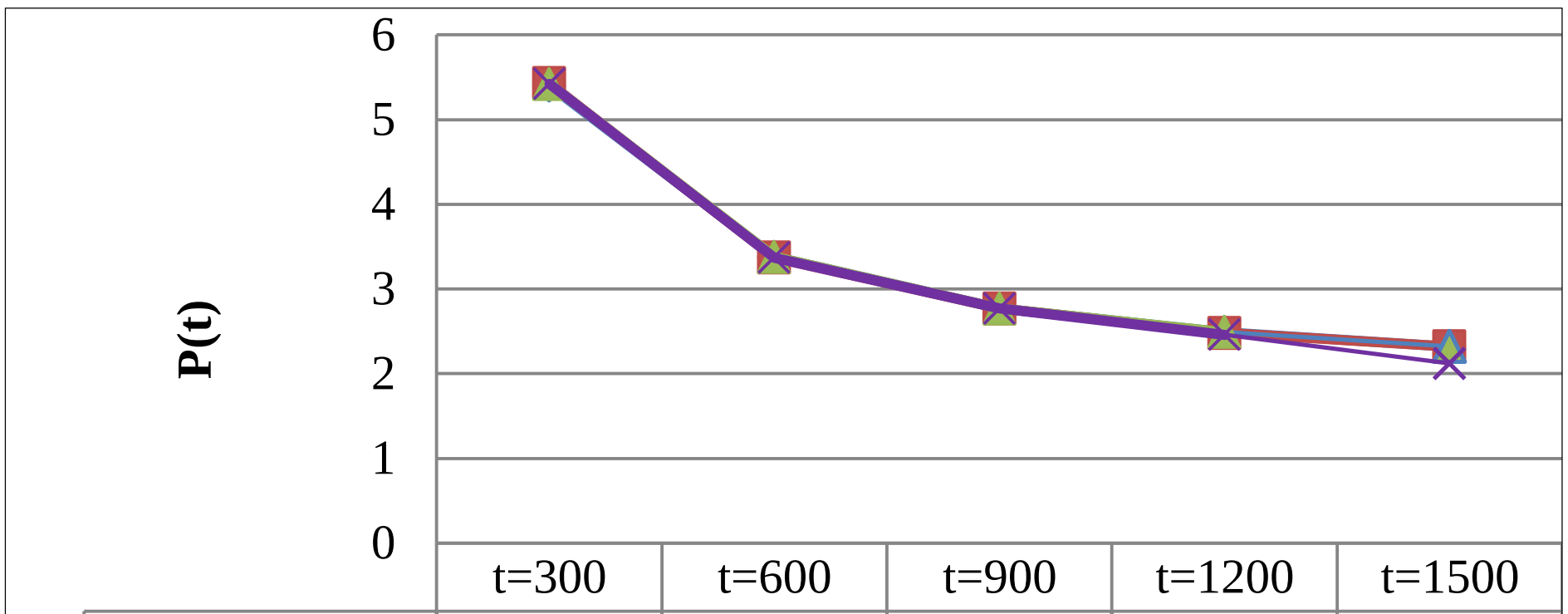

| | t=300 | t=600 | t=900 | t=1200 | t=1500 |
|---|---|---|---|---|---|
| Gut=0,4 Go | 5,4 | 3,38 | 2,77 | 2,49 | 2,32 |
| Gut=0,8 Go | 5,43 | 3,37 | 2,77 | 2,48 | 2,32 |
| Gut=1,2 Go | 5,42 | 3,38 | 2,77 | 2,49 | 2,32 |
| Gut=1,6 Go | 5,42 | 3,37 | 2,77 | 2,46 | 2,12 |

**Graph 1.3.2. The graph of the dependence of [P(x,0)-$P_2$ (x,t)]/[P(x,0)-$P_1$ (x,t)]=p(t) at the point x=25000 m of the gas pipeline on time t.**

The value of the pressure drop in the damaged pipeline is several times greater than the pressure drop in the undamaged one when the volume of leaked gas mass is $G_{ut}$ =0,8 $G_0$Thus, during the period of t = 300 seconds, the value of the pressure drop in the damaged pipeline (51,25- 39,07)= 14,16×$10^4$ Pa, while in the undamaged pipeline, this value will be (51,25-48,64) = 2,61 ×$10^4$ PaHowever, the ratio of the pressure drop in the damaged pipeline to the pressure drop in the undamaged pipeline becomes equal to the value when the gas leakage is $G_{ut}$ =0,4$G_0$, which is 5.4. Despite the gas leakage increasing by a factor of 2, the value of this ratio remains unchanged, only varying with time. As seen from Graph 1.3.2, the value of this ratio decreases only for $G_{ut}$ =1,6 $G_0$ at t = 1500 seconds. The reason for this is that, as shown in Graph 1.3.1, the pressure in the damaged pipeline at x = 25000 m approaches zero after 1200 seconds.

In other words, the process of the gas pipeline being almost completely emptied occurs in the damaged pipeline. Therefore, in the analyzed operating conditions of parallel gas pipelines [$P_1$(L,t) = $P_1$(L,t) = $P_s$ =constat], after t = 1200 seconds, the efficiency of gas supply to consumers is significantly compromised, meaning that most consumers cannot be supplied with gas. The proposed operating condition does not effectively account for uninterrupted gas supply to consumers in emergency regimes, thus rendering the operation of parallel gas pipelines in a single hydraulic regime ineffective within the specified boundary conditions.

Let's analyze the variant where the boundary condition of maintaining constant pressure at the end of the gas pipeline is $\ell_2$=75000 m. In this case, considering different values of the leaking gas mass every 300 seconds, we calculate the values of gas pressure at x = 75000 m and record them in the following Table 1.3.2. From the expression $P_i(x,0)=P_b-2aG_0x$, the value of gas pressure in the stationary regime at x = 75000 m will be 43,75 ×$10^4$ Pa.

| ℓ₂=75000 m | $P_1(75000,t)$, $10^4$ Pa | | | | $P_2(75000,t)$, $10^4$ Pa | | | |
|---|---|---|---|---|---|---|---|---|
| t, sec | $G_{ut}$ = 0,4 $G_0$ | $G_{ut}$ = 0,8 $G_0$ | $G_{ut}$ = 1,2 $G_0$ | $G_{ut}$ = 1,6 $G_0$ | $G_{ut}$ = 0,4 $G_0$ | $G_{ut}$ = 0,8 $G_0$ | $G_{ut}$ = 1,2 $G_0$ | $G_{ut}$ = 1,6 $G_0$ |
| 300 | 43,72 | 43,69 | 43,66 | 43,63 | 37,95 | 32,14 | 26,34 | 20,54 |
| 600 | 43,61 | 43,48 | 43,34 | 43,2 | 36,74 | 29,72 | 22,71 | 15,7 |
| 900 | 43,49 | 43,23 | 42,98 | 42,72 | 36,08 | 28,41 | 20,74 | 13,07 |
| 1200 | 43,38 | 43,01 | 42,64 | 42,27 | 35,65 | 27,55 | 19,45 | 11,36 |
| 1500 | 43,28 | 42,88 | 42,35 | 41,88 | 35,34 | 26,93 | 18,53 | 10,12 |

**Table 1.3.2. Pressure values at the leakage point x = 75000 m of the gas pipeline for different values of leaked gas mass (2nd variant).**

We'll construct the following graph using the data provided in Table 1.3.2 (Graph 1.3.3).

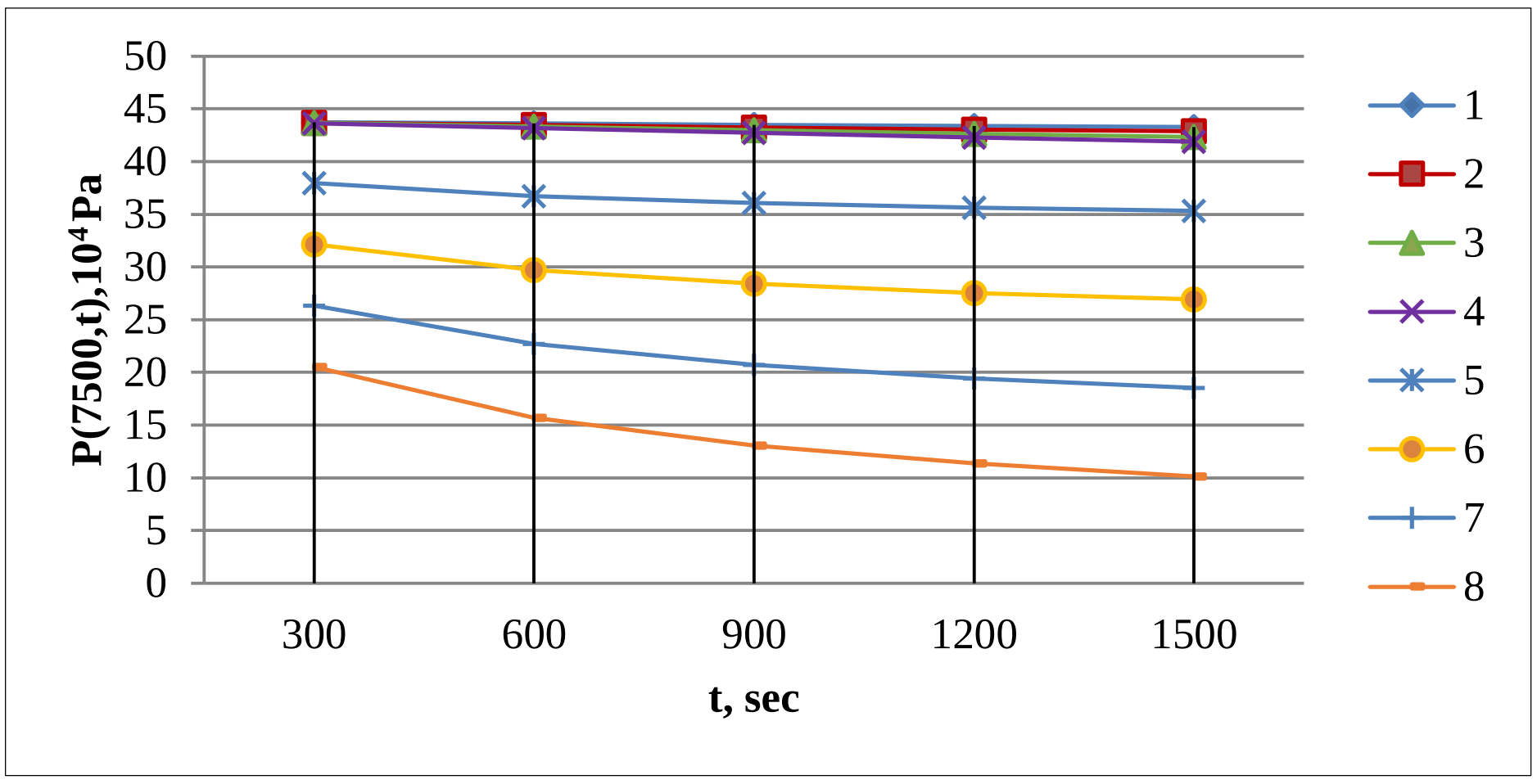

**Graph 1.3.3. The graph illustrates the variation of pressure over time at the point x = 75000 m for the second scenario. In this graph, curves 1, 2, 3, and 4 represent the pressure variation for the undamaged pipeline (1- $G_{ut}$= 0,4 $G_0$, 2-$G_{ut}$=0,8$G_0$, 3-$G_{ut}$=1,2 $G_0$, 4-$G_{ut}$=1,6 $G_0$ ); while curves 5, 6, 7, and 8 represent the pressure variation for the damaged pipeline (at the point of breach) (5-$G_{ut}$=0,4 $G_0$, 6-$G_{ut}$=0,8 $G_0$,7-$G_{ut}$=1,2 $G_0$, 8-$G_{ut}$ =1,6 $G_0$ ).**

Based on the analysis of Graph 1.3.3, it can first be said that the non-stationary process occurring in the damaged pipeline due to the integrity breach in the last part of the belt is not felt in the undamaged pipeline. Specifically, 300 seconds after the accident, the volume of the leaked gas mass is $G_{ut}$ =0,4 $G_0$ . At this time, the pressure drop in the damaged pipeline is (43,75 -37,95)= 5,8 ×$10^4$ Pa, while the pressure drop in the undamaged pipeline will be (43,75 -43,72)= 0,03×$10^4$ Pa. Thus, the value of the pressure drop in the damaged pipeline is 193 times greater than the pressure drop in the undamaged pipeline. This ratio changes over time depending on the steady value of the gas leakage, being 51 at t = 600 seconds, 29.5 at t = 900 seconds, 21.9 at t = 1200 seconds, and decreasing to 18.0 at t = 1500 seconds. Using the data from Table 2, we will construct the graph of the coefficient p(t), which characterizes the ratio of the pressure drop values in the damaged pipeline to the pressure drop values in the undamaged pipeline, depending on different values of gas leakage and time (Graph 1.3.4).

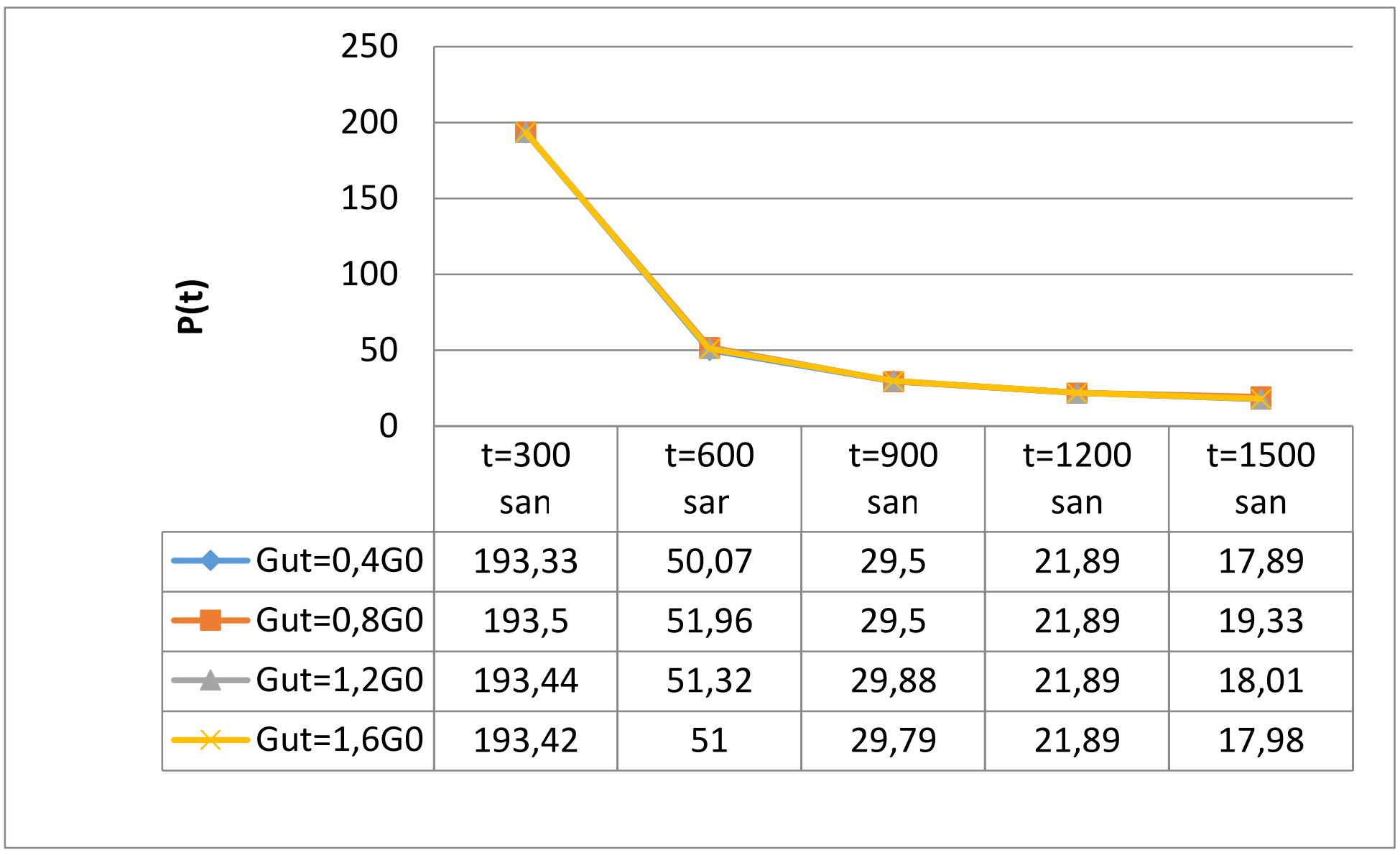

**Graph 1.3.4. Graph of the coefficient [P(x,0)-$P_2$ (x,t)]/[P(x,0)-$P_1$ (x,t)]=p(t) as a function of time *t*, at point x= 25000 on the gas pipeline**.

The analysis of Graph 1.3.4 shows that the values of p(t), which characterize the ratio of the pressure drop in the damaged pipeline to that in the undamaged pipeline, remain nearly constant across different gas leakage rates and over time. This constancy is attributed to the maintenance of stable pressure in the final segment of the operational gas pipeline. However, the mass of gas delivered to consumers decreases sharply as the value of the coefficient describing the ratio of pressure drops between the damaged and undamaged pipeline sections decreases from 193 to 18.

During the rupture of the second line at a distance of $\ell_2$ =25000m and at a gas leakage rate of $G_{ut}$=0.4 $G_0$, the values of pressure distribution in the damaged and undamaged pipeline are calculated at different times, every 25×$10^4$ m, to determine the conformity to the law of pressure distribution as follows (Table 1.3.3):

| $G_{ut}$=0,4 $G_0$ | $P_1$(x,t), $10^4$ Pa | $P_2$(x,t) , $10^4$ Pa | $P_1$(x,t) , $10^4$ Pa | $P_2$(x,t) , $10^4$ Pa | $P_1$(x,t) , $10^4$ Pa | $P_2$(x,t) , $10^4$ Pa |
|---|---|---|---|---|---|---|
| X, $10^4$ m | t=300 sec. | | t=900 sec. | | t=1500 sec. | |
| 0 | 51,95 | 51,95 | 47,96 | 47,96 | 45,67 | 45,67 |
| 2,5 | 49,94 | 44,17 | 39,66 | 41,68 | 45,22 | 37,28 |
| 5 | 46,95 | 44,54 | 41,04 | 42,47 | 43,9 | 39,17 |
| 7,5 | 43,55 | 42,67 | 40,87 | 41,62 | 42,07 | 39,86 |
| 10 | 40 | 40 | 40 | 40 | 40 | 40 |

**Table 1.3.3. For the case of the rupture of the first line of the gas pipeline and a gas leakage rate of $G_{ut}$=0,4 $G_0$, the pressure values along the length of the damaged and undamaged pipeline, depending on time, are as follows.**

Using the data from Table 3, we will create Graph 1.3.5

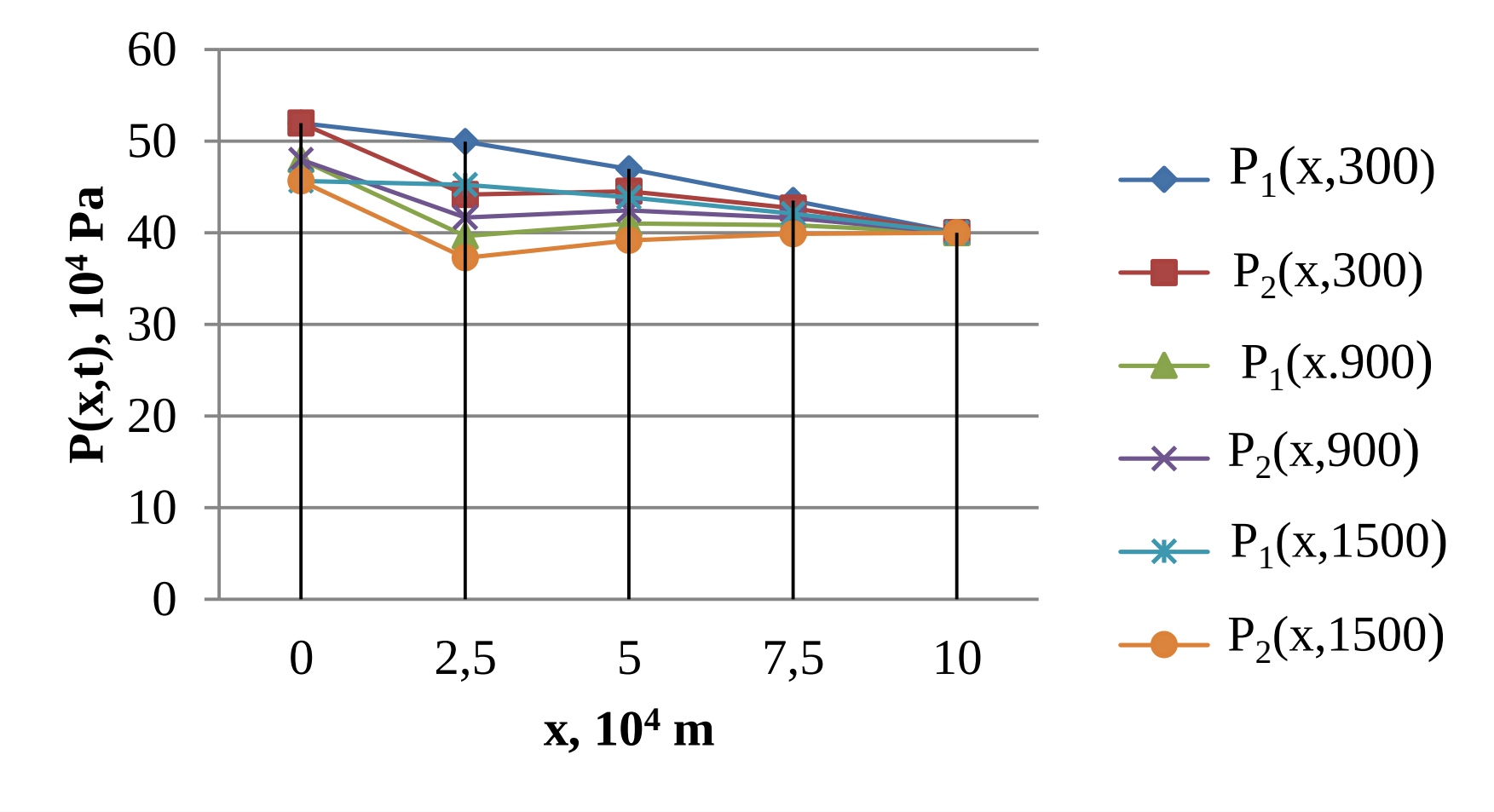

**Graph 1.3.5. The distribution of pressure along the length of the damaged and undamaged pipelines in the non-stationary regime due to the rupture of the parallel gas pipeline at a distance of $\ell_2$ =25000 m and a gas leakage rate of $G_{ut}$=0,4 $G_0$.**

The analysis of Graph 1.3.5 reveals that, except for the pressure at the final point of the gas pipeline, where the pressure remains constant under steady operating conditions, the pressure values along both the damaged and undamaged sections of the pipeline decrease over time due to the rupture of the parallel pipeline at a distance of $\ell_2$=25000 m. Additionally, the rate of pressure decrease increases along the direction from the final point to the starting point of the pipeline in both the damaged and undamaged sections. As a result, a process of gas discharge occurs at the beginning of the gas pipeline. By increasing the mass of leaking gas without changing the operating conditions of the gas pipeline, we can determine the maximum discharge time at which this process is possible.

In order to maintain the rupture point of the second line at a distance of $\ell_2$=25000 m and with the gas leakage rates of $G_{ut}$=0,8 $G_0$ and $G_{ut}$= 1,2 $G_0$ , we calculate the pressure values along the length of the damaged and undamaged pipelines at different times, every 2,5 ×$10^4$ m, and record them in the following tables (Table 1.3.4 and Table 1.3.5).

| $G_{ut}$=0,8 $G_0$ | $P_1(x,t)$, $10^4$ Pa | $P_2(x,t)$ , $10^4$ Pa | $P_1(x,t)$ , $10^4$ Pa | $P_2(x,t)$ , $10^4$ Pa | $P_1(x,t)$ , $10^4$ Pa | $P_2(x,t)$ , $10^4$ Pa |
|---|---|---|---|---|---|---|
| **X, $10^4$ m** | **t=300 sec.** | | **t=900 sec** | | **t=1500 sec.** | |
| 0 | 48,91 | 48,91 | 40,93 | 40,93 | 36,33 | 36,33 |
| 2,5 | 48,64 | 37,09 | 42,89 | 28,07 | 39,19 | 23,3 |
| 5 | 46,41 | 41,57 | 42,84 | 34,59 | 40,31 | 30,84 |
| 7,5 | 43,35 | 41,6 | 41,67 | 37,99 | 40,39 | 35,98 |
| 10 | 40 | 40 | 40 | 40 | 40 | 40 |

**Table 1.3.4. Pressure values along the length of the damaged and undamaged pipelines, depending on time, for the case of the rupture of the first line of the gas pipeline and a gas leakage rate of $G_{ut}$=0,8 $G_0$.**

| $G_{ut}$=1,2 $G_0$ | $P_1(x,t)$, $10^4$ Pa | $P_2(x,t)$ , $10^4$ Pa | $P_1(x,t)$ , $10^4$ Pa | $P_2(x,t)$ , $10^4$ Pa | $P_1(x,t)$ , $10^4$ Pa | $P_2(x,t)$ , $10^4$ Pa |
|---|---|---|---|---|---|---|
| X, $10^4$ m | t=300 sec. | | t=900 sec. | | t=1500 sec. | |
| 0 | 45,86 | 45,86 | 33,89 | 33,89 | 27 | 27 |
| 2,5 | 47,33 | 30,01 | 38,72 | 16,48 | 33,15 | 9,33 |
| 5 | 45,86 | 38,61 | 40,52 | 28,13 | 36,71 | 22,55 |
| 7,5 | 43,15 | 40,52 | 40,63 | 35,11 | 38,72 | 32,09 |
| 10 | 40 | 40 | 40 | 40 | 40 | 40 |

**Table 1.3.5. Pressure values along the length of the damaged and undamaged pipelines, depending on time, for the case of the rupture of the first line of the gas pipeline and a gas leakage rate of $G_{ut}$=0,4 $G_0$ .**

We will create the following graphs using the data from Table 1.3.4 and Table 1.3.5 (Graph 1.3.6 and Graph 1.3.7).

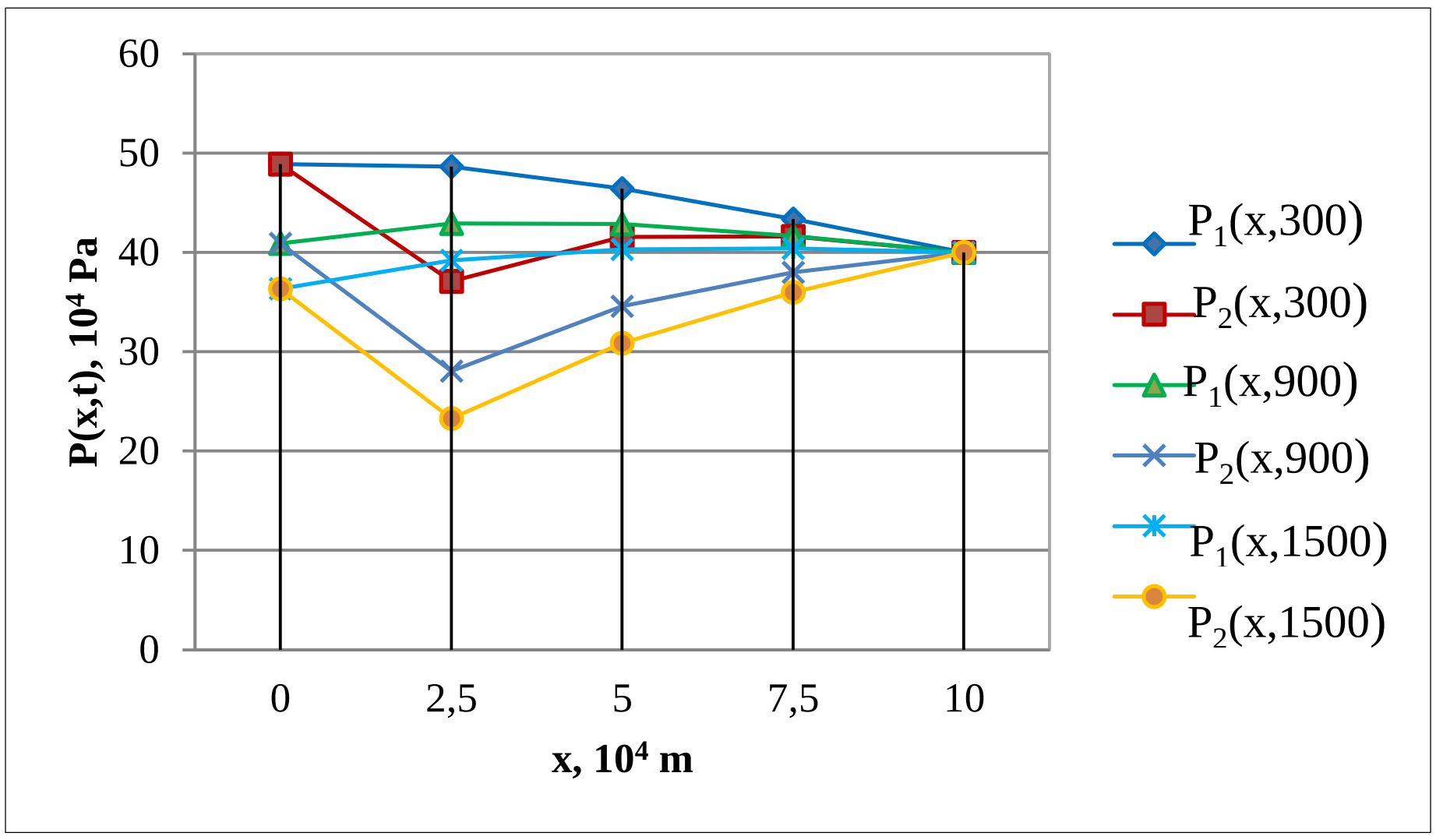

**Graph 1.3.6. The distribution of pressure along the length of the damaged and undamaged pipelines in the non-stationary regime due to the rupture of the parallel gas pipeline at a distance of $\ell_2$ =25000 m and a gas leakage rate of $G_{ut}$=0,8 $G_0$.**

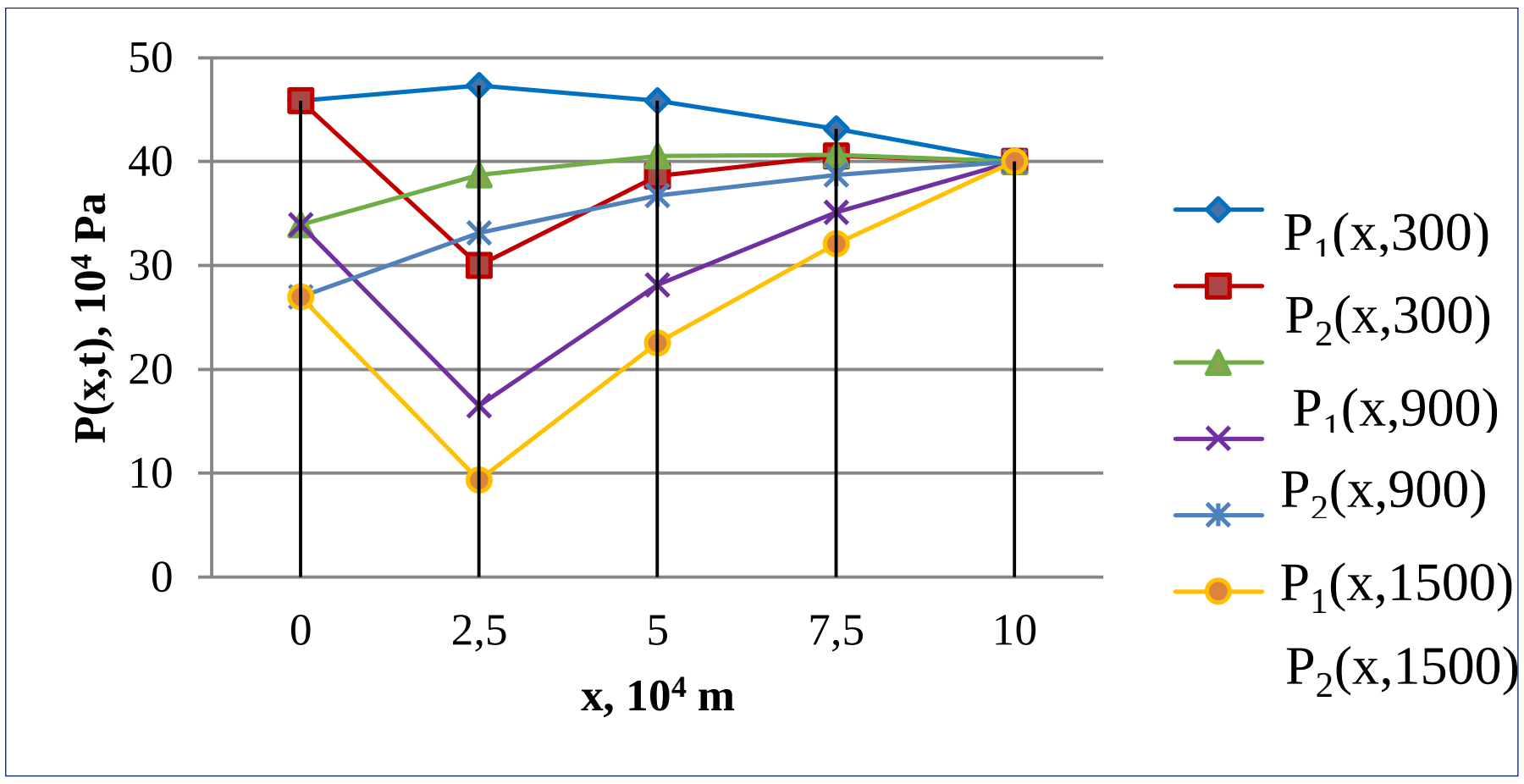

**Graph 1.3.7. Visualization of the longitudinal distribution of pressures in a non-stationary regime along the damaged and undamaged pipeline sections, when the integrity of a parallel gas pipeline is compromised at a distance of $\ell_2$ =25000 m and gas leakage occurs at $G_{ut}$=1,2$G_0$**

The analysis of Graph 1.3.6 and 1.3.7 indicates that, under steady-state operating conditions with constant pressure maintained at the endpoint of the gas pipeline, the profile of longitudinal pressure decreases along the damaged and undamaged pipeline sections in the direction from the endpoint to the starting point of the pipeline, which resembles the curve depicted in Graph 1.3.5. However, due to an increase in the mass flow rate of gas leakage, the extent of pressure reduction along the damaged pipeline increases sufficiently. Specifically, at the point xəttinin $\ell_2$ =25000 m along the damaged pipeline segment (where the rupture occurred), if the pressure decrease over a duration of t = 1500 seconds is (43,5- 37,28) ×$10^4$ =6,22×$10^4$ Pa for gas leakage at $G_{ut}$=0,4$G_0$, then for gas leakage at Gut=0.8G0, it would be (43,5- 23,3) ×$10^4$ =20,2×$10^4$ Pa, and for $G_{ut}$=1,2$G_0$ $_0$, it would be (43,5- 9,33) ×$10^4$ =34,17×$10^4$ Pa. Thus, as a result of the mass flow rate of gas leakage increasing by a factor of (1.2/0.4)=3, the pressure decrease at the point $\ell_2$ =25000 m along the damaged pipeline increases by a factor of (34.17/6.22)=5.5. This means that the pressure decrease accommodates nearly twice the increased magnitude of the mass flow rate of gas leakage. On the other hand, according to Graph 1.3.7, in the damaged pipeline, the minimum value of pressure decrease along the pipeline in the direction from the endpoint to the starting point is observed at the point $\ell_2$ =25000 m, reaching a value of (9,33·$10^4$ Pa), and this value increases gradually towards the starting section, reaching a value of 27×$10^4$ Pa at the starting point. At the same time, this value is equal to the pressure at the same point (x = 0) of the undamaged pipeline. The reason for this is the parallel operation of the gas pipelines in a unified hydraulic regime. In other words, the pipelines have the same starting and ending points. As a result, at the beginning section, the gas flow from the undamaged pipeline merges with the damaged pipeline. From this, we can conclude that in non-stationary regimes, the process of diverting the gas flow from the undamaged section to the damaged section and extracting it from the leakage point (Scheme 1.3.1) results in a discharge process at the beginning of both pipeline sections. As a result, in the section from the beginning of the gas pipeline to the point of rupture, the possibility of alternate emptying and filling of the damaged pipeline section is accompanied by a sharp decrease in pressure, and these pressure values are often significantly lower than those at the endpoint of the pipeline. This contradicts the requirement for normal gas flow in the pipeline for gas-dynamic processes. Therefore, as mentioned above, the exploitation process of parallel gas pipelines operating in a unified hydraulic regime is not effectively accounted for.

**The analysis of non-stationary gas dynamic processes under exploitation conditions with maintained constant pressure at the beginning of parallel pipeline systems.**

Let's analyze the specified exploitation conditions for the gas pipeline. Let's assume that at the beginning of the utilized parallel gas pipeline, the pressure is maintained constant. For the solution of the equations in mathematical model (2), we accept the initial stationary distribution of gas flow pressures in the calculation domain.

$t=0;\ P_i(x,0)=P_b-2aG_0x:\ i=1,2,3$

Using the derived general model and maintaining constant pressure at the beginning of the gas pipeline, we consider the pressure functions that describe the flow rate G(t) and its abrupt changes at the end of the pipeline. Thus, the following boundary conditions are adopted for the solution of equation (2):.

At the point x=0 $\begin{cases} P_1(x,t) = P_b \\ P_2(x,t) = P_b \end{cases}$

At the point $x = \ell_2$ $\begin{cases} P_2(x,t) = P_3(x,t) \\ \dfrac{\partial P_3(x,t)}{\partial X} - \dfrac{\partial P_2(x,t)}{\partial X} = -2aG_{ut}(t) \end{cases}$

At the point x= L $\begin{cases} P_1(x,t) = P_3(x,t) \\ \dfrac{\partial P_1(x,t)}{\partial X} + \dfrac{\partial P_2(x,t)}{\partial X} = -2aG_k(t) \end{cases}$

Using the Gauss method as outlined above, if we determine the coefficients specified in (3), (4), and (5) functions and substitute them into the equations, and after carrying out the uncomplicated operations, we obtain the following expressions to determine the distribution of pressure in transformed form:

$$P_1(x,S) = \frac{P_1(x,0)}{S} + \frac{\beta}{2}\left[\frac{G_0}{S} - G_s(S)\right]\frac{sh\lambda x}{ch\lambda L} - \\ -\beta G_{ut}(S)\frac{sh\lambda x sh\lambda \ell_2}{sh\lambda 2L} \qquad 0 \le x \le L \qquad (11)$$

$$P_2(x,S) = \frac{P_2(x,0)}{S} + \frac{\beta}{2}\left[\frac{G_0}{S} - G_s(S)\right]\frac{sh\lambda x}{ch\lambda L} - \\ -\beta G_{ut}(S)\frac{sh\lambda x sh\lambda(2L-\ell_2)}{sh\lambda 2L} \qquad 0 \le x \le \ell_2 \qquad (12)$$

$$P_3(x,S) = \frac{P_3(x,0)}{S} + \frac{\beta}{2}\left[\frac{G_0}{S} - G_s(S)\right]\frac{sh\lambda x}{ch\lambda L} - \\ -\beta G_{ut}(S)\frac{sh\lambda \ell_2 sh\lambda(2L-x)}{sh\lambda 2L} \qquad \ell_2 \le x \le L \qquad (13)$$

The original general solutions for functions (11), (12), and (13) will be written as follows:

$$P_1(x,t) = P_b - 8aG_oL\sum_{n=1}^{\infty}(-1)^{n-1}\frac{e^{-\alpha_1 t}}{[\pi(2n-1)]^2}Sin\frac{\pi(2n-1)x}{2L} - \\ -\frac{c^2}{L}\sum_{n=1}^{\infty}(-1)^{n-1}Sin\frac{\pi(2n-1)x}{2L}\int_0^t[G_s(\tau)]\cdot e^{-\alpha_1(t-\tau)}d\tau - \qquad 0 \le x \le L \qquad (14) \\ -\frac{c^2}{L}\sum_{n=1}^{\infty}(-1)^n Sin\frac{\pi n\ell_2}{2L}Sin\frac{\pi nx}{2L}\int_0^t G_{ut}(\tau)\cdot e^{-\alpha_2(t-\tau)}d\tau$$

$$P_{2,3}(x,t)=P_b-8aG_oL\sum_{n=1}^{\infty}(-1)^{n-1}\frac{e^{-\alpha_1 t}}{[\pi(2n-1)]^2}Sin\frac{\pi(2n-1)x}{2L}-$$

$$\frac{c^2}{L}\sum_{n=1}^{\infty}(-1)^{n-1}Sin\frac{\pi(2n-1)}{2L}\int_0^t[G_s(\tau)]\cdot e^{-\alpha_1(t-\tau)}d\tau- \qquad (15)$$

$$-\frac{c^2}{L}\sum_{n=1}^{\infty}(-1)^n\int_0^t G_{ut}(\tau)\cdot e^{-\alpha_2(t-\tau)}d\tau\begin{Bmatrix} Sin\frac{\pi n(2L-\ell_2)}{2L}Sin\frac{\pi nx}{L},0\le x\le\ell_2 \\ Sin\frac{\pi n\ell_2}{2L}Sin\frac{\pi n(2L-x)}{L},\ell_2\le x\le L\end{Bmatrix}$$

Here, $\alpha_1=\frac{\pi^2(2n-1)^2c^2}{8aL^2}$, $\alpha_2=\frac{\pi^2n^2c^2}{8aL^2}$

As noted, it is essential to develop methods for analyzing the dynamics of the main section of gas pipelines in parallel systems to better understand the operation of the systems. To achieve this, it is necessary to adopt simplification methods based on theoretical and experimental tests in relevant areas of application and to use calculation programs based on simple principles. Currently, such programs are available, and they are specialized simplifications for the operational management of pipeline networks [22].

$G_s(t)=G_0$, $G_{ut}(t)=G_{ut}=constat$, n=1,2,3…..10.

First, by maintaining the position of the rupture point of the second line at a distance of $\ell_2$ =25000 m along the gas pipeline and considering the variations of gas leakage at $G_{ut}$=0,8 $G_0$ and $G_{ut}$= 1,6 $G_0$, we calculate the values of pressure along the length of the damaged and undamaged pipeline for various times at intervals of 2,5 ×$10^4$ m, and we present the following tables (Table 1.3.6 and Table 1.3.7).

| $G_{ut}$=0,8 $G_0$ | $P_1(x,t)$, $10^4$ Pa | $P_2(x,t)$ , $10^4$ Pa | $P_1(x,t)$ ), $10^4$ Pa | $P_2(x,t)$ , $10^4$ Pa | $P_1(x,t)$ , $10^4$ Pa | $P_2(x,t)$ , $10^4$ Pa |
|---|---|---|---|---|---|---|
| X, $10^4$ m | t=300 sec. | | t=900 sec. | | t=1500 sec. | |
| 0 | 55 | 55 | 55 | 55 | 55 | 55 |
| 2,5 | 51,19 | 39,64 | 50,73 | 35,91 | 50,32 | 34,43 |
| 5 | 47,33 | 42,5 | 46,34 | 38,08 | 45,5 | 36,03 |
| 7,5 | 43,35 | 41,6 | 41,67 | 37,99 | 40,39 | 35,98 |
| 10 | 39,07 | 39,07 | 36,5 | 36,5 | 34,81 | 34,81 |

**Table 1.3.6. Pressure values along the length of the damaged and undamaged pipeline as a function of time for the rupture and gas leakage at $G_{ut}$=0,8 $G_0$ in the first variant of the gas pipeline.**

| $G_{ut}$=1,6 $G_0$ | $P_1(x,t)$, $10^4$ Pa | $P_2(x,t)$ , $10^4$ Pa | $P_1(x,t)$ , $10^4$ Pa | $P_2(x,t)$ , $10^4$ Pa | $P_1(x,t)$ , $10^4$ Pa | $P_2(x,t)$ , $10^4$ Pa |
|---|---|---|---|---|---|---|
| x, $10^4$ m | t=300 sec. | | t=900 sec. | | t=1500 sec. | |
| 0 | 55 | 55 | 55 | 55 | 55 | 55 |
| 2,5 | 51,13 | 28,04 | 50,22 | 20,57 | 49,38 | 17,62 |
| 5 | 47,17 | 37,5 | 45,18 | 28,67 | 43,5 | 24,56 |
| 7,5 | 42,96 | 39,44 | 39,59 | 32,23 | 37,04 | 28,2 |
| 10 | 38,15 | 38,15 | 33 | 33 | 29,61 | 29,61 |

**Table 1.3.7. Pressure values along the length of the damaged and undamaged pipeline as a function of time for the rupture and gas leakage at $G_{ut}$=1,6 $G_0$ in the first variant of the gas pipeline.**

Using the data provided in Tables 1.3.6 and 1.3.7, let's construct the following graphs (Graph 1.3.8 and Graph 1.3.9).

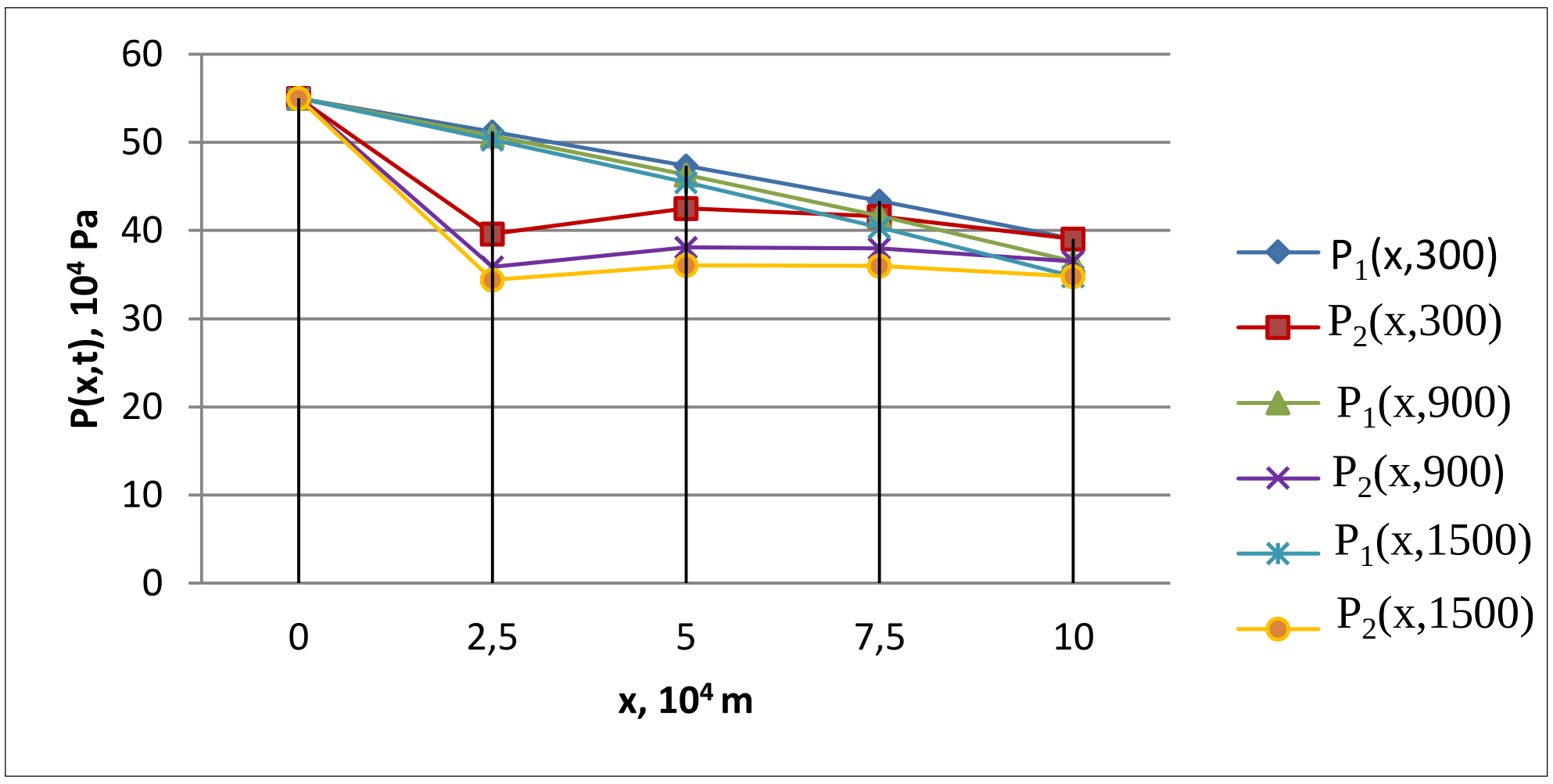

**Graph 1.3.8. The distribution of pressure along the length of the damaged and undamaged pipelines in the non-stationary regime due to the rupture and gas leakage at $G_{ut}$=0,8$G_0$ in the first variant of the parallel gas pipeline.**

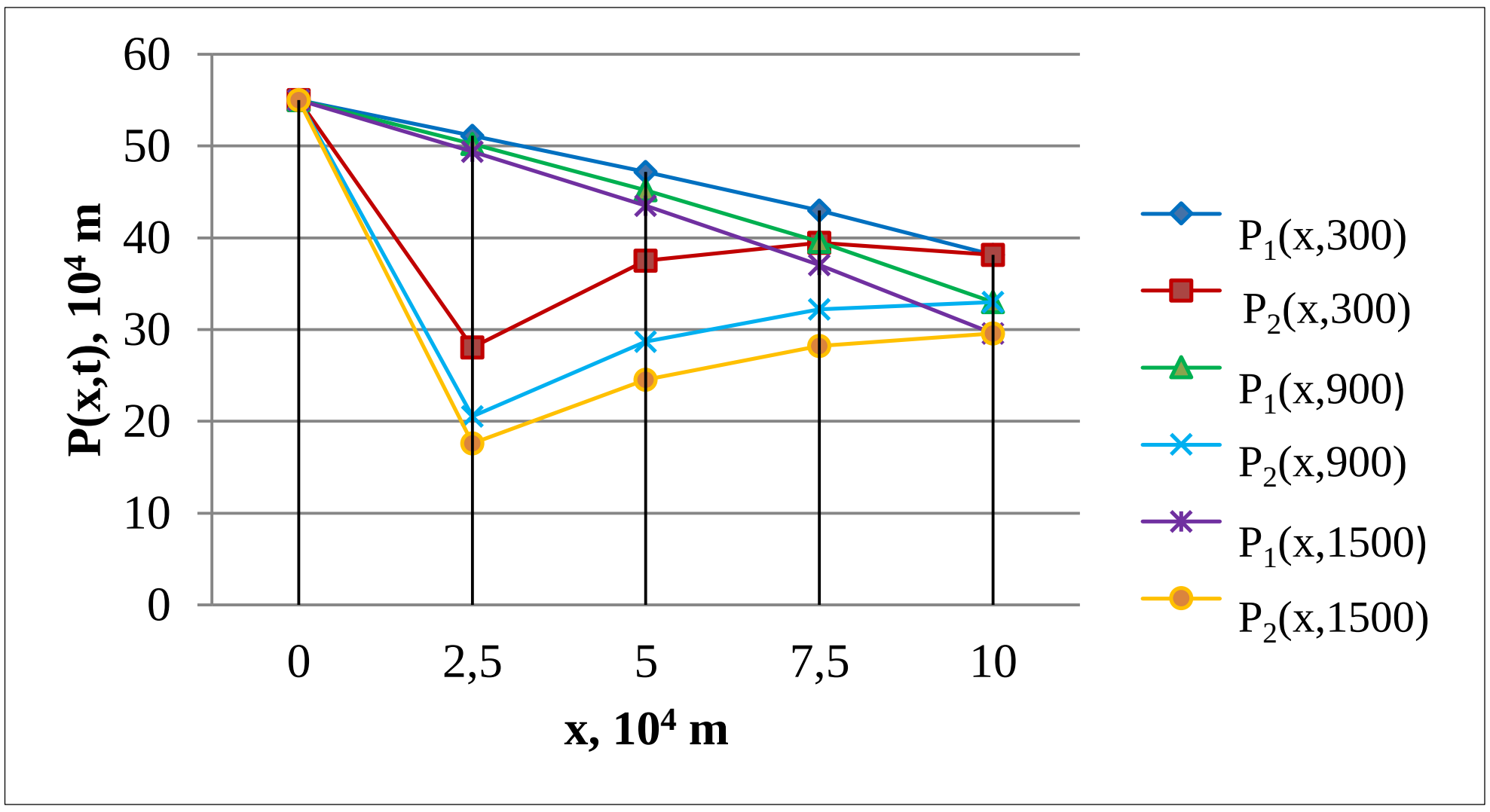

**Graph 1.3.9. The distribution of pressure along the length of the damaged and undamaged pipelines in the non-stationary regime due to the rupture and gas leakage at $G_{ut}$=1,6$G_0$ in the first variant of the parallel gas pipeline.**

The analysis of Graph 1.3.8 and Graph 1.3.9 reveals that, under conditions where pressure remains constant at the starting point of the gas pipeline, the curvature depicted by the declining pressure along the length of the pipeline from the endpoint towards the starting point is contrary to the view of the pipeline's endpoint under conditions of constant operational pressure.

In the case where the damaged section of the line is at $\ell_2$ =25000 m (the damaged section due to the rupture), in the Gut=0.8$G_0$ variant, at t=1500 seconds, if the pressure decrease due to gas leakage is (43,5- 34,43) $\times 10^4$ =9,07$\times 10^4$ Pa, then in the $G_{ut}$=1,6$G_0$ variant, the pressure decrease due to gas leakage would be (43,5- 17,62) $\times 10^4$ =25,88$\cdot 10^4$ Pa.

So, as a result of the gas leakage mass loss increasing by a factor of (1.6/0.8)=2, the pressure drop at the point $\ell_2$ =25000 m on the damaged pipeline increases by a factor of (25.88/9.07)=2.9. Thus, the pressure drop factor reflects the increased magnitude of gas leakage.

On the other hand, according to Graph 1.3.9, the pressure drop value along the pipeline in the damaged section, in the direction from the starting point of the pipeline to its endpoint, reaches its minimum value (17,62$\times 10^4$ Pa) at the point $\ell_2$ =25000 m and increases towards the endpoint, reaching a value of 29,61$\times 10^4$ Pa. At the same time, the pressure at the endpoint of the undamaged pipeline is the same. This pressure value is approximately 1.4 times lower than the values in the steady-state regime. The reason for this is the operation of the parallel gas coupling in a single hydraulic regime. In other words, it means that the pipelines have the same starting and ending points. As a result, at the endpoint, the gas flow from the undamaged pipeline converges with the damaged pipeline. From this, it can be inferred that in non-steady-state regimes, the process of diverting the gas flow from the undamaged to the damaged pipeline and removing it from the leakage point (Scheme 1.3.1) results in the emptying process occurring in both pipelines. As a result, the possibilities of alternative emptying and filling of the damaged pipeline section to the point of rupture are accompanied by a sharp decrease in pressure in that section, and these pressure values are several times lower than the values at the endpoint of the pipeline.

To conduct a more detailed analysis of the information regarding the non-normal gas yield compromising the reliability of the gas pipeline, we assume the location of the gas leakage point on the pipeline to be $\ell_2$ = 7,5$\times$ $10^4$ m and calculate the values of gas pressure in the pipeline at intervals of 300 seconds for different values of gas mass leakage at the point x = 75000 m. We then record these pressure values in the following Table 1.3.8.

| $G_{ut}$=1,2 $G_0$ | $P_1(x,t)$, $10^4$ Pa | $P_2(x,t)$ , $10^4$ Pa | $P_1(x,t)$ , $10^4$ Pa | $P_2(x,t)$ , $10^4$ Pa | $P_1(x,t)$ , $10^4$ Pa | $P_2(x,t)$ , $10^4$ Pa | $P_1(x,t)$ , $10^4$ Pa | $P_2(x,t)$ , $10^4$ Pa |
|---|---|---|---|---|---|---|---|---|
| X, $10^4$ m | t=300 sec. | | t=900 sec. | | t=1200 sec. | | t=1500 sec. | |
| 0 | 55 | 55 | 55 | 55 | 55 | 55 | 55 | 55 |
| 2,5 | 50,65 | 48,02 | 48,13 | 42,61 | 47,09 | 40,91 | 46,22 | 39,59 |
| 5 | 45,86 | 38,61 | 40,52 | 28,13 | 38,43 | 24,96 | 36,71 | 22,5 |
| 7,5 | 39,83 | 22,51 | 31,22 | 8,98 | 28,14 | 4,95 | 25,65 | 1,83 |
| 10 | 30,86 | 30,86 | 18,89 | 18,89 | 15,04 | 15,04 | 12 | 12 |

**Table 1.3.8. Pressure values along the length of the damaged and undamaged pipeline segments depending on time for the case of rupture of the gas coupling and gas leakage at $G_{ut}$=1,2 $G_0$ in the 2nd variant.**

Using the data provided in Table 8, let's construct the following graphs (Graph 1.3.10).

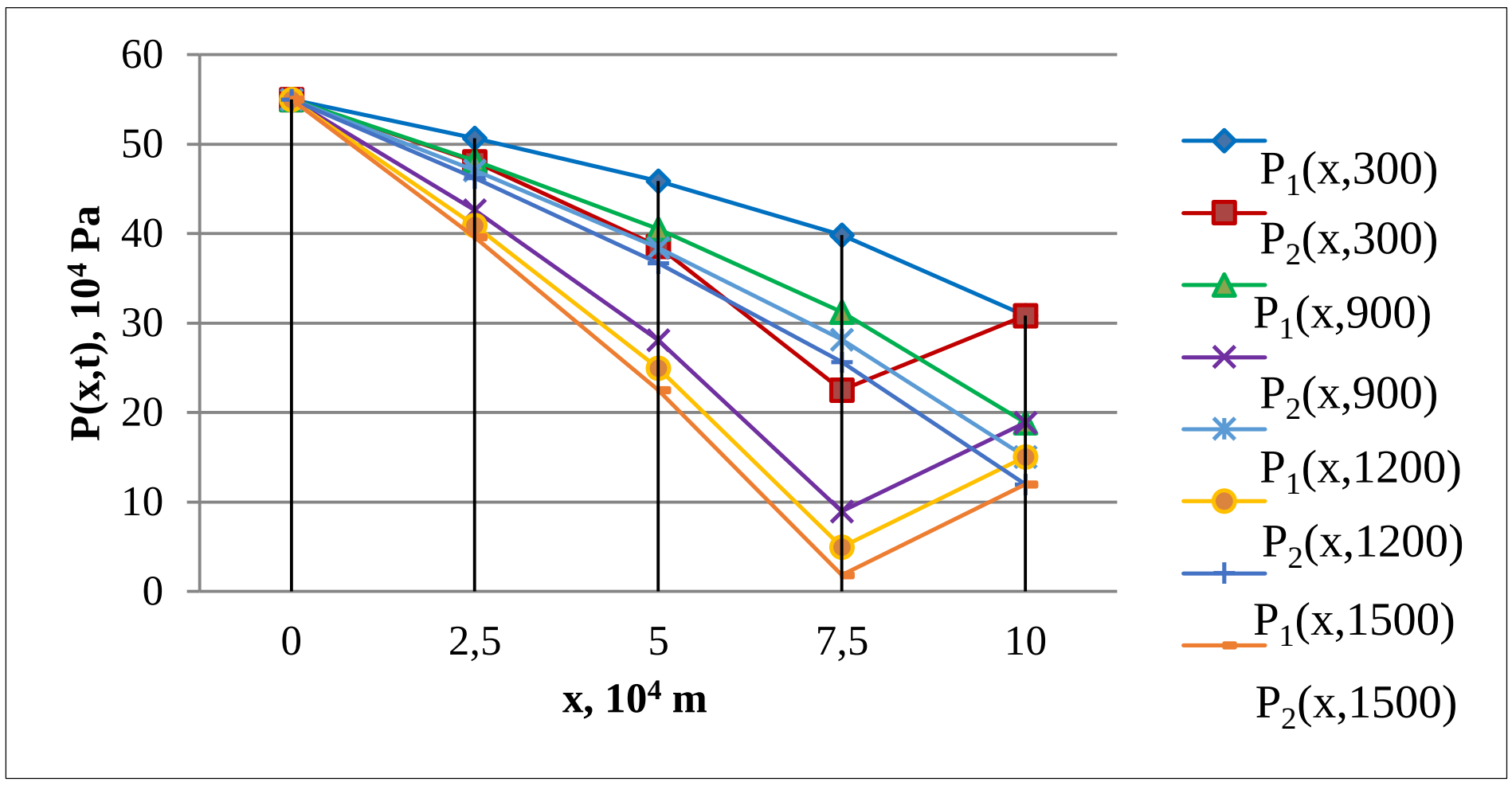

**Graph 1.3.10. The distribution of pressure along the length of the damaged and undamaged pipeline segments in non-steady-state conditions due to the rupture of the gas coupling at $\ell_2$ =75000 m distance and gas leakage at $G_{ut}$=1,2$G_0$.**

From Graph 1.3.10, it is evident that at the point $\ell_2$ =75000 m on the damaged pipeline (where the coupling is ruptured), the pressure drop due to gas leakage at $G_{ut}$=1,2$G_0$ 0 scenario is (43,5- 1,83)×$10^4$ =41,92×$10^4$ Pa at t = 1500 seconds. On the other hand, according to Graph 1.3.10, the pressure drop along the pipeline in the damaged section, in the direction from the starting point of the pipeline to its endpoint, reaches its minimum value (1,83×$10^4$ Pa) at the point $\ell_2$ =75000 m and increases towards the endpoint, reaching a value of 12×$10^4$ Pa at the endpoint. On the other hand, the pressure at the endpoint of the damaged pipeline is equal to that value. This value is 3.3 times lower than their values in the steady-state regime. Therefore, the gas efficiency for consumers cannot exceed 40%. This means that despite the ingress of gas leakage from the damaged pipeline into the undamaged pipeline, the gas mass emptied at t = 1500 seconds due to leakage surpasses the incoming gas mass, resulting in a sharp emptying of the final section.

We consider the point where the gas leak occurs in the gas pipeline to be $\ell_2$ = 7,5× $10^4$ m, and we calculate the values of the gas pressure in the pipeline at intervals of every t=300 seconds for different values of gas mass leakage at the point x=75000 m. We record these pressure values in the following Table 1.3.9.

| $\ell_2$=75000 m | $P_1$(75000,t), $10^4$ Pa | | | | $P_2$(75000,t), $10^4$ Pa | | | |
|---|---|---|---|---|---|---|---|---|
| t, sec | $G_{ut}$ = 0,4 $G_0$ | $G_{ut}$ = 0,8 $G_0$ | $G_{ut}$ = 1,2 $G_0$ | $G_{ut}$ = 1,6 $G_0$ | $G_{ut}$ = 0,4 $G_0$ | $G_{ut}$ = 0,8 $G_0$ | $G_{ut}$ = 1,2 $G_0$ | $G_{ut}$ = 1,6 $G_0$ |
| 300 | 42,44 | 41,14 | 39,83 | 38,53 | 36,67 | 29,59 | 22,51 | 15,44 |
| 600 | 40,86 | 37,96 | 35,07 | 32,17 | 33,98 | 24,21 | 14,44 | 4,67 |
| 900 | 39,57 | 35,39 | 31,22 | 27,04 | 32,16 | 20,57 | 8,98 | 0 |
| 1200 | 38,55 | 33,34 | 28,24 | 22,94 | 30,82 | 17,89 | 4,95 | 0 |
| 1500 | 37,72 | 31,69 | 25,65 | 19,62 | 29,78 | 15,8 | 1,83 | 0 |

**Table 1.3.9. Pressure values of the damaged and undamaged pipeline segments over time for various values of gas leakage in the 2nd variant of the gas pipeline with the rupture of the gas coupling and gas leakage.**

Using the data provided in Table 1.3.9, let's construct the following graphs (Graph 1.3.11)

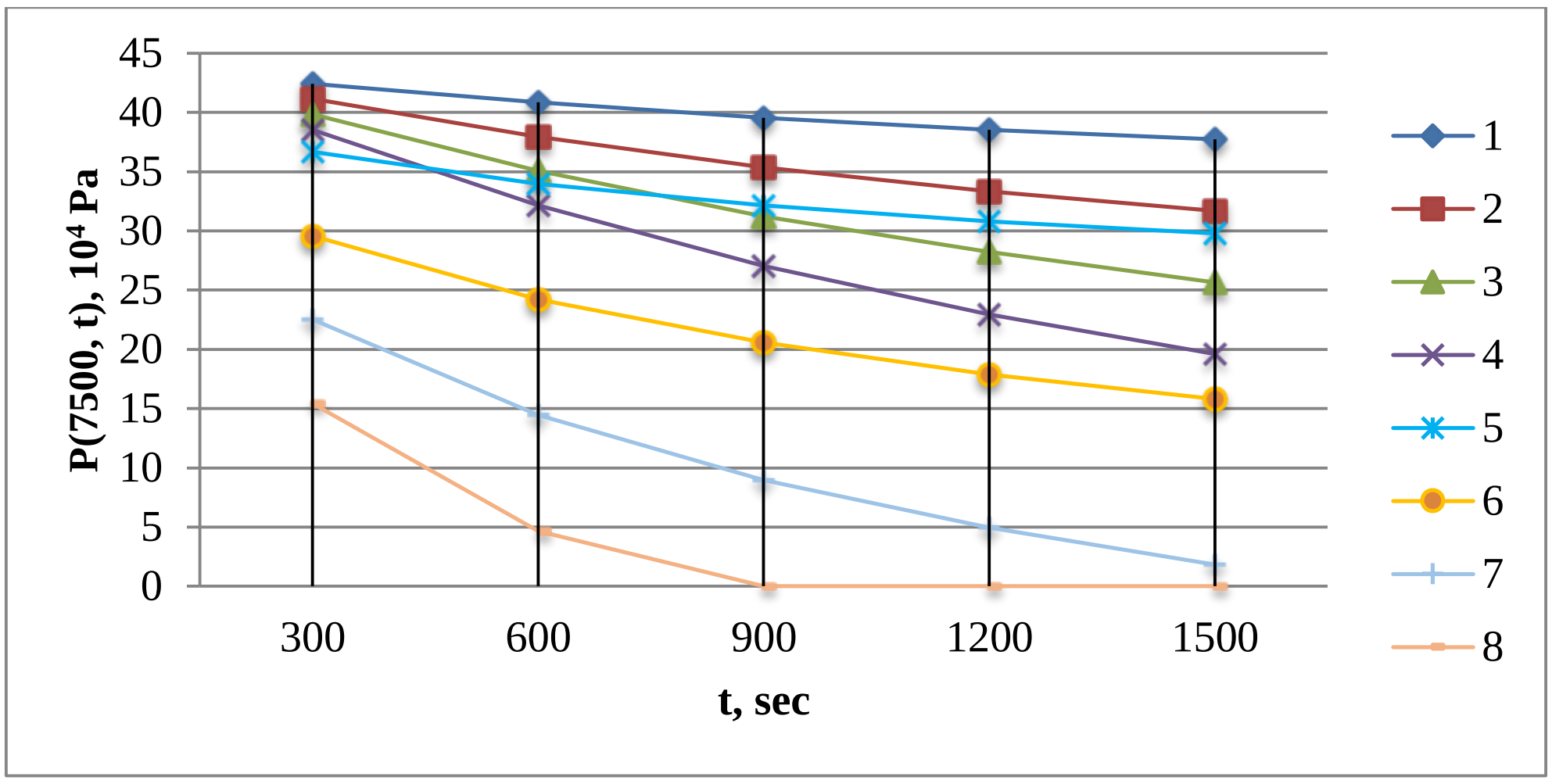

**Graph 1.3.11. The depiction of the variation of pressure over time at point x=75000 m in the gas pipeline: 1, 2, 3, and 4 - in the undamaged pipeline xəttində (1-$G_{ut}$=0,4 $G_0$, 2-$G_{ut}$=0,8 $G_0$, 3-$G_{ut}$=1,2 $G_0$, 4-$G_{ut}$ =1,6 $G_0$ ); 5, 6, 7, and 8 - in the damaged pipeline (at the point of coupling rupture): (5-$G_{ut}$=0,4 $G_0$, 6-$G_{ut}$=0,8 $G_0$, 7-$G_{ut}$=1,2 $G_0$, 8-$G_{ut}$ =1,6 $G_0$ ).**

From Table 1.3.9 and Graph 1.3.11, it is evident that the pressure drop in the damaged pipeline is multiple times greater than the drop in the undamaged pipeline. For instance, within a duration of t=300 seconds and with the volume of leaked gas being 40% of the gas consumption in the steady-state regime ($G_{ut}$=0,4 $G_0$), the pressure drop in the damaged pipeline segment [P(75000,0) -$P_2$(75000,300)] = (43,75- 42,44)×$10^4$= 1,31×$10^4$ Pa, whereas the pressure drop in the undamaged pipeline segment [P(75000,0) -$P_1$(25000,300)]= (43,75-36,67) ×$10^4$=7,08×$10^4$ Pa.

The pressure drop in the damaged pipeline segment is 5.4 times greater than the pressure drop in the undamaged pipeline segment. This ratio varies over time for a constant value of gas leakage: it reaches a maximum value at t=300 seconds and then decreases over time, reaching 3.38 at t=600 seconds, 2.77 at t=900 seconds, 2.49 at t=1200 seconds, and further decreases to 2.32 at t=1500 seconds. This indicates that the ratio reaches its maximum value at t=300 seconds for a constant value of gas leakage ($G_{ut}$ =0,4 $G_0$) and then decreases over time (Graph 1.3.12).

The ratio varies over time for a constant value of gas leakage at $G_{ut}$=0,8 $G_0$ it reaches 5.43 at t=300 seconds, 3.37 at t=600 seconds, 2.77 at t=900 seconds, 2.48 at t=1200 seconds, and further decreases to 2.31 at t=1500 seconds. This indicates that the ratio reaches its maximum value (5.43) at t=300 seconds for a constant value of gas leakage ($G_{ut}$=0,8 $G_0$) and then decreases over time, reaching a value of 2.31 at t=1500 seconds (Graph 1.3.12). On the other hand, it is apparent that, despite an increase in the volume of gas leakage, the ratio of the pressure drop in the damaged and undamaged sections remains constant until t=900 seconds.

However, for $G_{ut}$ =1,6 $G_0$, this value decreases to 2.1 at t=1200 seconds and to 1.81 at t=1500 seconds. The reason for this is that the pressure in the damaged pipeline is approaching zero at the x=75000 m point, as shown in Graph 1.3.11. In other words, it indicates that the gas pipe in the damaged pipeline segment is essentially fully emptied. Therefore, after 900 seconds, the gas supply regime for consumers is disrupted.

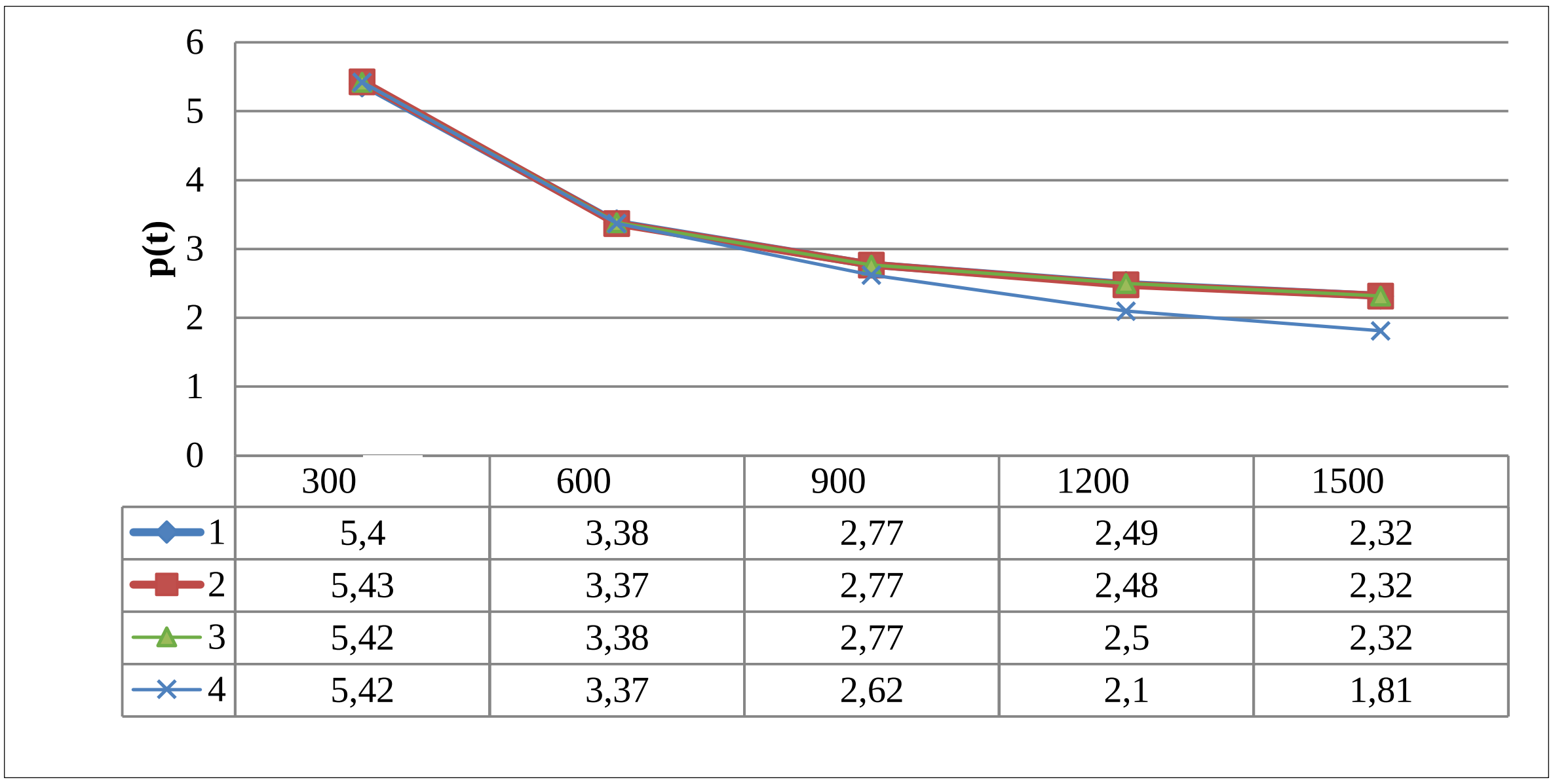

| | 300 | 600 | 900 | 1200 | 1500 |
|---|---|---|---|---|---|
| 1 | 5,4 | 3,38 | 2,77 | 2,49 | 2,32 |
| 2 | 5,43 | 3,37 | 2,77 | 2,48 | 2,32 |
| 3 | 5,42 | 3,38 | 2,77 | 2,5 | 2,32 |
| 4 | 5,42 | 3,37 | 2,62 | 2,1 | 1,81 |

**Graph 1.3.12. The graph of p(t) dependency on time (t) for the [Px, 0)-$P_2$ (x,t)]/[P(x, 0) -$P_1$ (x,t)]=p(t) at the x= 75000 m point of the gas pipeline (1-$G_{ut}$=0,4 $G_0$, 2-$G_{ut}$=0,8 $G_0$, 3-$G_{ut}$=1,2 $G_0$, 4-$G_{ut}$ =1,6 $G_0$).**

According to the mentioned conditions, it can be concluded that in the analyzed operational conditions of parallel gas pipelines [$P_1(0,t) = P_2(0,t) = P_s$ =constat], the gas supply regime to consumers is significantly disrupted after t = 1200 seconds. In other words, a large portion of consumers cannot be supplied with gas. Since the proposed operational condition does not comply with the continuous gas supply efficiency standards for consumers in accident scenarios, the effective utilization of parallel gas pipelines operating in a single hydraulic regime under the specified conditions is not feasible.

**The analysis of the non-stationary gas dynamic process considers the variation of pipeline losses at the beginning and end points of the gas pipeline under the given operating conditions.**

Let's analyze the following operating conditions for the gas pipeline: Suppose that the variation of mass flows at the beginning and end of the utilized parallel gas pipeline varies over time. For solving the mathematical model (2), we assume the initial condition of the gas flow pressures for the initial stationary distribution.

$t=0;\ P_i(x,0)=P_b-2aG_0x:\ i=1,2,3$

Using the general model obtained, let's consider the case where the variation of mass flow at the beginning and end of the gas pipeline is known. Thus, the following boundary conditions are assumed for solving equation (2):

At the point x=0 $$\begin{cases} P_1(x,t)=P_2(x,t) \\ \dfrac{\partial P_1(x,t)}{\partial X}+\dfrac{\partial P_2(x,t)}{\partial X}=-2aG_0(t) \end{cases}$$

At the point $x=\ell_2$ $$\begin{cases} P_2(x,t)=P_3(x,t) \\ \dfrac{\partial P_3(x,t)}{\partial X}-\dfrac{\partial P_2(x,t)}{\partial X}=-2aG_{ut}(t) \end{cases}$$

At the point x= L $\begin{cases} P_1(x,t) = P_3(x,t) \\ \dfrac{\partial P_1(x,t)}{\partial X} + \dfrac{\partial P_2(x,t)}{\partial X} = -2aG_k(t) \end{cases}$

Following the above sequence, by using the Gauss method, let's determine the coefficients stated in equations (3), (4), and (5) after substituting and performing straightforward operations. Then, we obtain the following expressions to determine the distribution of pressure in transformed form:

$$P_1(x,S) = \frac{P_1 - 2aG_0 x}{S} + \frac{\beta}{2}\frac{G_0}{S}\frac{sh\lambda(x - \frac{L}{2})}{ch\lambda\frac{L}{2}} + \frac{\beta}{2}G_0(S)\frac{ch\lambda(L-x)}{sh\lambda L} - \qquad 0 \le x \le L \qquad (16)$$
$$-\frac{\beta}{2}G_s(S)\frac{ch\lambda x}{sh\lambda L} - \frac{\beta}{2}G_{ut}(S)\frac{ch\lambda(L - \ell_2 - x)}{sh\lambda L}$$

$$P_2(x,S) = \frac{P_1 - 2aG_0 x}{S} + \frac{\beta}{2}\frac{G_0}{S}\frac{sh\lambda(x - \frac{L}{2})}{ch\lambda\frac{L}{2}} + \frac{\beta}{2}G_0(S)\frac{ch\lambda(L-x)}{sh\lambda L} - \qquad 0 \le x \le \ell_2 \quad (17)$$
$$-\frac{\beta}{2}G_s(S)\frac{ch\lambda x}{sh\lambda L} - \frac{\beta}{2}G_{ut}(S)\frac{ch\lambda(L - \ell_2 + x)}{sh\lambda L}$$

$$\ell_2 \le x \le L$$

$$P_3(x,S) = \frac{P_1 - 2aG_0 x}{S} + \frac{\beta}{2}\frac{G_0}{S}\frac{sh\lambda(x - \frac{L}{2})}{ch\lambda\frac{L}{2}} + \frac{\beta}{2}G_0(S)\frac{ch\lambda(L-x)}{sh\lambda L} - \qquad (18)$$
$$-\frac{\beta}{2}G_s(S)\frac{ch\lambda x}{sh\lambda L} - \frac{\beta}{2}G_{ut}(S)\frac{ch\lambda(L + \ell_2 - x)}{sh\lambda L}$$

The original general solution of equations (16), (17), and (18) will be written as follows:

$$P_1(x,t) = P_b - 2aG_0\frac{L}{4} + 4aG_oL\sum_{n=1}^{\infty}\frac{e^{-\alpha_1 t}}{[\pi(2n-1)]^2}Cos\frac{\pi(2n-1)x}{L} +$$
$$+\frac{c^2}{L}\sum_{n=1}^{\infty}Cos\frac{\pi nx}{L}\int_0^t\left[G_0(\tau) - (-1)^n G_s(\tau)\right]\cdot e^{-\alpha_2(t-\tau)}d\tau + \qquad 0 \le x \le L \quad (19)$$
$$+\frac{c^2}{2L}\left[G_0(\tau) - G_s(\tau) - G_{ut}(\tau)\right]d\tau - \frac{c^2}{L}\sum_{n=1}^{\infty}Cos\frac{\pi n(\ell_2 + x)}{L}\int_0^t G_{ut}(\tau)\cdot e^{-\alpha_2(t-\tau)}d\tau$$

$$P_{2,3}(x,t) = P_b - 2aG_0\frac{L}{4} + 4aG_oL\sum_{n=1}^{\infty}\frac{e^{-\alpha_1 t}}{[\pi(2n-1)]^2}Cos\frac{\pi(2n-1)x}{L} +$$
$$+\frac{c^2}{L}\sum_{n=1}^{\infty}Cos\frac{\pi nx}{L}\int_0^t\left[G_0(\tau) - (-1)^n G_s(\tau)\right]\cdot e^{-\alpha_2(t-\tau)}d\tau + \qquad (20)$$
$$+\frac{c^2}{2L}\left[G_0(\tau) - G_s(\tau) - G_{ut}(\tau)\right]d\tau - \frac{c^2}{L}\sum_{n=1}^{\infty}Cos\frac{\pi n(x - \ell_2)}{L}\int_0^t G_{ut}(\tau)\cdot e^{-\alpha_2(t-\tau)}d\tau \to \ell_2 \le x \le L$$

Here, $\alpha_1 = \dfrac{\pi^2(2n-1)^2 c^2}{8aL^2}$, $\alpha_2 = \dfrac{\pi^2 n^2 c^2}{8aL^2}$

An analysis of the investigated conditions [$P_1(0,t) = P_2(0,t) = P_b$ and $P_1(L,t) = P_2(L,t) = P_s$] of the parallel gas pipelines showed that under non-stationary regimes, the specified operating condition does not fulfill the function of providing normal gas supply to consumers for a certain period. For this purpose, we analyze the following operating conditions:

Firstly, in the first variant of the gas pipeline (at a distance of $\ell_2$ =25000 m), during the rupture of the second line and for the $G_{ut}$=0,4 $G_0$ gas leakage scenario, pressure values are calculated at different times for every 2,5×$10^4$ m meters along the length of the damaged and undamaged pipelines to determine the conformity of pressure distribution over time. The calculated pressure values at various times are recorded in the table below (Table 1.3.10).

| $G_{ut}$=0,4 $G_0$ | $P_1(x,t)$, $10^4$ Pa | $P_2(x,t)$, $10^4$ Pa | $P_1(x,t)$, $10^4$ Pa | $P_2(x,t)$, $10^4$ Pa | $P_1(x,t)$, $10^4$ Pa | $P_2(x,t)$, $10^4$ Pa |
|---|---|---|---|---|---|---|
| X, km | | t=300 sec. | | t=900 sec. | | t=1500 sec. |
| 0 | 53,45 | 53,47 | 50,92 | 51,34 | 48,65 | 49,85 |
| 2,5 | 50,57 | 47,69 | 48,53 | 45,08 | 46,43 | 43,28 |
| 5 | 47,16 | 45,97 | 45,48 | 43,84 | 43,52 | 42,35 |
| 7,5 | 43,5 | 43,11 | 41,95 | 41,74 | 40,03 | 40,68 |
| 10 | 39,66 | 39,77 | 37,98 | 39,09 | **36,02** | 38,57 |

**Table 1.3.10. Pressure Values over Time Dependent on the Length of the Damaged and Undamaged Pipeline in the First Variant of the Gas Pipeline for the Rupture of the Second Line and a Gas Leakage Scenario of $G_{ut}$=0,4 $G_0$:**

Based on the data provided in Table 1.3.10, we create the following graph (Graph 1.3.13)

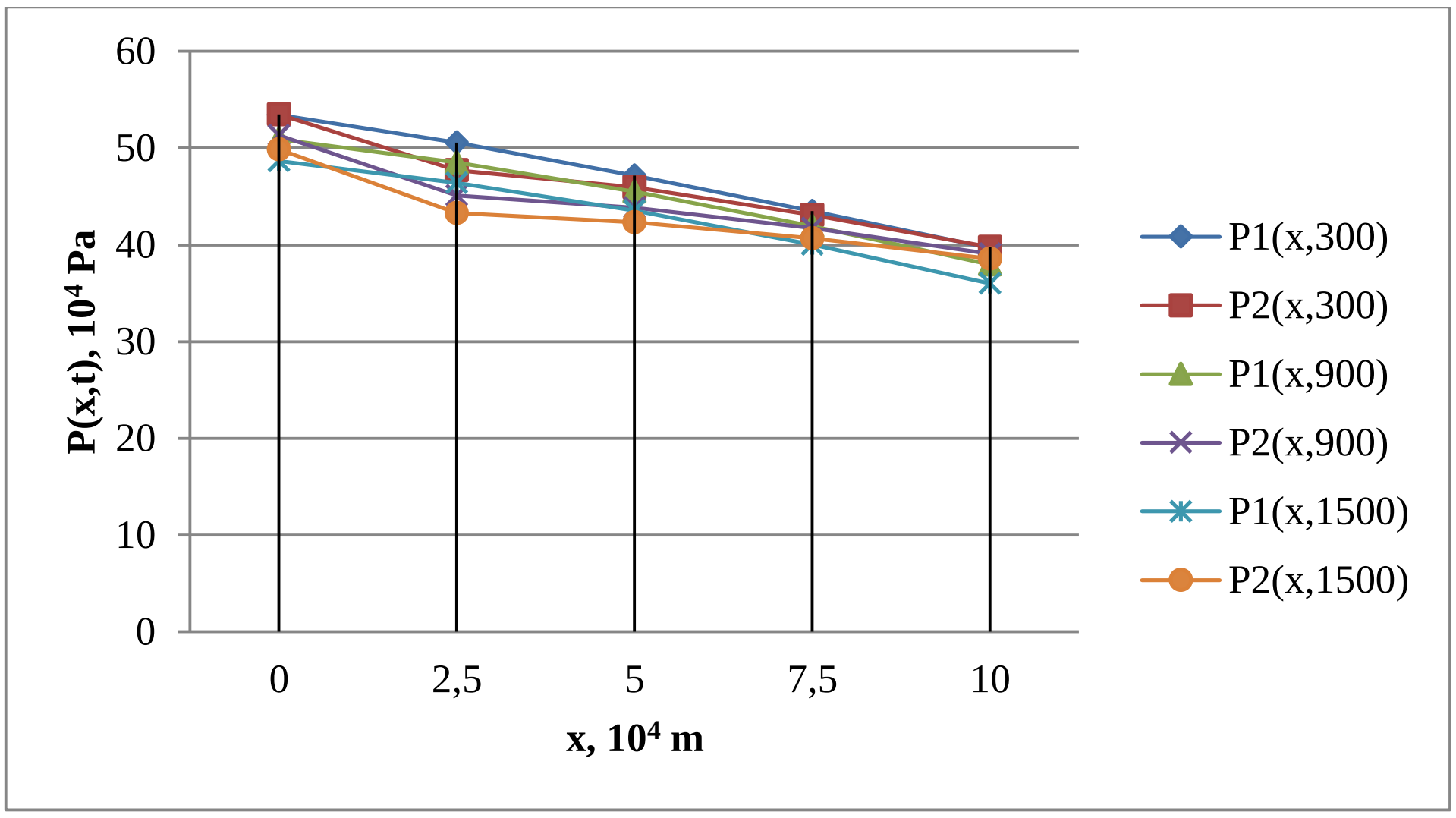

**Graph 1.3.13. Visualization of the non-stationary gas flow distribution along the length of the damaged and undamaged pipelines in the parallel gas pipeline at a distance of $\ell_2$ =25000 m, during the exposure of the 2nd line (damage) and gas leakage with $G_{ut}$=0,4$G_0$.**

It appears from the analysis of Graph 1.3.13 that, under the examined operating conditions of the gas pipeline, the pressure values of both the damaged and undamaged pipeline sections decrease over time due to the rupture of the 2nd line at a distance of $\ell_2$ =25000 m . On the other hand, at any given moment in the initial section of the gas pipeline, the pressure value of the undamaged pipeline section is greater than the pressure drop of the damaged pipeline section. Additionally, in the final section of the gas pipeline, at the examined moments, the pressure drop of the undamaged pipeline section is greater than the pressure drop of the damaged pipeline section.

This extreme condition of the gas dynamic process is explained by the fact that in the event of an accident caused by the rupture of one of the pipelines in a parallel gas pipeline system operating under unit hydraulic conditions, the main sections of the gas pipeline direct the flow of gas towards the damaged section. In other words, in operating conditions, the mass flows of the initial and final sections of the pipeline are directed towards the damaged section of the pipeline over time.As a result of the withdrawal of gas from both the initial and final sections of the undamaged pipeline and the continuous supply of gas to consumer facilities from its final section, the pressure and gas flow along the length of the undamaged pipeline vary both in coordinate and time due to non-uniform gas consumption. Thus, during the occurrence of a non-stationary flow regime resulting from the rupture of one of the pipelines in a parallel gas pipeline operating under unit hydraulic conditions, the non-stationary flow regime is transferred to the other undamaged pipelines, and when gas supply and withdrawal (leakage) are stopped simultaneously, the distribution of pressure in the pipeline segment can ensure the quality of the stationary regime.

Let's note that the process described above in the two operating conditions of the gas pipeline ($P_1(0,t)= P_2(0,t)=P_b$=const and $P_1(L,t)= P_2(L,t)=P_s$=const) has not been discovered. To illustrate our statements clearly, we will prepare the following table using the data from Table 1.3.10 (Table 1.3.11).

| t,sec | The pressure drop at the starting point of the gas pipeline | | The pressure drop at the leakage point of the gas pipeline | | The pressure drop at the end point of the gas pipeline | |
|---|---|---|---|---|---|---|
| | In the undamaged pipe line | In the damaged pipe line | In the undamaged pipe line | In the damaged pipe line | In the undamaged pipe line | In the damaged pipe line |
| | $[P_b-P_2(0,t)]$, $10^4$Pa | $P_b-P_2(0,t)]$, $10^4$Pa | $[P(\ell_2,0) - P_1(\ell_2,t)]$, $10^4$Pa | $[P(\ell_2,0) - P_2(\ell_2,t)]$, $10^4$Pa | $[P_s-P_1(L,t)]$, $10^4$Pa | $[P_s-P_2(L,t)]$, $10^4$Pa |
| 300 | 1,55 | 1,53 | 0,68 | 3,56 | 0,34 | 0,23 |
| 900 | 4,08 | 3,66 | 2,72 | 6,17 | 2,02 | 0,91 |
| 1500 | 6,35 | 5,15 | 4,82 | 7,97 | 3,98 | 1,43 |

**Table 1.3.11. Pressure Drop Values Over Time in the Initial, Leakage, and Final Sections of Damaged and Undamaged Pipeline Sections during the Rupture of the 2nd Line at $G_{ut}$=0,4 $G_0$**

Table 1.3.11 reveals that contrary to the process occurring at the initial and final points, the pressure drop value at the leakage point of the damaged pipeline section is greater than that of the undamaged pipeline section. This aligns with the physics of the gas dynamic process. In other words, the value of gas mass flow leaking into the unconfined environment at the point of damage in the parallel gas pipeline surpasses the value of gas mass flow entering the initial and final sections of the undamaged pipeline section.At the leakage point, the pressure drop values in both pipeline sections increase over time. As a result, a discharge process occurs in the leakage section of the damaged pipeline. Without altering the operating conditions of the gas pipeline, we can determine the maximum discharge time of the process by increasing the leaking gas mass.

Taking into account the rupture point of the 2nd line at a distance of $\ell_2$ =25000 m fixed, and considering gas leakage values of $G_{ut}$=0,8 $G_0$ and $G_{ut}$= 1,2 $G_0$, we calculate the pressure distribution along the length of the damaged and undamaged pipeline sections at various times, every 2,5×$10^4$ m, and record their values in the following tables (Table 1.3.12 and Table 1.3.13).

| $G_{ut}$=0,8 $G_0$ | $P_1(x,t)$, $10^4$ Pa | $P_2(x,t)$ , $10^4$ Pa | $P_1(x,t)$ , $10^4$ Pa | $P_2(x,t)$ , $10^4$ Pa | $P_1(x,t)$ , $10^4$ Pa | $P_2(x,t)$ , $10^4$ Pa |
|---|---|---|---|---|---|---|
| x, km | t=300 sec. | | t=900 sec. | | t=1500 sec. | |
| 0 | 51,91 | 51,95 | 46,83 | 47,67 | 42,3 | 44,7 |
| 2,5 | 49,88 | 44,12 | 45,8 | 38,81 | 41,61 | 35,32 |
| 5 | 46,82 | 44,45 | 43,46 | 38,91 | 39,54 | 37,2 |
| 7,5 | 43,26 | 42,47 | 40,15 | 40,17 | 36,32 | 37,61 |
| 10 | 39,32 | 39,54 | 35,96 | 38,18 | 32,04 | 37,13 |

**Table 1.3.12. Pressure Values Over Time Along the Length of Damaged and Undamaged Pipeline Sections during the Rupture of the 2nd Line in Variant 1 of the Gas Pipeline, with Gas Leakage at $G_{ut}$=0,8 $G_0$**

| $G_{ut}$=1,6 $G_0$ | $P_1(x,t)$, $10^4$ Pa | $P_2(x,t)$ , $10^4$ Pa | $P_1(x,t)$ , $10^4$ Pa | $P_2(x,t)$ , $10^4$ Pa | $P_1(x,t)$ , $10^4$ Pa | $P_2(x,t)$ , $10^4$ Pa |
|---|---|---|---|---|---|---|
| x, km | t=300 sec. | | t=900 sec. | | t=1500 sec. | |
| 0 | 48,81 | 48,89 | 38,66 | 40,34 | 29,6 | 34,4 |
| 2,5 | 48,51 | 37 | 40,36 | 26,57 | 31,97 | 19,38 |
| 5 | 46,15 | 41,39 | 39,42 | 32,84 | 31,57 | 26,9 |
| 7,5 | 42,77 | 41,19 | 36,54 | 35,58 | 28,89 | 31,47 |
| 10 | 38,65 | 39,07 | 31,92 | 36,36 | 24,07 | 34,27 |

**Table 1.3.13. Pressure Values Over Time Along the Length of Damaged and Undamaged Pipeline Sections during the Rupture of the 2nd Line in Variant 1 of the Gas Pipeline, with Gas Leakage at $G_{ut}$=1,6 $G_0$**

Using the data from Table 1.3.12 and Table 1.3.13, we create the following graphs (Graph 1.3.14 and Graph 1.3.15).

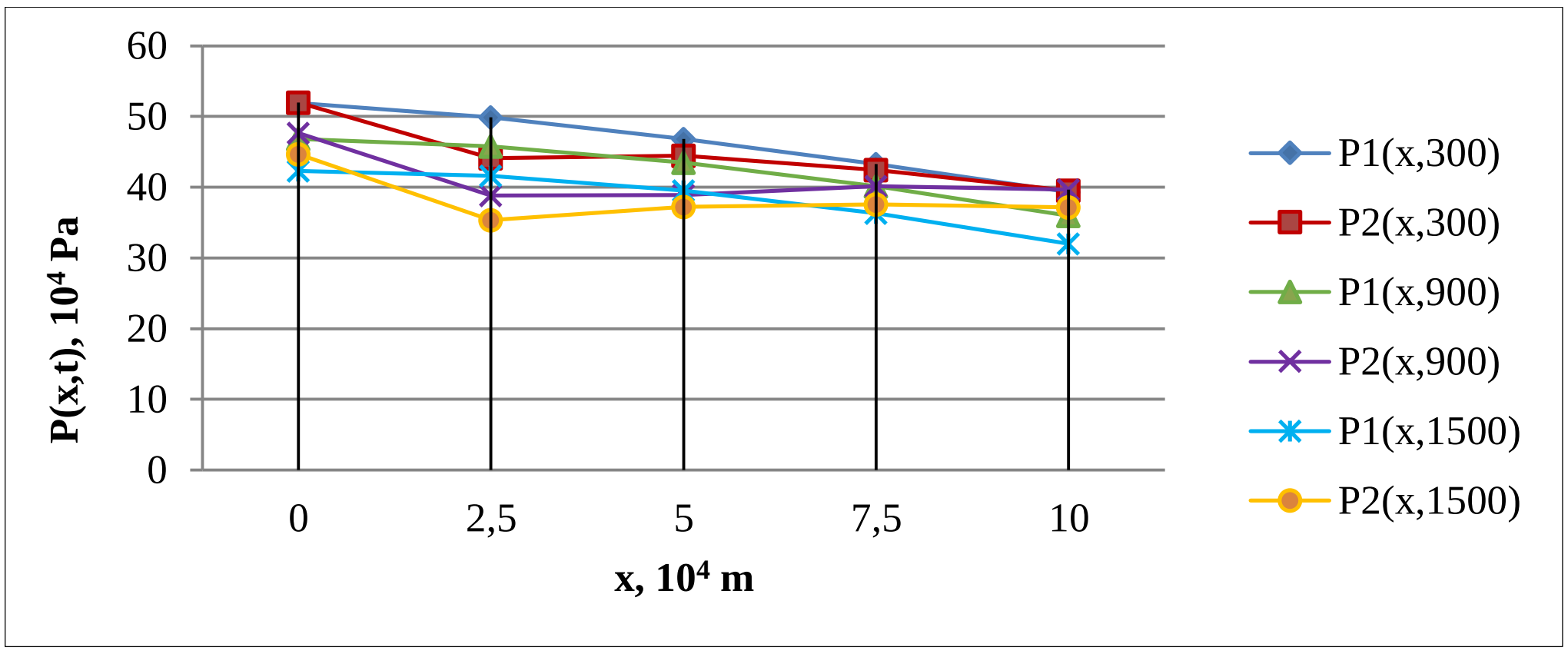

**Graph 1.3.14. Pressure Distribution Along the Length of Damaged and Undamaged Pipeline Sections in the Non-stationary Gas Flow Regime during the Rupture of the 2nd Line and Gas Leakage at Gut = 0.8 G0 in Variant 1 of the Parallel Gas Pipeline.**

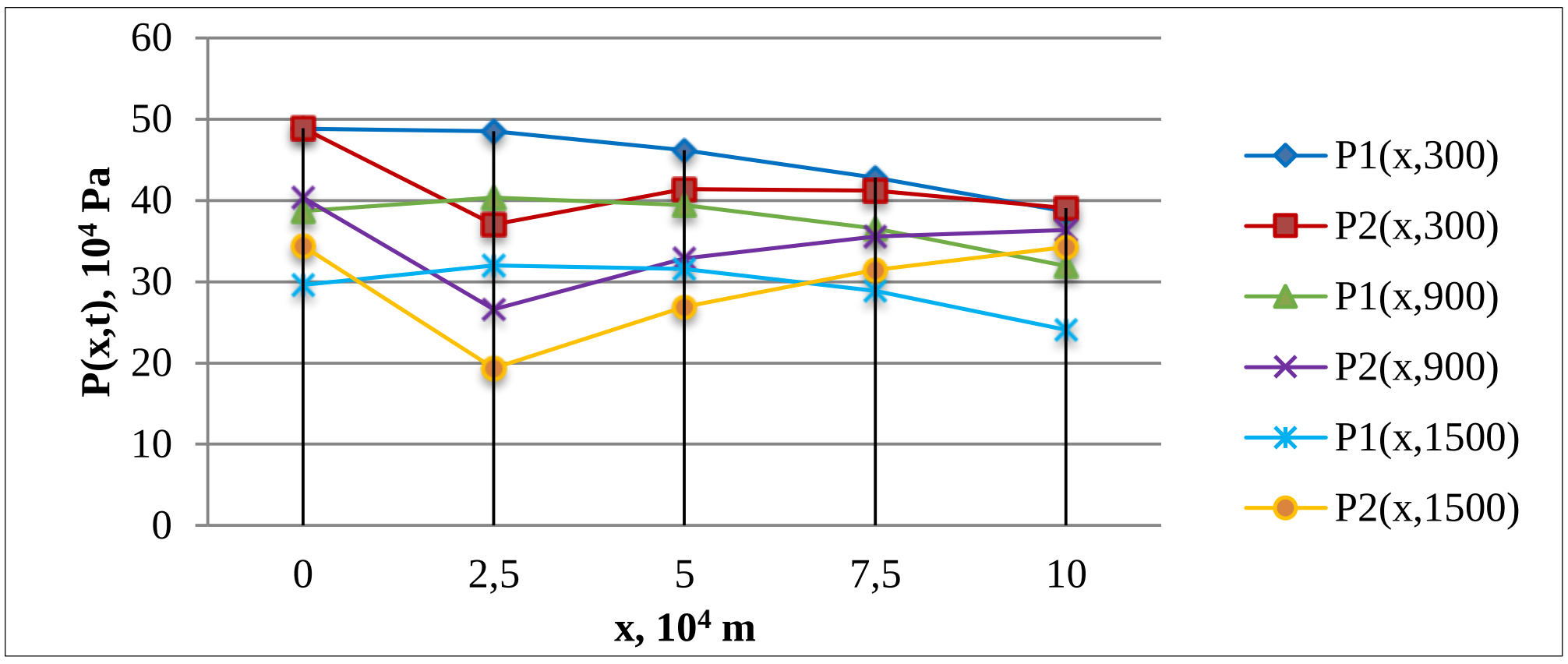

**Graph 1.3.15. Pressure Distribution Along the Length of Damaged and Undamaged Pipeline Sections in the Non-stationary Gas Flow Regime during the Rupture of the 2nd Line and Gas Leakage at Gut = 1.6 G0 in Variant 1 of the Parallel Gas Pipeline.**

From the analysis of Graphs 1.3.14 and 1.3.15, it is apparent that the distribution of pressure along the length of the pipeline does not change depending on the volume of leaked gas mass. Specifically, at any given moment, the pressure value of the undamaged pipeline section at the beginning and end of the gas pipeline is higher than the pressure value of the damaged pipeline section. However, increasing the mass flow of gas leakage leads to a proportional increase in the mentioned pressure ratios over time.

To illustrate our findings clearly, we created the following table using the data from Tables 1.3.12 and 1.3.13 (Table 1.3.14)

| t,sec | The pressure drop ratio of the pipeline at the starting point of the gas pipeline | | The pressure drop ratio of the pipeline at the leakage point of the gas pipeline | | The ratio of pressure drop values of the pipeline at the end point of the gas pipeline | |
|---|---|---|---|---|---|---|
| | $[P_b-P_2(0,t)]/[P_b-P_2(0,t)]$ | | $[P(\ell_2,0)-P_1(\ell_2,t)]/[P(\ell_2,0)-P_2(\ell_2,t)]$ | | $[P_s-P_1(L,t)]/[P_s-P_2(L,t)]$ | |
| | $G_{ut}=0,8\ G_0$ | $G_{ut}=1,6\ G_0$ | $G_{ut}=0,8\ G_0$ | $G_{ut}=1,6\ G_0$ | $G_{ut}=0,8\ G_0$ | $G_{ut}=1,6\ G_0$ |
| 300 | 1,01 | 1,01 | 5,20 | 5,20 | 1,48 | 1,45 |
| 900 | 1,11 | 1,11 | 2,28 | 2,27 | 2,22 | 2,22 |
| 1500 | 1,23 | 1,23 | 1,65 | 1,65 | 2,77 | 2,78 |

**Table 14. Pressure Ratio Over Time in the Initial, Leakage, and Final Sections of Damaged and Undamaged Pipeline Sections during the Rupture of the 2nd Line in Variant 1 of the Gas Pipeline, with Gas Leakage at $G_{ut}=0,8\ G_0$ and $G_{ut}=1,6\ G_0$**

From the analysis of Table 1.3.14, it can be concluded that the ratio of pressure values between the undamaged pipeline sections at the beginning, leakage, and final points does not change depending on the increase in the mass flow of leaked gas. However, the ratio of pressure values at the beginning and end points of the pipeline sections increases over time. Conversely, the ratio of pressure values at the leakage point decreases over time. This indicates that the discharge process of the damaged pipeline section is sluggish and does not allow for effective discharge. These indicators will positively influence the gas supply regime for consumers during non-stationary gas flow, ensuring uninterrupted gas delivery. Thus, the analyzed operating condition of the gas pipeline meets the criteria for seamless gas supply to consumers.

Certainly, in order to confirm our assertions, we need to analyze the distribution of pressure along the damaged and undamaged pipeline sections in the non-stationary gas flow regime, considering that the leakage point is near the end of the parallel gas pipeline (in the 2nd

variant).Let's assume that the 2nd line of the parallel gas pipeline ruptures at a distance of $\ell_2$=75000 m. For this scenario, we determine the variation of pressure over time at points near the beginning of the pipeline, at the point $\ell_2$=25000 m, and at the location of the rupture, and record the results for the different mass flow rates of the leaked gas in the following tables (Table 1.3.15 and Table 1.3.16).

| **$\ell_2$=75000 m** | **$P_1(25000,t)$, $10^4$ Pa** | | | | **$P_2(25000,t)$, $10^4$ Pa** | | | |
|---|---|---|---|---|---|---|---|---|
| t, sec | $G_{ut}$ = 0,4 $G_0$ | $G_{ut}$ = 0,8 $G_0$ | $G_{ut}$ = 1,2 $G_0$ | $G_{ut}$ = 1,6 $G_0$ | $G_{ut}$ = 0,4 $G_0$ | $G_{ut}$ = 0,8 $G_0$ | $G_{ut}$ = 1,2 $G_0$ | $G_{ut}$ = 1,6 $G_0$ |
| 300 | 51 | 50,76 | 50,51 | 50,27 | 50,61 | 49,97 | 49,33 | 48,69 |
| 600 | 50,31 | 49,37 | 48,43 | 47,49 | 49,85 | 48,45 | 47,04 | 45,64 |
| 900 | 49,45 | 47,65 | 45,84 | 44,04 | 49,21 | 47,17 | 45,12 | 43,08 |
| 1200 | 48,51 | 45,77 | 43,03 | 40,29 | 48,66 | 46,07 | 43,48 | 40,89 |
| 1500 | 47,53 | 43,82 | 40,1 | 36,39 | 48,18 | 45,11 | 42,04 | 38,97 |

**Table 1.3.15. Variation of Pressure Over Time at the Point $\ell_2$=25000 m in the Damaged and Undamaged Pipeline Sections during the Rupture of the 2nd Line in Variant 2 of the Gas Pipeline, Depending on Different Mass Flow Rates of Leaked Gas**

| **$\ell_2$=75000 m** | **$P_1(75000,t)$, $10^4$ Pa** | | | | **$P_2(75000,t)$, $10^4$ Pa** | | | |
|---|---|---|---|---|---|---|---|---|
| t, sec | $G_{ut}$ = 0,4 $G_0$ | $G_{ut}$ = 0,8 $G_0$ | $G_{ut}$ = 1,2 $G_0$ | $G_{ut}$ = 1,6 $G_0$ | $G_{ut}$ = 0,4 $G_0$ | $G_{ut}$ = 0,8 $G_0$ | $G_{ut}$ = 1,2 $G_0$ | $G_{ut}$ = 1,6 $G_0$ |
| 300 | 43,07 | 42,38 | 41,7 | 41,01 | 40,19 | 36,62 | 33,06 | 29,5 |
| 600 | 42,06 | 40,37 | 38,69 | 37 | 38,71 | 33,67 | 28,63 | 23,6 |
| 900 | 41,03 | 38,3 | 35,58 | 32,86 | 37,58 | 31,41 | 25,24 | 19,07 |
| 1200 | 39,98 | 36,21 | 32,44 | 28,67 | 36,62 | 29,5 | 22,37 | 15,25 |
| 1500 | 38,93 | 34,11 | 29,29 | 24,47 | 35,78 | 27,82 | 19,85 | 11,88 |

**Table 1.3.16. Variation of Pressure Over Time at the Point $\ell_2$=75000 m (Rupture Point) in the Damaged and Undamaged Pipeline Sections during the Rupture of the 2nd Line in Variant 2 of the Gas Pipeline, Depending on Different Mass Flow Rates of Leaked Gas.**

Using the data from Table 1.3.15 and Table 1.3.16, we can create the following graphs (Graph 1.3.16 and 1.3.17)

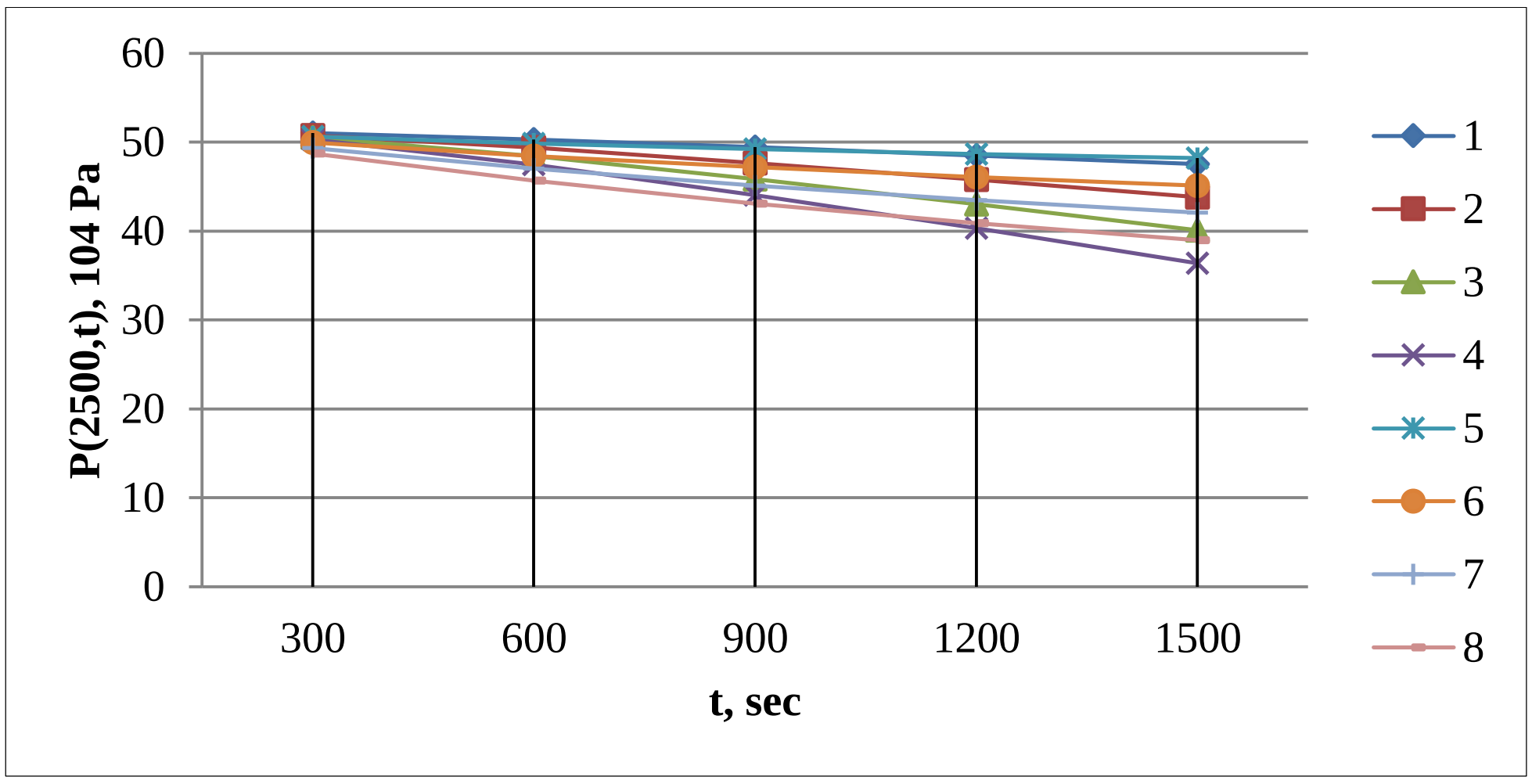

**Graph 1.3.16: It illustrates the variation of pressure over time at the point x=25000 m during the rupture of the 2nd line in Variant 2 of the gas pipeline for both the damaged and undamaged pipeline sections, considering different mass flow rates of leaked gas. Lines 1, 2, 3, and 4 represent the undamaged pipeline, while lines 5, 6, 7, and 8 represent the damaged pipeline.**

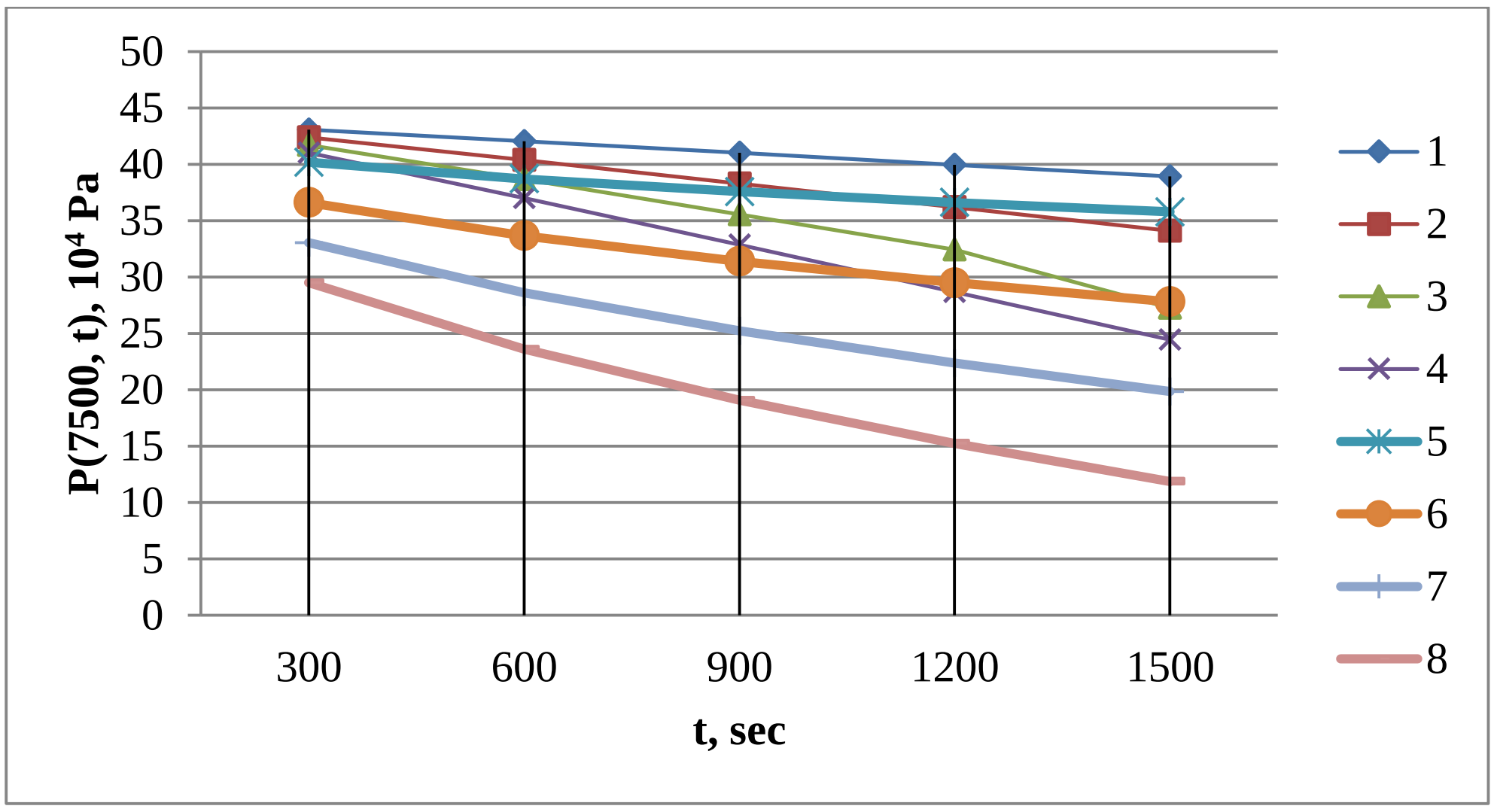

**Graph 1.3.17. The view of the pressure of the damaged and undamaged pipeline segments at x=75000 m point as a function of time, dependent on the values of the escaping gas mass. 1, 2, 3, and 4 - in the undamaged pipeline segment. 5, 6, 7, and 8 - in the damaged pipeline segment.**

From Graph 1.3.16, it is evident that the magnitude of pressure drop in the damaged pipeline segment is several times greater than that in the undamaged pipeline segment. For example, over a duration of t=300 seconds and with a volume of escaping gas mass ($G_{ut}$ =0,4 $G_0$), the pressure drop in the undamaged pipeline segment can be calculated as [$P_1$(25000,0) - $P_2$(25000,300)]= ( 51,25- 51) $\times 10^4$= 0,25$\times 10^4$ Pa, while the pressure drop in the damaged pipeline segment would be [P(25000,0) -$P_1$(25000,300)] = (51,25-50,61) $\times 10^4$=0,64$\times 10^4$ Pa.
So, the value of the pressure drop in the damaged pipeline segment is 2.6 times greater than that in the undamaged pipeline segment. This ratio varies over time for a constant gas leakage rate; at t=600 seconds, it is 1.49; at t=900 seconds, it reduces to 1.13; at t=1200 seconds, it further decreases to 0.95, and finally, at t=1500 seconds, it decreases to 0.83. It is evident that the value of this ratio reaches its maximum at t=300 seconds for a constant gas leakage rate ($G_{ut}$ =0,4 $G_0$) and then decreases over time (Graph 1.3.18).

The value of this ratio remains constant for different fixed gas leakage rates ($G_{ut}$=0,8 $G_0$, $G_{ut}$=1,2 $G_0$ and $G_{ut}$=1,6$G_0$), but it varies over time. At t=300 seconds, it is 2.6; at t=600 seconds, it is 1.49; at t=900 seconds, it decreases to 1.13; at t=1200 seconds, it further reduces to 0.95; and finally, at t=1500 seconds, it decreases to 0.83 (Graph 1.3.2). On the other hand, it appears that regardless of the increase in gas leakage volume, the ratio of the pressure drop in the damaged and undamaged pipeline segments remains constant. Therefore, it can be concluded that the process of complete gas discharge does not occur in the damaged pipeline segment.

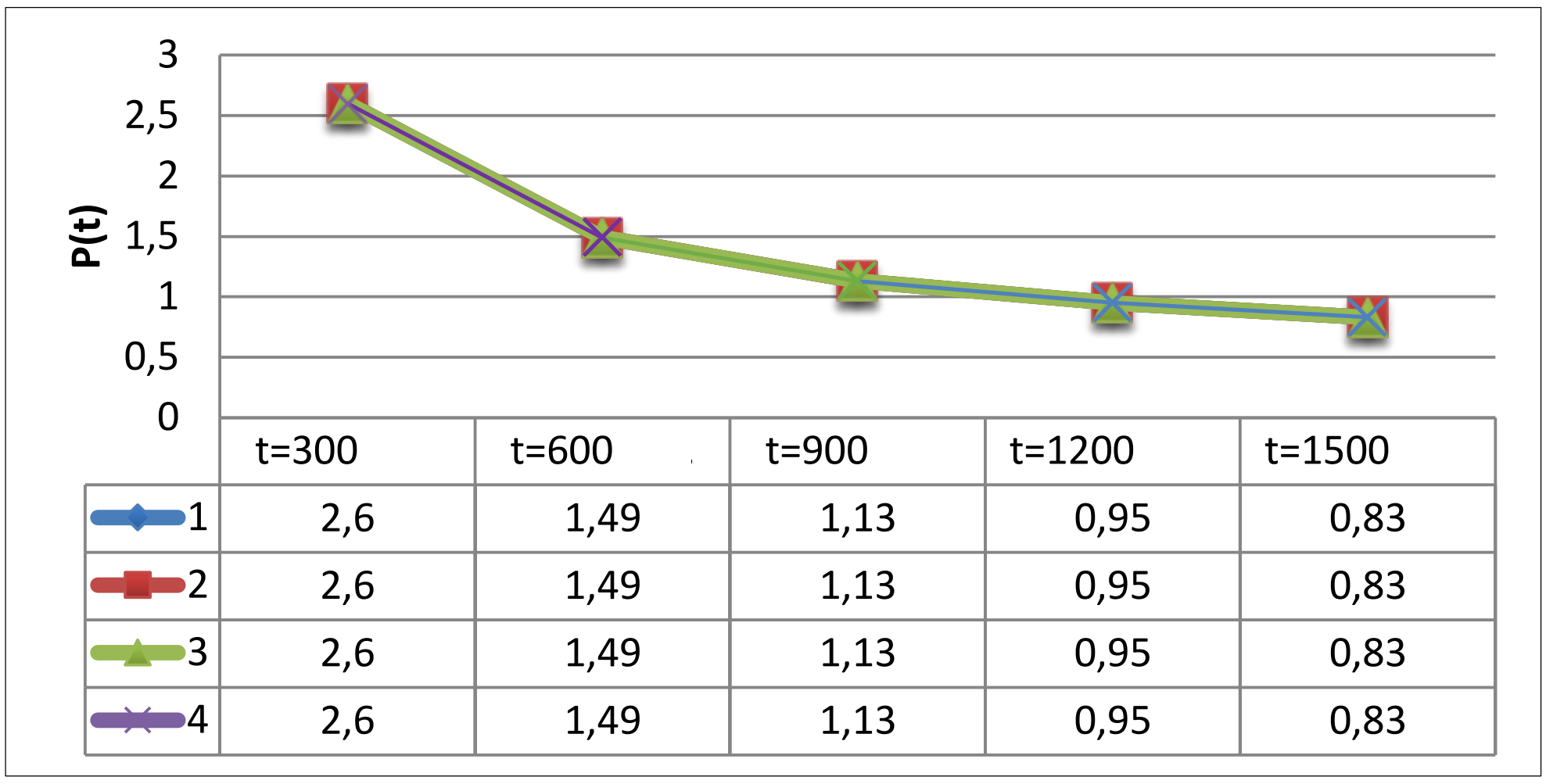

**Graph 1.3.18. Time-dependent graph of the expression [P(x,0)-P₂ (x,t)]/[P(x,0) -P₁ (x,t)]=p(t) at the point x=25000 m along the gas pipeline for different values of gas leakage rates (1-$G_{ut}$=0,4 $G_0$, 2-$G_{ut}$=0,8 $G_0$, 3-$G_{ut}$=1,2 $G_0$, 4-$G_{ut}$ =1,6 $G_0$).**

With the location of the damage to the pipeline at $\ell_2$=75000 m , let's construct the graph of the ratio p(t)=[P (75000,0) -$P_2$(75000,t)]/ [P(75000, 0) –$P_1$(75000,t)] based on the data provided in Table 1.3.16 and Graph 1.3.17.

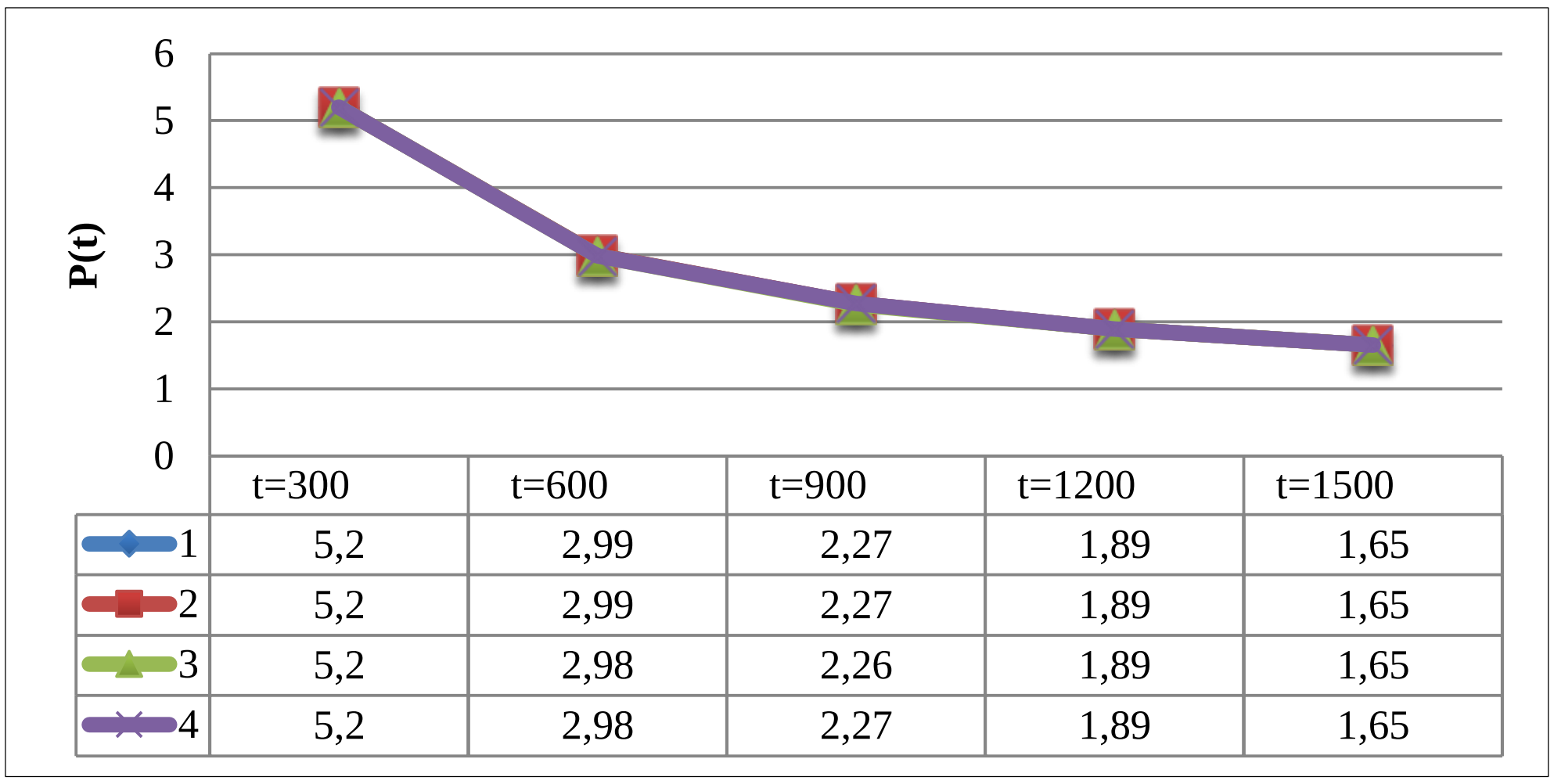

**Graph 1.3.19. The graph of the ratio [Px, 0)-$P_2$ (x,t)]/[P(x, 0) -$P_1$(x,t)]=p(t) as a function of t for the gas pipeline at x=75000 m, with different values of (1-$G_{ut}$=0,4 $G_0$, 2-$G_{ut}$=0,8 $G_0$, 3-$G_{ut}$=1,2 $G_0$, 4-$G_{ut}$ =1,6 $G_0$ ).**

Indeed, as seen from Graph 1.3.19, regardless of the location of the leakage point, the pressure of the gas in both the damaged and undamaged pipelines does not approach zero, indicating that the process of complete emptying does not occur.

The efficiency of the exploitation of main gas pipelines primarily depends on the actual parameters of non-stationary regimes, primarily on the transmission capacity of gas pipelines and the dynamic condition of pipelines affecting the gas supply regime of consumers. The results of the analyzed exploitation conditions show that the change in gas pressure along the length of the gas flow and its portion does not vary significantly in a wide range, as in other exploitation conditions. For example, in the non-stationary regimes of the first variant of the gas pipeline ($\ell_2$=75000 m), the difference in pressure at the start and leakage points of the damaged pipeline is as follows:

1. At the end of the gas pipeline, maintaining constant pressure in the pipelines under exploitation conditions: for $\mathbf{G_{ut}=0,8\ G_0}$ and t=1500 seconds (based on Table 1.3.4),
   $P_1(0,t)-P_1(\ell_2,t)=(36,33-23,3)\times10^4=13,03\times10^4$ Pa
2. At the beginning of the gas pipeline, maintaining constant pressure in the pipelines under exploitation conditions: for $G_{ut}=1,2\ G_0$ and t=1500 seconds (based on Table 1.3.6),
   $P_2(0,t)-P_2(\ell_2,t)=(55-34,43)\times10^4=20,57\times10^4$ Pa
3. At the beginning and end points of the gas pipeline, measuring the consumption of pipelines under exploitation conditions: for . $G_{ut}=1,2\ G_0$ and t=1500 seconds (based on Table 1.3.12),
4. $P_2(0,t)-P_2(\ell_2,t)=(44,7-35,32)\times10^4=9,5\times10^4$ Pa

Certainly, it is known that the wide range of pressure variations at the beginning and leakage points of the pipeline is not regulated by the ratio of the leakage amount to the initial consumption value. Therefore, based on the differences in pressure values, it can be said that the change in the consumption of pipelines at the beginning and end points of the gas pipeline under exploitation conditions results in at least 1.37 times less loss of gas to the environment compared to other exploitation conditions.

Thus, the change in the consumption of pipelines at the beginning and end points of the gas pipeline, measured in exploitation conditions, is superior to other exploitation conditions based on both the reliability of gas supply to consumers and the lower amount of gas leakage into the environment. Therefore, it is necessary to take these characteristics into account during the exploitation and reconstruction of gas pipelines.

**Results:**

1. Based on theoretical research, the mathematical models of non-stationary regimes in pipeline systems have been studied, and the methodology and algorithms for their application have been proposed. According to theoretical investigations, the laws of pressure variation over time along the length of gas transmission systems under different operating conditions during the non-stationary gas flow generated by the rupture of parallel gas pipelines have been determined.
2. Using the analysis and generalization methodology of existing materials, a methodology for calculating the non-stationary operating regimes of pipelines has been proposed. The proposed methodology allows for selecting optimal operating conditions by performing various calculations in accident regimes, ensuring reliable gas supply to consumers, and minimizing gas losses according to the criteria of gas efficiency and gas loss.
3. The conditions for selecting rational operating regimes based on the obtained analytical model of non-stationary gas flow determine new trends in the development of system reconstruction technology. These conditions contribute not only to the acceptance of fuel resulting from the reconstruction of parallel gas pipelines but also to reducing damage to consumer facilities, thereby increasing the reliability of gas supply.

# Chapter 2. Investigation of Reconstruction Methods for Gas Pipelines with Various Configurations to Improve Technological Processes

In +this work, the initial investigation focused on the methods for solving the system of differential equations that describe non-stationary gas flow along the pipeline with a constant cross-section, taking into account the relief of the pipeline’s linear section. A methodology for calculating non-stationary gas flow under steady relief conditions in pipelines was developed. The study examined the impact of the pipeline's linear profile and gravitational forces on flow characteristics, analyzing the effect of relief on the gas-dynamic process. The gas pressure along the length of the pipeline over time was calculated for various scenarios of gas leakage, and the results were analyzed. Consequently, conclusions were drawn on which sections of the pipeline should be prioritized for loop installation to enhance the productivity of the gas pipeline through reconstruction, and these conclusions were substantiated by theoretical investigation.

The goal of the study is to develop a probabilistic analytical method that will select the most efficient option for reconstructing ring gas pipelines. Existing methods for calculating the efficiency of ring networks are summarized, and methods for calculating the ring gas pipeline network under both steady-state and unsteady-state conditions are presented. As a result, a calculation scheme has been developed to determine the compliance of the ring gas networks' transmission capacity with the required reliability of the network. Additionally, the scheme justifies the acceptance of the hydraulic junction of semi-ring gas networks as an energy source for new consumers. The gas-dynamic condition of the entire length of the ring gas pipeline network with a flow rate-based zone along the road in various flow regimes, including unsteady and steady-state, has been analyzed according to the scheme's requirements.

One of the important problems remains to make the operating mode of consumer enterprises independent of the operating mode of the gas distribution pipeline feeding them with gas fuel, in other words, to ensure uninterrupted gas supply to these enterprises. One of the solutions to the problem of improving the reliability of gas pipelines is the use of new effective scientifically based technologies for the reconstruction of pipeline systems. Obviously, connecting fittings, lupings and connectors are one of the components of increasing the reliability of gas networks. To this end, a theoretically sound reporting scheme was developed for the design of gas networks, as well as for the development of new technological foundations of the gas pipeline network during its reconstruction and the prevention of possible gas losses due to modern equipment installed on these networks. This article explores the analytical expression for determining the economic length of interconnections between parallel gas pipelines and outlines a reconstruction variant for the effective methods of reconstructing multi-lines parallel gas pipelines.

As a result of generalizing existing modern methods, an algorithm has been developed to detect an efficient management of the operation of the circular gas pipeline in emergency conditions. A Calculation scheme has been developed to determine the compliance of the ring gas networks' transmission capacity with the required reliability of the network. The gas-dynamic condition of the entire length of the ring gas pipeline network with a flow rate-based zone along the road in various flow regimes, including unsteady and steady-state, has been analyzed in line with the scheme's requirements.

This paper presents innovative methodologies aimed at the effective operation, management, and optimization of gas pipelines and is of significant importance for solving modern engineering problems. The findings of the study present engineering solutions aimed at addressing the challenges in real-world applications by comparing the performance of various algorithms. Consequently, this work contributes to the advancement of cutting-edge approaches in the field of engineering and opens new perspectives for future research.

The management systems for the operational processes of existing parallel gas pipelines under non-stationary conditions are inefficient and technically complex. Additionally, the outdated elements, which are physically worn and do not meet the requirements of modern control technologies, fail to react quickly enough and only report problems after a significant amount of gas is lost to the environment. For this reason, it is crucial to develop new technological schemes that adopt modern control methods in the reconstruction phase of existing pipeline systems in order to eliminate negative outcomes.

In this work, based on the analysis of the mathematical model of non-stationary gas flow in parallel gas pipelines, an algorithm based on machine learning technology was developed for the information database of the control and management point of technological processes. The developed algorithm can be recommended for solving engineering technology problems related to the prediction of technical maintenance and management of emergency situations in pipelines. Additionally, in order to ensure the continuity of gas supply to consumers during emergency situations, an empirical formula for determining the locations of connecting pipes in real-time mode based on the analysis of data from pressure sensors was derived. The application of this formula at the central control point allows for the implementation of a remote activation system for the connecting pipes.

The analysis reveals that the pressure variations at the beginning and end of the parallel pipeline are directly dependent on the activation time of the automatic valves. This serves as the basis for obtaining information at the control point regarding the activation of the valves on the pipeline.

### 2.1. Study of Non-Stationary Gas Flow for Justifying Reconstruction Solutions of Terrain-Based Gas Pipeline Systems

**Introduction:** In accordance with the technological design standards for gas pipelines, the impact of the relief should be considered if there are points along the linear section of the gas pipeline that are more than 100 meters above or below the starting point. The hydraulic calculations for main pipelines are generally based on the condition of stationary gas flow. In practice, however, the movement of gas in pipelines is often non-stationary, meaning it varies over time. Pressure and flow rate vary both along the length of the pipeline and over time. Various factors can cause non-stationary gas movement in pipelines. These factors include the consumption of gas by consumers in volumes exceeding the required amount over a certain period, the startup and shutdown of equipment at compressor stations, the operation of shut-off valves, and accidental leaks in the pipeline. Undoubtedly, studying the non-stationary movement of gas is crucial in cases where the integrity of the relief pipeline fails (accidental leaks).

The topography of the pipeline is of great importance during the reconstruction of gas pipelines, as the elevation difference along the pipeline significantly affects the gas-dynamic characteristics of the flow. It is known that gravitational forces have a greater impact on the flow characteristics in sections of the gas pipeline where the gas density is high. The purpose of this work is to identify reconstruction solutions for relief gas pipelines. To this end, the impact of topography on non-stationary flows resulting from accidental leaks in gas pipelines with constant relief is investigated. Numerous local and international studies on these issues have been reviewed [35,37,38,45,50]. In these works, the authors have explored issues related to stationary regimes of gas pipelines, stationary and non-stationary transport of petroleum products, and non-stationary and non-isothermal regimes of offshore pipelines. The main distinction of our research is the establishment of mathematical expressions for the calculation and analysis of non-stationary regimes resulting from the loss of integrity in relief gas pipeline networks. Mathematical modeling is practically the primary method for studying non-stationary processes.

In many studies dedicated to modeling processes in trunk oil and gas pipelines [18,29,73], the inclination angle of the pipeline's horizontal surface is often not considered. This is because when the endpoint of the gas pipeline is located less than 100 meters above or below the starting point, this component is relatively small and does not create significant changes in the model. However, when dealing with significant gradients along the length of the gas pipeline, both optimal management problems and the necessity to consider the gravitational effect due to the high gas density make it especially important to take into account the inclination angle.

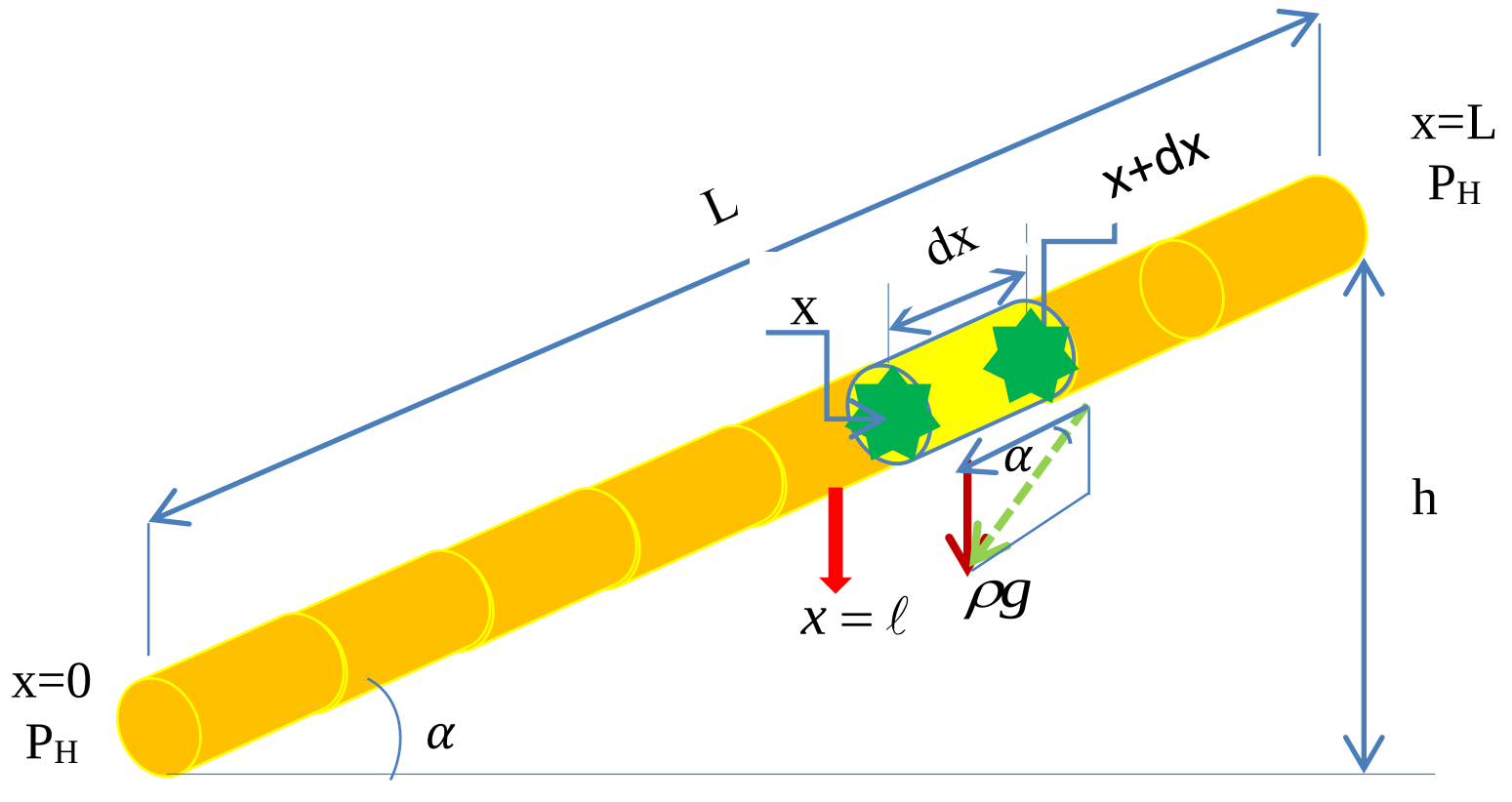

**Figure 2.1.1. Schematic diagram of a relief gas pipeline**

Here,

h: Vertical height difference between the starting and ending sections of the relief gas pipeline (h > 100 m).

α: Inclination angle of the ascending section of the gas pipeline's horizontal surface, degrees.

L: Length of the relief gas pipeline, meters.

x=$\ell$-Location of the breach in the gas pipeline network, meters (m)..

$P_H$ - Constant pressure along the length of the gas pipeline, in Pascals (Pa).

x=0 and x=L: Correspond to the starting and ending sections of the gas pipeline, respectively.

***Analysis of the Mathematical Model for Determining the Methodology of Calculating Non-Stationary Gas Flow in Relief Pipelines.***

Considering the topography of the linear section of the gas pipeline, developing mathematical models for calculating the non-stationary flow regime of the gas pipeline is relevant. To obtain these mathematical models, let's examine the investigated physical processes in more detail. For example, a gas leak with a flow rate of $G_{ut}$ occurs at the point x= $\ell$ of the pipeline. Such an accident is characterized by a change in pressure and the continuity of the flow (Figure 2.1.1). In this case, the isothermal non-stationary movement in pipes with a constant cross-section is described by the following system of non-linear partial differential equations:

- Continuity equations; - Flow motion equations.

$$\begin{cases} -\dfrac{\partial P}{\partial X} = \lambda \dfrac{\rho V^2}{2d} + \rho g Sin\alpha \\ -\dfrac{1}{c^2}\dfrac{\partial P}{\partial t} = \dfrac{\partial G}{\partial X} + 2aG_{ut}(t)\delta(x-\ell) \end{cases} \tag{1}$$

$$\rho = \frac{P}{zRT} = \frac{P}{c^2}$$

It is clear that the system of non-linear equations (1) cannot be solved directly, as such equations cannot be integrated. In engineering calculations, approximate methods are used. In other words, the equations are linearized. To solve the equations, some authors [39] effectively use the modeling of processes in hydraulic gas transport systems with intermediate elements based on I.A. Charny's motion equations. By appropriately simplifying the characteristics of these technological elements, it becomes possible to transition to a closed impulse system in a part of the non-linear system with varying parameters. The expression characterizing this linearization is as follows [78].

$$\begin{cases} -\dfrac{\partial P}{\partial X} = 2aG + \rho g Sin\alpha \\ -\dfrac{1}{c^2}\dfrac{\partial P_2}{\partial t} = \dfrac{\partial G_2}{\partial X} + 2aG_{ut}\delta(x-\ell) \end{cases} \tag{2}$$

Here,

$2a = \lambda \dfrac{V}{2d}$ - Charny's linearization coefficient;

G=ρv mass flow rate of the gas flow $\frac{Pa \cdot \sec}{m}$;

$G_{ut}$ - value of the mass flow rates measured at the leakage point of the pipeline at different times, $\frac{Pa \cdot \sec}{m}$;

d–internal diameter of the pipe m;

$\delta(x-\ell)$ - Dirac function
x –coordinate aligned with the pipe axis and directed along the gas flow m;
P –absolute average gas pressure at the cross-sections of the gas pipeline Pa;
v –average velocity of gas at the cross-section m/sec;
ρ –gas density, kq/$m^3$;
α– inclination angle of the pipe from the horizontal, degrees
t –time, sec;
λ –hydraulic resistance coefficient of the gas pipeline section, dimensionless
c –speed of sound in gas m/sec;
g –acceleration due to gravity, 9,81 m/$sec^2$

If we substitute x-from the first equation and assign $\frac{\partial G}{\partial X}$ from the second equation, then the equation will be as follows.

$$\begin{cases} -\frac{\partial^2 P}{\partial X^2} = 2a\frac{\partial G}{\partial X} + \frac{\partial P}{\partial X}\frac{g}{c^2}Sin\alpha \\ \frac{\partial G}{\partial X} = -\frac{1}{c^2}\frac{\partial P}{\partial t} - 2aG_{ut}\delta(x-\ell) \end{cases} \quad (3)$$

Then in equation (3), if we substitute the expression of $\frac{\partial G}{\partial X}$ from the first equation, we obtain the following equation for P(x,t) pressure, which is essential for solving the mathematical problem.

$$\frac{\partial^2 P}{\partial X^2} = \frac{2a}{c^2}\frac{\partial P}{\partial t} - \frac{\partial P}{\partial X}\frac{g}{c^2}Sin\alpha + 2aG_{ut}\delta(x-\ell) \quad (4)$$

Thus, ensuring the solution of the (4) tenancy reflecting the learning process, understanding how the non-stationary process will change, and also having information about the initial and boundary conditions are necessary to fully determine the given process. This means having information about the indicators of the initial state of the flow parameters (at t = 0 time) and what happens at the boundaries of the gas chamber (at t > 0, x = 0, and x = L divisions). Therefore, we accept the initial and boundary conditions (the values of the expenditures measured at the initial and final points are taken as constants, and the value of this constant is assumed to be equal to "0") as follows: Suppose that in the exploitation of the relief gas chamber, the pressure is kept constant along the length of the chamber. Then,

t=0; P(x,0)= $P_H$

The requirement of boundary conditions being equal to zero enables obtaining the comprehensive solution of the analytical model, but it does not restrict the applicability of the model for abrupt changes in gas flow or pressure. Therefore, for t > 0, the boundary conditions at the initial and final points of the gas chamber, under the condition of equality of expenditures measured, are as follows:.

$$\left.\frac{\partial P(x,t)}{\partial X}\right|_{x=0} = \left.\frac{\partial P(x,t)}{\partial X}\right|_{x=L} = 0$$

The (4) expression can be presented in a more comprehensive format

$$\frac{\partial^2 P}{\partial X^2} = k\frac{\partial P}{\partial t} - k_1\frac{\partial P}{\partial X} + 2aG_{ut}\delta(x-\ell) \qquad (5)$$

Here, k= $\frac{2a}{c^2}$ k1= $\frac{g}{c^2}Sin\alpha$

The problem statement encompasses the consideration of the non-stationary nature of a relief gas chamber along the length of a pipeline for resolving the reconstruction of the pipeline. For this purpose, we accept the following equality to convert the (5) tenancy into a heat transfer tenancy:

$$P(x,t) = U(x,t)e^{\lambda_1 t}e^{\lambda_2 x} \qquad (6)$$

Considering the (6) expressions in terms of (5), we define the heat transfer tenancy for analyzing the non-stationary nature of the linear section of gas chambers as follows:.

$$\frac{\partial^2 U}{\partial X^2} = k\frac{\partial U}{\partial t} + 2aG_{ut}\delta(x-\ell)e^{-\lambda_1 t}e^{-\lambda_2 x} \qquad (7)$$

Here, $\lambda_1 = -\frac{g^2 Sin^2\alpha}{2a4c^2}$; $\lambda_2 = -\frac{gSin\alpha}{2c^2} \rightarrow$ k1= - $2\lambda_2$

Therefore, the initial and boundary conditions will be determined as follows:

t=0; $U(x,t) = P_H e^{-\lambda_2 x}$

t>0; $\left.\frac{\partial U(x,t)}{\partial X} + \lambda_2 U(x,t)\right|_{x=0} = \left.\frac{\partial U(x,t)}{\partial X} + \lambda_2 U(x,t)\right|_{x=L} = 0$

It is understood that applying the Laplace transform [U(x,s)= $\int_0^\infty U(x,t)e^{-st}dt$ ] to equation (7) converts it into a second-order ordinary differential equation, and its general solutions will be as follows:.

$$\frac{d^2U(x,s)}{dx^2} - ksU(x,s) = 2aG_{ut}\delta(x-\ell_2)\frac{e^{-\lambda_2 x}}{s+\lambda_1} - kP_H e^{-\lambda_2 x}$$

$$\frac{d^2U(x,s)}{dx^2} - ksU(x,s)=0$$

$$U(x,S) = C_1 e^{\sqrt{kS}x} + C_2 e^{-\sqrt{kS}x}$$

Then,

$$U(x,s) = \frac{P_H e^{-\lambda_2 x}}{s+\lambda_1} - \frac{P_H}{s+\lambda_1}\left[ch\sqrt{ks}x - \frac{\lambda_2}{\sqrt{ks}} sh\sqrt{ks}x\right] + \\ + \frac{2a}{\sqrt{ks}}\frac{G_{ut}}{s+\lambda_1} e^{-\lambda_2 \ell} sh\sqrt{ks}(x-\ell) + C_1 e^{\sqrt{kS}x} + C_2 e^{-\sqrt{kS}x} \quad (8)$$

Applying the Laplace transform to the initial and boundary conditions, considering equation (8), we determine the constants $C_1$ and $C_2$

$$C_1 = \frac{kP_H}{2\left(\sqrt{ks}-\lambda_2\right)\sqrt{ks}} - \frac{2ae^{-\lambda_2 \ell} G_{ut}}{2(s+\lambda_1)} \frac{ch\sqrt{ks}(L-\ell) + \frac{\lambda_2}{\sqrt{ks}} sh\sqrt{ks}(L-\ell)}{\left(\sqrt{ks}-\lambda_2\right) sh\sqrt{ks}L} \quad (9)$$

$$C_2 = \frac{kP_H}{2\left(\sqrt{ks}+\lambda_2\right)\sqrt{ks}} - \frac{2ae^{-\lambda_2 \ell} G_{ut}}{2(s+\lambda_1)} \frac{ch\sqrt{ks}(L-\ell) + \frac{\lambda_2}{\sqrt{ks}} sh\sqrt{ks}(L-\ell)}{\left(\sqrt{ks}+\lambda_2\right) sh\sqrt{ks}L} \quad (10)$$

We obtain the values of equations (9) and (10) by considering equation (8) and the transformed tenancy of the dynamic state of the relief gas chamber, in other words, the distribution of the pressure profile along the length of the relief gas chamber.

$$U(x,s) = \frac{P_H e^{-\lambda_2 x}}{s+\lambda_1} + \frac{2a}{\sqrt{ks}}\frac{G_{ut}}{s+\lambda_1} e^{-\lambda_2 \ell} sh\sqrt{ks}(x-\ell) - \frac{2ae^{-\lambda_2 \ell}}{k(s+\lambda_1)^2} G_{ut} \frac{\lambda_2 sh\sqrt{ks}(L-\ell-x)}{sh\sqrt{ks}L} \\ - \frac{2ae^{-\lambda_2 \ell}}{k(s+\lambda_1)^2} G_{ut} \frac{\left[ch\sqrt{ks}(L-\ell+x) + ch\sqrt{ks}(L-\ell-x)\right]}{sh\sqrt{ks}L} \frac{\sqrt{ks}}{2} \\ + \frac{2ae^{-\lambda_2 \ell}}{k(s+\lambda_1)^2} G_{ut} \frac{\left[ch\sqrt{ks}(L-\ell+x) - ch\sqrt{ks}(L-\ell-x)\right]}{sh\sqrt{ks}L} \frac{\lambda_2^2}{2\sqrt{ks}} \quad (11)$$

Based on the solution to the problem under consideration, determining the mathematical expression of the non-stationary flow of gas in the gas chamber is more appropriate. Therefore, to find the original solutions of equations (11), we use the inverse Laplace transform rule for each boundary, obtaining the following original equation:.

$$U(x,t) = P_H e^{-\lambda_2 x} e^{-\lambda_1 t} - \frac{c^2 G_{ut}}{2L} e^{-\lambda_2 \ell} e^{-\lambda_1 t}\left[\frac{\left(e^{-\lambda_1 t}-1\right)}{\lambda_1} + 2\sum_{n=1}^{\infty} Cos\frac{\pi n(x-\ell)}{L}\frac{1-e^{-(\alpha-\lambda_1)t}}{\alpha-\lambda_1}\right] + \\ + \frac{c^2 e^{-\lambda_2 \ell}}{L} G_{ut} \sum_{n=1}^{\infty} Cos\frac{\pi n(x+\ell)}{L} e^{-\lambda_1 t}\left[\frac{t}{\alpha-\lambda_1} - \frac{1-e^{-(\alpha-\lambda_1)t}}{(\alpha-\lambda_1)^2}\right]\cdot\left[\alpha + \frac{\pi n c^2 \lambda_2}{2aL} - \frac{c^2\lambda_2^2}{2a}\right] + \\ + \frac{c^4 e^{-\lambda_2 \ell}\lambda_2^2}{2a2L} e^{-\lambda_1 t} G_{ut}\left[\frac{t}{\lambda_1} + \frac{e^{-\lambda_1 t}-1}{\lambda_1^2}\right] \quad (12)$$

If we consider equation (12) in terms of equation (6), according to the initial and boundary conditions, we obtain the following mathematical expression for the non-stationary gas flow pressure distribution along the length of the relief gas chamber as a function of time.

$$P(x,t)=P_H-\frac{c^2G_{ut}}{L}e^{-\lambda_2(\ell-x)}\cdot A \tag{13}$$

Here,

$$A=\left\{\begin{array}{l}\dfrac{c^2\lambda_2^2}{2a2\lambda_1}\left[t+\dfrac{e^{\lambda_1 t}-1}{\lambda_1}\right]+\left[\dfrac{\left(e^{-\lambda_1 t}-1\right)}{2\lambda_1}+\sum\limits_{n=1}^{\infty}Cos\dfrac{\pi n(x-\ell)}{L}\dfrac{1-e^{-(\alpha-\lambda_1)t}}{\alpha-\lambda_1}\right]+\\ +\sum\limits_{n=1}^{\infty}Cos\dfrac{\pi n(x+\ell)}{L}\left[\dfrac{t}{\alpha-\lambda_1}-\dfrac{1-e^{-(\alpha-\lambda_1)t}}{(\alpha-\lambda_1)^2}\right]\cdot\left[\alpha+\dfrac{\pi nc^2\lambda_2}{2aL}-\dfrac{c^2\lambda_2^2}{2a}\right]\end{array}\right\}$$

$$\alpha=\frac{\pi^2n^2c^2}{2aL^2};\qquad \lambda_1=-\frac{g^2Sin^2\alpha}{2a4c^2};\quad \lambda_2=-\frac{gSin\alpha}{2c^2}$$

***Analysis of the impact of relief parameters during gas pipeline leakage.***

Taking into account the topography of the linear section of the gas chamber, let's review the proposed method for determining the non-stationary operating regime of the gas chamber using equation (13), which is designed for identifying technological processes. Pressure control in technological processes is of utmost importance. For this purpose, we adopt the following parameters for the relief gas chamber:

$P_H$ = 14 ×10$^4$ Pa: $2a=0{,}1\dfrac{1}{\sec}$ ; $c=383{,}3\dfrac{m}{\sec}$; L = 3×10$^4$ m: d= 0,7 m: n=1,2,3,…..12.

In the example considered, the linear section of the gas pipeline is depicted with a constant relief from its starting to its ending section. The terrain covered by the pipeline spans a distance of 30,000 m from the starting point to the endpoint, with a vertical rise of 1,000 m (Figure 1).

Considering the data mentioned above in equation (13), we determine the lawful variation of gas pressure along the length and over time, depending on different locations of potential breaches (gas leaks) in the relief gas pipeline ($\ell$=5× 10$^3$ m; $\ell$=15 ×·10$^3$ m ; $\ell$=25 ×10$^3$ m and $G_{ut}$=10 Pa×sec/m )

Firstly, assuming the location of the gas pipeline leak at $\ell$=5×10$^3$ m, we record the values of gas pressure every 50 seconds and every 2.5 km as shown in the following table (Table 2.1.1).

| x, km | t,sec | | | | | | | | | | | |
|---|---|---|---|---|---|---|---|---|---|---|---|---|
| | **50** | **100** | **150** | **200** | **250** | **300** | **350** | **400** | **450** | **500** | **550** | **600** |
| 0 | 13,85 | 13,63 | 13,40 | 13,16 | 12,91 | 12,65 | 12,38 | 12,10 | 11,81 | 11,51 | 11,19 | 10,85 |
| 2500 | 13,84 | 13,60 | 13,35 | 13,09 | 12,83 | 12,55 | 12,26 | 11,96 | 11,64 | 11,31 | 10,97 | 10,61 |
| 5000 | 13,82 | 13,57 | 13,30 | 13,02 | 12,73 | 12,43 | 12,12 | 11,79 | 11,45 | 11,10 | 10,73 | 10,34 |
| 7500 | 13,81 | 13,53 | 13,24 | 12,94 | 12,63 | 12,31 | 11,97 | 11,62 | 11,25 | 10,87 | 10,47 | 10,05 |
| 10000 | 13,80 | 13,49 | 13,18 | 12,86 | 12,52 | 12,18 | 11,81 | 11,43 | 11,04 | 10,62 | 10,19 | 9,74 |
| 12500 | 13,78 | 13,46 | 13,12 | 12,77 | 12,41 | 12,03 | 11,64 | 11,23 | 10,80 | 10,36 | 9,89 | 9,40 |
| 15000 | 13,76 | 13,41 | 13,05 | 12,67 | 12,28 | 11,88 | 11,45 | 11,01 | 10,55 | 10,07 | 9,57 | 9,04 |
| 17500 | 13,74 | 13,37 | 12,97 | 12,57 | 12,15 | 11,71 | 11,25 | 10,78 | 10,28 | 9,76 | 9,22 | 8,65 |
| 20000 | 13,72 | 13,32 | 12,89 | 12,46 | 12,00 | 11,53 | 11,04 | 10,52 | 9,99 | 9,42 | 8,84 | 8,22 |
| 22500 | 13,70 | 13,26 | 12,81 | 12,33 | 11,84 | 11,33 | 10,80 | 10,25 | 9,67 | 9,06 | 8,43 | 7,77 |
| 25000 | 13,68 | 13,20 | 12,71 | 12,20 | 11,67 | 11,12 | 10,55 | 9,95 | 9,33 | 8,67 | 7,99 | 7,28 |
| 27500 | 13,65 | 13,14 | 12,61 | 12,06 | 11,49 | 10,90 | 10,28 | 9,63 | 8,96 | 8,25 | 7,52 | 6,75 |
| 30000 | 13,62 | 13,07 | 12,50 | 11,91 | 11,29 | 10,65 | 9,98 | 9,29 | 8,56 | 7,80 | 7,01 | 6,18 |

**Table 2.1.1. Pressure values over time along the length of the relief gas pipeline where the pipeline integrity is compromised ($\ell=5\cdot10^3$ m).**

Using Table 2.1.1, let's plot the distribution of pressure every 50 seconds along the relief gas pipeline at the starting, ending, and compromised sections (Graph 2.1.1).

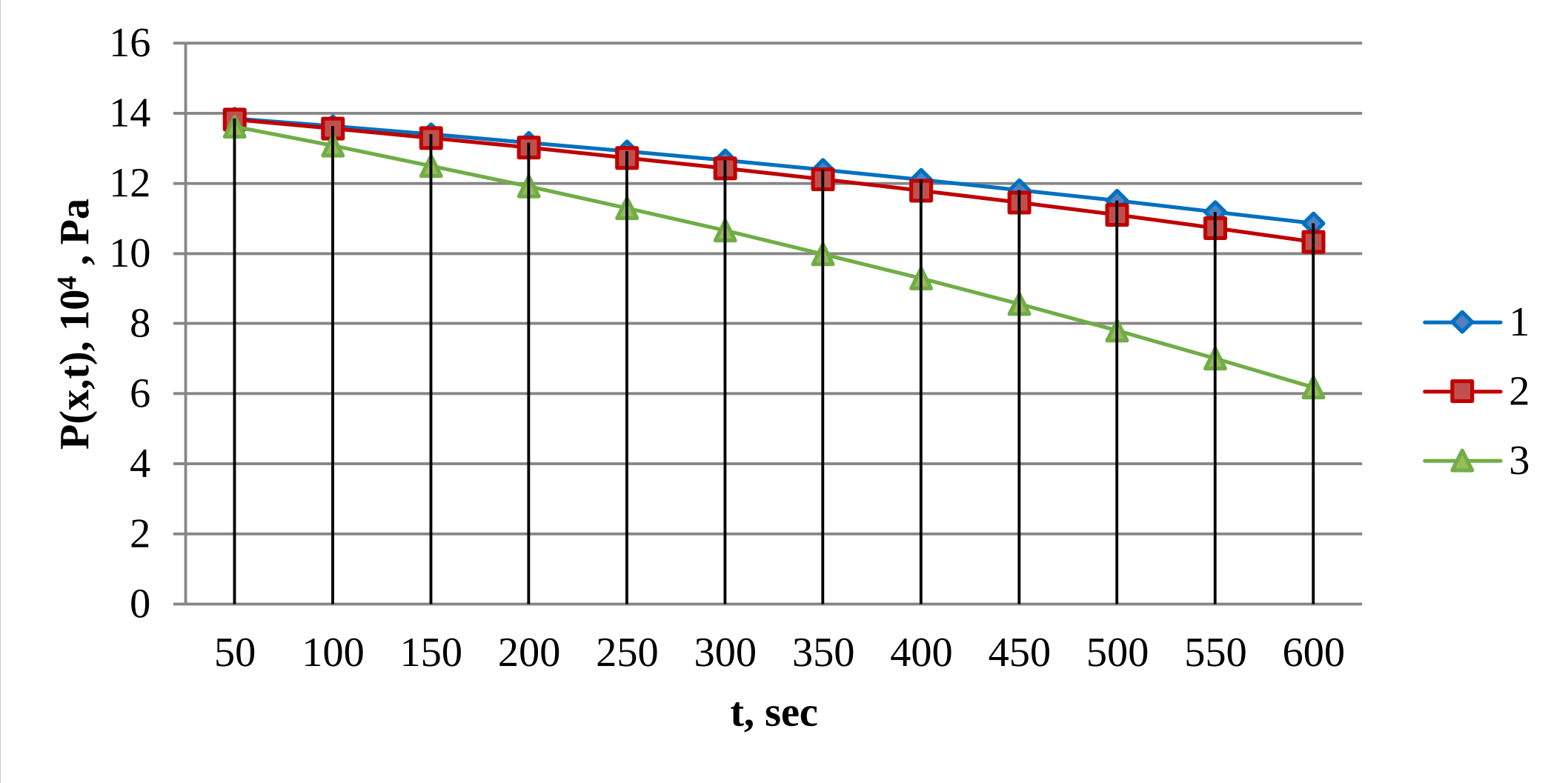

**Graph 2.1.1. Graphical representation of the pressure variation over time in different sections of the non-stationary gas flow in the relief gas pipeline where the pipeline integrity is compromised (at $\ell=5\cdot10^3$ m, 1 - at the starting section - P(0,t), 2 - at the section where the pipeline integrity is compromised - P($\ell$,t) and 3 - at the ending section - P(L,t)**

From the analysis of Graph 2.1.1, the first noticeable gas-dynamic process is that the rate of pressure drop over time at the end of the gas pipeline is greater than at the starting section. This indicates that the location of the pipeline integrity breach is closer to the starting section of the gas pipeline. As seen from Table 2.1.1, if the pressure at the beginning of the gas pipeline decreases from $14\cdot10^4$ Pa to $10{,}85\cdot10^4$ Pa after 600 seconds, the pressure drop at the end of the gas pipeline is from $14\cdot10^4$ Pa to ,18·$10^4$ Pa. Therefore, the pressure drop at the end of the gas pipeline is approximately 2.5 times greater than the pressure drop at the beginning of the gas pipeline.

This is in stark contrast to the gas-dynamic processes occurring in non-relief (linear) gas pipelines. That is, in non-relief gas pipelines, if the location of the pipeline integrity breach is near the beginning section, the rate of pressure drop over time at the beginning is many times greater than the rate of pressure drop at the end of the gas pipeline [7].
The reason for this process direction in the relief gas pipeline is the effect of gravitational force at the beginning of the pipeline. It is known that in a gas pipeline, the relief direction causing the gas flow to move upwards will increase the height difference (h) between the beginning and end sections, thereby increasing the gravitational force. The gas flow, with its weight increased, will tend to move in the opposite direction from top to bottom under the influence of gravitational force. On the other hand, based on the known characteristics of the gas-dynamic process, the overall gas flow in a compromised pipeline tends to move towards the leak area (Figure 2.1.1). Since the gas leak occurs near the beginning of the pipeline, the gas flow at the end of the pipeline will tend to move towards the leak, i.e., downward. Thus, in accordance with both known laws, the gas density at the end of the pipeline will decrease over time due to the counterflow at the end section. As a result of the decrease in gas density, the pressure ($P=\rho c^2$) at the end of the pipeline will also decrease accordingly. Therefore, the noted process is precisely due to the effect of the gas pipeline's relief. Consequently, the dynamics of the non-stationary gas flow resulting from pipeline integrity breaches in relief gas pipelines are entirely opposite to those in linear gas pipelines. This is because of the counterflow created by the gas leak near the beginning of the pipeline and the gravitational force arising from the relief at the start of the gas pipeline (Figure 2.1.1).

Now, to understand how the analyzed process occurs in the linear sections of the relief gas pipeline, let's use Table 2.1.1 to plot the dependence of pressure variation along the pipeline length every 2.5 km and over time every 50 seconds (Graph 2.1.2).

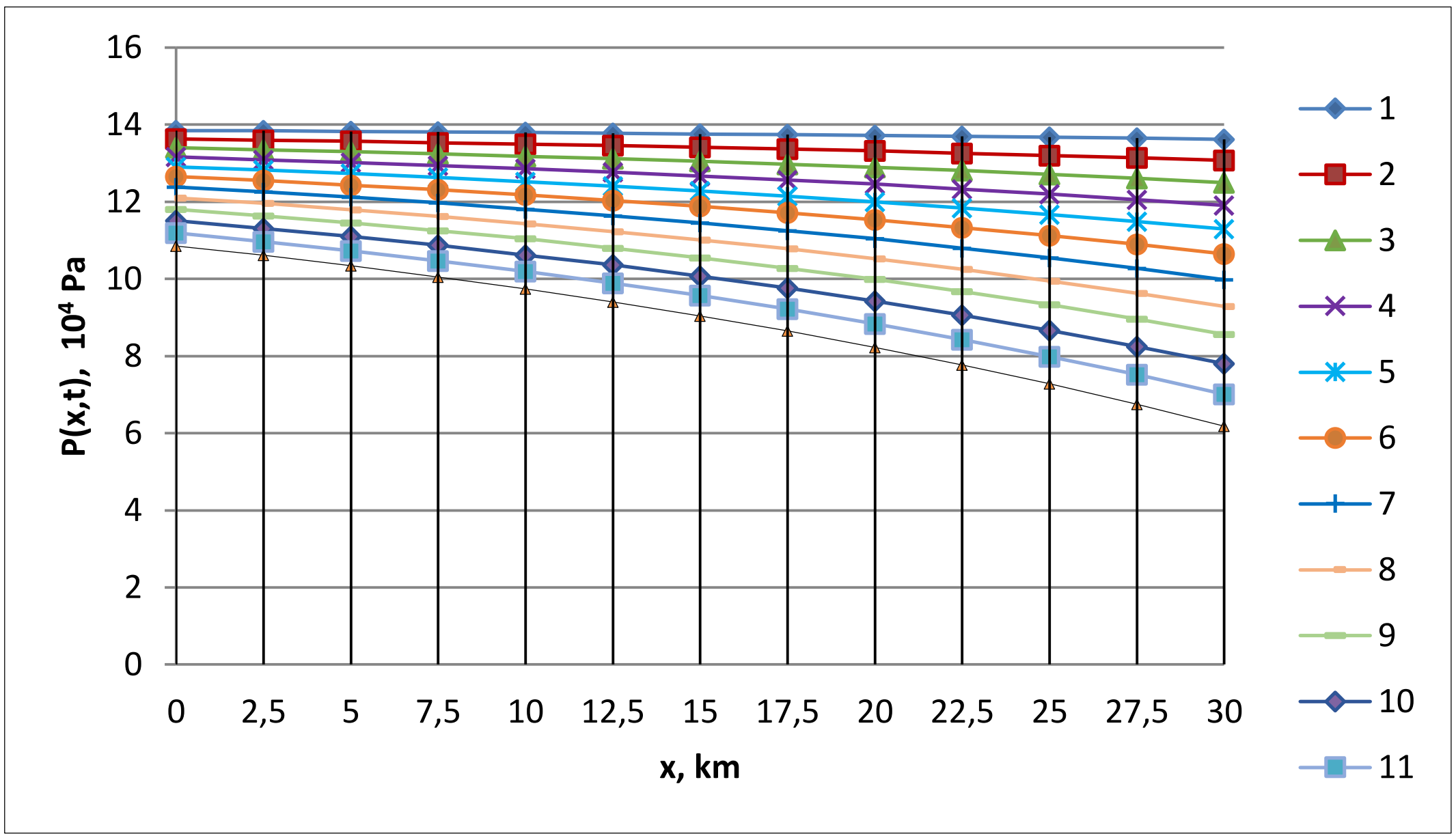

**Graph 2.1.2. The distribution of pressure along the length of the relief gas pipeline with compromised integrity in a non-stationary flow regime over time ($l=5\times10^3$m, 1 - t=50 sec, 2 - t=100 sec, 3 - t= 150 sec, 4 - t=200 sec, 5 - t=250 sec, 6 - t=300 sec, 7 - t=350 sec, 8 - t=400 sec, 9 - t=450 sec, 10 - t=500 sec, 11 - t=550 sec, 12 - t=600 sec).**

From the analysis of Graph 2.1.2, it can be concluded that, as a result of the pipeline integrity breach in the relief gas pipeline, the value of the pressure drop along the length increases from the starting point to the end point over time.

Based on Table 2.1.2, the value of the pressure drop at the beginning section of the pipeline after 50 seconds is (14−13.85)× $10^4$ =0.25×$10^4$Pa, after 100 seconds is (14−13.63)×$10^4$=0.37×$10^4$ Pa, after 150 seconds is (14−13.40)× $10^4$ =0.6×$10^4$ Pa, after 200 seconds is (14−13.16)×$10^4$=0.84×$10^4$ Pa, after 250 seconds is (14−12.91)× $10^4$ =1.09×$10^4$ Pa, after 300 seconds is (14−12.65)× $10^4$ =1.35×$10^4$ Pa, and thus, after 600 seconds, this value increases to (14−10.85)× $10^4$ =3.15×$10^4$ Pa. On the other hand, after t=300 seconds:

at x=0 point: (14−12.65)× $10^4$=1.35×$10^4$Pa,
at x=2500 m point: (14−12.55)× $10^4$=1.45×$10^4$ Pa,
at x=5000 m point: (14−12.43)× $10^4$=1.57×$10^4$ Pa,
at x=7500 m point: (14−12.31)× $10^4$=1.69×$10^4$ Pa,
at x=10000 m point: (14−12.18)× $10^4$=1.92×$10^4$ Pa,
at x=12500 m point: (14−12.03)× $10^4$=1.97×$10^4$ Pa,
at x=15000 m point: (14−11.88)× $10^4$=2.12×$10^4$ Pa,
at x=17500 m point: (14−11.71)× $10^4$=2.29×$10^4$ Pa,

thus, increasing to x=30000 m point: (14−10.65)× $10^4$=3.35×$10^4$ Pa at the end section of the gas pipeline

As a result of these analyses, it can be noted that the value of the pressure drop increases from the beginning section to the end section of the gas pipeline, and these values also change in an increasing direction over time. However, the ratio of the pressure drop value at the beginning section of the pipeline to the pressure drop value at any point along its length does not change over time (Graph 2.1.3).

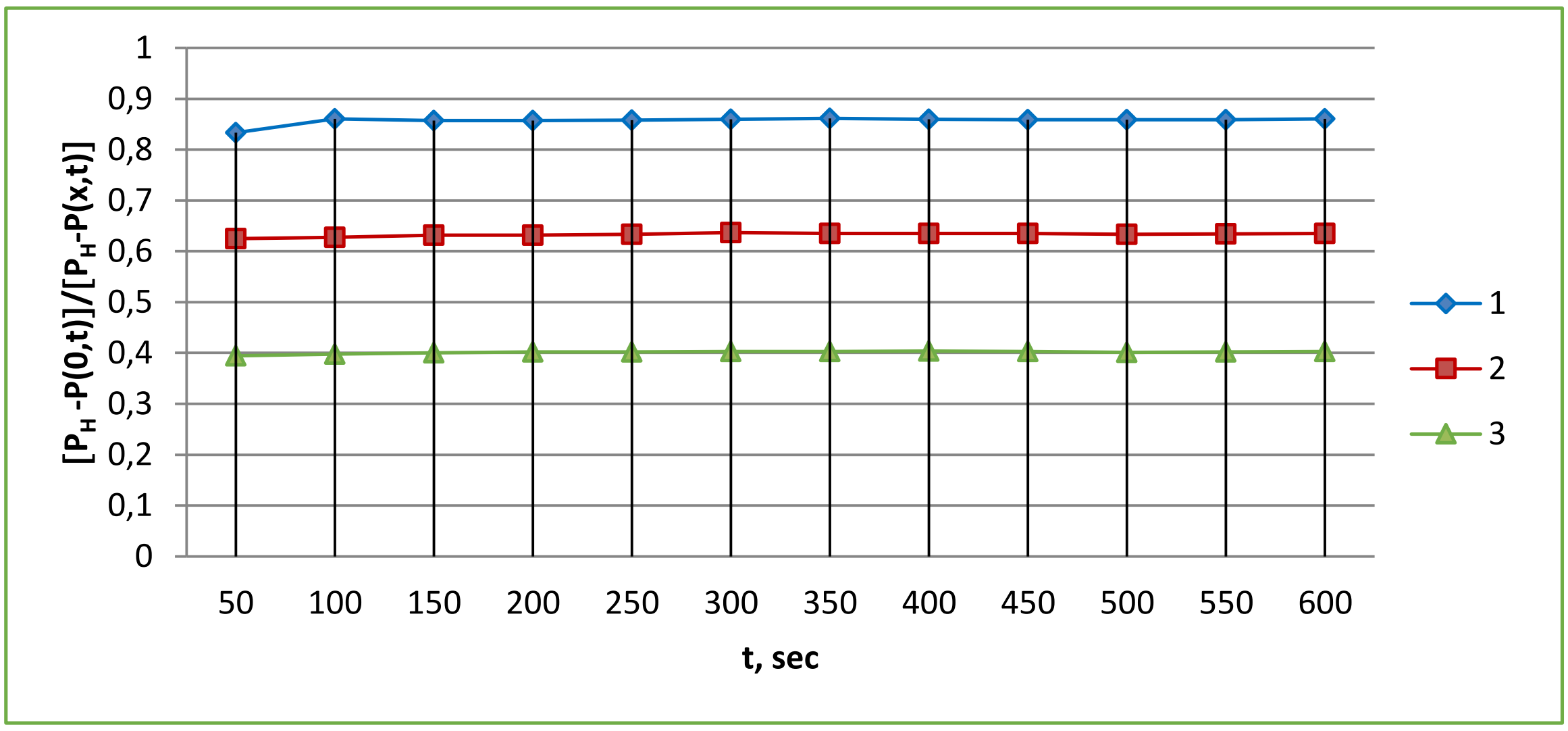

**Graph 2.1.3. The graph of the ratio of the pressure drop value at the beginning section of the relief gas pipeline with a breach to the pressure drop value at different sections over time (l=5×$10^3$m, 1 - at the beginning section - x=5 km, 2 - in the middle section of the gas pipeline - x=15 km, 3 - at the end section - x=30 km)**

From the analysis of Graph 2.1.3, it is clear that the ratio of the pressure drop value at the beginning section of the gas pipeline to the pressure drop values at the middle and end sections remains constant over time. However, according to Graph 2.1.3, the ratio of the pressure drop value at the beginning section of the gas pipeline to the pressure drop value at the leakage section reaches its minimum value (0.83×$10^4$ Pa) in the first 50 seconds, increases to (0.86×$10^4$ Pa) by t=100 seconds, and remains constant thereafter.

This variability may be due to the rounding of numbers during calculation. That is, the ratio of pressure drops at that section, and at that moment, it is constant. This can be explained by the fact that since the gradient of the relief along the length of the gas pipeline is constant, the effect of gravitational force along the pipeline is also constant. It should also be noted that in a relief-free (linear) gas pipeline, the ratio of the pressure drop value at the beginning section to the pressure drop values at different sections changes according to an extremum function, reaching maximum or minimum values depending on the location of the gas leak. However, if the gas leak occurs exactly in the middle section of the linear gas pipeline, only in this case does the ratio of the pressure drops at the beginning and end sections remain constant [7].

Now, by considering the leakage point in the relief gas pipeline to be at $\ell = 15\times10^3$ m, we calculate the pressure values of the gas pipeline every 50 seconds and at every 2.5 km, and record them in the following table (Table 2.1.2).

| x, km | t,sec | | | | | | | | | | | |
|---|---|---|---|---|---|---|---|---|---|---|---|---|
| | 50 | 100 | 150 | 200 | 250 | 300 | 350 | 400 | 450 | 500 | 550 | 600 |
| 0 | 13,89 | 13,72 | 13,56 | 13,38 | 13,20 | 13,01 | 12,81 | 12,60 | 12,39 | 12,16 | 11,92 | 11,68 |
| 2500 | 13,88 | 13,70 | 13,52 | 13,33 | 13,13 | 12,93 | 12,71 | 12,49 | 12,26 | 12,01 | 11,76 | 11,49 |
| 5000 | 13,87 | 13,68 | 13,48 | 13,28 | 13,06 | 12,84 | 12,61 | 12,37 | 12,12 | 11,86 | 11,58 | 11,30 |
| 7500 | 13,86 | 13,65 | 13,44 | 13,22 | 12,99 | 12,75 | 12,50 | 12,24 | 11,97 | 11,69 | 11,39 | 11,08 |
| 10000 | 13,85 | 13,63 | 13,40 | 13,16 | 12,91 | 12,65 | 12,38 | 12,10 | 11,81 | 11,51 | 11,19 | 10,85 |
| 12500 | 13,84 | 13,60 | 13,35 | 13,09 | 12,83 | 12,55 | 12,26 | 11,96 | 11,64 | 11,31 | 10,97 | 10,61 |
| 15000 | 13,82 | 13,57 | 13,30 | 13,02 | 12,73 | 12,43 | 12,12 | 11,79 | 11,45 | 11,10 | 10,73 | 10,34 |
| 17500 | 13,81 | 13,53 | 13,24 | 12,94 | 12,63 | 12,31 | 11,97 | 11,62 | 11,25 | 10,87 | 10,47 | 10,05 |
| 20000 | 13,80 | 13,49 | 13,18 | 12,86 | 12,52 | 12,18 | 11,81 | 11,43 | 11,04 | 10,62 | 10,19 | 9,74 |
| 22500 | 13,78 | 13,46 | 13,12 | 12,77 | 12,41 | 12,03 | 11,64 | 11,23 | 10,80 | 10,36 | 9,89 | 9,40 |
| 25000 | 13,76 | 13,41 | 13,05 | 12,67 | 12,28 | 11,88 | 11,45 | 11,01 | 10,55 | 10,07 | 9,57 | 9,04 |
| 27500 | 13,74 | 13,37 | 12,97 | 12,57 | 12,15 | 11,71 | 11,25 | 10,78 | 10,28 | 9,76 | 9,22 | 8,65 |
| 30000 | 13,72 | 13,32 | 12,89 | 12,46 | 12,00 | 11,53 | 11,04 | 10,52 | 9,99 | 9,42 | 8,84 | 8,22 |

**Table 2.1.2. Pressure values along the length of the gas pipeline with a breach at different times ($\ell$ = 15× $10^3$ m).**

Using Table 2.1.2, we construct a graph showing the distribution of pressure every 50 seconds at the beginning, end, and breached sections of the relief gas pipeline (Figure 2.1.4).

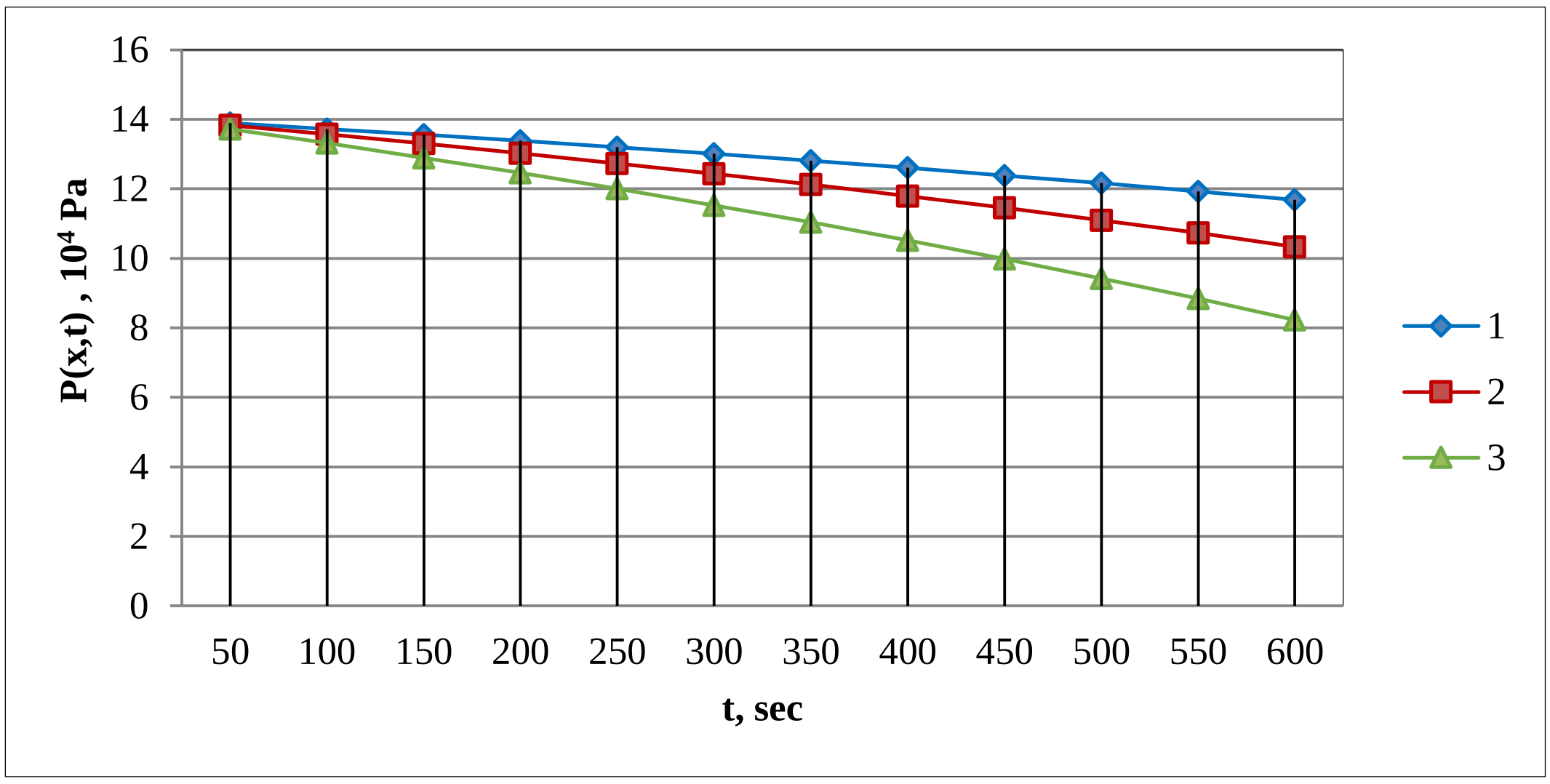

**Figure 2.1.4. Graph of pressure variation over time at different sections of the relief gas pipeline with a breach during non-stationary gas flow ($\ell = 15\cdot10^3$ m, 1 - at the beginning section - P(0,t), 2 - at the breach section - P($\ell$,t), and 3 - at the end section - P(L,t)).**

Analysis of Figure 2.1.4 shows that, despite the change in the location of the breach ($\ell = 15 \cdot 10^3$ m), the gas dynamic process occurring in the pipeline is similar to the process observed at the previous breach location ($\ell = 5\times10^3$ m). However, the pressure drop values at the beginning and end sections of the pipeline decrease relatively.

For instance, at a leakage distance of $\ell=5\times10^3$ meters, the pressure at the initial section drops from $14\times10^{44}$ Pa to $10{,}85\cdot10^4$ Pa (Table 2.1.1). At a leakage distance of $\ell=15\cdot10^3$meters, the initial pressure drops from $14\times10^4$ Pa to $11{,}68\times10^4$ Pa (Table 2.1.2). Similarly, when a leak occurs at a distance of $\ell=5\times10^3$ meters, the pressure at the end section of the gas pipeline drops from $14\times10^4$ Pa to $8{,}22\times10^4$ Pa (Table 2.1.1). If the leakage point is at a distance of $\ell=15\times10^3$ meters, the pressure decreases from $14\times10^4$ Pa to $8{,}22\times10^4$ Pa (Table 2.1.2).

Thus, as the leakage point moves from the initial section to the middle section of the gas pipeline, the pressure drop values decrease by approximately 10% at the initial section and by 30% at the end section of the gas pipeline. According to the unsteady gas flow equation (13) and the gas dynamic process, the pressure drop values at the leakage point do not decrease. This means that the pressure drop values do not change at these points because the flow rate of the leaking gas does not vary depending on the different locations of the leakage point.

Now, to understand how this gas dynamic process occurs in other sections of the uneven gas pipeline, let's use Table 2.1.2 to plot the pressure variation as a function of the pipeline length at every 2.5 km and over time at every 100 seconds (Graph 2.1.5).

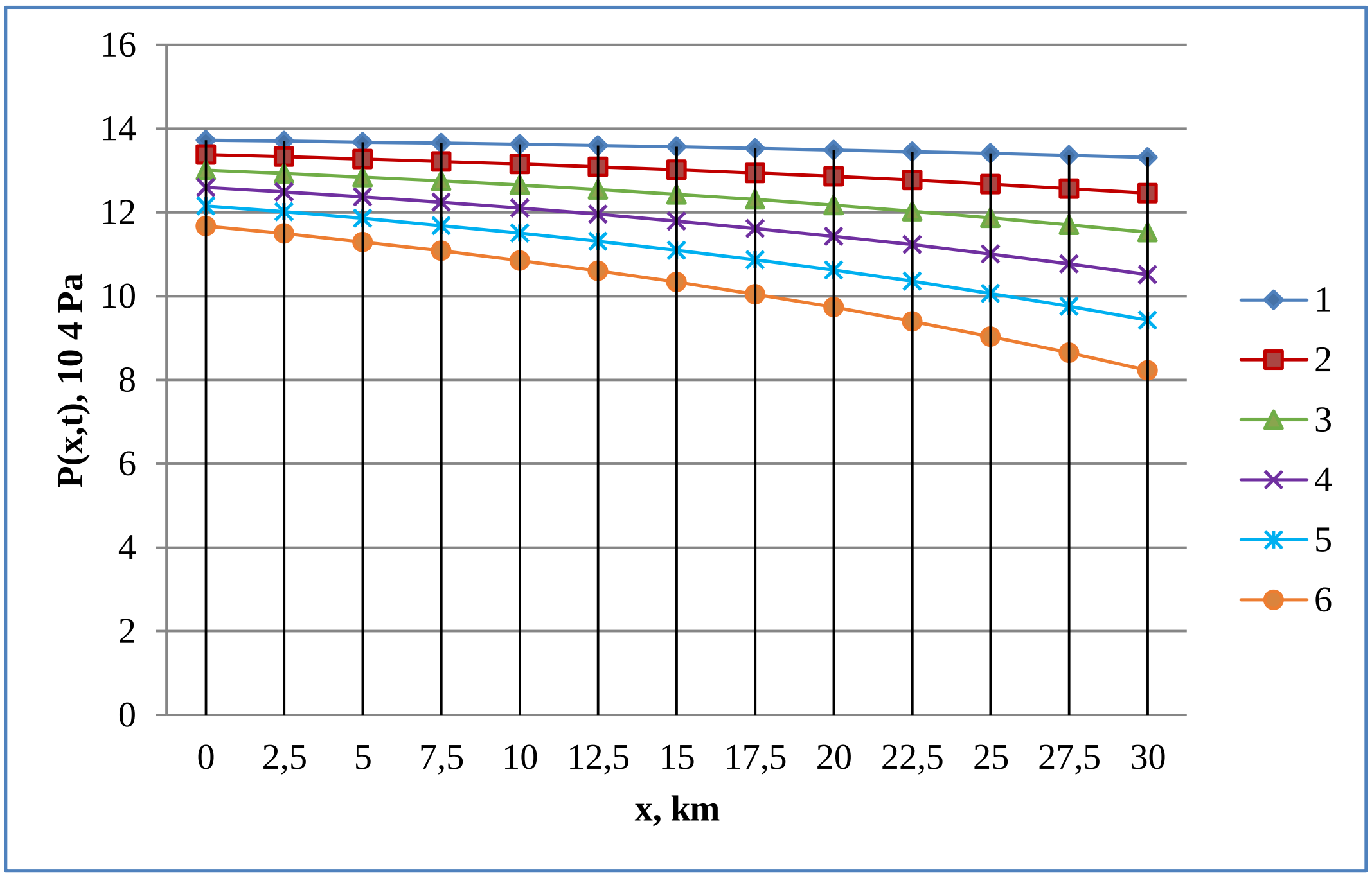

**Graph 2.1.5. The distribution of pressure along the length of the uneven gas pipeline with a compromised seal in an unsteady flow regime over time ($\ell=15\times10^3$ m, 1-t =100 seconds, 2- t=200 seconds, 3- t=300 seconds, 4- t=400 seconds, 5- t=500 seconds, 6- t=600 seconds).**

As a result of the analysis of Graph 2.1.5, it can be noted that the value of the pressure drop [($P_H$ –P(x,t)] increases from the initial section to the end section of the gas pipeline and also increases over time. However, the ratio of the pressure drop value at the initial section of the pipeline to the pressure drop value at any point along the length of the pipeline does not change over time (Graph 2.1.6).

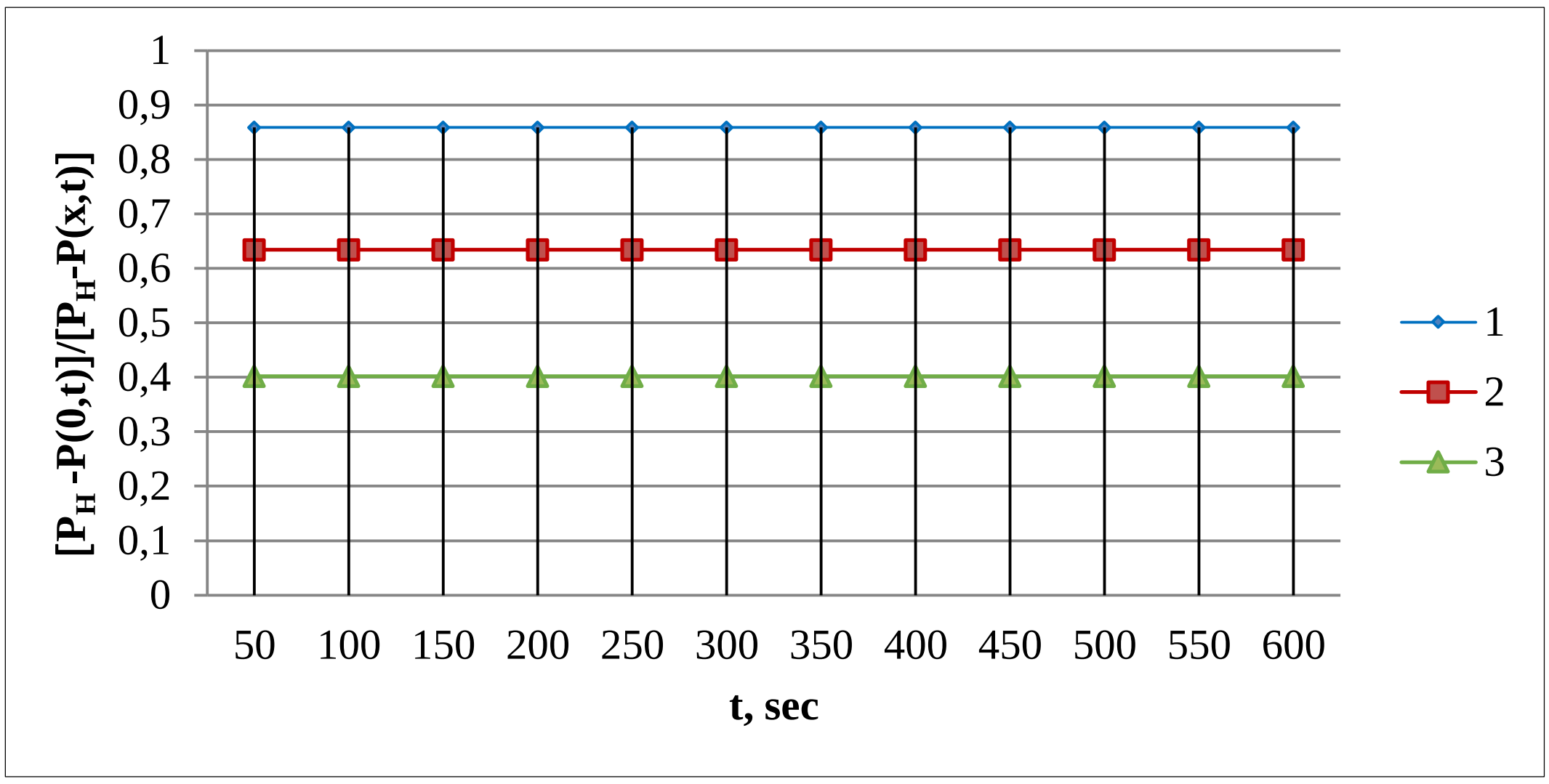

**Graph 2.1.6. The time dependence graph of the ratio of the pressure drop value at the initial section of the uneven gas pipeline with a compromised seal to the pressure drop values at different sections ($\ell$ =15×103 m, 1- initial section at x=5 km, 2- middle section of the gas pipeline at x=15 km, 3- end section at x=30 km).**

From the analysis of Graph 2.1.6, it is clear that, similar to Graph 2.1.3, despite the change in the leakage point of the uneven pipeline, the ratio of the pressure drop value at the initial section of the gas pipeline to the pressure drop values at different sections remains constant over time.

Consequently, by assuming the leakage point in the gas pipeline to be at $\ell = 25 \cdot 10^3$ m, we calculate the pressure values of the gas pipeline at every 50 seconds and every 2.5 km and record them in the following table (Table 2.1.3).

| x, km | t,sec | | | | | | | | | | | |
|---|---|---|---|---|---|---|---|---|---|---|---|---|
| | 50 | 100 | 150 | 200 | 250 | 300 | 350 | 400 | 450 | 500 | 550 | 600 |
| 0 | 13,92 | 13,80 | 13,67 | 13,54 | 13,41 | 13,27 | 13,12 | 12,97 | 12,81 | 12,64 | 12,47 | 12,29 |
| 2500 | 13,91 | 13,78 | 13,65 | 13,51 | 13,36 | 13,21 | 13,05 | 12,89 | 12,71 | 12,53 | 12,35 | 12,15 |
| 5000 | 13,90 | 13,76 | 13,62 | 13,47 | 13,31 | 13,15 | 12,98 | 12,80 | 12,61 | 12,42 | 12,22 | 12,00 |
| 7500 | 13,90 | 13,74 | 13,59 | 13,42 | 13,25 | 13,08 | 12,89 | 12,70 | 12,50 | 12,29 | 12,08 | 11,85 |
| 10000 | 13,89 | 13,72 | 13,56 | 13,38 | 13,20 | 13,01 | 12,81 | 12,60 | 12,39 | 12,16 | 11,92 | 11,68 |
| 12500 | 13,88 | 13,70 | 13,52 | 13,33 | 13,13 | 12,93 | 12,71 | 12,49 | 12,26 | 12,01 | 11,76 | 11,49 |
| 15000 | 13,87 | 13,68 | 13,48 | 13,28 | 13,06 | 12,84 | 12,61 | 12,37 | 12,12 | 11,86 | 11,58 | 11,30 |
| 17500 | 13,86 | 13,65 | 13,44 | 13,22 | 12,99 | 12,75 | 12,50 | 12,24 | 11,97 | 11,69 | 11,39 | 11,08 |
| 20000 | 13,85 | 13,63 | 13,40 | 13,16 | 12,91 | 12,65 | 12,38 | 12,10 | 11,81 | 11,51 | 11,19 | 10,85 |
| 22500 | 13,84 | 13,60 | 13,35 | 13,09 | 12,83 | 12,55 | 12,26 | 11,96 | 11,64 | 11,31 | 10,97 | 10,61 |
| 25000 | 13,82 | 13,57 | 13,30 | 13,02 | 12,73 | 12,43 | 12,12 | 11,79 | 11,45 | 11,10 | 10,73 | 10,34 |
| 27500 | 13,81 | 13,53 | 13,24 | 12,94 | 12,63 | 12,31 | 11,97 | 11,62 | 11,25 | 10,87 | 10,47 | 10,05 |
| 30000 | 13,80 | 13,49 | 13,18 | 12,86 | 12,52 | 12,18 | 11,81 | 11,43 | 11,04 | 10,62 | 10,19 | 9,74 |

Similarly, using Table 2.1.3, we plot the graph of the pressure distribution in the initial, end, and compromised sections of the uneven gas pipeline at every 50 seconds (Graph 2.1.7).

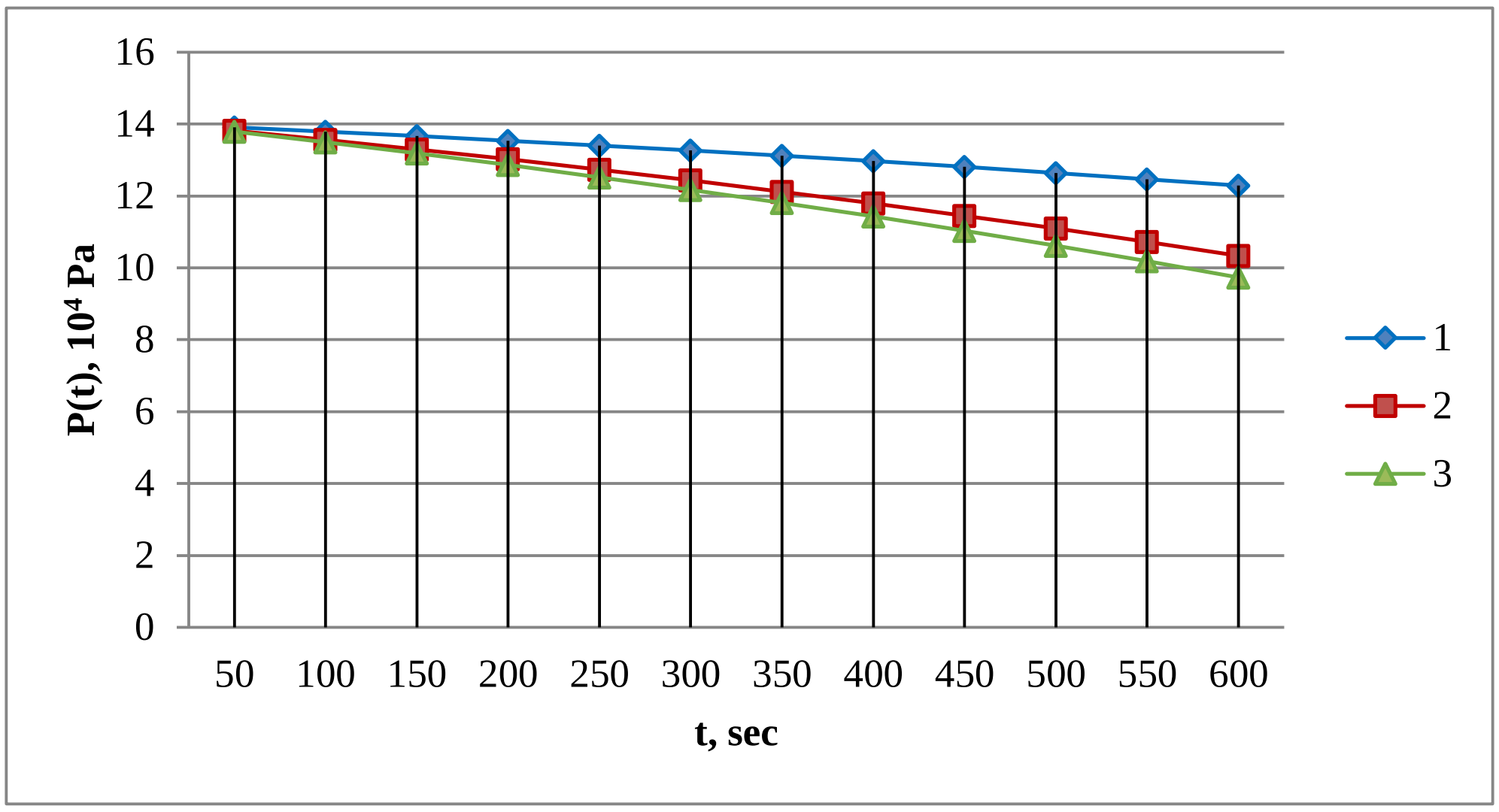

**Graph 2.1.7. The graph of pressure variation over time in the unsteady gas flow of the uneven gas pipeline with a compromised seal ($\ell$ = 5×10$^3$ m, 1- initial section - P(0,t), 2- compromised section - P($\ell$,t), and 3- end section - P(L,t)).**

From the analysis of Graph 2.1.7, it can be noted that, similar to Graph 2.1.1 and Graph 2.1.4, the pressure drop values at the compromised section (leakage point) remain the same. This means that the pressure drop values at these points do not change because the flow rate of the gas leaking into the environment does not vary depending on the different locations of the leakage point. However, the pressure drop values at the initial and end sections of the gas pipeline change depending on the different locations of the leakage point.

Thus, when the location of the leakage point is $\ell$ = 5×10$^3$ m, the pressure drop value [($P_H$ –P(0,t)]=0,15×10$^4$. When the location is $\ell$ =15×10$^3$ m, the pressure drop value [($P_H$ – P(0,t)]=0,11·10$^4$. If the compromised location of the gas pipeline is at a distance of $\ell$ = 25×10$^3$ m, the pressure drop value decreases to 0,08·10$^4$.

Thus, as the location of the compromised seal in the gas pipeline moves from the initial section to the end section, the recorded pressure values at x=0 at t=50 seconds (13,85×10$^4$ Pa, 13,89×·10$^4$ Pa, 13,92·×10$^4$ Pa) increase, and correspondingly, the pressure drop values decrease (Table 2.1.1, Table 2.1.2, and Table 2.1.3).

This pattern corresponds to the process occurring in straight (linear) gas pipelines. That is, in linear gas pipelines, as the location of the compromised seal moves from the initial section to the end section, the pressure drop value at the initial section decreases. Continuing the analysis of Graph 2.1.7, it can also be noted that the recorded pressure values at x=L increase similarly to the initial section(13,62 ×· 10$^4$ Pa, 13,72 × 10$^4$ Pa, 13,8 × 10$^4$ Pa), and correspondingly, the pressure drop values at the end section of the gas pipeline decrease (Table 2.1.1, Table 2.1.2, and Table 2.1.3).

This is the exact opposite of the gas dynamic process occurring in straight (linear) gas pipelines. In linear gas pipelines, as the location of the compromised seal moves closer from the initial section to the end section, the pressure drop value at the end section should increase over time [7]. As mentioned above, during the compromise of the seal in uneven gas pipelines, the dynamics of the unsteady gas flow at the end section are the exact opposite of the dynamics at the end section of linear gas pipelines. The reason for this is that as the location of the gas leak moves closer from the initial section to the end section of the gas pipeline, the height difference resulting from the inclination differences between the initial and end sections increases, and correspondingly, the value of the gravitational force increases. As a result, the mass flow rate of the gas leak increases as it approaches the end section of the gas pipeline, influenced by the terrain (Figure 2.1.1). As the mass flow rate of the reverse flow increases, the gas density at the end section of the gas pipeline decreases, leading to a decrease in pressure in that section, according to the mathematical model of unsteady gas flow.

To obtain information about the pressure drop values along the length of the uneven gas pipeline, let's use Table 2.1.3 to plot the graph of pressure variation along the pipeline length at every 2.5 km and over time at every 100 seconds (Graph 2.1.8).

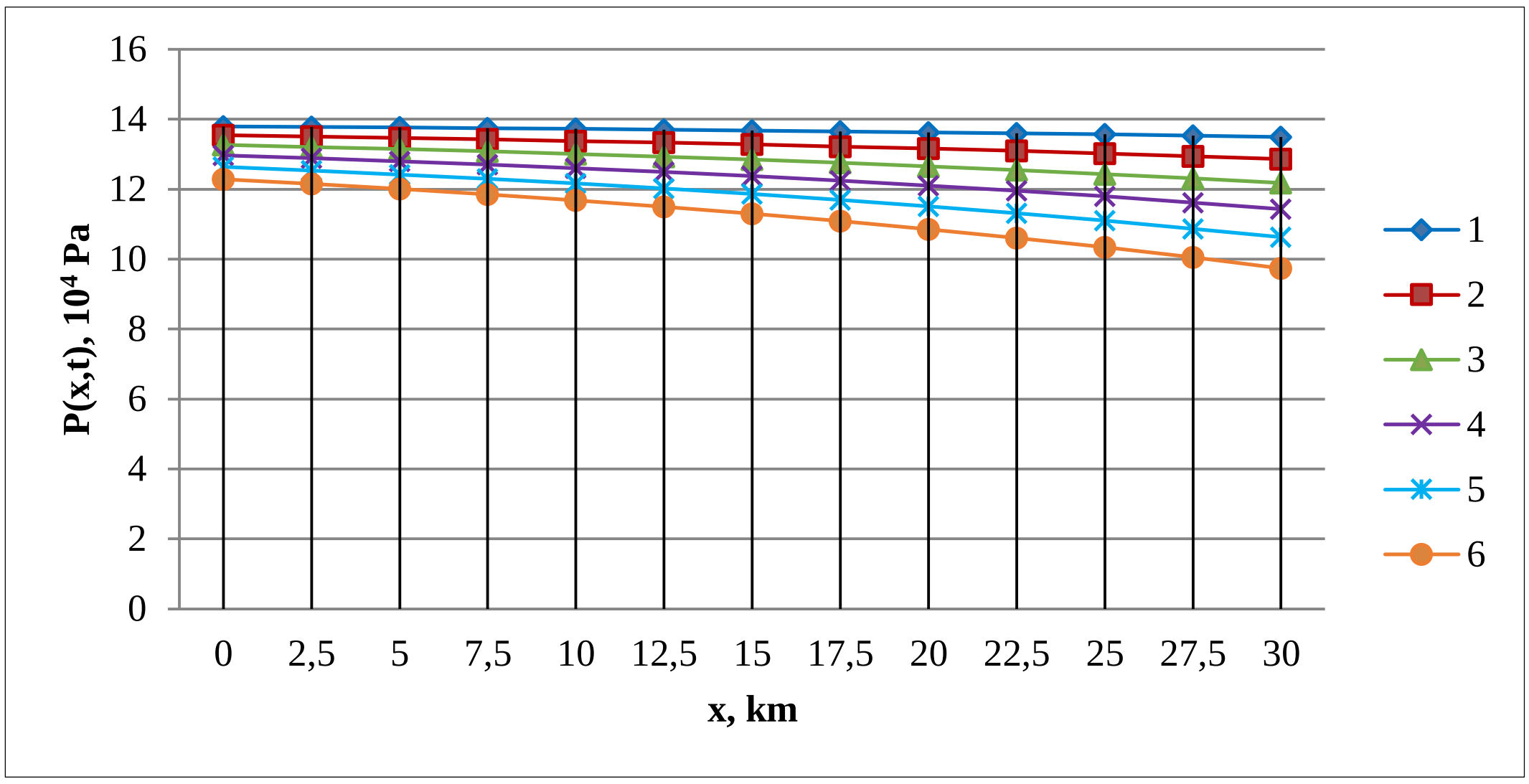

**Graph 2.1.8. The distribution of pressure along the length of the uneven gas pipeline with a compromised seal in an unsteady flow regime over time ( $\ell$ =25·10$^3$ m 1- t=100 seconds, 2- t=200 seconds, 3- t=300 seconds, 4- t=400 seconds, 5- t=500 seconds, 6- t=600 seconds).**

From the analysis of Graph 8, it can be noted that, regardless of the location of the compromised seal, the pressure drop value [($P_H$ –P(x,t)] increases from the initial section to the end section of the gas pipeline and also increases over time. However, the ratio of the pressure drop value at the initial section of the pipeline to the pressure drop value at any point along the length of the pipeline does not change over time (Graph 2.1.9).

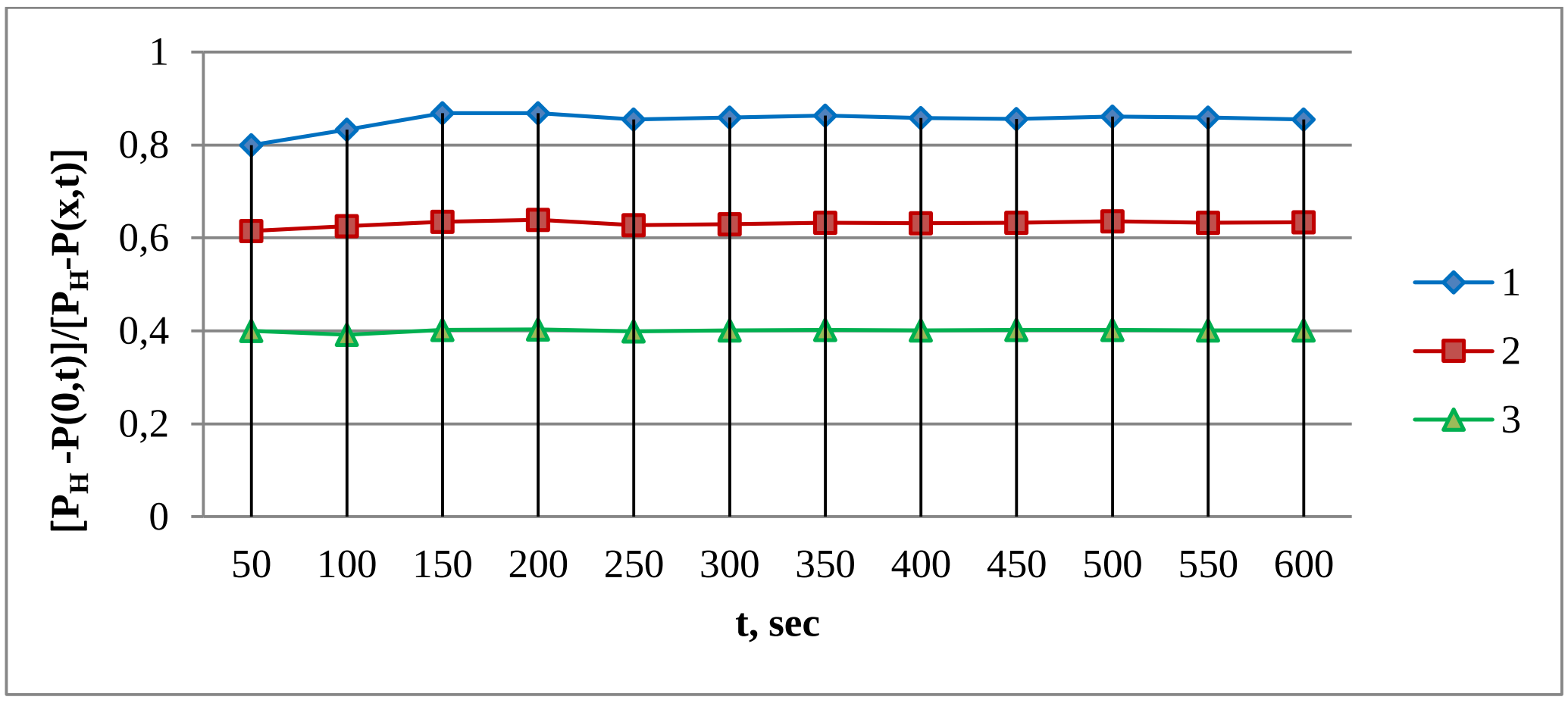

**Graph 2.1.9. The time dependence graph of the ratio of the pressure drop value at the initial section of the uneven gas pipeline with a compromised seal to the pressure drop values at different sections ( $\ell$ =25·10$^3$ m, 1- initial section at x=5 km, 2- middle section of the gas pipeline at x=15 km, 3- end section at x=30 km)**

From the analysis of Graph 2.1.9, it is clear that, similar to Graphs 2.1.3 and 2.1.6, despite the change in the leakage point of the uneven pipeline, the ratio of the pressure drop value at the initial section of the gas pipeline to the pressure drop values at different sections remains constant over time. Additionally, it can be stated that this ratio remains stable over time regardless of the location of the leakage point and that the values of this ratio do not change depending on the location of the leakage point. For example, when the compromised seal location is at $\ell = 5 \times 10^3$ m, the ratio of the pressure drop value at the initial section to the pressure drop value at the middle sections is 0.63 according to Graph 2.1.3. When the compromised seal location is at $\ell = 15 \times \cdot 10^3$ m and $\ell = 30 \cdot 10^3$ m, this ratio remains 0.63, meaning it does not change (Graph 2.1.6 and Graph 2.1.9).

**Justification for the Reconstruction Solution of the Operated Uneven Gas Pipeline.**

The main finding from the conducted analyses is that, due to the influence of the unsteady gas flow regime in the uneven gas pipeline, the pressure drop value at the end section of the pipeline decreases more significantly compared to the initial section, regardless of the location of the gas leak. This is because the pressure decreases along the length of the gas pipeline, and the forces of weight and friction acting along the x-axis move in the opposite direction to the x-axis. This negatively impacts the reliability of the gas pipelines. According to the equation of motion, the pressure drop in the pipeline is directly proportional to the flow rate of the gas stream passing through the cross-section of the pipeline. Therefore, the decrease in pressure drop at the end section of the gas pipeline will result in a decrease in gas flow rate in that section. The reduction in gas supply volume to consumers fed from the end section of the gas pipeline will disrupt the continuous operation of many demanders. Consequently, the decrease in the amount of gas supplied to consumers will reduce the reliability level of the gas pipeline. It is essential to consider this characteristic in the reconstruction of uneven gas pipelines, as one of the key indicators of reconstruction is to increase their reliability.

As a result of the analyses conducted above, we have clarified the reasons for the unusual and opposite gas dynamic process in uneven gas pipelines compared to linear gas pipelines based on the known regularities of the gas dynamic process. Now, let's substantiate this process from a mathematical perspective. Numerous studies on the unsteady gas flow through pipelines confirm that the mathematical modeling of motion is indeed the primary method for studying unsteady processes. It is clear that the first equation of system (1) expresses the equation of motion for the flow of a compressible medium, derived from the law of change of momentum. The first part on the right-hand side of this equation determines the intensity of change in the momentum along the length of the gas pipeline, characterized by the difference between the momentum brought through section x and the momentum carried out through section x+dx (Figure 2.1.1). If we accept the volume of momentum between sections x and x+dx in the cross-section of the pipe (F) as Fdx (x+dx-x=dx), then the rate of change of momentum in volume Fdx with time is characterized as the second term, indicating the unsteady nature of the process.

The second equation in system (1) on the right-hand side shows the difference between the amount of gas passing through sections x+dx and x. The accumulation of gas in the elementary volume Fdx at the same time is given on the left-hand side of the equation. If more gas exits than enters through the cross-section x+dx, the accumulation in the elementary volume Fdx will be negative. This explains the negative sign in equation (1). The conditions appearing in the equation on the right-hand side of system (1) show the forces acting along the xxx-axis on the elementary mass of gas Fdx.

These forces include the gravitational force [ρgSinαFdx], friction [$-\lambda \frac{\rho V^2}{2d} Fdx$], and the pressure gradient [$-\frac{\partial P}{\partial X} Fdx$ ]. The negative sign of these forces is predetermined by the fact that the pressure gradient $\frac{\partial P}{\partial X}$ is negative, which causes the pressure to decrease along the length of the gas pipeline. On the other hand, considering the x-axis and the friction force, the direction of the gravitational force is noted (Figure 2.1.1). The friction force acts in the opposite direction of the x-axis. Thus, the reason for the pressure drop at the end of the uneven gas pipeline always being higher compared to the beginning, regardless of the gas leak location, is due to the combined effects of friction and gravitational forces, as well as the backflow caused by the location and volume of the gas leak. It should be noted that the breach in the pipeline's integrity also creates a backflow at the beginning of the gas pipeline (Graphs 2.1.3, 2.1.6, and 2.1.9). However, theoretical analysis has determined that if the gas leak occurs at a distance close to the beginning of the pipeline, the pressure drop at the beginning is 2.5 times less than at the end (Graphs 2.1.1 and 2.1.2).

One of the options for the reconstruction of gas pipelines is to increase the capacity of the existing pipeline to supply new consumers connected to the current pipeline with the necessary gas. This requirement can be met either by increasing the power or number of gas-pumping equipment at the beginning of the pipeline (such as in a compressor station) or by connecting an optimally parameterized looping system based on the technical and economic analysis of the reconstruction indicators [81].

To increase the capacity of uneven gas pipelines, various factors must be considered when selecting methods for detecting reconstruction solutions for the looping system, such as pipeline parameters, linear profile, characteristics of the transported product, and reliability. Considering that it is appropriate to apply the looping system option to increase the capacity of the existing uneven gas pipeline during the reconstruction phase, it is recommended to install the looping at the beginning of the pipeline to ensure continuous gas supply to additional new consumers fed from the gas pipeline in both stationary and non-stationary modes. This solution has been substantiated by the above-mentioned theoretical analyses.

**Conclusion**

1. A solution has been obtained for the system of differential equations describing unsteady gas flow along a pipeline with a constant cross-section, considering the pipeline's terrain. The inclusion of boundary conditions set to zero in the solution allows for a compact notation of the analytical model, yet it does not limit the model's applicability to cases of sharp changes in gas flow or pressure. Thus, a methodology for calculating unsteady gas flow under constant terrain conditions in pipelines has been developed.

2. The impact of the gravitational force resulting from the slope of the pipeline profile on the flow characteristics of the linear section of the pipe has been studied. In the analysis of the unsteady flow, the influence of the terrain on the gas dynamic process for gas transport with the same elevation difference has been investigated. The results of pressure calculations for various locations of gas leaks have been analyzed. The comparison of calculation examples and their theoretical analysis allows for the assessment of the accuracy of the proposed mathematical expressions for calculating the pressure in uneven gas pipelines depending on both length and time.
3. When selecting methods for identifying reconstruction solutions for an uneven gas pipeline system, various factors such as pipeline parameters, linear profile, characteristics of the transported product, and reliability have been considered. Based on this, conclusions have been drawn regarding the preferred location for installing the looping in the pipeline for its reconstruction, and these conclusions have been substantiated through theoretical investigation.

## 2.2. Technological Foundations of Reconstruction Solutions for Looping Gas Pipeline Systems

**Introduction.** With the increasing population in the country, new buildings and structures of various purposes are being constructed in the prospective areas of settlements. This leads to a growing energy balance and an increase in the demand for gas. Consequently, existing gas networks are typically used to supply new buildings with gas fuel. If the configuration of the existing gas pipeline is a ring network, the consumption (load) of the existing ring network is increased to accommodate the volume of new consumers. Thus, there is a need for the reconstruction of the existing ring gas pipeline [1,2].

Reconstruction of the ring gas pipeline presents several modeling problems for managing complex interconnected pipeline systems under both steady-state and unsteady-state conditions. Modeling of pipeline systems for gas transportation must be carried out by altering the structure of gas flows and modifying the consumption and supply volumes of gas through the use of existing ring pipeline systems. The difficulty in modeling unsteady flow regimes for such systems is related to the complexity and variability of the dynamic state of gas flow in pipelines, as well as the need to account for many different factors.

Our goal is to develop and implement solutions for the reconstruction of the ring distribution network in accordance with the requirements for efficiency and reliability indicators. To achieve this, it is necessary to recommend the hydraulic junction node as an energy source for new consumers and provide a visual schematic of this junction. The mathematical model and solution method for the gas-dynamic state of the ring gas pipeline network, which includes a flow rate-based zone along the road under various flow regimes, are presented according to the requirements of this schematic.

When studying the gas-dynamic issues of ring networks, it should be noted that in both distribution ring networks diverging from the starting point, gas flows in opposite directions and intersects at any point, creating equilibrium at the intersection. That is, at the starting point of the gas pipeline, the pressures in the hydraulic junction sections of both distribution branches equalize. Thus, the flow streams of each distribution branch intersect at any point $x=x_1$ (see Scheme 2.2.1). This may cause problems in the gas supply regime for consumers. The solution to this problem is based on two conditions:

The first condition is based on the equilibrium of the aforementioned hydraulic junction node. If the gas flow into the hydraulic junction is considered positive and the gas flow exiting the junction and within the junction itself is considered negative, then the algebraic sum of the calculated gas flow rates at the junction should be zero.

The second condition is that the sum of pressure losses in one semi-ring should be equal to the sum of pressure losses in the other semi-ring. If we assume the clockwise direction of gas pressure as positive and the opposite direction as negative, then the difference in the algebraic sum of pressure losses in both rings should be zero (see Scheme 2.2.1). In practice, this condition is not always possible to meet, and it is assumed that the discrepancy in pressure losses between semi-rings should not exceed 10%.

Numerous local and international studies have addressed these issues [21,77,80,88,89]. These works present methods of hydraulic schematic theory, mathematical modeling, numerical methods for solving problems related to the gas-dynamic state of pipelines, optimization of pipeline networks, and algorithms for calculating the gas-dynamic parameters of ring gas pipelines. However, the reviewed materials often do not extensively cover the analysis of unsteady gas flow parameters throughout the entire length of the ring pipeline network.

The main distinctive feature of our research is the analysis and examination of functions that enable the precise solution of the unsteady flow regime problem along the entire length of the ring gas pipeline, corresponding to the number and volume of road-based gas flows connected to the ring pipeline. Additionally, it is recommended that the hydraulic junction node ($x=x_1$, see Scheme 2.2.1) be considered as an energy supply source for new consumers in the reconstruction solution of existing ring gas pipelines.

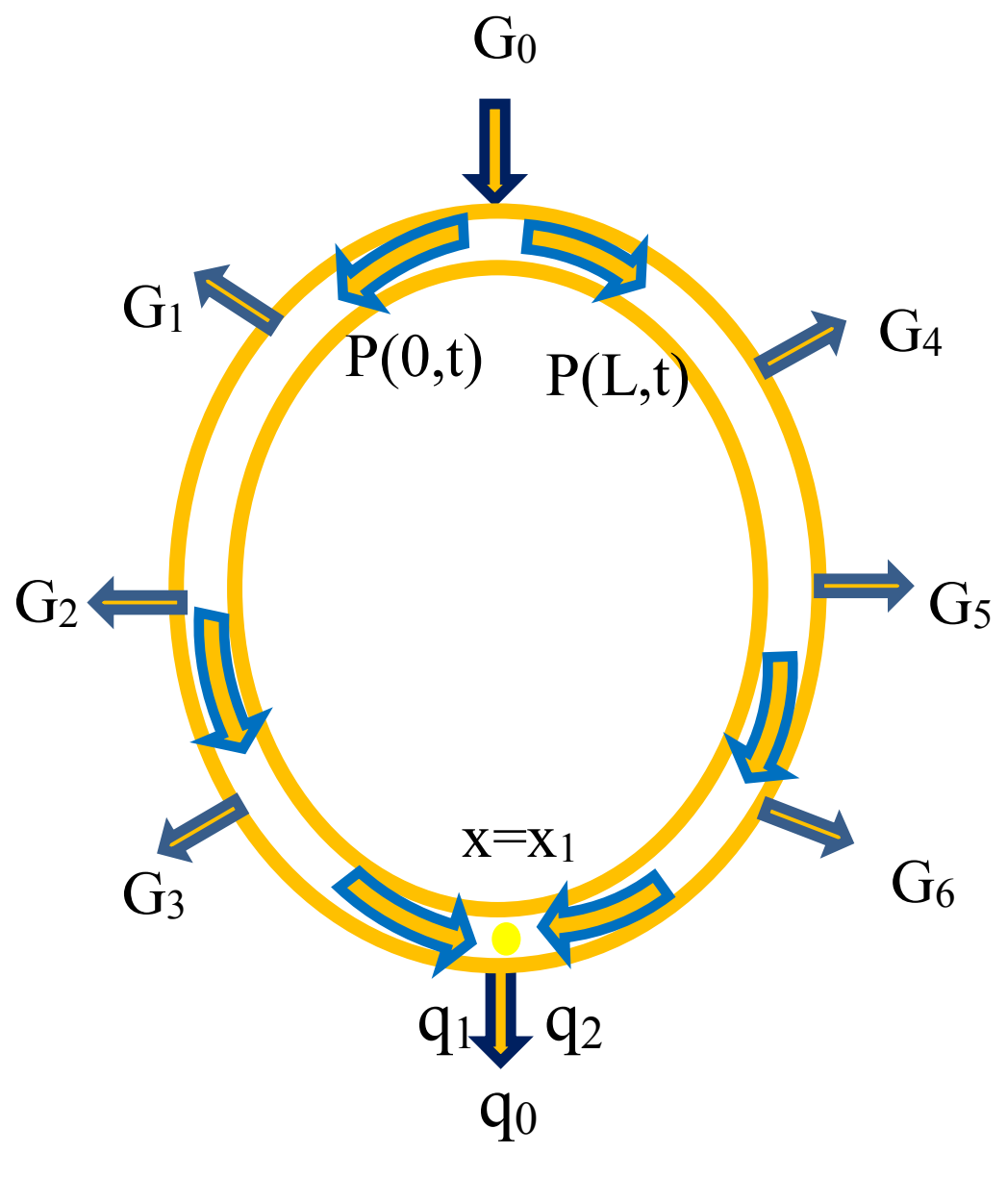

**Scheme 2.2.1. Visual Diagram of the Ring Gas Supply System**

As shown in Scheme 2.2.1, a characteristic feature of ring gas supply systems is that the starting and ending points of the network converge at a single point. In a steady-state regime, the starting ($P_b$) and ending ($P_s$) pressures are equal ($P_b = P_s = P_1$). Since in ring networks the direction of gas flow is opposite from the start and end points, these flow streams meet at some point $x=x_1$. This meeting point is referred to as the hydraulic junction of the ring network, and in some cases, it is called a "dead" point.

From this perspective, during the operation of ring gas supply systems, the pressures of the flow streams entering the hydraulic junction must be equal, and the total flow rates of the incoming streams from both branches must equal the road-based consumption (q0) exiting the junction ($q_1 + q_2 = q_0$; see Scheme 2.2. 1). Otherwise, the utilization of gas consumption by the consumers supplied by the network may be disrupted due to varying pressure flows, which is unacceptable.

The main advantage of ring gas supply systems is that even if a leak occurs at any point in the network or if a section needs to be removed, gas supply to consumers continues uninterrupted despite the deviation from normal operating conditions.

### *Mathematical Model Analysis for Unsteady Gas Flow in Ring Pipelines*

When addressing solutions to the considered problem, it is advisable to obtain solvable mathematical expressions for studying unsteady processes and managing the occurring physical processes. Therefore, using unsteady-state models is appropriate. These models enable the resolution of various issues, including those related to the reconstruction of ring networks. Consequently, models and algorithms for calculating hydraulic schematics of ring gas pipelines need to be refined.

For unsteady gas flow, I.A. Charnynin's differential equations system relates gas flow parameters, such as pressure P(x,t) and gas mass flow G(x,t), as functions of time t and distance x. Based on Diagram 1, the system of specific derivative non-linear equations for the problem under consideration can be expressed as follows [1,2].

$$\begin{cases} -\dfrac{\partial P}{\partial X} = \lambda \dfrac{\rho V^2}{2d} \\ -\dfrac{1}{c^2}\dfrac{\partial P}{\partial t} = \dfrac{\partial G}{\partial X} + \sum_{i=1}^{n} Gi(t)\delta(x - x_i) \end{cases} \tag{1}$$

Here, $G=\rho v$; $P=\rho \cdot c^2$ ; $\rho = \dfrac{P}{zRT}$

T – Average temperature of the gas.
z=z (p,T) – Compressibility factor.
c – Speed of sound in the gas for an isothermal process, m/sec.
P – Pressure in the cross-sections of the gas pipeline, Pa.
x – Coordinate along the axis of the gas pipeline, m.
t – Time coordinate, sec.
λ – Hydraulic resistance coefficient.
ρ – Average density of the gas, $kg/m^3$.
v – Average velocity of the gas flow, m/sec.
d – Diameter of the gas pipeline, m.
G – Mass flow rate of the gas.
$Gi(t)$ – Road-based mass flow rates from various points of the ring gas pipeline. $\dfrac{Pa \cdot \sec}{m}$
$\delta(x - x_i)$ – Dirac delta function.

It is evident that the nonlinear system of equations (1) cannot be solved directly, as such equations are not integrable. In engineering calculations, approximate methods are employed. In other words, the equations are linearized. Linearization of equations, as proposed by I.A. Charnynin, is more practical for solving these problems. The expression characterizing this linearization is given as follows [78].

$$2a = \lambda \frac{V}{2d} \qquad (2)$$

Here, 2a is the linearization coefficient according to I.A. Charnynin, then the system of equations (1) and (2) will be transformed into the following form:

$$\begin{cases} -\frac{\partial P}{\partial X} = 2aG \\ -\frac{1}{c^2}\frac{\partial P}{\partial t} = \frac{\partial G}{\partial X} + \sum_{i=1}^{n} Gi(t)\delta(x - x_i) \end{cases} \qquad (3)$$

When differentiating the first expression of the system of equations (3) with respect to x and solving for $\frac{\partial G}{\partial X}$ from the second expression, the equation will be in the following form:

$$\begin{cases} -\frac{\partial^2 P}{\partial X^2} = 2a\frac{\partial G}{\partial X} \\ \frac{\partial G}{\partial X} = -\frac{1}{c^2}\frac{\partial P}{\partial t} - \sum_{i=1}^{n} Gi(t)\delta(x - x_i) \end{cases} \qquad (4)$$

Then, by substituting the expression for $\frac{\partial G}{\partial X}$ from the first equation into the expression (4), we obtain the heat transfer equation for the mathematical solution of the problem.

$$\frac{\partial^2 P}{\partial X^2} = \frac{2a}{c^2}\frac{\partial P}{\partial t} + 2a\sum_{i=1}^{n} Gi(t)\delta(x - x_i) \qquad (5)$$

Thus, for the uniqueness of the solution to equation (5) that reflects the studied process and for the complete determination of the given process, the initial and boundary conditions must be correctly specified.

**Initial Conditions:** When inertia forces are not considered, the initial distribution of the sought function (pressure) along the pipeline at t=0 is provided.

**Boundary Conditions:** The boundary conditions ensure that the pressure and flow rate at the endpoints of the pipeline change according to a predetermined law during the process under consideration.

A correct specification of the boundary conditions completes the mathematical model of the process and allows for a comprehensive and accurate study of the physical phenomena occurring within the system. The sum of the intensities of the concentrated road-based flow rates at the pipeline junctions is considered equal to the initial flow rate in steady-state regimes. The diameter of the investigated gas pipeline is uniform throughout. The supply intensity and the total mass flow rate of the gas is $G_0$, and the pressure at the starting point is $P_1$.Then,

$t=0$; $P(x,0)=P_1-2aG_0x$: $G_0x=\sum G_ix_i$, $i=1,2,3,...,n$

For the main sections of the gas pipeline, we assume the following boundary conditions for the process under consideration (continuity of flow at the initial and final points):

**Initial Point:** The flow rate and pressure at the beginning of the pipeline are set according to the given conditions.

**Final Point:** The flow rate and pressure at the end of the pipeline are adjusted to ensure continuity and adherence to the predetermined law during the process.

$$x=0,L \quad \text{point} \begin{cases} P(0,t)=P(L,t) \\ \dfrac{\partial P(0,t)}{\partial X}=\dfrac{\partial P(L,t)}{\partial X} \end{cases}$$

It is clear that by applying the Laplace transform [U(x,s)= $\int_0^\infty U(x,t)e^{-st}dt$ ], equation (5) is transformed into second-order ordinary differential equations. The general solutions for these equations will be as follows:

$$P(x,s)=\sqrt{\frac{A}{s}}c^2\sum_{n=1}^{n}G_i(s)\int_0^x sh\sqrt{As}(x-y)\delta(y-x_i)dy- \\ -\sqrt{\frac{A}{s}}\int_0^x sh\sqrt{ks}(x-y)(P_1-2aG_0y)dy+C_1sh\sqrt{As}x+C_2ch\sqrt{As}x \tag{6}$$

Here $A=\dfrac{2a}{c^2}$

By applying the Laplace transform to the initial and boundary conditions, we determine the constants $C_1$ and $C_2$ in equation (8).

$$C_1=\sqrt{\frac{A}{s}}c^2\sum_{n=1}^{n}G_i(s)\frac{sh\sqrt{As}\left(x_i-\frac{L}{2}\right)}{2sh\sqrt{As}\frac{L}{2}}+\frac{2aG_0L}{2S}\frac{ch\sqrt{As}\frac{L}{2}}{2sh\sqrt{As}\frac{L}{2}}-\frac{2aG_0}{s\sqrt{As}} \tag{7}$$

$$C_2=-\sqrt{\frac{A}{s}}c^2\sum_{n=1}^{n}G_i(s)\frac{ch\sqrt{As}\left(\frac{L}{2}-x_i\right)}{2sh\sqrt{As}\frac{L}{2}}+\frac{P_1}{S}-\frac{2aG_0L}{2S} \tag{8}$$

By considering the values of expressions (7) and (8) in equation (6), we obtain the transformed equation for the dynamic state of the process, or in other words, for the ring-shaped gas pipeline. The resulting equation will represent the transformed distribution of pressure along the ring-shaped gas pipeline in a non-stationary flow regime.

$$P(x,S)=\frac{P_1-2aG_0x}{S}-\frac{2aG_0L}{S}\frac{sh\sqrt{As}\left(\frac{L}{2}-x\right)}{2sh\sqrt{As}\frac{L}{2}}-\sqrt{\frac{A}{s}}c^2\sum_{n=1}^{n}G_i(s)\frac{ch\sqrt{As}\left(\frac{L}{2}-x+x_i\right)}{2sh\sqrt{As}\frac{L}{2}} \tag{9}$$

Based on the solution of the problem, since determining the mathematical expression for the physical process of the non-stationary flow of gas in the pipeline is more suitable, we use the inverse Laplace transform for each boundary to find the original solution of equation (9). This process results in the following original equation:

$$P(x,t)=P_1-2aG_0x-2aG_0\left(\left(\frac{L}{2}-x\right)+L\sum_{n=1}^{\infty}(-1)^n\frac{Sin(L-x)\frac{\pi n}{L}}{\pi n}\int_0^t e^{-n^2\alpha(t-\tau)}d\tau\right)-$$

$$-\frac{c^2}{L}\int_0^t\sum_{i=1}^{n}G_i(\tau)\left[1+2\sum_{n=1}^{\infty}(-1)^n Cos[L-2(x-x_i)]\frac{\pi n}{L}e^{-n^2\alpha(t-\tau)}\right]d\tau \qquad (10)$$

Here, $G_0x=\sum_{i=1}^{n}G_ix_i$, $\alpha=\frac{4\pi^2c^2}{2aL^2}$

If we assume that the road-based gas consumption used by consumers is constant ($G_i$(t) $=G_i$ =const) then through several mathematical operations, equation (10) will be simplified to the following form:

$$P(x,t)=P_1-2aG_0\frac{L}{2}+2aG_0L\sum_{n=1}^{\infty}Sin\frac{\pi nx}{L}\left(\frac{1-e^{-n^2\alpha t}}{\pi n^3}\right)$$

$$-\frac{c^2t}{L}\sum_{i=1}^{n}G_i-\frac{2c^2}{L}\sum_{i=1}^{n}G_i\sum_{n=1}^{\infty}Cos\frac{2\pi n(x-x_i)}{L}\left(\frac{1-e^{-n^2\alpha t}}{\alpha n^2}\right) \qquad (11)$$

To determine the law of variation of pressure in a ring-shaped network, we accept the following data:

$P_1=14\times10^4$ Pa : $G_0$= 10 Pa ×sec/m : $2a=0,1\frac{1}{\text{sec}}$ ; $c=383,3\frac{m}{\text{sec}}$;d= 0,7 m: L = $3\times10^4$ m

First, assuming that the number of road-based flows connected to the ring-shaped gas pipeline is three, with the following parameters:

($G_1$=3 Pa ·sec/m ; $G_2$=4 Pa·sec/m; $G_3$=3 Pa ·sec/m)

And that the junction points are located at distances:

$x_1=0,3\cdot10^4$ m; $x_2=1,5\cdot10^4$ m: $x_3=2,7\cdot10^4$ m

We calculate the pressure values of the gas pipeline as a function of time and distance, and record them in the following table (Table 2.2.1).

| x,m | P(x,t), Pa | | |
|---|---|---|---|
| | t=50 sec | t=300 sec | t=900 sec |
| 0 | 122717,5 | 110478,6 | 81094,8 |
| 1000 | 146064,4 | 134445,6 | 105061,8 |
| 3000 | 186412,1 | 176006,9 | 146623,1 |
| 6000 | 229710,7 | 220960,7 | 191576,9 |
| 9000 | 254914,2 | 247478,2 | 218094,4 |
| 12000 | 264319,1 | 257726,5 | 228342,7 |
| 15000 | 260046,5 | 253743,6 | 224359,8 |
| 17000 | 251713,7 | 245279,8 | 215896 |
| 20000 | 231199,6 | 224097,4 | 194713,6 |
| 23000 | 203327 | 195052,9 | 165669,2 |
| 27000 | 158388,6 | 147974,3 | 118590,6 |
| 30000 | 122357,1 | 110108,7 | 80724,91 |

**Table 2.2.1. Values of the Pressure Variation Along the Length of the Ring-shaped Network as a Function of Time ($G_1+G_2+G_3=G_0$)**

Using Table 2.2.1, we plot the graph of pressure distribution along the length of the ring-shaped gas pipeline as a function of time (Figure 2.2.1).

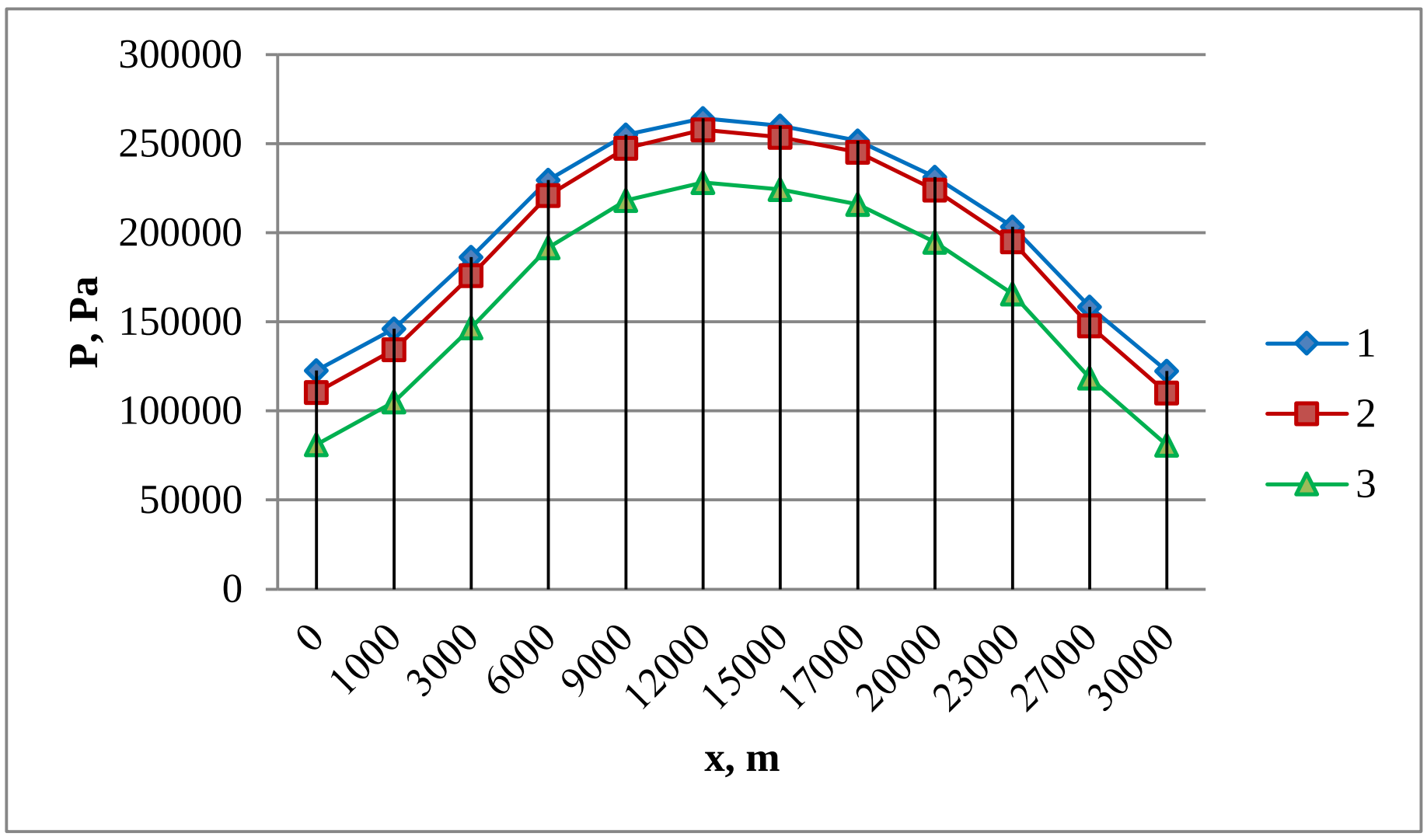

**Figure 2.2.1. Pressure Distribution Along the Length of the Road-based Gas Pipeline as a Function of Time (1-t=50 sec, 2-t=300 sec, 3- t=900 sec).**

The analysis of Figure 1 reveals the following:

**Pressure Increase and Maximum:** The pressure increases from the initial part of the pipeline towards the middle (counterclockwise direction) and reaches its maximum value near the hydraulic junction in the middle section of the pipeline.

**Pressure Decrease:** From the middle section towards the end of the pipeline (clockwise direction), the pressure decreases and eventually stabilizes to match the initial pressure at the starting point, i.e., [P(0,t)=P(L,t)].

**Reason for Pressure Increase:** The increase in pressure in the middle sections of the ring-shaped pipeline is due to the incremental extraction of the gas mass flow ($G_0$) from various points along the pipeline ($G_i x_i$ ,i=1,2,3,...,n)

**Pressure Decrease Over Time:** At the start, the pressure is $14\cdot10^4$ Pa, and after 50 seconds, it decreases to $12.3\times10^4$ Pa. Over time, this value gradually reduces, reaching $11.05\times10^4$ Pa at 300 seconds and $8.1\times10^4$ Pa at 900 seconds (refer to Table 2.2.1). This decrease continues until a new stationary state is achieved.

**Cause of Pressure Decrease:** The decrease in pressure at the starting and ending points is due to the counterflow of the gas and hydraulic losses. According to gas dynamic processes, hydraulic losses are directly proportional to the gas mass flow.

To provide further clarity on the results of the analysis, we assume that the total sum of the flow rates is kept constant while increasing the number of flow points to four, with the following values:

($G_1$=2 Pa ·sec/m; $G_2$=3 Pa ·sec/m; $G_3$=2 Pa·sec/m; $G_4$=3 Pa ·sec/m; $x_1=0,3\times10^4$ m; $x_2=1,5\times10^4$ m: $x_3=2,7\times10^4$ m; $x_4=2,0\times10^4$ m)

By calculating the pressure values of the gas pipeline depending on time and length, we record them in the following table (Table 2.2.2).

| x,m | P(x,t), Pa | | |
|---|---|---|---|
| | t=50 s | t=300 s | t=900 s |
| 0 | 122717,5 | 110478,6 | 81094,8 |
| 1000 | 146064,4 | 134445,6 | 105061,8 |
| 3000 | 186412,1 | 176006,9 | 146623,1 |
| 6000 | 229710,7 | 220960,7 | 191576,9 |
| 9000 | 254914,2 | 247478,2 | 218094,4 |
| 12000 | 264319,1 | 257726,5 | 228342,7 |
| 15000 | 260046,5 | 253743,6 | 224359,8 |
| 17000 | 251713,7 | 245279,8 | 215896 |
| 20000 | 231199,6 | 224097,4 | 194713,6 |
| 23000 | 203327 | 195052,9 | 165669,2 |
| 27000 | 158388,6 | 147974,3 | 118590,6 |
| 30000 | 122357,1 | 110108,7 | 80724,91 |

**Table 2.2.2. Values of the pressure variation along the length of the gas pipeline over time in a ring network with side ( $G_1$ +$G_2$+$G_3$+$G_4$=$G_0$ ).**

Using Table 2.2.2, we construct the graph of the pressure distribution along the length of the ring gas pipeline over time (Figure 2.2.2).

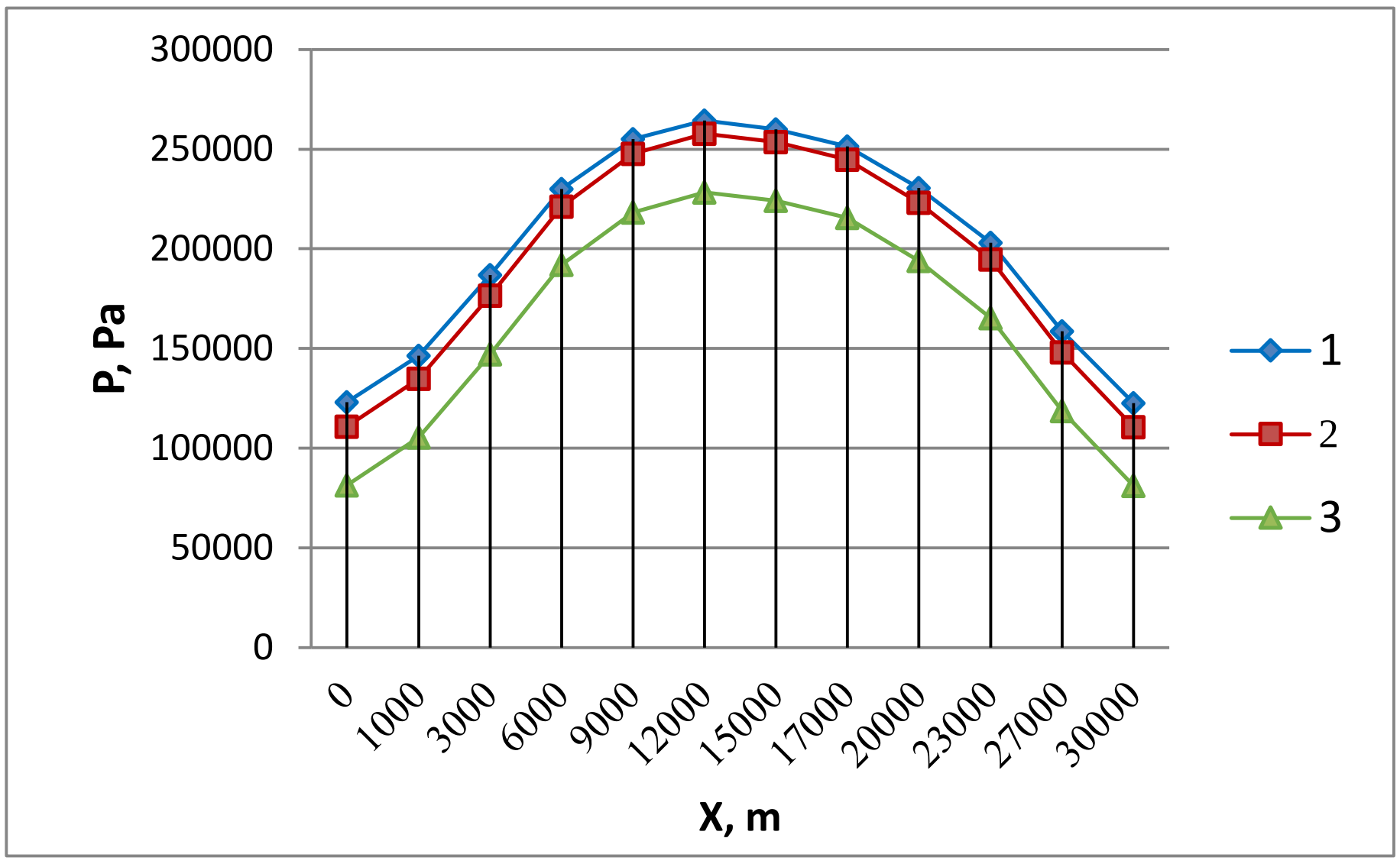

**Figure 2.2.2. Distribution of pressure along the length of the ring gas pipeline with side flows over time. (1-t=50 sec, 2-t=300 sec, 3- t=900 sec).**

Analysis of Figure 2.2.2 reveals that in ring pipelines, increasing the number of side flow junctions while keeping the total side flow constant (as noted in Figure 2.2.1) does not affect the pressure distribution over time or along the length of the pipeline. This is because the gasodynamic calculations for the ring are performed according to Kirchhoff's first law at all network junctions. In this case, the distribution of flows can vary depending on the diameters of the pipeline sections. The research was conducted under the condition that the diameter of the existing ring pipeline remained constant. Thus, the loss of pressure balance in closed loops is not possible. The restoration of pressure balance in all closed loops depends not on the number of network flows but on the mass flow rate. To verify this, we compare the pressure values of the gas flow noted in Table 2.2.1 with those in Table 2.2.2 by creating the following table (Table 2.2.3).

Table 2.2.3. Comparison of Pressure Values Over Coordinates and Time in Ring Pipelines with Constant Total Side Flows.

| x, m | P(x,t), $10^4$ Pa | | | | | |
|---|---|---|---|---|---|---|
| | t=50 sec | | t=300 sec | | t=900 sec | |
| | $\sum_{i=1}^{3} Gi = 10$, Pa ×sec/m | $\sum_{i=1}^{4} Gi = 10$, Pa×sec/m | $\sum_{i=1}^{3} Gi = 10$, Pa×sec/m | $\sum_{i=1}^{4} Gi = 10$, Pa ×sec/m | $\sum_{i=1}^{3} Gi = 10$, Pa ×sec/m | $\sum_{i=1}^{4} Gi = 10$, Pa×sec/m |
| 0 | 12,2 | 12,3 | 11,0 | 11,0 | 8,1 | 8,1 |
| 1000 | 14,6 | 14,6 | 13,4 | 13,5 | 10,5 | 10,5 |
| 3000 | 18,6 | 18,6 | 17,6 | 17,6 | 14,7 | 14,7 |
| 6000 | 23,0 | 23,0 | 22,1 | 22,1 | 19,2 | 19,2 |
| 9000 | 25,5 | 25,5 | 24,7 | 24,8 | 21,8 | 21,8 |
| 12000 | 26,4 | 26,4 | 25,8 | 25,8 | 22,8 | 22,8 |
| 15000 | 26,0 | 26,0 | 25,4 | 25,4 | 22,4 | 22,4 |
| 17000 | 25,2 | 25,1 | 24,5 | 24,5 | 21,5 | 21,5 |
| 20000 | 23,1 | 23,0 | 22,4 | 22,3 | 19,5 | 19,4 |
| 23000 | 20,3 | 20,3 | 19,5 | 19,5 | 16,6 | 16,5 |
| 27000 | 15,8 | 15,9 | 14,8 | 14,8 | 11,9 | 11,9 |
| 30000 | 12,2 | 12,3 | 11,0 | 11,0 | 8,1 | 8,1 |

**Table 2.2.3 clearly shows that the pressure values for a ring pipeline with three mass flows are identical to those for a ring pipeline with four mass flows at the same time and coordinates (length) when the total mass flow is the same ($\sum_{i=1}^{3} Gi = \sum_{i=1}^{4} Gi$).**

Now, we will consider the case where the number of junctions for the side branches in ring gas pipelines ($x_1$=0,3×$10^4$ m; $x_2$=1,5×$10^4$ m: $x_3$=2,7×$10^4$ m; $x_4$=2,0×$10^4$ m) is kept constant at four, but their total flows are increased by 50%. Thus, we assume $G_1$=3 Pa ·sec/m; $G_2$=5 Pa ·s/m; $G_3$=3 Pa ·sec/m; $G_4$=4 Pa ·sec/m. We will calculate the pressure values of the gas pipeline over time and length and record these values in the table below for comparison with the values in Table 2.2.3 (Table 2.2.4).

| x, m | P(x,t), $10^4$ Pa | | | | | |
|---|---|---|---|---|---|---|
| | t=50 sec | | t=300 sec | | t=900 sec | |
| | $\sum_{i=1}^{4} Gi$, Pa ·sec/m | $1{,}5\sum_{i=1}^{4} Gi$, Pa ·sec/m | $\sum_{i=1}^{4} Gi$, Pa ·sec/m | $1{,}5\sum_{i=1}^{4} Gi$, Pa ·sec/m | $\sum_{i=1}^{4} Gi$, Pa ·sec/m | $1{,}5\sum_{i=1}^{4} Gi$, Pa ·sec/m |
| 0 | 12,3 | 11,4 | 11,0 | 9,6 | 8,1 | 5,2 |
| 1000 | 14,6 | 15,0 | 13,5 | 13,2 | 10,5 | 8,8 |
| 3000 | 18,6 | 21,0 | 17,6 | 19,5 | 14,7 | 15,0 |
| 6000 | 23,0 | 27,5 | 22,1 | 26,2 | 19,2 | 21,8 |
| 9000 | 25,5 | 31,3 | 24,8 | 30,1 | 21,8 | 25,7 |
| 12000 | 26,4 | 32,6 | 25,8 | 31,7 | 22,8 | 27,2 |
| 15000 | 26,0 | 32,0 | 25,4 | 31,0 | 22,4 | 26,6 |
| 17000 | 25,1 | 30,7 | 24,5 | 29,7 | 21,5 | 25,3 |
| 20000 | 23,0 | 27,6 | 22,3 | 26,5 | 19,4 | 22,1 |
| 23000 | 20,3 | 23,4 | 19,5 | 22,2 | 16,5 | 17,8 |
| 27000 | 15,9 | 16,8 | 14,8 | 15,2 | 11,9 | 10,82 |
| 30000 | 12,3 | 11,4 | 11,0 | 9,6 | 8,1 | 5,2 |

**Table 2.2.4. Comparison of Pressure Values over Time and Length in Ring Pipelines with Different Total Flows of Side Branches.**

Using the total side flow values $1{,}5\sum_{i=1}^{4} Gi$ from Table 2.2.4, we construct the graph of pressure distribution over time along the length of the pipeline (Figure 2.2.3).

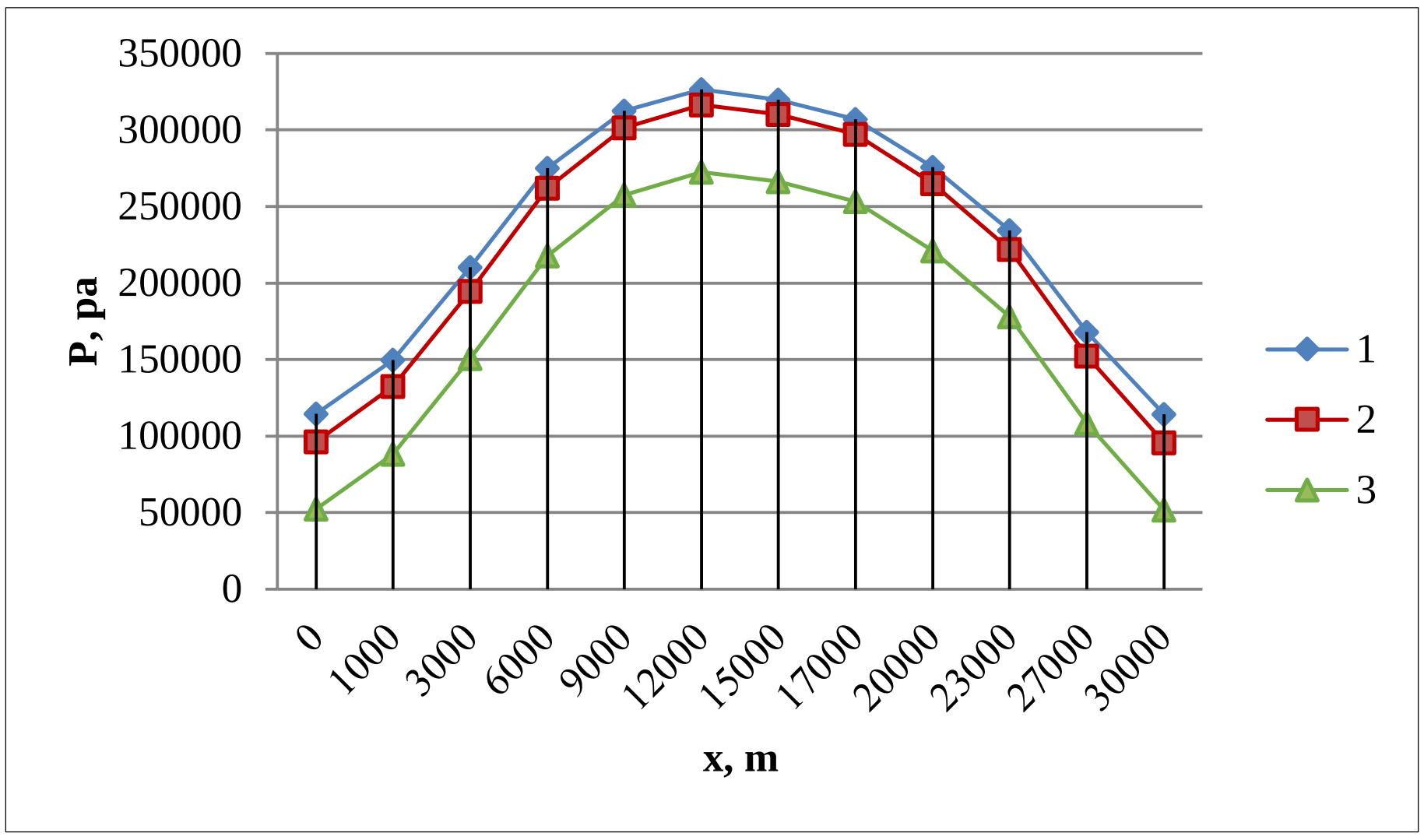

**Figure 2.2.3. Distribution of pressure over time along the length of the pipeline for a ring gas pipeline with side flows. (1-t=50 sec, 2-t=300 sec, 3- t=900 sec).**

Analysis of Table 2.2.4 and Figure 2.2.3 shows that increasing the total volume of side flows by 50% while keeping the number and location of side flow junctions constant has a significant impact on pressure values in closed loops. Specifically, when the total side gas flow is $\sum_{i=1}^{4} Gi$ = 10 Pa×sec/m, the pressure at the start of the gas pipeline decreases from $14\times10^4$ Pa to $12.3\times10^4$ Pa after 50 seconds. When the total side gas flow increases to $\sum_{i=1}^{4} Gi$ = 15 Pa×sec/m, the pressure at the start decreases from $14\times10^4$ Pa to $11.4\times10^4$ Pa, resulting in a difference of $0.9\times10^4$ Pa in the pressure drop.

After 300 seconds, the difference in the pressure drop at the start of the gas pipeline is $1.4\times10^4$ Pa, and after 900 seconds, it increases to $2.9\times10^4$ Pa. On the other hand, moving from the start towards the middle of the pipeline, the pressure increases and reaches a maximum value near the middle. At 50 seconds, this maximum pressure is $26\times10^4$ Pa, whereas with a total side gas flow of 15 Pa·s/m, the maximum pressure increases to $32\times10^4$ Pa. The increase in the maximum pressure due to the increase in side gas flow is $6\times10^4$ Pa. After 300 seconds, the increase in maximum pressure at the midpoint of the gas pipeline is $5.6\times10^4$ Pa, and after 900 seconds, this increase decreases to $4.2\times10^4$ Pa.

Therefore, as the side gas flow increases, the pressure drop at the start of the ring pipeline increases over time, while the increase in the maximum pressure at the midpoint of the pipeline decreases (Figures 2.2.1 and 2.2.3). The specific feature of the gasodynamic process in ring networks is that as the volume of side flows increases, the pressure reaches its maximum value at the hydraulic stability point of the pipeline ($x = x_1$), while the pressure reserve at the start decreases. This characteristic should be taken into account when reconstructing ring gas pipelines.

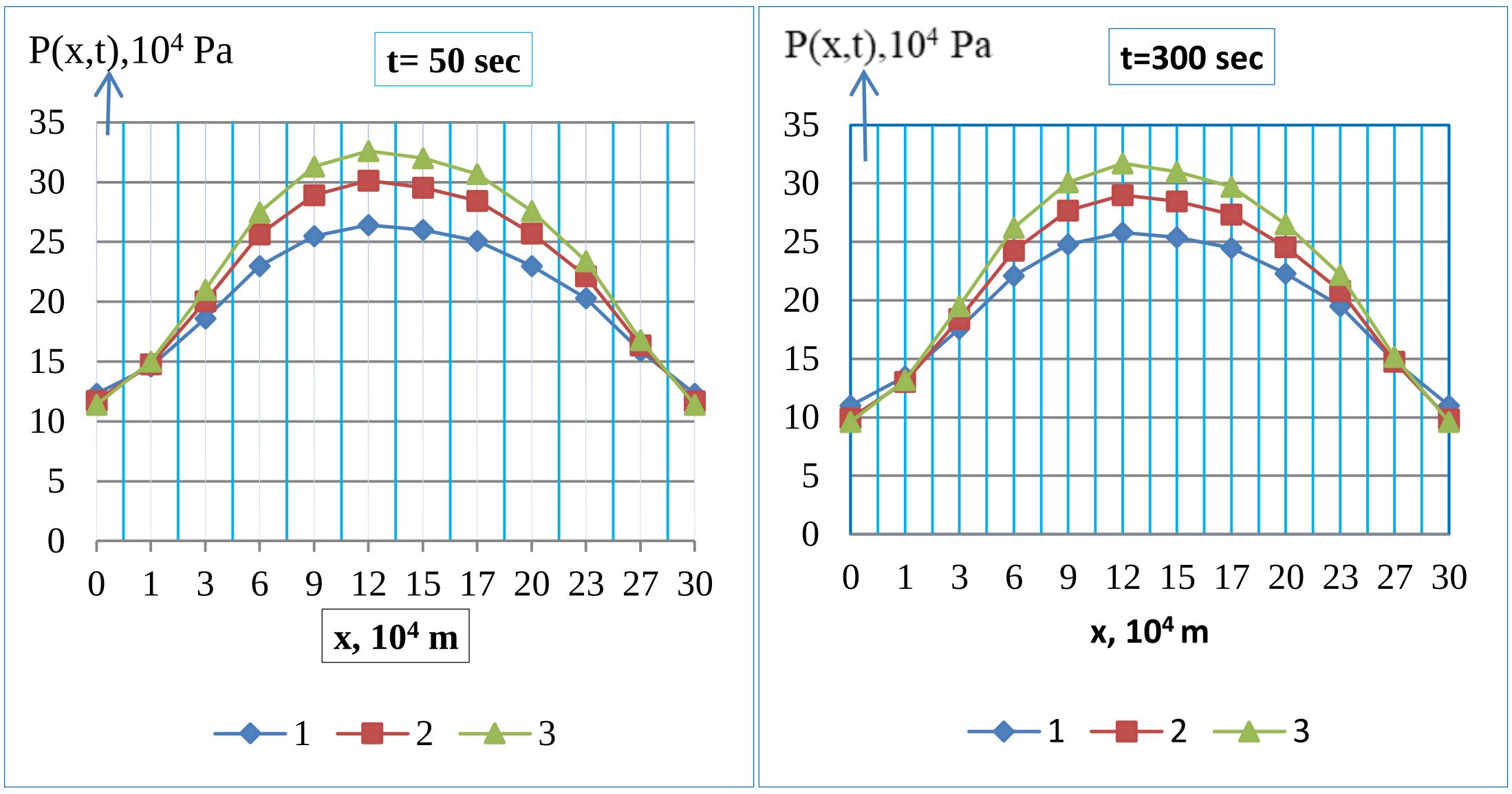

**Figure 2.2.4. Distribution of pressure along the length of the ring gas pipeline over time with varying total side flows mənzərəsi ( 1- $\sum G_i$=10 Pa ·sec/m, 2- $\sum G_i$ = 13 Pa ·sec/m, 3- $\sum G_i$ =15 Pa ·sec/m).**

Analysis of Figure 2.2.4 reveals that increasing the total side flows has a minor impact on the pressure distribution along the entire length of the ring network. Specifically, a 30% increase in total flows results in a 14% increase in the maximum pressure at the hydraulic junction after 50 seconds (from $26.4 \times 10^4$ Pa to $30.1 \times 10^4$ Pa). However, the pressure at the start and end points of the ring decreases by 5% (from $12.3 \times 10^4$ Pa to $11.7 \times 10^4$ Pa). After 300 seconds, the maximum pressure at the hydraulic junction increases by 12.4% (from $25.8 \times 10^4$ Pa to $29 \times 10^4$ Pa), while the pressure at the ring's start and end points decreases by 11% (from $11 \times 10^4$ Pa to $9.9 \times 10^4$Pa). Similarly, a 50% increase in total side flows results in a 23.5% increase in the maximum pressure at the hydraulic junction after 50 seconds. The pressure at the start and end points of the ring decreases by up to 10.8%. After 300 seconds, the maximum pressure at the hydraulic junction increases by 22.9%, while the pressure at the start and end points decreases by up to 14.6%.

Based on the analysis results, we can conclude that increasing the initial flow rates of the existing network by up to 30% can be considered an efficient option to meet the required reliability indicator for the ring gas distribution network. This is because an increase beyond this volume results in a pressure drop of more than 11% at the start and end points of the ring, leading to a gradual decrease in pressure reserve over time.

Assuming that the pressure difference for all consumers is fully utilized, a decrease in the initial pressure will also cause a decrease in pressure at point x1, and the pressure of the gas supplied to the consumers fed from this point will drop. Simultaneously, the decrease in pressure will affect the operation of equipment and devices using the gas. Devices and equipment that operate at maximum pressure will immediately feel the drop in pressure, which will reduce the pressure they experience. In such conditions, only devices and equipment that can operate close to the nominal pressure will be functional. However, it will not be possible to sufficiently lower the hydraulic resistance of the equipment. As the pressure in the devices decreases, the duration of their effective operation will increase. This is explained by the fact that as the operation time of the devices increases, the number of devices operating simultaneously will also increase. Consequently, the gas supply regime for consumers will be disrupted. Therefore, during the reconstruction phase of the ring gas pipelines, it is crucial to develop an analytical calculation scheme for optimizing the volume of gas consumption by new consumers and determining the location of hydraulic junctions.

## Justification of the Reconstruction Solution for the Operating Ring Gas Pipeline

The principle of reconstructing ring gas supply systems is based on the analyses we have reviewed above. Specifically, based on these analyses, the goal is to identify an efficient reconstruction option by utilizing the existing ring network to supply gas to new consumers. The application of this option results in the hydraulic connection point ($x=x_1$) of the ring gas supply system being used as the energy supply source for new consumers connected to the existing network. In this way, new consumers can be supplied with gas without any dismantling or assembly work on the existing network. From this perspective, the reconstruction option is considered economically viable.

It is clear that the choice of this option can depend on the volume of new consumers connecting to the ring network (i.e., the volume of increased gas consumption in the existing network). In other words, the connection of new consumers to the existing network should not lead to a decrease in the operating quality of equipment and devices or in the performance of the gas supply system. Given that the pressure at the hydraulic connection (hydraulic junction) point will not decrease due to the reduction in initial consumption, the reduction in the equivalent resistance faced by the equipment and devices should be such that the amount of gas supplied to consumers remains unchanged. Thus, the reconstruction option is conducted by determining the network's required reliability based on its transmission capacity.

Ring gas pipelines are typically used to supply gas to consumers in large residential areas because, as noted, these networks have high reliability and accumulation capacity. Therefore, such networks are classified as low, medium, and high pressure. For consumers supplied by low and medium pressure ring networks, the gas consumption regime can be regulated by changing the pressure at the junction. This is explained by the fact that there are no devices or equipment regulating gas pressure or consumption between the regulating point and the consumers. In this case, the amount of gas consumption G0 used at the beginning of the network can be increased to a certain extent without significantly altering the initial pressure. However, as the consumption increases, the initial pressure will decrease, and this decrease is not linear.

Based on the research results, one notable feature of the ring gas pipeline is that the variation in flow parameters under non-stationary conditions does not depend on the number and location of the connected in-line losses. This is because the pressure is a continuous function along the entire length of the ring network, and the gas-dynamic calculations of the ring are governed by Kirchhoff's first law at all network junctions. Therefore, it only changes according to the total amount of in-line losses.

Based on the above, we can assume the equalities $\sum_{i=1}^{n} G_i = G_1$ and $x_i = x_1$ in equation (12). Then, equation (12) can be expressed as follows:

$$P(x,t) = P_1 - 2aG_0\frac{L}{2} + 2aG_0L\sum_{n=1}^{\infty} Sin\frac{\pi nx}{L}\left(\frac{1-e^{-n^2\alpha t}}{\pi n^3}\right) - \frac{c^2t}{L}G_1 - \frac{2c^2}{L}G_1\sum_{n=1}^{\infty} Cos\frac{2\pi n(x-x_1)}{L}\left(\frac{1-e^{-n^2\alpha t}}{\alpha n^2}\right) \qquad (12)$$

Here, $G_1$ is the increased consumption of the existing ring gas pipeline corresponding to the volume of gas consumption required by new consumers ($G_1=G_0+q_0$), in Pa · sec/m.

The use of the ring pipeline aims to improve the reliability of gas supply [5]. Continuous delivery of gas to consumers does not fully solve the reliability problem. Additionally, the integration of new consumers into the existing ring network must ensure that the required gas consumption volume is met with reserve road consumption and inlet pressure according to the demand. This introduces new requirements for the methodology of calculating the dynamic state of gas flow along the ring pipeline.

From the analyzed graphs, it is clear that the ring gas pipeline achieves its maximum value at x= $x_1$ regardless of the total road consumption and the number of junctions (Graphs 2.2.1 to 2.2.5). The main goal of the reconstruction solution for ring gas pipelines is to develop a calculation scheme for the pressure at the connection point for new consumers. It is clear that to ensure uninterrupted delivery of gas to consumers, it is appropriate to select the section of the ring where the gas flow pressure achieves its maximum value.

Thus, to obtain the stated results, we take the derivative of equation (12) with respect to x. The resulting equation will be as follows:

$$\frac{dP(x,t)}{dx} = 2aG_0 \sum_{n=1}^{\infty} Cos\frac{\pi nx}{L}\left(\frac{1-e^{-n^2\alpha t}}{n^2}\right) - \frac{4\pi c^2}{L^2} G_1 \sum_{n=1}^{\infty} Sin\frac{2\pi n(x-x_1)}{L}\left(\frac{1-e^{-n^2\alpha t}}{\alpha n}\right) \quad (13)$$

To determine the rate of change of pressure at the point x=$x_1$ using equation (14), you would typically follow these steps:

$$\frac{dP(x_1,t)}{dx} = 2aG_0\left[\left(\frac{\pi^2 x_1^2}{4L^2} - \frac{\pi^2 x_1}{2L} + \frac{\pi^2}{6}\right) - \sum_{n=1}^{\infty} Cos\frac{\pi n x_1}{L}\left(\frac{e^{-n^2\alpha t}}{n^2}\right)\right] \quad (14)$$

To perform engineering calculations and determine the location of the hydraulic junction node ($x_1$) in the ring gas pipeline, we accept the simplification n=2 for the application of expression (14) within its known domains, based on theoretical and experimental tests [28]. Then,

$$x_1 = L\frac{1-\sqrt{1-\frac{4Z}{\pi^2}\left[\left(\frac{2\pi^2-3e^{-4\alpha t}}{12}\right)+\frac{G_1}{G_0}\right]}}{Z}, \quad m \quad (15)$$

Here, $Z = (1+8e^{-4\alpha t})$, $\alpha = \frac{4\pi^2 c^2}{2aL^2}$

As mentioned above, in both stationary and non-stationary regimes, the pressure function across the entire length of the ring network reaches its maximum value at the point x=$x_1$. In this case, the location of the hydraulic junction node $x_1$ is determined by formula (15).

To confirm our statements and for engineering calculations, we accept the known data and the parameters obtained below.

$P_1 = 14 \cdot 10^4$ Pa : $G_0$= 10 Pa ·sec/m : $2a = 0{,}1\frac{1}{sec}$ ; $c = 383{,}3\frac{m}{sec}$ ;d= 0,7 m: L = $3 \cdot 10^4$ m

The location of the hydraulic junction node x1 is determined every 50 seconds (at t=50,150,200,250,…,700 seconds) by considering the given data in formula (15) for different values of the ratio between the increased flow rate corresponding to the new consumers' gas demand and the existing initial flow rate of the gas pipeline ($G_1/G_0$=1;1,05;1,1;1,15;1,2;1,25; 1,3;1,35;1,4;1,45;1,5;1,55). The results of the calculations are presented in the following graph (Graph 2.2.5).

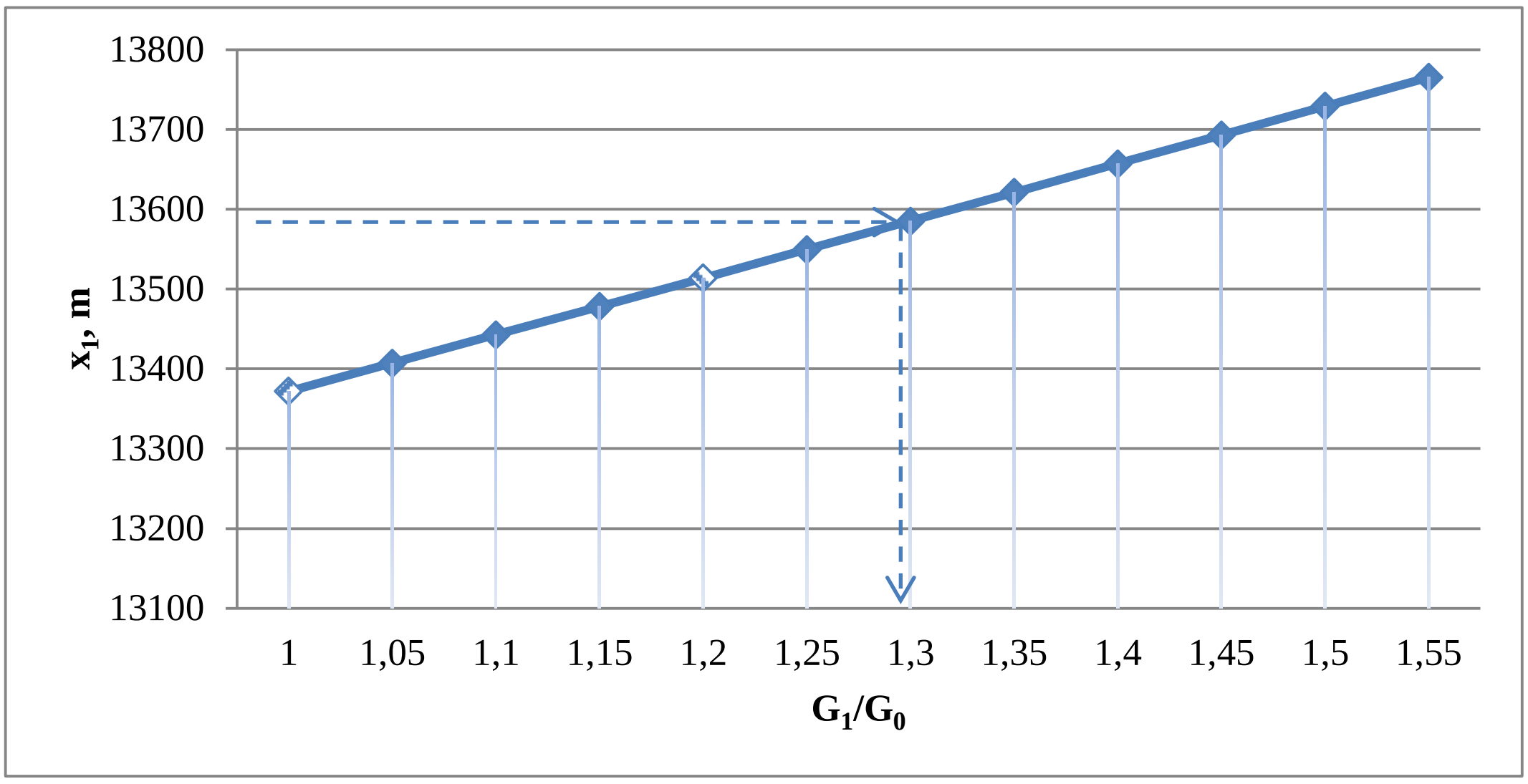

**Graph 2.2.5: Dependence graph of the hydraulic junction node of the ring gas pipeline on the ratio of the increased flow to the existing initial flow rate for the reconstruction.**

From the analysis of Graph 2.2.5, it becomes clear that the dependence of the hydraulic junction node x1, which characterizes the hydraulic junction, on the ratio of the increased flow rate to the initial flow rate during the reconstruction stage remains constant over time. However, the value of the flow rate ratio is linearly dependent on the value of x1. Specifically, as the flow rate increases, the value of the hydraulic junction node also increases.

First, we determine the actual value of the maximum hydraulic junction node for different scenarios of the total and number of passing flows using equation (12). The results of the calculation are presented in the following graph (Graph 2.2.6).

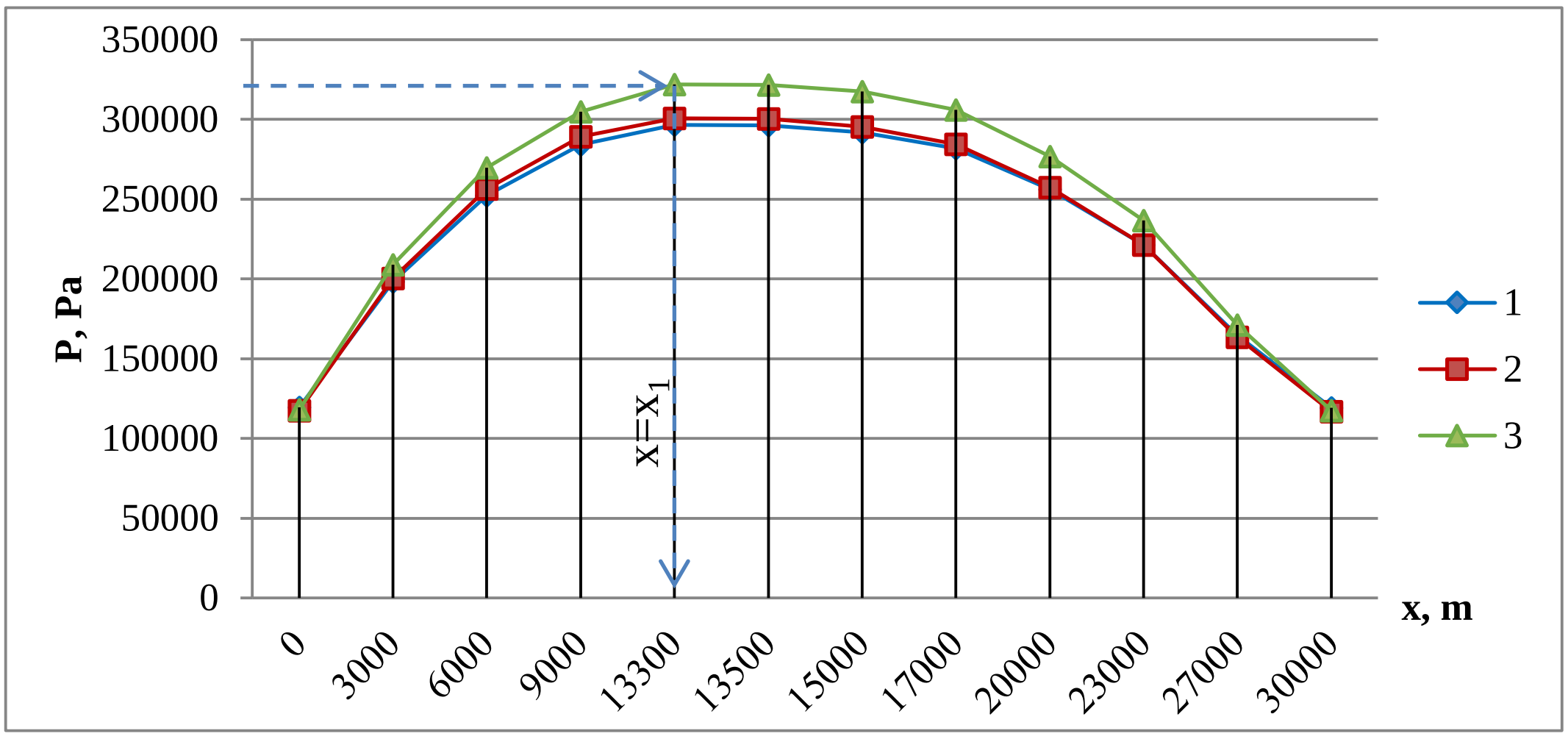

**Graph 2.2.6: Pressure distribution profile along the length of the ring gas pipeline, depending on the total value of the passing flows and the number of junctions.**

**(1- $\sum_{i=1}^{3} G_i$ =12 Pa×sec/m, 2- $\sum_{i=1}^{4} G_i$ =13 Pa ×sec/m , 3- $G_1$ -14Pa ×sec/m, t=50 sec).**

Based on the analysis of Graph 2.2.6, we can conclude that regardless of the different values of the total passing flows and the number of junctions, the ring gas pipeline always reaches its maximum value at the point x = $x_1$. Therefore, considering the findings from the previous analyses, it can be stated that the location of the hydraulic junction point in a ring pipeline with passing flows depends solely on the design parameters of the gas pipeline. Consequently, the location of the hydraulic junction points will vary at different distances depending on the value of the parameters of each ring pipeline with passing flows. From this perspective, it is crucial to analytically determine the location of the hydraulic junction point during the reconstruction phase of ring gas pipelines.

Next, we compare the actual maximum value of the hydraulic junction point with the calculated value. According to Graph 2.2.6, the actual maximum value of the hydraulic junction point is $x_1^f = 13300m$. The calculated maximum values of the hydraulic junction point are determined using equation (15) depending on the ratio of the initial flow rate of the pipeline, increased to meet the gas consumption needs of new consumers, and time. The comparison of the results for both scenarios is presented in the following table (Table 2.2.5).

| t, sec | $G_1/G_0$ | **Actual value** $x_1^f$, m | **Calculated value** $x_1^h$, m | **Relative error** $\frac{x_1^h - x_1^f}{x_1^h}$ |
|---|---|---|---|---|
| 50 | 1 | 13300 | 13372,0 | 0,005 |
| 100 | 1,05 | 13300 | 13407,3 | 0,008 |
| 150 | 1,1 | 13300 | 13442,8 | 0,011 |
| 250 | 1,15 | 13300 | 13478,3 | 0,013 |
| 300 | 1,2 | 13300 | 13513,9 | 0,016 |
| 400 | 1,25 | 13300 | 13549,6 | 0,018 |
| 450 | 1,3 | 13300 | 13585,3 | 0,021 |
| 500 | 1,35 | 13300 | 13621,2 | 0,024 |
| 550 | 1,4 | 13300 | 13657,1 | 0,026 |
| 600 | 1,45 | 13300 | 13693,1 | 0,029 |
| 650 | 1,5 | 13300 | 13729,2 | 0,031 |
| 700 | 1,55 | 13300 | 13765,3 | 0,034 |

**Table 2.2.5. Comparison of the calculated and actual values of the hydraulic junction point ($x_1$) in different scenarios of total passing flows and the number of junctions in ring gas pipelines over time.**

The slight relative error observed is due to the linearization of nonlinear equation systems and the simplification of mathematical expressions. However, this is acceptable for engineering calculations and for the reconstruction solution of ring pipelines. The mathematical expressions derived and the analytical formula obtained are suitable for practical use in studying physical processes.

Based on the data from Graphs 2.2.5 and 2.2.6, as well as the analysis of Graph 2.2.4, it can be noted that increasing the initial flow rate by up to 30% to match the maximum value of the hydraulic junction point of the studied ring gas pipeline does not cause significant changes in the gas dynamic processes of the pipeline. Therefore, if the gas supply for new construction within the area served by the existing ring gas pipeline amounts to up to 30% of the initial flow rate of the ring pipeline, accepting the hydraulic junction node of the ring pipeline as the energy supply source for new consumers is technologically and economically viable. Considering this option during the reconstruction phase of existing ring gas pipelines can lead to a significant increase in the budget revenues of the gas industry.

Thus, the adoption of the option of using the hydraulic connection node of the ring pipeline as an energy supply source for new consumers in the reconstruction of ring gas transmission systems is scientifically and technologically justified. To make the reconstruction even more effective, the location of the hydraulic connection node is determined using formula (15), and this node is accepted as the energy source for new consumers.

If the volume of gas supplied to new consumers exceeds 30% of the initial flow rate, this reconstruction option is often not economically or technologically feasible. Therefore, in the reconstruction of ring systems, two additional options are considered for comparison [1,94]. In the first option, if the gas consumption required for new consumers exceeds 30%, dismantling and installation work on the existing network pipelines is necessary to avoid abnormal gas dynamics caused by the connection to the ring network. In other words, the diameter of the ring pipeline is increased according to the gas consumption. In the second option, an independent gas line is designed and constructed for new consumers from the source of the ring gas pipeline system (compressor station, gas distribution station, etc.). In this case, the existing network supplies the existing consumers, while the independent network supplies the new consumers separately with gas. On the other hand, if the existing gas pipeline is severely physically and morally worn out, the technical-economic indicators of both the option of using the existing gas pipeline and the option of designing an independent new gas pipeline may not meet the reconstruction requirements due to the system's unfitness for operation. In such cases, an additional third option is also considered for comparison in the reconstruction of ring gas pipelines. In this option, the existing network sections are dismantled, and the unusable parts are removed. As a result, the distribution scheme of the system is renovated, taking into account the new consumers. Modern equipment and devices are applied in this option.

## CONCLUSION

1. The solution to the differential equations system describing the non-stationary gas flow in a ring pipeline with a constant cross-sectional area has been obtained. The requirements for the equality of boundary conditions included in the solution allow for a compact representation of the analytical model, but do not restrict the model's applicability to sharp changes in gas flow or pressure. Thus, the methodology for calculating non-stationary gas flow in ring pipelines has been developed.
2. The impact of the flow characteristics of the ring gas pipeline on the gasodynamic process has been investigated. The results of gas pressure calculations for various volumes of line gas flows and connection points have been analyzed and evaluated. The comparison of computational examples and their theoretical analyses allows for judgment on the accuracy of the proposed mathematical expressions for pressure calculations in ring gas pipelines, both in terms of length and time.
3. An analytical method has been developed based on the technological principles of ring gas pipeline system reconstruction solutions. The proposed analytical method can be generalized for the ring network of pipelines with different gas flow intensities in various zones of the ring, facilitating the reconstruction process and allowing for accurate determination of the connection point for new consumers in both stationary and non-stationary modes.

## 2.3. Advanced Engineering Approaches for Efficient Ring Gas Pipeline Operation and Emergency Management

**Introduction.** The operation, management, and optimization of ring gas pipelines is of significant importance for solving modern engineering problems. The recent application of advanced technologies, particularly machine learning facilitates the resolution of engineering problems. Gas pipeline systems are complex systems characterized by numerous elements, and managing them requires a large amount of information corresponding to the number of possible elements. Managing such complex gas pipelines is impossible without a systematic approach based on the combined consideration of complex concepts of control theory, such as system, information, goal-orientation, and feedback. New approaches to operation methods aimed at ensuring a higher level of quality are required in managing these systems.

The theory of selecting justified parameters for the operation of ring gas pipelines requires the calculation of non-stationary processes resulting from any emergency [85,87, 100,103].

The difficulty in modeling unsteady flow regimes for such systems is related to the complexity and variability of the dynamic state of gas flow in pipelines, as well as the need to account for many different factors. Our goal is the operation of the ring distribution network in accordance with the requirements for efficiency and reliability indicators.

When studying the gas-dynamic issues of ring networks, it should be noted that in both distribution rings diverging from the starting point, gas flows in opposite directions and intersects at any point, creating equilibrium at the intersection. At the starting point of the gas pipeline, the pressures in the hydraulic junction sections of both distribution branches equalize. If the integrity of the gas pipeline is violated in the x=ℓ section, then the flow streams of each distribution branch intersect at the leak point x=ℓ (Figure 2.3.1).

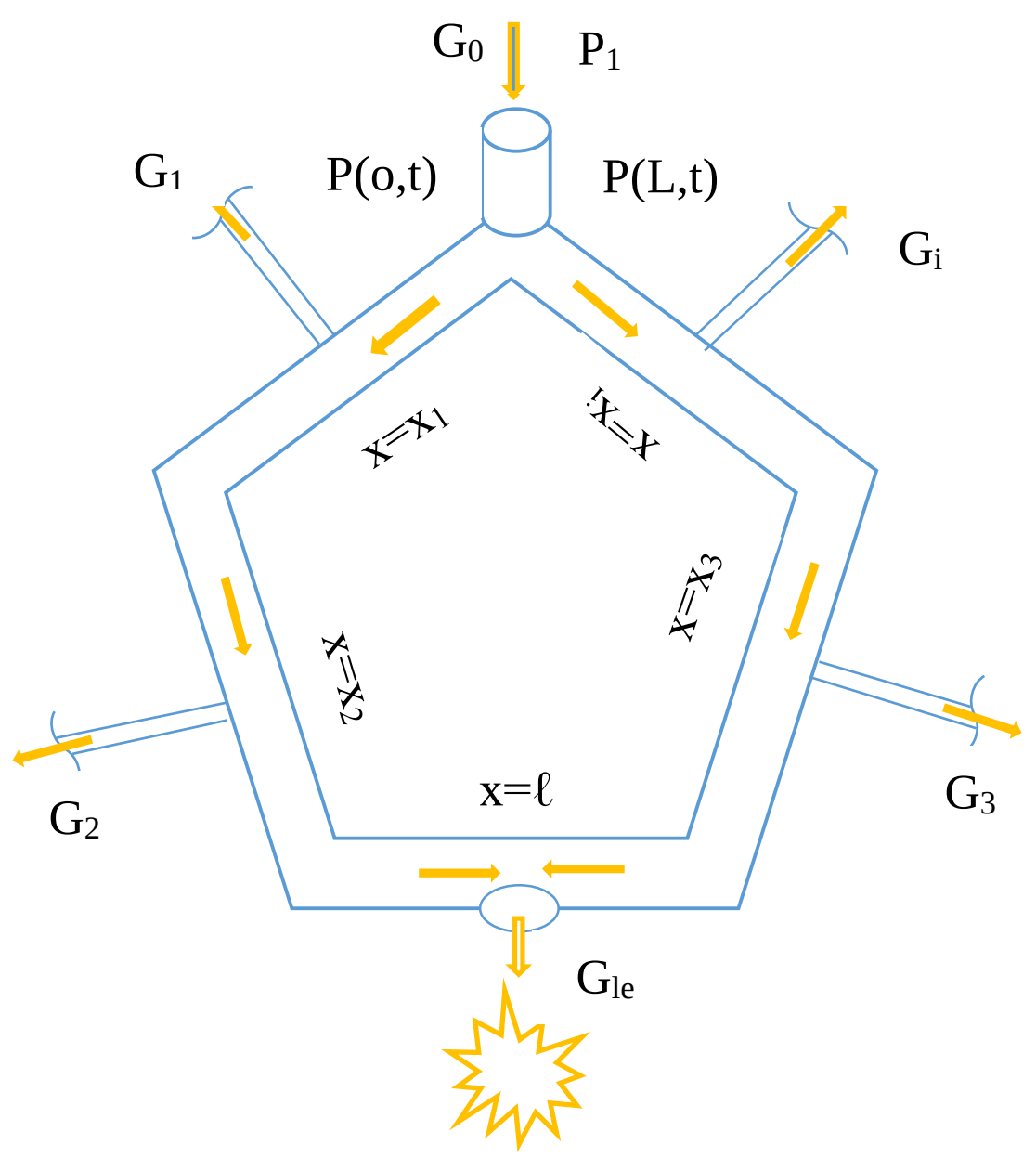

**Figure 2.3.1. Visual scheme of a circular gas supply system.**

Leaks in pipeline networks remain one of the primary causes of countless losses to pipeline stakeholders and the surrounding environment. Pipeline failures can lead to severe environmental disasters, human casualties, and financial losses. To prevent such hazards and maintain a safe and reliable pipeline infrastructure, it is essential to implement early warning and leakage detection technology for pipeline safety. The pipeline leakage detection technologies are discussed in [84, 98, 113], summarizing the latest advancements. Various leak detection and localization methods in pipeline systems are reviewed, with their strengths and weaknesses highlighted. A comprehensive analysis of the current state of machine learning and artificial intelligence in solving gas supply issues has been conducted by several researchers [ 98, 89, 111,114]. These studies also discuss various types of machine learning and artificial intelligence methods that can be used to process and interpret data in different areas of the oil and gas industry.Comparative performance analysis is carried out to provide guidance on which leak detection method is most suitable for specific operational parameters. Furthermore, research gaps and open issues for the development of reliable pipeline leakage detection systems are discussed.

Due to their complexity and scale, pipelines are susceptible to leaks that can result in severe environmental, economic, and safety consequences. Effective leak detection (LD) program management is necessary for mitigating these risks, but it is not enough to have an LD system, it must be designed and installed properly.

Existing leak detection and localization methods are known to be applied to gas pipelines.

The network architecture and hierarchical model for detecting and localizing leaks are described as follows [83]:

Natural gas pipeline monitoring sensor networks are made up of a great number of sensor nodes, the sinks, and a control and management center. The whole network adopt a sort of clustering structure, which is described in Figure 2.3.2 The sensor nodes, installed along the pipeline, are in charge of signal collecting and data preprocessing. Then the processing results are sent to the sink via a multi-hop route. The monitoring networks take an in-tube communication mode for message exchange between sensor nodes. The sink, as the cluster leader, is responsible for managing sensor nodes in its cluster,combining preprocessing results for final decision and localizing leak position if a leak has occurred. Finally, the diagnosis result is sent to the control and management center, which judges whether a leak warning alarm should be announced.

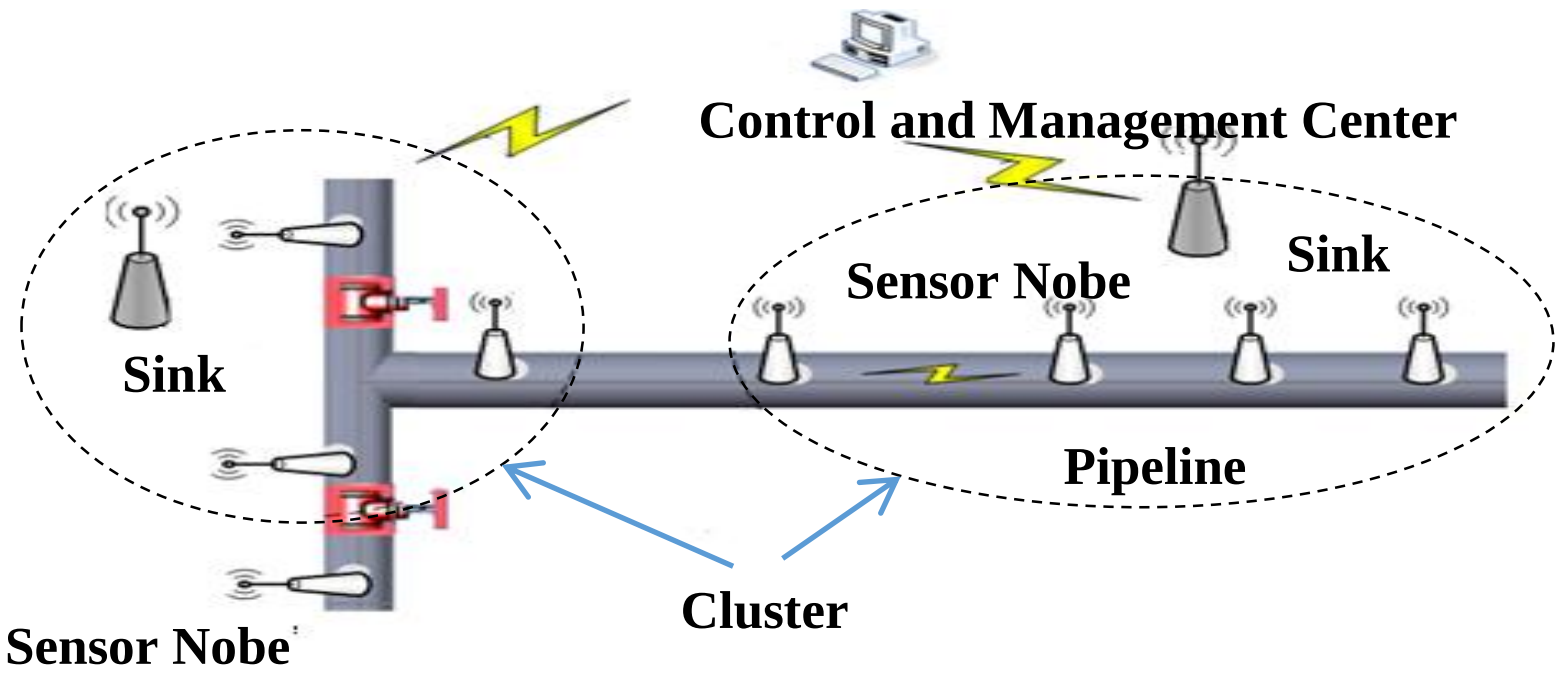

**Figure 2.3.2 Architecture of pipeline monitoring sensor networks [83].**

The localization method based on time difference of arrival (TDOA) is described as follows: [83].

Several sensor nodes are placed at fixed points to make up the sensor array in accordance with some geometric relationship. These nodes measure and calculate the relative time difference of signal waves transmitting to each sensor node. The time differences are substituted into the equations which satisfy the geometric relationship of sensor array for solving the leak point coordinate. Furthermore, leak signals will lose some energy during the process of propagation, that is to say, signal intensity will decline with the increase of the distance between leak point and sensor position. According to this feature, the leak source region can be roughly determined via the signal intensity received by sensor nodes and leak points can be accurately localized by analyzing the signal attenuation characteristics. Figure 2.3.3 shows the principle of TDOA localization through two sensor nodes.

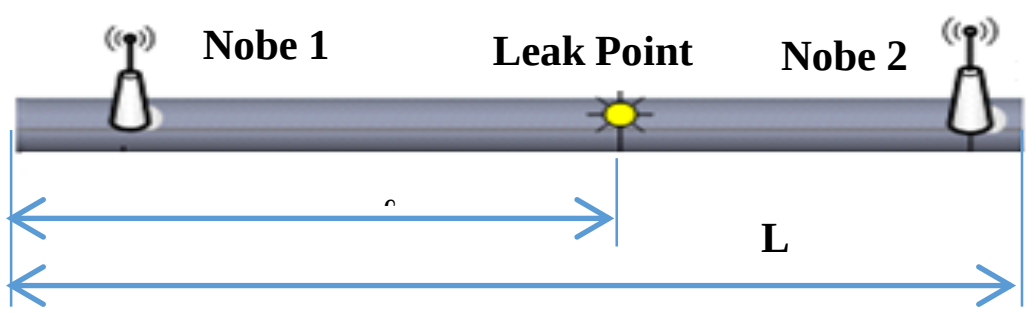

**Figure 2.3.3. Principle of time difference of arrival localization [83].**

There are also different types of machine learning and artificial intelligence techniques that can be used to process and interpret data in different areas of the gas industry. Achievements and developments promise the advantages of machine learning and artificial intelligence methods in terms of large storage capabilities and high efficiency of numerical calculations. Using these latest technologies, development and progress have become smarter, making the decision-making procedure simple and straightforward.

It is known that recently, machine learning and artificial intelligence systems have been present in various fields, including the management of gas supply.

There are three methods for detecting leaks in gas pipelines: qualitative, quantitative, and analytical. The practical examples given mainly concern the quantitative method.

The purpose of our research was focused on the method of detecting cases of violation of the integrity (tightness) of circular gas pipelines, that is, emergency situations, based on an analytical database. In particular, it involved the development of a reporting indicator for the effective management of gas dynamic processes in non-stationary modes.

For ring gas pipeline systems, new approaches have been investigated for modernizing existing control process monitoring systems. These approaches are based on modern advancements in control theory and information technology, aimed at selecting emergency and technological modes. The diagram showing the methodological process is in the figure below.

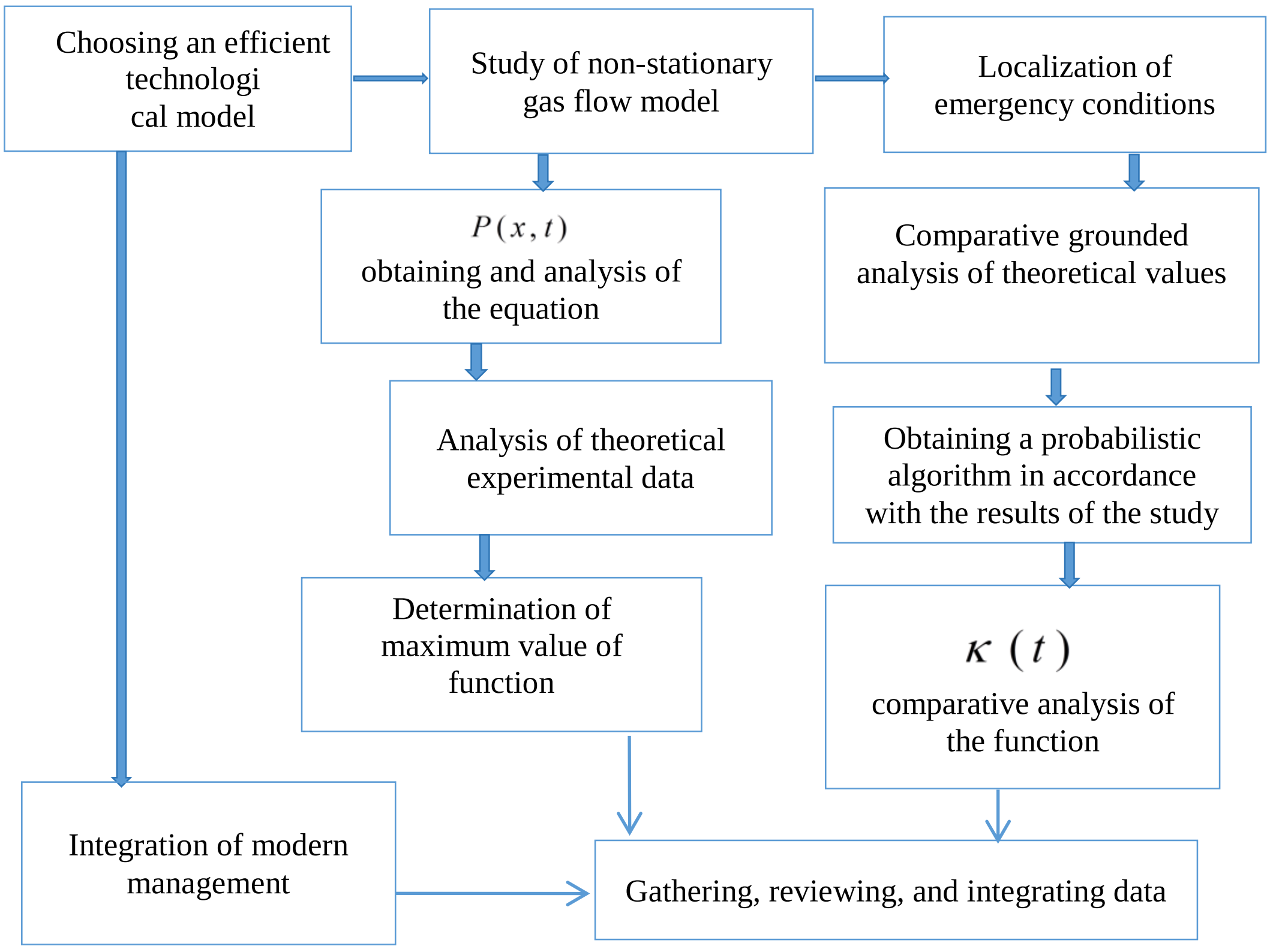

Analysis of mathematical model for identifying methodology for unsteady gas flow in ring pipelines. When addressing solutions to the considered problem, it is advisable to obtain solvable mathematical expressions for studying unsteable processes and managing the occurring physical processes. Therefore, using unsteady-state models is appropriate.

These models enable the resolution of various issues, including those related to the operation of ring networks. Consequently, models and algorithms for calculating hydraulic schematics of ring gas pipelines need to be refined. For unsteady gas flow, I.A. Charnynin's differential equations system relates gas flow parameters, such as pressure P(x,t) and gas mass flow G(x,t), as functions of time t and distance x. Based on Figure 2.3.1, the system of specific derivative non-linear equations for the problem under consideration can be expressed as follows [87,88].

$$\begin{cases} -\dfrac{\partial P}{\partial x} = \lambda \dfrac{\rho V^2}{2d} \\ -\dfrac{1}{c^2}\dfrac{\partial P}{\partial t} = \dfrac{\partial G}{\partial x} + \sum_{i=1}^{n} G_i(t)\delta(x - x_i) + G_{le}\delta(x - \ell) \end{cases} \tag{1}$$

Here, $G = \rho v$; $P = \rho c^2$ ; $\rho = \dfrac{P}{zRT}$

T – Average temperature of the gas.

z=z (p,T) – Compressibility factor.

c – Speed of sound in the gas for an isothermal process, m/sec.

P – Pressure in the cross-sections of the gas pipeline, Pa.

x – Coordinate along the axis of the gas pipeline, m.

t – Time coordinate, sec.

λ – Hydraulic resistance coefficient.

ρ – Average density of the gas, kg/m$^3$.

v – Average velocity of the gas flow, m/sec.

d – Diameter of the gas pipeline, m.

G – Mass flow rate of the gas.

$\sum_{i=1}^{n} Gi(t)$ –the sum of the consumption of the mass of gas allocated to consumers from different points of the ring gas pipeline. $\dfrac{\text{Pa} \times \text{sec}}{\text{m}}$

$G_{le}(t)$ – the flow rate of the mass of gas and the mass of gas flowing from the point $x = \ell$, $\dfrac{\text{Pa} \times \text{sec}}{\text{m}}$

$\delta(x - x_i), \delta(x - \ell)$ – Dirac delta function.

It is evident that the nonlinear system of equations (1) cannot be solved directly, as such equations are not integrable. In engineering calculations, approximate methods are employed. In other words, the equations are linearized. Linearization of equations, as proposed by I.A. Charnynin is more practical for solving these problems. The expression characterizing this linearization is given as follows [85].

$$2a = \lambda \frac{\text{v}}{2\text{d}} \tag{2}$$

After linearization, Equation (1) becomes the heat transfer equation for the mathematical solution of the problem.

$$\frac{\partial P^2}{\partial x} = \frac{2a}{c^2}\frac{\partial G}{\partial x} + 2a\sum_{i=1}^{n} G_{le}(t)\delta(x - x_i) + 2aG_{le}\delta(x - \ell) \tag{3}$$

Thus, for the uniqueness of the solution to equation (3) that reflects the studied process and for the complete determination of the given process, the initial and boundary conditions must be correctly specified. The sum of the intensities of the concentrated road-based flow rates at the pipeline junctions is considered equal to the initial flow rate in steady-state regimes. The diameter of the investigated gas pipeline is uniform throughout. The total mass flow rate of the gas is $G_0$, and the pressure at the starting point is $P_1$.

Then, t=0; $P(x,0)=P_1-2aG_0x$: $G_0x=G_ix_i$, i=1,2,3,...,n

For the main sections of the gas pipeline, we assume the following boundary conditions for the process under consideration (continuity of flow at the initial and final points):

1.Initial Point: The flow rate and pressure at the beginning of the pipeline are set according to the given conditions.

2.Final Point: The flow rate and pressure at the end of the pipeline are adjusted to ensure continuity and adherence to the predetermined law during the process.

$$x=0,L \quad \text{point} \quad \begin{cases} P(0,t)=P(L,t) \\ \dfrac{\partial P(0,t)}{\partial X}=\dfrac{\partial P(L,t)}{\partial X} \end{cases}$$

It is clear that by applying the Laplace transform $\left[P(x,s)=\int_0^\infty P(x,t)e^{-st}dt\right]$, equation (3) is transformed into second-order ordinary differential equations. The general solutions for these equations will be as follows:

$$P(x,s)=\sqrt{\frac{A}{s}}c^2\left[\sum_{i=1}^{n}G_i(t)\int_0^x sh\sqrt{As}(x-y)\delta(y-x_i)dy+G_{le}(s)\int_0^x sh\sqrt{As}(x-y)\delta(y-\ell)dy\right]-$$
$$-\sqrt{\frac{A}{s}}\int_0^x sh\sqrt{As}(x-y)(P_1-2aG_0y)dy+C_1sh\sqrt{As}x+C_2ch\sqrt{As}x \qquad (4)$$

Here, $A=\dfrac{2a}{c^2}$

By applying the Laplace transform to the initial and boundary conditions, we determine the constants $C_1$ and $C_2$ in equation (4). After By considering the values of expressions $C_1$ and $C_2$ in equation (4), we obtain the transformed equation for the dynamic state of the process, or in other words, for the ring-shaped gas pipeline. The resulting equation will represent the transformed distribution of pressure along the ring-shaped gas pipeline in a non-stationary flow regime.

$$P(x,S)=\frac{P_1-2aG_0x}{S}-\frac{2aG_0L}{S}\frac{sh\sqrt{As}\left(\frac{L}{2}-x\right)}{2sh\sqrt{As}\frac{L}{2}}-$$
$$\sqrt{\frac{A}{s}}c^2\sum_{n=1}^{n}G_i(s)\frac{ch\sqrt{As}\left(\frac{L}{2}-x+x_i\right)}{2sh\sqrt{As}\frac{L}{2}}-\sqrt{\frac{A}{s}}c^2G_{le}\frac{ch\sqrt{As}\left(\frac{L}{2}-x+\ell\right)}{2sh\sqrt{As}\frac{L}{2}} \qquad (5)$$

`

Based on the solution of the problem, since determining the mathematical expression for the physical process of the non-stationary flow of gas in the pipeline is more suitable, we use the inverse Laplace transform for each boundary to find the original solution of equation (5). This process results in the following original equation.

$$P(x,t)=P_1-2aG_0x-2aG_0\left(\left(\frac{L}{2}-x\right)+L\sum_{n=1}^{\infty}(-1)^n\frac{Sin(L-x)\frac{\pi n}{L}}{\pi n}\int_0^t e^{-n^2\alpha(t-\tau)}d\tau\right)-$$

$$-\frac{c^2}{L}\int_0^t\sum_{i=1}^{n}G_i(\tau)\times\left[1+2\sum_{n=1}^{\infty}(-1)^n Cos[L-2(x-x_i)]\frac{\pi n}{L}e^{-n^2\alpha(t-\tau)}\right]d\tau- \quad (6)$$

$$-\frac{c^2}{L}\int_o^t G_{le}\left[1+2\sum_{n=1}^{\infty}(-1)^n Cos[L-2(x-\ell)]\frac{\pi n}{L}e^{-n^2\alpha(t-\tau)}\right]d\tau$$

Here, $G_0x=\sum_{i=1}^{n}G_ix_i$, $\alpha=\frac{4\pi^2c^2}{2aL^2}$

If we assume that the road-based gas consumption used by consumers is constant ($G_i$ (t) $=G_i$ =const and $G_{le}(t)=G_{le}=\text{const}$ ) then through several mathematical operations, equation (6) will be simplified to the following form:

$$P(x,t)=P_1-2aG_0\left[\frac{L}{2}-L\sum_{n=1}^{\infty}Sin\frac{\pi nx}{L}\left(\frac{1-e^{-n^2\alpha t}}{\pi n^3}\right)\right]-\frac{c^2t}{L}\left[\sum_{i=1}^{n}G_i+G_{le}\right]-\frac{2c^2}{L}\left(\frac{1-e^{-n^2\alpha t}}{\alpha n^2}\right)\times$$

$$\times\left[\sum_{i=1}^{n}G_i\sum_{n=1}^{\infty}Cos\frac{2\pi n(x-x_i)}{L}+G_{le}\sum_{n=1}^{\infty}Cos\frac{2\pi n(x-\ell)}{L}\right] \quad (7)$$

To determine the law of variation of pressure in a ring-shaped network, we accept the following data:

$P_1 = 14\times10^4$ Pa : $G_0$= 10 Pa $\times$ sec/m $2a=0.1\frac{1}{\text{sec}}$ ; $c=383.3\frac{m}{san}$; d= 0.7 m: L = $3\times10^4$ m;

First, we look at the case when the integrity of the gas pipeline is not violated ($G_{le}$=o), and let's assume that the number of connected pipes for gas yield from the circular gas pipeline to consumers is three($G_1$ ,$G_2$, $G_3$) with the following parameters:($G_1$=3 Pa$\times$ sec/m ; G2=4 Pa$\times$ sec/m; $G_3$=3 Pa$\times$ sec/m)

And that the junction points are located at distances:

$x_1=0.3\times10^4$ m; $x_2=1.5\times10^4$ m: $x_3=2.7\times10^4$ m.

We calculate the pressure values of the gas pipeline as a function of time and distance, and record them in the following table (Table 2.3.1).

| x,m | P(x,t), Pa | | |
|---|---|---|---|
| | t=50 sec | t=300 sec | t=900 sec |
| 0 | 122717.5 | 110478.6 | 81094.8 |
| 1000 | 146064.4 | 134445.6 | 105061.8 |
| 3000 | 186412.1 | 176006.9 | 146623.1 |
| 6000 | 229710.7 | 220960.7 | 191576.9 |
| 9000 | 254914.2 | 247478.2 | 218094.4 |
| 12000 | 260046.5 | 253743.6 | 224359.8 |
| 15000 | 264319.1 | 257726.5 | 228342.7 |
| 17000 | 251713.7 | 245279.8 | 215896 |
| 20000 | 231199.6 | 224097.4 | 194713.6 |
| 23000 | 203327 | 195052.9 | 165669.2 |
| 27000 | 158388.6 | 147974.3 | 118590.6 |
| 30000 | 122357.1 | 110108.7 | 80724.91 |

**Table 2.3.1. Values of the Pressure Variation Along the Length of the Ring-shaped Network as a Function of Time ($G_1+G_2+G_3=G_0$).**

Using Table 2.3.1, we plot the graph of pressure distribution along the length of the ring-shaped gas pipeline as a function of time (Figure 2.3.4).

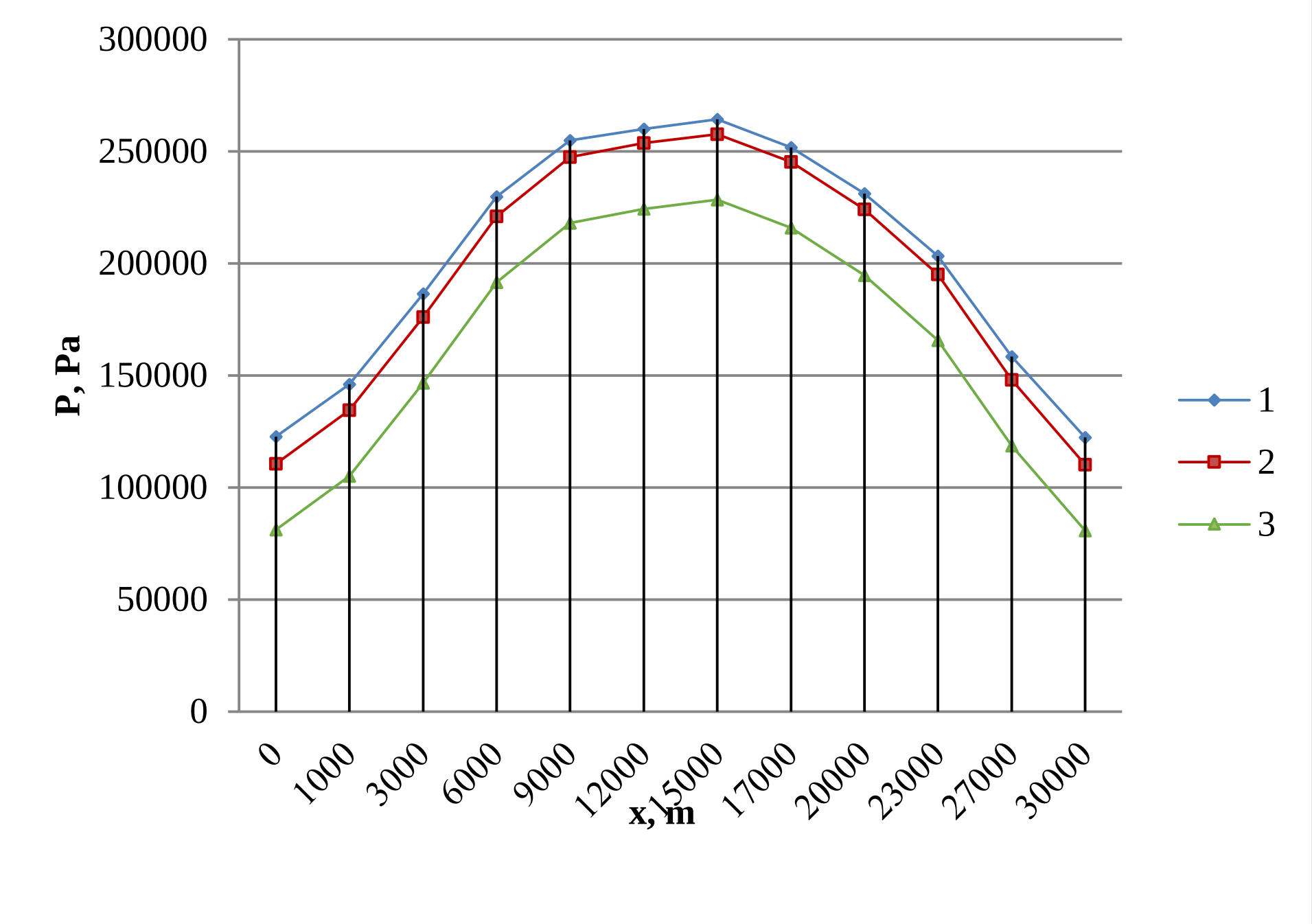

**Figure 2.3.4. Pressure Distribution Along the Length of the Road-based Gas Pipeline as a Function of Time (1-t=50 sec, 2-t=300 sec, 3- t=900 sec).**

The analysis of Figure 2.3.4 reveals the following:

**Pressure Increase and Maximum:** The pressure increases from the initial part of the pipeline towards the middle (counterclockwise direction) and reaches its maximum value near the hydraulic junction in the middle section of the pipeline.

Pressure Decrease: From the middle section towards the end of the pipeline (clockwise direction), the pressure decreases and eventually stabilizes to match the initial pressure at the starting point, i.e., [P(0,t)=P(L,t)].

Reason for Pressure Increase: The increase in pressure in the middle sections of the ring-shaped pipeline is due to the incremental extraction of the gas mass flow ($G_0$) from various points along the pipeline ($G_i x_i$ i=1,2,3,...,n)

Pressure Decrease Over Time: At the start, the pressure is $14\times 10^4$ Pa, and after 50 seconds, it decreases to $12.3\times 10^4$ Pa. Over time, this value gradually reduces, reaching $11.\times$ Pa at 300 seconds and $8.1\times 10^4$ Pa at 900 seconds (refer to Table 2.3.1). This decrease continues until a new stationary state is achieved.

Cause of Pressure Decrease: The decrease in pressure at the starting and ending points is due to the counterflow of the gas and hydraulic losses. According to gas dynamic processes, hydraulic losses are directly proportional to the gas mass flow.

Let us now analyze the conditions in which the integrity of the pipeline is violated at a distance of $\ell=2.0\times 10^4$ m , but the accident mode in which the gas consumption and locations of the ring gas pipelines allocated to the consumers remain constant ($x_1=0.3\times 10^4$ m; $x_2=1.5\times 10^4$ m: $x_3=2.7\times 10^4$ m;). We will look at a situation where the volume of gas mass leaked into the environment is equal to 50% of its total consumption ($G_{le}=0.5g$). Thus, we assume that $G_1=3$ Pa $\times$ sec/m; $G_2=4$ Pa$\times$ sec/m; $G_3=3$ Pa$\times$ sec/m and $G_{le}$ = 5 Pa$\times$ sec/m gas. We will calculate the pressure values in the pipeline in time and length under accident conditions and enter these values into the table below to compare them with the values in Table 2.3.1 (Table 2.3.2).

| x, m | P(x,t), $10^4$ Pa | | | | | |
|---|---|---|---|---|---|---|
| | t=50 sec | | t=300 sec | | t=900 sec | |
| | $\sum_{i=1}^{3} Gi$ | $\sum_{i=1}^{3} G_i + G_{le}$ | $\sum_{i=1}^{3} Gi$ | $\sum_{i=1}^{3} G_i + G_{le}$ | $\sum_{i=1}^{3} Gi$ | $\sum_{i=1}^{3} G_i + G_{le}$ |
| | Pa× sec/m | | Pa× sec/m | | Pa× sec/m | |
| 0 | 12.3 | 11.4 | 11.0 | 9.6 | 8.1 | 5.2 |
| 1000 | 14.6 | 15.0 | 13.5 | 13.2 | 10.5 | 8.8 |
| 3000 | 18.6 | 21.0 | 17.6 | 19.5 | 14.7 | 15.0 |
| 6000 | 23.0 | 27.5 | 22.1 | 26.2 | 19.2 | 21.8 |
| 9000 | 25.5 | 31.3 | 24.8 | 30.1 | 21.8 | 25.7 |
| 12000 | 26.0 | 32.0 | 25.4 | 31.0 | 22.4 | 26.6 |
| 15000 | 26.4 | 32.6 | 25.8 | 31.7 | 22.8 | 27.2 |
| 17000 | 25.1 | 30.7 | 24.5 | 29.7 | 21.5 | 25.3 |
| 20000 | 23.0 | 27.6 | 22.3 | 26.5 | 19.4 | 22.1 |
| 23000 | 20.3 | 23.4 | 19.5 | 22.2 | 16.5 | 17.8 |
| 27000 | 15.9 | 16.8 | 14.8 | 15.2 | 11.9 | 10.82 |

**Table 2.3.2. Comparison of pressure values depending on the length and time in circular pipelines in technological and emergency conditions.**

Using the values in the breakdown condition $\left(\sum_{i=1}^{3} G_i + G_{le}\right)$ from Table 2.3.2, we construct a table of pressure distribution over time along the length of the pipeline (fig. 2.3.5).

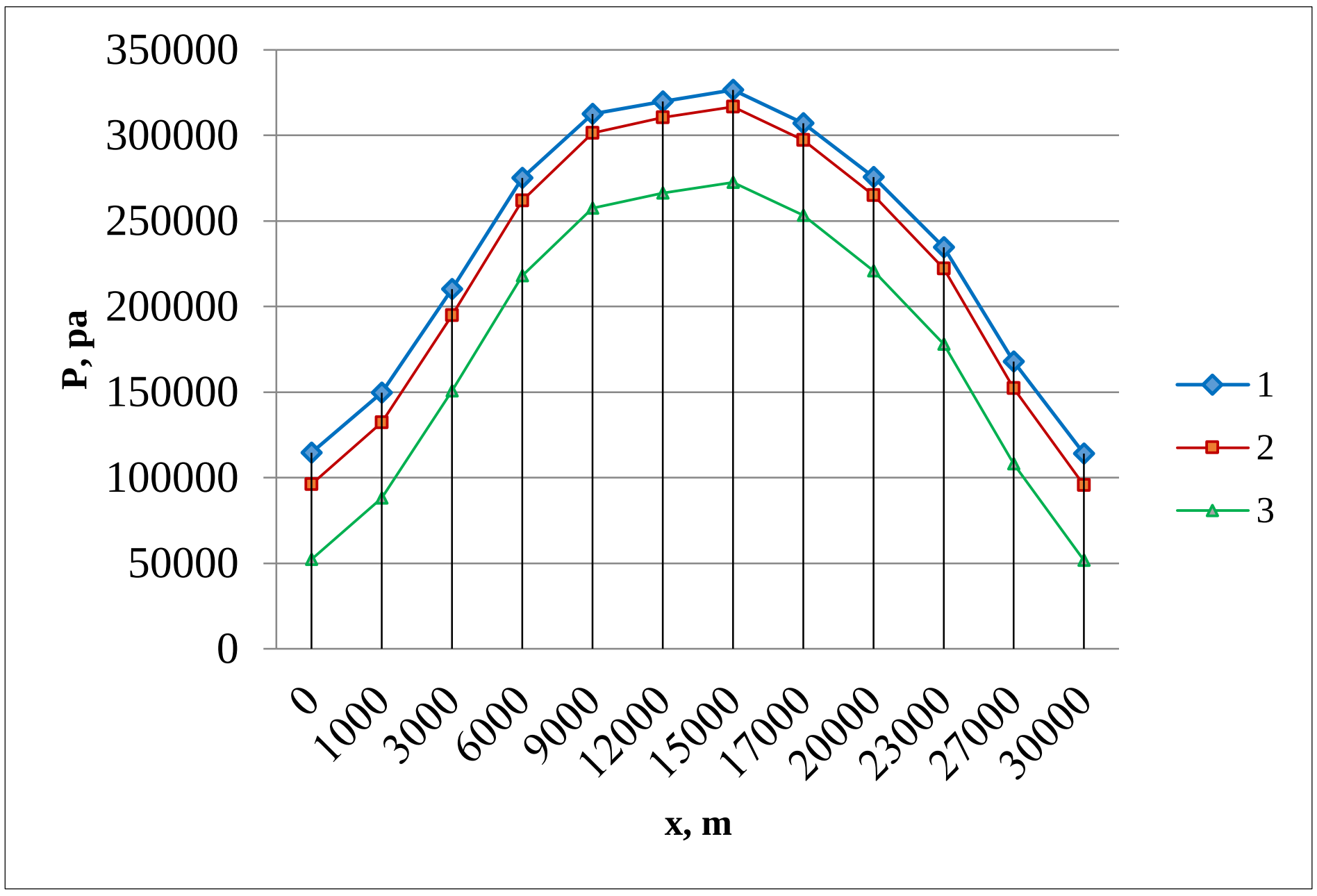

**Figure 2.3.5. Distribution of its pressure depending on the time along the length of the annular gas pipeline under emergency conditions. (1-t=50s, 2 - t=300s, 3-t=900s).**

Analysis of Tables 2.3.2 and Figure 2.3.5 shows that emergency situations have a significant impact on pressure values in closed circuits. Especially if the pressure at the beginning of the gas pipeline drops from $14 \times 10^4$ Pa to $12.3 \times 10^4$ Pa after 50 seconds in normal technological mode. As a result of the violation of the integrity of the pipeline (in emergency cases), the starting pressure of the gas pipeline decreases from $14 \times 10^4$ Pa to $11.4 \times 10^4$ Pa. So, in the event of an accident, the price of the pressure drop in the starting section of the gas pipeline is $0.9 \times 10^4$ Pa more than in the normal technological mode. After 300 seconds, this difference increases to 1.4 $\times$ $10^4$ Pa, after 900 seconds to $2.9\times10^4$ Pa. On the other hand, as shown in Figure 2.3.5, as the pipeline moves from its initial section to its middle section, the gas flow pressure increases and reaches its maximum value in the middle section. After 50 seconds, the maximum pressure value increases to $26\times10^4$ Pa, and in emergency cases, the maximum pressure value increases to $32\times10^4$ pa. Thus, in the event of an accident, the increase in maximum pressure due to gas leakage exceeds its value in normal technological mode by 6 $\times10^4$ Pa. After 300 seconds, the maximum pressure increases by 5.6 x $10^4$ Pa, after 900 seconds this increase is 4.2 $\times$ $10^4$ Pa. So, in emergency situations, the pressure difference in the initial section of the annular pipeline increases over time, and the difference in the maximum pressure values at the midpoint of the pipeline decreases (fig. 2.3.4 and 2.3.5). Thus, both in the usual technological mode and in emergency conditions, the pressure of the gas flow reaches its maximum value at the point of hydraulic stability of the annular pipeline ( $x=\frac{L}{2}$ ). And this is a common feature of the dynamic process of ring networks. This is due to the fact that during the operation of the gas pipeline in both technological and emergency modes, there is a reverse flow and a decrease in the pressure Reserve in the starting Section. It is this feature that must be taken into account when managing the operation of circular gas pipelines.

The gasodynamic processes mentioned in the emergency conditions of gas pipelines of other configurations have not been observed. For example, studies conducted under conditions of violation of the integrity of a linear and parallel gas pipeline say that the quality of pressure changes varies significantly depending on the location of the leakage point [85,87]. This variability is most pronounced in the pressures at the initial and final points of the gas pipeline. So, if an accident occurs at the beginning of the pipeline , then the pressure drop at the beginning of the pipeline will be approximately 600 seconds $8{,}62\times10^4$ Pa. But the pressure drop at the end of the pipeline is smaller, about $0.28\times10^4$ Pa, and it differs from the initial pressure drop by about 30 times. If an accident occurs at the end of the pipeline section, then after 600 seconds at the start, the pressure drop will be $0{,}16\times10^4$ Pa. However, the pressure drop at the end of the pipeline will be significantly larger, about $8.62\times10^4$ Pa, which is about 53 times the initial pressure drop. The result of the analysis indicates that determining the location of leakage in pipelines analytically, based on the ratio of pressure drops occurring at the beginning and end of the pipeline, is more suitable for the intended purpose.

Based on the analysis carried out, it can be noted that it is impossible to determine the place of leakage in the annular gas pipeline by analytical method using Equation (7) based on existing reporting methods. Because, as can be seen from Figure 2.3.5, the prices of pressure drops in the initial and final sections of the annular gas pipeline in the event of an accident are always equal to each other, regardless of the place of leakage. As a result of the research, it can be said that one of the distinctive features of the gasodynamic process of the circular gas pipeline is that the change in the values of gas flow parameters (pressue) depending on the time in unstable conditions does not depend on the number and location of pipes connected to the network. This is due to the fact that the pressure is a continuous function along the entire length of the ring network (7), and the gas dynamic calculations of the ring are based on Kirchhoff's law at all nodes of the network.From the analysis of exerimental reports, on the other hand, it turned out that under emergency conditions it receives the maximum value of the pressure of the pipeline along its length exactly in the middle part of the ring (Figure 2.3.4, Figure 2.3.5).

Thus, we need to develop a reliable method through gas-dynamic analysis of the circular gas pipeline and processing of the obtained data. In this case, we have obtained an algorithm for making a valid decision in accordance with the results of an experimental study of mathematical expressions and thus for solving conflicting results. The data of the algorithm can be considered as a new database for detecting and localizing emergency situations in the annular pipeline.

$$\frac{P(\frac{L}{2},t)-P_1}{P(0,t)} = \kappa(t) \quad (8)$$

It has been established that if $\kappa(t) \succ 1$ meets the conditions of inequality, then the change in pressure depending on the time characterizes the accident conditions in the gas pipeline, and if $\kappa(t) \prec 1$, this is the result of technological processes. Thus, the function $\kappa(t)$ is a reliable parameter that allows you to make the right decision to efficiently manage the annular gas pipeline.

The application of pressure sensors in gas pipelines to find and analyze the causes of non-stationary pressure in the pipeline is considered. The experimental results show that the middle part of the annular pipelines is a characteristic part for managing the processes of emergency situations for gas transport. In this regard, in order to ensure the reliability of monitoring systems in the process of managing ring gas pipeline systems, it is recommended to install two pressure sensors in the initial and middle parts of the pipeline (Figure 2.3.6). Thus, a highly reliable and effective technological model of managing accident processes using the Internet of Things (IoT) platform is presented [82].

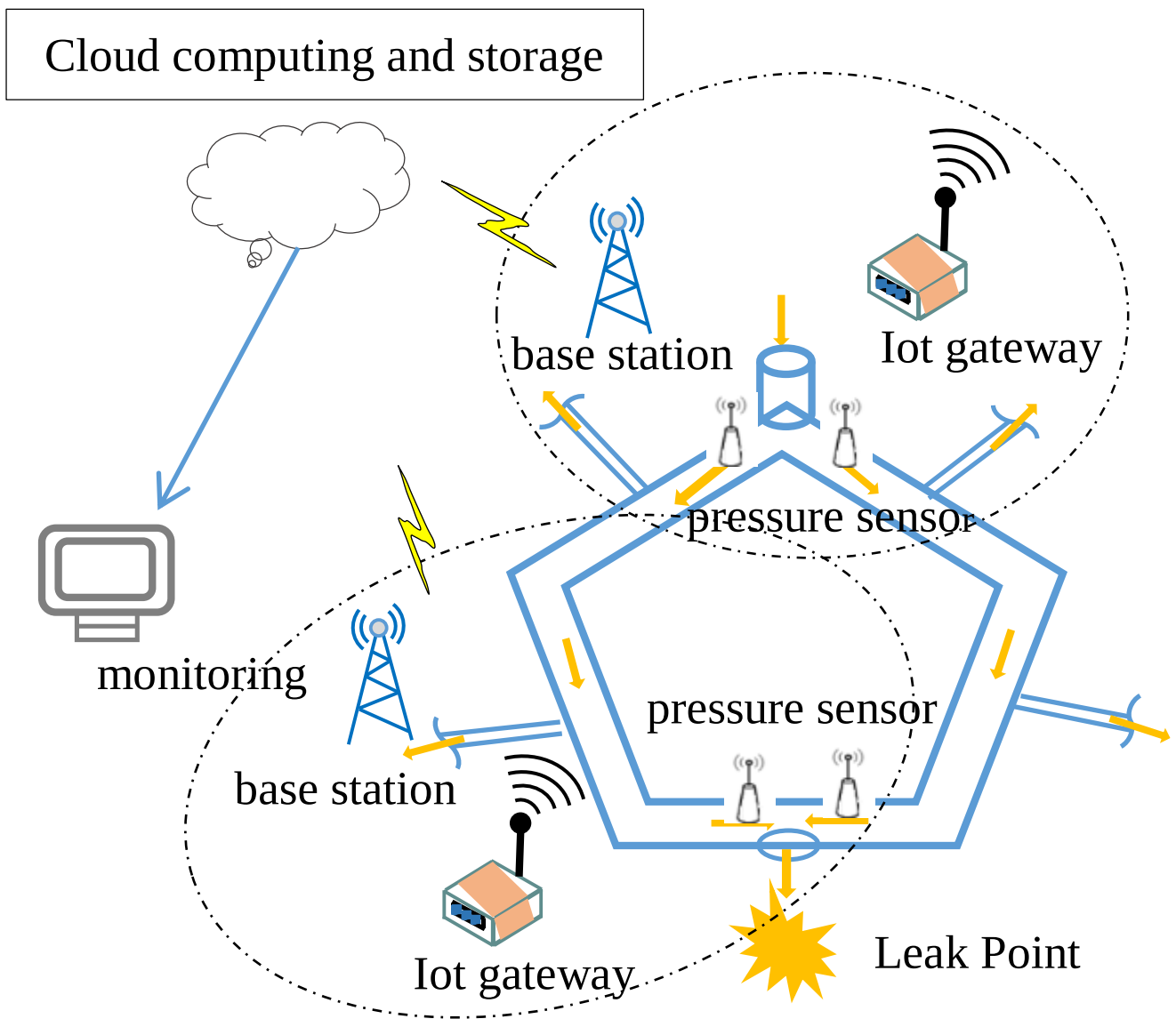

**Figure 2.3.6. Scheme of sensor networks for monitoring annular pipelines.**

The values of the function  are automatically calculated sequentially in the control centers and it is determined whether the value obtained is greater or less than one unit.

To solve it, computer processing of the data (pressure ) from the sensors installed in the gas pipeline is necessary, as can be seen from Figure 2.3.6. Data from the wireless pressure sensor is supplied to the Iot - gateway operating in master Mode. Thus, the information in its network enters the data processing center through the cloud computing and data storage network.

Finally, the Information Center makes a decision on whether there is an emergency in the gas pipeline.

To ensure comprehensive control of the integrity of the system, it is important that the data center accurately detects whether pressure fluctuations in the starting and middle sections of the gas pipeline are caused by accidents or technological processes and is able to distinguish the relationships between these processes. In this context, modern control methods, particularly those enhanced by artificial intelligence (AI) technologies, should be prioritized. For instance, AI technologies, particularly hybrid models, provide a more resilient and effective approach by offering accurate predictions even with limited data. Specifically, indicators that assess algorithms' capabilities in real-time incident detection and data augmentation for handling data insufficiencies have been presented. This approach lays a stronger foundation for addressing existing challenges in pipeline engineering.

Furthermore, key tasks include tracking the evolution of incidents, using data augmentation to address data gaps, enhancing model interpretability, and improving the management of pipeline integrity. Future research should focus on the optimization and integration of AI models to support informed decision-making based on a reliable information base. Additionally, AI enables the analysis of large volumes of data from various sources to predict incident risks, optimize management processes, and provide critical insights that support engineering decisions. Consequently, prioritizing the broader application of hybrid models in real-time monitoring and predictive technologies is essential for advancing this field.

## CONCLUSIONS

The solution to the differential equations system describing the non-stationary gas flow in a ring pipeline with a constant cross-sectional area has been obtained. The requirements for the equality of boundary conditions included in the solution allow for a compact representation of the analytical model but do not restrict the model's applicability to sharp changes in gas flow or pressure. Thus, the methodology for calculating non-stationary gas flow in ring pipelines has been developed.

To ensure the unconditional resolution of efficient management issues, a function characterized by the ratio of pressure drops has been defined and tested. An analysis of the $\kappa(t)$ function and parameters, which can perform several operations in a fixed time, showed that accounting for this indicator at the operational stage can optimize the minimization of information transfer time to central control points. The established connection allows you to distinguish between technological and emergency modes of operation of gas pipelines based on indicators of changes in the parameters of the final and middle sections of the gas flow.

A potential approach to increasing the efficiency, accuracy, and scale of the management of ring-shaped tube networks has been defined as the integration of IoT and cloud technologies in emergency detection systems. When a leak occurs in the gas pipeline, data analysis and pressure sensor data are made possible by the combination of these technologies. Therefore, future research in this field should give priority to the creation of plans that can satisfactorily solve these issues and guarantee the effectiveness and safety of these systems. Thus, it is recommended to introduce a system using IoT and cloud technologies for detecting emergency situations in circular gas networks.

## 2.4. Technological Foundations of Parallel Gas Pipeline Reconstruction

The improvement of the technical and economic indicators of the gas industry is mainly related to two areas. Improving the equipment and technology of gas transmission and distribution networks, aimed primarily at reducing the capital intensity and operating costs of facilities. And then improving the dispatching control of gas supply, optimizing operating modes and increasing the reliability of gas supply and saving energy resources spent on their own needs.

One of the solutions to the problem of improving the reliability of gas pipelines is the use of new effective scientifically based technologies for the reconstruction of pipeline systems. Connecting fittings, loopings and connectors are one of the components of increasing the reliability of gas networks. Loops are designed both to increase pressure at the end point of the gas pipeline or to increase throughput, and to reduce pressure at the starting point. The purpose of the calculation is to determine the length of the looping, at which the necessary effect will be provided. A luping gas pipeline is one of the simplest types of complex gas pipelines.

Obviously, as the number of lupings increases, the effectiveness of binders decreases. In accordance with the decrease in the diameters of the points, a decrease in the efficiency of the connections is also evident. If the diameters of the parallel spans do not change in length and there is no luping, then using couplings as a means of increasing the capacity of the gas pipeline is useless. But the connectors allow you to reduce throughput during scheduled pipeline repairs. In other words, it is necessary to adopt such a scheme of parallel gas pipelines so that in the event of an accident, the connectors perform the luping function.

According to the regulations, parallel gas pipelines are connected to each other through interconnectors and operate in a unified hydraulic mode (Figure 2.4.1). Parallel gas pipelines are widely used in practice due to their high reliability compared to linear pipelines and their ability to provide consumers with an uninterrupted energy carrier. When an accident occurs in one of the levels or repair work is carried out, the damaged section is isolated from the main gas pipeline using automatic valves. Subsequently, the valves on the interconnectors near the accident site are opened. As a result, consumers are not deprived of gas supply (Figure 2.4.2). In other words, the gas from the damaged level is transferred to the undamaged one through interconnectors, and this process continues until the damaged pipeline is repaired [1]. From this perspective, the construction and reconstruction of gas pipelines are aimed at creating a gas supply system that is independent of the pipeline's technological operating mode and various repair works.

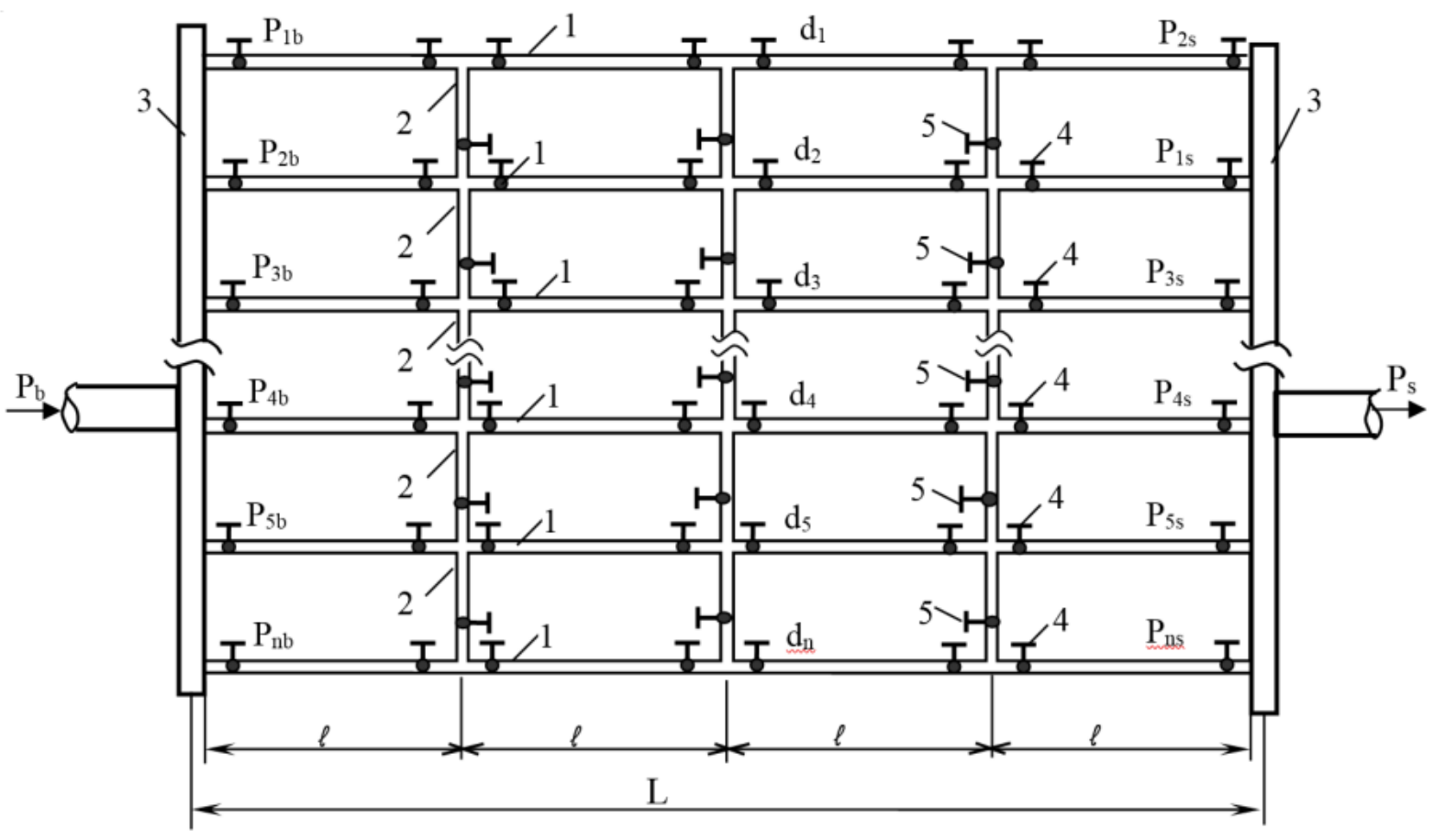

**Figure 2.4.1. Schematic diagram of a multiline parallel gas pipeline.**

1 - Pipeline sections on the first line; 2 - Interconnectors; 3 - Collectors; 4 - Automatic valves on pipeline sections; 5 - Automatic valves on interconnectors; L - Length of the gas pipeline;  - Distance between interconnectors; $P_i$ - Initial pressure of the pipeline; $P_f$ - Final pressure of the pipeline; $d_1$, $d_2$, $d_3$, ..., dn - Diameters of the pipeline sections on each line.

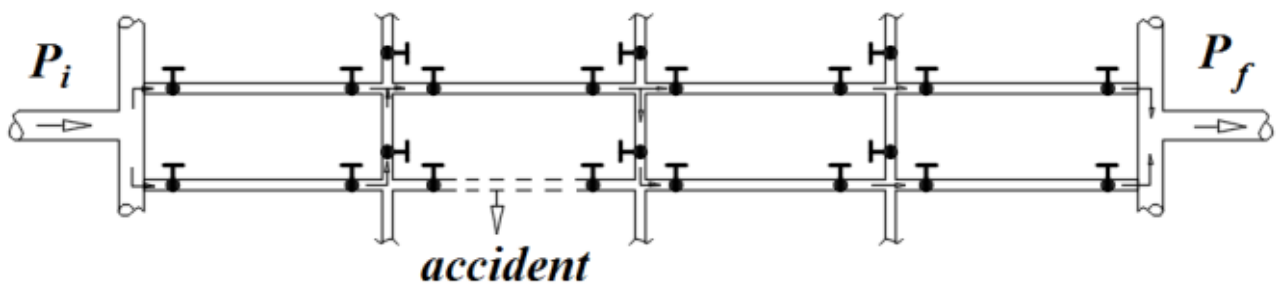

**Figure 2.4.2. Schematic of interconnecting multiline gas pipelines in emergency mode.**

During the reconstruction of parallel gas pipelines, it is necessary to select an optimal distance between interconnectors to minimize the cumulative costs while generating an alternative income from providing consumers with uninterrupted gas supply. As previously mentioned, depending on the different technological operating modes of the pipeline and the need to enhance the capacity of the system and each line, interconnectors are put into operation. It is clear that only through this process can consumers be continuously supplied with gas. This process is carried out as follows:

1.The damaged or to-be-repaired section of a parallel pipeline is identified.

2.The connecting fittings previously placed on the pipeline on either side of the identified section are activated (Figure 2.4.2). At this point, the damaged section of the pipeline is separated from the main system and prepared for repair.

3.To ensure the uninterrupted flow of gas in the damaged section, the gas valves on the connectors located on the right and left sides of the damaged section are activated (opened) [1]. This allows the gas flow to be directed to the undamaged section. As a result, consumers' gas supply remains uninterrupted (Figure 2.4.2).

Assuming that the number of line of the parallel gas pipeline shown in Figure 2.4.1 is three (n = 3), and if the last, third line is damaged, the implementation of the transition processes mentioned above will result in a three-level scheme for the pipeline. As seen in Figure 2.4.2, after activating the connectors in the damaged section, it acts as a bypass (looping). It is known that looping is used to increase the capacity (efficiency) of gas pipelines. Research has shown that the productivity of the pipeline changes depending on the length of the installed looping. In other words, as the length of the looping increases, the productivity of the pipeline will increase accordingly [33]. On the other hand, increasing the length of the installed looping will increase the capital investment in the system.

From this perspective, it is essential to determine the optimal parameters of the installed looping (i.e., the diameter and length) to ensure that the capacity of the multiline parallel gas pipeline increases significantly while minimizing the cost of installing the looping. As shown in Figure 2.4.1, the length of the looping depends on the distance between the steps of the connectors on the parallel gas pipeline. Therefore, working to determine the optimal distance between interconnectors is one of the fundamental methods for the reconstruction of multiline parallel gas pipelines. Therefore, determining the optimal length of the looping between interconnectors during the reconstruction of parallel gas pipelines is crucial, as the income obtained as a result of increasing the pipeline's productivity allows for the cost of the system's reconstruction to be covered.

We can represent this relationship in the following analytical form:

$$Z = \phi(\ell). \tag{1}$$

To determine this function, we must first define the productivity (throughput capacity) of the multiline parallel line gas pipelines as a function of the distances between interconnectors.

Our goal is to analyze the efficiency of gas supply to consumers based on distances between interconnectors, gas flow modes, and repair processes. First, we establish the existing productivity of the multiline parallel gas pipeline before any damage occurs. Using Figure 2.4.1, we describe the steady-state behavior of gas flow for each line of the gas pipeline:

For the 1st line: . $$P_{1i}^2 - P_{1f}^2 = c^2 \frac{\lambda L}{d} Q_1^2 \quad (2)$$

For the 2nd line: $$P_{2i}^2 - P_{2f}^2 = c^2 \frac{\lambda L}{d} Q_2^2. \quad (3)$$

For the 3rd line: $$P_{3i}^2 - P_{3f}^2 = c^2 \frac{\lambda L}{d} Q_3^2. \quad (4)$$

And so on, with a similar equation applicable to the n-th line:

$$P_{ni}^2 - P_{nf}^2 = c^2 \frac{\lambda L}{d} Q_n^2. \quad (5)$$

These equations assume that the pipeline consists of pipes with the same diameter. Here, the variables λ and d represent the hydraulic resistance factor and diameter of the pipes in the existing gas pipeline lines, respectively. Variable c- is the speed of sound in the gas [85], and

$c = \sqrt{zRT}.$

In this equation,

$R$ - represents the universal gas constant;

$T$ - denotes the absolute temperature;

$Z$ - represents the gas compressibility factor.

$P_{1i}, P_{2i}, P_{3i}, \ldots, P_{ni}$ - are the initial pressures at the respective lines of the parallel gas pipeline.

$P_{1f}, P_{2f}, P_{3f}, \ldots, P_{nf}$ - are the pressures at the end of the lines of the parallel gas pipeline.

Since these multiline parallel gas pipelines operate in hydraulic mode, their pressures are equal:

$P_{1i} = P_{2i} = P_{3i} = \ldots = P_{ni} = P_i$ ,

$P_{1f} = P_{2f} = P_{3f} = \ldots = P_{nf} = P_f$

Where $P_i$ and $P_f$ represent the initial and final pressures of the multiline parallel gas pipeline, as shown in Figure 2.4.1.

To determine the total productivity of the multiline parallel gas pipeline, we sum equations (2), (3), (4), and (5). We get:

$$\Pi\left(P_i^2 - P_f^2\right) = c^2 L\left(\frac{Q_1^2 + Q_2^2 + Q_3^2 + \ldots + Q_n^2}{d}\right)\lambda. \quad (6)$$

It is important to note that for a steady-state mode, the sum of losses in the lines is equal to the total pipeline losses ($Q_0$). So:

$$Q_1^2 + Q_2^2 + Q_3^2 + \ldots + Q_n^2 = Q_0^2. \quad (7)$$

Now, using equation (7) in equation (6), we can calculate the required flow rate as follows:

$$Q_0 = \sqrt{\frac{\left(P_i^2 - P_f^2\right) n \cdot d}{c^2 \lambda}} \cdot \frac{1}{\sqrt{L}}. \quad (8)$$

Next, when one line of the parallel gas pipeline is damaged, that is, when gas is redirected from the damaged section to the undamaged section using connectors, we can determine the total productivity (efficiency) of the pipeline for the transition process.

For research purposes, we divide the damaged section into three parts (as shown in Figure 2.4.2):

- The distance from the start of the pipeline to the first connector (looping system) $0 \le x \le \ell$;
- The damaged section between two connectors $x \le \ell \le (L - x - \ell)$;
- The distance from the last connector to the end of the $\ell \le (L - x - \ell) \le L$.

We then define the following variables:

$\tilde{Q}$ - the changed productivity of the pipeline due to gas flow redirection from the damaged section to the undamaged section (looping system).

$\ell$ - the unknown distance between connectors.

$P_1$ and $P_2$ - the pressures at the start and end of the damaged section. Using these variables, we can write the equations for the steady-state gas flow in each of these sections:

For the looping system $0 \le x \le \ell$:

$$P_i^2 - P_1^2 = \lambda c^2 \frac{Q_\ell^2}{2d} x. \tag{9}$$

For the damaged section $x \le \ell \le (L - x - \ell)$:

$$P_1^2 - P_2^2 = \lambda c^2 \frac{Q_\ell^2}{d} \ell. \tag{10}$$

For the looping system at the end $\ell \le (L - x - \ell) \le L$:

$$P_2^2 - P_f^2 = \lambda c^2 \frac{Q_\ell^2}{2d} (L - x - \ell). \tag{11}$$

Combining equations (9), (10), and (11), we can determine the efficiency (consumption) of the second and third lines (branch system) of a three-line parallel chamber.

For the second and third lines of the three-line parallel chamber, based on the power formula:

$$Q_\ell = \sqrt{\frac{(P_i^2 - P_f^2)}{\lambda c^2}} \cdot \frac{\sqrt{2}}{\sqrt{\ell + L}}. \tag{12}$$

For the first line of the three-line parallel chamber (n = 1), according to equation (8), we have:

$$Q_\ell = \sqrt{\frac{(P_i^2 - P_f^2)}{\lambda c^2}} \cdot \frac{1}{\sqrt{L}}. \tag{13}$$

It is clear that, for the transition process of a three-line parallel chamber, we can calculate the total efficiency (consumption) by summing equations (12) and (13):

$$\tilde{Q} = \sqrt{\frac{(P_i^2 - P_f^2)}{\lambda c^2}} \cdot \left[ \frac{\sqrt{2}}{\ell + L} + \frac{1}{\sqrt{L}} \right]. \tag{14}$$

Here, $\tilde{Q} = Q_\ell + Q_1$.

So, when the connectors are located at the beginning and end of the gas pipeline, the efficiency increases. However, for the other sections of the pipeline, the length $\ell$ between the connectors is essential. Therefore, we need to determine the optimal length of $\ell$ so that the productivity remains close to its original line before any damage to the gas pipeline (until the repair period). To find the required expression for point (1), we assume n=3 in expression (8) (for a three-line parallel pipeline) and determine the ratios of expression (8) with formula (14) as follows:

$$\phi(\ell)=\frac{\tilde{Q}}{Q_0}=\left[\frac{\sqrt{\frac{2L}{\ell+L}}+1}{\sqrt{3}}\right]. \tag{15}$$

From equation (15), it is clear that to find the optimal distance between the connectors in a multiline parallel gas pipeline, you need to find the value of l that makes $\phi(\ell)$ approach unity $[\phi(\ell)\rightarrow 1]$ . However, under this condition, it is possible to provide consumers with a continuous supply of gas, regardless of damage to the pipeline and its repair. In other words, the majority of consumers can use gas normally, and the system's reliability will be high.

On the other hand, as previously mentioned, it is necessary to determine the value of $\ell$ in such a way that the cost of the system's reconstruction is not significantly higher than the revenue obtained from the uninterrupted supply of gas to consumers. In other words:

$$\phi(\ell)\cdot S_g \geq S_{rec}. \tag{16}$$

If we replace the "≥" symbol with "=" in equation (16), we get the following equation:

$$\frac{2L}{\ell+L}+2\sqrt{\frac{2L}{\ell+L}}+1-\frac{\sqrt{3}\cdot S_{rec}}{S_g}=0. \tag{17}$$

Here, $S_{rec}$ -is the total incurred cost of the reconstruction of the multi-layer parallel gas pipeline, and

$S_g$ -is the profit obtained as a result of optimal placement of the connectors.

By solving equation (17), you can determine the optimal length of the step between connectors in the reconstruction of multiline paralel gas pipelines.

$$\ell=L\left[1-\frac{2}{\left(\sqrt{\beta}-1\right)^2}\right]. \tag{18}$$

Where, $\beta=\frac{\sqrt{3}\cdot S_{rec}}{S_g}$.

However, determining the optimal length of connectors leads to the need to consider many reconstruction variants, as the value of β ($\beta \in ]0:\infty[$ ) can be infinite. As a result, equation (18) can take an infinite number of values. Thus, the number of reconstruction variants according to equation (18) will be sufficiently large.

To avoid the complexity mentioned above, the issue of determining the lengths of connectors based on the complex nature of the system can be considered. As shown in Figure 2.4.1, if we assume the length between connectors as $\ell$, then the number of steps (number of connectors) m will be $m=\frac{L}{\ell}$ . On the other hand, it is clear that for one step, there will be approximately 8 automatic valves (connector fittings) and 3 connectors installed for a three-layer parallel gas pipeline. Thus, when the number of connectors is m, the total cost of the incurred expenses for the reconstruction of gas pipelines can be determined as follows [44, 79]:

$$S_{rec} = m\cdot Z \tag{19}$$

here, Z -is the value of expenses for each step of pipeline reconstruction. Z is given by:

$$Z = E_n\left(8K_{a.v} + 3K_{con}\right) + 8C_{a.v} + 3C_{con}. \qquad (20)$$

Where $K_{a.v}$, $K_{con}$ -are the capital costs for the installation of an automatic valve and connector during reconstruction.

$C_{a.v}$, $C_{con}$ – are operating expenses for automatic valves and connectors.

$E_n$- is the normative factor for the comparative efficiency of capital investment for gas pipelines.

Thus, to determine the reconstruction variant, we consider equation (19) in inequality (16). Then we have:

$\phi(\ell)S_g \geq mZ$ or

$$\left[\frac{\sqrt{\frac{2L}{\ell+L}}+1}{3}\right]\cdot S_g \geq \frac{L}{\ell}Z \qquad (21)$$

The left side of equation (21) in the inequality represents the profit obtained by consumers during the reconstruction of the connectors, and the right side represents the expenses incurred in the reconstruction of the system. It is clear that when the condition of inequality (21) is met, the length of the specified connector ensures the economically efficient reconstruction of the system [46, 79].

When one of the layers of the parallel gas pipeline is damaged, or during the operation of the connectors for different reasons, the profit obtained for each 1 $m^3$ of gas supplied to consumers can be determined analytically as follows:

$$S_g = \frac{Q_0 \cdot t \cdot e}{n} \cdot \omega \cdot L. \qquad (22)$$

Where $Q_0$ - is the productivity of multi-layer gas pipelines, $m^3$/hour,

n -is the number of layers in the gas pipeline,

t -is the period of repairing the gas pipelines or the period during which consumers are supplied with gas without interruption,

ω -is the gas pipeline's accident rate, *1/km· year,*

L- is the length of the gas pipeline, *km,*

e -is the income obtained for each 1 $m^3$ of gas supplied to consumers with different categories.

Thus, the analysis of inequality (21) shows that during the reconstruction of parallel gas pipelines, both the profit obtained and the cost of the connectors are dependent on the distance between the connectors. Then, if we assume that the left and right sides of inequality (21) are equal, we get:

$$\left(\sqrt{\frac{2L}{3\ell+3L}}+\frac{1}{\sqrt{3}}\right)\cdot \wp = \frac{L}{\ell}. \qquad (23)$$

Here, $\wp$ is defined as follows:

$$\wp = \frac{S_g}{Z}.$$

From equation (23), we obtain the following equation.

$$Wl^3 + \left( \frac{2\sqrt{3}L}{\wp} - L \right) \cdot \ell^2 + \left( \frac{2\sqrt{3}L}{\wp} + \frac{3L^2}{\wp^2} \right) . \ell + \frac{3L^3}{\wp^2} = 0 \qquad (24)$$

Solving equation (24) allows us to determine the unknown length of the step between connectors for the economically efficient reconstruction of parallel gas pipelines.

$$\ell = \frac{2L}{\sqrt{3}} \left( \mu sh \frac{4}{3} - \frac{1}{\wp} + \frac{1}{2\sqrt{3}} \right). \qquad (25)$$

Here,

$$\mu = \sqrt{\left( \frac{1}{\wp} - \frac{5}{\sqrt{3}} \right)^2 - 8};$$

$$\eta = \frac{9\sqrt{3}}{\wp} \left( \frac{3}{\wp^2} + 1 \right) - \left( \frac{2\sqrt{3}}{\wp} - 1 \right)^3; \quad \xi = \frac{\eta}{3\sqrt{3}\mu};$$

$$\phi = \ln\left[ \xi + \sqrt{1 + \xi^2} \right].$$

Thus, for engineering calculations in the reconstruction of parallel gas pipelines, you can confidently use the formula (25). To do this, we accept the following initial data of a three-line parallel gas pipeline:

-the initial pressure of the gas pipeline is $P_i$ = 2 MPa;
-the final pressure of the gas pipeline is $P_f$ = 0.85 MPa;
-the length of the gas pipeline is L = 40 km;
-the diameter of the gas pipeline threads is D = 0.5 m;
-the speed of sound sliding in the gas is c = 383.3 m/sec;
-the number of lines is n = 3.

For gas pipelines, the value of the fracture intensity $\lambda = 2.5 \cdot 10^3$ Mpa 1/km·year is assumed, and the value of the coefficient of friction $\lambda = 0.03$.

To begin with, let's determine the output of the parallel gas pipeline using the formula (8).

$Q_0 = 0{,}17 \cdot 10^{-3}$ *MPa·sec/m* = 161028 *$m^3$/h.*

The amount of damage for each 1 $m^3$ of gas that is not presented to various categories of consumers, based on statistical data studied by Professor O.M.Ivanchov, we take e = 1.2 €/$m^3$.

The repair period of gas pipelines according to building standards T = 6÷24 hours, based on which we assume T = 6 hours. Then according to formula (22),

$$S_g = 618348 \ \frac{€}{year}$$

it will be.

To determine the value (Z) of the total reduced costs incurred in the system during the reconstruction of the gas pipeline, we take the following initial data.

On the website of the agrochemiyainvest society on the internet (http://zadvizhki.narod.ru/) it is indicated that the cost of an automatic crane (electric drive) $K_{a.k}$ = 1500 €. Accordingly, the operating costs of the crane will amount to $C_{a.k}$ = 145.5 €/year. If the cost of installing one contactor is $K_{con}$ = 103 €, then its operating costs will amount to 10 €/year.First, using expression (23), we find the value of Z (we take $E_n$ = 0.12).

Z = 2659 €/year

So, $\wp = 232{,}55$; $\mu = 0{,}555$; $\eta = 1{,}023$; $\xi = 0{,}355$ taking into account the calculated prices of the coefficients in formula (25),

$$\ell = 0{,}39 \cdot L = 0{,}39 \cdot 40 = 15{,}6 \ km$$

Thus, for this particular example, the most convenient length for the reconstruction option is the placement of $\ell = 15{,}6 \ km$ distance between connectors that will be installed on a three-circuit gas pipeline with a length of $L = 40 \ km$. But if the parameters of the gas flow in the pipeline differ from the example considered, then the length of the lupin that will be installed for the gas pipeline at this length may change.

In other words, the length between the connectors may vary depending on the geometric and technical parameters of the parallel lines location, gas flow parameters, as well as the category of consumers. For this reason, determining the optimal distance between connectors is of great importance in the method of reconstruction of multicomponent gas pipelines.

## Conclusion

A method for the reconstruction of complex gas transportation networks has been developed with the aim of conserving energy resources and ensuring the reliability of gas supply to consumers. The advantage of this method is characterized by the practical significance of technological calculations. Timely implementation of measures allows for an increase in the pipeline's productivity and a reduction in the energy costs for gas transportation. Based on the results of the research, a function has been developed that characterizes the reliable operation of parallel gas pipeline systems in stationary and emergency modes. Using this function enables the determination and enhancement of the gas transmission capacity of the looping system that arises due to the operation of connectors during pipeline emergencies.

By employing the theoretically and technically grounded formula, it is possible to efficiently carry out the reconstruction of existing multi-line parallel gas pipelines based on economic and technical principles. The implementation of the proposed calculation scheme for the economically efficient placement of connectors in existing and newly commissioned parallel gas pipelines will result in more favorable economic indicators of the reconstruction compared to other alternatives.

### 2.5. Technological Schemes and Control Methods in the Reconstruction of Parallel Gas Pipeline Systems Under Non-Stationary Conditions

One of the promising directions for the modern development of the fuel and energy sector is the technical and economic justification of the reconstruction methods for existing gas pipeline systems to ensure continuous energy supply to consumers. This is directly related to the continuous provision of energy resources to new consumers emerging as a result of the development of the industrial infrastructure in the country, both in the near and distant future. Parallel gas supply systems are complex structures consisting of numerous elements. Due to the physical wear and aging of these elements, existing gas pipeline systems no longer meet the requirements of modern technology. As a result, the reduction in the operational reliability of these pipelines hinders the future development of their management systems. To this end, the modernization of elements that have undergone physical and moral wear is required to increase the efficiency of gas pipeline system operation. On the other hand, the issue of gas supply in the required volumes to new consumers must also be addressed. One of the methods for the reconstruction of parallel gas pipeline systems is based on these principles[87].

To address the shortcomings in the operation of parallel gas pipelines in a timely manner using known methods of reconstruction, the following measures are planned: increasing the gas supply volume in the network by expanding the diameter of the gas pipeline; adjusting the distance between connecting pipes installed to transfer gas flow from one pipeline configuration to another in order to provide the required volume of gas to consumers supplied by parallel gas pipelines during emergency situations; installing additional new connecting pipe networks when needed; improving the schemes of parallel gas pipeline systems during emergency or planned maintenance to prevent significant gas losses; and defining modern control methods, among others.

One of the important issues in improving the management systems of non-stationary gas flow is the replacement of outdated elements with more modern components in order to increase the operational reliability of the parallel gas pipeline system. When modeling the operation of gas distribution networks, it is essential to place these elements and separation devices (valves) in parallel gas pipeline networks in such a way that, during emergency situations, the damage to the environment caused by gas leaks and the reduction in gas supply consumption are minimized, thereby minimizing the damage to consumers. It should also be noted that, for the technical and economic efficiency of the reconstruction, the installation of separation valves and elements, which are planned to be added, must take into account the construction, installation, and operational costs when integrating them into the structure of the parallel pipeline network.

In the context of reconstruction, the rational placement of separation devices (valves) and elements requires determining the key parameters of the gas distribution systems, along with identifying the categories of consumer enterprises [123]. At this stage, in selecting the management method, it is necessary to solve the problem of choosing the most optimal and rational approach for ensuring the reliable and safe operation of the parallel gas network [104].

The operation and reconstruction of parallel gas pipeline systems involve making management decisions that create various technical and economic challenges. In many cases, it is impossible to resolve these issues without implementing an automated control system. To address this, the present study proposes methods and mathematical models for monitoring and managing the non-stationary operating modes of the linear section of the parallel gas pipeline. These models, based on pressure functions, enable the analysis and calculation of the dynamic characteristics of technological processes resulting from the disruption of the integrity of any section of the parallel gas pipeline.

The methodological basis of this work is the development of technological schemes and algorithms for improving automated control systems of non-stationary gas transportation processes during the reconstruction of parallel gas pipelines. This can address key issues such as increasing the efficiency, reliability, and safety of gas transportation systems at the required level.

Thus, for the reconstruction of existing parallel gas pipelines, the use of modern separation devices and elements is proposed, along with the reliable management of non-stationary gas flows. The following structural scheme is proposed for this purpose (Figure 1).

In the case of an emergency situation in parallel pipelines, the diameter of the connecting pipe, which is installed to transfer the gas flow from one pipeline configuration to another, is taken to be equal to the diameter of the gas pipeline. At the same time, from the moment the gas pipeline is put into operation, along with its initial parameters ($P_1$, $P_2$, $G_0$, L, d), the distance between the connecting pipes ( $\ell$ ) is also recorded in the technical passport of the gas pipeline.

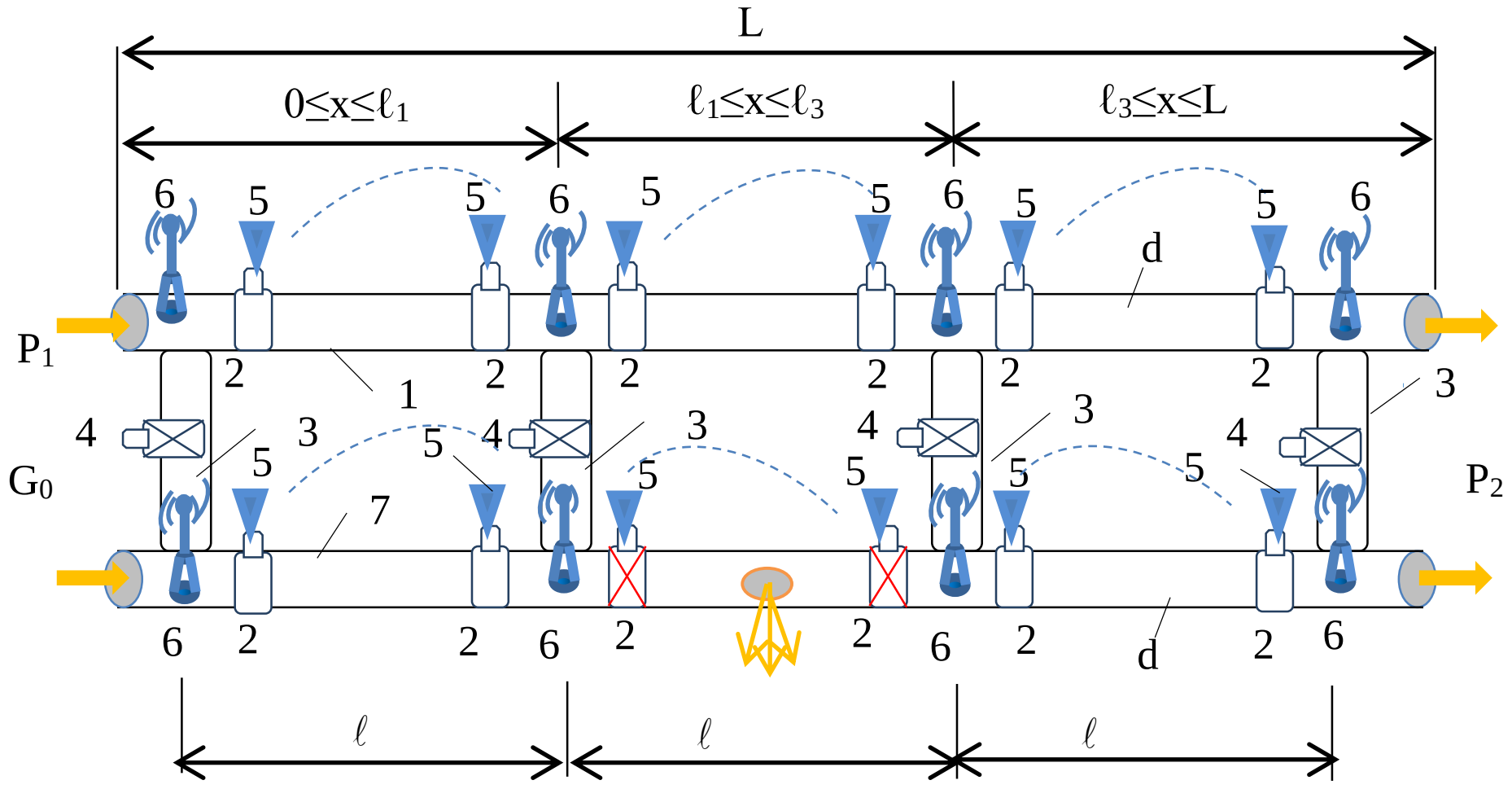

**Figure 2.5.1. Proposed structural diagram for the reconstruction phase of parallel gas pipelines. 1- undamaged pipeline; 2- pipeline isolation valves; 3- pipeline connecting the parallel pipelines (connecting pipes); 4- isolation valves of the connecting pipes; 5- valve position monitoring sensor; 6- wireless pressure sensor; 7- damaged pipeline.**

**L - length of the gas pipeline, m; $P_1$ - pressure at the beginning of the gas pipeline in stationary mode, Pa; $P_2$ - pressure at the end of the gas pipeline in stationary mode, Pa; $G_0$ - mass flow rate at the beginning of the gas pipeline, $\frac{Pa \times \sec}{m}$ ; d - diameter of the pipeline, m; $\ell$- distance between the connecting pipes, m.**

According to Figure 2.5.1, the installation of automatic shut-off valves on the parallel gas pipeline is planned as follows. They are to be installed after the connecting pipe at the beginning of the gas pipeline and before the connecting pipe at the end of the gas pipeline. In the other sections of the parallel gas pipeline, automatic shut-off valves are to be installed on both the right and left sides of each pipeline in the network.

For the efficient management of gas pipelines in emergency situations, the automatic activation of valves has proven to be effective in practice. Previously, valves calibrated based on the rate of pressure drop in the gas pipeline would automatically activate each valve separately if the pressure drop rate exceeded 0.1 MPa/min during emergency conditions [96]. Recently, however, with the application of valve position monitoring sensors, the automatic activation of valves has been successfully tested in practice and has been confirmed as a reliable technology [100].

The main advantage of the scheme we propose is that, in the event of a gas leak, the calibration (adjustment) of the valves is based not only on the rate of pressure drop but also on the activation of the valve position monitoring sensor. This improves the reliability of valve operation in emergency conditions. Specifically, when the integrity of one of the parallel gas pipeline pipes is compromised, the valve of the nearest shut-off valve to the leak location will close when the pressure drop reaches the calibrated value (the pressure drop rate reaches 0.2-0.5 MPa per minute [118]). The sensor monitoring the position of the closed valve automatically activates, which also activates the sensor monitoring the valve on the other side of the leak. The activated sensor automatically opens the valve it is monitoring and closes the shut-off valve. Therefore, both the rate of pressure drop and the activation of valve sensors will simultaneously trigger and close the valves on both sides of the leak.As shown in Figure 2.5.1, each of the automatic shut-off valves installed on the parallel gas pipeline is equipped with a position monitoring sensor. The sensors of the automatic shut-off valves installed on the pipeline configurations between the connecting pipes communicate with each other through the technical operation capabilities of the Internet of Things gateway. This allows the activation of one sensor to automatically trigger the sensor it is connected to. If the automatic shut-off valve on the left side of the leak point closes due to the limited pressure drop rate (Figure 2.5.1), the sensors designed to activate the valves will automatically close the valve on the right side of the leak. In other words, when the valve on the left side of the leak is automatically closed due to the pressure drop rate, the sensor monitoring the position of that valve activates, which triggers the sensor on the right side of the leak, connected via wireless communication channels. This will automatically close the valve on the right side of the leak. Thus, the automatic shut-off valves on both the right and left sides of the leak point will be closed simultaneously.

For the remote activation of pneumatic valves installed on the connecting pipes of the parallel gas pipeline, wireless pressure sensors are proposed [121]. As shown in Figure 2.5.1, it is suggested that the pressure sensors be installed on the section of the pipeline where the connecting pipes of the gas pipelines are located. These sensors allow for the monitoring of pressure variations in the pipeline configurations, enabling the control center to track the data via the internet.

It is known that the following technologies are available for predicting the technical operation of pipelines and managing emergency situations: the Internet of Things (IoT) [87]; machine learning and artificial intelligence [92,118]; advanced communication protocols such as 5G for real-time data transmission [93,95]; Arduino Uno controller - cloud storage for data collection and analysis [7,15,27]; hybrid models [128]; real-time monitoring or artificial intelligence [95,101]. These technologies are playing an increasingly important role in complex gas transportation systems and should be considered for the efficient operation of gas pipelines. It is well known that trunk and distribution gas pipelines are mainly designed in remote and inaccessible areas of residential settlements. Therefore, previously, due to the difficulty of monitoring the technological regime of these gas pipelines, their operational costs were correspondingly high. However, now, with the application of innovation, changes occurring in the gas pipelines can be monitored and tracked remotely via the internet through devices installed on them. As an example, the application of the Internet of Things (IoT) technology for monitoring and managing the non-stationary operating modes of the linear section of the gas pipeline can be highlighted [129]. The Internet of Things (IoT) technology is a system that enables devices containing hardware, software, and networks to interact, cooperate, and share information. Currently, IoT technology sensors remotely monitor the technological and emergency modes of these pipelines in real-time via the internet and send the collected data directly to control centers. Thus, if a leak is detected, the system automatically activates (closes) the valves to isolate the damaged section of the pipeline from the main section. Numerous local and international studies have addressed these issues [94,98,108,112,116].

In these works, the authors describe how the integration of the mentioned technologies and methodologies into gas pipelines will enhance the efficiency, safety, and reliability of gas pipeline systems. The materials presented, however, do not extensively address the methods for detecting leakage points in the linear pipeline or the analysis of non-stationary gas flow parameters along the entire length of the network as a result of this process. The key difference in our research is that it analyzes various transition processes of non-stationary gas flow in parallel gas pipelines. Specifically, in the event of a breach in the integrity of one of the pipeline sections of parallel gas pipelines, the process of activating the valves located to the left and right of the damaged section, as well as the remote activation of pneumatic valves installed on the connecting pipes for uninterrupted gas supply to consumers, are analyzed together, as shown in Figure 2.5.1. From the moment the valves on the pipeline are closed ($t=t_1$), the analysis includes the non-stationary gas flow parameters in three different sections of the damaged pipeline until the activation of the connecting pipes (change in pressure values) ($t=t_2$). Additionally, the technological scheme for accepting the control method of gasodynamic processes of these transition processes through a control center is developed. In this research, the use of pressure sensors for monitoring the transition processes, alongside the valve position sensor, is proposed for the localization of the damaged section in the emergency condition of parallel gas pipelines.

The methodology of the proposed control system for the parallel gas pipeline is based on the following scheme (Figure 2.5.2).

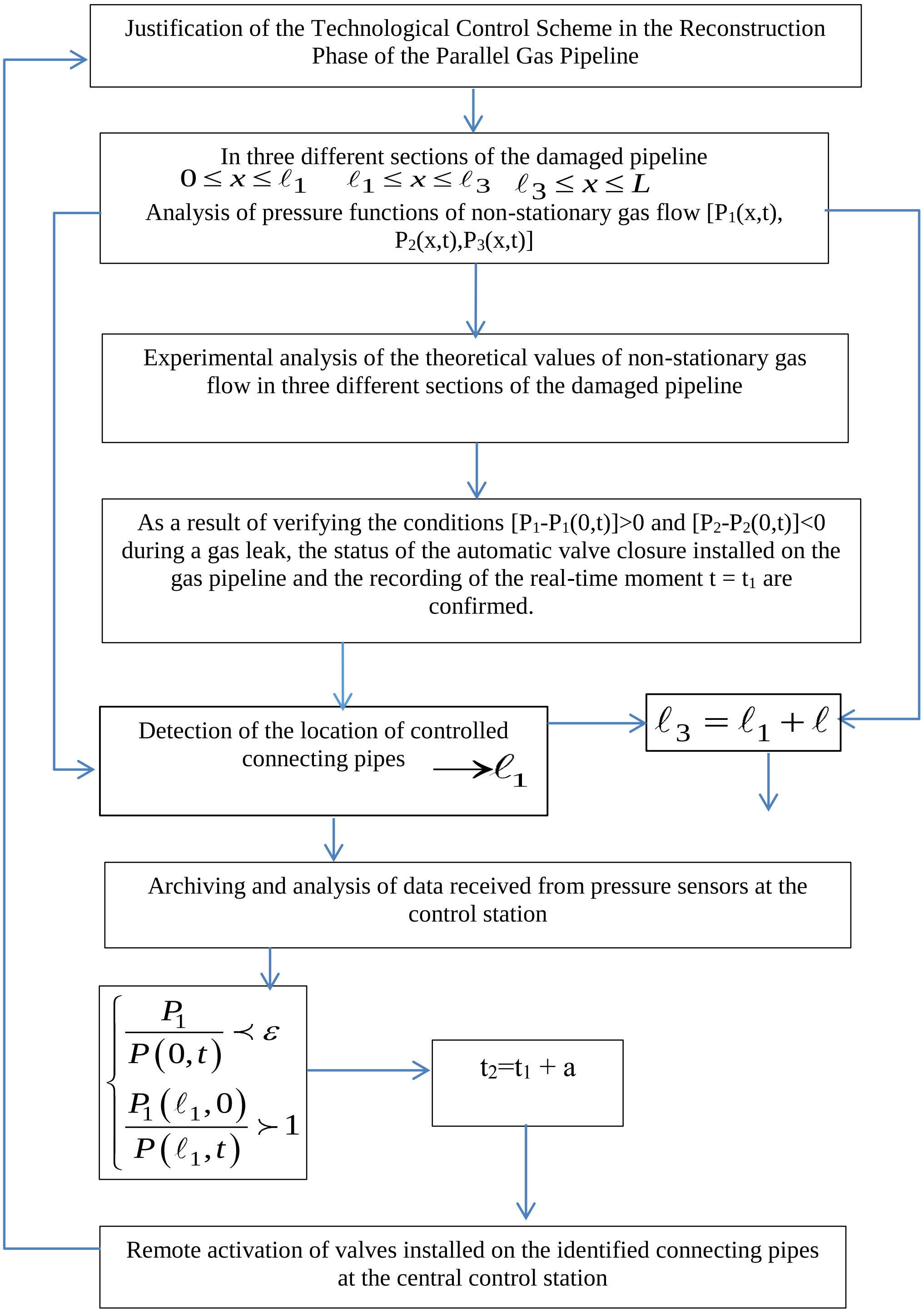

**Figure 2.5.2. Methodological diagram of the parallel gas pipeline control system.**

The analysis of the methodological scheme of the control system for the parallel gas pipeline (Figure 2.5.2) shows that data collection through the Internet of Things (IoT), as well as data storage and analysis using cloud storage networks, is proposed to be integrated into a control center, and a reporting scheme is suggested for an efficient wireless control system method. Thus, as a result of the breach of the integrity of any pipeline in the parallel gas pipeline, the valves close the valves located to the left and right of the leakage point at $t=t_1$ due to the pressure drop rate reaching its limited value and the activation of the sensors installed on them, effectively isolating the damaged section from the main structure of the parallel gas pipeline system.

Based on Figure 2.5.2., the mathematical models of the gas dynamic process occurring in three different sections of the damaged pipeline ($0 \le x \le \ell_1$, $\ell_1 \le x \le \ell_3$ , $\ell_3 \le x \le L$) from the moment the valves are closed have been thoroughly analyzed in the study(Chapter 4.1) we investigated, and the pressure functions of the non-stationary gas flow for the three sections are represented by the following expressions [88].

$$0 \le x \le \ell_1$$

$$P_1(x,t) = P_1(x,t_1) + \frac{c^2}{\ell_1}\int_{t_1}^{t} G_0(\tau)d\tau + \frac{2c^2}{\ell_1}\sum_{n=1}^{\infty}\mathrm{Cos}\frac{\pi n x}{\ell_1}\int_{t_1}^{t} G_0(\tau)\cdot e^{-\alpha_3(t-\tau)}d\tau + \\ + \frac{c^2}{2a\ell_1}\frac{dP_1(0,t_1)}{dx}\int_{t_1}^{t}\left[1+2\sum_{n=1}^{\infty}\mathrm{Cos}\frac{\pi n x}{\ell_1}e^{-\alpha_3\tau}\right]d\tau \tag{1}$$

$$\ell_1 \le x \le \ell_3$$

$$P_2(x,t) = P_2(x,t_1) - \frac{c^2}{(\ell_3-\ell_1)}\int_{t_1}^{t} G_{ut}(\tau)\times\left[\begin{array}{l} 1+2\sum_{n=1}^{\infty}(-1)^n\,\mathrm{Cos}\frac{\pi n(\ell_2-\ell_3)}{\ell_3-\ell_1}\times \\ \times\mathrm{Cos}\frac{\pi n(x-\ell_1)}{\ell_3-\ell_1}e^{-\alpha_4(t-\tau)}d\tau \end{array}\right] \tag{2}$$

$$\ell_3 \le x \le L$$

$$P_3(x,t) = P_3(x,t_1) - \frac{c^2}{(L-\ell_3)}\int_{t_1}^{t} G_s(\tau)\left[1+2\sum_{n=1}^{\infty}(-1)^n\,\mathrm{Cos}\frac{\pi n(x-\ell_3)}{L-\ell_3}\cdot e^{-\alpha_5(t-\tau)}\right]d\tau - \\ - \frac{c^2}{2a(L-\ell_3)}\frac{dP_3(L,t_1)}{dx}\int_{t_1}^{t}\left[1+2\sum_{n=1}^{\infty}(-1)^n\,\mathrm{Cos}\frac{\pi n(x-\ell_3)}{L-\ell_3}\cdot e^{-\alpha_5\tau}\right]d\tau \tag{3}$$

Here,

$$\alpha_3 = \frac{\pi^2 n^2 c^2}{2a\ell_1^{\,2}},\ \alpha_4 = \frac{\pi^2 n^2 c^2}{2a(\ell_3-\ell_1)^2},\ \alpha_5 = \frac{\pi^2 n^2 c^2}{2a(L-\ell_3)^2}$$

Let us assume that for the purpose of localizing the gas leak, the valves installed on the pipelines are automatically closed via the valves at $t=t_1=300$ seconds. Then, for each of the three different sections of the parallel gas pipeline, we accept the following given parameters of the specific gas pipeline in order to determine the regularity of the change in non-stationary gas flow for the considered case [88].

$$\mathrm{P}_1(0,t_1)=13.36\times10^4\,\mathrm{Pa};\ \mathrm{P}_1(5000,t_1)=12.82\times10^4\,\mathrm{Pa};\ \mathrm{P}_1(10000,t_1)=\mathrm{P}_2(10000,t_1)=12.19\times10^4\,\mathrm{Pa};$$

$$\mathrm{P}_2(14500,t_1)=11.56\times10^4\,\mathrm{Pa};\quad \mathrm{G}_0=10\,\frac{\mathrm{Pa}\times\mathrm{san}}{\mathrm{m}};\ \mathrm{P}_2(20000,t_1)=\mathrm{P}_3(20000,t_1)=11.24\times10^4\,\mathrm{Pa};$$

$$\mathrm{P}_3(25000,t_1)=10.86\times10^4\,\mathrm{Pa};\ \mathrm{P}_2=11\times10^4\,\mathrm{Pa};\ \mathrm{P}_3(30000,t_1)=10.4\times10^4\,\mathrm{Pa};\quad \mathrm{P}_1=14\times10^4\,\mathrm{Pa};\ \mathrm{L}=3\times10^4\,\mathrm{m};$$

$$\ell_1=1\times10^4\,\mathrm{m};\ \ell_2=1.45\times10^4\,\mathrm{m};\ \ell_3=2\times10^4\,\mathrm{m};\ \mathrm{c}=383.3\frac{\mathrm{m}}{\mathrm{sec}};\ 2a=0.1\frac{1}{\mathrm{sec}}$$

For the considered period, it can be assumed that $G_0(t)$=constant [88]. Then, using equation (1), we calculate the pressure values of the non-stationary gas flow along the length of the damaged pipeline section 1 every 120 seconds, and record them in the table below (Figure 2.5.1).

| x. km | $P_1(x.t)$, $10^4$ Pa | | | | | |
|---|---|---|---|---|---|---|
| | t=$t_1$+a. sec | | | | | |
| | a=0 | a=120. | a=240 | a=360 | a=480 | a =600 |
| 0 | 13.36 | 14.58 | 15.46 | 16.35 | 17.24 | 18.13 |
| 5 | 12.82 | 13.67 | 14.55 | 15.44 | 16.33 | 17.22 |
| 10 | 12.19 | 12.91 | 13.8 | 14.69 | 15.57 | 16.46 |

**Figure 2.5.1. Pressure variation values along the length ( $0 \le x \le \ell_1$ ) of the damaged gas pipeline section 1 as a function of time.**

It can be seen from Figure 2.5.1 that after the automatic valves on the damaged pipe are closed, the pressure along the length of section 1 of the gas pipeline increases over time. Specifically, at x = 0, during the period from a = 0 to 600 seconds, the pressure increases from $13.36 \times 10^4$ Pa to $18.13 \times 10^4$ Pa. At x = 5 km, the pressure at a = 0 seconds is $12.82 \times 10^4$ Pa, and after 600 seconds, it increases to $17.22 \times 10^4$ Pa. Similarly, at the final point x = 10 km, the pressure at a = 0 seconds is $12.19 \times 10^4$ Pa, and at a = 600 seconds, it increases to $16.46 \times 10^4$ Pa.

That is, during the period from a = 0 to 600 seconds, the rate of increase along the length is 1.35. Accordingly, the pressure values of the non-stationary gas flow along the length of the second section of the damaged pipeline are calculated every 120 seconds and recorded in the following table (Table 2.5.2).

| x. km | $P_2(x.t)$, $10^4$ Pa | | | | | |
|---|---|---|---|---|---|---|
| | t=$t_1$+a. sec | | | | | |
| | a=0 | a=120 | a=240 | a=360 | a=480 | a =600 |
| 10 | 12.19 | 10.59 | 9.7 | 8.81 | 7.93 | 7.04 |
| 14.5 | 11.56 | 10.52 | 9.63 | 8.74 | 7.86 | 6.97 |
| 20 | 11.24 | 10.01 | 9.13 | 8.24 | 7.35 | 6.46 |

**Table 2.5.2. Pressure variation along the length and with respect to time in the second section ( $\ell_1 \le x \le \ell_3$ ) of the damaged gas pipeline.**

From Table 2.5.2, it is observed that in the second section of the damaged pipeline, the pressure decreases as you move from the starting point to the endpoint. Specifically, at the starting point (x=10 km), the pressure decrease rate is 1.73 over a period of 600 seconds. At x=14.5 km, the rate of decrease is 1.65, and at the final point (x=20 km), the decrease rate is again 1.73. Based on the analysis, we can conclude that the pressure decrease rates at the starting and ending points of the second section are equal, while at the leak point, the rate of decrease increases.

Following this, the pressure values along the length of the third section of the damaged pipeline will be calculated and noted in Table 2.5.3 every 120 seconds.

| x. km | $P_2(x.t)$, $10^4$ Pa | | | | | |
|---|---|---|---|---|---|---|
| | t=$t_1$+a. sec | | | | | |
| | a=0 | a=120 | a=240 | a=360 | a=480 | a =600 |
| 20 | 11.24 | 10.52 | 9.63 | 8.74 | 7.86 | 6.97 |
| 25 | 10.86 | 10.01 | 9.13 | 8.24 | 7.35 | 6.46 |
| 30 | 10.4 | 9.19 | 8.3 | 7.41 | 6.52 | 5.63 |

**Table 2.5.3. Pressure variation values along the length ( $\ell_1 \le x \le \ell_3$ ) of the third section of the damaged gas pipeline and over time.**

From Table 2.5.3, it can be seen that in the third section, the pressure decreases along the direction from the starting point to the end point over time. For instance, at the starting point (x = 20 km), at t = 0 seconds, the pressure is $11.24 \times 10^4$ Pa. At t = 240 seconds, it decreases to $9.52 \times 10^4$ Pa, and at t = 600 seconds, it further decreases to $6.97 \times 10^4$ Pa. Thus, the rate of decrease at t = 600 seconds is 1.6. At the x = 25 km point, at t = 0 seconds, the pressure is $10.86 \times 10^4$ Pa, and at t = 600 seconds, the pressure decreases to $6.46 \times 10^4$ Pa, giving a rate of decrease of 1.68. At the x = 30 km endpoint, at t = 600 seconds, the pressure decreases to $5.63 \times 10^4$ Pa, resulting in a decrease rate of 1.85. Based on the analysis, we can conclude that the rate of pressure decrease increases from the starting point to the end point of the third section. As a result of the sharp decrease in pressure in the final section of the damaged gas pipeline, the volume of gas supplied to consumers will also decrease.

Using Tables 1, 2, and 3, we can construct a graph (Figure 2.5.3) illustrating the distribution of pressure over time along the length of the damaged pipeline in the 1st, 2nd, and 3rd sections of the complex gas pipeline.

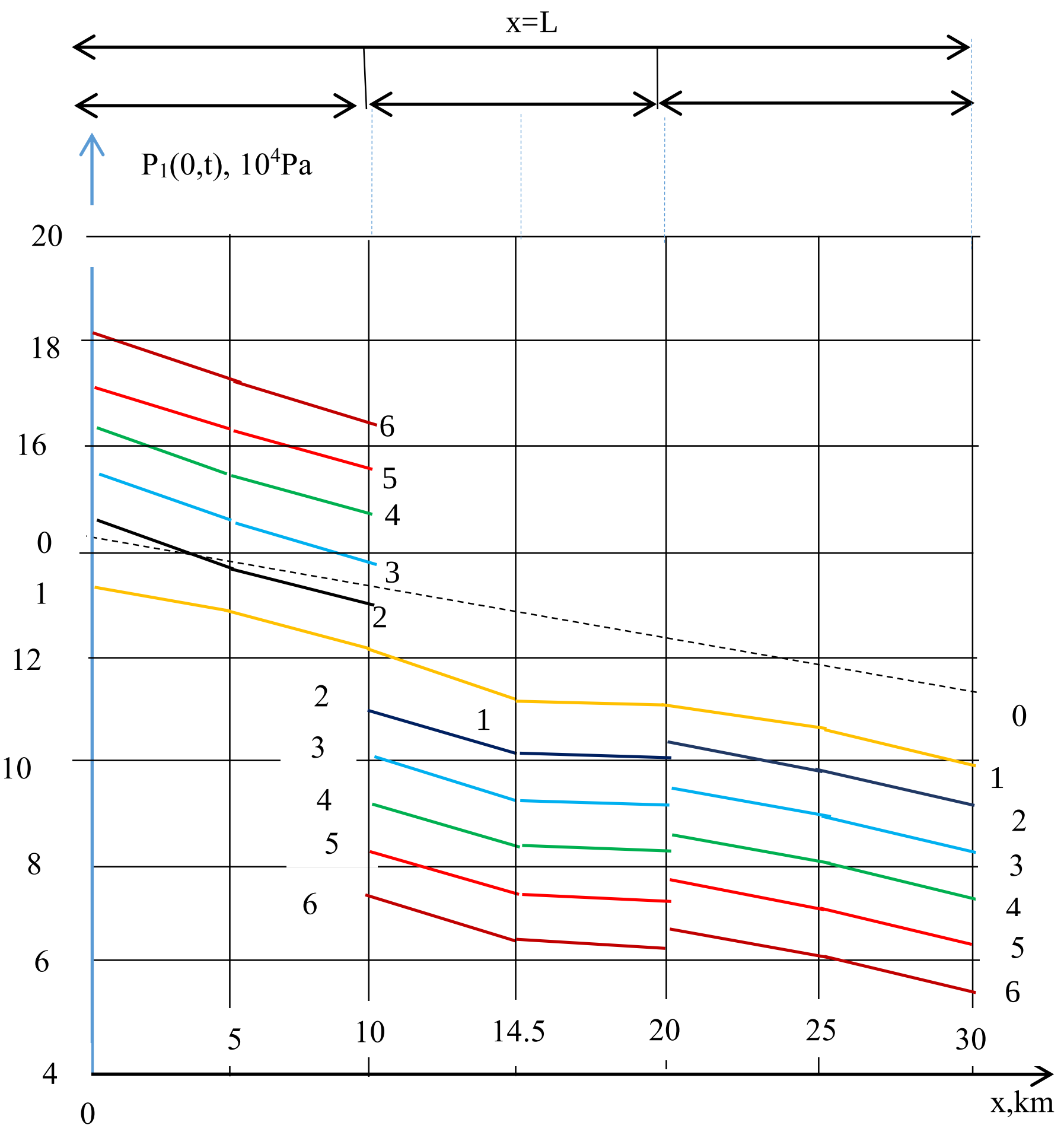

**Figure 2.5.3. The graph of pressure variation over time along the length of the gas pipeline after the valves at $x=\ell_1$ and $x=\ell_3$ are closed at t=t1. (0 - t=0 stationary state, 1 - $t=t_1$, 2 - $t=t_1+120$ sec, 3 - $t=t_1+240$ sec, 4 - $t=t_1+360$ sec, 5 - $t=t_1+480$ sec, 6 - $t=t_1+600$ sec).**

As shown in Figure 2.5.3, in the case where the integrity of the parallel gas pipeline arrangements is not violated, i.e., in a stationary regime, the pressure values along the length of $0 \le x \le \ell_1$ $\ell_1 \le x \le \ell_3$ both pipelines are equal and do not change over time (0-0 line). When the tightness of one of the pipelines is compromised (according to Figure 2.5.3, the 7th arrangement), no change occurs in the gas flow dynamics in the first pipeline, meaning it continues to operate in a stationary regime (0-0 line). In the pipeline with compromised integrity, however, the pressure values along the length decrease over time (as shown by the 1-1 line in Figure 2.5.3). At $t=t_1=300$ seconds, the pressure at the starting point x=0 changes from $P_1=14\times10^4$ Pa to $13.36\times10^4$ Pa, at x=10 km it decreases from $P(\ell_1,0)=13\times10^4$ Pa to $12.19\times10^4$ Pa, at the leak point (x=15 km) from $P(\ell_2,0)=12.6\times10^4$ Pa to $11.6\times10^4$ Pa, at x=20 km from $P(\ell_3,0)=12\times10^4$ Pa to $11.24\times10^4$ Pa, and at the final point x=L the pressure decreases from $P_2=11\times10^4$ Pa to $10.4\times10^4$ Pa.

Based on the analysis, it can be concluded that the pressure decreases along the length of the damaged pipeline, and the rate of decrease at the leak point is higher compared to the other sections.

As noted, in the case of an accident, at $t=t_1$, the valves in the sections $x=\ell_1$ and $x=\ell_3$ of the pipeline are automatically closed.

Based on the theoretical analysis of the gas pipeline's non-stationary flow regimes, it can be concluded that the opening of pneumatic valves on the connecting pipes at $t=t_1$ in the damaged sections is not advisable according to the gas dynamic process requirements. This is because, in the damaged section at $x=\ell_1$, the pressure is lower ($12.19\times10^4$ Pa < $13\times10^4$ Pa) compared to the pressure in the undamaged pipeline. As a result, the gas flow will occur in the reverse direction: gas will flow from the undamaged section into the damaged section, which will disrupt the continuity and stability of the gas supply to consumers.

Therefore, it is more appropriate to wait until the connecting pipes are opened at a later time $t=t_2$ , ensuring that the pressure in the damaged section exceeds the pressure in the undamaged section. According to Figure 2.5.3, during this time, the increase in pressure at the beginning of the damaged section at x=0 surpasses the pressure increase at $x=\ell_1$. This period is crucial for ensuring the protection of the main equipment of the gas pipeline, particularly the compressor stations. To avoid potential damage, the compression ratio ε\varepsilonε of the compressor's driving aggregates must be limited, with the value of $\varepsilon\leq1.3$ being the accepted standard to prevent overcompression, represented as $[P_1/P(0,t_2)\leq1.3]$ [88].

On the other hand, as shown in Figure 2.5.3, after the valves installed on the pipeline are automatically closed, the volume of gas supplied to consumers gradually decreases over time. To ensure continuous supply of the required amount of gas to consumers, it is essential to remotely activate the valves installed on the connecting pipeline sections at $x=\ell_1$ and $x=\ell_3$. For this purpose, the control center must first receive information about the activation of the valves and quickly determine their location within the pipeline.

The experimental analysis of the theoretical values of the non-stationary gas flow pressure shows that, after the valves are closed, the pressure in the starting section of the pipeline increases over time. However, in the final section, the pressure decreases with time. The analysis indicates that this particular gas dynamic behavior can only occur in the case of valve closure. Therefore, by analyzing data from the pressure sensors installed at the start and end points of the gas pipeline at the control center, this characteristic can be identified, which will signal the valve closure process.

To determine the location of these valves, we use the equation from (1) to express the pressure change at the starting point of the pipeline as follows:

The experimental analysis of the theoretical values of the non-stationary gas flow pressure indicates that, after the valves are closed, the pressure in the initial section of the pipeline increases over time. However, in the final section, the pressure decreases with time. This behavior suggests that this particular gas dynamic phenomenon occurs only in the case of valve closure. Therefore, by analyzing the data from the pressure sensors installed at the start and end points of the pipeline at the control center, this characteristic can be detected, which will indicate the valve closure process.

To determine the location of these valves, we use the equation from (1) to express the pressure variation at the starting point of the pipeline as follows:

$$P_1(0,t) = P_1(0,t_1) - \frac{2a\ell_1}{3}G_0 + \frac{2c^2}{\ell_1}G_0\sum_{n=1}^{\infty}\frac{e^{-\alpha n^2(t-t_1)} + e^{-\alpha n^2 t} - e^{-\alpha n^2 t_1}}{\alpha n^2} \quad (4)$$

Here , $\alpha = \frac{\pi^2 c^2}{2a\ell_1^{\,2}}$

The simplified form of equation (4) will be as follows:

$$P_1(0,t) = P_1(0,t_1) + \frac{2c^2 G_0}{\ell_1}(t-t_1)(1-C) - 2aG_0\ell_1\left(\frac{1}{3}+\frac{2}{\pi^2}\right) \quad (5)$$

The constant "C" is the Euler constant, and its value is taken as C = 0.577215.
After performing several mathematical operations in equation (5), we obtain the following quadratic equation:

$$2aG_0\left(\frac{1}{3}+\frac{1}{\pi^2}\right)\ell_1^2 + \left[P_1(0,t) - P_1(0,t_1)\right]\ell_1 - 2c^2G_0(1-C)(t-t_1) = 0 \quad (6)$$

Thus, by accepting only the positive root of the quadratic equation (6), we obtain the following analytical formula for determining the location of the automatic valves installed on the connecting pipes at the points

x= $\ell_1$ and x= $\ell_3$

$$\ell_1 = \frac{\sqrt{Z^2 + 8c^2G_0(1-C)(t-t_1)Z_1} - Z}{2Z_1} \quad (7)$$

Here, $Z = \left[P_1(0,t) - P_1(0,t_1)\right]$; $Z_1 = 2aG_0\left(\frac{1}{3}+\frac{2}{\pi^2}\right)$

As mentioned, once the parallel gas pipelines are commissioned and the value of $\ell$, representing the step length between the connecting pipes, is known at the control center, the following analytical expression is obtained to determine the location of the automatic valves installed at the x= $\ell_3$

$$\ell_3 = \ell_1 + \ell \quad (8)$$

As a result of the analysis conducted above, it can be noted that the technological determination of the activation time for the valves installed in the connecting pipeline should be based on the pressure values at x=0 and x= $\ell_1$ points of the damaged pipeline at t=$t_2$. Therefore, for the transmission of information to the control center and the reliable management of the gas pipeline, it is crucial to initially receive data from the sensors. To analyze and verify the data obtained at the control center, the following system equation is adopted for the creation of an algorithm based on machine learning mechanisms.

$$\begin{cases} \dfrac{P_1}{P(0,t)} \prec \varepsilon \\ \dfrac{P_1(\ell_1,0)}{P(\ell_1,t)} \succ 1 \end{cases} \quad (9)$$

In the initial stage, based on the changes in pressure values at both the starting and ending sections of the gas pipeline, the control center can record the valve activation status. Using equations (7) and (8), the location of the connecting pipelines to the right and left of the leak point is determined. Subsequently, the analysis of the data from the pressure sensors at the control center confirms condition (9). This allows for the detection of real-time operation, i.e., the time $t=t_2$. At this moment, the valves on the connecting pipeline are activated from the control center. As a result, consumers of different categories supplied by the parallel gas pipeline will be continuously provided with gas.

As seen in Figure 2.5.1 and Figure 2.5.2, during the time period $t=t_1+120$ seconds, the pressure at the damaged section of $x=\ell_1$ is greater than the pressure in the undamaged pipeline ($12.91\times10^4 > 12.19\times10^4$), and at this moment, ε=14.58/14=1,05<1.3. Therefore, condition (9) is satisfied. This means that, for the given accident scenario and the specific gas pipeline, the control center can activate the valves installed on the connecting pipelines at $t_2=t_1+120$ seconds. For the purpose of verifying the above findings and the derived empirical formula, we refer to the theoretical experimental data.

According to the data from Figure 2.5.1 and the specific gas pipeline information, using equations (7) and (8), we can determine the location of the valve to be activated on the connecting pipeline, using the pressure value at the starting point of the gas pipeline at $t=t_1+120$ seconds.

$$Z=\left[P_1(0,t)-P_1(0,t_1)\right]=14.58-13.36=1.22\,; \quad Z_1=2aG_0\left(\frac{1}{3}+\frac{2}{\pi^2}\right)=0.1\times10\times(0.333+0.203)=0.536$$

$$\ell_1=\frac{\sqrt{Z^2+8c^2G_0(1-C)(t-t_1)Z_1}-Z}{2Z_1}=0.9\times10^4\,\mathrm{m}\,; \quad \ell_3=0.9\times10^4+1\times10^4=1.9\times10^4\,\mathrm{m}$$

To obtain theoretical experimental results, the actual distance between the connecting pipes in the specific parallel gas pipeline was accepted as $\ell^f_1=1\times10^4$ m and $\ell^f_3=2\times10^4$ m. Thus, the relative error between the actual value of the valve installed on the connecting pipeline of the specific parallel gas pipeline and the value determined using equation (7) is 0.1. This value of the relative error is based on the linearization of the non-linear differential equation system of non-stationary gas flow in pipeline systems [89] and the simplification of equation (1). However, this indicator is acceptable for determining the location of the connecting pipeline because the distance between the connecting pipelines installed on parallel gas pipelines is sufficiently large from both a technological and economic perspective. For a sufficiently long distance, the relative error mentioned above does not interfere with identifying the location of the connecting pipelines. For instance, if the length of the parallel gas pipeline is at least 10 km, the distance between the connecting pipes will be at least 1 km from a technological and economic standpoint [123]. Therefore, the relative error value does not hinder the identification of the location of the connecting pipes. Finally, the established formulas (7) and (8) are useful for engineering calculations and for the central control center.

After determining the location of the connecting pipes, the data obtained from the pressure sensors installed at that location are used to determine inequality (9). If the condition is met, the pneumatic valves installed at those points are closed by the control center at $t=t_2$. Thus, the consumers supplied by the parallel gas pipeline will be continuously provided with gas. To manage the aforementioned processes, we adopt the following technological system based on Figure 2.5.2. (Figure 2.5.4).

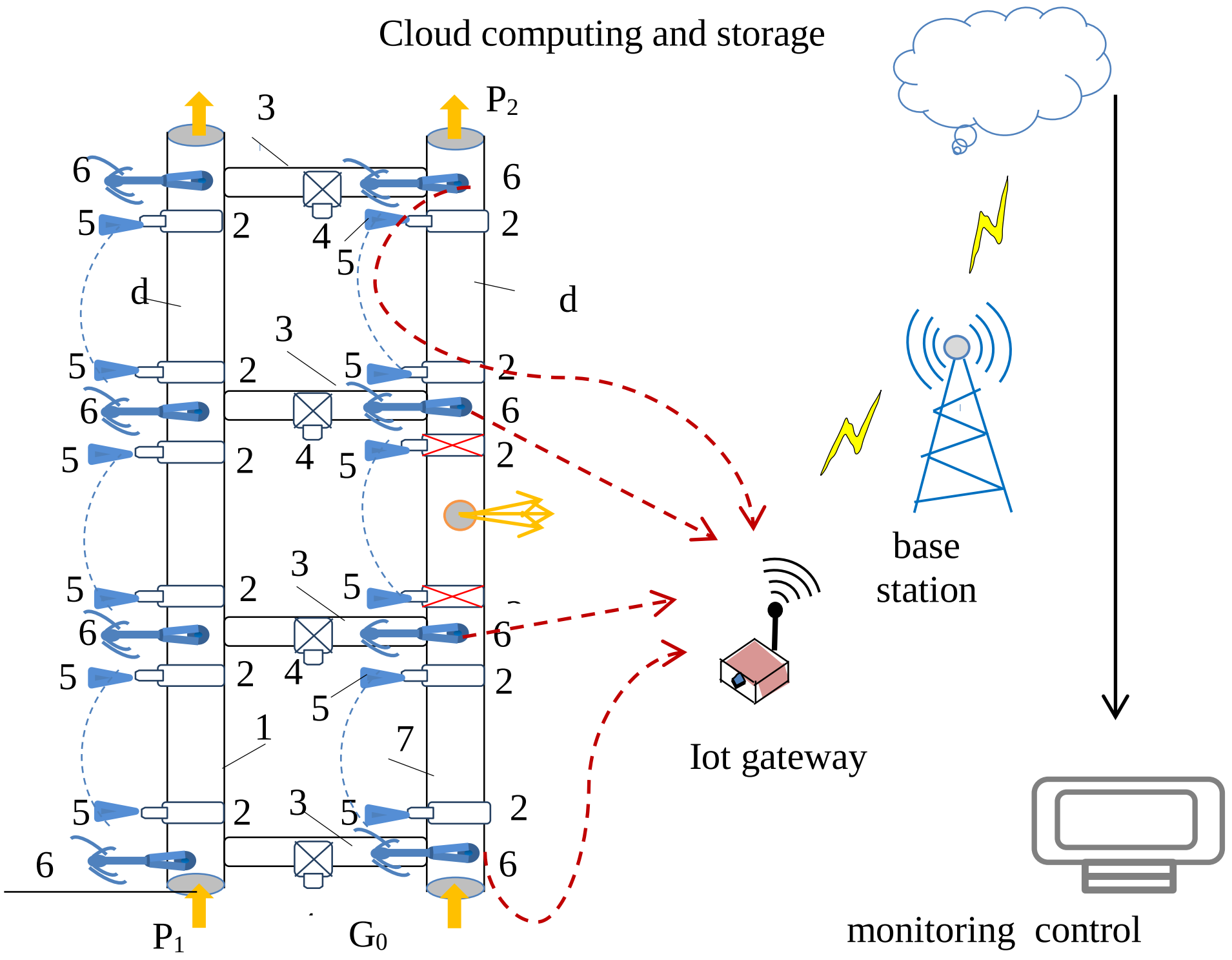

**Figure 2.5.4. Technological scheme of the IoT and Cloud computing and data storage system for managing non-stationary gas flow in a parallel gas pipeline.**

For monitoring the reliability of managing the non-stationary gas flow in the parallel gas pipeline, the implementation of the IoT and cloud computing and data storage system at the control center will be carried out in the following sequence, based on Figure 2.5.4.

As seen in Figure 2.5.4, when the integrity (tightness) of the 7th pipeline layout is compromised in the 8th section, the automatic valves on both sides of that section close simultaneously at time $t=t_1$ based on the control sensor's impulse. The central control center records this event based on the data from the pressure sensors located at the ends of the pipeline, $[P_1(0,t)-P_1>0$ and $P_2(L,t)-P_2<0]$, indicating that the valves are closing. This is because after that point, the pressure at the beginning of the pipeline increases, while the pressure at the final section decreases. Such a gas dynamic process can only occur after the automatic valves have closed.

At the control center, using the equations (7) and (8), the locations of the connecting pipes on both sides of the damaged section are determined. The data from the pressure sensors installed at those locations is transferred to the Internet of Things (IoT) gateway. From the IoT gateway, the data is forwarded to the cloud storage system for storage and analysis through the base station equipment. The data received from the control center's monitoring and control devices is analyzed and calculated. The calculated parameter values are determined, and if the condition of the system equation (9) is satisfied, the valves installed on the connecting pipes ($x=\ell_1$ and $x=\ell_3$) are opened by the personnel at the control center at time $t=t_2$. Thus, consumers fed by the parallel gas pipeline will be continuously supplied with gas.

The damaged section of the parallel gas pipeline is repaired by the personnel of a special emergency brigade according to technical-operational requirements for a specific period. After the repairs, the parallel gas pipeline is put back into operation in its previous stationary mode.

## CONCLUSION

As a result of the analysis of the mathematical model of non-stationary gas flow reflecting the control processes of parallel gas pipelines, an empirical algorithm based on machine learning technology has been defined for the information base of the control center. The application of the developed algorithm to the technological control system will increase the efficiency of automated decision-making processes.

The algorithms designed based on the developed schematic can be recommended for the reconstruction of existing parallel gas pipeline systems to solve engineering problems. This will enable the improvement of the central control system for the efficient operation of parallel gas pipelines.

The analysis shows that the pressure changes at the start and end sections of the parallel pipeline are directly dependent on the valve operation time. In the event of a failure, the rapid closure of the valves through the pressure sensors in the damaged sections of the parallel gas pipeline leads to the isolation of the damaged pipeline from the main network, significantly reducing gas loss to the environment and minimizing technological losses.

Based on the analysis of the data received from the pressure sensors in the event of a failure, an empirical formula has been derived for determining the location of the remotely activated connecting valves in real-time. The application of this formula at the central control center ensures the continuous supply of gas to consumers, thereby improving the reliability and safety of the parallel gas pipeline.

A technological scheme ensuring the application of modern control systems during the reconstruction phase of existing parallel gas pipelines has been proposed to guarantee the reliability and safety of parallel gas pipelines. To further develop this technology and implement it in operational gas pipelines, continuous scientific and practical research in this field is necessary. As a result, the implementation of an efficient management system using IoT and cloud technologies in parallel gas networks is recommended.

# Chapter 3. Investigation of Optimization Methods for the Parameters of Reconstructed Gas Pipelines

One of the important tasks in gas supply systems is the independence of the operating mode of consumer enterprises from the operating mode of the gas distribution pipeline providing gas fuel, in other words, ensuring uninterrupted supply of these enterprises with gas. That is, it is important, not to make the work of consumer enterprises dependent on the operating mode of the distribution pipeline supplying them with gas fuel, and to organize an uninterrupted gas supply to these consumers. For this purpose, it is important to develop new technological bases for the operation of the network and to develop a theoretically based computational scheme to prevent gas losses due to modern equipment installed in gas transport networks. One of the solutions to the problem of improving the reliability of gas pipelines is the use of effective scientifically based technologies. In the article, the efficient placement of connecting fittings to prevent gas losses in the operation of parallel gas pipelines was investigated, and proposals were made for theoretical and technical bases of the operation of parallel gas pipelines.

The principle of reconstructing gas pipelines is carried out in accordance with modern technological project standards and optimality criteria. The conducted research is dedicated to the comparative analysis of determining the parameters of looping installed to ensure the required gas flow to consumers. Thus, with the help of models, the optimal parameters of looping were determined to ensure an increase in the productivity of linear gas pipelines. At the same time, reconstruction options for the practical implementation of this principle were confirmed based on theoretical calculations. The construction of loops from pipes with optimal parameters on the gas pipeline has significant economic effects. The scale of this effect will increase in line with the rising share of gas consumption in the country's fuel balance.

## 3.1. Study of Efficient Placement of Connectors in Gas Pipelines

In stationary modes, the dynamic state of gas networks remains stable along the entire length of the gas pipeline, regardless of time, and therefore does not cause problems in the operation of control stations. However, in non-stationary modes, that is, in emergency modes, due to changes in the dynamic state of gas pipelines, targeted control of the operating mode is of great importance. At this time, to increase the reliability of networks and prevent gas losses, it is necessary to install additional fittings and equipment. In other words, a new computational scheme must be developed for the uninterrupted supply of gas to consumers and the amount of gas lost into the environment.

**Determination of the distance between the connectors in order to avoid losses in the operation of gas pipelines:**

It should be noted that the amount of gas lost to the environment when an accident occurs in pipelines also depends on the distance between the connecting fittings. Thus, when accidents occur on pipelines, the damaged part is separated from the main part of the gas pipeline by the connectors. The main purpose of the separation is to create conditions for repairing the damaged part of the pipe. Since maintenance work on the network in operation is a dangerous ones, it should be performed according to the appropriate rules and instructions. Failure to comply with these requirements can lead to accidents such as burns, explosions, and poisoning. Therefore, the executor is given a separate task for each repair. In these tasks, the composition of the brigade and the rules of safe working conditions, the permissible pressure for carrying out repairs should be specified.

**PROBLEM STATEMENT**

The normal pressure that can be left on the tape for repair is considered to be between 200 and 1200 Pa. Reducing the pressure to "0" at low pressures can create a risk of explosion in the gas pipeline due to the air mixture there. At low pressures, dropping the pressure to "0" can create a risk of explosion in the gas pipeline due to the air mixture in the pipeline. If the pressure is higher than the norm, it will practically not allow welding. Therefore, repairs to gas pipelines take a long time. During this period, the operating regime of gas supply to consumers is disrupted, and large amounts of gas are also released into the environment as a raw material. The amount of gas released into the environment is the sum of the losses resulting from two processes:

1) Gas lost into the environment during the period up to the closure of the automatic taps;

2) In terms of maintenance safety, the emptying of the pipeline between closed taps ensures that the gas is mixed with the environment (Figure 3.1.1).

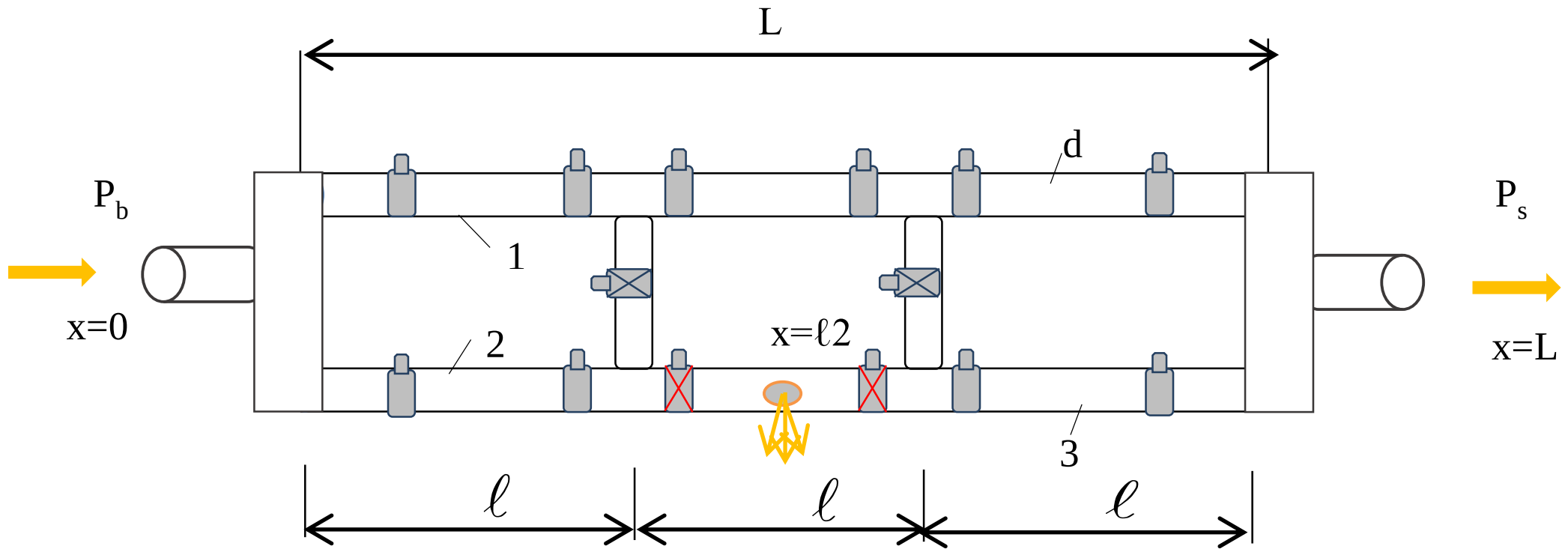

**Figure 3.1.1. Scheme of connecting the sections in the accident mode of parallel gas pipelines.**

The amount of gas lost in the initial operation is not dependent on the distance resulting from the placement of the automatic taps. In the second operation, the amount of gas lost is directly dependent on the distance between the automatic taps. So the longer the distance between the taps, the greater the amount of gas lost into the environment. On the other hand, as can be seen from the diagram (Fig. 1), the smaller the distance between the automatic taps tubes, the higher the cost of the pipeline. That is why we need to determine the distance between the automatic taps tubes economically. So the existing gas pipeline needs to have a number of automatic taps that do not exceed the cost of the pipeline due to the damage caused by the gas lost to the environment. That is to say,

$$S_{p.c.} \leq S_l \qquad (1)$$

Here,

$S_{p.c.}$ - is the cost incurred in the pipeline as a result of the installation of automatic taps on the pipeline and is determined depending on the operation mode and type of tap.

$S_l$ - damage caused by the release of gas into the environment in emergency mode is defined below:

$$S_l = (Q_1 + Q_2) \cdot C_{gas}, \text{ €} \qquad (2)$$

$Q_1$ - the amount of gas lost in the period from the moment of the accident to the time when the automatic taps are closed, m$^3$.

$Q_2$ - the amount of gas lost as a result of the emptying of the band for repair of the damaged part after the automatic taps are started, m$^3$.

$C_{gas}$ - is the cost of 1 m$^3$ of gas, €/m$^3$.

First of all, let's determine the amount of gas that has been lost in the period until the automatic tubes are shut down. For this purpose, in order for the gas pipeline shown in Figure 1 to be operational, the dynamic state of all processes must be determined in accordance with the operating principle when an accident occurs. In other words, by determining the mathematical solution of the occurring physical processes, we need to find the time-dependent distribution of the gas flow parameters along the pipeline. When parallel gas transmission systems operate in emergency mode, the system of specially derived nonlinear equations that mathematically represent the physical processes of the non-stable flow of a gas flow is as follows [1].

$$\begin{cases} -\dfrac{\partial P}{\partial X} = \lambda \dfrac{\rho V^2}{2d} \\ -\dfrac{1}{c^2}\dfrac{\partial P}{\partial t} = \dfrac{\partial G}{\partial X} \end{cases} \tag{3}$$

Here, c - is the speed of sound transmission of the gas for isothermal processing, m/sec.
P - the pressure in the sections of the gas pipeline, Pa.

x - coordinates along the axis of the gas pipeline, km.
t - time coordinates, sec.
λ - hydraulic resistance coefficient,
ρ - average gas density, kg/m$^3$
V - average gas flow rate, m/sec
d - gas pipeline diameter, m
G - gas flow mass consumption, $\frac{Pa \cdot \sec}{m}$

It is clear that the equation system (3) cannot be solved nonlinearly. Because these kinds of equations are not integrable at the moment. Approximate methods are used in engineering calculations. In other words, the equations are corrected. They use the following rationalization techniques to solve the equations: During the calculations, I.A.Charnov's method of rectangulation was shown to be more appropriate when the belts are out of order. The expression that characterizes this righteousness is as follows [78]:

$$2a = \lambda \frac{V}{2d} \tag{4}$$

Here, 2a - is Charnov's rationalization coefficient, then the equation system (1) becomes the following.

$$\begin{cases} \dfrac{\partial P}{\partial X} = 2aG \\ -\dfrac{1}{c^2}\dfrac{\partial P}{\partial t} = \dfrac{\partial G}{\partial X} \end{cases} \tag{5}$$

Here, G=ρv; P=ρ·c$^2$

(5) If we derive the equation system from the 1st expression and the 2nd expression according to x If we set $\frac{\partial G}{\partial X}$, then

$$\begin{cases} -\frac{\partial^2 P}{\partial X^2} = 2a\frac{\partial G}{\partial X} \\ \frac{\partial G}{\partial X} = -\frac{1}{c^2}\frac{\partial P}{\partial t} \end{cases} \quad (6)$$

If we replace the $\frac{\partial G}{\partial X}$ expression in the (6) expression 1 of the equation system, we get the heat transfer equation for the mathematical solution of the problem.

$$\frac{\partial^2 P}{\partial X^2} = \frac{2a}{c^2}\frac{\partial P}{\partial t} \quad (7)$$

Solution of the equation (7) to be the only solution and for the given process to be properly determined.

The distribution of the requested function (pressure) at t=0 is given when the starting conditions and inactivity forces are not taken into account. The boundary conditions will ensure a pre-determined regular change in pressure and flow at the end points of the pipeline during the time of the process considered [6]. Providing the correct boundary conditions completes the mathematical model of the process and enables a thorough and accurate study of occurring physical events.

$$t=0; \quad P_i(x,0)=P_b-2aG_0x: \; i=1, 2$$

In the process of treating the three points, we accept the boundary conditions (the flow is measured at the starting and ending points) as follows [44].

$$\text{at } x=0 \; \begin{cases} P_1(x,t) = P_2(x,t) \\ \frac{\partial P_1(x,t)}{\partial X} + \frac{\partial P_2(x,t)}{\partial X} = -2aG_0(t) \end{cases}$$

$$\text{at } x=L \; \begin{cases} P_1(x,t) = P_3(x,t) \\ \frac{\partial P_1(x,t)}{\partial X} + \frac{\partial P_3(x,t)}{\partial X} = -2aG_S(t) \end{cases}$$

$$\text{at } x=l_2 \; \begin{cases} P_2(x,t) = P_3(x,t) \\ \frac{\partial P_2(x,t)}{\partial X} - \frac{\partial P_3(x,t)}{\partial X} = -2aG_{ut}(t) \end{cases}$$

where $G_0(t)$, $G_s(t)$ and $G_{ut}(t)$ - the value of the mass flow measured at different times at the leak point for the start [9], end, and emergency modes of the gas pipeline, $\frac{Pa \cdot \sec}{m}$.

When the mode of gas pipelines is broken, the Laplace conversion method is preferred in practice. Laplace transformation means

$P(x,S) = \int_0^{\infty y} P(x,t)e^{-st}dt$ integral. Here S=a+ib is called a complex number, $i=\sqrt{-1}$ a virtual number [3]. Thus, if we apply the Laplace conversion method to heat transfer equations (7), we get a two-shape equation.

$$\frac{d^2P(x,S)}{dx^2} = \frac{2a}{c^2}\left[SP(x,S) - P(x,0)\right] \quad (8)$$

It is clear that the general solutions of the second-level differential equations (8) will be as follows.

$$P_i(x,S)=\frac{P_i(x,o)}{S}+A_1Sh\sqrt{\frac{2as}{c^2}}x+B_1Ch\sqrt{\frac{2as}{c^2}}x \qquad (9)$$

Applying the Laplace transformation to the starting and boundary conditions, we find the Ai and Bi coefficients in equation (9). Taking into account the values of the coefficients in statement (9), we get the transformed equation of the process in question, in other words the dynamic state of the gas pipeline. These equations will represent the reverse distribution of pressure throughout the band for the sections considered.

$$P_1(x,S)=\frac{P_b-2aG_1x}{S}+Z_2-\frac{2aG_{ut}(S)}{2\sqrt{\mu s}}\frac{ch\sqrt{\mu s\left(L-\ell_2-x\right)}}{sh\sqrt{\mu s}2L} \qquad (10)$$

$$P_2(x,S)=\frac{P_b-2aG_2x}{S}+Z_2-\frac{2aG_{ut}(S)}{2\sqrt{\mu s}}\frac{ch\sqrt{\mu s\left(L-\ell_2+x\right)}}{sh\sqrt{\mu s}L} \qquad (11)$$

$$P_3(x,S)=\frac{P_b-2aG_2x}{S}+Z_2-\frac{2aG_{ut}(S)}{2\sqrt{\mu s}}\frac{ch\sqrt{\mu s\left(L+\ell_2-x\right)}}{sh\sqrt{\mu s}L} \qquad (12)$$

here, $\mu=\frac{2a}{c^2}$,

$$Z_2=\frac{2aG_0}{2S\sqrt{\mu s}}\frac{sh\sqrt{\mu s}\left(x-\frac{L}{2}\right)}{ch\sqrt{\mu s}\frac{L}{2}}+\frac{2aG_0(S)}{2\sqrt{\mu s}}\frac{ch\sqrt{\mu s}(L-x)}{sh\sqrt{\mu s}L}-\,-\frac{2aG_k(S)}{2\sqrt{\mu s}}\frac{ch\sqrt{\mu s}x}{sh\sqrt{\mu s}L}$$

Since it is more appropriate to determine the mathematical expression of the physical process of the non-stable flow of a gas flow based on the solution of the problem in an undamaged configuration, we will use the inverse Laplace conversion rule of each limit. (8÷10) to find the original solution of the equations. So, according to the starting and boundary conditions, we get the following mathematical expression of the distribution of pressure along the belt axis as a function of time.

$$P_1(x,t)=Z_2-\frac{c^2}{L}\sum_{n=1}^{\infty}\text{Cos}\frac{\pi n\left(\ell_2+x\right)}{L}\int_0^t G_{ut}(\tau)\cdot e^{-\alpha_2(t-\tau)}d\tau \qquad (13)$$

$$P_2(x,t)=Z_2-\frac{c^2}{L}\sum_{n=1}^{\infty}\text{Cos}\frac{\pi n\left(x-\ell_2\right)}{L}\int_0^t G_{ut}(\tau)\cdot e^{-\alpha_2(t-\tau)}d\tau \qquad (14)$$

$$P_3(x,t)=P_2(x,t) \qquad \ell_2\le x\le L \quad \text{for this condition} \qquad (15)$$

Here, $\alpha_1=\frac{\pi^2(2n-1)^2c^2}{8aL^2}, \alpha_2=\frac{\pi^2n^2c^2}{8aL^2}$

$$Z_2=P_b-2aG_o\frac{L}{4}+4aG_oL\sum_{n=1}^{\infty}\frac{e^{-\alpha_1t}}{\left[\pi(2n-1)\right]^2}\text{Cos}\frac{\pi(2n-1)x}{L}+\frac{c^2}{L}\sum_{n=1}^{\infty}\text{Cos}\frac{\pi nx}{L}\int_o^t\left[G_o(\tau)-(-1)^nG_k(\tau)\right]\cdot e^{-\alpha_2(t-\tau)}d\tau+$$
$$+\frac{c^2}{2L}\int_o^t\left[G_o(\tau)-G_k(\tau)-G_{ut}(\tau)\right]d\tau$$

For the emergency mode of a parallel gas pipeline operating in a smooth hydraulic mode, we accept the following expressions, (13), (14), (15) and depending on the different gas leakage points along the line of damaged and undamaged structures: Data for determining the regularity of gas pressure changes.

$P_b = 55\times10^4$ MPa: $P_s = 40\times10^4$ MPa:$G_0$ = 30 Pa×sec/m = 3 кg·sec/m$^3$ : $2a = 0{,}1\frac{1}{\text{sec}}$ ;

$c = 383{,}3\frac{m}{\text{sec}}$ L = $10^5$ m $\ell_2 = 0{,}25\times10^4$ m; $\ell_2 = 0{,}5\times10^4$ m ; $\ell_2 = 0{,}75\times10^4$ m : d= 0,7 m:

First, we take the $\ell_2 = 0{,}25\times10^4$ m; $\ell_2 = 0{,}5\times10^4$ m; $\ell_2 = 0{,}75\times10^4$ m leakage points in the gas pipeline and write them down in Table 1 below, calculating the pressure values in the gas pipeline at t = 300 seconds every 5 km for each of its various locations.

Table 3.1.1.

| x, km | $\ell_2 = 0{,}25 \cdot 10^4$ m; | | | $\ell_2 = 0{,}5 \cdot 10^4$ m; | | | $\ell_2 = 0{,}75 \cdot 10^4$ m; | | |
|---|---|---|---|---|---|---|---|---|---|
| | $P_1(x,t)$, $10^4$ Pa | $P_2(x,t)$, $10^4$ Pa | $P_3(x,t)$, $10^4$ Pa | $P_1(x,t)$,$10^4$ Pa | $P_2(x,t)$, $10^4$ Pa | $P_3(x,t)$, $10^4$ Pa | $P_1(x,t)$,$10^4$ Pa | $P_2(x,t)$, $10^4$ Pa | $P_3(x,t)$, $10^4$ Pa |
| 0 | 49,1 | 49,1 | - | 52,09 | 52,09 | - | 53,81 | 53,81 | - |
| 5 | 48,83 | 47,49 | - | 51,23 | 50,88 | - | 52,31 | 52,78 | - |
| 10 | 48,7 | 45,78 | - | 50,78 | 49,61 | - | 51,67 | 51,72 | - |
| 15 | 48,5 | 43,92 | - | 50,29 | 48,27 | - | 50,99 | 50,63 | - |
| 20 | 48,22 | 41,92 | - | 49,75 | 46,85 | - | 50,28 | 49,51 | - |
| 25 | 47,88 | 39,75 | 39,75 | 49,17 | 45,35 | - | 49,55 | 48,34 | - |
| 30 | 47,48 | - | 40,42 | 48,56 | 43,74 | - | 48,78 | 47,13 | - |
| 35 | 47,03 | - | 40,92 | 47,92 | 42,03 | - | 47,99 | 45,86 | - |
| 40 | 46,54 | - | 41,28 | 47,24 | 40,17 | - | 47,17 | 44,52 | - |
| 45 | 46 | - | 41,49 | 46,53 | 38,17 | - | 46,31 | 43,1 | - |
| 50 | 45,42 | - | 41,6 | 45,8 | 36 | 36 | 45,42 | 41,6 | - |
| 55 | 44,81 | - | 41,6 | 45,03 | - | 36,67 | 44,5 | 39,99 | - |
| 60 | 44,17 | - | 41,52 | 44,24 | - | 37,17 | 43,54 | 38,28 | - |
| 65 | 43,49 | - | 41,36 | 43,42 | - | 37,53 | 42,53 | 36,42 | - |
| 70 | 42,78 | - | 41,13 | 42,56 | - | 37,74 | 41,48 | 34,42 | - |
| 75 | 42,05 | - | 40,84 | 41,67 | - | 37,85 | 40,38 | 32,25 | 32,25 |
| 80 | 41,28 | - | 40,51 | 40,75 | - | 37,85 | 39,22 | - | 32,92 |
| 85 | 40,49 | - | 40,13 | 39,79 | - | 37,77 | 38 | - | 33,42 |
| 90 | 39,67 | - | 39,72 | 38,78 | - | 37,61 | 36,7 | - | 33,78 |
| 95 | 39,31 | - | 39,29 | 37,73 | - | 37,38 | 35,33 | - | 33,99 |
| 100 | 38,81 | - | 38,81 | 37,09 | - | 37,09 | 34,1 | - | 34,1 |

Let's plot the distribution of the instantaneous pressure t = 300 sec along the pipeline line, depending on the various leak sites in the pipeline using Tables 3.1.1 and 3 (see Graph 3.1.1).

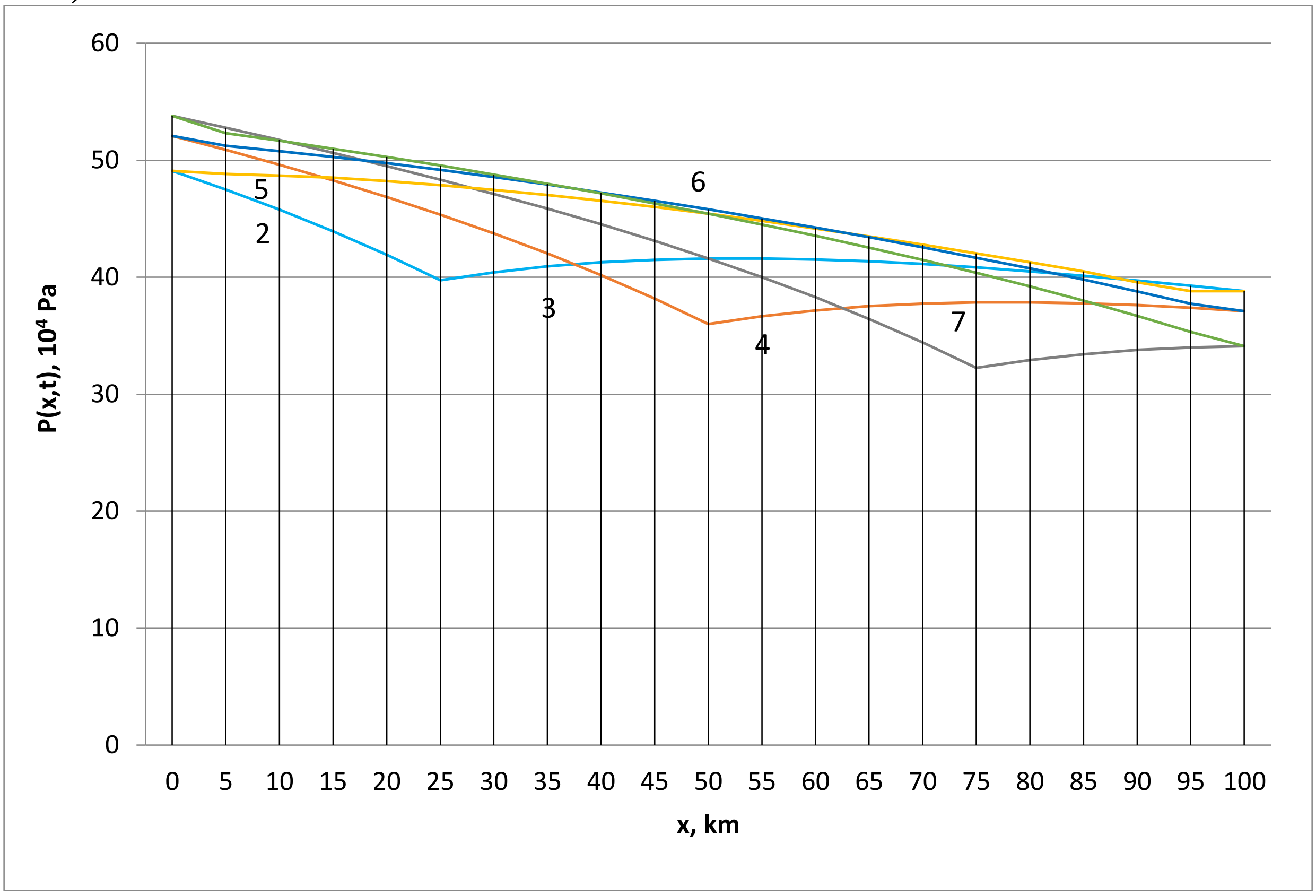

**Graph 3.1.1. Distribution of pressure along a faulty parallel gas pipeline line depending on the different locations of the leak.**

**1-fixed mode pressure change; 2-3-4 - 1 (harmless) change of pressure in adjustment; 5-6-7 - 2 nd order pressure change (damaged); 5-2- $\ell_2$ = 0,25 ×10⁴ m ; 6-3- $\ell_2$ = 0,5× 10⁴ m ; 7-4- $\ell_2$ = 0,75 × 10⁴ m ;**

The arrangements on parallel gas pipelines operating in single hydraulic mode have a common starting and ending point. Therefore, when one of the systems is damaged, the total amount of gas is directed to the damaged part of the gas pipeline. Graph 3.1.1 shows that the location of leakage points does not alter the quality of the image.

Using expression (11), first define expression $P_1(0,t)$, characterizing the pressure change at the starting point of the pipeline (x=0), and expression $P_2(L,t)$, characterizing the pressure change at the end point of the pipeline, point (x=L), using equations (12). After we accept simplifications $G_0(\tau)=G_k(\tau), G_{ut}(\tau)=const$ for the case n=2 using himin expressions, we define the difference of functions

$$\left[P_1(0,t)=P_b(0,t)\right]-\left[P_2(L,t)=P_s(L,t)\right]$$

as follows.

$$P_1(0,t)-P_2(L,t)=P_b-P_s-2aG_{ut}\frac{L}{2}+2aG_{ut}l_2$$

Here, taking into account that

$$P_s=P_b-2aG_0L \quad \text{we get}$$

$$P_b(0,t) - P_s(L,t) = 2aG_0L + 2aG_{ut}\left[\ell_2 - \frac{L}{2}\right] \tag{16}$$

We use equation (16) to determine the mass consumption ($G_{ut}$) of the gas lost into the environment in emergency mode.

$$G_{ut} = \frac{P_b(0,t) - P_s\ (L,t) - 2aG_0L}{2a\left[\ell_2 - \frac{L}{2}\right]} \tag{17}$$

As we know, the duration of the location of the leakage point and the working time of the automatic taps attached to the pipeline are taken equally [57, 79]. That is, the location of the leak point is determined, and the position of the pre-assembled connectors is known in the pipe network. If the time to locate the leak point is t = $t_1$, and the dependency between mass consumption and volume consumption is

$$Q = G_{ut}\frac{Fg}{\rho}$$

then the formula (15) will be expressed as follows:

$$Q_1 = \frac{Fgt_1}{2a\rho} \cdot \frac{P_s(0,t) - P_s(L,t) - 2aG_0L}{\left[\ell_2 - \frac{L}{2}\right]} \tag{18}$$

Obviously, to prevent the gas from leaking into the environment, they turn off the automatic taps brackets on the right and left of the damaged part at $t = t_1$. We obtained formula (18) by taking the integral of formula (17) to determine the total amount of gas lost from the leak point into the environment $\left(0 \le t \le t_1\right)$,

$$\int_0^{t_1} G_{ut}\, dt = G_{ut}\, t_1$$

Here, $F$ - is the cross-sectional area of the gas pipeline, m$^2$

$g$ - is acceleration of gravity, m/sec$^2$

ρ - is the average density of gas in the cross-sectional area, kg/m$^3$

*2a* - is the value of the smoothing factor, l/sec

$P_b$ (0,t) - the pressure at the starting point of the gas pipeline at time $t = t_1$ during the accident mode, kg/m$^2$

$P_s$ (L,t) - the pressure at the end point of the gas pipeline at time $t = t_1$ during the accident mode, kg/m$^2$

$G_0$ - mass flow rate of gas in stationary mode, kg·sec/m$^3$

$L$ – length of the gas pipeline, m.

$\ell_2$ – distance from the starting point of the gas pipeline to the leakage point, *m*.

Now, let's determine the amount of gas lost after the emergency shutdown of the automatic taps. The amount of gas lost during this period will be determined by the following expression, depending on the distance between the taps:

$$Q_2 = F \cdot \ell \cdot \frac{P_m}{P_0} \cdot \frac{T_m}{T_0}, \quad m^3 \tag{17}$$

Here, $P_m$ - represents the average pressure in the gas pipeline and is determined as follows:

$$P_m = \frac{2}{3}\left( P_b + \frac{P_s^2}{P_b + P_s} \right), \quad kg/m^2$$

$P_0$ - represents the value of atmospheric pressure, with $P_0 = 10^4 kg/m^2$

$T_m$ - top is the average temperature of gas in the pipeline, assumed to be within the range of 298 to 323K.

$T_0$ - is the absolute temperature of the gas at $0^0$C, with a value of $T_o = 293K$

Considering equations (16) and (17) in the inequality of equation (1), then

$$S_{p.c.} \le \left[ \frac{Fgt_1}{2a\rho} \frac{P_b(0,t) - P_s(L,t) - 2aG_0 \cdot L}{\left[ \ell_2 - \frac{L}{2} \right]} + F \cdot \ell \frac{P_m}{P_0} \frac{T_m}{T_0} \right] \cdot C_{gas} \tag{18}$$

By replacing the '≤' sign with the '=' sign in inequality (18), we can economically determine the distance between automatic taps.

From here,

$$\ell = \frac{P_o}{F \cdot P_{or}} \frac{T_o}{T_{or}} \cdot \left[ \frac{S_{k.x}}{C_{qaz}} - \frac{Fgt_1}{2a\rho} \frac{P_b(0,t) - P_s(L,t) - 2aG_0 L}{\left[ \ell_2 - \frac{L}{2} \right]} \right] \tag{19}$$

By substituting the values of the constants into equation (19) and making a few straight forward replacements, we obtain the distance between automatic taps as follows:

$$\rho = 0.73 \ kg/m^3, \ 2a = 0.1 \ l/\sec \ \ and \ \ T_m = 323K$$

are assumed.

$$\ell = \left[ \frac{10767}{P_{or}} \cdot \frac{C_{k.x.}}{C_{qaz} \cdot d^2} - 113{,}6 \cdot \frac{t_1}{P_{or}} \frac{P_b(0,t) - P_s(L,t) - 0{,}1 \cdot G_0 \cdot L}{\ell_2 - \frac{L}{2}} \right] \tag{20}$$

To determine the length of the effective distance between the cranes to be installed on the gas pipeline according to the formula (20), we take the above initial data of the gas pipeline:

$P_b = 55\times10^4$ MPa : $P_s = 40\times10^4$ MPa : $G_0 = 30$ Pa×sec/m= 3 кg×sec/m$^3$:

$2a = \frac{1}{\sec}$ ; $c = 383.3 \frac{m}{\sec}$; L = $10^5$ m: d= 0,7 m:

Using Table 3.1.1, we take the values $P_1(0,t_1)$ and $P_2(L,t_1)$ of the pressures at the beginning and end of the belt for option $l_2$=0,75·$10^4$ m of the leakage point for $t_1$=300 seconds.

$P_1(0,t_1) = P_b(0,t) = 53{,}81\times10^4$ Pa= $5.381 \cdot 10^4$ кg/m$^2$ :

$P_2(L,t_1) = P_s(L,t) = 34{,}1 \cdot 10^4$ Pa =$3.41 \cdot 10^4$ кg/m$^2$

The cost of installing an automatic crane with an electric drive with a diameter of d=0.7 m (steel DN 700 En, PN 25 C) $S_{p.c.}$ = 150000 € can be accepted.

Considering the current cost of 1m$^3$ of gas as $C_{gas}$=300 €/m$^3$, we can determine the economically optimal distance between the s using the principle outlined in equation (20). First, we determine the value of the average pressure in the gas pipeline.

if $P_{or} = 4.79 \cdot 10^4 kg / m^2, \ell = 7344m.$

Therefore, in order to avoid losses in the operation of gas pipelines, it is technologically and economically efficient that the distance between the connectors should be $\ell$ =7344 m according to the initial data above. The analysis of the diagram in Figure 1 suggests that after the automatic tubes have been shut down, the connectors placed on the right and left of the leakage point of the band should be operated simultaneously.

According to the rules of technical design of main gas pipelines, emergency automatic valves must ensure that the valves are turned off if the pressure in the gas pipeline drops to 10-15% of the operating pressure within 1-3 minutes. Linear connection fittings must be remotely controlled in accordance with the standards of the technological project. СНиП 2.05.06-85* - in accordance with clause 12.10, the distance between the linear connecting fittings installed on the pipeline should not exceed 10 km. On the other hand, when laying main gas pipelines in parallel, connectors must be provided: - for gas pipelines with the same pressure - with connecting fittings. The connectors are located on linear cranes (before and after cranes) at a distance of at least 40 km and not more than 60 km from each other.

Based on our theoretical research, it was found that the distance between the connectors can vary depending on the parameters of the parallel gas pipeline, the type of consumers and the cost of the equipment to be installed.

Therefore, the proposed new reporting scheme will have a significant impact on the effective placement of links.

From the above-mentioned considerations, it can be concluded that the picture of changes in the dynamic state of the pipeline, as well as the level of gas supply to consumers, directly depends on the distance between automatic cranes. Therefore, during the operation of pipelines,determining the distance between automatic taps installed on them is one of the main principles of gas pipeline reconstruction.

## CONCLUSIONS

Mathematical modeling based on the numerical solution of a complete system of equations describing the unsteady gas flow in a gas pipeline is used as the initial basis. Since the object of the study is not only the network itself, but also the connecting fittings and connectors, a complete system of gas dynamic equations was taken into account to determine its dynamic state in all processes, in accordance with the principle of operation in the event of an emergency event. During the analysis, the distribution of gas pressures along the length of the gas pipeline at various points of the leakage point was studied.

Violation of the tightness of one of the lines directly leads to a change in the operating mode of the other. It is established that the qualitative picture of the dynamic state of the damaged line is similar to that of an intact line. However, regardless of the operating conditions of the gas pipeline, the pressure drop rate at the leak point of the damaged pipeline is about 2-3 times higher than the pressure of the intact gas pipeline at this point.

The application of the proposed new design scheme for uninterrupted gas supply to consumers for various purposes, powered by a pipeline network during the operation of parallel gas pipelines, not only avoids losses during the operation of gas pipelines, but at the same time is acceptable for determining the optimal parameters of connecting fittings.

It can be concluded from the above that the dynamic state of the pipeline and the level of gas supply to consumers depend directly on the distance between the automatic taps. Therefore, the determination of the distance between the automatic cranes to be installed in the pipeline operation process is one of the basic principles of the design of gas pipelines.

## 3.2. Fundamental Principles of Reconstruction for the Linear Section of Gas Pipeline Systems

**Introduction.** The research is dedicated to the comparative analysis of modern technical solutions used to ensure the required transmission capacity of existing main gas pipelines during the reconstruction phase. However, despite numerous studies conducted in this field, there are significant shortcomings in the accuracy of determining the technological parameters of pipeline operation during mathematical modeling in the operation and reconstruction of pipeline systems. Unlike previous studies, different approaches have been investigated to improve the accuracy of mathematical models in solving the existing problem. In this regard, the development of new mathematical models and the improvement of their accuracy during the design of looping in the reconstruction phase of pipeline systems is one of the urgent tasks. According to the technological scheme of the linear section of the gas pipeline, networks are designed in two types:

1. Pipe networks with a constant diameter (medium diameter) along the length of the gas pipeline.
2. Telescopic networks.

In medium-diameter gas transmission systems, the diameter of the gas pipeline remains constant along the length of the linear section. In other words, the diameter is the same from the starting point to the endpoint of the gas pipeline. The linear section of a telescopic gas pipeline is designed in a scheme where the diameter of the gas pipeline increases and decreases sequentially at different points. This characteristic of the telescopic scheme prevents losses during its reconstruction due to the dismantling of existing networks. Thus, the variability in the diameters of the gas pipeline allows for its dismantling at one point and reapplication at another. Therefore, the application of the disassembly and assembly method in the reconstruction of telescopic linear gas transmission systems is carried out based on scientific and technical principles. During this process, to make the reconstruction more efficient, existing parts are dismantled, improved, and reassembled in the necessary sections of the system after thorough repair.

If the physical wear of the existing system's networks is significant, the aforementioned reconstruction option may not be economically efficient in some cases. Therefore, a second option is also considered for comparison in the reconstruction of telescopic systems. In this option, an independent network is designed from the system's source to the new consumers. In other words, the increase in the required gas flow corresponding to the volume of new consumers is not added to the existing network, thus eliminating the need for dismantling and assembly work on the existing network. In this case, the existing network supplies its previous consumers, while the independent network supplies the newly emerged consumers separately with the energy carrier. On the other hand, if the existing gas pipeline is significantly worn both physically and morally, the option of using the existing gas pipeline and the option of designing an independent gas pipeline to supply the new consumers with the energy carrier may not meet the technical and economic parameters required for reconstruction due to the system's unserviceability. In such cases, an additional third option is considered for comparison in the reconstruction of telescopic networks. In this option, the existing network parts are dismantled, and unserviceable working parts are removed. As a result, the distribution scheme of the system is reconstructed, taking into account the new consumers. This option also considers the application of modern equipment and devices.

### 3.3. Determination of Efficient Volume for Reuse of Existing Gas Pipeline Networks

The application of the three variants mentioned below for reconstruction purposes is appropriate for the purpose due to the inclusion of new consumers to existing telescopic gas pipeline systems.

1. Variant of dismantling existing networks and reusing them after improving the existing networks for exploitation;
2. Variant of designing additional gas pipelines from the source of gas supply to new consumers;
3. Variant of reconstructing the distributive scheme from the point of view of enhancing the technical condition of the existing system.

The second variant, as noted, requires more resources and entails long-term commitments, potentially leading to economic losses. Therefore, conducting a precise technical-economic analysis based on the parameters of the first variant is essential for the acceptance of the second variant. In this regard, proposing the following calculation scheme (Scheme 3.3.1) for the dismantling and reusing variant is crucial for the decision-making process.

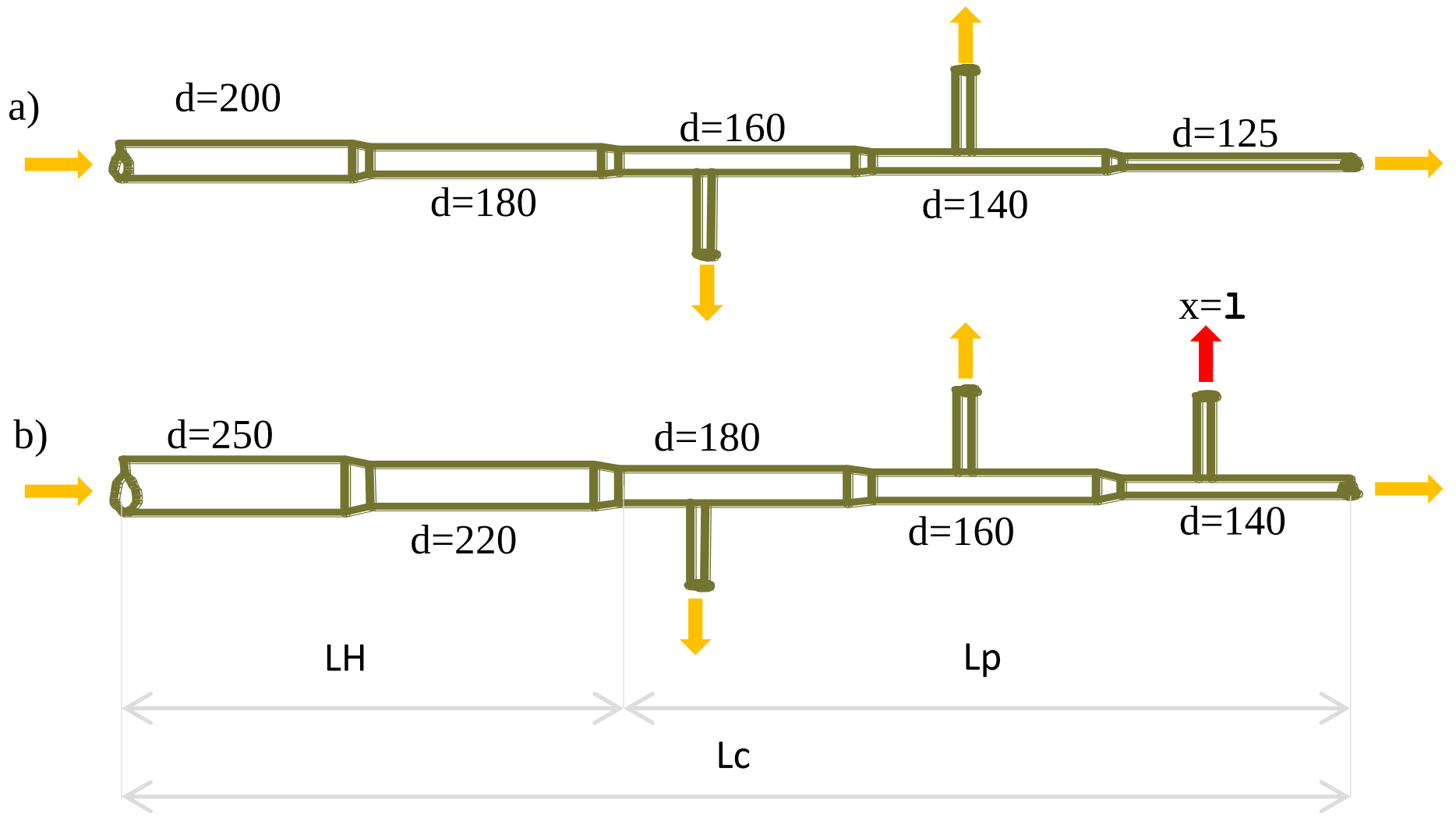

**Scheme 3.3.1. Diagram of the reconstructed section of the telescopic gas pipeline network.**
**a) Existing scheme before reconstruction,**
**b) Proposed scheme during reconstruction:**
**Lc- length of the gas pipeline, m; $L_H$ -Total length of new pipes installed for the reconstruction of the gas pipeline, m ; $L_p$- Total length of pipes reused for the dismantling and reinstallation of existing networks, m**

Scheme 3.3.1 proposes the reconstruction of the existing gas pipeline to accommodate the connection of a new pipeline segment at point x= l on the network (Scheme 3.3.1, b). This adjustment aims to increase the capacity of the gas pipeline at its origin to adequately supply gas to the newly connected consumers on the network. As a result of the increased gas consumption, the hydraulic resistance of the existing pipeline network increases, disrupting the operational regime of the gas dynamic process. Consequently, due to the significant pressure losses in the network, the efficiency of gas delivery to consumers decreases, leading to a decrease in the reliability of the gas pipeline. Therefore, increasing the hydraulic efficiency of the gas pipelines should be considered to eliminate frictional forces and local hydraulic resistance losses.To achieve this, the hydraulic efficiency of the network should be adjusted so that its characteristic intersects with the flow (consumption) rate. In this regard, to enhance the hydraulic efficiency of the existing telescopic gas pipeline, the cross-sectional area (diameter) of the pipeline is increased along the length of the network (Scheme 3.3.1, b). This reconstruction principle is based on this concept.

The technical decision for the reconstruction project should be based on the achievements of scientific and technical progress in the construction and design of gas transportation systems. Reconstruction involves implementing appropriate measures to improve both the system as a whole and the technical level of each individual element. As a rule, reconstructing gas transportation systems with a very low technical level is not suitable for the intended purpose. The reconstruction of engineering systems and equipment becomes purposeful when the efficiency of the additional expenses incurred for its implementation is not lower than the efficiency of other methods for intensifying technical processes during operation. The economic criterion is based on ensuring the minimum incurred expenses to achieve the effective improvement of gas transportation systems.

For this purpose, using Scheme 3.3.1, we accept the following inequality:.

$$S_p L_p \leq S_H L_H \alpha_t (1-F)(1-M) \tag{1}$$

Here,

$S_p$ -represents the direct expenses incurred for the reconstruction of the gas pipeline in the first variant due to the reuse (reinstallation) of dismantled existing pipelines of various diameters required for the gas pipeline's sections.

$S_H$ -represents the direct expenses incurred for the reconstruction of the gas pipeline in the first variant due to the installation of new pipes in the dismantled sections of the gas pipeline.

M -represents an equivalent reflecting the moral wear indicator of the reconstructed existing gas pipeline.

F- represents an equivalent reflecting the physical wear indicator of the reconstructed existing gas pipeline.

$\alpha_t$ -represents an equivalent reflecting the remaining service life of the reused existing pipe compared to the service life of the new pipe, determined according to the following guidelines.

$$\alpha_t = \frac{t_{\text{remai}}}{t_{x.m}}$$

$t_{\text{remai}}$ -is determined according to the following guidelines and represents the remaining service life of the reused pipe line of the existing gas pipeline.

$$t_{\text{remai}} = t_{x.m} \cdot (1-F), \text{ il}$$

$t_{x.m}$ -represents the service life of the system resulting from the application of new pipes, which varies depending on the purpose and material of the system. For example, for cast iron pipes, a service life of 40 years is accepted.

The reliability indicator of the gas pipeline must not be overlooked during the reconstruction process. Considering the reliability parameter is crucial for determining the reconstruction variant, as it lays the groundwork for ensuring the supply of energy carriers to consumers. One of the main parameters of this provision is the average accident intensity. In this regard, to assess the safety of gas pipeline operation, it is necessary to take into account both the accident intensity index of the reused (reinstalled) pipe and the accident intensity indices of the new pipes. This will further refine the technical performance of the reconstruction variant. Consequently, the probability of accidents occurring for each pipeline with different technical conditions regarding the accident intensity parameter will be as follows:.

$$\lambda_p = \omega_p L_p \quad ; \qquad \lambda_H = \omega_H L_H$$

If we consider these equalities (1) inequality, the following expression is obtained:.

$$S_p \omega_p L_p^2 \le S_H \omega_H L_H^2 \alpha_t (1-F)(1-M) \tag{2}$$

If we consider the equation $L_H=(L_c-L_p)$ based on Scheme 1 and the substitution of $\alpha_\lambda = \frac{\omega_H}{\omega_p}$, the inequality (2) would be as follows.

$$S_p L_p^2 \le S_H (L_c - L_p)^2 \alpha_\lambda \alpha_t (1-F)(1-M) \tag{3}$$

The parameters indicated on the left side of inequality (3) result from the utilization of existing pipes, while the expression on the right side characterizes the expenses incurred for the reconstruction of the system due to the use of new pipes. It is evident that when the condition of inequality (3) is satisfied, it is appropriate to reuse the existing pipes after dismantling for restoration, in line with the purpose of the reconstruction variant. In this case, the total cost incurred for the gas pipeline's mentioned reconstruction variant will be equal to the sum of the expenses on both sides of inequality (3), and we accept that expression as equivalent to the function with argument Lp.

$$S_p L_p^2 + S_H (L_c - L_p)^2 \mathrm{K} = f(L_p) \tag{4}$$

Here,

K= $\alpha_\lambda \alpha_u \alpha_t (1-M)$ represents the technical condition of the reconstructed existing gas pipeline.

If we analyze equation (4), we see that the cost of the reconstruction work on the system is quadratically dependent on the length Lp of the reused pipes. To determine this dependence, we accept the following initial data:

We consider a gas pipeline with a diameter of d = 219 mm and a length of the network of Lc = 10000 m. According to information available in the database, the economic indicators of a cast iron gas pipeline are as follows:

If the approximate cost of installing 1 meter of a gas pipeline with a diameter of d = 219 mm is 100$/m and the approximate cost of installing 1 meter of a cast iron pipe with the same diameter is 50$/m, then:

$S_H$ = 100 + 50 = $150/m = 150 $\cdot 10^3$ $ /km.

Let's assume the cost incurred for the dismantling and reinstallation of the existing pipes of the gas pipeline is 50$/m.

Thus, $S_p$ = 50 + 30 = 80$/m = 80 $\cdot 10^3$ $/km.

For calculating equation (4), we accept three variants:

1. The technical condition of the existing pipes of the gas pipeline is at a low level, characterized by the factor K = 0.25.
2. The technical condition of the existing pipes in the gas pipeline is at a medium level, characterized by the factor K = 0.5.
3. The technical condition of the existing pipes of the gas pipeline is at a good level, characterized by the factor K = 0.8.

We compile the following table (Table 3.3.1) using the results of the calculations performed:

| $L_p$, km | $\alpha \cdot S_H \cdot (L_s - L_p)^2$, M $ | | | $S_p L_p^2$, M $ | | | $f(L_p)$, MN $ | | |
|---|---|---|---|---|---|---|---|---|---|
| | K=0,2 | K=0,5 | K=0,8 | K=0,2 | K=0,5 | K=0,8 | K=0,2 | K=0,5 | K=0,8 |
| 0 | 3000 | 7500 | 12000 | 0 | 0 | 0 | 3000 | 7500 | 12000 |
| 1 | 2430 | 6075 | 9720 | 80 | 80 | 80 | 2510 | 6155 | 9800 |
| 2 | 1920 | 4800 | 7680 | 320 | 320 | 320 | 2240 | 5120 | 8000 |
| 3 | 1470 | 3675 | 5880 | 720 | 720 | 720 | 2190 | 4395 | 6600 |
| 4 | 1080 | 2700 | 4320 | 1280 | 1280 | 1280 | 2360 | 3980 | 5600 |
| 5 | 750 | 1875 | 3000 | 2000 | 2000 | 2000 | 2750 | 3875 | 5000 |
| 6 | 480 | 1200 | 1920 | 2880 | 2880 | 2880 | 3360 | 4080 | 4800 |
| 7 | 270 | 675 | 1080 | 3920 | 3920 | 3920 | 4190 | 4595 | 5000 |
| 8 | 120 | 300 | 480 | 5120 | 5120 | 5120 | 5240 | 5420 | 5600 |
| 9 | 30 | 75 | 120 | 6480 | 6480 | 6480 | 6510 | 6555 | 6600 |
| 10 | 0 | 0 | 0 | 8000 | 8000 | 8000 | 8000 | 8000 | 8000 |

Using Table 3.3.1, let's plot the graph of F($L_p$) as a function of Lp (Graph 3.3.1).

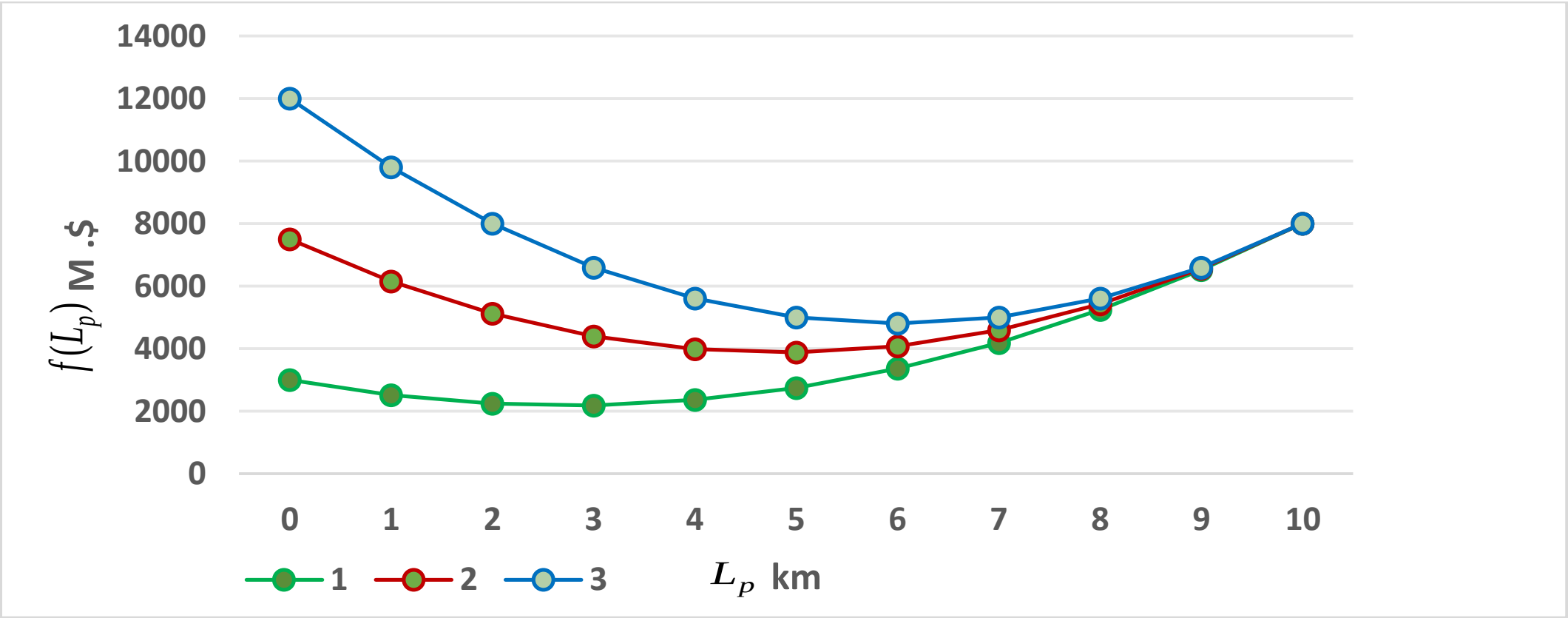

**Graph 3.3.1. Dependence of F(Lp) function on the values of different technical indicators and lengths of reused pipes (1- K = 0.2; 2- K = 0.5; 3- K = 0.8).**

In Graph 3.3.1, the cost of the reconstruction works, as indicated, reaches its minimum value depending on the length of the reused (dismantling-reinstallation) pipe. The minimum value of the function varies depending on the technical condition indicator of the existing pipes. Specifically:

In Variant 1 with a technical condition indicator of K = 0.2, the minimum value is attained at $L_p^* = 3$ km distance.

In Variant 2 with a technical condition indicator of K = 0.5, the minimum value is attained at $L_p^* = 5$ km distance.

In Variant 3 with a technical condition indicator of K = 0.8, the minimum value is attained at $L_p^* = 6$ km distance.

Exactly, as the technical condition indicator of the existing pipes increases, the cost of their reuse also increases. According to the horizontal axis of Graph 1, the length of utilization of the existing pipes decreases from Lp = 0 to Lp = $L_p^*$ and then begins to increase for values greater than that. Therefore, Lp = $L_p^*$ is the optimal length of utilization of the dismantled existing pipe line from a technical-economic perspective.

Additionally, from the analysis of Graph 1, it can be concluded that the function (4) varies extremally with respect to Lp. Considering this point, we take the derivative of (4) with respect to Lp and after setting it equal to zero, we determine the optimal length of the reused pipe for the effectiveness of gas pipeline reconstruction ($L_p = L_p^*$).

$$L_p^* = \frac{L_c}{1 + \frac{S_p}{2S_H K}} \quad (5)$$

Using principle (5), it's possible to effectively conduct the reconstruction of telescopic gas pipelines from a technical-economic perspective. In other words, based on the technical condition of the line sections of the reconstructed existing gas pipeline (including physical wear, remaining service life for operation, accident intensity index, and spiritual wear), the optimal volume of reused pipes to be used for the demounting-remounting operation is determined.

### 3.4. Investigation of Optimal Parameters for Looping Installed in the Reconstruction of Linear Gas Pipelines

Unlike telescopic networks, the acceptance of reconstruction variants for medium-diameter linear gas pipelines is impractical. This is because when new consumer facilities are connected to medium-diameter linear gas pipelines, the narrowest section of the pipeline diameter along its length increases according to the volume of consumers, while the length of the pipeline remains constant. This does not allow for the reuse of the decommissioned operational section. In such cases, designing a separate gas pipeline for new consumers may always be more effective for the reconstruction variant. The distance of new consumers from the source significantly increases the capital investment required for reconstruction. Therefore, alternative ways of supplying gas to new consumers from the existing network should be explored. Connecting a loop to the linear gas pipeline can facilitate using the existing network. The connection of the loop to the linear section of the gas pipeline can occur in three cases depending on geographical conditions and relief (Scheme 3.4.1).

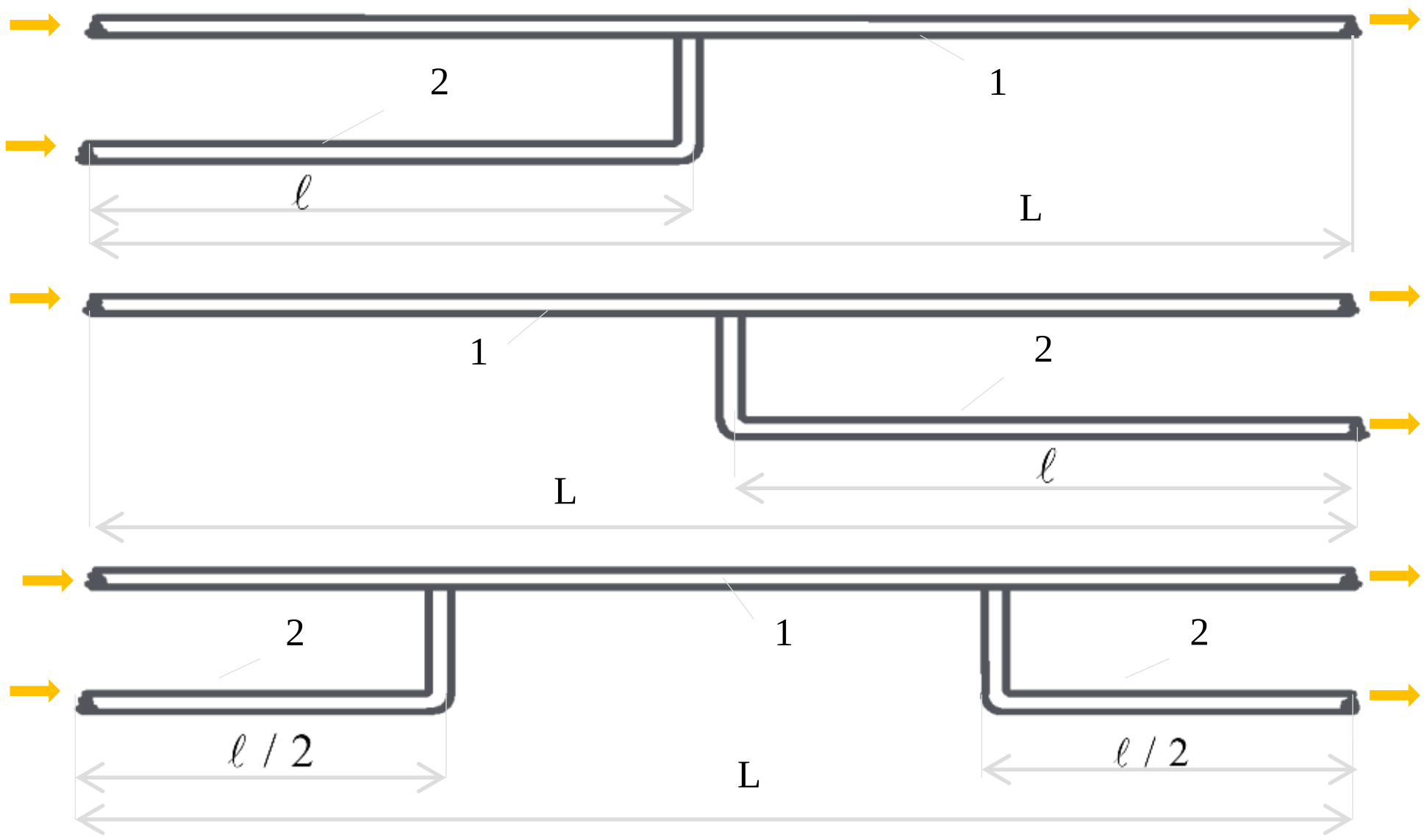

**Scheme 3.4.1. Scheme for connecting the loop to the linear section of the gas pipeline.**

**a) Connected to the beginning section of the loop gas pipeline; b) Connected to the end section of the loop gas pipeline; c) Connected to both the beginning and end sections of the loop gas pipeline.**
**1- Linear section of the gas pipeline; 2-Loop; L- Length of the gas pipeline; $\ell$ - Length of the loop.**

Thus, based on the indicators of physical and moral wear of medium-diameter networks, three variants are appropriate for the application of the reconstruction method resulting from the connection of new consumers:

1. Application of the loop system to increase the gas transmission capacity of the existing gas network for the supply of new consumers with gas fuel;
2. Project designing an additional gas pipeline from the source of the system to the new consumers;
3. Renovation variant of the scheme from the perspective of improving the technical condition of the existing system.

The integrity of the reconstruction should be approached from a point of view where the parameters of technical and economic indicators are determined for all three variants, and the minimum should be accepted. The gas transmission capacity of the pipeline section for gas flow under normal regimes is calculated according to the following principles without taking into account the length of the gas pipeline topography. [58].

$$Q_0 = F\sqrt{\frac{(P_H^2 - P_k^2)\cdot D}{\lambda z \mathrm{RT}\, L}} \qquad (6)$$

The theoretical research has determined that the length of the loop to be connected to the existing gas network to accommodate new consumers' demand is determined according to the following guidelines [1]:

$$l = \frac{4\cdot L\cdot \left[b^2 - 1\right]}{3\cdot b^2}, m$$

Here,

$$b = 1 + \frac{Y}{100}$$

"Y" represents the percentage that takes into account the consumption volume of new consumers fed from the existing gas network.

**Example:** For a pipeline with a length of L=50000 m, the volume of new consumers connected to the pipeline has increased by Y=20%. Correspondingly, a 20% increase in the gas transmission capacity (efficiency) of the pipeline is required. The addition of these new consumers to the existing engineering system results in $b = 1 + \frac{20}{100} = 1{,}2$.

Then,

$$\ell = \frac{4\cdot 50000\cdot \left[(1{,}2)^2 - 1\right]}{3\cdot (1{,}2)^2} = 20370m$$

So, in order to supply the new consumers with natural gas, a pipeline of the same diameter with a length of 20370 m should be added to the existing network. The economic feasibility of adding the pipeline at the specified length needs to be clarified. Determining the parameters of the pipeline to be added to increase the productivity of oil and gas pipelines is a topic discussed in many works by various authors. [63,66,70]The main difference in our research is that when determining the length of the pipeline, the technical indicators of the gas pipeline system and the technological requirements of the reconstruction are taken into account. On the other hand, the determination of the optimal diameter of the pipeline is based on its dependence on the length parameter obtained from technical and economic aspects.

If we assume that the diameter of the pipeline is equal to the diameter of the gas pipeline's main section and perform calculations based on the provided diagrams (Scheme 3.4.1), it becomes clear that regardless of the different installation schemes for the pipelines, the gas transmission capacity (productivity) for all three schemes is calculated using the same known principles and is as follows:

$$Q_\ell = 2F\sqrt{\frac{(P_H^2 - P_k^2)\cdot D}{\lambda z RT\,(L+\ell)}} \tag{7}$$

$P_H$ and Pk - respectively the pressure at the beginning and end of the gas pipeline section; D - diameter of the gas pipeline;

F - minimum cross-sectional area of the gas pipeline;

λ - coefficient of hydraulic resistance; z - compressibility factor of the gas;

T - absolute temperature; R - gas constant;

L - length of the pipeline;

If we use equations (6) and (7) to find the ratio $Q_0 / Q_\ell$,

$$\alpha = \frac{Q_0}{Q_\ell} = \sqrt{\frac{L+\ell}{L}} \tag{8}$$

$\alpha$ - is the coefficient representing the increase in gas transmission capacity (productivity) after the installation of the bypass pipeline to the gas pipeline.

Let's calculate the values for $\ell$=5, 10,15, … ,50 km using equation (8), where L=100 km. Then we can plot the graph of $\alpha$ as a function of $\ell$. Let's proceed with the calculations (**Graph 3.4.1**).

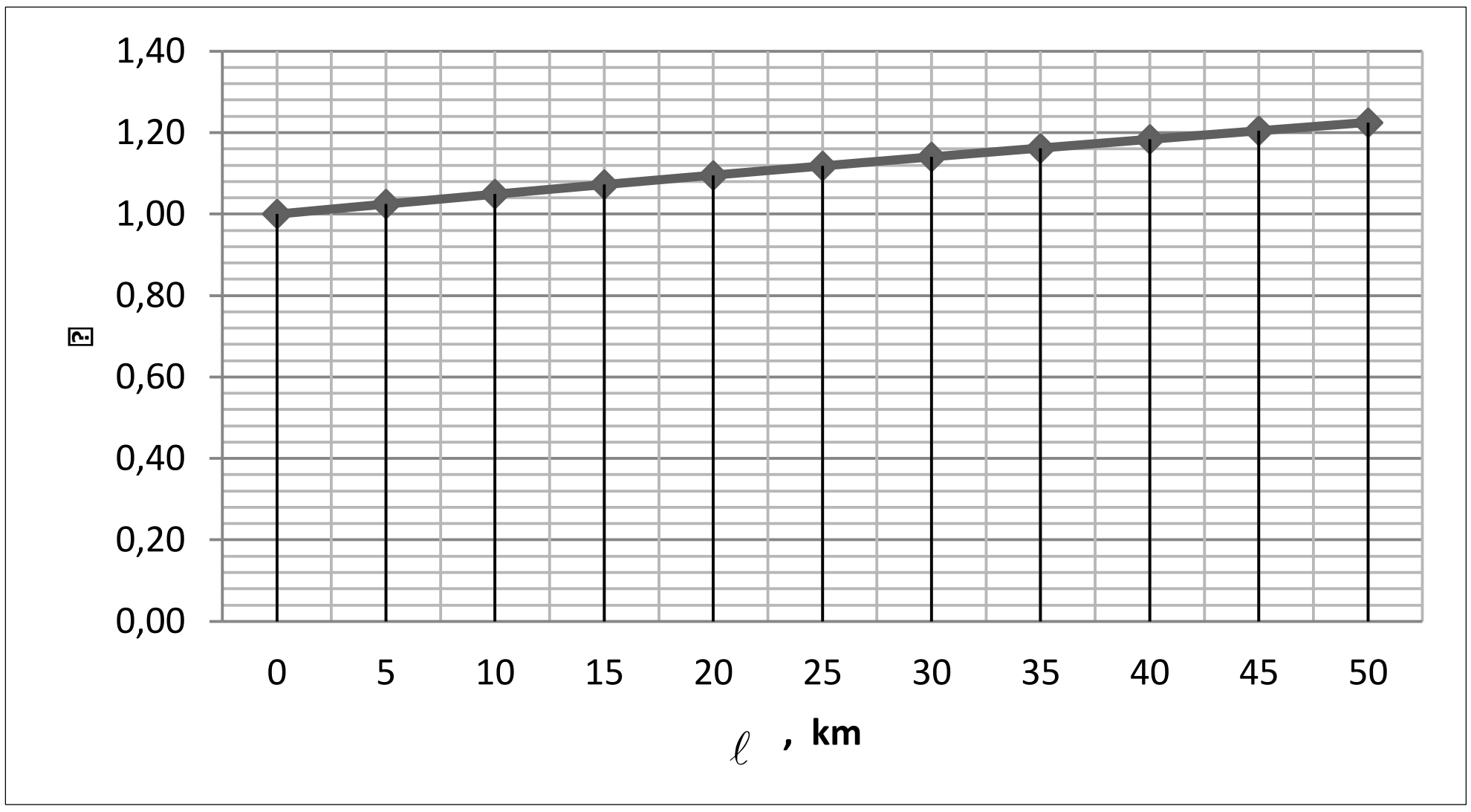

**Graph 3.4.1. The dependence of the gas transmission capacity coefficient during the reconstruction stage of gas pipelines on the length of the installed looping.**

Graph 3.4.1 shows that as the length of the looping increases, the value of the coefficient α also increases. For instance, for a pipeline with a length of 100 km, installing a looping with a length of 20 km is required to increase the gas transmission capacity by 10% ($\alpha$ =1,1). However, to increase the gas transmission capacity by 20% ($\alpha$ =1,2), a looping with a length of 45 km is needed. Therefore, it is evident that as the length of the looping increases, the value of the gas pipeline's gas transmission capacity also increases. This implies that by installing longer loopings along the gas pipeline, we can provide gas to consumers (new consumers) to a greater extent.

As the length of the loopings increases, the amount of expenditure incurred on the pipeline will also increase. Therefore, in the economic aspect of the reconstruction phase of the gas pipelines, we need to find such a length for the looping that the expenditure incurred due to its installation does not exceed the income obtained from the continuous supply of energy carriers to consumer institutions (increasing gas transmission capacity). In other words,

If condition (9) is not met, then the option of designing an additional gas pipeline from the source of the gas pipeline system to serve new consumers is appropriate.

$S_{com}$ -depending on the installed length of the loopline, the income obtained as a result of the provision of the necessary amount of uninterrupted gas fuel to consumer facilities due to the increase in the efficiency of the gas pipeline system is determined as follows:.

$$S_{com.} = (\alpha - 1) \cdot Q_0 \cdot e, \quad \$ \tag{10}$$

Here, the revenue obtained per cubic meter of gas consumption for e-consumers is expressed in dollars per cubic meter ($/m3).

$S_{exp.}$- The value of the additional cost incurred for the reconstruction of the system, depending on the length of the loopline installed on the main section of the gas pipeline and the technical condition of the system, is determined by the following expression.

$$S_{exp.} = C_b(1-F)(1-M) \cdot \frac{\ell}{L}, \ \$ \tag{11}$$

Here,

$C_{b}$- is the restoration cost of the existing gas pipeline during reconstruction, in $;

F is the coefficient reflecting the physical wear and tear of the existing gas pipeline during reconstruction;

M- is the coefficient reflecting the moral depreciation of the existing gas pipeline during reconstruction.

From the analysis of equation (11), it is clear that if the looping length is $\ell = 0$, no additional cost is required for the gas pipeline system, in other words, $S_{exp.} = 0$. On the other hand, if $\ell = L$, an additional cost equal to the construction cost of the gas pipeline is required, which means it is equivalent to the cost of a two-line parallel gas pipeline system.

Taking into account equalities (10) and (11) in equality (9), we obtain the following inequality.

$$\left[ \sqrt{\frac{L+\ell}{L}} - 1 \right] \cdot Q_0 \cdot e \geq C_b(1-F)(1-M) \cdot \frac{\ell}{L} \tag{12}$$

In inequality (12), replace the « ≥» sign with the « = » sign and determine $\ell$.

$$\ell = L \cdot \left[ \frac{2\varphi - 1}{\varphi^2} \right] \tag{13}$$

Here,

φ is the coefficient characterizing the ratio of the cost of reconstructing gas pipelines to the revenue obtained from the increased gas transmission capacity of the system, and it is defined as follows.

$$\varphi = \frac{C_b \cdot (1-F)(1-M)}{Q_0 \cdot e}$$

Thus, to determine the economic and technical efficiency of the length of the installed looping, the following key initial data should be known:

- The initial productivity of the gas pipeline at the reconstruction stage;
- The cost per cubic meter ($m^{3)}$ of gas supplied to new consumers, depending on their designation;
- The restoration cost, physical wear, and moral depreciation values of the gas pipeline at the reconstruction stage.

Once the aforementioned initial data are known, we can determine the economically efficient length of the looping to be installed in the reconstruction of medium-diameter gas pipelines using formula (13).

If φ=1, according to formula (13), $\ell = L$. This means that if the cost incurred for the reconstruction of medium-diameter gas pipelines, depending on the technical condition of the system, equals the revenue obtained from increasing the gas transmission capacity of the system, then the installation of looping along the length of the gas pipeline is economically justified. If φ=0,5, then $\ell = 0{,}938L$, meaning that the installation of looping along 93.8% of the gas pipeline length is technically and economically justified. If φ=0,5, then, $\ell = 0$, meaning that the installation of looping is not justified from a technical and economic point of view. Accordingly, we construct the following graph (Graph 3.4.2) based on the analysis conducted.

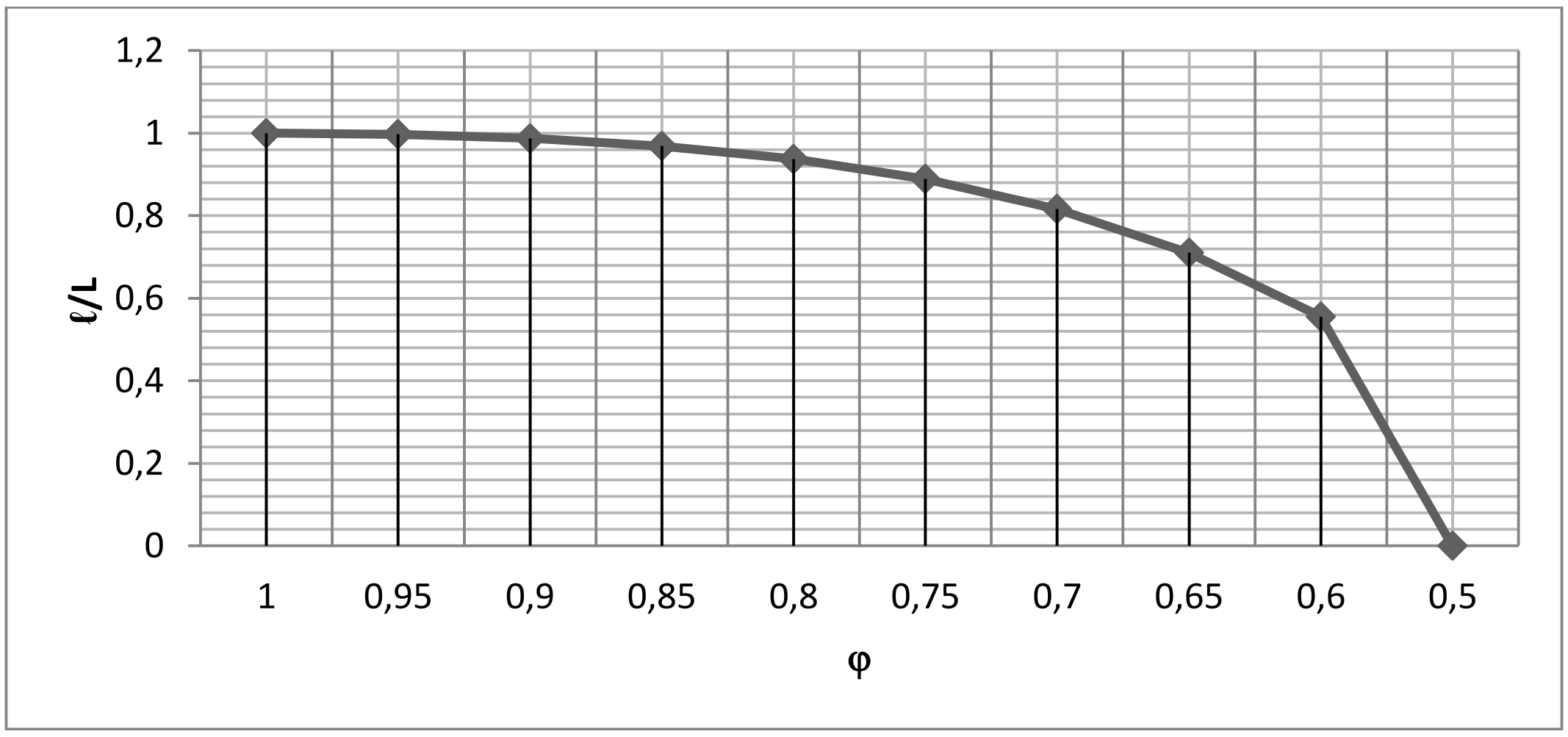

**Graph 3.4.2: Dependence of the economic efficiency indicator ϕ of medium-diameter gas pipelines during the reconstruction stage on the length of the installed looping.**

From the analysis of Graph 3, it is evident that the economic efficiency indicator ϕ\phiϕ lies within the range [1,0.5), indicating that the application of the looping system to increase the gas transmission capacity of the existing gas network for the purpose of supplying gas fuel to new consumers connected to the existing medium-diameter gas pipeline is economically justified. For example, if ϕ=0.65, according to Graph 4, when the length of the gas pipeline is 100 km, $\ell$=28 km.

According to Graph 3.4.2, for $\ell$=28 km, the value of α\alphaα is taken as 1.17. This means that to increase the productivity of the existing gas pipeline by 17% of its volume for the uninterrupted supply of gas to newly connected consumers, looping with a length of 28 km is required. In this case, the effectiveness of the reconstruction will mainly depend on the technical condition of the existing gas pipeline, denoted as [ $C_b \cdot (1-F)(1-M)$ ]. Determining the coefficient ϕ from the formula indicates that the quantity $Q_0$ and its value e of the increased gas consumption, depending on the designation of newly connected consumers to the gas pipelines, can be assumed constant during the reconstruction stage.

In other words, if the technical condition of the existing reconstructed gas pipeline (physical and moral depreciation, remaining service life, etc.) is very poor, installing looping on the existing gas pipeline to supply new consumers with gas may not be technically and economically feasible. This is because connecting new consumers to a gas pipeline that has only 5 years of service life remaining due to its physical depreciation does not meet the technical operational requirements of the system.

Therefore, for values of $\varphi<0.5$, it is appropriate to design an additional gas pipeline from the source of the system to new consumers. For values of physical depreciation $F>0.6$ (more than 60%) of the existing morally worn-out gas pipeline, the renovation option of the scheme should be preferred from the perspective of improving the technical condition of the existing gas pipeline. Considering the new gas consumers, both the linear sections of the worn-out gas pipelines and their distribution networks should be updated to meet modern requirements. For example, currently, there is new equipment and connecting fittings that successfully execute the management operations of gas pipelines in line with modern technological developments. In this case, the increased volume of gas consumption required by the new consumers connected to the existing gas pipeline is achieved by increasing the power or number of the equipment at the compressor station.

Thus, rational reconstruction options based on technical-economic indicators are justified through reports to ensure the necessary (increased) productivity of linear gas pipelines. Rationalization is a calculation method that allows for an evaluative analysis according to the technical-economic criteria of the studied options. The report scheme established the minimal length of loopings required to increase the throughput capacity of various sections of main gas pipelines. As a continuation of the research, let's determine which other parameters of the proposed looping need to be optimized. Specifically, the diameter of the looping that raises the productivity of the gas pipeline section to the required level should be determined relative to the diameter of the existing pipeline to further minimize reconstruction costs. For this purpose, the diameter of the loopings installed on the existing gas pipeline should be optimized.

It should be noted that the diameter of the looping, whose efficient length is determined by formula (13), is assumed to be equal to the diameter of the reconstructed gas pipeline. The fundamental results of the proposed formula are not dependent on the diameter of the looping, meaning the value of $\varphi$ remains unchanged regardless of the chosen diameter. For this purpose, effective calculations need to be carried out to enable the analysis of reconstruction options by applying the optimal diameter of the looping according to technical-economic criteria. For the analysis, we divide the linear section of the gas pipeline along its length into three parts: the first part before the looping is installed (x), the second part with the installed looping ($\ell$), and the third part after the looping is installed (L-x-$\ell$) (Scheme 3.4.2).

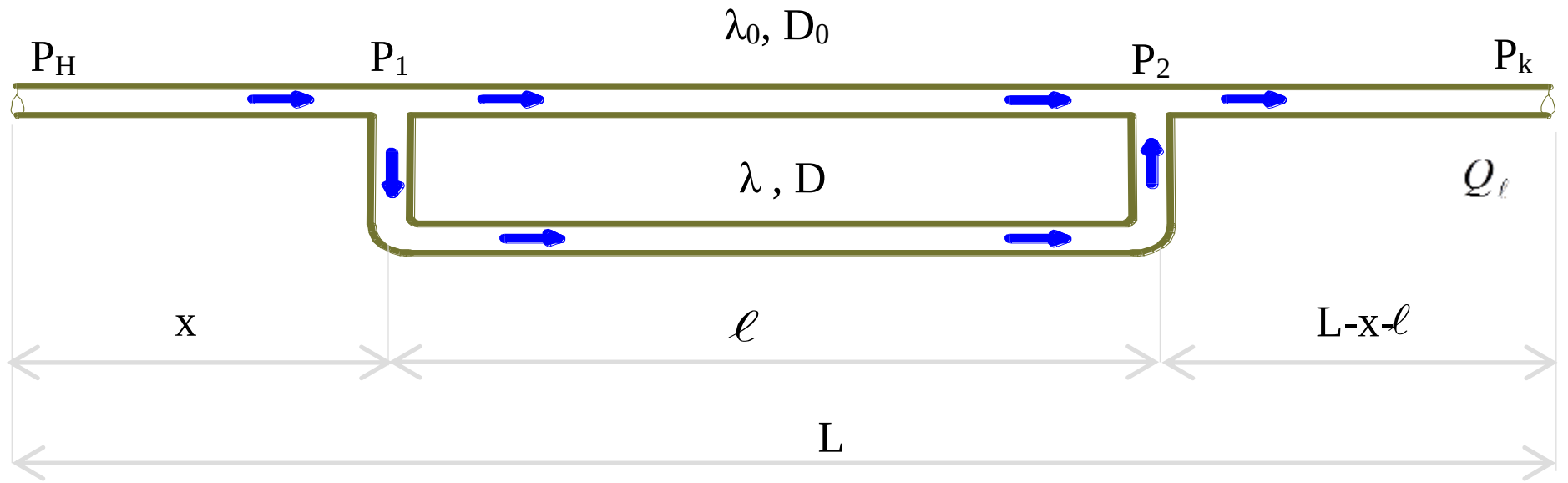

**Scheme 3. 4.2. Diagram of the linear pipeline section with looping Let us present the following notes:**

$Q_\ell$ – the increase in the productivity of the gas pipeline section due to the installation of the looping;

$\lambda_0$ and $D_0$ − the resistance coefficient and diameter of the existing pipe of the operational gas pipeline;

$\lambda$ and D − the unknown values of the resistance coefficient and diameter of the installed looping;

x – the length of the first section of the pipe;

$\ell$ − the unknown length of the looping;

L – the total length of the reconstructed section of the gas pipeline.

Based on Scheme 3.4.2, let's write the equations of stationary gas movement for the first, second, and third sections along the length of the gas pipeline [8].

$$P_1^2 - P_H^2 = c^2 \frac{\lambda_0 x}{D_0^5} Q^2 \tag{14}$$

$$P_2^2 - P_1^2 = c^2 \frac{\lambda_0 \ell}{D_0^5} (1+\beta^2) Q^2 \tag{15}$$

$$P_k^2 - P_2^2 = c^2 \frac{\lambda_0 (L - \ell - x)}{D_0^5} Q^2 \tag{16}$$

Here, $P_н$, $P_к$ – pressure at the beginning and end of the gas pipeline section;
$P_1$, $P_2$ - pressure at the beginning and end of the looping.

$$c = \sqrt{zRT} \qquad\qquad \beta = \sqrt{\left(\frac{\lambda D_0^5}{\lambda_0 D_\ell^{\,5}}\right)}$$

By combining equations (14) to (16), we obtain the following expression.

$$P_k^2 - P_н^2 = \frac{c\lambda_0}{D_0^5} \left[\left(L + \ell\beta^2\right)\right] Q_\ell^{\,2} \tag{17}$$

Based on equation (6), the expression for the stationary gas flow in a straight pipeline without looping is determined by the following statement.

$$P_k^2 - P_н^2 = \frac{c\lambda_0}{D_0^5} L Q_o^{\,2} \tag{18}$$

If we divide equation (17) by equation (18), the ratio can be expressed in a more compact form.

$$\frac{Q_\ell^{\,2}}{Q_0^2} = \frac{L}{L + \ell\beta^2} \tag{19}$$

It is known that the resistance coefficients of the existing operational gas pipeline and the looping installed during the reconstruction phase may adhere to the same expression adopted during the project design. For example, if the operational section of the existing gas pipeline is fouled but subject to square gas flow conditions, the resistance coefficients of the existing and new pipeline sections can be considered equal. Therefore, β can be expressed as follows [3]:

$$\beta^2 = \left(\frac{D_0}{D_\ell}\right)^5 \tag{20}$$

In this case, the expression for the throughput capacity of the linear pipeline section will take the following form:

$$Q = \beta Q_0 = \beta \cdot F \sqrt{\frac{\left(P_н^2 - P_к^2\right)D}{\lambda_0 zRTL}} \tag{21}$$

When calculating gas pipelines with loopings, it is convenient to use consumption coefficients. The gas pipeline shown in Scheme 3 consists of three sections: the first section (x), the second section ($\ell$), and the third section (L-x-$\ell$). The first and third sections are single-line with diameter D0 and their consumption coefficient is $\beta_0$. The second section is a double-line gas pipeline with diameters D0 and D, and consumption coefficients $\beta_0$ and $\beta$, respectively. All three sections of the gas pipeline are connected sequentially. Therefore, according to the scheme, we establish the following equality.

$$\frac{L}{\beta^2} = \frac{\ell}{(\beta_0 + \beta)^2} + \frac{L - x - \ell}{\beta_0^2} + \frac{x}{\beta_0^2} \quad \text{or}$$

$$\frac{L}{\beta^2} = \frac{\ell}{(\beta_0 + \beta)^2} + \frac{L - \ell}{\beta_0^2} \tag{22}$$

According to equation (22), the consumption coefficient β will be determined by the following relationship.

$$\beta = \frac{\beta_0}{\sqrt{1 - \frac{\ell}{L}\left[1 - \left(\frac{\beta_0}{\beta_0 + \beta}\right)^2\right]}} \tag{23}$$

Now let's compare the throughput capacity of the gas pipeline with looping ($Q_\ell$) to the throughput capacity of the gas pipeline without looping ($Q_0$). Based on equation (21) and under the condition that the pressures $P_H$ and $P_k$ are the same before and after installing the looping, we can write that the ratio of throughput capacities is $Q_\ell/Q_0=\beta_0/\beta$. Taking equation (23) into account, we obtain the following expression:

$$\frac{Q_\ell}{Q_0} = \sqrt{1 - \frac{\ell}{L}\left[1 - \left(\frac{\beta_0}{\beta_0 + \beta}\right)^2\right]} \tag{24}$$

Using the method described earlier [71], we determine the gas flow regime in the segmented sections of the gas pipeline. If it is not possible to determine the regime, then the quadratic zone of the turbulent regime is assumed. In that case, equation (24) is established as follows:

$$\frac{Q_\ell^{\ 2}}{Q_0^2} = 1 - \frac{3}{4}\frac{L}{\ell} \qquad (25)$$

Using equations (19) and (25), the coefficient of flow resistance β can be determined.

$$\beta = \frac{L}{\ell}\left(1 - \frac{1}{1 - \frac{3}{4}\frac{L}{\ell}}\right) \qquad (26)$$

If we consider equations (20) and (26), we can derive the following principle for determining the optimal (efficient) diameter of the looping installed during the reconstruction stage of the gas pipeline's linear section:.

$$D_\ell = D\frac{1}{\left[\frac{L}{\ell}\left(1 - \frac{1}{1 - \frac{3}{4}\frac{L}{\ell}}\right)\right]^{0,2}} \qquad (27)$$

The parameter β, which is the main parameter of the principle (13), is determined by equation (27). It is evident from equation (27) that the optimal diameter of the looping can vary significantly depending on its length. However, specific characteristics mentioned in the literature have not been taken into account thus far. Based on the results and norms, it has been proposed to equate the optimal diameter of the looping to 0.8 times the diameter of the existing gas pipeline (0.8D). This proposal is based on the technical condition of the reconstructed existing gas pipeline, which affects the length of the installed looping. To determine this dependency, we use equation (27) to calculate the ratios at different diameters of the looping and plot the corresponding dependency graph (Graph 3.4.3) on the $\ell$ -axis.

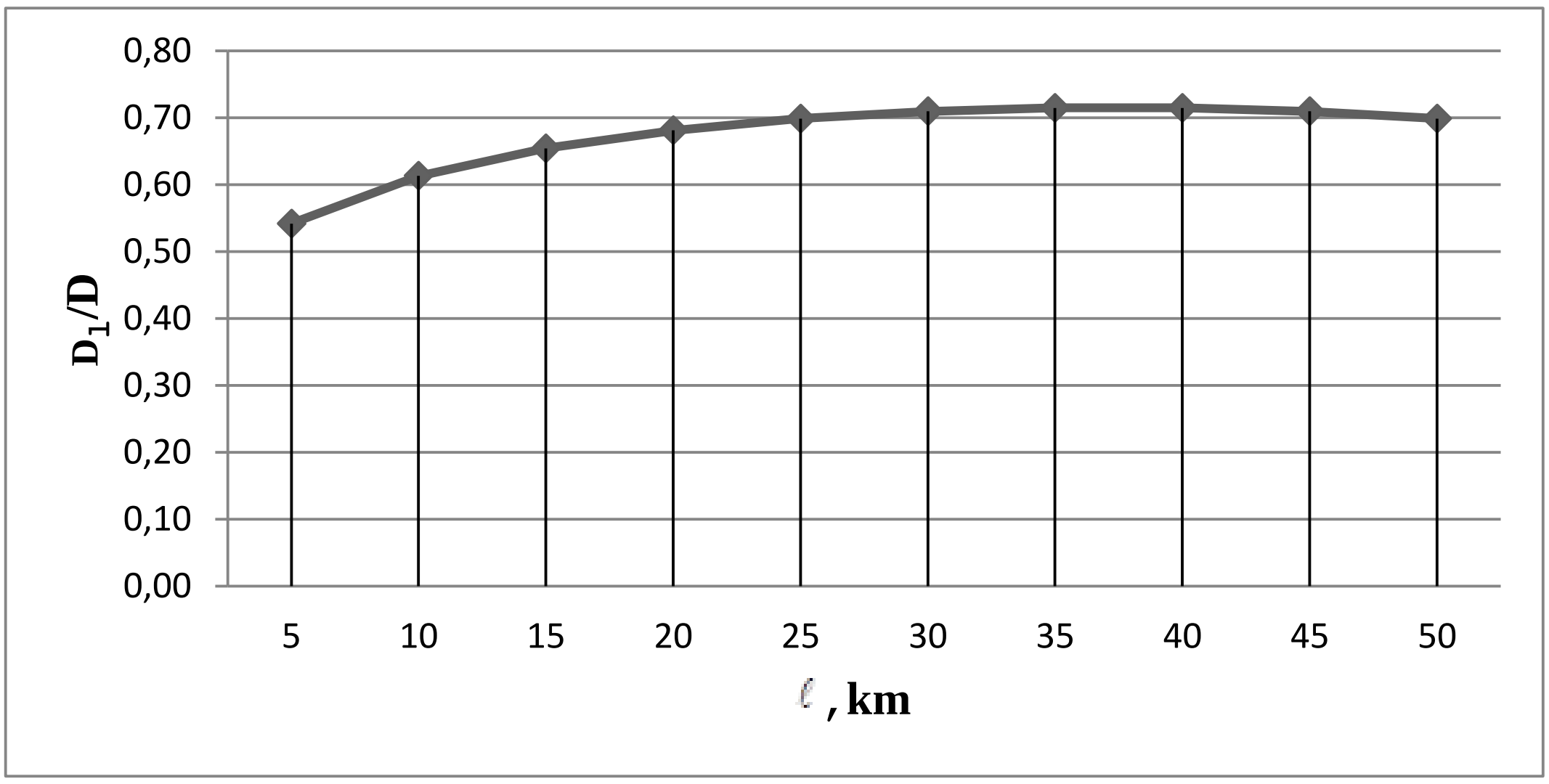

**Graph 3.4.3. Dependence of the ratio of the gas pipeline's linear section diameter to the diameter of the installed looping on the length of the looping.**

It is evident from the analysis of Graph 3.4.3 that the diameter of the looping varies depending on its length. Specifically, for a gas pipeline's linear section of 100 km, the diameter of the installed looping increases from values between 5 km to 25 km in length. However, for lengths $\ell \geq 25$ km, the diameter size appears to stabilize and remains approximately constant. Referring back to Graph 2, it is evident that for the linear section of the existing gas pipeline with a length of 100 km, to increase its productivity by 2%, installation of a looping with a length of $l$=5 km is required.

On the other hand, to increase the productivity of the existing gas pipeline by 12%, it is necessary to install a 25 km long sleeve. Therefore, to achieve a 2% increase in the productivity of the gas pipeline, it is economically viable to install a sleeve with a diameter equal to 0.55 times the diameter of the existing gas pipeline and a length of 5 km.

Analysis of Graph 3 reveals that the economic efficiency factor (ϕ) is 0.51, indicating that the economic benefit exceeds the reconstruction cost approximately twice. To increase the productivity of the gas pipeline by 12%, it is cost-effective to install a sleeve of 25 km in length with a diameter equal to 0.7 times the diameter of the gas pipeline. The economic efficiency factor (ϕ) in this case is 0.54, suggesting an estimated increase in reconstruction costs of 5%.

Installing a sleeve of 50 km in length with a diameter equal to 0.7D of the gas pipeline results in an increase in reconstruction costs of 14% ($\phi = 0.59$). The increase in the length of the installed sleeve contributes to the increase in reconstruction costs.

## Conclusion

1. The report outlines disparities that reflect the comparison of technical and economic parameters of reconstructed gas pipelines. By optimizing these disparities, an assessment scheme has been developed to ensure the cost-effectiveness of reusing existing gas pipeline networks through the installation of telescopic gas sleeves. This defined approach can effectively execute the reconstruction of telescopic gas pipelines from a technical and economic perspective.
2. Theoretical investigations into the principles of reconstructing the linear sections of gas pipeline systems have been conducted and prepared. Research has been dedicated to comparative analysis aimed at determining the optimal parameters of lupines used to increase the efficiency of existing medium-diameter gas pipelines through reconstruction. Based on these analyses, guidelines have been proposed for calculating the optimization parameters of installed lupines to meet consumer demand for gas flow. Theoretical reports have also confirmed the acceptance of reconstruction variants for practical implementation.
3. Therefore, for the effective implementation of the reconstruction of linear sections of main gas pipelines, interdependent relationships between economic, technical, and technological indicators have been identified and comprehensively considered. Simultaneously, to address one of the most critical issues of enhancing the reliability of gas transportation systems, economically viable variants for the reconstruction of linear sections of accepted main gas pipelines have been determined based on their technical and economic conditions.

# Chapter 4. Technological Foundations of Management Decisions in Gas Pipeline Reconstruction

Based on the results of theoretical research, a computational scheme has been developed that allows the establishment of normative-technical documents for the organization and technological decision-making in the exploitation of main gas pipelines. The study extensively explores the dynamic behavior of processes for the supply of modern automatic devices for gas pipelines and the use of an efficient automated control system. The analytical determination of the optimal time for transition processes has been widely applied to ensure the high reliability of pipeline transportation and the most favorable operating conditions for the system. Methods for calculating complex transient processes in main gas pipelines, from non-stationary to stationary modes, have been utilized when gas ingress occurs in a non-stationary regime. A comparison of mathematical expressions for calculating transient processes in complex main gas pipelines has been conducted through theoretical sources.

Based on the results of theoretical research, an optimal strategy for the reconstruction of complex gas pipelines has been selected. The strategy revolves around ensuring the efficient operation regime of the gas pipeline and minimizing the potential damage from accidents to an acceptable level by enhancing the reliability of the elements installed on the pipeline's critical section under operating conditions. A new method for determining the minimum time $t=t_1$ for fixing has been employed. To address such issues, comprehensive calculation schemes have been applied, and modern tools have been utilized. In other words, operational data enabling the analysis of non-stationary operating conditions has been used to prepare calculation schemes aimed at determining the location and timing of pipeline segment ruptures as well as the operational sequence and timing of the installed elements.

In this study, boundary conditions were analyzed to describe mass transfer during the entry and exit of non-uniform gas flow through the connecting pipeline when the main gas pipeline operates in normal (stationary) mode. Additionally, a new calculation scheme has been employed to analyze the dynamics issues in complex main gas pipelines, identifying previously unknown relationships and legal requirements, allowing for a more comprehensive analysis of the problem and providing substantial findings and recommendations. Thus, mathematical expressions have been employed to analytically calculate complex transition processes in main pipeline chambers as a result of the simultaneous entry and exit processes of gas slugs from roadside pipes with varying dynamic parameters into the main pipeline chamber operating in a non-stationary regime. As a result, a cost-effective variant of the sequence of installation of outgoing and incoming pipelines has been identified to enhance the reliability of the pipeline operation.

A highly reliable and efficient technological scheme has been proposed for managing emergency processes in gas transportation based on the principles of the reconstruction phase. For complex gas pipeline systems, new approaches have been investigated for the modernization of existing control process monitoring systems. These approaches are based on modern achievements in control theory and information technology, aiming to select emergency and technological modes. One of the pressing issues is to develop a method to minimize the transmission time of measured and controlled data on non-stationary flow parameters of gas networks to dispatcher control centers. Therefore, the reporting schemes obtained for creating a reliable information base for dispatcher centers using modern methods to efficiently manage the gas dynamic processes of non-stationary modes are of particular importance.

### 4.1. Analysis of Transition Processes in Various Operating Modes for Selecting an Efficient Technological Scheme for Gas Transport Systems

**Introduction:** Currently, researchers in the practice of designing and reconstructing complex gas pipelines typically rely on fundamental principles. Consequently, the motion of the gas flow is assumed to be constant, and based on this principle, the operational mode and other parameters of the pipelines are determined. Subsequently, these parameters are adjusted during operation, taking into account possible errors in the initial data. It is noteworthy that, as indicated by the analyses, the operating modes of main gas pipelines tend to be variable during exploitation.

The accidents occurring in the exploited main gas pipelines cause significant losses to the country's economy and lead to environmental pollution. Therefore, the development and preparation of effective methods for the analysis of non-stationary regimes of natural gas transportation in main gas pipelines undoubtedly arouses scientific and practical interest. Thus, in the conditions of intensive development of gas pipeline transportation, the need for engineering methods to calculate transient processes in complex main gas pipelines becomes particularly relevant. However, as the analysis shows, the solution of dynamic problems in complex main pipelines requires new approaches and solution methods. At the same time, solving dynamic problems in complex main gas pipelines will allow discovering previously unknown relationships and legalities, conducting a more thorough analysis of the problem, and providing substantial results and recommendations. The relevance of this problem has increased in recent years, especially in connection with the ever-increasing requirements for the accuracy of calculations for automated control systems of technological processes. Calculation of transient processes in the main gas pipeline section is necessary for determining the project and exploitation time pressure limit and for developing and preparing automatic protection means in case of exceeding this limit [26,72].

During the analysis of transient processes, calculating the pressure surges of the gas flow in the compressor stations will allow determining the maximum values of pressure rise and the minimum values of pressure drop relative to their steady-state values at compressor stations (CS). Indicators of pressure increase or decrease in the compressor station become known at the moment of activating the automatic protection system and regulating it until the release of protective devices necessary for protection. Furthermore, studying transient processes in complex gas pipelines during complex situations, such as pressure increase or decrease, is crucial for identifying and resolving accidents in a timely manner in the main gas pipelines. This is especially important for the operational dispatch control of automated technological management systems in the automated technological management systems of main gas pipelines [13].

However, one of the challenges in calculating complex transient processes in main gas pipelines during the non-stationary regime, when gas flow transitions from a non-stationary regime to a stationary regime in the pipeline (according to the proposed scheme), is the diverse parameters of the interconnected sections of the pipelines at the moment of gas flow transition. The nature of the challenges lies in the excessive complexity of the mathematical representation of dynamic processes in these systems. The difficulty arises because the distance between connectors or operating valves is identified after the location of the leakage is known, and these processes are described by heat conduction equations. In other words, all processes are described by non-linear differential equations. Computational methods are available for determining the lengths of the steps between connectors, and they are installed depending on the length of the main gas pipeline and the designation of gas consumers [8]. However, the quantity of gas flow transferred from the damaged section to the undamaged section of the pipeline and the law of change are not known in advance and should be determined according to the location of the gas pipeline leakage. Therefore, the methods for calculating transient processes in complex main gas pipelines should be sufficiently developed to allow for both their design and utilization processes, as well as the development and preparation of automatic control means.

In the practice of engineering calculations, the analysis of transient processes in complex main gas pipelines is frequently encountered when the condition of installed valves at various points changes (when opening or closing) [75]. The reason for re-examining this issue is to address a series of problems for the automated control systems of technological processes in pipeline chambers, which result from determining the course of the technological process and the sequence of operation of technological devices. Considering this point, when connecting a parallel gas pipeline to exploitation, the calculation method should be applied to determine the sequence of operation (opening or closing) of valves installed on connectors and pipelines, as well as to determine the pressure value at characteristic points of the parallel main gas pipeline. For this purpose, a new calculation scheme has been proposed for the design and reconstruction of the parallel main gas pipeline that perfectly accommodates automated control systems for technological processes (Scheme 4.1.1).

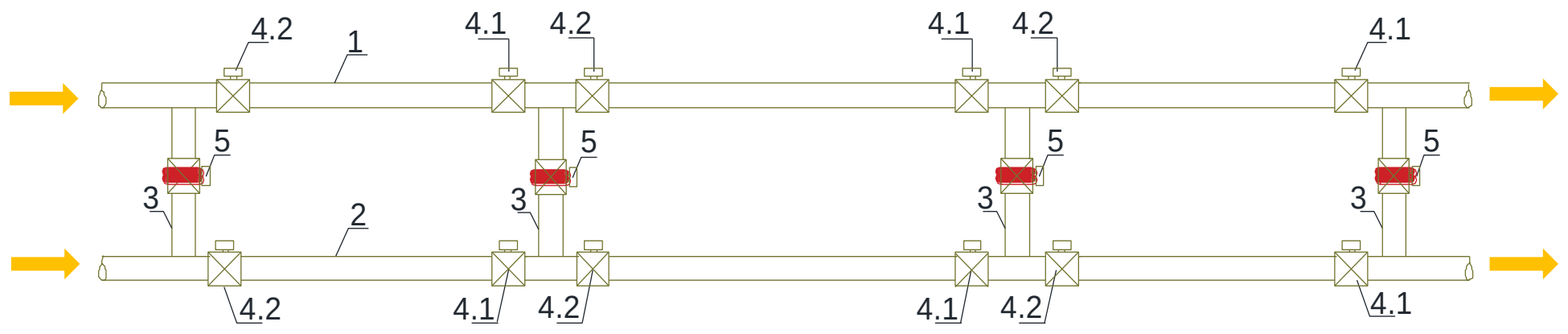

***Scheme 4.1.1. Proposed schematic diagram of the parallel main gas pipeline.***

***1,2 - pipeline chambers on level 1 and level 2; 3 - connectors; 4.1 and 4.2 - automatic valves on pipelines; 5 - automatic valves on connectors.***

It is clear that fittings (valves, drawers, automatic taps, control valves, etc.) are installed on gas pipelines in order to control the flow of gas transported through pipelines, to separate one part of the pipeline from another, and to ensure the operation modes of technological equipment.

The technological aspects of the proposed reporting scheme that distinguish from the scheme of existing complex gas pipelines are as follows:

1. The distance (step) between pipes connecting the parallel pipeline networks of the complex gas pipeline is determined in advance, and the automatic valves installed on these connectors are in a closed position in the normal operating mode (stationary mode) of the gas pipeline. Only when there is a disruption in the integrity of the pipeline networks or during repairs are the automatic valves installed on the connectors activated, i.e., opened, and after the gas pipeline returns to normal operating mode, they are closed again.
2. In the pipeline networks, automatic valves are installed on the right side after the first connector and on the left side before the last connector. It is proposed to install valves on both the left and right sides of the other connectors (Scheme 4.1.1). Thus, in emergency modes, the activation of valves installed on connectors and pipeline networks creates a loop system, resulting in uninterrupted gas supply to consumers from the gas pipeline (Scheme 4.1.2).

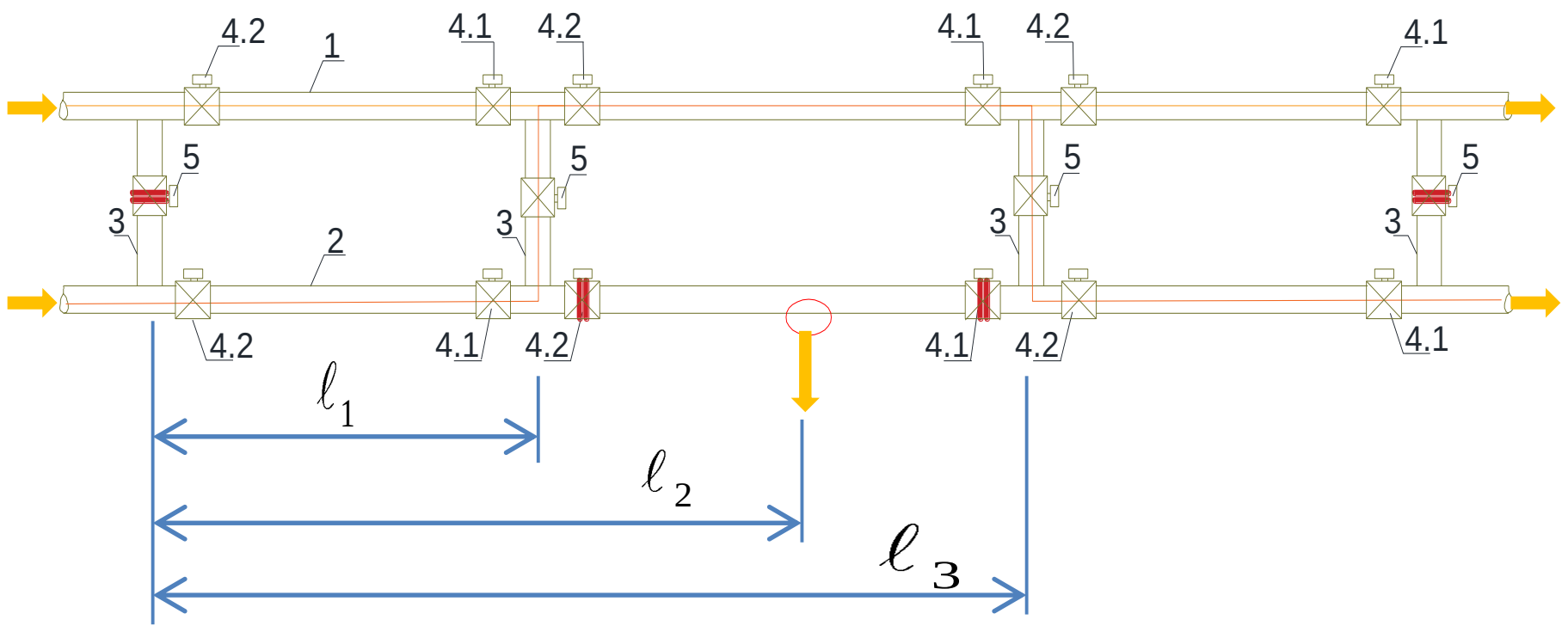

**Scheme 4.1.2.** ***Principal scheme of the proposed technological operating mode for the complex gas pipeline in emergency mode.***

From Scheme 4.1.2, it is evident that in the event of a breach at point x=$\ell_2$ in the gas pipeline, the sequential operation of the technological operating mode for automated control systems of technological processes is implemented as follows. The scheme shows that two automatic valves are placed on each side of a connector, both on the left and right sides. However, in the event of an emergency, only one of these valves should operate. Therefore, the operational principle of these valves should be controlled not based on pressure drop but through the application of automated control systems in the emergency dispatching centers. To achieve this, during an emergency, the location of the automatic valves that will operate in the control center must be known first. Initially, the location of the pipeline integrity breach (leakage point x=$\ell_2$) is identified in the emergency dispatching center. As mentioned, the length of the step between the valves becomes known from the moment the gas pipeline is put into operation. Assuming we know this distance, the following calculation method should be provided to determine the location of the valves ($\ell_3$ and $\ell_1$) installed on the left -4.2 and right -4.1 parts of point x=$\ell_2$ in the control center. To first determine the location of the valve situated on the left side of the emergency point, the following statement is designated:.

$$\ell_3 = \left( \frac{\ell_2}{\ell} + a \right) \cdot \ell \text{, m}$$

Here, $a \in \left]0:1\right[$, İn other words $\left( \frac{\ell_2}{\ell} + a \right)$=n: n=1,2,3,4......

Then, $\ell_1 = \ell_3 - \ell$

Thus, as the leakage point is identified in the emergency dispatch center, the locations of the automatic valves installed on the right and left sides of the damaged pipeline are determined and closed. As a result, the damaged (integrity compromised) section of the pipeline is separated from the main section of the gas pipeline. As a result, the valves on the connectors with known positions next to the 4.2 and 4.1 automatic valves are sequentially opened, allowing the gas flow from the damaged section of the gas pipeline to be directed to the undamaged section. This ensures uninterrupted gas supply to consumers relying on the gas pipeline.

After the closure of the 4.2 and 4.1 automatic valves, the transition process in the gas pipeline is divided into three sections. In the first section ($0 \le x \le \ell_1$), a filling process occurs in the pipeline, and the pressure value along the pipeline increases over time. The pressure increase affects the compression ratio of the gas compressors at the compressor station and can lead to a disruption of their operating mode (Scheme 4.1.3).

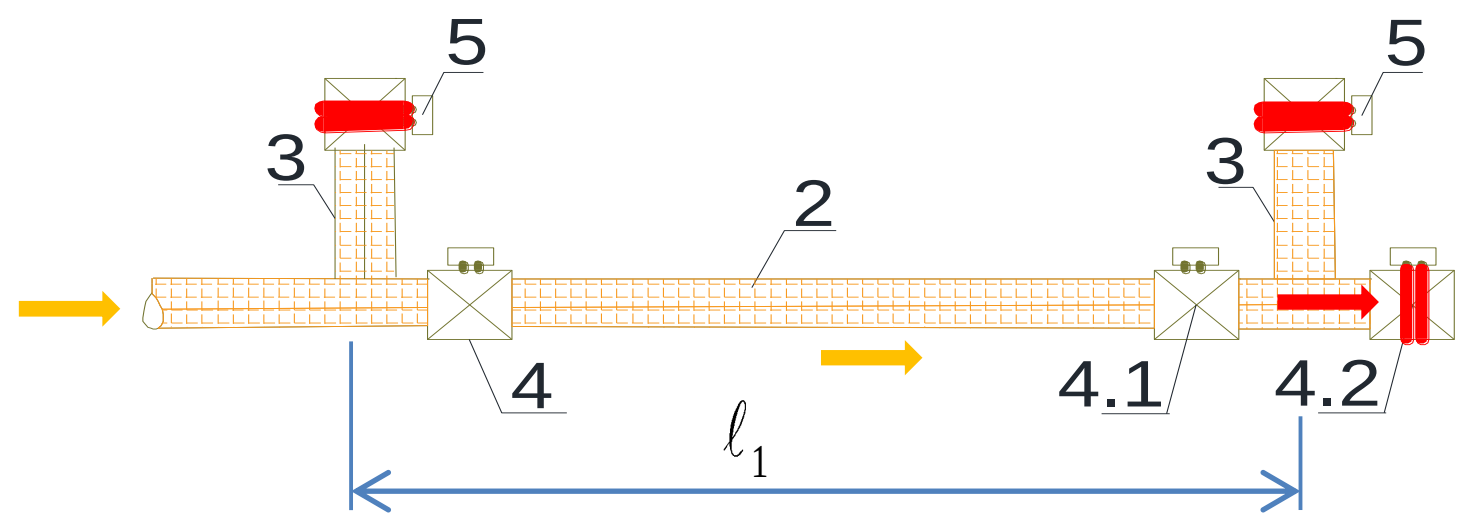

**Scheme 4.1.3. Principal scheme of the technological operating mode in the first section of the complex gas pipeline during the accident regime.**

In order to prevent the progression of this process, it is essential to channel the compressed gas flow to the undamaged pipeline. For this purpose, the displacements installed on the connectors should be activated (opened) at the specified time t=t2, and the gas flow from the damaged pipeline should be redirected to the parallel pipeline's undamaged pipeline. Determining the time t=t2 becomes a critical issue at this point. In the second part ($\ell_1 \leq x \leq \ell_3$), a venting process occurs at the point of rupture of the pipeline's integrity (x=$\ell_2$) due to the gas escaping into the atmosphere (Scheme 4.1.4).

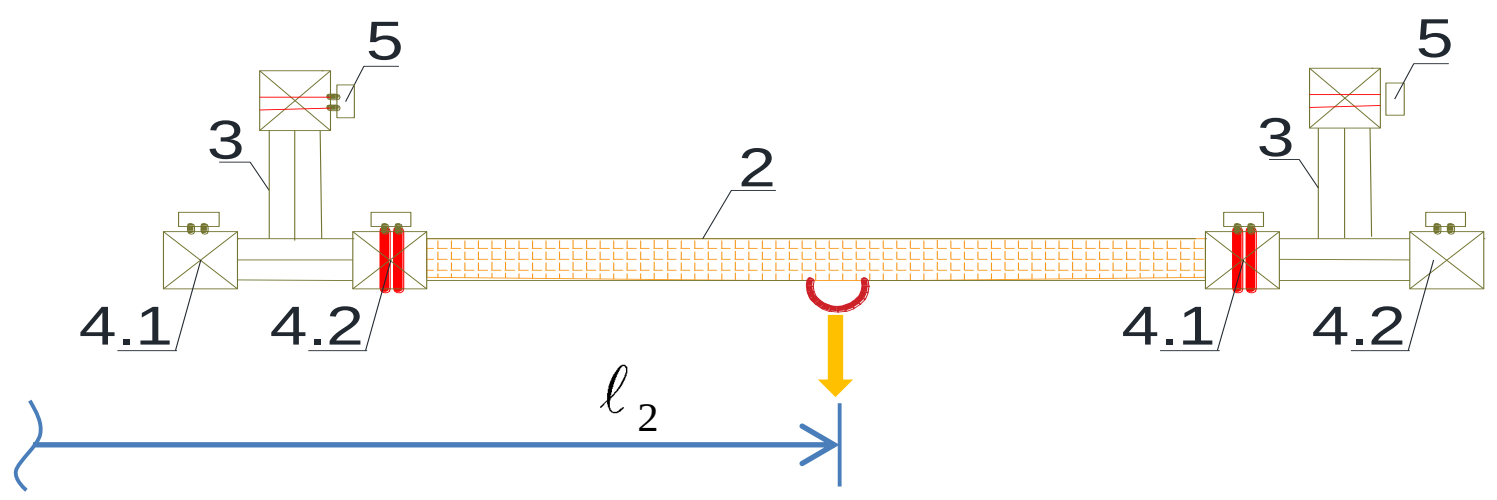

**Scheme 4.1.4. The principal scheme of the technological operation mode in the second part of the complex gas pipeline in the accident mode.**

In this case, it is not possible to prevent the loss of the escaped gas. This is because both valves (sleeves) at the starting and ending points of that part (4.2 and 4.1) are closed. Consequently, due to the closure of these sleeves, the damaged section separates from the main part of the pipeline. At the same time, repair (welding) of the section with a damaged seal involves hazardous gas operations, so it must be carried out in accordance with the relevant regulations and instructions.

In the 3rd section of the pipeline ($\ell_3 \leq x \leq L$), the condensation process occurs due to gas consumption by consumers (Scheme 4.1.5). As seen in the scheme, the closure of the valve at the x= $\ell_3$ point of the damaged pipeline and the valves on the connections -3 result in a decrease in the pressure along the section over time, leading to a reduction in the gas stream's purity.

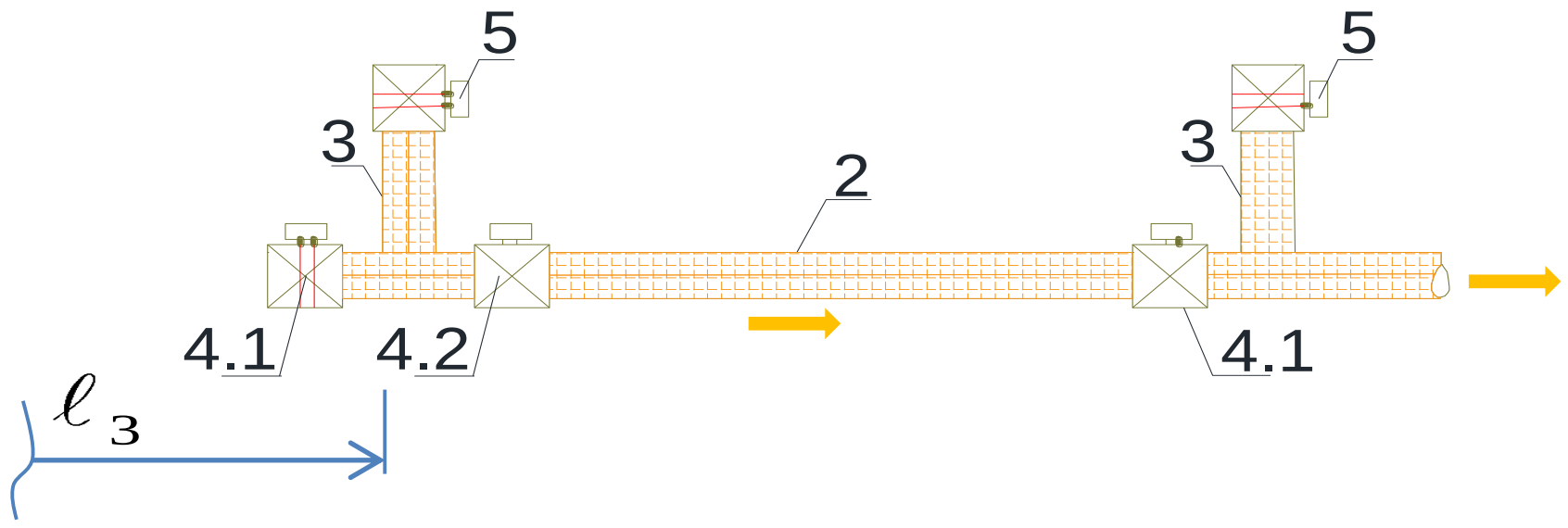

**Scheme 5. Principal Scheme of the Technological Operating regime in the 3rd section of the Complex Gas Pipeline in emergency mode.**

The decrease in gas consumption at the end of the pipeline leads to a decrease in the reliability of gas supply to consumers. In other words, consumers cannot be supplied with the necessary gas pressure. After a certain time, the gas supply regime for consumers can be completely disrupted. At this point, the drawers (taps) at x= $\ell_1$ and x= $\ell_3$ should be activated at the -3 position, and the automatic valves on the -5 should be opened, allowing the required portion of gas from the undamaged pipe belt to flow into the undamaged third part of the damaged pipe where the accident occurred. Thus, consumers in the pipeline are continuously supplied with gas. At t = $t_2$, it is essential to activate the drawers (taps) installed on the connectors.

To assess the feasibility of the proposed scheme in the event of an accident, it is necessary to determine its dynamic state in all processes according to the operating principle. In other words, providing a mathematical solution for the upcoming physical processes involves determining the pressure of gas along the pipeline over time. For the mathematical solution of the processes, we use the following transfer functions [26,60].

For 1st part ($0 \le x \le \ell_1$).

$$\frac{\partial^2 \tilde{P}_1}{\partial X^2} = \frac{2a}{c^2}\frac{\partial \tilde{P}_1}{\partial t} \tag{1}$$

For 2nd part ($\ell_1 \le x \le \ell_3$).

$$\frac{\partial^2 \tilde{P}_2}{\partial X^2} = \frac{2a}{c^2}\frac{\partial \tilde{P}_2}{\partial t} + 2a\tilde{G}_{ut}(t)\delta(x-\ell_2) \tag{2}$$

For 3rd part ($\ell_3 \le x \le L$).

$$\frac{\partial^2 \tilde{P}_3}{\partial X^2} = \frac{2a}{c^2}\frac{\partial \tilde{P}_3}{\partial t} \tag{3}$$

Here, $\tilde{G}$ =ρv; $\tilde{P}$ =ρ·c$^2$

The mass flow measurements at the leakage point at various times represent the values of mass expenditures for the failure mode of the $\tilde{G}_s(t)$ pipeline belt, $\frac{Pa \times \sec}{m}$.

$\delta(x-\ell_2)$ - Dirac's Function

c- For an isothermal process, the speed of sound propagation of gas is expressed in m/s

$\tilde{P}$ - Pressure at the narrowest points of the gas pipeline, Pa.
x-coordinate of the gas pipeline along its length, in kilometers.
t - time coordinate, seconds.
$\lambda$ - hydraulic resistance coefficient.
$\rho$ - average density of gas, kg/m$^3$.
v - average velocity of the gas flow, m/s.
d - diameter of the gas pipeline, meters.
The distribution of the function at the initial conditions for the physical process under consideration is given at the starting time t = 0. Since the process we are examining is a transient process, we will consider the connection of automatic valves installed at the -4.2 and -4.1 points on both sides of the x=$\ell_2$ point as the initial condition. If these valves are connected at time t=$t_1$, then the distribution of the sought function (pressure) at t=$t_1$ is given without taking into account inertial forces.

For t=$t_1$

$$\tilde{P}_i(x,t_1) = P_i(x,t_1) \rightarrow i = 1,2,3$$

The boundary conditions provide the initially known law-compliant variations in expenditure at four points along the pipeline during the considered time of the process. Proper specification of the boundary conditions completes the mathematical model of the process and enables a comprehensive and precise investigation of the occurring physical events

at the point of $$\frac{\partial \tilde{P}_1(x,t_1)}{\partial x} = -2a\tilde{G}_0(t)$$

at the point of $$\begin{cases} \dfrac{\partial \tilde{P}_1(x,t_1)}{\partial x} = 0 \\ \dfrac{\partial \tilde{P}_2(x,t_1)}{\partial x} = 0 \end{cases}$$ at the point of $$\begin{cases} \dfrac{\partial \tilde{P}_2(x,t_1)}{\partial x} = 0 \\ \dfrac{\partial \tilde{P}_3(x,t_1)}{\partial x} = 0 \end{cases}$$

at the point of $$\frac{\partial \tilde{P}_3(x,t_1)}{\partial x} = -2a\tilde{G}_s(t)$$

Here, $\tilde{G}_0(t)$ and $\tilde{G}_s(t)$ represent the values of mass flow measured at different times at the initial and final points of the gas pipeline, $\frac{Pa \cdot san}{m}$.

When the stability of gas pipelines is compromised in experiments, preference is given to the Laplace transformation method. When referring to Laplace transformation P(x,S)= $\int_0^\infty P(x,t)e^{-st}dt$ . S=a+ib represents a complex number, where , i= $\sqrt{-1}$ is an imaginary number. It is clear that through the application of Laplace transformation, equations (1), (2), and (3) are transformed into second-order differential equations, and their general solutions will be as follows.

$$\tilde{P}_1(x,s) = \frac{P_1(x,t_1)}{S} + c_1 Sh\lambda x + c_2 Ch\lambda x \qquad 0 \le x \le \ell_1 \qquad (4)$$

$$\tilde{P}_2(x,S)=\frac{P_2(x,t_1)}{S}+c_3Sh\lambda x+c_4Ch\lambda x-\beta\tilde{G}_{ut}(S)\int_{\ell_1}^{x}\delta(y-\ell_2)sh\lambda(y-x)dy \qquad \ell_1\le x\le\ell_3 \quad (5)$$

$$\tilde{P}_3(x,S)=\frac{P_3(x,t_1)}{S}+c_5Sh\lambda x+c_6Ch\lambda x \qquad \ell_3\le x\le L \quad (6)$$

Here, $\lambda=\sqrt{\frac{2as}{c^2}}$ , $\beta=\sqrt{\frac{2ac^2}{s}}$

Applying the Laplace transformation to the initial and boundary conditions, considering equations (4), (5), and (6), we can determine the constants $C_1$, $C_2$, $C_3$, $C_4$, $C_5$, and $C_6$. By substituting the values of these constants into equations (4), (5), and (6), we obtain the transformed form of the dynamic behavior of the considered process, in other words, the Laplace transform of the pressure distribution along the pipeline for the studied sections.

$$\tilde{P}_1(x,S)=\frac{P_1(x,t_1)}{S}+\beta\tilde{G}_0(S)\frac{ch\lambda(x-\ell_1)}{sh\lambda\ell_1}+\frac{dP_1(0,t_1)}{dx}\frac{ch\lambda(x-\ell_1)}{S\lambda sh\lambda\ell_1} \qquad 0\le x\le\ell_1 \quad (7)$$

$$\tilde{P}_2(x,S)=\frac{P_2(x,t_1)}{S}-\beta\tilde{G}_{ut}(S)\frac{ch\lambda(\ell_2-\ell_3)ch\lambda(x-\ell_1)}{sh\lambda(\ell_3-\ell_1)}-\begin{cases}0\to\ell_1\le x\le\ell_2\\ \beta\tilde{G}_{ut}(S)sh\lambda(x-\ell_2)\to\ell_2\le x\le\ell_3\end{cases} \quad (8)$$

$$\tilde{P}_3(x,S)=\frac{P_3(x,t_1)}{S}-\beta\tilde{G}_s(S)\frac{ch\lambda(x-\ell_3)}{sh\lambda(L-\ell_3)}-\frac{dP_3(L,t_1)}{dx}\frac{ch\lambda(x-\ell_3)}{S\lambda sh\lambda(L-\ell_3)} \quad \ell_3\le x\le L \quad (9)$$

Based on the solution of the considered problem, since it is more convenient to determine the mathematical expression for the non-stationary flow of the gas stream in each of the three sections, we use the reverse Laplace transformation procedure for equations (7), (8), and (9) to find the original solution for each boundary. Thus, according to the initial and boundary conditions, we obtain the mathematical expressions for the distribution of pressure along the pipeline over time.

$$0\le x\le\ell_1$$

$$\tilde{P}_1(x,t)=P_1(x,t_1)+\frac{c^2}{\ell_1}\int_{t_1}^{t}\tilde{G}_0(\tau)d\tau+\frac{2c^2}{\ell_1}\sum_{n=1}^{\infty}\text{Cos}\frac{\pi nx}{\ell_1}\int_{t_1}^{t}\tilde{G}_0(\tau)\cdot e^{-\alpha_3(t-\tau)}d\tau+$$
$$+\frac{c^2}{2a\ell_1}\frac{dP_1(0,t_1)}{dx}\int_{t_1}^{t}\left[1+2\sum_{n=1}^{\infty}\text{Cos}\frac{\pi nx}{\ell_1}e^{-\alpha_3\tau}\right]d\tau \quad (10)$$

$$\ell_1\le x\le\ell_3$$

$$\tilde{P}_2(x,t)=P_2(x,t_1)-\frac{c^2}{(\ell_3-\ell_1)}\int_{t_1}^{t}\tilde{G}_{ut}(\tau)\left[1+2\sum_{n=1}^{\infty}(-1)^n\text{Cos}\frac{\pi n(\ell_2-\ell_3)}{\ell_3-\ell_1}\cdot\text{Cos}\frac{\pi n(x-\ell_1)}{\ell_3-\ell_1}e^{-\alpha_4(t-\tau)}d\tau\right] \quad (11)$$

$$\ell_3\le x\le L$$

$$\tilde{P}_3(x,t)=P_3(x,t_1)-\frac{c^2}{(L-\ell_3)}\int_{t_1}^{t}\tilde{G}_s(\tau)\left[1+2\sum_{n=1}^{\infty}(-1)^n \operatorname{Cos}\frac{\pi n(x-\ell_3)}{L-\ell_3}\cdot e^{-\alpha_5(t-\tau)}\right]d\tau-$$
$$-\frac{c^2}{2a(L-\ell_3)}\frac{dP_3(L,t_1)}{dx}\int_{t_1}^{t}\left[1+2\sum_{n=1}^{\infty}(-1)^n \operatorname{Cos}\frac{\pi n(x-\ell_3)}{L-\ell_3}\cdot e^{-\alpha_5\tau}\right]d\tau \quad (12)$$

Here, $\alpha_3=\frac{\pi^2 n^2 c^2}{2a\ell_1^{\,2}}$, $\alpha_4=\frac{\pi^2 n^2 c^2}{2a(\ell_3-\ell_1)^2}$, $\alpha_5=\frac{\pi^2 n^2 c^2}{2a(L-\ell_3)^2}$

For $t=t_2$ duration, it is possible to assume $\tilde{G}_0(\tau)=\tilde{G}_0(\tau)=$ constant. In this case, we accept the laws expressed in equations (10), (11), and (12) as well as the given conditions to determine the regularity of the change in gas pressure for each of the three different sections of the complex pipeline's technological process.

$P_b = 14\cdot 10^{-2}$ MPa : $P_s = 11\cdot 10^{-2}$ MPa : $G_0 = 10\ Pa\times sec/m$ .: 2a=0.1 1/sec;
c=383.3 m/sec; L =3 $\times 10^4$ m: $\ell_1 = 1\times 10^4$ m ; $\ell_2 = 1{,}45\times 10^4$ m ; $\ell_3 = 2\times 10^4$ m :

Firstly, assuming the location of the gas pipeline leakage point is $\ell_2 = 1{,}45\times 10^4$ m and the closing time of automatic valves is $t=t_1=300$ seconds, we accept (22) and (23) equations. By using these equations, we calculate values for $P_1(x,t_1)$, $P_2(x,t_1)$ and $P_3(x,t_1)$ at every 5000 m and note them in the table below (Table 4.1.1).

| X, km | 0 | 5 | 10 | 14,5 | 20 | 25 | 30 |
|---|---|---|---|---|---|---|---|
| $P_1$ (x,300) $10^{-2}$ MPa | 13,36 | 12,82 | 12,19 | 11,56 | 11,24 | 10,86 | 10,40 |
| $P_2$ (x,300) $10^{-2}$ MPa | - | - | 12,19 | 11,56 | 11,24 | - | - |
| $P_3$ (x,300) $10^{-2}$ MPa | - | - | - | - | 11,24 | 10,86 | 10,40 |

Using the data from Table 4.1.1, we calculate the values of the $\tilde{P}_1(x,t)$ expression for every 60 seconds and every 5 km for the first section of the complex gas pipeline, noting them in the table below, **Table 4.1.2.**

| x , km | $\tilde{P}_1(x,t)$, $10^{-2}$ Mpa | | | | | | | | | | |
|---|---|---|---|---|---|---|---|---|---|---|---|
| | t=0 | t=60 | t=120 | t=180 | t=240 | t=300 | t=360 | t=420 | t=480 | t=540 | t=600 |
| 0 | 13,36 | 14,13 | 14,58 | 15,02 | 15,46 | 15,91 | 16,35 | 16,79 | 17,24 | 17,68 | 18,13 |
| 5 | 12,82 | 13,22 | 13,67 | 14,11 | 14,55 | 15,0 | 15,44 | 15,89 | 16,33 | 16,77 | 17,22 |
| 10 | 12,19 | 12,47 | 12,91 | 13,36 | 13,8 | 14,24 | 14,69 | 15,13 | 15,57 | 16,02 | 16,46 |

Using Table 4.1.2, we construct the graph of the pressure distribution over time for the starting (x=0), middle (x=5 km), and end (x=10 km) points of the first section of the complex gas pipeline (Graph 4.1.1).

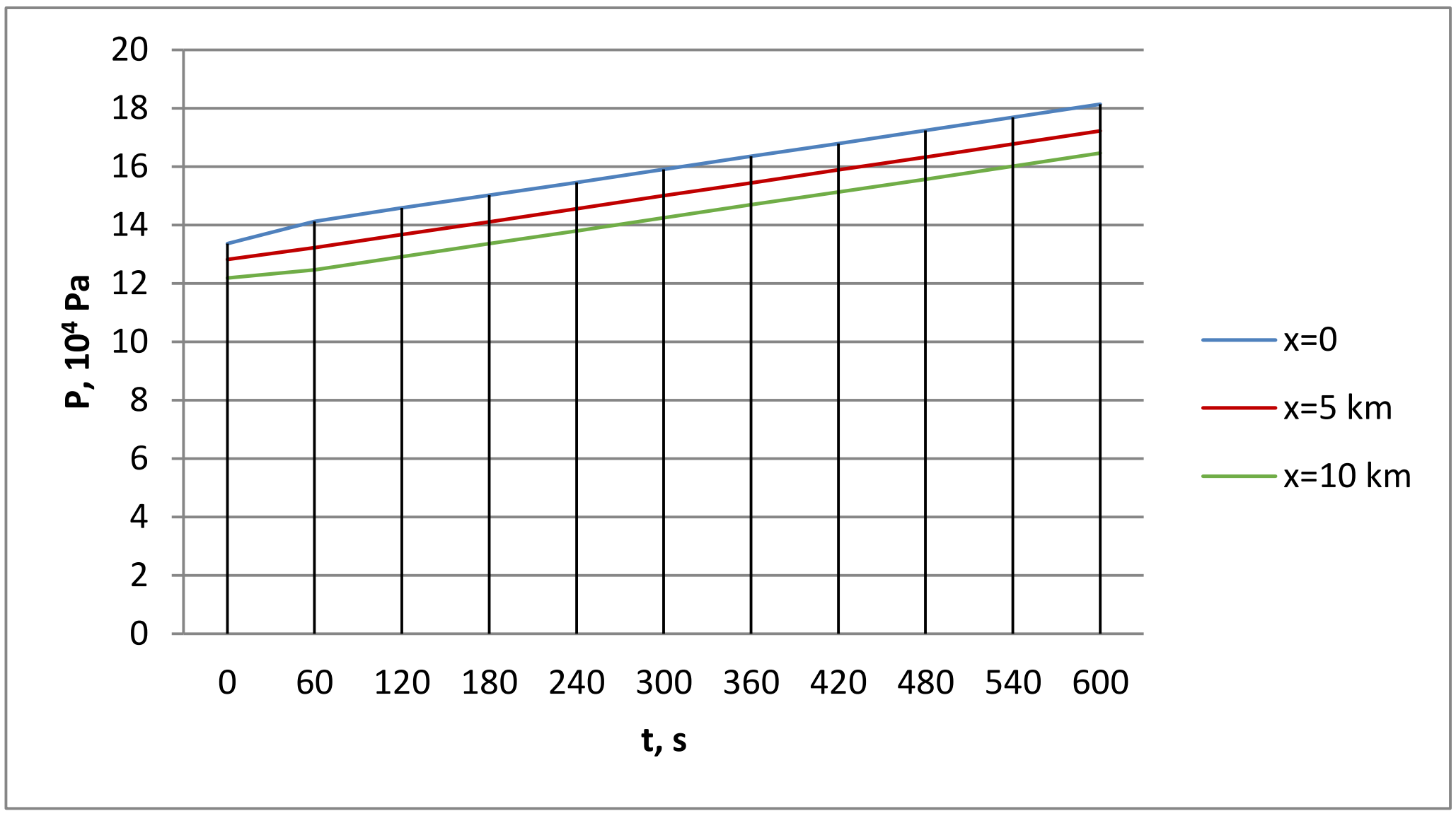

From **Graph 4.1.1**, it is evident that at time t=$t_1$, after the closing of the valves located at x= $\ell_1$ and x= $\ell_3$ points of the compromised pipeline, the pressure along the line x= $\ell_1$ experiences an increase over time. The rate of increase in pressure is calculated at these coordinates. For example, at the starting point x=0, the recorded pressure value is 14.13 × $10^4$Pa at t=60 seconds. This value increases to 15.91 × $10^4$Pa at t=300 seconds and further to 18.13 × $10^4$Pa at t=600 seconds. Thus, the ratio of the increase at t=600 seconds is 1.28. At x=5 km, the pressure recorded at t=60 seconds is 13.22 × $10^4$Pa, and it increases to 17.22 × $10^4$Pa at t=600 seconds, resulting in a ratio of 1.3. Similarly, at x=10 km, the pressure increases from 12.47 × $10^4$Pa at t=60 seconds to 16.46 × $10^4$Pa at t=600 seconds, with a ratio of increase being 1.32. Based on the analysis, it can be observed that the ratio of pressure increase along the axis from the starting point to the endpoint of the first section is increasing.

As mentioned, the compressors used at the beginning of the gas pipeline must be activated at points x=$\ell_1$ and x=$\ell_3$ to protect against the impact of pressure increase. This activation should occur precisely at t=$t_2$, meaning the duration of this time should be determined and justified theoretically.

Continuously, using the data from Table 4.1.1, we calculate the values of the $\tilde{P}_2(x,t)$ expression for the starting point (x=10 km), the leakage point (x=14.5 km), and the endpoint (x=20 km) of the second section of the complex gas pipeline every 60 seconds, and record them in Table 4.1.3.

| x, km | $\tilde{P}_2(x,t)$, $10^{-2}$ Mpa | | | | | | | | | | |
|---|---|---|---|---|---|---|---|---|---|---|---|
| | t=0 | t=60 | t=120 | t=180 | t=240 | t=300 | t=360 | t=420 | t=480 | t=540 | t=600 |
| 10 | 12,19 | 11,77 | 11,32 | 10,87 | 10,43 | 9,98 | 9,54 | 9,1 | 8,65 | 8,21 | 7,77 |
| 14,5 | 11,56 | 11,03 | 10,59 | 10,15 | 9,7 | 9,26 | 8,81 | 8,37 | 7,93 | 7,48 | 7,04 |
| 20 | 11,24 | 10,86 | 10,42 | 9,97 | 9,53 | 9,09 | 8,64 | 8,2 | 7,75 | 7,31 | 6,87 |

Using Table 4.1.3, we create the graph of the distribution of pressure over time for the starting point (x=10 km), leakage point (x=14.5 km), and endpoint (x=20 km) of the second section of the complex gas pipeline (Graph 4.1.2).

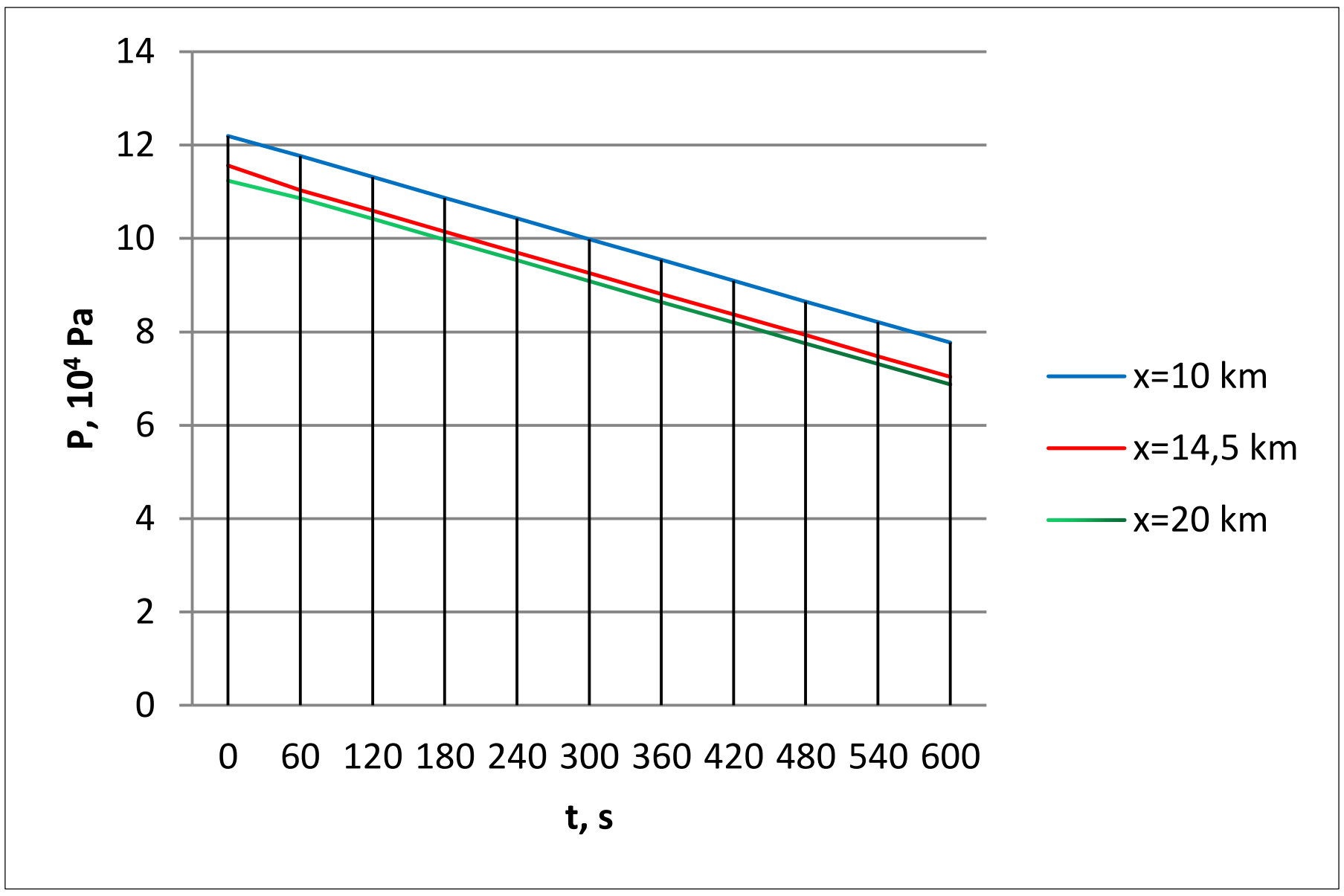

From Graph 4.1.2, it is evident that at t=t1, after the closing of the control valves located at x= $\ell_1$ and x= $\ell_3$, the pressure along the line decreases over time for the second section of the gas pipeline ($\ell_1 \le x \le \ell_3$)The rate of this decrease is higher at the point of valve closure. In this case, at the starting point x=10 km, the pressure, noted at t=60 seconds, decreases from 11.77·10^4 Pa to $9.98\times10^4$ Pa at t=300 seconds, and further decreases to $7.77\times10^4$ Pa at t=600 seconds. In other words, the reduction ratio at t=600 seconds is 1.51 for x=14.5 km. At t=60 seconds, the recorded pressure value is $11.03\times10^4$ Pa, and at t=600 seconds, the pressure value decreases to $7.04\times10^4$ Pa. In other words, the reduction ratio is 1.57 for x=20 km. At t=60 seconds, the recorded pressure value is $10.86\times10^4$ Pa, and at t=600 seconds, the pressure value decreases to $6.87\times10^4$ Pa. It means the reduction ratio is 1.58. Analysis suggests that the reduction ratio of pressure increases in the direction from the starting point to the final point of the second section.

As noted, it is impossible to release gas into the environment as a result of leakage. To determine the mass loss of the gas due to leakage, we perform the following consecutive operations based on equation (8):.

$$\tilde{P}_2(x,S) - P_2(x,t_1) - \frac{c^2}{\ell_3-\ell_1}\tilde{G}_{ut}(S)\frac{ch\sqrt{\frac{2as}{c^2}}(\ell_2-\ell_3)\,ch\sqrt{\frac{2as}{c^2}}(x-\ell_1)}{\frac{sh\sqrt{\frac{2as}{c^2}}(\ell_3-\ell_1)}{\sqrt{\frac{2as}{c^2}}(\ell_3-\ell_1)}} \qquad (13)$$

Here, $\lambda = \sqrt{\frac{2as}{c^2}}$

If we look at the expression in the (13) statement in terms $S \to 0$ (expressed in the transformed form), then we adopt the following rule.

$$\int_{t_1}^{\infty} \frac{\partial \tilde{P}_2(x,t)}{\partial t} dt = -\frac{c^2}{\ell_3 - \ell_1} \int_{t_1}^{\infty} \tilde{G}_{ut}(t) dt$$

Here,

$$\int_{t_1}^{\infty} \tilde{G}_{ut}(t) dt = -\frac{\ell_3 - \ell_1}{c^2} \int_{t_1}^{\infty} \frac{d\tilde{P}_2(\ell_3,t)}{dt} dt \qquad \text{or,} \quad \int_{t_1}^{\infty} \tilde{G}_{ut}(t) dt = -\frac{\ell_3 - \ell_1}{c^2} \tilde{P}_2(\ell_3,t) \qquad (14)$$

Therefore, this gas loss is inevitable. From equation (14), it can be seen that the amount of this loss depends on the distance between displacements. As the distance between them increases, the amount of lost gas loss will also increase.

In conclusion, using the data from Table 4.1.1, we calculate the values of the expression $\tilde{P}_3(x,t)$ at the starting point (x=20 km), the midpoint (x=25 km), and the endpoint (x=30 km) of the third segment of the complex gas pipeline, at intervals of 60 seconds. We record these values in the following( Table 4.1.4).

| x ,km | $\tilde{P}_3(x,t)$, $10^{-2}$ Mpa | | | | | | | | | | |
|---|---|---|---|---|---|---|---|---|---|---|---|
| | t=0 | t=60 | t=120 | t=180 | t=240 | t=300 | t=360 | t=420 | t=480 | t=540 | t=600 |
| 20 | 11,24 | 10,96 | 10,52 | 10,08 | 9,63 | 9,19 | 8,74 | 8,3 | 7,86 | 7,41 | 6,97 |
| 25 | 10,86 | 10,46 | 10,01 | 9,57 | 9,13 | 8,68 | 8,24 | 7,8 | 7,35 | 6,91 | 6,46 |
| 30 | 10,4 | 9,63 | 9,19 | 8,74 | 8,3 | 7,85 | 7,41 | 6,97 | 6,52 | 6,08 | 5,63 |

Using Table 4.1.4, we create a graph depicting the distribution of pressure over time at the starting point (x=20 km), midpoint (x=25 km), and endpoint (x=30 km) of the third segment of the complex gas pipeline (Graph 4.1.3).

From Graph 4.1.3, it is apparent that at time t=$t_1$, after connecting the automatic valves located at points x=$\ell_1$ and x=$\ell_3$ on the depressurized pipeline, the pressure along the length of the gas pipeline's third segment ($\ell_3 \le x \le L$) decreases over time.

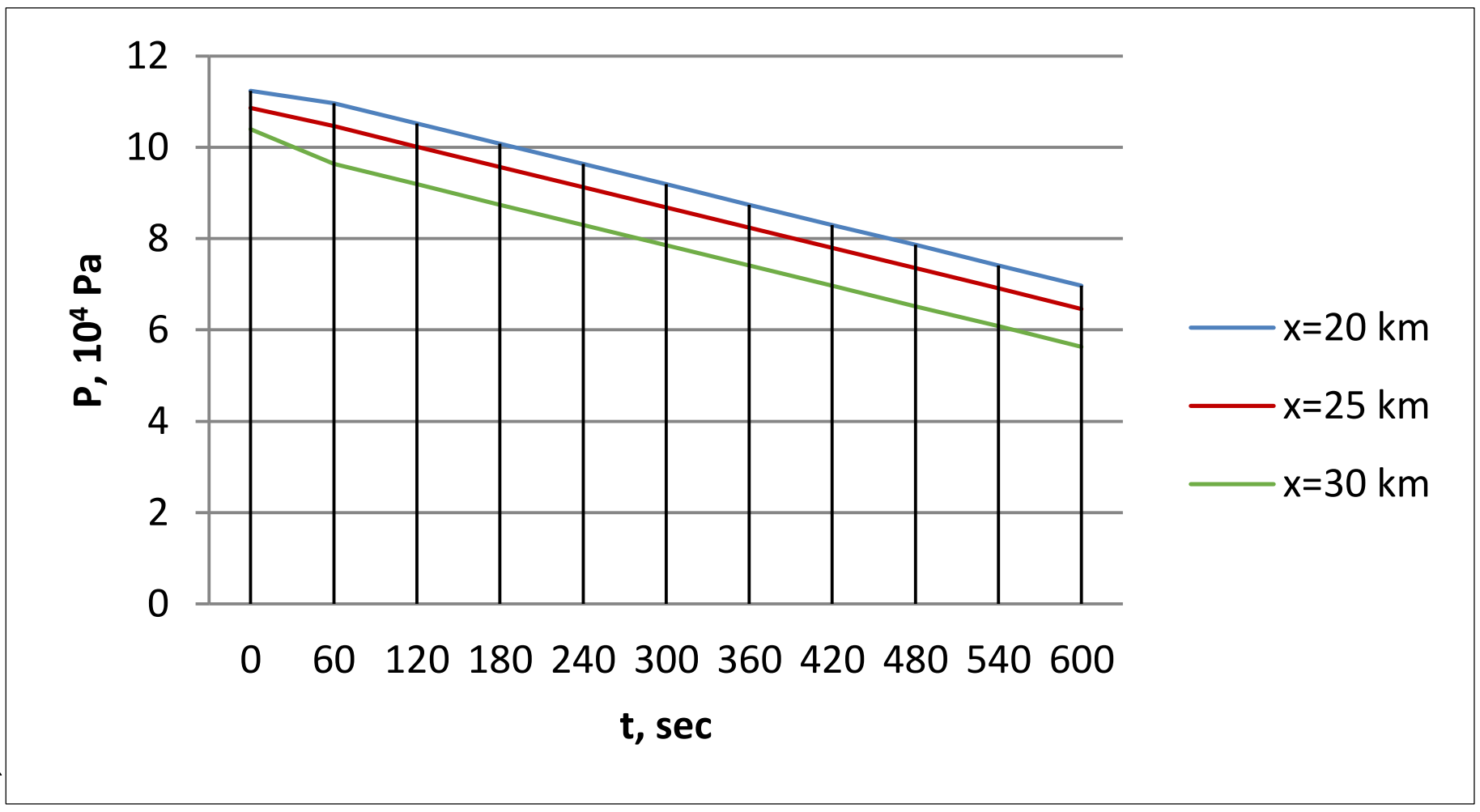

The rate of this decrease is calculated at specific coordinates. Specifically, at the starting point x=20 km, the pressure value decreases from 10.96×$10^4$ Pa at t=60 seconds to 9.19×$10^4$ Pa at t=300 seconds, and further decreases to 6.97×$10^4$ Pa at t=600 seconds. Exactly, at the point

x=25 km, the pressure decreases from $9.63\times10^4$ Pa at t=60 seconds to $6.46\times10^4$ Pa at t=600 seconds, indicating a relative decrease of 1.57 at t=600 seconds. Indeed, at the endpoint x=30 km, the pressure decreases from $10.86\times10^4$ Pa at t=60 seconds to $5.63\times10^4$ Pa at t=600 seconds, resulting in a relative decrease of 1.93. Based on the analysis, we can observe that the relative decrease in pressure increases from the starting point to the endpoint in the direction of the third segment. As noted, during this time, the volume of gas supplied to consumers also decreases. In other words, after t=600 seconds, the volume of gas supplied to consumers will decrease by an average of 48 percent.

In order to improve the gas supply to consumers and restore its continuity, it is essential to activate the automatic valves installed on connectors at points x= $\ell_1$ and x= $\ell_3$ at the specified time t= $t_2$ . The determination and justification of this duration rely on utilizing expression (10) to identify the lawfulness of the pressure change at the starting point of the gas pipeline.

$$\tilde{P}_1(0,t) = P_1(0,t_1) - \frac{2a\ell_1}{3}\tilde{G}_0 + \frac{4a\ell_1}{\pi^2}\tilde{G}_0\sum_{n=1}^{\infty}\frac{e^{-\alpha n^2(t-t_1)} + e^{-\alpha n^2 t} - e^{-\alpha n^2 t_1}}{n^2} \quad (15)$$

The simplified expression of formula (15) is as follows.

$$\tilde{P}_1(0,t) = P_1(0,t_1) + \frac{4a\ell_1}{\pi^2}\tilde{G}_0\left[e^{-\alpha(t-t_1)} + e^{-\alpha t} - e^{-\alpha t_1} - \alpha(t-t_1)(1-C)\right] - \frac{2a\ell_1}{3}\tilde{G}_0 \quad (16)$$

Here, $\alpha = \frac{\pi^2 c^2}{2a\ell_1^2}$ ; "C" is a Euler's formula, and its value is taken as $C = 0,577215$ .

It is evident that one of the general problems in optimizing the efficiency of the gas supply system is the individual and complex analysis of the parameters of pipeline sections during the investigation. One of the optimized indicators during the research of the gas compressor units is the specific power q, which is a unified parameter. If the increase in system efficiency is considered as the $\bar{Q}$ parameter, then the determination of the optimal number n of additional units operating in the compressor station is determined using the following principles $\left(n = \frac{\bar{Q}}{q}\right)$

.

$$F_{\alpha.\varepsilon}(n) = (S_0 + S_1\varepsilon)n^{\frac{2}{3}} + 0{,}5S_2 \cdot Q^{-\frac{1}{2}}(n-\alpha) = 0 \quad (17)$$

As it appears from this equation, the optimal number of operating units depends on the number of spare (reserve) equipment denoted by α and the provided compression ratio denoted by ε. $S_0$, $S_1$, and $S_2$ are parameters for determining the costs incurred by the gas compressor units, and the approximate accuracy, estimated using the variables q and ε, is as follows:

$$\bar{S} = S_0 + S_1\varepsilon + S_2 \cdot q^{-\frac{1}{2}}$$

The analyzed investigation of expression (17) reveals that, if the following condition is satisfied,

$$\varepsilon \geq \frac{S_0}{S_1} \quad (18)$$

In that case, this equation has only one positive root. Analyses of the processed technical and economic parameters of gas compressor units have shown that when ε approaches 1.35 (ε ⟩ 1,35), there is only one optimal value satisfying this condition. As a result, it is clear that the basis of the normal operating mode of the compressor station is the ratio ε, which reflects the number of gas compressor units and their compression degree. Generally, since the main equipment in the exploitation of complex main gas pipelines is the compressor station, significant attention should be paid to the limited value of ε for the protection of their crucial compressor units. As a conclusion of the analyses conducted above, it can be stated that the theoretical determination of the duration $t=t_2$ should be based on the ε ratio. In this case, for the determination of $t_2$, the condition represented by the change in the initial pressure of the gas pipeline we analyzed in stationary mode - $P_b$ and the gas pipeline's first section in the accident mode (after the automatic valves are connected) denoted by $\tilde{P}_1(o,t)$ should be taken into account.

$$\varepsilon \geq \frac{\tilde{P}_1(o,t)}{P_b} \tag{19}$$

If we replace the "≥" symbol with "=", we obtain the following equality for $t= t_2$ in expression (19).

$$\tilde{P}_1(o,t) = \varepsilon P_b$$

In that case, considering expression (16) in relation to expression (19), we derive the following analytical rule for determining the time of activating the automatic valves installed on connectors at points $x=\ell_1$ and $x=\ell_3$.

$$t_2 = t_1 + \frac{\ell_1}{\tilde{G}_0} \frac{\varepsilon P_b - P(0,t_1) - 2a\ell_1\left(\frac{1}{3} - \frac{1}{\pi^2}\right)}{(2-C)} \tag{20}$$

Since the process we are examining is a transition process, we have assumed the initial condition as $t= t_1$.It is evident that, when the valve is opened at point $x=\ell_2$ in the complex gas pipeline, the equations expressing the distribution of pressure for all three segments until the main connection of automatic valves at the left-4.2 and right-4.1 sections are as follows [7].

$$P_{1,2}(x,t) = P_1 - 2aG_0x + \frac{8aL}{\pi^2}G_o\sum_{n=1}^{\infty}\frac{e^{-\alpha_1(2n-1)^2t}}{\left[(2n-1)\right]^2}\text{Cos}\frac{\pi(2n-1)x}{L} - \frac{4aL}{\pi^2}G_o\sum_{n=1}^{\infty}\left[1-(-1)^n\right]\frac{e^{-\alpha_1n^2t}}{n^2}\text{Cos}\frac{\pi nx}{L} -$$
$$-\frac{c^2t_1}{L}G_{ut} - 2aG_{ut}\left[\frac{x^2}{2L} + \frac{\ell_2^2}{2L} + \frac{L}{3} - \ell_2\right] + \frac{4aL}{\pi^2}G_{ut}\sum_{n=1}^{\infty}\text{Cos}\frac{\pi nx}{L}\text{Cos}\frac{\pi n\ell_2}{L}\cdot e^{-\alpha_1n^2t} \tag{21}$$

The solution of the expressions for the infinite path in rule (21) is as follows.

$$\sum_{n=1}^{\infty}\left[\frac{e^{-\alpha_1(2n-1)^2t}}{\left[(2n-1)\right]^2}\text{Cos}\frac{\pi(2n-1)x}{L} - \left[1-(-1)^n\right]\frac{e^{-\alpha_1n^2t}}{n^2}\text{Cos}\frac{\pi nx}{L}\right] = 0$$

If we consider this simplification in expression (21), it becomes as follows.

$$P_{1,2}(x,t) = P_1 - 2aG_0x - \frac{c^2t_1}{L}G_{ut} - 2aG_{ut}\left[\frac{x^2}{2L} + \frac{\ell_2^2}{2L} + \frac{L}{3} - \ell_2\right] +$$
$$+\frac{4aL}{\pi^2}G_{ut}\sum_{n=1}^{\infty}\text{Cos}\frac{\pi nx}{L}\text{Cos}\frac{\pi n\ell_2}{L}\cdot\frac{e^{-\alpha_1n^2t}}{n^2} \tag{22}$$

$$P_3(x,t) = P_1(x,t) - 2aG_{ut}\ell_2 + 2aG_{ut}x \quad (23)$$

Here, $\alpha_1 = \dfrac{\pi^2 nc^2}{2aL^2}$

By using expression (22), we can determine the value of the expression. Based on Table 4.1.1, the value of the expression at t=300 seconds is =13.36×10$^4$Pa.
In that case, considering ε=1.35 and using the provided values and expression (20), we can determine the time for activating the automatic valves installed on connectors.

$$t_2 = t_1 + \frac{\ell_1}{c^2\tilde{G}_0}\frac{\varepsilon P_b - P(0,t_1) - 2a\ell_1\tilde{G}_0\left(\frac{1}{3} - \frac{1}{\pi^2}\right)}{(2-C)} \text{=300+255=555 sec}$$

So, if the proposed scheme and the given parameters are used, and if the valves installed on the damaged pipeline are closed at $t_1$ = 300 seconds when the gas pipeline is depressurized, then the valves on the connectors should be opened at $t_2$ = 255 seconds. If we consider the moment of the accident, this duration $t_2$ will be 555 seconds.

**Results:**

Taking into account the transition duration of the motion equation based on classical methods, the solution of the system of differential equations describing the non-stationary gas flow in a pipeline with a constant minimum diameter has been achieved. The solution, introduced at the initial time of the solution, allows meeting the boundary conditions' requirement of zero equality, enabling obtaining the comprehensive notation of the analytical model. However, it does not limit the scope of the model for sharp changes in gas flow or pressure.

The methodology for calculating transient processes described both by arbitrary boundary conditions (closing of valves) and heat transfer equations in main gas pipelines has been provided. The method for calculating changes in pressure along the pipeline for various operating modes of the studied section of a complex gas pipeline has been developed and prepared.

The proposed methods can be widely used in the design, reconstruction, and creation of automated technological control systems for main gas pipelines. An analytical approach has been employed to determine the time of transition from one process to another in complex gas pipelines. The adopted principle is suitable for the application of automated control systems for technological processes in pipelines and is convenient for engineering calculations.

Examples of calculations and their comparison with theoretical sources allow evaluating the high accuracy of the proposed methods for calculating transient processes in complex main gas pipelines.

## 4.2. Optimization of Non-Stationary Flow Parameters of Reconstructed Complex Gas Pipelines

**İntroduction:** As known from experience, one of the main problems in the operation of gas transportation system facilities is faults that lead to significant gas losses. The consequences of natural gas leakage can include the following:

Contamination of the environment, damage to buildings and structures, injuries to individuals due to ignition or explosion of gas-air mixtures, interruption of gas supply to consumers, and economic losses. For underground trunk pipelines, this is more dangerous because the operational process can continue without visible disturbances until cracks reach a critical size. Subsequently, these cracks enlarge, resulting in the catastrophic failure of the gas pipeline. It is essential to timely detect, eliminate, and minimize the amount of natural gas loss by solving the problems associated with the detection and elimination of gas leaks to ensure the efficient and safe operation of complex gas pipelines.

Therefore, the development and implementation of effective measures for preventing and eliminating accidents in gas distribution networks, including the use of mathematical expressions allowing for accurate calculation of gas losses during accidents, are essential tasks in ensuring the uninterrupted operation of these systems. Accidents in gas pipelines occur as a result of damage to the pipeline and installed elements, corrosion effects, and failure of welded joints, as well as natural and man-made events affecting the composition and movement of the soil. These accidents are accompanied by the appearance of leaks in the pipes, deterioration of joint connections, and the formation of cracks. Gas leakage occurs in this situation, often resulting in ignition.

The formation and development of small cracks eventually lead to material failure. Consequently, the operational lifespan of underground gas pipelines is significantly reduced under pressure. The occurrence and progression of microcracks are influenced by the chemical composition of the soil where the gas pipeline is laid. Aggressive environments that infiltrate the volume of structural elements in gas pipelines lead to changes in their mechanical properties, resulting in alterations in stress-strain conditions and, concurrently, issues with energy supply as well as changes in the functioning of service organs. Gas networks operating under pressure are also at risk of exposure to aggressive environments, which in some cases lead to rapid corrosion and subsequent failure of pipeline sections [19]. Accident prevention primarily relies on adhering to industrial safety requirements and the operational regulations of gas distribution networks. It is crucial to exercise strict control over the integrity of the pipelines and the absence of cracks used for the ecological and economic aspects of the gas pipeline system, ensuring no formation of cracks. Special attention must be paid to ensuring the reliability and safety of pipeline transportation systems during their operation.

Ensuring the reliable and trouble-free operation of complex gas pipelines, as well as meeting the increasing demands for their long-term use and ecological safety, is one of the most important tasks. Developing and implementing effective methods to prevent accidents, reduce their frequency, and eliminate accidents are crucial objectives.

The task of identifying and determining the damage inflicted on pipeline sections as a result of accidents, as well as determining their location, is crucial. Therefore, high demands are placed on gas leakage detection systems for various operating conditions.

- Efficiency (speed);
- High sensitivity;
- Accuracy in determining the location of gas leakage;
- Reliability and accuracy of automatic detection in online mode;

There is no interference with detection modes.

***The methods for determining the location of the rupture point in a pipeline.***

Methods for Detecting Pipeline Leakage Sites: The methods for detecting pipeline leaks are divided into programs and hardware systems. Program systems typically collect information from sensors commonly used in pipeline operation (such as pressure, temperature, and flow sensors) to detect and localize potential leaks based on program algorithms. Hardware-based leakage monitoring systems, on the other hand, utilize sensors that are not typically associated with normal pipeline operations.

There are four main different approaches to solving the leakage control problem [47, 74]:

1. The method of capturing the controlled part of the pressure wave that occurs at the moment of gas leakage is sufficiently sensitive. For instance, the registration period for a change in the nominal flow rate of 1.5–3% is possible within 5–10 minutes, and when determining the lower gas leakage point, the average error in a 100 km pipeline is approximately 500–1000 m. However, this method requires continuous monitoring. There may be a difference between the registration period and the extinguishing time of pressure waves.

2. The hydraulic correlation method for locating gas leaks [52]. The essence of this method lies in the restoration of piezometric levels in the controlled section of the pipeline according to pressure measuring indicators. This allows for the constant registration of leaks at 5–10% of the nominal flow rate and an average error of approximately 1500–2000 meters for a 100-kilometer pipeline section when determining the location of the leak. In periodic mode, it allows for querying pressure sensors one by one, thereby reducing the system's idle value. However, this method requires sufficiently functioning sensors, and during free flow, it may produce incorrect signals due to changes in the rheological properties of the product and the movement of foreign bodies within the pipeline with product flow. Despite this, the correlation principle, due to the presence of reliable hardware and software support with good accuracy in determining the location of gas leaks, is widely used in conjunction with other methods in all modern gas leakage diagnostic technologies [82].

3. Acoustic method for detecting gas leaks. This method is based on the detection of characteristic sound noise in the ultrasonic range generated during the flow of gas under pressure. The acoustic method is sufficiently accurate (1-5 m) and expensive for locating gas leaks. Expensive sensors are installed every 100-300 m along the length of the pipeline, and the equipment for information processing is also costly. Its application has been found during the monitoring of underwater passages in modern projects in the United States [87].

4. The differential method of balancing flow rates in sections equipped with flow meters at the ends of the pipeline. This method allows for the detection of "small" leaks (less than 1%) and is the most accurate of all existing methods [62].

In [18], based on the numerical solution of the mathematical model of pressure changes in the pipeline in the form of a telegraph equation, a system of equations has been obtained that allows for the detection of small leaks under non-stationary flow conditions. Indeed, we also prefer this method, but there has been a different approach to solving the problem. The main distinguishing feature of our research is the use of a reporting scheme that allows for the detection of gas leaks at a fixed moment, i.e., in a short period of time.

When preparing solutions for the given problem, it is advisable to study non-stationary processes and obtain mathematical expressions that can manage the occurring physical processes. To achieve this goal, it is necessary to use non-stationary models. To maintain a compromise between adequacy and simplicity, the requirements for mathematical models compel us to closely monitor the physical processes of gas flow through pipes and necessitate the preparation of specific models or models for each transition process.

Thus, the provision of gas to consumers by complex gas pipelines relies on the productive capacity of gas distribution networks. Taking into account non-stationarity allows for the use of additional connection fittings and connectors in networks designed for maximum duration productivity. Therefore, optimizing the non-stationary flow parameters of reconstructed complex gas pipelines is crucial. This is because one of the components contributing to the reliability of complex gas networks is the connection fittings, junctions, and connectors [10]. Connector fittings are designed to completely close or open the flow of the working environment in the pipeline, depending on the technological requirements. For this purpose, we propose the following scheme for the reconstruction of complex gas networks: (Scheme 4.2.1).

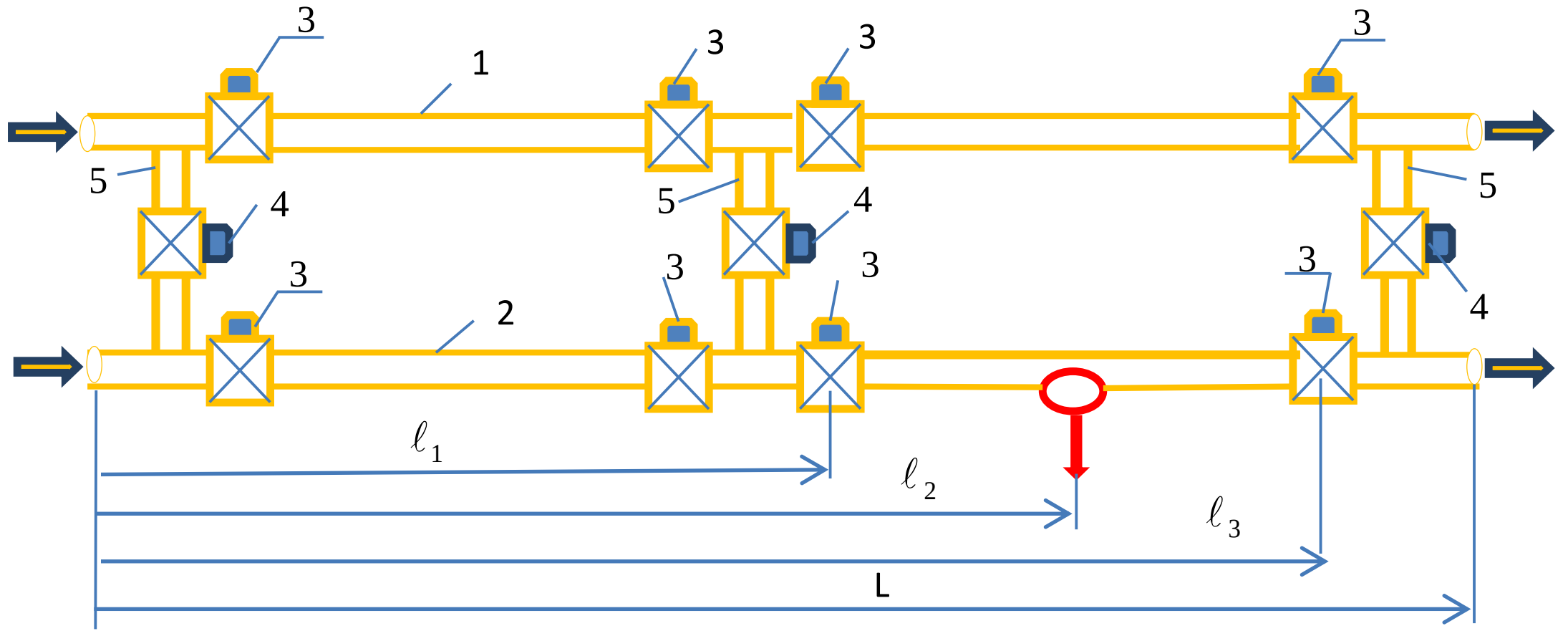

**Scheme 4.2.1: Proposed scheme for increasing the reliability of complex gas networks during reconstruction.**

**1, 2—parallel gas pipeline lines. 3- connection fittings (fittings) serving the pipeline.**
**4-connection fittings serving the connecting pipes. 5- connecting pipes.**
**L- is the length of the parallel gas network.**
**$\ell_2$ -Distance from the starting point of the gas pipeline to the leakage point.**

**$\ell_1, \ell_3$ –Distance from the leakage point to the nearest connection fittings on the left and right, respectively.**

For the exploitation of the proposed gas pipeline scheme in accident mode, it is necessary to determine the transition processes in accordance with the operating principle in a non-stationary regime. For this purpose, we divide the transition processes into three main sections:

The first section is from the starting point of the first gas pipeline to the connection fitting on the left side of the leakage point ( $0 \le x \le \ell_1$ ) , $P_1(x,t)$.

The second section is the damaged portion of the pipeline, between the connection fittings on the left and right sides of the leakage point ( $\ell_1 \le x \le \ell_3$ ), $P_2(x,t)$.

The third section is from the connection fitting on the right side of the leakage point to the end point of the pipeline ( $\ell_3 \le x \le L$ ) , $P_3(x,t)$.

A similar issue was addressed in [3]. In this work, the legality of changes in gas pressure was determined depending on the leakage points of the damaged and undamaged sections of a single hydraulic regime-operated parallel gas pipeline. In parallel gas pipelines operating under a single hydraulic regime, the pipelines extend from common starting to end points. Therefore, when one section of the pipeline is damaged, the other pipelines experience the same incident, and the total gas volume is precisely directed towards the damaged section of the pipeline. However, in the case we are considering, the complex (parallel) gas pipeline operates under independent hydraulic regimes, meaning the operating conditions of the pipelines are not dependent on each other. In this scenario, when an incident occurs in one pipeline, the resulting non-stationary process does not affect the other pipelines.

***Analysis of the mathematical modeling based on the numerical solution of the special equations describing the non-stationary gas flow in the gas pipeline.***

For the non-stationary gas flow, I.A. Charny represented the system of partial differential equations related to gas flow parameters, such as pressure P(x, t) and gas mass flow rate G(x, t), as functions of two variables: time t measured along the axis and distance x. Therefore, for the specific non-linear non-stationary problem under consideration, the specialized system of partial differential equations will be as follows: [69,78 ]

For the first and third sections,

$$\begin{cases} -\dfrac{\partial P}{\partial X} = \lambda \dfrac{\rho V^2}{2d} \\ -\dfrac{1}{c^2}\dfrac{\partial P}{\partial t} = \dfrac{\partial G}{\partial X} \end{cases} \tag{1}$$

For second section

$$\begin{cases} -\dfrac{\partial P_2}{\partial X} = \lambda \dfrac{\rho V^2}{2d} \\ -\dfrac{1}{c^2}\dfrac{\partial P_2}{\partial t} = \dfrac{\partial G_2}{\partial X} + 2aG_{ut}(t)\delta(x-\ell_2) \end{cases} \tag{2}$$

$$\rho = \frac{P}{zRT}$$

Here,

T - average temperature of the gas

z = z(p, T) - compressibility coefficient

c - speed of sound propagation of the gas for an isothermal process, m/second

P - pressure at the ends of the gas pipeline, Pa (Pascal)

x - longitudinal coordinate of the gas pipeline, m (meter)

t - time coordinate, second

λ - hydraulic resistance coefficient

ρ - average density of the gas, kg/ $m^3$ (kilogram per cubic meter)

V - average velocity of the gas flow, m/second

d - diameter of the gas pipeline, m (meter)

G - mass flow rate of the gas

$\delta(x-\ell_2)$ - Dirac function

It is clear that solving the non-linear system of equations (1) and (2) is not possible because these types of equations cannot be integrated. Approximate methods are used in engineering calculations. In other words, equations are linearized. The following linearization methods are used to solve the equations: It has been shown during calculations that Chernov's linearization is more favorable when the pipelines are damaged. The expression characterizing this linearization is as follows: [78].

$$2a = \lambda \frac{V}{2d} \tag{3}$$

Here, 2a represents the Chernov linearization factor, in which the system of equations (1) and (2) will be as follows:

$$\begin{cases} -\dfrac{\partial P_{1,3}}{\partial X} = 2aG_{1,3} \\ -\dfrac{1}{c^2}\dfrac{\partial P_{1,3}}{\partial t} = \dfrac{\partial G_{1,3}}{\partial X} \end{cases} \tag{4}$$

$$\begin{cases} -\dfrac{\partial P_2}{\partial X} = 2aG_2 \\ -\dfrac{1}{c^2}\dfrac{\partial P_2}{\partial t} = \dfrac{\partial G_2}{\partial X} + 2aG_{ut}(t)\delta(x-\ell_2) \end{cases} \tag{5}$$

Here, $G=\rho v$; $P=\rho \cdot c^2$

If we differentiate equation (4) with respect to $x$ and equation (5) with respect to $\dfrac{\partial G}{\partial X}$, the equations will be as follows:

$$\begin{cases} -\dfrac{\partial^2 P_{1,3}}{\partial X^2} = 2a\dfrac{\partial G_{1,3}}{\partial X} \\ \dfrac{\partial G_{1,3}}{\partial X} = -\dfrac{1}{c^2}\dfrac{\partial P_{1,3}}{\partial t} \end{cases} \tag{6}$$

$$\begin{cases} -\dfrac{\partial^2 P_2}{\partial X^2} = 2a\dfrac{\partial G_2}{\partial X} \\ \dfrac{\partial G_2}{\partial X} = -\dfrac{1}{c^2}\dfrac{\partial P_2}{\partial t} - 2aG_{ut}(t)\delta(x-\ell_2) \end{cases} \tag{7}$$

After substituting the expression for $\dfrac{\partial G}{\partial X}$ into the first equation of the system (6) and (7), we obtain the heat transfer equation for solving the mathematical problem.

$$\frac{\partial^2 P_{1,3}}{\partial X^2} = \frac{2a}{c^2}\frac{\partial P_{1,3}}{\partial t} \tag{8}$$

$$\frac{\partial^2 P_2}{\partial X^2} = \frac{2a}{c^2}\frac{\partial P_2}{\partial t} + 2aG_{ut}(t)\delta(x-\ell_2) \tag{9}$$

Thus, to fully determine the given process, the solution of equations (8) and (9) should be sufficient, provided that the initial and boundary conditions are correctly specified. In the initial conditions, the distribution of the sought function (pressure) along the length of the pipe at t=0 is given when inertia forces are not considered. The boundary conditions are provided to ensure the prescribed law-abiding changes in pressure and flow at the endpoints of each segment of the pipe throughout the considered time. Properly specifying the boundary conditions completes the mathematical model of the process and allows for a comprehensive and accurate investigation of the physical events occurring.

$t=0$; $P_i(x,0)=P_b-2aG_0x$: $i=1,2,3$

For the main sections of the gas pipeline considered in the process, we assume the following boundary conditions (where flows are measured at the initial and final points):

At x=0 point $\frac{\partial P_1(x,t)}{\partial X} = -2aG_0(t)$

At $x = \ell_1$ point $\begin{cases} P_1(x,t) = P_2(x,t) \\ \frac{\partial P_1(x,t)}{\partial X} = \frac{\partial P_2(x,t)}{\partial X} \end{cases}$

*At* $X = \ell_3$ point $\begin{cases} P_2(x,t) = P_3(x,t) \\ \frac{\partial P_2(x,t)}{\partial X} = \frac{\partial P_3(x,t)}{\partial X} \end{cases}$

At x= L point $\frac{\partial P_3(x,t)}{\partial X} = -2aG_S(t)$

Here, $G_0(t)$, $G_s(t)$ and $G_{ut}(t)$ respectively represent the mass flows measured at the beginning, end, and leakage point of the pipeline at different times, $\frac{Pa \cdot san}{m}$.

In experiments, the Laplace transformation method is preferred when the integrity of gas pipelines is compromised. When referring to Laplace transformation, P(x,S)= $\int_0^\infty P(x,t)e^{-st}dt$ is considered. Here, S=a+ib represents a complex number, where , i= $\sqrt{-1}$is called the imaginary unit. It is clear that the application of Laplace transformation converts equations (8) and (9) into second-order ordinary differential equations, and their general solutions will be as follows.

$$P_1(x,s) = \frac{P_i - 2aG_0x}{S} + c_1Sh\lambda x + c_2Ch\lambda x \qquad 0 \le x \le \ell_1 \tag{10}$$

$$P_2(x,S) = \frac{P_i - 2aG_0x}{S} + c_3Sh\lambda x + c_4Ch\lambda x - \beta G_{ut}(S)\int_{\ell_1}^{x} \delta(y - \ell_2)sh\lambda(y - x)dy \;, \ell_1 \le x \le \ell_3 \tag{11}$$

$$P_3(x,S) = \frac{P_i - 2aG_0x}{S} + c_5Sh\lambda x + c_6Ch\lambda x \qquad \ell_3 \le x \le L \tag{12}$$

Here, $\lambda = \sqrt{\frac{2as}{c^2}}$, $\beta = \sqrt{\frac{2ac^2}{s}}$

Applying the Laplace transform to the initial and boundary conditions, we can determine the constants $C_1$, $C_2$, $C_3$, $C_4$, $C_5$ and $C_6$ from equations (10), (11), and (12). These constants' values are obtained by considering equations (10), (11), and (12) and analyzing the dynamic state of the gas pipeline, or in other words, the transformed equations represent the distribution of pressure along the pipeline in its altered form for the considered sections.

$$P_1(x,S)=\frac{P_1-2aG_0x}{S}+\frac{\beta G_0}{S}\frac{sh\lambda\left(x-\frac{L}{2}\right)}{ch\lambda\frac{L}{2}}+\beta G_0(S)\frac{ch\lambda(L-x)}{sh\lambda L}-\beta G_s(S)\frac{ch\lambda x}{sh\lambda L}-$$
$$-\beta G_{ut}(S)\frac{ch\lambda(L-\ell_2)ch\lambda x}{sh\lambda L}\qquad 0\le x\le \ell_1 \tag{13}$$

$$P_2(x,S)=\frac{P_1-2aG_0x}{S}+\frac{\beta G_0}{S}\frac{sh\lambda\left(x-\frac{L}{2}\right)}{ch\lambda\frac{L}{2}}+\beta G_0(S)\frac{ch\lambda(L-x)}{sh\lambda L}-\beta G_s(S)\frac{ch\lambda x}{sh\lambda L}-$$
$$-\beta G_{ut}(S)\frac{ch\lambda(L-\ell_2)ch\lambda x}{sh\lambda L}-\begin{cases}0\to \ell_1\le x\le \ell_2\\ \beta G_{ut}(S)sh\lambda(\ell_2-x)\to \ell_2\le x\le \ell_3\end{cases} \tag{14}$$

$$P_3(x,S)=\frac{P_1-2aG_0x}{S}+\frac{\beta G_0}{S}\frac{sh\lambda\left(x-\frac{L}{2}\right)}{ch\lambda\frac{L}{2}}+\beta G_0(S)\frac{ch\lambda(L-x)}{sh\lambda L}-\beta G_s(S)\frac{ch\lambda x}{sh\lambda L}-$$
$$-\beta G_{ut}(S)\frac{ch\lambda\ell_2 ch\lambda(L-x)}{sh\lambda L}\qquad \ell_3\le x\le L \tag{15}$$

Based on the solution of the problem under consideration, since determining the mathematical expression of the physical process of the non-stationary flow of gas in each of the three sections is more appropriate (13), (14), and (15), we use the inverse Laplace transform rule to find the original solution of each equation. Thus, according to the initial and boundary conditions, the mathematical expression of the distribution of pressure along the length of the pipe over time is obtained as follows:

$$P_1(x,t)=P_1-2aG_0\frac{L}{2}+8aG_oL\sum_{n=1}^{\infty}\frac{e^{-\alpha_1 t}}{[\pi(2n-1)]^2}Cos\frac{\pi(2n-1)x}{L}+$$
$$+\frac{2c^2}{L}\sum_{n=1}^{\infty}Cos\frac{\pi nx}{L}\int_0^t\left[G_0(\tau)-(-1)^nG_s(\tau)\right]\cdot e^{-\alpha_2(t-\tau)}d\tau+\frac{c^2}{L}\int_0^t\left[G_0(\tau)-G_s(\tau)-G_{ut}(\tau)\right]d\tau-\qquad 0\le x\le \ell_1$$
$$-\frac{2c^2}{L}\sum_{n=1}^{\infty}Cos\frac{\pi nx}{L}Cos\frac{\pi n(L-\ell_2)}{L}\int_0^tG_{ut}(\tau)\cdot e^{-\alpha_2(t-\tau)}d\tau \tag{16}$$

$$P_{2,3}(x,t)=P_1-2aG_0\frac{L}{2}+8aG_oL\sum_{n=1}^{\infty}\frac{e^{-\alpha_1 t}}{[\pi(2n-1)]^2}Cos\frac{\pi(2n-1)x}{L}+$$
$$+\frac{2c^2}{L}\sum_{n=1}^{\infty}Cos\frac{\pi nx}{L}\int_0^t\left[G_0(\tau)-(-1)^nG_s(\tau)\right]\cdot e^{-\alpha_2(t-\tau)}d\tau+\frac{c^2}{L}\int_0^t\left[G_0(\tau)-G_s(\tau)-G_{ut}(\tau)\right]d\tau- \tag{17}$$
$$-\frac{2c^2}{L}\sum_{n=1}^{\infty}\int_0^tG_{ut}(\tau)\cdot e^{-\alpha_2(t-\tau)}d\tau\cdot\begin{cases}Cos\frac{\pi nx}{L}Cos\frac{\pi n(L-\ell_2)}{L}\to \ell_1\le x\le \ell_2\\ Cos\frac{\pi n\ell_2}{L}Cos\frac{\pi n(L-x)}{L}\to \ell_2\le x\le L\end{cases}$$

$$\text{Here,}\quad \alpha_1=\frac{\pi^2(2n-1)^2c^2}{2aL^2},\qquad \alpha_2=\frac{\pi^2n^2c^2}{2aL^2}$$

### *Analysis and Calculation of Transient Processes in Complex Gas Pipelines*

During the determination of the calculation time for transient processes during accidents in main pipeline networks, it should be considered that the simplification process does not have a significant impact, in other words, a short period is required. For example, the simplification $G_0(t)=G_s(t)=G_0$ is correct for a certain period of time, but as the duration of time increases, this equality is violated and continues until the creation of a new steady-state regime.

Reducing the calculation time also enhances the efficiency of dispatcher stations. This is because promptly identifying the location of a pipeline rupture allows for the activation of shut-off valves, thereby localizing the damaged section of the gas pipeline and reducing gas loss. To achieve this, it is necessary to use calculation programs based on simple principles, adopting simplification methods that have undergone theoretical and experimental tests in known areas of application. Currently, such programs are available [28].

$$G_0(t)=G_s(t)=G_0,\quad G_{ut}(t)=G_{ut}=constat$$

When accepting these simplifications, equations (16) and (17) will be as follows:

$$P_1(x,t)=P_1-2aG_0x+\frac{8aLG_0}{\pi^2}\sum_{n=1}^{\infty}Cos\frac{\pi nx}{L}\frac{e^{-(2n-1)^2\alpha_3 t}}{(2n-1)^2}-\frac{4aL}{\pi^2}G_0\sum_{n=1}^{\infty}((1-(-1)^n)Cos\frac{\pi nx}{L}\frac{e^{-n^2\alpha_3 t}}{n^2}-$$
$$-2aG_{ut}\left(\frac{x^2}{2L}+\frac{\ell_2^2}{2L}+\frac{L}{3}-\ell_2\right)+\frac{4aL}{\pi^2}G_{ut}\sum_{n=1}^{\infty}Cos\frac{\pi nx}{L}Cos\frac{\pi n\ell_2}{L}\frac{e^{-n^2\alpha_3 t}}{n^2} \tag{18}$$

$$P_{2,3}(x,t)=P_1(x,t)-2aG_{ut}(x-\ell_2) \tag{19}$$

Here, $\alpha_3=\dfrac{\pi^2c^2}{2aL^2}$

Using the equations (18) and (19) mentioned above, we accept the following data to determine the lawfulness of the gas pressure variation along the length of the complex pipeline, depending on different locations of gas leakage ( $\ell_2$ =0,5 ·$10^4$ m ; $\ell_2$ =5 ·$10^4$ m ; $\ell_2$ = 9,5·$10^4$ m).

$P_b$ = 55 ×$10^4$ Pa; $P_s$ = 25×$10^4$ Pa; $G_0$= 30 Pa ×sec/m ; $2a=\dfrac{1}{san}$ ; $c=383{,}3\dfrac{m}{\sec}$

L = $10^5$ m: d= 0,7 m: n=1,2,3,…..12.

First, by assuming the location of the gas leakage point in the gas pipeline as $\ell_2$ =0,5 ×$10^4$ m, we calculate the values of the gas pressure in the pipeline every 100 seconds and every 12.5 km (including at the gas leakage point and at the 95th km), and record them in the table below (Table 4.2.1).

| x, $10^3$m | P(x,t), $10^4$Pa | | | | | | | | |
|---|---|---|---|---|---|---|---|---|---|
| | t=100,sec | t=200,sec | t=300,sec | t=400, sec | t=500, sec | t=600, sec | t=700, sec | t=800, sec | t=900, sec |
| 0 | 52,23 | 50,58 | 49,30 | 48,21 | 47,25 | 46,38 | 45,58 | 44,83 | 44,13 |
| 5 | 52,10 | 50,50 | 49,24 | 48,16 | 47,21 | 46,35 | 45,55 | 44,81 | 44,11 |
| 12,5 | 49,72 | 48,34 | 47,18 | 46,18 | 45,27 | 44,44 | 43,67 | 42,95 | 42,27 |
| 25 | 47,11 | 46,29 | 45,45 | 44,64 | 43,88 | 43,16 | 42,48 | 41,83 | 41,21 |
| 37,5 | 43,69 | 43,33 | 42,83 | 42,28 | 42,01 | 41,14 | 40,58 | 40,04 | 39,51 |
| 50 | 40,00 | 39,88 | 39,64 | 39,31 | 38,93 | 38,53 | 38,10 | 37,67 | 37,24 |
| 62,5 | 36,25 | 36,22 | 36,13 | 35,96 | 35,73 | 35,47 | 35,17 | 34,85 | 34,51 |
| 75 | 32,50 | 32,49 | 32,46 | 32,39 | 32,27 | 32,10 | 31,90 | 31,67 | 31,41 |
| 87,5 | 28,75 | 28,75 | 28,74 | 28,71 | 28,64 | 28,54 | 28,40 | 28,22 | 28,01 |
| 95 | 26,50 | 26,50 | 26,49 | 26,47 | 26,42 | 26,34 | 26,21 | 26,05 | 25,86 |
| 100 | 25,00 | 25,00 | 24,99 | 24,97 | 24,93 | 24,84 | 24,72 | 24,57 | 24,37 |

**Table 4.2.1. Values of the pressure of the damaged density gas pipeline along the length depending on time ( $\ell_2 = 0.5 \times 10^4$ m).**

Using Table 4.2.1, let's plot the graph of pressure distribution along the length of the complex gas pipeline as a function of every 200 seconds (Graph 4.2.1).

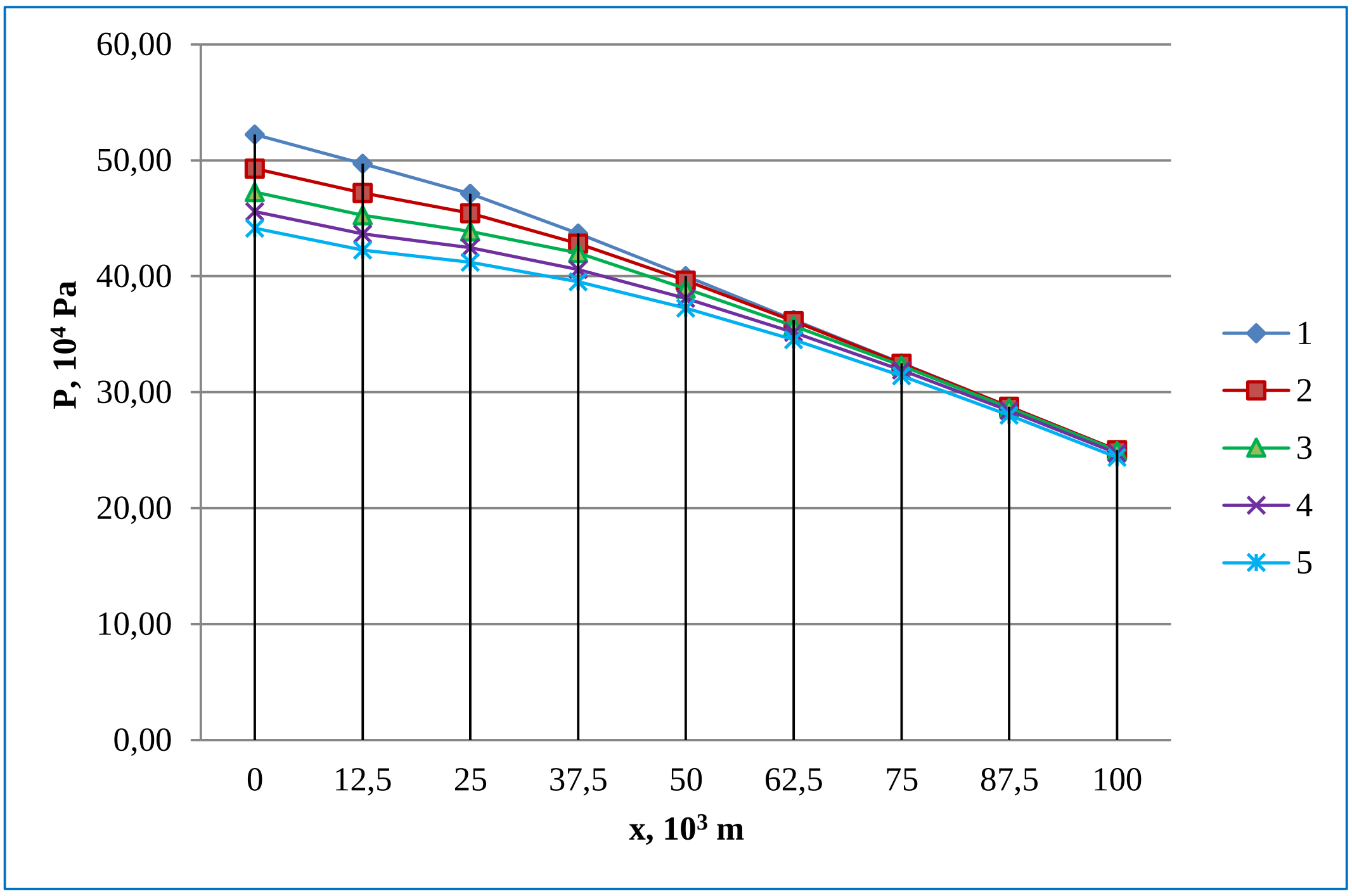

**Graph 4.2.1. Pressure distribution along the length of the complex gas pipeline as a function of time ( $\ell_2$ = 0,5 ×· $10^4$ m, 1-t=100 sec, 2-t=300 sec, 3- t=500 sec, 4- t=700 sec, 5- t=900 sec).**

Sequentially assuming the gas leakage point in the gas pipeline to be $\ell_2 = 5 \cdot 10^4$ meters, the pressure values of the gas pipeline are calculated at intervals of every 100 seconds and every 12.5 kilometers, including the 5th and 95th kilometers, and recorded in the following table (Table 4.2.2).

| x, $10^3$m | P(x,t), $10^4$Pa | | | | | | | | |
|---|---|---|---|---|---|---|---|---|---|
| | t=100, sec | t=200, sec | t=300, sec | t=400, sec | t=500, sec | t=600, sec | t=700, sec | t=800, sec | t=900, sec |
| 0 | 55 | 54,9 | 54,66 | 54,34 | 53,96 | 53,56 | 53,14 | 52,71 | 52,28 |
| 5 | 53,5 | 53,38 | 53,14 | 52,81 | 52,43 | 52,03 | 51,6 | 51,17 | 50,74 |
| 12,5 | 51,23 | 51,05 | 50,76 | 50,4 | 50 | 49,59 | 49,16 | 48,73 | 48,29 |
| 25 | 47,34 | 46,93 | 46,49 | 46,05 | 45,61 | 45,17 | 44,73 | 44,29 | 43,85 |
| 37,5 | 43,05 | 42,35 | 41,76 | 41,23 | 40,75 | 40,28 | 39,82 | 39,38 | 38,93 |
| 50 | 37,95 | 37,1 | 36,45 | 35,89 | 35,38 | 34,9 | 34,44 | 33,99 | 33,54 |
| 62,5 | 35,55 | 34,85 | 34,26 | 33,73 | 33,25 | 32,78 | 32,32 | 31,87 | 31,43 |
| 75 | 32,33 | 31,93 | 31,49 | 31,05 | 30,61 | 30,17 | 29,73 | 29,29 | 28,84 |
| 87,5 | 28,72 | 28,54 | 28,25 | 27,89 | 27,50 | 27,09 | 26,66 | 26,23 | 25,79 |
| 95 | 26,49 | 26,38 | 26,13 | 25,80 | 25,02 | 25,02 | 24,60 | 24,17 | 23,74 |
| 100 | 24,99 | 24,89 | 24,66 | 24,33 | 23,96 | 23,49 | 23,14 | 22,71 | 22,27 |

**Table 4.2.2. Pressure values along the length of the gas pipeline over time ( $\ell_2 = 5\times 10^4$ m).**

Using Table 4.2.2, let's create a graph illustrating the distribution of pressure along the length of the complex gas pipeline every 200 seconds (Graph 4.2.2).

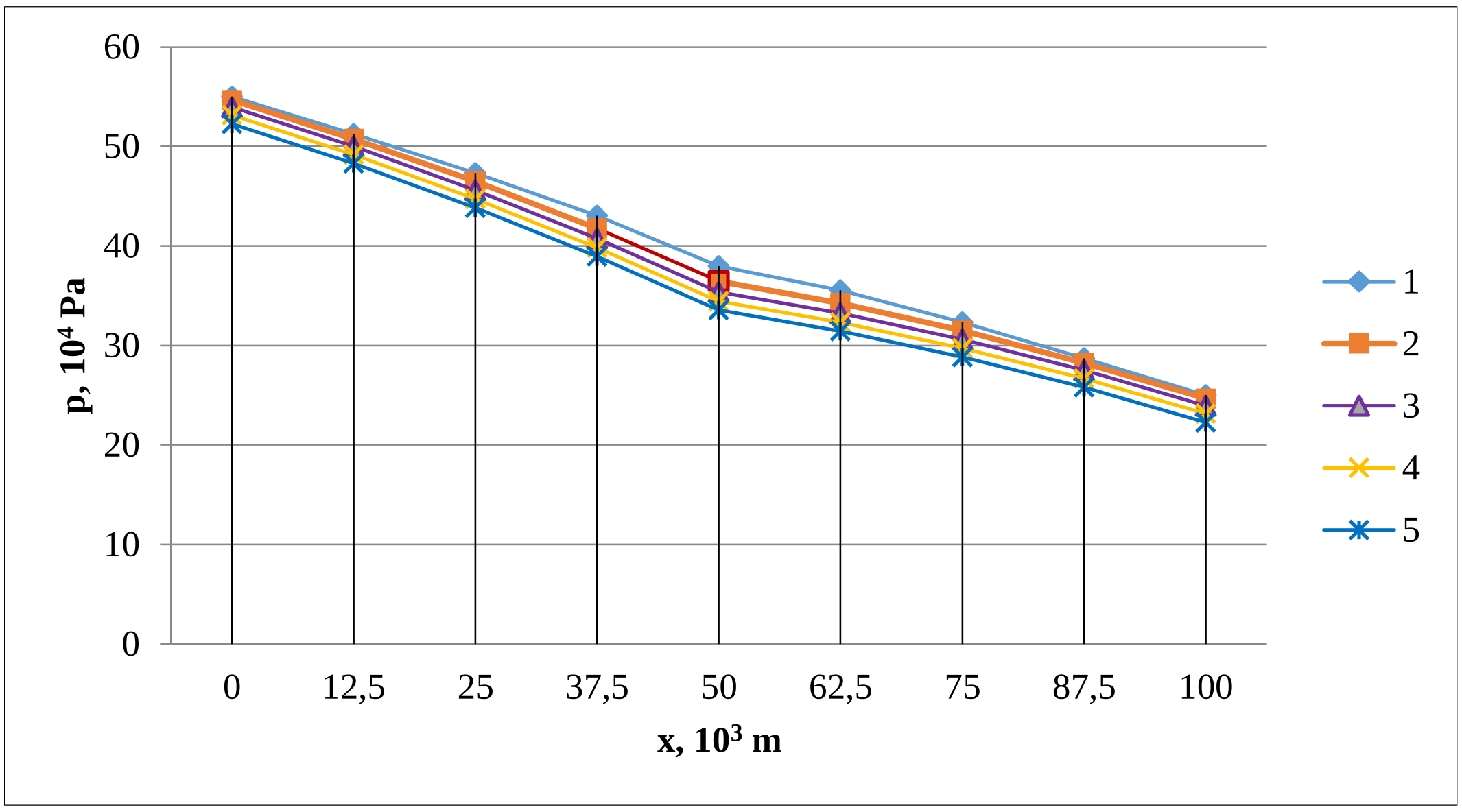

**Graph 4.2.2. Visualization of the distribution of pressure along the length of the complex gas pipeline over time, dependent on the damaged density. ( $\ell_2 = 5\times10^4$ m, 1- t=100 sec, 2- t=300 sec, 3- t=500 sec, 4- t=700 sec, 5- t=900 sec).**

Finally, assuming the location of the gas leakage point at $\ell_2 = 9{,}5 \cdot 10^4$ m, we calculate the pressure values of the gas pipeline every 300 seconds and every 12.5 kilometers (as well as at the gas leakage point and at the 5th kilometer), and record them in the following table (Table 4.2.3).

| x, $10^3$m | P(x,t), $10^4$Pa | | | | | | | | |
|---|---|---|---|---|---|---|---|---|---|
| | t=100,sec | t=200, sec | t=300,sec | t=400, sec | t=500, sec | t=600, sec | t=700, sec | t=800, sec | t=900, sec |
| 0 | 55,00 | 55,00 | 54,99 | 54,97 | 54,93 | 54,84 | 54,72 | 54,57 | 54,37 |
| 5 | 53,50 | 53,50 | 53,49 | 53,47 | 53,42 | 53,34 | 53,21 | 53,05 | 52,86 |
| 12,5 | 51,25 | 51,25 | 51,24 | 51,21 | 51,14 | 51,04 | 50,9 | 50,72 | 50,51 |
| 25 | 47,50 | 47,49 | 47,46 | 47,38 | 47,26 | 47,10 | 46,90 | 46,67 | 46,41 |
| 37,5 | 43,75 | 43,72 | 43,62 | 43,45 | 43,23 | 42,96 | 42,67 | 42,35 | 42,01 |
| 50 | 39,99 | 39,88 | 39,63 | 39,30 | 38,93 | 38,52 | 38,10 | 37,67 | 37,24 |
| 62,5 | 36,18 | 35,83 | 35,32 | 34,77 | 34,20 | 33,64 | 33,08 | 32,54 | 32,00 |
| 75 | 32,10 | 31,28 | 30,44 | 29,64 | 28,88 | 28,16 | 27,48 | 26,83 | 26,21 |
| 87,5 | 27,22 | 25,83 | 24,68 | 23,67 | 22,77 | 21,94 | 21,17 | 20,45 | 19,77 |
| 95 | 23,56 | 21,96 | 20,70 | 19,62 | 18,67 | 17,81 | 17,01 | 16,27 | 15,57 |
| 100 | 22,23 | 20,58 | 19,30 | 18,21 | 17,25 | 16,38 | 15,58 | 14,83 | 14,13 |

**Table 4.2.3. Pressure values of the gas pipeline as a function of time, depending on the length of the damaged gas pipeline. ( $\ell_2$ = 9,5×$10^4$ m).**

Using Table 4.2.3, let's create a graph illustrating the distribution of pressure along the length of the complex gas pipeline every 200 seconds (Graph 4.2.3).

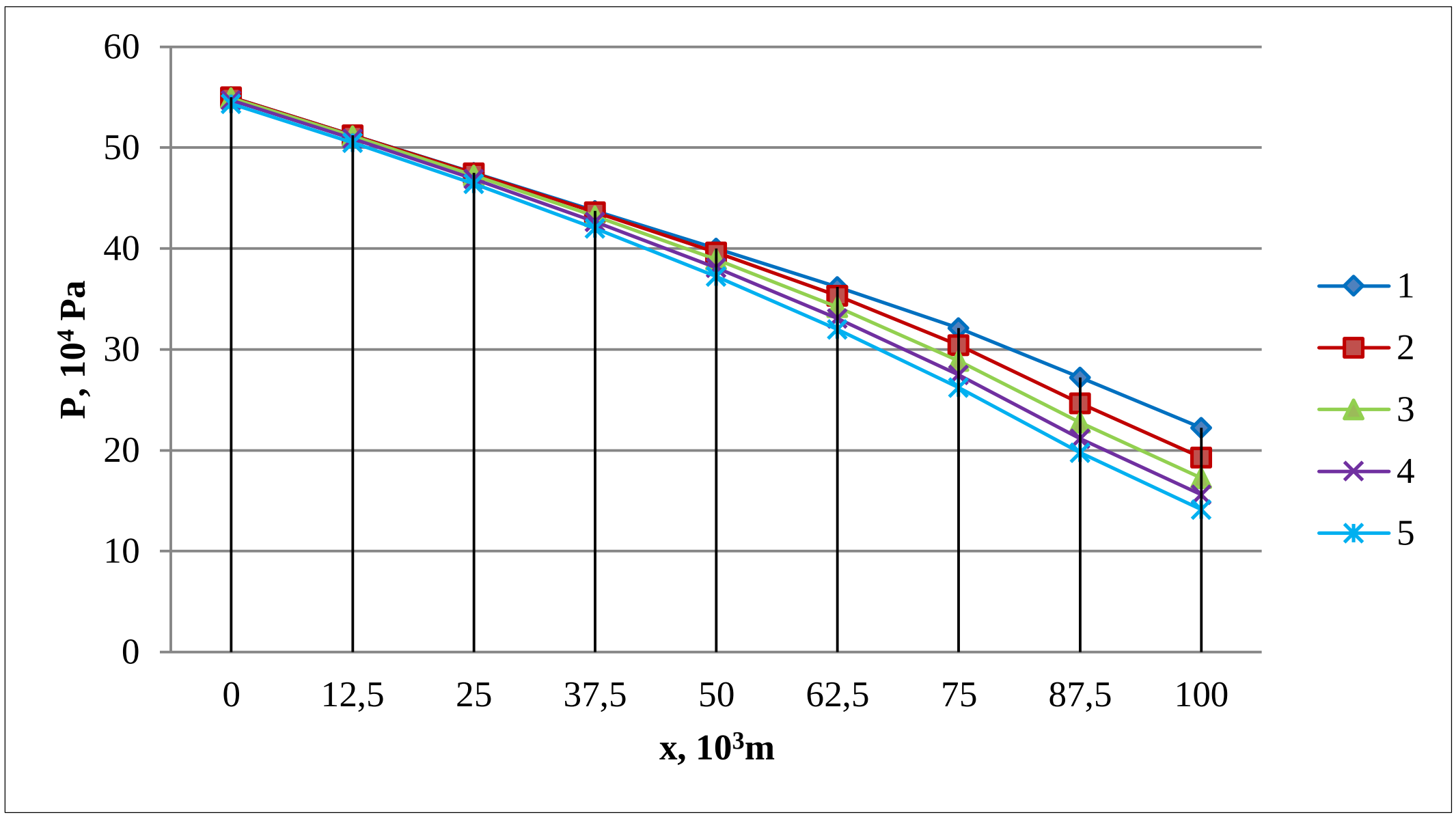

**Graph 4.2.3. The view of the distribution of pressure along the length of the damaged complex gas pipeline over time ( $\ell_2$ = 5 ×$10^4$ m, 1-t=100 sec, 2-t=300 sec, 3-t=500 sec, 4-t=700 sec, t=900 sec).**

Pressure is a key parameter in the operation of many engineering facilities, including gas supply systems. Control over this parameter is not only important during operation but also crucial during emergencies. The efficiency of gas pipelines' operation depends significantly on the control and regulation of this parameter. It should be noted that the more precise the control over this parameter, the higher the efficiency of gas pipelines' operation, which is essential for gas consumers as well. In other words, continuous control and regulation are necessary to ensure the proper operation of gas supply systems.

The analysis of the graphs and tables reveals that the quality of the longitudinal pressure variation over time varies significantly depending on the location of the leakage point. This variability is most pronounced in the pressures at the initial and final points. It appears that if an accident occurs at the beginning of the pipeline ($\ell_2 = 0{,}5 \times 10^4$ m), then at this time, the pressure drop at the beginning of the pipeline will be approximately $8.62 \times 10^4$ Pa after 600 seconds [($P_1$-$P_1$(0,t)] =(55-46,38) = $8{,}62\times10^4$ Pa. However, at the end of the pipeline, the pressure drop will be much smaller, approximately [($P_2$-$P_3$(L,t)]= (25-24,72) =$0{,}28\times10^4$ Pa, and it differs from the initial pressure drop by approximately 30 times (Graph 4.2.1 and Table 4.2.1).

If an accident (rupture) occurs at the end of the pipeline section ($\ell_2 = 9{,}5\times10^4$ m), then at the beginning, the pressure drop will be $0{,}16\times10^4$ Pa after 600 seconds (55-54,84)=$0{,}16\times10^4$ Pa However, at the end of the pipeline, the pressure drop will be much larger, approximately (25-16,38)=$8{,}62\times10^4$ Pa, which is about 53 times greater than the initial pressure drop (Graph 4.2.3 and Table 4.2.3).

Now let's consider the scenario where the accident occurs in the middle section of the pipeline ($\ell_2 = 5 \times10^4$ m). In this case, the pressure drop at the beginning will be $1{,}44\times10^4$ Pa after 600 seconds (55-53,56)=$1{,}44\times10^4$ Pa. However, at the end of the pipeline, the pressure drop will also be $1{,}44\times10^4$ Pa after 600 seconds (25-23,56) =$1{,}44\times10^4$ Pa, which is equal to the initial pressure drop (Graph 4.2.2 and Table 4.2.2).

The difference in pressure drops [($P_1$-$P_1$(0,t)]- [($P_2$-$P_3$(L,t)] varies over time depending on the constant value of gas leakage. For instance, at a distance of $\ell_2 = 0{,}5\times10^4$ m from the leakage point, it is 2.77 at t=100 seconds, 5.67 at t=300 seconds, and 7.59 at t=500 seconds (Graph 4.2.4). This means that near the beginning of the gas pipeline, if a gas leak occurs, the value of the difference in pres increases over time, while near the end, it decreases over time. However, in the middle section of the pipeline, if a gas leak occurs, the value of the difference in pressure drops remains constant over time (Graph 4.2.4).

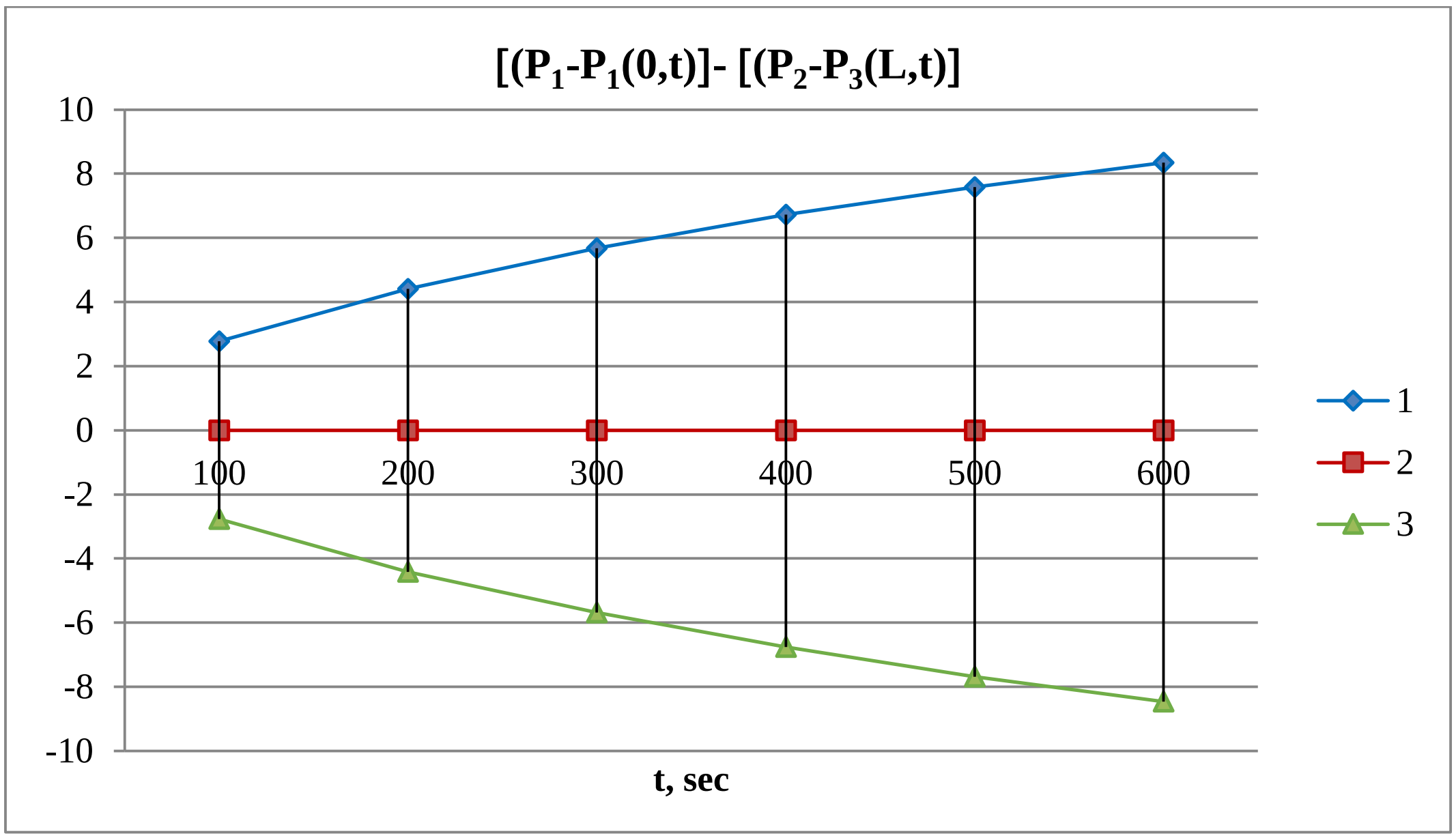

**Graph 4.2.4 The view of the distribution of the difference in pressure drops between the beginning and end points of the gas pipeline over time and distance from the location of the pipeline's rupture (1- $\ell_2$ = $0{,}5\times 10^4$ m, 2- $\ell_2 = 5 \times 10^4$ m, 3- $\ell_2 = 9{,}5\times10^4$ m).**

***Optimizing parameters of non-stationary flow.***

Similarly, it can be noted that in the vicinity of the beginning of the gas leakage point of the pipeline ($\ell_2$ = 0,5×10$^4$ m), after t=300 seconds, the value of the initial pressure drop is approximately $\left\{\frac{(55-49,30)}{(25-24,97)}\right\}$=190 times greater than the pressure drop generated near the end of the pipeline (Graph 4.2.1 and Table 4.2.1).

If a leakage occurs near the end of the pipeline ($\ell_2$ = 9,5×10$^4$ m), after t = 300 seconds, the pressure drop at the beginning of the pipeline is approximately $\left\{\frac{(55-54,99)}{(25-19,3)}\right\}$=0,002 times smaller than the pressure drop generated at the beginning of the pipeline (Graph 4.2.3 and Table 4.2.3). If the leakage occurs in the middle section of the pipeline, then after t = 300 seconds, the pressure drop at the beginning of the pipeline equals the pressure drop generated at the end of the pipeline, $\left\{\frac{(55-54,53)}{(25-24,37)}\right\}$=1 (Graph 4.2.4, Table 4.2.4, and Graph 4.2.5).

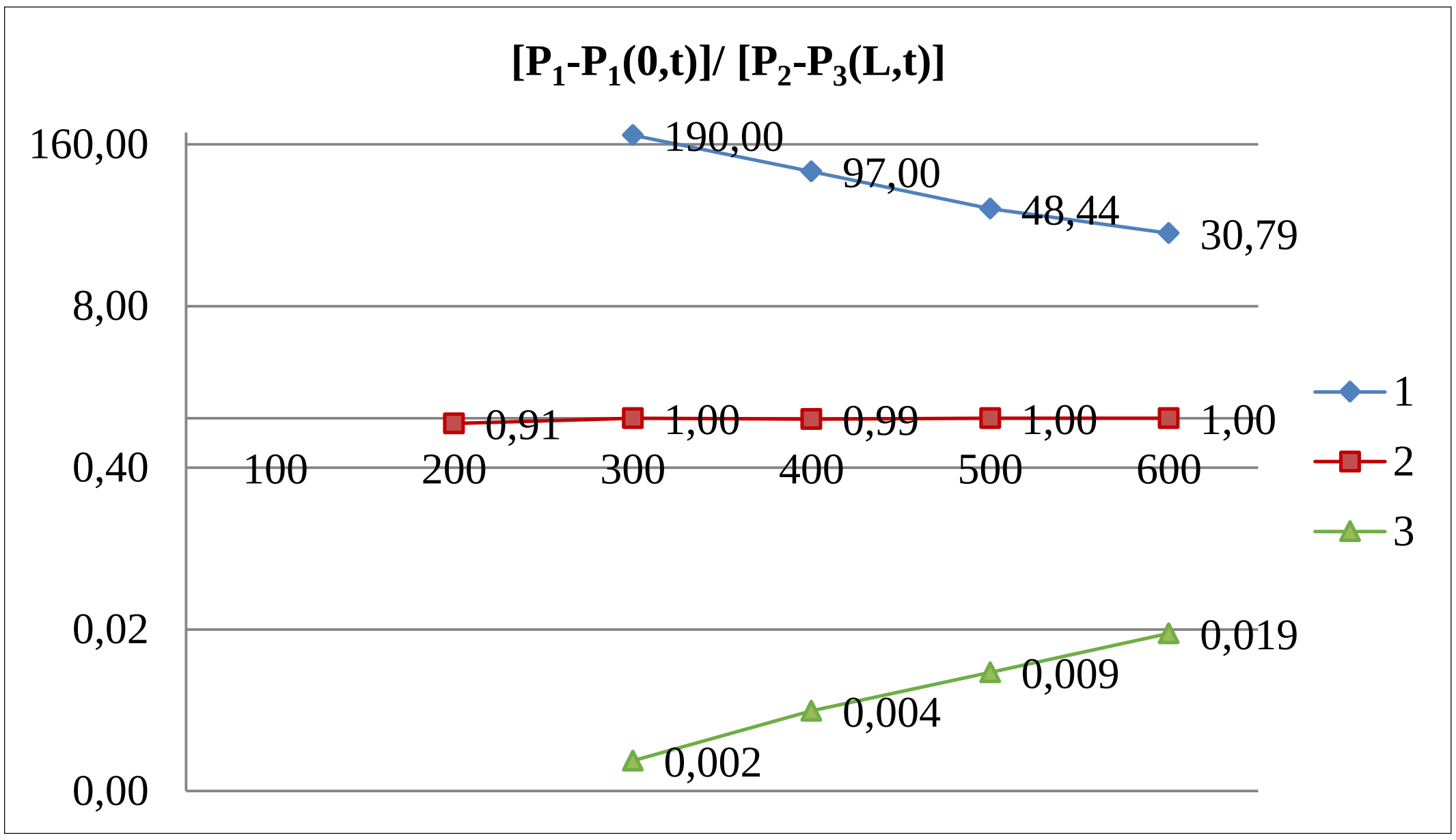

**Graph 4.2.5. The view of the distribution of the ratio of pressure drops at the beginning and end points of the gas pipeline over time, depending on the location and time of the pipeline's leakage. (1- $\ell_2$ = 0,5 × 10$^4$ m, 2- $\ell_2$ = 5× 10$^4$ m , 3- $\ell_2$ = 9,5 ×10$^4$ m).**

From the analysis of Graph 4.2.5, it can be concluded that the ratios of the pressure drops occurring at the beginning and end of the pipeline vary significantly depending on the location of the gas leakage point. However, as the duration of the incident increases, these ratios either decrease or increase. For instance, in the vicinity of the end of the pipeline ($\ell_2$ =9,5×10$^4$ m), after t = 300 seconds, the ratio of the pressure drop at the beginning of the pipeline to that at the end is 0.002. However, at t = 400 seconds, this ratio increases to 0.004, at t = 500 seconds to 0.009, and at t = 600 seconds, it further increases to 0.019.

The further analysis of Graph 4.2.5 reveals that the value of this ratio reaches its maximum or minimum at t=300 seconds for a constant value of gas leakage, and these values vary over time. Specifically, as the duration of the incident increases, the ratio of pressure drops due to the leakage near the beginning of the pipeline decreases. However, in the middle section of the pipeline, the ratio of pressure drops due to leakage stabilizes after t=300 seconds. In subsequent times, this value remains unchanged. Thus, by optimizing the parameters of complex gas pipeline systems during incident scenarios, we determine the fixed time $t=t_1$ for non-stationary flow parameters.

***The determination of the location of leakage in complex gas pipelines is carried out using analytical methods.***

Non-stationary gas flow models offer the advantage of being able to relate the analytical determination of the leakage location and the variation of pressure over time when valves on the pipeline are opened (Scheme 4.2.1). However, the difficulties arising in this situation make it impossible to accurately pinpoint the exact location where leaks or cracks occur. The superiority of stationary gas flow models lies in the simplicity of calculations and the precise determination of damaged locations. The disadvantage is that, after pressure drops at certain points in the gas pipeline, a certain period is required to identify the new stationary movement state of the gas flow dynamics. During this period, a large amount of gas is emitted into the surrounding environment.

The occurrence of a gas leak in a gas pipeline results in the formation of a non-stationary gas flow. The duration of the non-stationary regime is determined by the transition time from one stationary state to another, and it depends on the speed of propagation of pressure fluctuations along the length of the gas pipeline. Pressure fluctuation waves in the gas flow propagate at the speed of sound. Therefore, with a gas pipeline length of 100 km, the movement duration of the wave will be approximately 5 minutes. Thus, it is not difficult to locate the leak by observing the pressure drop in the network. In this case, there is no need to wait for the formation of a new stationary regime in the gas flow.

In the fixed time interval $t=t_1$, it is possible to observe a deviation in the pressure change from its value in the operational regime. The ability to detect small gas leaks will depend entirely on the accuracy of measurement instruments. In an automated dispatch control system, calculations are made in real-time for any changes in gas pressure at selected points. In this case, the actual and calculated pressure values for each section should be compared quickly. Therefore, determining the minimum time required to deliver reliable information to the dispatch center is a crucial issue.

The result of the analysis indicates that determining the location of leakage in pipelines analytically, based on the ratio of pressure drops occurring at the beginning and end of the pipeline, is more suitable for the intended purpose. Additionally, determining the location of the leakage at the time of fixation is crucial since the ratio of pressure drops occurring at the beginning and end of the pipeline varies significantly over time from the moment of the incident. Therefore, considering the factors mentioned, such as the simplicity and accuracy of the mathematical expressions obtained as a result of optimizing non-stationary flow parameters, will ensure reliability in determining the location of leakage in complex gas pipelines through analytical methods.

Based on Graph 4.2.2, it is evident that if a leakage (incident) occurs near the middle section of the pipeline (for the specific gas pipeline considered, where $\ell_2 = 5 \cdot 10^4$ m), then the pressure drop values at the end of the pipeline over time are approximately equal to the initial values. In this case, we can assume the equality $G_0(t) = G_s(t)$.

However, based on the analysis of the pressure variation along the pipeline over time for other instances of the location of the breach. If the location of the breach occurs near the beginning of the pipeline (considering a specific gas pipeline where $\ell_2$ =0,5·10$^4$ m), then it can be observed from Graph 4.2.5 that the value of the pressure drop at the end of the pipeline remains indeterminate until t= 300 seconds. This is because, as you mentioned, it takes about 5 minutes for the non-stationary gas flow wave to reach and be recorded at the final point at the moment of the accident.

On the other hand, if the location of the leak occurs near the end of the pipeline (considering a specific gas pipeline where $\ell_2$ = 9,5×10$^4$ m), then it also takes $t$=300 seconds for the gas flow wave to reach and be recorded at the starting point. It can be seen from Graph 4.2.5 that after $t$=300 seconds, the pressure drops at the beginning and end points of the pipeline take different values. It is understood that the change in gas consumption at the beginning and end points of the gas pipelines directly depends on the values of the pressure drops in those sections [26,50,78,92].

So, by using equations (18) and (19), we determine the following formula to identify the location of damage in non-stationary flow parameters in complex gas networks at the fixed time $t=t_1$

$$\ell_2 = \frac{1}{2}L + \varphi \cdot L \cdot \left[\frac{1-\phi(t)}{1+\phi(t)}\right] \qquad (20)$$

Here, $\phi(t) = \dfrac{P_1 - P_1(0,t_1)}{P_2 - P_2(L,t_1)}$; $\varphi = \dfrac{2}{3} + \dfrac{e^{-2\alpha t_1} - 4e^{-\alpha t_1}}{L}$ $\qquad \alpha = \dfrac{\pi^2 c^2}{2aL^2}$

As mentioned above, in accident regimes, the function $\phi(t)$ takes extreme values at $t=t_1$. In this case, equation (20) is used to determine the location of the gas leakage point and the location of the operating valves. To confirm our statements and for engineering calculations, we accept the known data and parameters obtained from calculations (Tables 4.2.1-4.2.3 and Graphs 4.2.1-4.2.5).

$P_b$ = 55× 10$^4$ Pa; $P_s$ = 25× 10$^4$ Pa; $G_0$= 30 Pa×sec/m ; $2a = 0.1\dfrac{1}{san}$ ; c=383.3 m/sec

L = 10$^5$ m: d=0.7m: $t_1$=300 sec: $\ell_2^f$ = 50000 m: $\ell_2^f$ = 5000 m: $\ell_2^f$ = 95000 m.

Using equation (20) and considering different values of the leakage point, we determine the leakage location and compare the calculated value with the actual value. We record the results of the calculation in the following table (Table 4.2.4).

| $P_1(0,t)$, $10^4$ Pa | $P_3(L,t)$, $10^4$ Pa | t, sec | $\varphi$ | $\phi(t)$ | $\ell_2^h$,m | $\dfrac{\ell_2^f - \ell_2^h}{\ell_2^f}$ |
|---|---|---|---|---|---|---|
| | $\ell_2^f$ = 50000 m | | | | | |
| 55,00 | 24,99 | 100,00 | 0,27 | 0,00 | 76915,02 | 0,54 |
| 54,90 | 24,89 | 200,00 | 0,36 | 0,91 | 51737,90 | 0,03 |
| 54,66 | 24,66 | 300,00 | 0,45 | 1,00 | 50000,00 | 0,00 |
| 54,34 | 24,33 | 400,00 | 0,53 | 0,99 | 50394,81 | 0,01 |
| 53,96 | 23,96 | 500,00 | 0,59 | 1,00 | 50000,00 | 0,00 |
| 53,56 | 23,56 | 600,00 | 0,66 | 1,00 | 50000,00 | 0,00 |
| | $\ell_2^f$ = 5000 m | | | | | |
| 52,23 | 25 | 100,00 | 0,27 | - | - | - |
| 50,58 | 25 | 200,00 | 0,36 | 442,00 | 13668,81 | 0,63 |
| 49,3 | 24,97 | 300,00 | 0,45 | 190,00 | 5512,41 | 0,09 |
| 48,21 | 24,93 | 400,00 | 0,53 | 97,00 | -1437,48 | 4,48 |
| 47,25 | 24,84 | 500,00 | 0,59 | 48,44 | -6913,52 | 1,72 |
| 46,38 | 24,72 | 600,00 | 0,66 | 30,79 | -11380,80 | 1,44 |
| | $\ell_2^f$ = 95000 m | | | | | |
| 55 | 22,23 | 100,00 | 0,27 | 0,000 | 76915,02 | -0,19 |
| 55 | 20,58 | 200,00 | 0,36 | 0,000 | 86495,95 | -0,09 |
| 54,99 | 19,3 | 300,00 | 0,45 | 0,002 | 94800,88 | 0,00 |
| 54,97 | 18,21 | 400,00 | 0,53 | 0,004 | 102047,13 | 0,07 |
| 54,93 | 17,25 | 500,00 | 0,59 | 0,009 | 108251,16 | 0,14 |
| 54,84 | 16,38 | 600,00 | 0,66 | 0,019 | 113114,96 | 0,19 |

**Table 4.2.4. Comparison of the calculated and actual values of the gas leakage point for different leakage point values.**

Table 4.2.4 shows that the calculated leakage point location is accurately determined only at $t=t1=300$ seconds, and its relative error is minimized. The reason for the relative error lies in the simplification of the non-linear system of equations and the simplification of mathematical expressions. However, it is accepted for engineering calculations and dispatcher stations. Indeed, the principles obtained through the calculation of mathematical expressions derived from the study of physical processes and the application of the new calculation scheme are suitable for practical use in experience.

**Result:**

1. The mathematical modeling based on the numerical solution describing the non-stationary gas flow in the gas pipeline was used as the initial basis. The research object included not only the network itself but also the connecting fittings and fixtures; thus, a complete gas-dynamic system consisting of mass conservation and momentum change equations was considered. During the analysis, the distribution of gas pressures along the length of the gas pipeline was investigated at three different locations of the leakage point, leading to the identification of new, previously unknown effects, and the existing calculation methods were tested and improved accordingly.
2. The investigation employed a calculation scheme to determine the minimum time fixed at $t=t_1$ for capturing non-stationary flow parameters as a result of optimizing the parameters in the accident regime of complex gas pipelines. The nature of the developed method lies in selecting optimal strategies for the reconstruction of complex gas pipelines.
3. The investigation employed a calculation scheme to determine the minimum time fixed at $t=t_1$ for capturing non-stationary flow parameters as a result of optimizing the parameters in the accident regime of complex gas pipelines. The nature of the developed method lies in selecting optimal strategies for the reconstruction of complex gas pipelines.

4. The investigation employed a calculation scheme to determine the minimum time fixed at $t=t_1$ for capturing non-stationary flow parameters as a result of optimizing the parameters in the accident regime of complex gas pipelines. The nature of the developed method lies in selecting optimal strategies for the reconstruction of complex gas pipelines.

### 4.3. Development of a Method for Calculating Dynamics Resulting from Gas Influx and Discharge Processes in Pipeline Operation

**Introduction:** The main gas pipelines are the backbone of the gas transportation system. The pipeline network accounts for the majority of energy resource expenditures in gas transportation. Safety issues related to the operation of main gas transportation facilities are currently relevant both during the design of new gas pipelines and during the reconstruction of existing gas transmission networks.

The characteristics of the pipeline network, along with the imposed technological requirements, primarily determine the operating regime of the remaining equipment in compressor stations. Therefore, ensuring the capability to accurately calculate the operating regime of the main gas pipeline is crucial. In this regard, all significant factors influencing the operation of the main gas pipeline should be taken into account. The calculation method should avoid complexities that cannot be resolved mathematically.

The most effective way to reliably predict the dynamics of transition processes and their outcomes is through the mathematical (numerical) modeling of the observed physical processes. In this study, the proposed method for analyzing gas pipelines relies on utilizing a collection of mathematical models that describe mass transfer processes occurring when non-uniform gas flow enters and exits different sections of the pipeline operating in normal (stationary) mode. These mathematical models are used to numerically model various accident scenarios. This process is being investigated by us for the first time. Of course, these processes have been studied [21, 64]. These studies have examined the resolution of problems by gradually altering the number of connections, disconnections, and related entries. Previous researchers have separately examined the processes of non-uniform gas flow entering the gas pipeline and the separation (discharge) of non-uniform gas flow from the pipeline. In our research, however, we analyze the differential equation system describing the non-stationary isothermal flow occurring in the gas pipeline as a result of both events happening simultaneously. We undertake the analytical solution of this system and analyze the dynamics of the process to ensure its compliance with the laws.

The physical consistency of transition processes involves evaluating the results of processes related to the modeling of the non-stationary operating regime of a high-pressure gas pipeline as a result of the simultaneous entry and exit (discharge) of non-uniform gas flow from different locations of the pipeline operating in normal mode. This entails modeling the distribution of the gas pipeline along its length, requiring information on the varying pressure and flow parameters of the gas slug over time.

Upon reviewing the existing research, it is evident that all universal models, despite their diversity, adhere to general formulation principles. For example, two perspectives can be clearly observed in the establishment of these models. In the first perspective (stage), models are developed for each type of elementary object within the gas transportation system. For all these models, the solution types that emerge should ensure their mutual compatibility and the possibility of subsequent joint solutions using existing mathematical methods when constructing the model of a complex pipeline system. In the second stage, a scheme for the complex pipeline system under consideration is developed using this or that mathematical apparatus. For this purpose, we adopt the scheme of the process under consideration as follows: (Scheme 4.3.1).

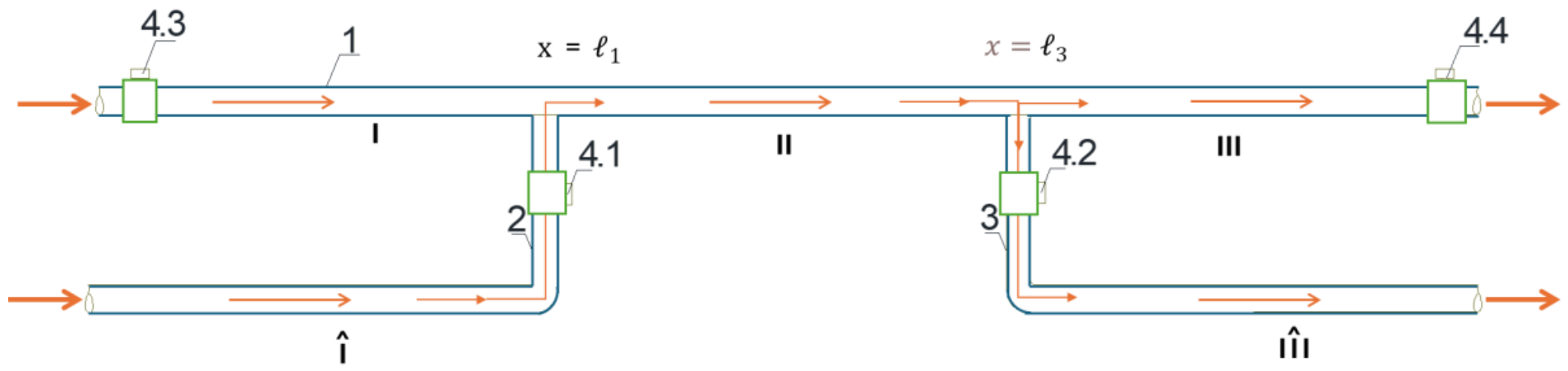

**Scheme 4.3.1: Technological scheme of the process of simultaneous entry and exit (discharge) of gas slugs from different locations of the gas pipeline at the same time.**

1. Investigated main gas pipeline; 2. Roadside pipeline entering the main gas pipeline at point x= $\ell_1$ 3. Roadside pipeline exiting the main gas pipeline at point x= $\ell_3$ , I, II, III - Separate technological processes of the main gas pipeline consisting of the 1st, 2nd, and 3rd sections. ; $\hat{I}$- Roadside pipeline entering section; $\hat{I}$-Roadside pipeline exiting section; 4.1 and 4.2 - Valves installed at the end of the roadside pipeline entering the main gas pipeline and at the beginning of the separated pipeline respectively. 4.3 and 4.4 - Valves installed at the beginning and end of the main gas pipeline respectively.

The accurate investigation of the analyzed scheme and the integrity of the obtained results require that the physical processes occurring in the roadside pipelines adhere to the requirements of the scheme's specifications.

Indeed, in the , $\hat{I}$ section, at point x= $\ell_1$ (in the entering section), the pressure of the gas should be higher than the pressure in the corresponding section of the main gas pipeline. In the . $\hat{III}$ section, at point x= $\ell_3$, the opposite should be true; Thus, the physical processes of these sections are adjusted to fit the transition process of the parallel gas pipelines analyzed in [3]. The conformity of the pressure change along the line and over time in the entering (connecting) pipeline of the main gas pipeline is assumed to be the same as the dynamics of the first section of the damaged pipeline in the parallel gas pipeline. In that section, a filling process occurs in the pipeline, causing the pressure along the line to increase over time. The conformity of the pressure change along the line and over time in the separated (discharged) pipeline is assumed to be the same as the dynamics of the third section of the damaged pipeline in the parallel gas pipeline. In this section, a depletion process occurs in the pipeline due to the consumption of gas by consumers, resulting in a decrease in the pressure along the line over time. Special models are developed for the calculation of non-stationary conditions in various sections. Analytical solutions that can practically be applied in real calculations for each considered scenario are presented, taking into account their characteristics, and recommendations are provided. Additionally, considering their length and time-dependent distribution, an isothermal gas flow with physical parameters is considered [3].

The non-stationary isothermal movement of gas in gas pipelines, which have steady-state sections, is described by the following system of non-linear differential equations [26,60]:

For the 1st section of the main gas pipeline.

$$\begin{cases} -\dfrac{\partial P^I}{\partial X} = \lambda \dfrac{\rho V^2}{2d} \\ -\dfrac{1}{c^2}\dfrac{\partial P^I}{\partial t} = \dfrac{\partial G^{II}}{\partial X} \end{cases} \qquad (1)$$

For the 2ndsection of the main gas pipeline.

$$\begin{cases} -\frac{\partial P^{II}}{\partial X} = \lambda \frac{\rho V^2}{2d} \\ -\frac{1}{c^2}\frac{\partial P^{II}}{\partial t} = \frac{\partial G^{II}}{\partial X} \end{cases} \quad (2)$$

For the 3rd section of the main gas pipeline.

$$\begin{cases} -\frac{\partial P^{III}}{\partial X} = \lambda \frac{\rho V^2}{2d} \\ -\frac{1}{c^2}\frac{\partial P^{III}}{\partial t} = \frac{\partial G^{III}}{\partial X} \end{cases} \quad (3)$$

For the main gas pipeline, for the section where the roadside pipeline enters at point x= $\ell_1$

$$\begin{cases} -\frac{\partial \hat{P}^{I}}{\partial X} = \lambda \frac{\rho V^2}{2d} \\ -\frac{1}{c^2}\frac{\partial \hat{P}^{I}}{\partial t} = \frac{\partial \hat{G}^{I}}{\partial X} \end{cases} \quad (4)$$

For the section of the main gas pipeline where the roadside pipeline enters at point x= $\ell_3$ .

$$\begin{cases} -\frac{\partial \hat{P}^{III}}{\partial X} = \lambda \frac{\rho V^2}{2d} \\ -\frac{1}{c^2}\frac{\partial \hat{P}^{III}}{\partial t} = \frac{\partial \hat{G}^{III}}{\partial X} \end{cases} \quad (5)$$

Here, G=ρv; P=ρ·$c^2$; c- speed of sound for the isothermal process in the gas flow, m/s., d - internal diameter of the pipe, m; x - coordinate aligned with the pipe axis and directed along the gas flow, m; P - absolute average gas pressure in the sections, Pa; v - average velocity of gas at the section, m/s; ρ - average gas density, kg/$m^3$; t - time, seconds; λ - dimensionless hydraulic resistance coefficient of the gas pipeline section; G - mass flow rate of the gas flow, (Pa·s)/m

The system of differential equations (1–5) is nonlinear and cannot be accurately integrated for any real boundary conditions. However, even with consideration of the peculiarities of gas pipeline operation, solving this simplified system of equations still presents certain difficulties.

It is noted that the accurate analytical solution of systems of equations describing the motion of real fluids in pipes [31] has not yet been obtained; therefore, it is necessary to use numerical methods, approximate analytics, or methods to transform non-linear equations into linear ones. In other words, various regularization methods should be employed. In solving the system of differential equations shown above with complex boundary conditions relevant to the problem under consideration, the regularization method proposed by I. A. Charnov was employed.

$$2a = \lambda \frac{V}{2d}$$

The main idea of the method is to transform the system of non-linear equations into a heat transfer equation system for the mathematical solution of the problem.

For the 1st section of the main gas pipeline.

$$\frac{\partial^2 P^I}{\partial X^2} = \frac{2a}{c^2}\frac{\partial P^I}{\partial t} \qquad 0 \le x \le \ell_1 \quad (6)$$

For the 2nd section of the main gas pipeline.

$$\frac{\partial^2 P^{II}}{\partial X^2} = \frac{2a}{c^2}\frac{\partial P^{II}}{\partial t} \qquad \ell_1 \le x \le \ell_3 \qquad (7)$$

For the 3rd section of the main gas pipeline.

$$\frac{\partial^2 P^{III}}{\partial X^2} = \frac{2a}{c^2}\frac{\partial P^{III}}{\partial t} \qquad \ell_1 \le x \le L \qquad (8)$$

For the roadside pipeline entering at point x= $\ell_1$ on the main gas pipeline.

$$\frac{\partial^2 \acute{P}^{I}}{\partial X^2} = \frac{2a}{c^2}\frac{\partial \acute{P}^{I}}{\partial t} \qquad 0 \le x \le \ell_1 \qquad (9)$$

For the roadside pipeline leaving at point x= $\ell_3$ on the main gas pipeline

$$\frac{\partial^2 P^{III}}{\partial X^2} = \frac{2a}{c^2}\frac{\partial P^{III}}{\partial t} \qquad \ell_1 \le x \le L \qquad (10)$$

The initial distribution of the function at the starting time t=0 for the physical process under consideration is given. We can assume this for the main gas pipeline.

$$P_i(x,0) = P_1 - 2aG_0x: \quad i= 1,2,3.$$

At the point x= $\ell_1$ , based on the technology of the transition process, the initial condition of the i-th section at the time of entry of the roadside pipeline into the main gas pipeline relies on the non-zero value of the initial condition being equal to the pressure of the roadside pipeline $\tilde{P}_1(x, t_2)$ According to Scheme 4.3.1, the pressure value $[\tilde{P}_1(x, t_2)]$ is the value at the moment of transition from the closed to open position of the valve at the end of the pipeline. At point x= $\ell_3$ , the transition process occurs in reverse, such that the gas flow from the 3rd section of the main gas pipeline flows into the roadside pipeline. In this case, the initial condition of the separated pipeline is assumed to be equal to the pressure distribution $\tilde{P}_3(x, t_2)$ from the starting point of the separated 3rd section to its end. Thus, the distribution of the function reflecting the initial conditions of the 1st and 3rd sections of the roadside gas pipelines accompanies different variations, meaning that during the connection period, a filling process occurs in the incoming pipeline, while an emptying process occurs in the separated pipeline.

Therefore, the non-zero initial conditions at x= $\ell_1$ correspond to the function $\tilde{P}_1(x, t_2)$ of the transition process in the first section at t=t2, and at x= $\ell_3$, the function $\tilde{P}_3(x, t_2)$ of the third section during the transition process is considered. As a result, for the initial conditions, we use the moment t = t2, meaning at that moment, the entry and exit of the overhead pipeline into the main gas pipeline occur simultaneously (at the same time). As shown in Scheme 4.3.1, the conformity reflecting the transition process occurs at the moments of activating the valves installed at points x= $\ell_1$ and x= $\ell_3$. The moment of their opening can be accepted as the initial condition. If these valves are opened at t=t2, then the distribution of the sought function (pressure) at t=t2 is given.

For t=t2

$$\acute{P}_i(x,0) = \tilde{P}_i(x, t_2): \quad i= 1,3.$$

The mathematical expression for solving equations (6-10) is determined based on the type of boundary conditions. Proper specification of the boundary conditions completes the mathematical model of the process and allows for comprehensive and accurate investigation of the physical events occurring. The law-abiding changes of known flow expenses at the beginning of the main gas pipeline and the pipeline, as well as at the end of the main pipe and the separated pipe, are provided. At points x= $\ell_1$ and x= $\ell_3$, the law of continuity is used to provide the law-abiding changes of known pressure and flow expenses.

At the point x= 0
$$\begin{cases} \frac{\partial P^I(x,t)}{\partial x} = -2aG_0(t) \\ \frac{\partial \hat{P}^I(x,t)}{\partial x} = -2a\hat{G}_0(t) \end{cases}$$

At the point x= $\ell_1$
$$\begin{cases} \frac{\partial P^{II}(x,t)}{\partial x} - \frac{\partial P^I(x,t)}{\partial x} = \frac{\partial \hat{P}^I(x,t)}{\partial x} \\ P^I(x,t) = P^{II}(x,t) \\ P^I(x,t) = \hat{P}^I(x,t) \end{cases}$$

At the point x= $\ell_3$
$$\begin{cases} \frac{\partial P^{II}(x,t)}{\partial x} - \frac{\partial P^{III}(x,t)}{\partial x} = \frac{\partial \hat{P}^{III}(x,t)}{\partial x} \\ P^{II}(x,t) = P^{II}(x,t) \\ P^{II}(x,t) = \hat{P}^{IIi}(x,t) \end{cases}$$

At the point x= $L$
$$\begin{cases} \frac{\partial P^{III}(x,t)}{\partial x} = -2aG_k(t) \\ \frac{\partial \hat{P}^{III}(x,t)}{\partial x} = -2a\hat{G}_K(t) \end{cases}$$

It is clear that the application of the Laplace transform converts equations (6-10) into second-order ordinary differential equations, and their general solutions will be as follows.

$$P^I(x,s) = \frac{P_i(x,0)}{S} + c_1 Sh\lambda x + c_2 Ch\lambda x \qquad 0 \le x \le \ell_1 \tag{11}$$

$$P^{II}(x,s) = \frac{P_i(x,0)}{S} + c_3 Sh\lambda x + c_4 Ch\lambda x \qquad \ell_1 \le x \le \ell_3 \tag{12}$$

$$P^{III}(x,s) = \frac{P_i(x,0)}{S} + c_5 Sh\lambda x + c_6 Ch\lambda x \qquad \ell_3 \le x \le L \tag{13}$$

$$\hat{P}^I(x,s) = \frac{\tilde{P}_1(x,0)}{S} + \hat{c}_1 Sh\lambda x + \hat{c}_2 Ch\lambda x \qquad 0 \le x \le \ell_1 \tag{14}$$

$$\hat{P}^{III}(x,s) = \frac{\tilde{P}_3(x,0)}{S} + \hat{c}_5 Sh\lambda x + \hat{c}_6 Ch\lambda x \qquad \ell_3 \le x \le L \tag{15}$$

Here, $\lambda = \sqrt{\frac{2as}{c^2}}$, $\beta = \sqrt{\frac{2ac^2}{s}}$

Applying the Laplace transform to the initial and boundary conditions, considering equations (11-15), we can find the constants $C_1$, $C_2$, $C_3$, $C_4$, $C_5$, $C_6$, $\hat{C}_1$, $\hat{C}_2$, $\hat{C}_5$ and $\hat{C}_6$.

$$c_1 = -B$$

$$\hat{c}_1 = -\hat{B}$$

$$c_2 = 4A_2\frac{ch\lambda(L-\ell_3)}{Z} - \frac{4(\Gamma+\hat{\Gamma})}{Z} + 4D_2\frac{sh\lambda(L-\ell_3)}{Z} - \frac{4g}{Zsh\lambda 2\ell_1}(B+\hat{B}) +$$

$$+2A_1\frac{g}{Zsh\lambda\ell_1} - 2D_1\frac{g}{Zch\lambda\ell_1} - \frac{A_1}{sh\lambda\ell_1} + \frac{2}{sh\lambda 2\ell_1}(B+\hat{B}) + \frac{D_1}{ch\lambda\ell_1} + B\frac{sh\lambda\ell_1}{ch\lambda\ell_1}$$

$$\hat{c}_2 = 4A_2\frac{ch\lambda(L-\ell_3)}{Z} - \frac{4(\Gamma+\hat{\Gamma})}{Z} + 4D_2\frac{sh\lambda(L-\ell_3)}{Z} - \frac{4g}{Zsh\lambda 2\ell_1}(B+\hat{B}) +$$

$$+2A_1\frac{g}{Zsh\lambda\ell_1} - 2D_1\frac{g}{Zch\lambda\ell_1} - \frac{A_1}{sh\lambda\ell_1} + \frac{2}{sh\lambda 2\ell_1}(B+\hat{B}) + \hat{B}\frac{sh\lambda\ell_1}{ch\lambda\ell_1}$$

$$c_3 = 2A_2\frac{sh\lambda 2\ell_1 ch\lambda(L-\ell_3)}{Z} - \frac{2sh\lambda 2\ell_1(\Gamma+\hat{\Gamma})}{Z} + 2D_2\frac{sh\lambda(L-\ell_3)sh\lambda 2\ell_1}{Z} +$$

$$+A_1\frac{g}{Z}\frac{sh\lambda 2\ell_1}{sh\lambda\ell_1} - \frac{2g}{Z}(B+\hat{B}) - D_1\frac{g}{Z}\frac{sh\lambda 2\ell_1}{ch\lambda\ell_1}$$

$$c_4 = 4A_2\frac{(1-sh^2\lambda\ell_1)ch\lambda(L-\ell_3)}{Z} - \frac{4(1-sh^2\lambda\ell_1)(\Gamma+\hat{\Gamma})}{Z} + 4D_2\frac{sh\lambda(L-\ell_3)(1-sh^2\lambda\ell_1)}{Z} +$$

$$+2A_1\frac{g}{Z}\frac{(1-sh^2\lambda\ell_1)}{sh\lambda\ell_1} - \frac{4g(1-sh^2\lambda\ell_1)}{Zsh\lambda 2\ell_1}(B+\hat{B}) - 2D_1\frac{g}{Z}\frac{(1-sh^2\lambda\ell_1)}{ch\lambda\ell_1} -$$

$$-\frac{A_1}{sh\lambda\ell_1} + \frac{2}{sh\lambda 2\ell_1}(B+\hat{B}) + \frac{D_1}{ch\lambda\ell_1}$$

$$c_5 = -4A_2\frac{sh\lambda L[sh\lambda\ell_1 sh\lambda(\ell_3-\ell_1) + ch\lambda\ell_3]}{Z} + \frac{4sh\lambda L[sh\lambda\ell_1 sh\lambda(\ell_3-\ell_1) + ch\lambda\ell_3](\Gamma+\hat{\Gamma})}{Zch\lambda(L-\ell_3)} +$$

$$+\frac{4gsh\lambda L[sh\lambda\ell_1 sh\lambda(\ell_3-\ell_1) + ch\lambda\ell_3]}{Zsh\lambda 2\ell_1 ch\lambda(L-\ell_3)}(B+\hat{B}) - 2A_1\frac{gsh\lambda L[sh\lambda\ell_1 sh\lambda(\ell_3-\ell_1) + ch\lambda\ell_3]}{Zsh\lambda\ell_1 ch\lambda(L-\ell_3)} -$$

$$-4D_2\frac{sh\lambda L sh\lambda(L-\ell_3)[sh\lambda\ell_1 sh\lambda(\ell_3-\ell_1) + ch\lambda\ell_3]}{Zch\lambda(L-\ell_3)} + 2D_1\frac{gsh\lambda L[sh\lambda\ell_1 sh\lambda(\ell_3-\ell_1) + ch\lambda\ell_3]}{Zch\lambda\ell_1 ch\lambda(L-\ell_3)} +$$

$$+A_1\frac{sh\lambda L ch\lambda\ell_3}{sh\lambda\ell_1 ch\lambda(L-\ell_3)} - \frac{2sh\lambda L ch\lambda\ell_3}{sh\lambda 2\ell_1 ch\lambda(L-\ell_3)}(B+\hat{B}) - D_1\frac{sh\lambda L ch\lambda\ell_3}{ch\lambda\ell_1 ch\lambda(L-\ell_3)} - \Gamma\frac{ch\lambda\ell_3}{ch\lambda(L-\ell_3)}$$

$$\hat{c}_5 = -4A_2\frac{sh\lambda L[sh\lambda\ell_1 sh\lambda(\ell_3-\ell_1) + ch\lambda\ell_3]}{Z} + \frac{4sh\lambda L[sh\lambda\ell_1 sh\lambda(\ell_3-\ell_1) + ch\lambda\ell_3](\Gamma+\hat{\Gamma})}{Zch\lambda(L-\ell_3)} +$$

$$+\frac{4gsh\lambda L[sh\lambda\ell_1 sh\lambda(\ell_3-\ell_1) + ch\lambda\ell_3]}{Zsh\lambda 2\ell_1 ch\lambda(L-\ell_3)}(B+\hat{B}) - 2A_1\frac{gsh\lambda L[sh\lambda\ell_1 sh\lambda(\ell_3-\ell_1) + ch\lambda\ell_3]}{Zsh\lambda\ell_1 ch\lambda(L-\ell_3)} -$$

$$-4D_2\frac{sh\lambda L sh\lambda(L-\ell_3)[sh\lambda\ell_1 sh\lambda(\ell_3-\ell_1) + ch\lambda\ell_3]}{Zch\lambda(L-\ell_3)} + 2D_1\frac{gsh\lambda L[sh\lambda\ell_1 sh\lambda(\ell_3-\ell_1) + ch\lambda\ell_3]}{Zch\lambda\ell_1 ch\lambda(L-\ell_3)} +$$

$$+A_1\frac{sh\lambda L ch\lambda\ell_3}{sh\lambda\ell_1 ch\lambda(L-\ell_3)} - \frac{2sh\lambda L ch\lambda\ell_3}{sh\lambda 2\ell_1 ch\lambda(L-\ell_3)}(B+\hat{B}) - D_1\frac{sh\lambda L ch\lambda\ell_3}{ch\lambda\ell_1 ch\lambda(L-\ell_3)} - \hat{\Gamma}\frac{ch\lambda\ell_3}{ch\lambda(L-\ell_3)} +$$

$$+D_2\frac{sh\lambda L}{ch\lambda(L-\ell_3)}$$

$$c_6 = 4A_2 \frac{ch\lambda L[sh\lambda\ell_1 sh\lambda(\ell_3-\ell_1)+ch\lambda\ell_3]}{Z} - \frac{4ch\lambda L[sh\lambda\ell_1 sh\lambda(\ell_3-\ell_1)+ch\lambda\ell_3]\left(\Gamma+\hat{\Gamma}\right)}{Zch\lambda(L-\ell_3)} -$$

$$-\frac{4gch\lambda L[sh\lambda\ell_1 sh\lambda(\ell_3-\ell_1)+ch\lambda\ell_3]}{Zsh\lambda 2\ell_1 ch\lambda(L-\ell_3)}\left(B+\hat{B}\right) + 2A_1 \frac{gch\lambda L[sh\lambda\ell_1 sh\lambda(\ell_3-\ell_1)+ch\lambda\ell_3]}{Zsh\lambda\ell_1 ch\lambda(L-\ell_3)} +$$

$$+4D_2 \frac{ch\lambda L sh\lambda(L-\ell_3)[sh\lambda\ell_1 sh\lambda(\ell_3-\ell_1)+ch\lambda\ell_3]}{Zch\lambda(L-\ell_3)} - 2D_1 \frac{gch\lambda L[sh\lambda\ell_1 sh\lambda(\ell_3-\ell_1)+ch\lambda\ell_3]}{Zch\lambda\ell_1 ch\lambda(L-\ell_3)} -$$

$$-A_1 \frac{ch\lambda L ch\lambda\ell_3}{sh\lambda\ell_1 ch\lambda(L-\ell_3)} + \frac{2ch\lambda L ch\lambda\ell_3}{sh\lambda 2\ell_1 ch\lambda(L-\ell_3)}\left(B+\hat{B}\right) + D_1 \frac{ch\lambda L ch\lambda\ell_3}{ch\lambda\ell_1 ch\lambda(L-\ell_3)} + \Gamma \frac{sh\lambda\ell_3}{ch\lambda(L-\ell_3)}$$

$$\hat{c}_6 = 4A_2 \frac{ch\lambda L[sh\lambda\ell_1 sh\lambda(\ell_3-\ell_1)+ch\lambda\ell_3]}{Z} - \frac{4ch\lambda L[sh\lambda\ell_1 sh\lambda(\ell_3-\ell_1)+ch\lambda\ell_3]\left(\Gamma+\hat{\Gamma}\right)}{Zch\lambda(L-\ell_3)} -$$

$$-\frac{4gch\lambda L[sh\lambda\ell_1 sh\lambda(\ell_3-\ell_1)+ch\lambda\ell_3]}{Zsh\lambda 2\ell_1 ch\lambda(L-\ell_3)}\left(B+\hat{B}\right) + 2A_1 \frac{gch\lambda L[sh\lambda\ell_1 sh\lambda(\ell_3-\ell_1)+ch\lambda\ell_3]}{Zsh\lambda\ell_1 ch\lambda(L-\ell_3)} +$$

$$+4D_2 \frac{ch\lambda L sh\lambda(L-\ell_3)[sh\lambda\ell_1 sh\lambda(\ell_3-\ell_1)+ch\lambda\ell_3]}{Zch\lambda(L-\ell_3)} - 2D_1 \frac{gch\lambda L[sh\lambda\ell_1 sh\lambda(\ell_3-\ell_1)+ch\lambda\ell_3]}{Zch\lambda\ell_1 ch\lambda(L-\ell_3)} -$$

$$-A_1 \frac{ch\lambda L ch\lambda\ell_3}{sh\lambda\ell_1 ch\lambda(L-\ell_3)} + \frac{2ch\lambda L ch\lambda\ell_3}{sh\lambda 2\ell_1 ch\lambda(L-\ell_3)}\left(B+\hat{B}\right) + D_1 \frac{ch\lambda L ch\lambda\ell_3}{ch\lambda\ell_1 ch\lambda(L-\ell_3)} - D_2 \frac{ch\lambda L}{ch\lambda(L-\ell_3)} +$$

$$+\Gamma \frac{sh\lambda\ell_3}{ch\lambda(L-\ell_3)}$$

We obtain the values of these constants by considering equations (11-15) and deriving the transformed form of the equation, representing the dynamic state of the gas pipeline in question. These equations will express the distribution of pressure along the pipeline in transformed form for the considered sections.

$$P^{I}(x,s) = \frac{P_i - 2aG_0x}{s} + 4A_2 \frac{ch\lambda(L-\ell_3)ch\lambda x}{Z} - \frac{4ch\lambda x(\Gamma+\hat{\Gamma})}{Z} + 4D_2 \frac{sh\lambda(L-\ell_3)ch\lambda x}{Z} -$$
$$-\frac{4\left(gch\lambda x - \frac{Z}{2}ch\lambda x\right)}{Zsh\lambda 2\ell_1}\left(B+\hat{B}\right) - 4B \frac{\left(\frac{Z}{2}sh\lambda(x-\ell_1)sh\lambda\ell_1\right)}{Zsh\lambda 2\ell_1} + 2A_1 \frac{gch\lambda x - \frac{Z}{2}ch\lambda x}{Zsh\lambda\ell_1} - \quad (16)$$
$$-2D_1 \frac{gch\lambda x - \frac{Z}{2}ch\lambda x}{Zch\lambda\ell_1}$$

$$P^{II}(x,s) = \frac{P_i - 2aG_0x}{s} + 4A_2 \frac{ch\lambda(L-\ell_3)[sh\lambda\ell_1 sh\lambda(x-\ell_1)+ch\lambda x]}{Z} - \frac{4[sh\lambda\ell_1 sh\lambda(x-\ell_1)+ch\lambda x](\Gamma+\hat{\Gamma})}{Z} +$$
$$+4D_2 \frac{sh\lambda(L-\ell_3)[sh\lambda\ell_1 sh\lambda(x-\ell_1)+ch\lambda x]}{Z} + 2A_1 \frac{g[sh\lambda\ell_1 sh\lambda(x-\ell_1)+ch\lambda x] - \frac{Z}{2}ch\lambda x}{Zsh\lambda\ell_1} - \quad (17)$$
$$-\frac{4\left(g[sh\lambda\ell_1 sh\lambda(x-\ell_1)+ch\lambda x] - \frac{Z}{2}ch\lambda x\right)}{Zsh\lambda 2\ell_1}\left(B+\hat{B}\right) - 2D_1 \frac{g[sh\lambda\ell_1 sh\lambda(x-\ell_1)+ch\lambda x] - \frac{Z}{2}ch\lambda x}{Zch\lambda\ell_1}$$

$$P^{III}(x,s) = \frac{P_i - 2aG_0x}{s} + 4A_2 \frac{ch\lambda(L-x)[sh\lambda\ell_1 sh\lambda(\ell_3-\ell_1)+ch\lambda x]}{Z} -$$
$$-\frac{4ch\lambda(L-x)[sh\lambda\ell_1 sh\lambda(\ell_3-\ell_1)+ch\lambda x](\Gamma+\hat{\Gamma})}{Zch\lambda(L-\ell_3)} - -\frac{\Gamma sh\lambda(x-\ell_3)}{ch\lambda(L-\ell_3)}$$
$$+4D_2 \frac{sh\lambda(L-\ell_3)ch\lambda(L-x)[sh\lambda\ell_1 sh\lambda(\ell_3-\ell_1)+ch\lambda x]}{Zch\lambda(L-\ell_3)} + 2A_1 \frac{ch\lambda(L-x)\left[g[sh\lambda\ell_1 sh\lambda(\ell_3-\ell_1)+ch\lambda x] - \frac{Z}{2}ch\lambda\ell_3\right]}{Zch\lambda(L-\ell_3)sh\lambda\ell_1} -$$
$$-\frac{4ch\lambda(L-x)\left(g[sh\lambda\ell_1 sh\lambda(\ell_3-\ell_1)+ch\lambda x] - \frac{Z}{2}ch\lambda\ell_3\right)}{Zsh\lambda 2\ell_1 ch\lambda(L-\ell_3)}\left(B+\hat{B}\right) - 2D_1 \frac{ch\lambda(L-x)\left\{g[sh\lambda\ell_1 sh\lambda(x-\ell_1)+ch\lambda x] - \frac{Z}{2}ch\lambda\ell_3\right\}}{Zch\lambda\ell_1 ch\lambda(L-\ell_3)}$$

(18)

$$\hat{P}^{I}(x,s)=\frac{\tilde{P}^{I}(x,s)}{s}+4A_2\frac{ch\lambda(L-\ell_3)ch\lambda x}{Z}-\frac{4ch\lambda x(\Gamma+\hat{\Gamma})}{Z}+4D_2\frac{sh\lambda(L-\ell_3)ch\lambda x}{Z}-$$
$$-\frac{4\left(gch\lambda x-\frac{Z}{2}ch\lambda x\right)}{Zsh\lambda 2\ell_1}\left(B+\hat{B}\right)-4\hat{B}\frac{\left(\frac{Z}{2}sh\lambda(x-\ell_1)sh\lambda\ell_1\right)}{Zsh\lambda 2\ell_1}+2A_1\frac{gch\lambda x-\frac{Z}{2}ch\lambda x}{Zsh\lambda\ell_1}-$$
$$-2D_1\frac{gch\lambda x}{Zch\lambda\ell_1}$$

(19)

$$\hat{P}^{III}(x,s)=\frac{\tilde{P}^{III}(x,t_2)}{s}+4A_2\frac{ch\lambda(L-x)[sh\lambda\ell_1 sh\lambda(\ell_3-\ell_1)+ch\lambda x]}{Z}-\frac{4ch\lambda(L-x)[sh\lambda\ell_1 sh\lambda(\ell_3-\ell_1)+ch\lambda x](\Gamma+\hat{\Gamma})}{Zch\lambda(L-\ell_3)}-$$
$$-\frac{\hat{\Gamma}sh\lambda(x-\ell_3)}{ch\lambda(L-\ell_3)}+4D_2\frac{ch\lambda(L-x)\left\{[sh\lambda\ell_1 sh\lambda(\ell_3-\ell_1)+ch\lambda x]sh\lambda(L-\ell_3)-\frac{Z}{4}\right\}}{Zch\lambda(L-\ell_3)}+$$

4

$$+2A_1\frac{ch\lambda(L-x)\left[g[sh\lambda\ell_1 sh\lambda(\ell_3-\ell_1)+ch\lambda x]-\frac{Z}{2}ch\lambda\ell_3\right]}{Zch\lambda(L-\ell_3)sh\lambda\ell_1}-\frac{ch\lambda(L-x)\left(g[sh\lambda\ell_1 sh\lambda(\ell_3-\ell_1)+ch\lambda x]-\frac{Z}{2}ch\lambda\ell_3\right)}{Zsh\lambda 2\ell_1 ch\lambda(L-\ell_3)}\text{x}$$
$$\text{x}\left(B+\hat{B}\right)-2D_1\frac{ch\lambda(L-x)\left\{g[sh\lambda\ell_1 sh\lambda(x-\ell_1)+ch\lambda x]-\frac{Z}{2}ch\lambda\ell_3\right\}}{Zch\lambda\ell_1 ch\lambda(L-\ell_3)}$$

.

(20)

Here,

$$A_1=\frac{d\tilde{P}^{I}(\ell_1,t_2)}{s\lambda dx};\ A_2=\frac{d\tilde{P}^{III}(\ell_3,t_2)}{s\lambda dx};\ B=\left[\frac{2aG_0(t)}{\lambda}-\frac{2aG_0}{s\lambda}\right];\ \hat{B}=\left[\frac{2a\hat{G}_0(t)}{\lambda}+\frac{d\tilde{P}(0,t_2)}{s\lambda dx}\right]$$

$$D_1=\frac{\tilde{P}^{I}(\ell_1,t_2)-P^{I}(\ell_1,0)}{s};\ \Gamma=\left[\frac{2aG_k(t)}{\lambda}-\frac{2aG_0}{s\lambda}\right];\ \ \hat{\Gamma}=\left[\frac{2a\hat{G}_k(t)}{\lambda}+\frac{d\tilde{P}^{III}(L,t_2)}{s\lambda dx}\right]$$

$$D_2=\frac{\tilde{P}^{III}(\ell_3,t_2)-P^{III}(\ell_3,0)}{s};\ g=[3sh\lambda L-sh\lambda(L-2\ell_3)];$$

$$\mathrm{Z}=\mathrm{sh}\lambda 2\ell_1[3\mathrm{ch}\lambda\mathrm{L}-\mathrm{ch}\lambda(\mathrm{L}-2\ell_3)]+2\mathrm{g}(1-\mathrm{sh}^2\lambda\ell_1).$$

Since the sums of equations (16-20) do not have an analytical solution due to the transcendental nature of the current source, it is not possible to directly revert to the original form here. To do this, numerous mathematical operations are required to adapt the summation to the Laplace transformation table. On the other hand, the original mathematical expressions obtained as a result of these operations will consist of numerous series and complex integrals, which are not considered efficient for engineering calculations. To mitigate this difficulty, it is more appropriate to consider the formation of a new steady-state regime in the gas pipeline to sufficiently accurately calculate the operating regime of the pipeline. In addition, the analysis of the new steady-state regime associated with the completion of the transition process is simpler and clearer. As noted in [27], the transition period in the gas pipeline is approximately 33.3 minutes. Thus, reaching the new steady-state regime will take approximately half an hour. On the other hand, it is known that accidents and planned repair works are the main reasons for the occurrence of transition processes. In both cases, the duration of repair and restoration works is envisaged to be within 6 to 24 hours according to regulations. So, within the analyzed half-hour period, complex processes resulting in dangerous increases or decreases in pressure in the pipelines, disruptions in gas supply to consumers, and similar situations will not cause an excessive load on the dynamic limit of the pipeline section.

In stationary regimes, the pressure remains constant at the beginning and end sections of the gas pipeline and is determined as follows.

$$\hat{G}_0(0)=G_o(0)=\hat{G}_k(0)=G_o(0)=\alpha G_0$$

In this case, the distribution of pressure along the pipeline for each section of the newly established stationary regime will be as follows.

$$P^{I}(x) = P_i - 2aG_0x - 2aG_0(\alpha-1)(x-\ell_1) + 4aG_0(\alpha-1)\frac{(L-\ell_3)(\ell_3-\ell_1)}{(2L+\ell_1-\ell_3)} + \\ +\left[\tilde{P}(\ell_3,t_2) - P(\ell_3,0)\right]\frac{L-\ell_3}{(2L+\ell_1-\ell_3)} + \frac{\ell_1[\tilde{P}(\ell_1,t_2)-P(\ell_1,0)]}{(2L+\ell_1-\ell_3)} + 2aG_0(\alpha-1)\frac{(L-\ell_1)^2-\ell^2{}_1}{(2L+\ell_1-\ell_3)} \quad (21)$$

$$\acute{P}^{I}(x) = P^{I}(x) - P_1 + 2aG_0x + \tilde{P}^{I}(x,t_2) - \acute{P}^{I}(\ell_1,t_2) + P^{I}(\ell_3,0) \quad (22)$$

$$P^{II}(x) = P_i - 2aG_0x - \frac{4aG_0(\alpha-1)}{(2L-\ell_1-\ell_3)}\left[\ell_1(x-\ell_1) - (L-\ell_3)(\ell_3-x) + \frac{x^2}{2} - \frac{(L-x)^2}{2}\right] + \\ + \frac{\ell_1[\tilde{P}(\ell_1,t_2)-P(\ell_1,0)]}{(2L+\ell_1-\ell_3)} + \frac{(L-\ell_3)[\tilde{P}(\ell_3,t_2)-P(\ell_3,0)]}{(2L+\ell_1-\ell_3)} \quad (23)$$

$$P^{III}(x) = P_i - 2aG_0x - 2aG_0(\alpha-1)(x-\ell_3) - \frac{4aG_0(\alpha-1)}{(2L+\ell_1-\ell_3)}\left[\frac{\ell_1(\ell_3-\ell_1)}{2} + \frac{\ell_3{}^2}{2} - \frac{(L-\ell_3)^2}{2}\right] + \\ + \frac{\ell_1[\tilde{P}(\ell_1,t_2)-P(\ell_1,0)]}{(2L+\ell_1-\ell_3)} + \left[\tilde{P}(\ell_3,t_2) - P(\ell_3,0)\right]\frac{L-\ell_3}{(2L+\ell_1-\ell_3)} \quad (24)$$

$$P^{III}(x) = P^{III}(x) - P_1 + 2aG_0x + \tilde{P}^{III}(x,t_2) - P^{III}(\ell_3,t_2) + P^{III}(\ell_3,0) \quad (25)$$

Here, $\alpha = \sqrt{\frac{L}{L+3\ell}}$

L- is the length of the main gas pipeline, in meters.

$\ell$ is the distance from the point $x = \ell_3$ to the end of the main gas pipeline, in meters.

In equations (21-25), $\tilde{P}(\ell_1,t_2)$ and $\tilde{P}(\ell_3,t_2)$ represent the functions characterizing the variation of gas pressure along the main gas pipeline and the feeder pipelines entering and exiting the main gas pipeline, respectively, in terms of both spatial and temporal changes (under non-stationary conditions) at the specified coordinates ($\ell_1$, $\ell_3$) and at $t= t_2$. As seen in Scheme 1, these feeder pipelines can perform the function of the connectors of the parallel-planed gas pipelines. As noted in [3], when the integrity of the gas pipelines is compromised, automatic valves installed on the damaged pipeline separate the damaged section from the main section of the gas pipeline. As a result, a new transition process occurs in the gas pipeline. During this process, a filling process occurs in the section of the pipe from the beginning of the damaged section to the valve at the beginning of the damaged area, and as a result, the pressure of the gas in that section increases over time. If we consider this process to be the same as the process of the roadside pipe entering the main gas pipeline at the point $x= \ell_1$, and if we assume that the valves installed on the damaged pipeline are closed at $t=t_1$, then the mathematical expression characterizing the process can be as follows [8].

$$\tilde{P}_1(x,t) = P_1(x,t_1) + \frac{c^2}{\ell_1}\int_{t_1}^{t}\tilde{G}_0(\tau)d\tau + \frac{2c^2}{\ell_1}\sum_{n=1}^{\infty}\text{Cos}\frac{\pi n x}{\ell_1}\int_{t_1}^{t}\tilde{G}_0(\tau)\cdot e^{-\alpha_3(t-\tau)}d\tau + \\ + \frac{c^2}{2a\ell_1}\frac{dP_1(0,t_1)}{dx}\int_{t_1}^{t}\left[1+2\sum_{n=1}^{\infty}\text{Cos}\frac{\pi n x}{\ell_1}e^{-\alpha_3\tau}\right]d\tau \qquad 0 \le x \le \ell_1 \quad (26)$$

On the other hand, as a result of the mentioned process, after the point where the damaged section is disconnected and the subsequent valve is closed, there is an emptying process in the section from the beginning of the pipeline to the end, and as a result, the pressure of gas along that section decreases over time. If we consider this process to be the same as the dynamics of the roadside pipe separated from the main gas pipeline at the point x= $\ell_3$, then the mathematical expression characterizing the process would be as follows [9].

$$\tilde{P}_3(x,t) = P_3(x,t_1) - \frac{c^2}{(L-\ell_3)}\int_{t_1}^{t}\tilde{G}_s(\tau)\left[1+2\sum_{n=1}^{\infty}(-1)^n \operatorname{Cos}\frac{\pi n(x-\ell_3)}{L-\ell_3}\cdot e^{-\alpha_5(t-\tau)}\right]d\tau - $$
$$- \frac{c^2}{2a(L-\ell_3)}\frac{dP_3(L,t_1)}{dx}\int_{t_1}^{t}\left[1+2\sum_{n=1}^{\infty}(-1)^n \operatorname{Cos}\frac{\pi n(x-\ell_3)}{L-\ell_3}\cdot e^{-\alpha_5\tau}\right]d\tau \qquad \ell_3 \le x \le L \quad (27)$$

If we assume t=t₂ duration to have $\tilde{G}_0(\tau)=\tilde{G}_0(\tau)$=constant, and accept the given data for the calculation of equations (26) and (27) [9].

$P_b = 14\times10^{-2}$ MPa : $P_s = 11\times10^{-2}$ MPa : $G_0 = 10\,\frac{Pa\times sec}{m}$: $2a = 0{,}1\,\frac{1}{sec}$ ;

$L = 3\times10^4$ m: $\ell_1 = 1\times 10^4$ m ; $\ell_3 = 2\times 10^4$ m : $c = 383{,}3\,\frac{m}{sec}$

If we assume $t=t_2$=600 seconds and determine the values of expressions $\tilde{P}(\ell_1,t_2)$ and $\tilde{P}(\ell_3,t_2)$ from equations (26) and (27) using calculation software (for example, Mathcad).

$\tilde{P}(\ell_1,t_2)$ = 14,24×10⁻² MPa: $\tilde{P}(\ell_3,t_2)$ = 9,1× 10⁻² Mpa

We calculate the values of boyu $P^{I}(x)$, $P^{II}(x)$ and $P^{III}(x)$ for every 2500 m, taking into account the determined values and the values provided above in equations (21), (23), and (25) for the newly formed steady-state regime of the process. We record the results in the table below (Table 4.3.1)

.

| x, km | $P^{I}(x)$, $10^4$ Pa | $P^{II}(x)$ , $10^4$ Pa | $P^{III}(x)$ , $10^4$ Pa | P(x,0) , $10^4$ Pa |
|---|---|---|---|---|
| 0 | 13,09 | - | - | 14 |
| 2,5 | 12,91 | - | - | 13,75 |
| 5 | 12,73 | - | - | 13,5 |
| 7,5 | 12,55 | - | - | 13,25 |
| 10 | 12,38 | 12,38 | - | 13 |
| 12,5 | - | 12,27 | - | 12,75 |
| 15 | - | 12,17 | - | 12,5 |
| 17,5 | - | 12,06 | - | 12,25 |
| 20 | - | 11,96 | 11,96 | 12 |
| 22,5 | - | - | 11,74 | 11,75 |
| 25 | - | - | 11,53 | 11,5 |
| 27,5 | - | - | 11,31 | 11,25 |
| 30 | - | - | 11,1 | 11 |

Here, P(x,0) represents the distribution of pressure depending on the coordinate in the initial steady-state regime of the main gas pipeline before the transition process, in Pascal (Pa). Using Table 4.3.1, we construct the graph of the distribution of pressure depending on the coordinate for the newly formed steady-state regime in all three sections of the main gas pipeline during the complex transition process (Graph 4.3.1).

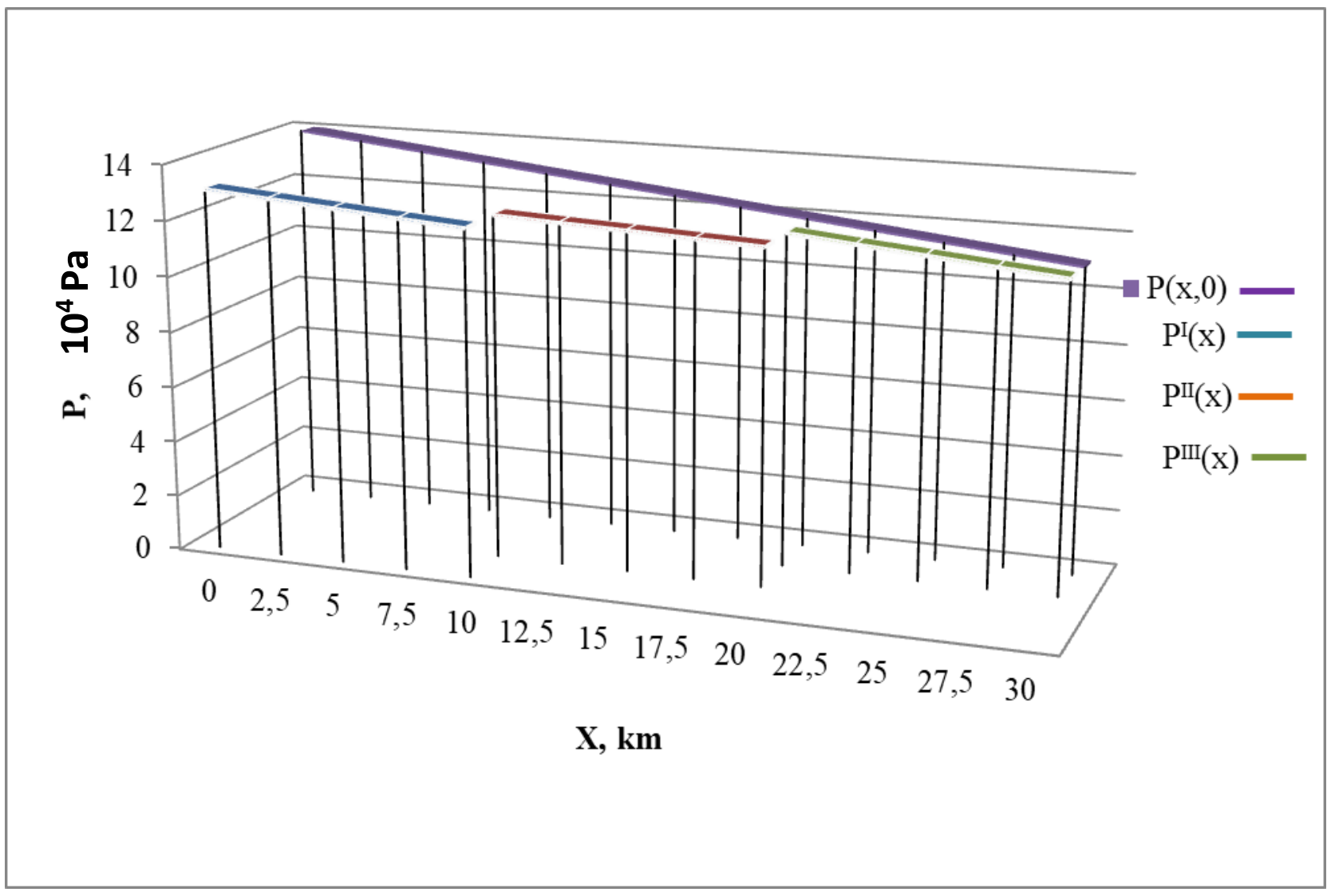

**Graph 4.3.1. Graph of pressure distribution depending on the coordinate in the newly formed steady-state regime.**

The analysis of Graph 4.3.1 reveals that the qualitative view of the newly formed steady-state regime resulting from the processes of gas ingress and egress (discharge) from different locations of the gas pipeline is consistent with the view of the previous steady-state regime. In other words, in both cases, the functions characterizing the distribution of pressure along the length of the pipeline are similar. However, as seen from Table 4.3.1, the values of pressure in the first section of the pipeline in the newly formed steady-state regime decrease by 6.5% to 4.8% compared to the values in the previous steady-state regime along the coordinate. This implies that at the point x= $\ell_1$ , where the gas enters the main pipeline through the lateral pipeline, the increase in hydraulic resistance in that section results in an increase in pressure losses along the length of the pipeline. A similar process occurs in the second section of the pipeline, but in this section, the decrease in pressure along the coordinate is reduced from 4.8% to 0.3%. Conversely, at the end of the third section of the pipeline, the opposite process occurs, such that the pressure at the end of that section in the newly formed steady-state regime exceeds the value in the previous regime by 1%. This implies that at the point x= $\ell_3$ , where the lateral pipeline connects to the main gas pipeline, the lateral pipeline plays the role of a booster. As we know, a booster is a parallel pipeline installed at certain sections of trunk gas pipelines to increase their productivity (pressure). Therefore, installing the lateral pipeline after the main pipeline serves the purpose of increasing the efficiency of gas pipeline exploitation.

## Conclusion

1. Gas pipeline movement along different sections, connected by an intermediary pipeline, and the process of non-uniform gas flow entering and exiting, based on the classical Laplace method, for solving the differential equations system describing the non-uniform isothermal gas flow in a pipeline with a constant cross-section, the solutions were given for the pressure and flow rate variations at the ends of the gas pipeline as well as at the junction points.
2. In non-stationary regimes, based on the transformed equations of the dynamic state of gas transportation processes in different sections of the main gas pipelines, the determination and analysis of pressure distribution in the new steady-state regime based on coordinates have been developed and successfully tested in theoretical calculations.
3. New mathematical expressions capable of analyzing the changes in pressure along the pipeline in the complex transition process of gas pipelines have been developed. Through specific calculations, these expressions have demonstrated sufficient adequacy and can be utilized to conduct more comprehensive analyses of the processes and provide substantial findings and recommendations.

**4.** To depict the non-steady-state flow of the studied processes, an analysis and qualitative overview of the obtained equations, considering the locations and sequence of separation and merging, were initially conducted. As a result, it is proposed to install the dedicated pipeline after the main pipeline inlet to enhance the reliability of gas pipeline operation.

### 4.4. Creation of an Information Base for Efficient Management of Reconstructed Complex Gas Pipeline Systems

**Introduction:** The reconstruction, management, and optimization of gas pipelines is of significant importance for solving modern engineering problems. This paper presents innovative methodologies aimed at the effective reconstruction of gas pipelines under unstable conditions. The research encompasses the application of machine learning and optimization algorithms, targeting the enhancement of system reliability and the optimization of interventions during emergencies. The findings of the study present engineering solutions aimed at addressing the challenges in real-world applications by comparing the performance of various algorithms. Consequently, this work contributes to the advancement of cutting-edge approaches in the field of engineering and opens new perspectives for future research.

The increasing consumption volumes and resource usage to meet the gas demand in the country can be ensured not only through the design and construction of new networks and structures but also through the reconstruction of gas pipelines and efficient management by utilizing potential reserves. Gas pipeline systems are complex systems characterized by numerous elements, and managing them requires a large amount of information corresponding to the number of possible elements. Managing such complex gas pipelines is impossible without a systematic approach based on the combined consideration of complex concepts of control theory, such as system, information, goal-oriented, and feedback. New approaches to reconstruction methods aimed at ensuring a higher level of quality are required in managing these systems.

The theory of selecting justified parameters for the operation and reconstruction of main gas pipelines requires the calculation of non-stationary processes resulting from any emergency [89]. This task should be considered one of the most pressing issues. The operating conditions of a reconstructed gas pipeline must ensure its ability to perform specified functions reliably for a fixed period.

he main purpose of the report scheme is to prepare a database in accordance with the requirements for the reconstruction and optimization technology of gas pipelines. This database can be applied to modern smart systems for managing gas supply using various methods and processes. It is known that recently, machine learning and artificial intelligence systems have been present in various fields, including the management of gas supply systems. A comprehensive analysis of the current state of machine learning and artificial intelligence in solving gas supply issues has been conducted by several researchers [92, 95, 124, 127]. These studies also discuss various types of machine learning and artificial intelligence methods that can be used to process and interpret data in different areas of the oil and gas industry. The use of these modern technologies can facilitate the decision-making process in the management of gas supply system reconstructions.

The novelty of the work lies in the development of reporting schemes for solving several important new problems related to the management and regulation of technological processes in gas transportation based on optimization methods and mathematical modeling. For this purpose, an algorithm for the operational management of gas transportation systems has been developed. This algorithm is based on the information obtained about the emergency state of pipelines by monitoring changes in the technological parameters of gas flow at the peripheral parts of the pipelines during emergency modes. This significantly enhances the reliability of gas supply management and the efficiency of transportation processes during the reconstruction phase of complex gas pipelines.

At the current stage, the issues we face differ fundamentally from those already addressed [101] in that they involve imposing specific requirements and time constraints during the development phase of the reconstruction concept for complex gas pipeline systems. This ensures the acquisition of structural, technical, and technological bases within the framework of future operation of complex gas pipeline systems, the clarification and elimination of technical and technological contradictions in the process of gas transportation system development and reconstruction, and the efficient management of the gas transportation system under modern conditions.

***Justification of an effective gas transportation scheme for the reconstruction of the existing parallel gas pipeline.***

The essence of the studied method is to justify a targeted technological variant for reconstruction by increasing the efficiency of gas transportation under non-stationary conditions. The reconstruction of gas transportation systems is characterized by the reliability of multi-level information bases transmitted to the dispatcher control center for complexes of local management systems of individual technological processes and objects.

The diagram showing the methodological process is in the figure below.

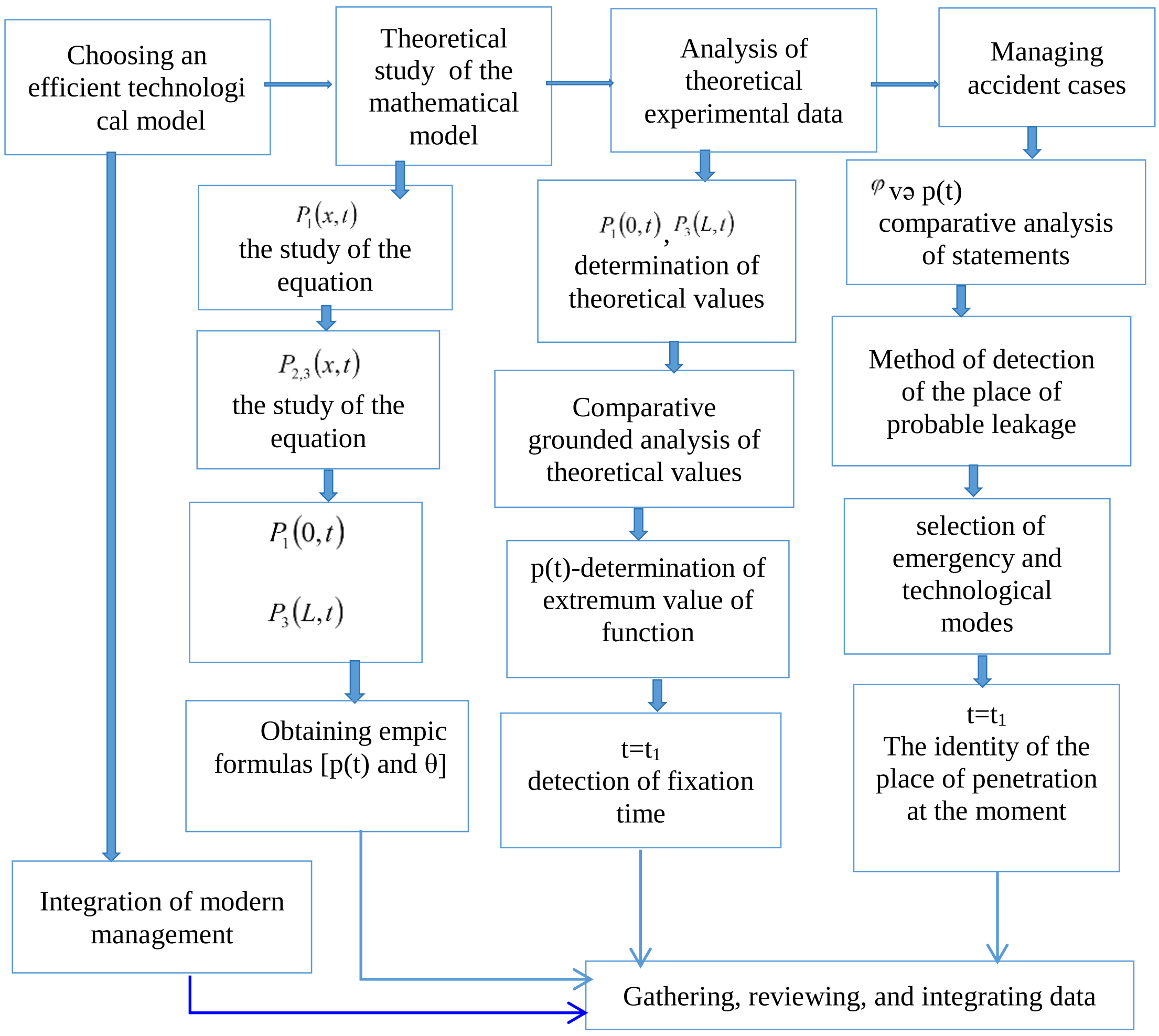

The reconstruction phase is based on the following key principles:

The advantages of an efficient technological scheme for complex gas pipeline systems;

The uninterrupted transportation of planned gas flows to consumers;

The level of reliability indicators for complex gas pipeline systems;

The safety of the operation and environmental aspects of complex gas pipeline systems;

Maximum conservation of energy and resources;

The enhancement of technological and economic efficiency levels of complex gas pipeline systems.

The "Main Pipeline Design" standards state that parallel gas pipelines are complex systems consisting of two or more pipeline lines. These lines have the same pressure at the beginning and end sections and are equipped with at least two connectors and connecting fittings to ensure the transfer of gas flow between the pipelines during technical operations and emergencies. These systems, having the same initial and final pressure and operating in a unified hydraulic regime, create challenges for the application of modern control systems in dispatcher centers during emergency regimes. In other words, the information obtained from changes in gas flow parameters at the end sections of the system due to the failure of one of the pipeline lines in parallel gas pipelines is not reliable. According to the research conducted in reference [90], it was confirmed that depending on the operating conditions of gas pipelines, the pressure values at the beginning and end sections of damaged (ruptured) and undamaged pipeline lines were equal. Moreover, it was even observed that, at any given time during operation, the change in pressure at the beginning and end of the gas pipeline was greater for the undamaged pipeline line than for the damaged one. In such conditions, identifying the damaged section in the dispatcher center becomes impossible. Indeed, the information obtained from the non-stationary nature of gas dynamic processes makes it unreliable for the effective management of gas pipeline systems.

Therefore, the development of effective gas transportation schemes for the reconstruction of existing parallel gas pipelines has become a pressing issue. This is crucial for ensuring the long-term utilization and environmental safety of complex gas pipeline transportation, as well as for achieving reliable and trouble-free operation and efficient management of emergency situations.

To prepare the technological-based reporting scheme for this purpose, the following issues need to be addressed:

1) Optimization of the operating regime of complex gas pipeline systems;

2) Implementation of effective strategies for managing parallel laid distributor pipeline systems during normal technological operating regimes and in case of emergencies;

3)Utilization of new methods with high quality for analyzing and optimizing the management of emergency and technological conditions and the selection of rational operating regimes;

4) Prevention and minimization of the consequences of emergency situations;

5) Ensuring the accuracy of damage location determination, automatic detection in real-time mode, and minimizing the impact of interventions in detection modes.

It is known that one of the main requirements for the reconstruction of gas pipelines is the application of modern technologies to the existing system. The China National Petroleum Corporation (CNPC) has independently developed 27 unique technologies across five series that have been successfully applied in the construction and management of various gas pipelines. These technologies have been implemented in the Lanzhou-Chengdu Chongqing Gas Pipeline, as well as pipeline projects in Sudan, Libya, India, Russia, and Central Asia. Some of these technologies include:Oil and Gas Pipeline Flow Assurance Technology;Oil and Gas Pipeline Simulation and Optimization Technology; Pipeline Integrity Management System;Early-warning and leakage detection technology for pipeline safety; Pipeline Emergency Repair Technology, and so on.

However, leaks in pipeline networks remain one of the primary causes of countless losses to pipeline stakeholders and the surrounding environment. Pipeline failures can lead to severe environmental disasters, human casualties, and financial losses. To prevent such hazards and maintain a safe and reliable pipeline infrastructure, it is essential to implement early warning and leakage detection technology for pipeline safety. The pipeline leakage detection technologies are discussed in [101], summarizing the latest advancements. Various leak detection and localization methods in pipeline systems are reviewed, with their strengths and weaknesses highlighted. Comparative performance analysis is carried out to provide guidance on which leak detection method is most suitable for specific operational parameters. Furthermore, research gaps and open issues for the development of reliable pipeline leakage detection systems are discussed.

It is known that one of the main requirements for the reconstruction of gas pipelines is the application of modern technologies to the existing system. The China National Petroleum Corporation (CNPC) has independently developed 27 unique technologies across five series that have been successfully applied in the construction and management of various gas pipelines. These technologies have been implemented in the Lanzhou-Chengdu Chongqing Gas Pipeline, as well as pipeline projects in Sudan, Libya, India, Russia, and Central Asia. Some of these technologies include:

Oil and Gas Pipeline Flow Assurance Technology; Oil and Gas Pipeline Simulation and Optimization Technology; Pipeline Integrity Management System; Early-warning and leakage detection technology for pipeline safety; Pipeline Emergency Repair Technology, and so on.

As an example, in [104, 106], we can mention a method based on a device that uses ultrasonic waves for radiation to manage pipelines. In the study being examined, acoustic sensors installed outside the pipeline detect the internal noise level and establish a baseline with specific characteristics. The selfsimilarity of this signal is continuously analyzed by the acoustic sensors. When a leak occurs, the resulting low-frequency acoustic signal is detected and investigated. If the properties of this signal differ from the original, an emergency signal is triggered. The received signal is stronger near the leak location, allowing for the localization of the leak. When a leak is detected, waves are reflected, and the location of the leak is determined by the emission speed of the wave and the return time of the scattered signals. Ripple conversion is used to process these signals and locate the leak. In general, this method is based on the detection of noise generated in the pipeline when a leak occurs [113].

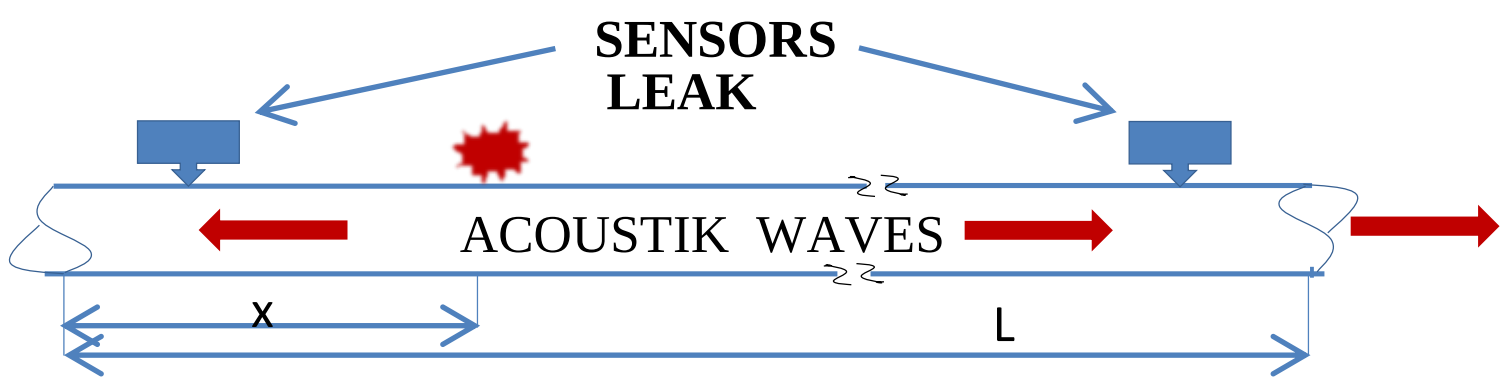

**Figure 4.4.1. Technological Figure of acoustic leak detection.**

Sensors are activated on the right and left of the point where the leak occurs, at any point x on a pipeline with length L, as shown in Figure1.

In this study [112], technologies for extracting the characteristics of acoustic signals have been summarized. Experiments highlighted the characteristics in the time, frequency, and frequency-time domains by measuring acoustic leaks and interference signals. Results indicate that external interference can be effectively eliminated through time domain, frequency domain, and frequency-time domain characteristics. It can be concluded that the acoustic leak detection method can be applied to natural gas pipelines and that its features can help reduce false alarms and missed alarms. Therefore, it can be concluded that the acoustic leak detection method can be applied to natural gas pipelines, helping to reduce false and missed alarm signals. Recently, many authors have extensively studied a system using an Internet of Things (IoT) analytical platform service to simulate real-time monitoring of pipelines and locate leaks within the pipeline [95]. In this study, pressure pulses based on the vibration principle of pipes are used for pipeline monitoring. The principle of time delay between pulse receptions at sensor positions is also applied in this research. The primary principle of the study was to create a wireless communication device that combines, programs, and interacts with the ThingSpeak IoT analytics platform through an Arduino and a Wi-Fi module. A total of five channels were created to collect data from five sensors used in experimental testing on the platform (Figure 4.4.2). These signal data are collected every 15 seconds,and all channels are updated every 2 minutes. ThingSpeak provides immediate visualization of data transmitted by the wireless communication device. Online analysis and processing of data were performed as it was received [82].

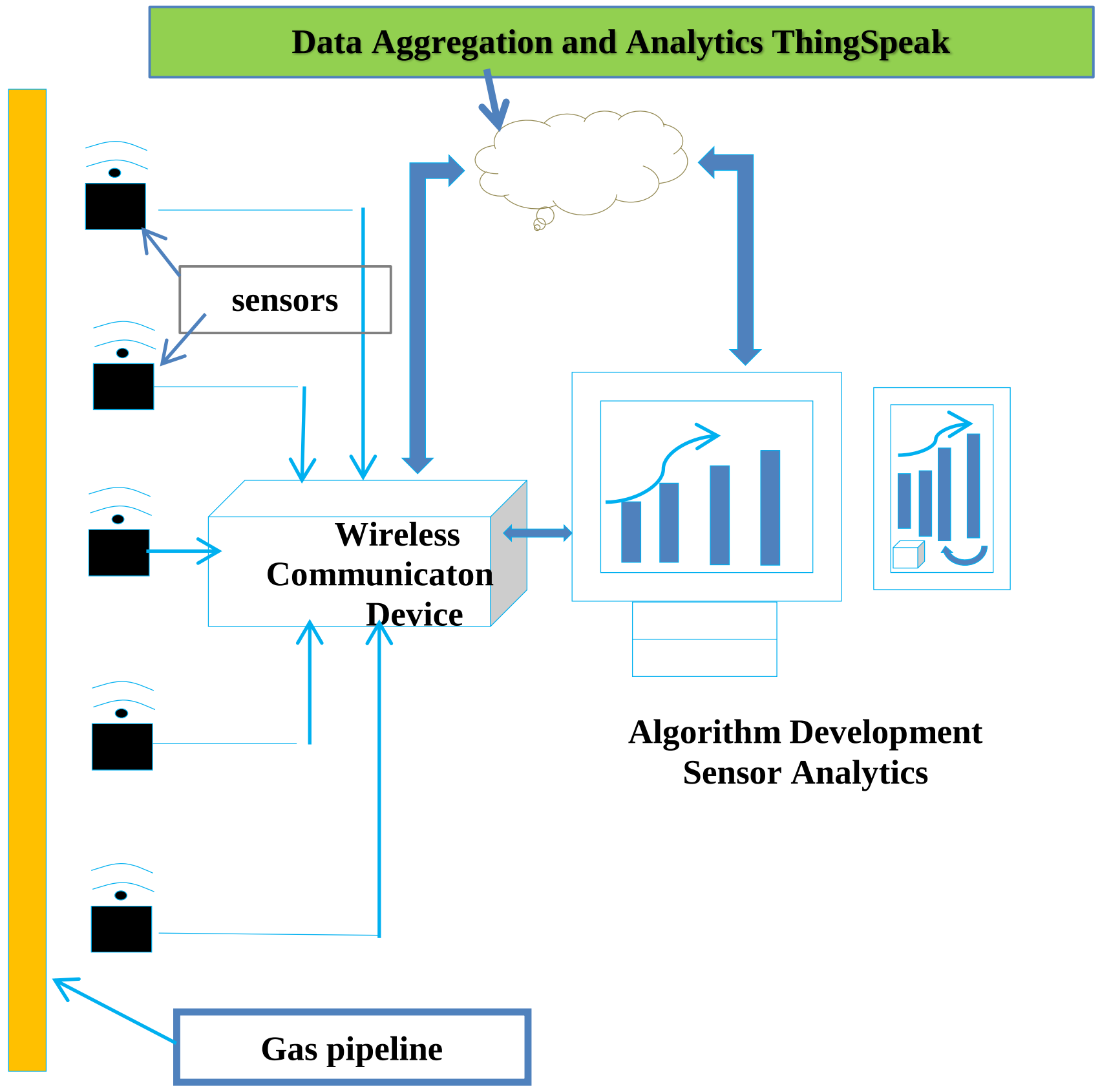

**Figure 4.4.2. IoT system used for pipeline monitoring.**

Using a data logging device and air as the transport fluid, the measured propagation speed of vibration was 355 m/s, ranging from 4.243 m to 4.246 m. This yielded only a 3 mm deviation when comparing the estimated leak location with the actual measured location of 4.23 m. The experimental results demonstrated that the performance of the wireless communication device compares favorably with that of the data logger and can determine the location of damage in actual pipelines when used for real-time monitoring. With this communication device and analytical platform, it is possible to monitor pipelines in real time from anywhere in the world on any device with internet access.

İintegrity management programs are also available today that can locate the damage. A thorough, organised, and unified integrity management programme provides the means to improve the safety of pipeline systems. Such a programme provides the information for a pipeline system operator to allocate resources for prevention, detection, and mitigation actions that will result in improved safety and a reduction in the number of incidents.See Figure 4.4.3 [22].

It covers performance and prescriptive based integrity management The performance-based integrity management programme employs more data and broad risk analysis, which permit pipeline system operators to comply with the requirements of this programme in the areas of inspection intervals, tools used, and mitigation techniques used.Inspection, prevention, detection, and mitigation are all part of the prescriptive process necessary to produce an integrity management programme.

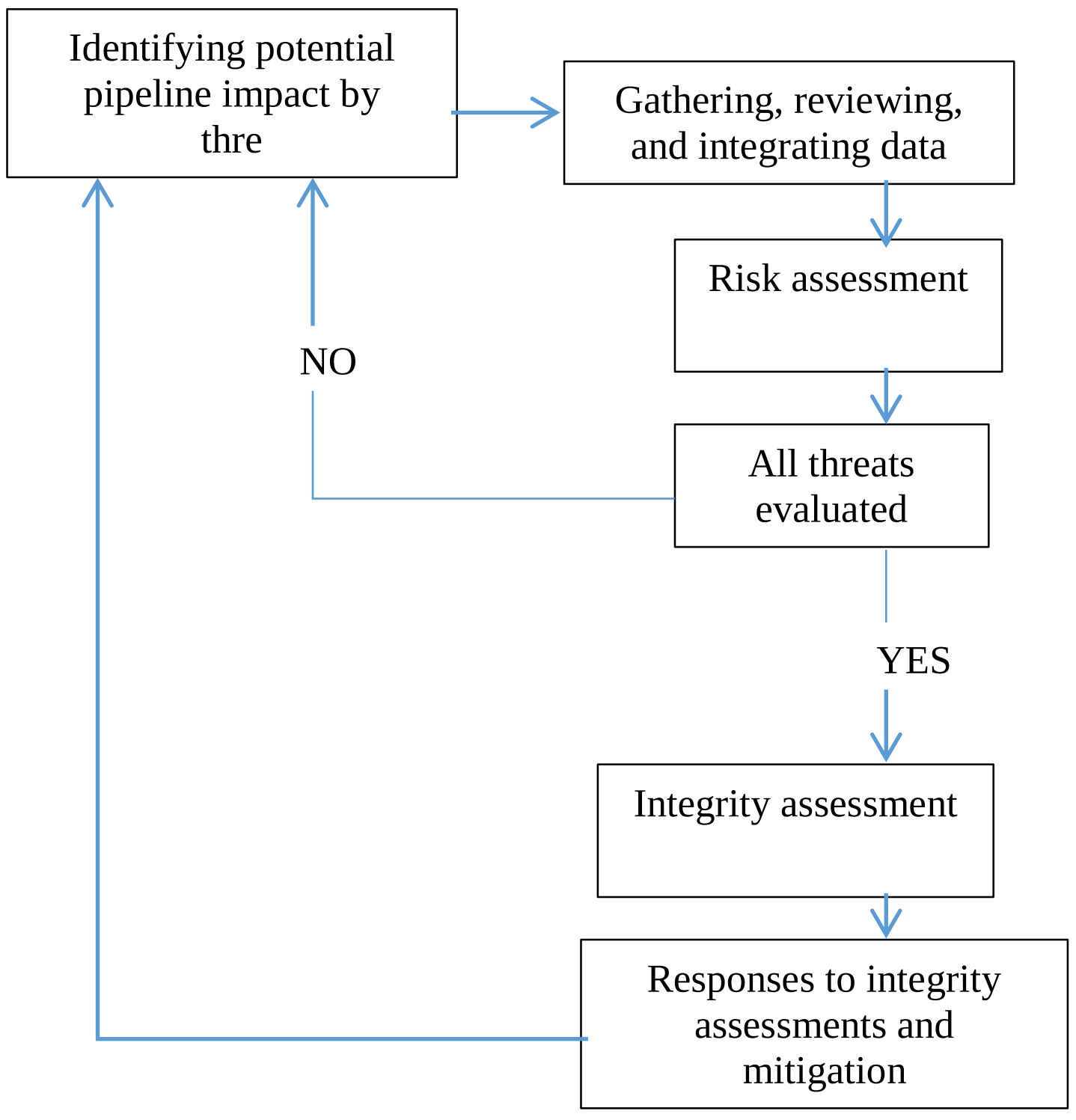

**Figure 4.4.3. Integrity management plan process flow diagram.**

Prior to proceeding with the performance-based integrity programme a pipeline system operator shall first ensure sufficient inspections are performed that provide the information on the pipeline condition required by the prescriptive based programme. The level of assurance of a performance-based programme shall meet or exceed that of a prescriptive programme [91].

It is known that there are three methods for locating leaks in gas pipelines: qualitative, quantitative, and analytical. The practical examples mentioned mainly pertain to the quantitative method. In our research, however, the focus was on the analytical method. Specifically, it involved developing a reporting Figure to establish a reliable information base for dispatch centers to effectively manage the gas dynamic processes in non-stationary modes. As mentioned above, accidents occurring in parallel gas pipelines operating in a unified hydraulic regime cause a decrease in pressure throughout all sections of both pipeline networks over time. The decrease in pressure at the end section of the pipeline deteriorates the process of gas supply to consumers.

If we assume that the pressure difference is fully utilized by all consumers, then as the consumption of the initial gas mass in the network increases, the pressure at its end point will decrease again, leading to various devices and appliances experiencing different levels of disruption. Devices and appliances operating at maximum pressure will immediately sense the drop in gas pressure, reducing their numbers in operation simultaneously. In such a regime, only devices and appliances close to the nominal pressure will be able to operate. Therefore, these issues should be addressed in the reconstruction of existing gas pipelines. Based on the principles of the reconstruction phase, selecting an efficient technological Figure that perfectly manages the technological and emergency processes of gas transportation is essential. As we know, the reliability indicator of main gas pipelines directly depends on providing consumers with the required uninterrupted gas supply [125]. Among the complex gas pipeline systems that excel at fulfilling this function is the parallel gas pipeline system.

As a result of accidents, we adopt the following technological scheme for the efficient management of gas pipelines, the detection of leakage points, and the provision of uninterrupted gas supply to consumers in complex gas pipeline systems [90].

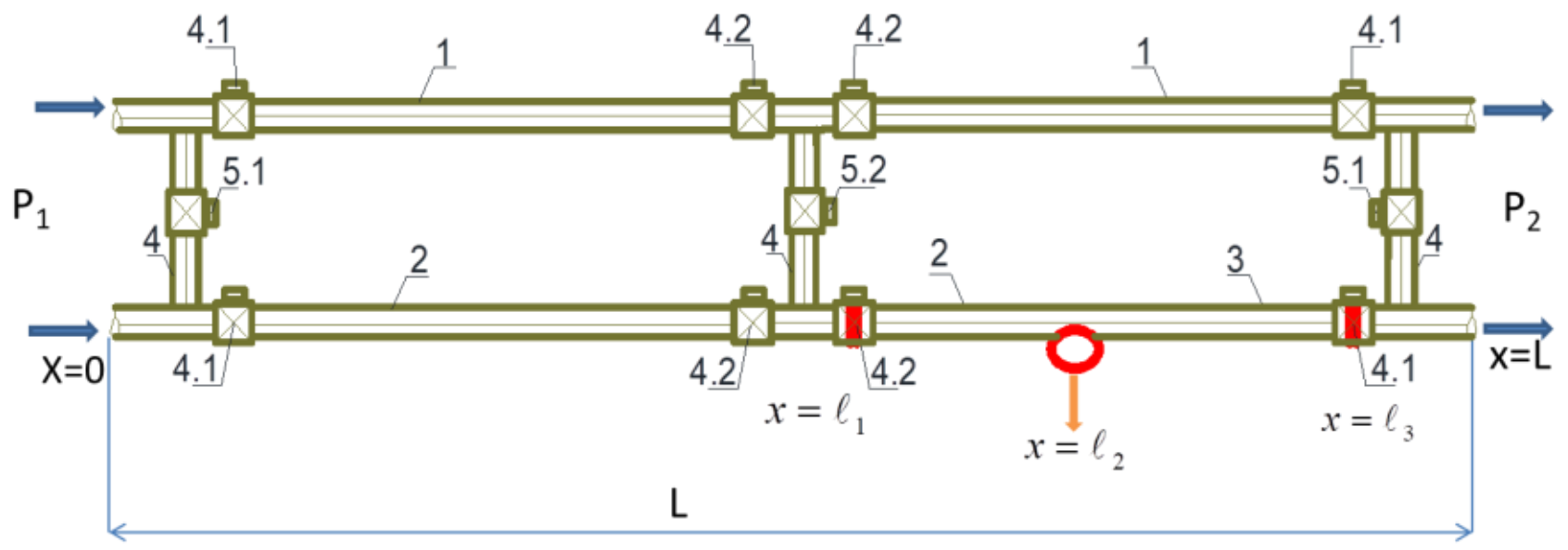

**Scheme 4.4.1. Principle diagram of a complex gas pipeline during the reconstruction phase.**

**1 - Undamaged pipeline ($0 \leq x \leq L$);**
**2 - Section from the starting point to the leakage point of the damaged pipeline ($0 \leq x \leq \ell_2$),**
**3 - Section from the leakage point to the end point of the damaged pipeline ($\ell_2 \leq x \leq L$);**
**4- Connectors;**
**4.1 and 4.2 - Automatic valves on the pipeline;**
**5.1 and 5.2 - Automatic valves on the connectors;**
**P1 - Pressure at the starting section of the gas pipeline, Pa;**
**P2 - Pressure at the end section of the gas pipeline, Pa;**
**L - Length of the gas pipeline, m.**
$\ell_1, \ell_3$ **– Distance from the starting point of the damaged pipeline to the nearest shut-off fittings to the left and right of the leakage point, respectively.**

In the proposed gas transportation scheme, at the fixed moment $t=t_1$, technological measures are implemented to isolate the damaged section from the main part of the pipeline. Specifically, the location of the damaged section is detected in the shortest possible time ($\ell_2$), and the positions of the valves installed to the right ($\ell_1$)and left ($\ell_3$) of the leakage point are identified and closed by the dispatcher center at that moment.

Based on the principle of managing technological measures, it is recommended to install automatic valves at the starting section of the pipeline network to the right of the first connector (4.1) and to the left of the last connector (4.2).

In the normal (stationary) operating regime of gas pipelines, the valves on the connectors are closed, while the valves on the pipelines are open. In emergency situations (when the pipeline's integrity is compromised), the valves to the right and left of the leakage point on the pipeline are closed ( $\ell_1, \ell_3$), isolating the damaged section from the main part of the gas pipeline. To ensure an adequate gas supply to consumers, the valves (5.1 and 5.2) located to the right and left of the leakage point on the connecting pipeline are opened under the control of dispatcher center management panels.

The main distinguishing feature of the proposed gas pipeline compared to the existing one is its operation in an independent hydraulic regime. Another key difference of the reconstructed gas pipeline is the independence of the pipeline dynamics from each other during emergencies and technological regimes. The proposed parallel gas pipeline is a complex system characterized by a compressor station serving each pipeline at the starting section and numerous elements such as connectors and valves.

***Establishing an information database for the efficient management of the reconstructed gas pipeline system.***

To establish an information database for the efficient management of the reconstructed gas pipeline system, the reporting scheme was implemented in the following sequence. First, let's assume that the integrity of the second pipeline in the parallel gas pipeline has been compromised, resulting in a large amount of gas leaking from the damaged section ( $x = \ell_2$) into the environment (Scheme 4.4.1).

At this time, the undamaged pipeline continues to operate in its previous stationary regime. In other words, the non-stationary regime that arises in the damaged pipeline does not affect the undamaged pipeline. The non-stationary regime of gas flow in the damaged gas pipeline due to the accident is divided into three different technological sections:

The section from the starting point of the pipeline to the shut-off valve to the left of the leakage point ( $0 \le x \le \ell_1$),

The damaged section of the pipeline, between the shut-off valves to the left and right of the leakage point. ( $\ell_1 \le x \le \ell_3$),

The section from the shut-off valve to the right of the leakage point to the end point of the pipeline ( $\ell_3 \le x \le L$ ).

The hydraulic models of all three sections are different. The hydraulic model of the pipeline system is constructed from the models of its components: multi-line pipelines with loops and connectors, shut-off valves, compressor stations, etc. The models of these components, in turn, consist of the computational expressions for the elements. The solution to this problem necessitates considering the entire system as a whole rather than just a separate section of the pipeline. Initially, there is a significant need for analytical expressions to solve the equations describing non-stationary events in parallel gas pipelines and for mathematical expressions to solve these equations that reflect real events in the pipeline. The non-stationary operating conditions of pipelines lead to significant pressure changes and disruptions in the normal operation of gas delivery to consumers. Studying these processes to efficiently manage the pipelines and considering the results of the obtained analytical solution methods at dispatcher control points will significantly enhance the reliability of gas supply and the efficiency of transportation processes.

For this purpose, it is essential to obtain the analytic expressions reflecting the non-stationary events in the two-line parallel gas pipelines and ensure the management of physical processes through the solution of these equations. The mentioned analytic expressions have been determined based on the resolution of the problem addressed in our published work [88], which specifically deals with the non-stationary flow of gas in each of the three sections. These expressions are as follows:

$$P_1(x,t)=P_1-2aG_0x+\frac{8aLG_0}{\pi^2}\sum_{n=1}^{\infty}Cos\frac{\pi nx}{L}\frac{e^{-(2n-1)^2\alpha_3 t}}{(2n-1)^2}-\frac{4aL}{\pi^2}G_0\sum_{n=1}^{\infty}((1-(-1)^n)Cos\frac{\pi nx}{L}*$$
$$*\frac{e^{-n^2\alpha_3 t}}{n^2}-2aG_{ut}\left(\frac{x^2}{2L}+\frac{\ell_2^2}{2L}+\frac{L}{3}-\ell_2\right)+\frac{4aL}{\pi^2}G_{ut}\sum_{n=1}^{\infty}Cos\frac{\pi nx}{L}Cos\frac{\pi n\ell_2}{L}\frac{e^{-n^2\alpha_3 t}}{n^2} \quad (1)$$

$$P_{2,3}(x,t)=P_1(x,t)-2aG_{ut}(x-\ell_2) \quad (2)$$

Here, $\alpha_3=\dfrac{\pi^2c^2}{2aL^2}$

By considering x=0 in equation (1) and x=L in equation (2), we establish the conformity of the pressure distribution at the beginning and end points of the gas pipeline over time.

$$P_1(0,t)=P_1-\frac{c^2t}{L}G_{ut}+\frac{8aLG_0}{\pi^2}\sum_{n=1}^{\infty}\frac{e^{-(2n-1)^2\alpha_3 t}}{(2n-1)^2}-2aG_{ut}\left(\frac{\ell_2^2}{2L}+\frac{L}{3}-\ell_2\right)+$$
$$+\frac{4aL}{\pi^2}G_{ut}\sum_{n=1}^{\infty}Cos\frac{\pi n\ell_2}{L}\frac{e^{-\alpha_2 t}}{n^2} \quad (3)$$

$$P_3(L,t)=P_2-\frac{c^2t}{L}G_{ut}+\frac{8aLG_0}{\pi^2}\sum_{n=1}^{\infty}(-1)^n\frac{e^{-(2n-1)^2\alpha_3 t}}{(2n-1)^2}-2aG_{ut}\left(\frac{L}{2}+\frac{\ell_2^2}{2L}+\frac{L}{3}-\ell_2\right)+$$
$$+\frac{4aL}{\pi^2}G_{ut}\sum_{n=1}^{\infty}(-1)^nCos\frac{\pi n\ell_2}{L}\frac{e^{-\alpha_2 t}}{n^2}-2aG_{ut}(\ell_2-L) \quad (4)$$

Based on equations (3) and (4) and the results of the analysis conducted in [1], we obtain the following expression:

$$\theta=\frac{1}{2}+\left[\frac{1-p(t)}{1+p(t)}\right]\cdot\left[\frac{2}{3}+\frac{e^{-2\alpha_2 t}-4e^{-\alpha_2 t}}{\pi^2}\right] \quad (5)$$

Here, $\alpha_2=\dfrac{\pi^2c^2}{2aL^2}$ ; p(t)= $\dfrac{P_1-P_1(0,t)}{P_2-P_2(L,t)}$; $\theta=\dfrac{\ell_2}{L}$

$G_{ut}$- Mass flow rate at the leakage point for the gas pipeline accident regime, Pa·sec/m
$G_0$- Mass flow rate at the beginning of the gas pipeline in the stationary regime, Pa·sec/m
c - Speed of sound propagation for gas in an isothermal process, m/sec
P - Pressure at the ends of the gas pipeline, Pa

x - Coordinate along the length of the gas pipeline, m
t - Time coordinate, sec
2a - Charney linearization, 1/sec

It should be noted that the function p(t) characterizes the ratio of pressure values at the beginning and end points of the gas pipeline. Based on the analysis of the parameters involved in equation (5), it becomes evident that there is a functional dependence between θ and p(t). Based on the conducted research, it was determined that if p(t) = 1, then θ takes a value of 0.5. This indicates that the gas leakage occurs in the exact middle section of the pipeline. When p(t) takes values greater than 1, θ takes values less than 0.5, confirming that the gas leakage occurs between the initial and middle sections of the pipeline. Conversely, when p(t) takes values less than 1, θ takes values greater than 0.5, confirming that the gas leakage is situated between the middle and final sections of the pipeline.

To confirm the above ideas, the expressions (3) and (4) were examined, and the following parameters of the gas pipeline were considered to determine the law of pressure change along the gas axis at different characteristic points of gas leakage location (at the beginning, $\ell_2 = 0{,}5\times 10^4$ m ; in the middle, $\ell_2 = 5\times 10^4$ m, and at the end, $\ell_2 = 9{,}5\times 10^4$ m).

$P_1 = 55\times 10^4$ Pa: $P_2 = 25\times 10^4$ Pa: $G_0$= 30 Pa×sec/m ; $c = 383{,}3\dfrac{m}{\text{sec}}$ ; L = $10\times 10^4$ m:

Currently, the calculation of gas flow regimes is implemented through the use of modern software. Firstly, the pressure values of the gas pipeline at the beginning and end of the pipeline are calculated for each of the three characteristic points of the gas leakage every 100 seconds. The results are recorded in Table 4.4.1 below.

| **t, sec** | **The location of the gas pipeline failure.** | | | | | |
|---|---|---|---|---|---|---|
| | $\ell_2 = 0{,}5\cdot 10^4$ m | | $\ell_2 = 5\cdot 10^4$ m | | $\ell_2 = 9{,}5\cdot 10^4$ m | |
| | The pressure in the starting part, $P_1(0,t)$, $10^4$ Pa | Pressure in the last part, $P_3(L,t)$, $10^4$ Pa | The pressure in the starting part, $P_1(0,t)$, $10^4$ Pa | Pressure in the last part, $P_3(L,t)$, $10^4$ Pa | The pressure in the starting part, $P_1(0,t)$, $10^4$ Pa | Pressure in the last part, $P_3(L,t)$, $10^4$ Pa |
| 100 | 52,23 | 25 | 55 | 24,99 | 55 | 22,23 |
| 200 | 50,58 | 25 | 54,9 | 24,89 | 55 | 20,58 |
| 300 | 49,3 | 24,99 | 54,66 | 24,66 | 54,99 | 19,3 |
| 400 | 48,21 | 24,97 | 54,34 | 24,33 | 54,97 | 18,21 |
| 500 | 47,25 | 24,93 | 53,96 | 23,96 | 54,925 | 17,25 |
| 600 | 46,38 | 24,84 | 53,56 | 23,49 | 54,84 | 16,38 |
| 700 | 45,58 | 24,72 | 53,14 | 23,14 | 54,72 | 15,58 |
| 800 | 44,83 | 24,57 | 52,71 | 22,71 | 54,565 | 14,84 |
| 900 | 44,13 | 24,37 | 52,28 | 22,27 | 54,37 | 14,14 |

**Table 4.4.1. Pressure values at the point of rupture of the gas pipeline depending on time.**

Using the data provided in Table 1, we plot the graph of p(t) as a function of time for each of the three characteristic points.(Graph 4.4.1)

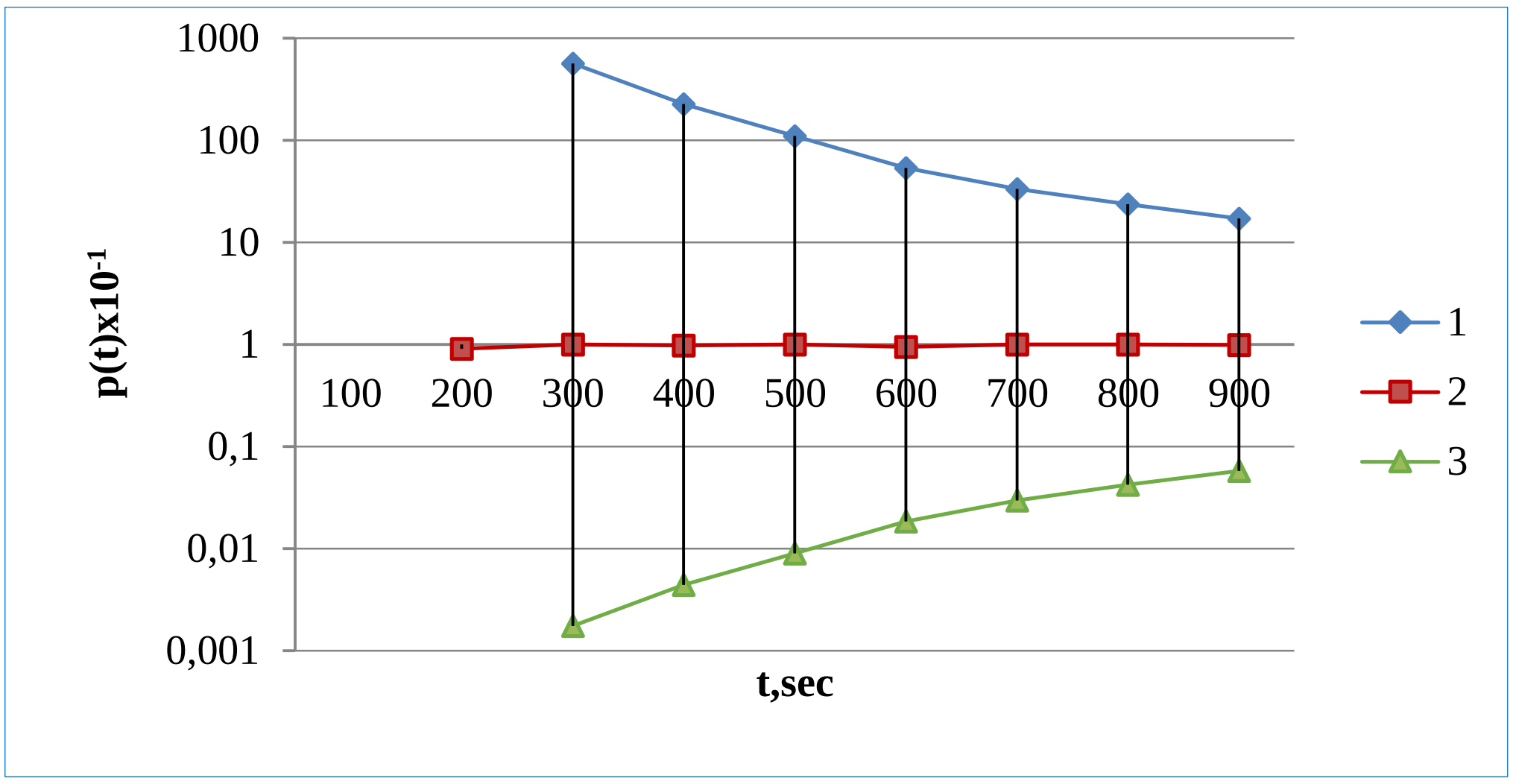

**Graph 4.4.1. Variation pattern of p(t)p(t)p(t) function over time, dependent on the location of the gas pipeline rupture. (1- $\ell_2$ =0,5×10$^4$ m; 2- $\ell_2$ =5×10$^4$ m; 3- $\ell_2$ =9,5×10$^4$ m).**

The analysis of Graph 4.4.1 reveals that, regardless of the location of the gas pipeline rupture based on the parameters mentioned above, the p(t) function attains maximum or minimum values within a duration of t=300 seconds. Therefore, the fixed time observed at dispatcher stations for the examined gas pipeline is $t_1$= 300 seconds.

Minimizing the transmission time of information to dispatcher control centers is crucial. Therefore, the adopted value of $t_1$= 300 seconds ensures the safe operation of gas pipelines. On the other hand, the p(t) function takes different values depending on the dependence on θ, which determines the location of the leakage. Specifically, when θ > 0.5, the p(t) function takes its minimum value (3rd line in Graph 1), while for θ < 0.5, the p(t) function takes its maximum value (1st line in Graph 1). However, when θ is approximately equal to 0.5, p(t) approaches 1, indicating a scenario where the leakage point is in the middle section of the pipeline ($\ell_2$ =5·10$^4$m).

The analysis of the graph leads to the following conclusion: when p(t) belongs to the interval ]0;1[, it characterizes the location of gas pipeline rupture occurrence between $0 < \ell_2 < \frac{L}{2}$. When p(t) is greater than 1, it confirms the occurrence of gas pipeline rupture between $\frac{L}{2} < \ell_2 < L$. When p(t) = 1, it confirms that the location of gas pipeline rupture is at the exact midpoint ($\ell_2 = \frac{L}{2}$). Clearly, these simple and understandable pieces of information are very important for the gas emergency control center.

Based on the analysis conducted, it was determined that the fixation period $t=t_1$ remains constant regardless of the location of accidents on the pipeline. To determine whether this time remains constant or varies depending on the parameters of the pipeline, we consider different parallel gas pipelines with the following parameters particular importance.

$$P_1 = 14 \times 10^4 \text{ Pa}: \ P_2 = 11\times 10^4 \text{ Pa}: \ G_0 = 10 \text{ Pa} \times \text{sec/m} ; \ c = 383{,}3\frac{m}{\text{sec}} ; \ L = 3\times 10^4 \text{ m}:$$

It is evident from the given information that the length of the newly adopted parallel gas pipeline is at least 3 times shorter than the length of the previously analyzed gas pipeline, and also that the values of other parameters (pressure and consumption) are different. We consider the length of the examined complex parallel pipeline at three characteristic points where the gas leakage occurs, similar to the previous pipeline: near the beginning of the pipeline $\ell_2$ = 0,5×10$^4$ m, at the middle $\ell_2$ =1,5×10$^4$ m, and near the end $\ell_2$ = 2,5× 10$^4$ m.

We use equations (3) and (4) to calculate the pressure values of the gas pipeline's beginning and end sections every 60 seconds for each of the three characteristic points where the gas leakage occurs. The results are recorded in the following table (Table 4.4.2).

| t, sec | The location of the gas pipeline failure. | | | | | |
|---|---|---|---|---|---|---|
| | $\ell_2 = 0,5 \cdot 10^4$ m | | $\ell_2 = 1,5 \cdot 10^4$ m | | $\ell_2 = 2,5 \cdot 10^4$ m | |
| | The pressure in the starting part, $P_1(0,t)$, $10^4$ Pa | The pressure in the last part, $P_2(L,t)$, $10^4$ Pa | The pressure in the starting part, $P_1(0,t)$, $10^4$ Pa | The pressure in the last part, $P_2(L,t)$, $10^4$ Pa | The pressure in the starting part, $P_1(0,t)$, $10^4$ Pa | The pressure in the last part, $P_2(L,t)$, $10^4$ Pa |
| 60 | 13,37 | 10,97 | 13,83 | 10,83 | 13,97 | 10,37 |
| 120 | 12,95 | 10,79 | 13,54 | 10,54 | 13,79 | 9,95 |
| 180 | 12,61 | 10,55 | 13,24 | 10,24 | 13,55 | 9,61 |
| 240 | 12,29 | 10,27 | 12,95 | 9,95 | 13,27 | 9,29 |
| 300 | 11,99 | 9,98 | 12,66 | 9,66 | 12,98 | 8,99 |
| 360 | 11,70 | 9,69 | 12,36 | 9,36 | 12,69 | 8,70 |
| 420 | 11,40 | 9,40 | 12,07 | 9,07 | 12,40 | 8,40 |
| 480 | 11,11 | 9,11 | 11,77 | 8,77 | 12,11 | 8,11 |
| 540 | 10,81 | 8,81 | 11,48 | 8,48 | 11,81 | 7,81 |
| 600 | 10,52 | 8,52 | 11,19 | 8,19 | 11,52 | 7,52 |

**Table 2. Pressure values of the gas pipeline depending on the location of the leakage and time.**

Using the data from Table 4.4.2, we plot the graph of p(t) as a function of time for each of the three characteristic points. (Graph 4.4.2).

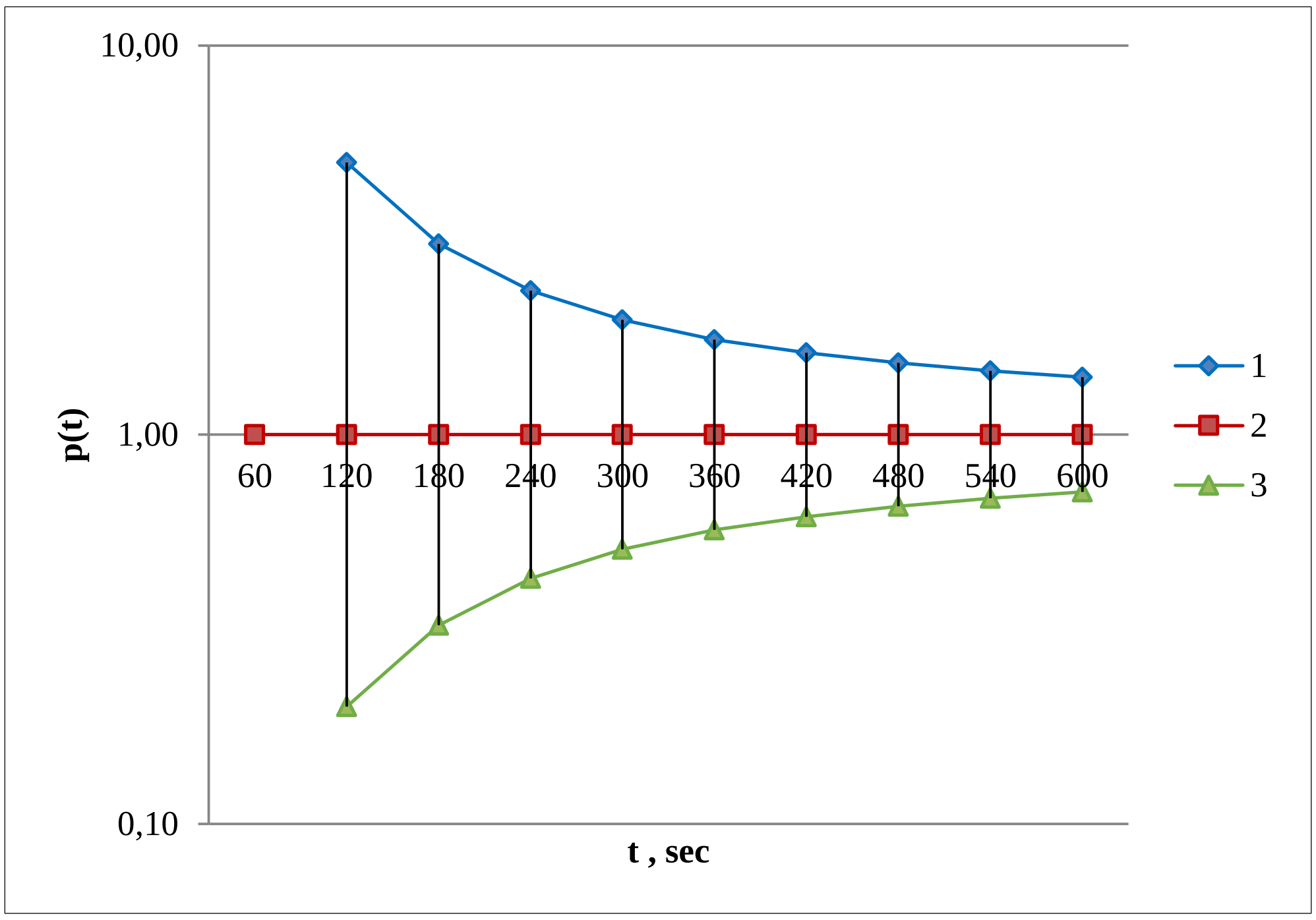

**Graph 4.4.2. The variation pattern of p(t)p(t)p(t) function depending on the time and the location of the gas pipeline's leakage. (1- $\ell_2 = 0,5\times10^4$ m; 2- $\ell_2 = 1,5\times10^4$ m; 3- $\ell_2 = 2,5\times10^4$ m).**

From the analysis of Graph 2, it can be noted that the variation of pressure over time in the outer sections depending on the location of the gas pipeline's leakage is similar to that observed in Graph 4.4.1. Specifically, when θ>0.5, the p(t) function takes its minimum value (indicated by the 3rd line in Graph 4.4.2), while for θ<0.5, the p(t) function takes its maximum value (indicated by the 1st line in Graph 4.4.2).

However, regardless of the location of the gas pipeline's leakage, the p(t) function takes its maximum or minimum value over a duration of t=120 seconds. This means that the fixed duration $t_1$=120 seconds is applicable for this gas pipeline. It can be concluded that the fixed duration varies depending on the length of the gas pipeline and the parameters of the gas flow. Thus, as the length of the pipeline segment decreases, this duration decreases accordingly. This is because the time required for the parameters of the non-stationary gas flow during the accident to reach the starting point of the gas pipeline is related to the change in the parameters. Considering the short length of the analyzed gas pipeline segment and the reliability of the information transmitted to the dispatcher station, it is necessary to accept a minimum duration of fixation that is sufficient to record the change in pressure of the non-stationary gas flow at the dispatcher station.

In other words, even if the speed of the non-stationary gas flow in the pipeline matches its sound propagation speed, the delivery time of information from the last three points of a 100,000-meter pipeline segment to the starting point cannot be faster than 261 seconds. It should be noted that considering this indicator during the reconstruction stage of gas pipelines for efficient management allows for the optimization of the issue of minimizing the transfer time of information and control system data to dispatcher control centers.

We conduct theoretical research to examine the consistency of the method for determining the fixing time $t=t_1$, which is one of the measures taken to reduce the workload of the dispatcher station. For this purpose, we perform engineering calculations based on theoretical research. Using equation (5), we determine the location of the pipeline's instability point at each of the three characteristic points, depending on time. First, we calculate based on the data of the gas pipelines under consideration, with a length of L=100,000 m, and record the calculation results in the following table (Table 4.4.3).

| t, sec | Actual values of failure of the pipeline with a length of 100000 m. | | | | | |
|---|---|---|---|---|---|---|
| | $\ell_2 = 0{,}5 \cdot 10^4$ m | | $\ell_2 = 5 \cdot 10^4$ m | | $\ell_2 = 9{,}5 \cdot 10^4$ m | |
| | Calculated by formula (5) price, | Relative error of actual and calculated prices | Calculated by formula (5) price, | Relative error of actual and calculated prices | Calculated by formula (5) price, | Relative error of actual and calculated prices |
| 100 | - | - | 5,00 | 0,00 | 5,00 | 0,47 |
| 200 | - | - | 5,20 | 0,02 | 9,20 | 0,03 |
| 300 | 0,55 | 0,04 | 5,00 | 0,00 | 9,45 | 0,01 |
| 400 | 0,33 | 3,44 | 5,04 | 0,00 | 9,67 | 0,02 |
| 500 | 0,15 | 1,61 | 5,00 | 0,00 | 9,84 | 0,04 |
| 600 | 0,04 | 1,38 | 5,12 | 0,01 | 9,96 | 0,05 |

**Table 4.4.3. Relative Error of the Actual and Calculated Values of the Gas Instability Point Location Depending on Time for a Pipeline with a Length of 100,000 meter:**

Sequentially, based on the given data for the next gas pipeline with a length of L = 30,000 m, we calculate the locations of instability points at each of the three characteristic points using equation (5) and record the results in the following table (Table 4.4.4).

| t, sec | Actual values of failure of the pipeline with a length of 30000 m. | | | | | |
|---|---|---|---|---|---|---|
| | Actual values of pipeline failure. | | | | | |
| | $\ell_2 = 0,5 \cdot 10^4$ m, | | $\ell_2 = 1,5 \cdot 10^4$ m | | $\ell_2 = 2,5 \cdot 10^4$ m | |
| | Calculated by formula (5) price, | Relative error of actual and calculated prices | Calculated by formula (5) price, | Relative error of actual and calculated prices | Calculated by formula (5) price, | Relative error of actual and calculated prices |
| 60 | 0,06 | 0,87 | 1,50 | 0,0 | 2,94 | 0,175 |
| 120 | 0,28 | 0,44 | 1,50 | 0,0 | 2,72 | 0,088 |
| 180 | 0,51 | 0,02 | 1,50 | 0,0 | 2,49 | 0,005 |
| 240 | 0,71 | 0,41 | 1,50 | 0,0 | 2,29 | 0,083 |
| 300 | 0,85 | 0,70 | 1,50 | 0,0 | 2,15 | 0,140 |
| 360 | 0,95 | 0,91 | 1,50 | 0,0 | 2,05 | 0,181 |

**Table 4.4.4. The values of the relative error of the actual and calculated positions of instability of the gas pipeline with a length of 30,000 m, depending on time.**

Based on the analysis of Tables 4.4.3 and 4.4.4, it can be noted that the prescribed equation (5) is useful for engineering calculations. The results of both analyzed gas pipelines' calculations showed that only at the fixed time $t=t_1$ is the location of instability accurately determined, with its relative error value not exceeding 0.04, or in other words, less than 4%.

Therefore, the prescribed calculation scheme is a suitable method for verifying the compliance of gas pipelines with the requirements for safe and efficient utilization, determining the parameters of gas flow in non-stationary regimes, as well as selecting appropriate parameters for the reconstruction of gas pipelines.

***The determination of the analytical expression based on which accidents and technological regimes are selected for dispatcher stations.***

The variation in pressure in a gas pipeline can result from both changes in the amount of gas supplied to consumers, the opening and closing of valves in the pipeline, and the operation of various elements within the pipeline. As a result, technological processes cause changes in pressure at the beginning and end sections of the gas pipeline. One of the most important tasks of operational control in gas pipeline technological regimes is to identify critical conditions in gas supply systems at dispatcher stations. For instance, pressure variations at the end sections of the pipeline can occur as a result of changes in technological regimes. Therefore, it is essential to prioritize understanding the pressure changes resulting from technological processes in the event of accidents.

The initial model of events and processes occurring in gas transmission systems is not always accurately represented with high precision. Therefore, to solve the problems of operational management of gas transmission systems, it is necessary to use analytical expressions to assess the real situation of gas transmission systems. Having a mathematical expression allows for selecting the parameters and control structures of the pipeline, determining optimality criteria and constraints, assessing accuracy, selecting appropriate technical control devices, and so on.

To solve the problems of managing trunk gas pipelines, it is necessary to understand the non-stationary dynamic characteristics of these systems. It is appropriate to utilize the formalization of technological and accident processes for the transportation of gas in pipeline systems. By employing the method of describing accident processes, it is possible to obtain unit characteristics. Based on theoretical research, the following inequality has been determined:.

$$\frac{\varphi-0,5}{\varphi+0,5} < p(t) < \frac{\varphi+0,5}{\varphi-0,5} \quad (6)$$

Here, $\varphi = \frac{2}{3} + \frac{e^{-2\alpha_2 t} - 4e^{-2\alpha_2 t}}{\pi^2}$

It has been determined that if the function p(t) satisfies the conditions of the inequality mentioned above, the variation of pressure over time corresponds to the lawfulness characterizing accidents in the gas pipeline, otherwise, it is the result of technological processes. Therefore, the function p(t) is a reliable parameter that enables dispatchers to make timely decisions for the efficient management of the gas pipeline.

At the point where time $t=t_1$ is fixed, the value of θ\thetaθ and the obtained values of the function p(t) are reliable information for the dispatcher station. This is because the accuracy of determining the leakage location is characterized by the moment $t=t_1$. To visualize the methodology, we use the flow diagram below.

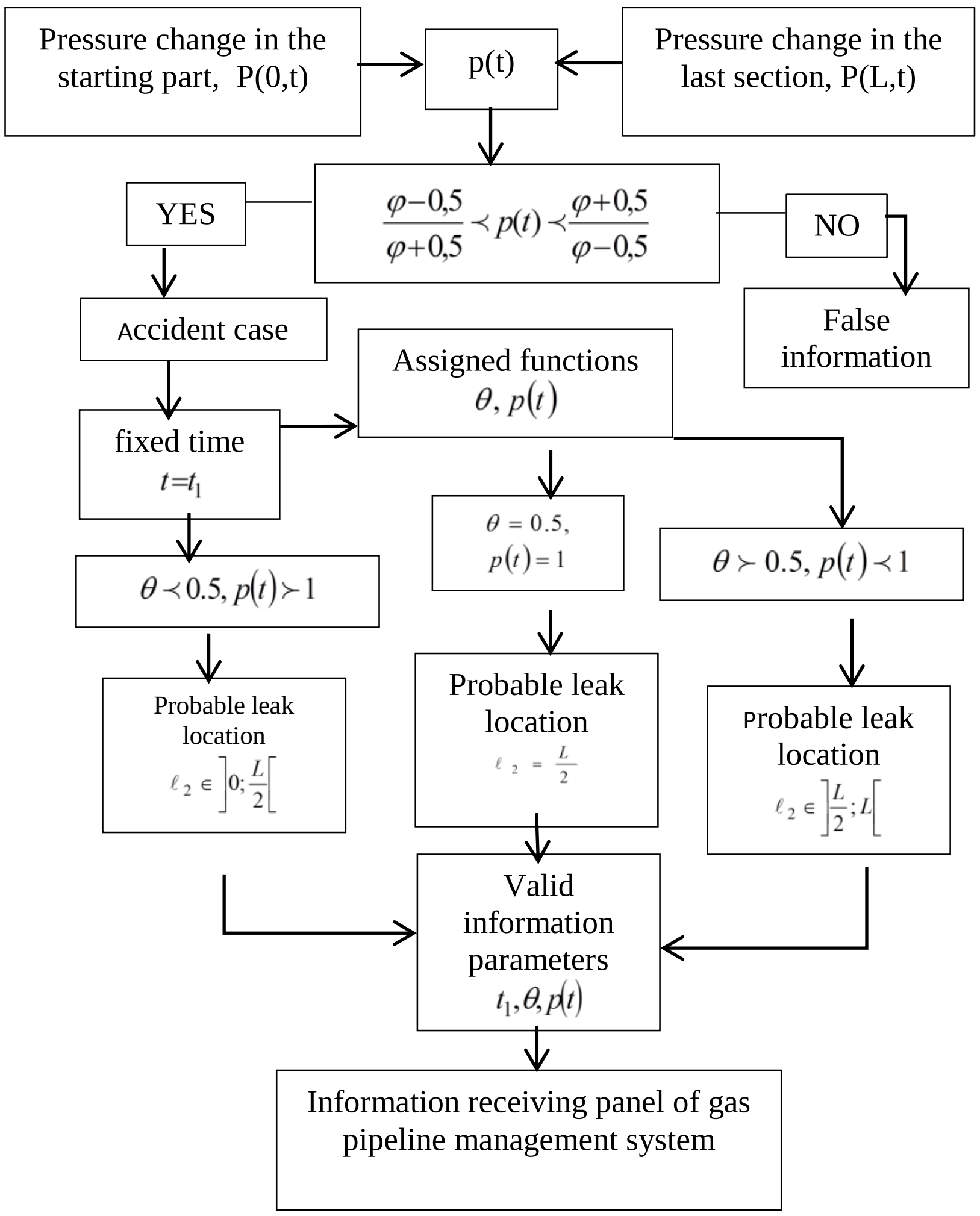

**Figure 4.4.5. A visual diagram of the methodology.**

On the other hand, for the reliability and effectiveness of automatic detection of the accident location in online mode, these parameters are useful. For instance, if pressure changes are observed in the initial and final sections of the gas pipeline, the dispatcher station first determines whether this change occurred as a result of the accident based on the conditions of inequality (6). Sequentially, the values of the function p(t) are automatically calculated, and its maximum or minimum value is determined.

The obtained values are assumed to be equal to $t=t_1$. Based on the value of p(t) at $t=t_1$ calculated according to the inequality (5), the value of θ is determined. If this value is less than 0.5 and p(t)>1, it is determined that the location of the accident is between the starting and middle sections of the gas pipeline, and based on the value of θ, the location of the gas leak ( $\theta = \frac{\ell_2}{L}$ ) is automatically identified. Also, the valves ( $\ell_1, \ell_3$ ), whose locations are already known, are activated to isolate the damaged section of the pipeline from the main section of the gas pipeline.

Based on the results of the research, it can be stated that using expressions (5) and (6), it is possible to distinguish between accidental and technological regimes. On the other hand, from the moment the occurrence of an accident is confirmed, it is possible to determine the location of gas leaks and the sequence of valve activations. Thus, it can be concluded that the values obtained at the fixed time $t=t_1$ and from these expressions provide a reliable information base for monitoring and efficiently managing the technological regimes of gas distribution pipelines under extreme conditions.

To ensure the comprehensive management of system integrity, it is essential that emergency dispatch center personnel accurately understand whether pressure fluctuations in the final segments of the gas pipeline are due to incidents or technological processes and can discern the relationships between these processes. In this context, modern control methods, particularly those enhanced by artificial intelligence (AI) technologies, should be prioritized. For instance, AI technologies, particularly hybrid models, provide a more resilient and effective approach by offering accurate predictions even with limited data.

In this study, we compare the performance of algorithms by presenting tables that highlight the accuracy, efficiency, and other relevant metrics of each algorithm (Table 4.4.3, Table 4.4.4). Specifically, indicators that assess algorithms' capabilities in real-time incident detection and data augmentation for handling data insufficiencies have been presented. This approach lays a stronger foundation for addressing existing challenges in pipeline engineering. Furthermore, key tasks include tracking the evolution of incidents, using data augmentation to address data gaps, enhancing model interpretability, and improving the management of pipeline integrity. Future research should focus on the optimization and integration of AI models to support informed decision-making based on a reliable information base. Additionally, AI enables the analysis of large volumes of data from various sources to predict incident risks, optimize management processes, and provide critical insights that support engineering decisions. Consequently, prioritizing the broader application of hybrid models in real-time monitoring and predictive technologies is essential for advancing this field.

## Conclusion

A highly efficient technological scheme with a high reliability indicator has been proposed, based on the principles of the reconstruction phase, for the management of technological processes in gas transportation. The application of modern technologies and equipment to the proposed gas transportation scheme is of particular importance for optimizing the gasification capabilities of consumers.

Several mathematical models have been developed and an effective computational scheme has been devised to address a number of important new problems related to the management and regulation of technological processes in gas transportation. This will significantly enhance the reliability and efficiency of gas supply management during the reconstruction phase of complex gas pipelines. The calculation results showed that the error in comparing the obtained algorithm values does not exceed 4%. The developed computational scheme can be used to create an information database for the efficient management of gas transportation systems.

To ensure the unconditional resolution of efficient management issues, a function characterized by the ratio of pressure drops has been defined and tested. The analysis of the p(t) function and parameters, which can perform multiple operations at the fixed time $t=t_1$ , demonstrated that considering this indicator during the reconstruction phase can optimize the minimization of the data transmission time from information and control systems to dispatcher control points. The established relationship enables the differentiation of technological (stationary) and non-stationary operating modes of gas pipelines based on the indicators of changes in the parameters of the gas flow's final sections.

Based on the results of the analyses, individual problems related to the reconstruction of parallel main gas pipelines have been solved using the developed algorithms and methods, taking into account the non-stationary regime. Consequently, targeted methods for the reconstruction processes of an efficiently managed gas transport system have been identified as the subject of the research.

A technological scheme has been proposed that adopts the integration of new technologies and methodologies to ensure efficiency, safety, and reliability for the effective management of the reconstructed gas pipeline system. A method for the rapid detection of faults, which should become an integral part of the automated dispatch control system of the gas pipeline, is proposed. This method involves recording the gas flow pressure at the start and end of the pipeline at a fixed moment $t=t_1$.

The applied optimization algorithms not only enhance system reliability but also optimize interventions during emergencies. Future research could further advance by exploring the application of these methods in other engineering fields and integrating new technologies. This work offers new perspectives for improving the efficiency of gas pipeline systems and lays the groundwork for future investigations.

The proposed technological model provides a comprehensive solution for managing complex gas transportation systems. Future work should focus on further improving the model by integrating additional predictive technologies and enhancing data transmission methods. Furthermore, continuous advancements in control systems and emergency response strategies are necessary to address the increasing complexity of global gas transportation networks.